\documentclass[english]{smfbook}

\usepackage{amsmath,amscd,amssymb}

\author{Takuro Mochizuki}
\address{Research Institute for Mathematical Sciences,
Kyoto University, Kyoto 606-8502, Japan}
\email{takuro@kurims.kyoto-u.ac.jp}
\title{Mixed twistor $D$-modules}
\alttitle{}

\newcommand{\nbiga}{\mathcal{A}}
\newcommand{\nbigb}{\mathcal{B}}
\newcommand{\nbigc}{\mathcal{C}}
\newcommand{\nbigd}{\mathcal{D}}
\newcommand{\nbige}{\mathcal{E}}
\newcommand{\nbigf}{\mathcal{F}}
\newcommand{\nbigg}{\mathcal{G}}
\newcommand{\nbigh}{\mathcal{H}}
\newcommand{\nbigi}{\mathcal{I}}
\newcommand{\nbigj}{\mathcal{J}}
\newcommand{\nbigk}{\mathcal{K}}
\newcommand{\nbigl}{\mathcal{L}}
\newcommand{\nbigm}{\mathcal{M}}
\newcommand{\nbign}{\mathcal{N}}
\newcommand{\nbigo}{\mathcal{O}}
\newcommand{\nbigp}{\mathcal{P}}
\newcommand{\nbigq}{\mathcal{Q}}
\newcommand{\nbigr}{\mathcal{R}}
\newcommand{\nbigs}{\mathcal{S}}
\newcommand{\nbigt}{\mathcal{T}}
\newcommand{\nbigu}{\mathcal{U}}
\newcommand{\nbigv}{\mathcal{V}}
\newcommand{\nbigw}{\mathcal{W}}
\newcommand{\nbigx}{\mathcal{X}}
\newcommand{\nbigy}{\mathcal{Y}}
\newcommand{\nbigz}{\mathcal{Z}}

\newcommand{\proj}{\mathbb{P}}
\newcommand{\seisuu}{{\mathbb Z}}
\newcommand{\rnum}{{\mathbb Q}}

\newcommand{\cnum}{{\mathbb C}}
\newcommand{\real}{{\mathbb R}}

\newcommand{\Tate}{\mathbb{T}}

\newcommand{\DD}{\mathbb{D}}

\newcommand{\II}{\mathbb{I}}

\newcommand{\gbigc}{\mathfrak C}
\newcommand{\gbigd}{\mathfrak D}

\newcommand{\gbigf}{\mathfrak F}
\newcommand{\gbigg}{\mathfrak G}

\newcommand{\gbigi}{\mathfrak I}

\newcommand{\gbigp}{\mathfrak P}

\newcommand{\gbigs}{\mathfrak S}
\newcommand{\gbigt}{\mathfrak T}

\newcommand{\gminia}{\mathfrak a}
\newcommand{\gminib}{\mathfrak b}
\newcommand{\gminic}{\mathfrak c}

\newcommand{\gminie}{\mathfrak e}

\newcommand{\gminik}{\mathfrak k}

\newcommand{\gminim}{\mathfrak m}

\newcommand{\gminip}{\mathfrak p}

\newcommand{\vecxi}{{\pmb \xi}}

\newcommand{\vece}{{\pmb e}}

\newcommand{\vecu}{{\pmb u}}
\newcommand{\vecw}{{\pmb w}}

\newcommand{\veczero}{{\pmb 0}}

\newcommand{\veca}{{\pmb a}}
\newcommand{\vecb}{{\pmb b}}
\newcommand{\vecbeta}{{\pmb \beta}}
\newcommand{\vecdelta}{{\pmb \delta}}

\newcommand{\vecc}{{\pmb c}}
\newcommand{\vecd}{{\pmb d}}

\newcommand{\veck}{{\pmb k}}
\newcommand{\vecm}{{\pmb m}}

\newcommand{\vecN}{{\pmb N}}
\newcommand{\vecI}{{\pmb I}}

\newcommand{\vecx}{{\pmb x}}

\newcommand{\vecF}{{\pmb F}}

\newcommand{\vecp}{{\pmb p}}
\newcommand{\vecq}{{\pmb q}}

\newcommand{\vecz}{{\pmb z}}
\newcommand{\vecC}{{\pmb C}}

\newcommand{\vecS}{\pmb{S}}
\newcommand{\vecT}{{\pmb T}}
\newcommand{\vecL}{{\pmb L}}
\newcommand{\vecvarphi}{{\pmb \varphi}}
\newcommand{\vecSigma}{{\pmb{\Sigma}}}
\newcommand{\vecX}{{\pmb{X}}}

\newcommand{\larr}{\leftarrow}
\newcommand{\llarr}{\longleftarrow}

\newcommand{\rarr}{\rightarrow}
\newcommand{\lrarr}{\longrightarrow}
\newcommand{\darr}{\downarrow}

\newcommand{\pf}{{\bf Proof}\hspace{.1in}}

\def\Hom{\mathop{\rm Hom}\nolimits}

\def\End{\mathop{\rm End}\nolimits}
\def\Ext{\mathop{\rm Ext}\nolimits}

\def\Cok{\mathop{\rm Cok}\nolimits}

\def\Image{\mathop{\rm Im}\nolimits}

\def\Re{\mathop{\rm Re}\nolimits}

\def\Gr{\mathop{\rm Gr}\nolimits}

\def\Tot{\mathop{\rm Tot}\nolimits}

\def\rank{\mathop{\rm rank}\nolimits}

\def\Ker{\mathop{\rm Ker}\nolimits}

\def\Gr{\mathop{\rm Gr}\nolimits}

\def\Res{\mathop{\rm Res}\nolimits}

\def\ord{\mathop{\rm ord}\nolimits}

\def\tr{\mathop{\rm tr}\nolimits}

\def\Diff{\mathop{\rm Diff}\nolimits}

\def\can{\mathop{\rm can}\nolimits}
\def\var{\mathop{\rm var}\nolimits}
\def\id{\mathop{\rm id}\nolimits}

\def\codim{\mathop{\rm codim}\nolimits}

\def\Supp{\mathop{\rm Supp}\nolimits}

\def\Irr{\mathop{\rm Irr}\nolimits}

\newcommand{\del}{\partial}
\newcommand{\delbar}{\overline{\del}}

\newcommand{\deltatilde}{\widetilde{\delta}}

\newcommand{\nhom}{{\mathcal Hom}}

\newcommand{\next}{{\mathcal Ext}}

\newcommand{\nbar}{\underline{n}}

\newcommand{\mbar}{\underline{m}}

\newcommand{\pbar}{\underline{p}}

\newcommand{\nibar}{\underline{2}}

\newcommand{\harmonicbundle}{(E,\delbar_E,\theta,h)}

\newcommand{\barz}{\overline{z}}
\newcommand{\zbar}{\barz}

\newcommand{\baralpha}{\overline{\alpha}}
\newcommand{\alphabar}{\baralpha}
\newcommand{\barlambda}{\overline{\lambda}}
\newcommand{\lambdabar}{\barlambda}
\newcommand{\fbar}{\overline{f}}

\newcommand{\etabar}{\overline{\eta}}
\newcommand{\xibar}{\overline{\xi}}
\newcommand{\omegabar}{\overline{\omega}}

\newcommand{\tbar}{\overline{t}}

\newcommand{\sbar}{\overline{s}}

\newcommand{\Abar}{\overline{A}}
\newcommand{\bbar}{\overline{b}}

\newcommand{\nbigxlambda}{\nbigx^{\lambda}}

\newcommand{\KMS}{{\rm{K\!M\!S}}}

\newcommand{\kmsmap}{\gminik}
\newcommand{\paramap}{\gminip}
\newcommand{\eigenmap}{\gminie}

\newcommand{\lefttop}[1]{{}^{#1}\!}

\def\irr{\mathop{\rm irr}\nolimits}
\def\reg{\mathop{\rm reg}\nolimits}
\def\Def{\mathop{\rm Def}\nolimits}

\def\nil{\mathop{\rm nil}\nolimits}

\newcommand{\Fzero}{F^{(\lambda_0)}}

\newcommand{\Vzero}{V^{(\lambda_0)}}

\newcommand{\psizero}{\psi^{(\lambda_0)}}

\newcommand{\tildepsi}{\widetilde{\psi}}
\newcommand{\psitilde}{\tildepsi}

\newcommand{\psitildezero}{\psitilde^{(\lambda_0)}}

\newcommand{\nbigqzero}{\nbigq^{(\lambda_0)}}

\newcommand{\Uzero}{U^{(\lambda_0)}}

\newcommand{\deldel}{\eth}
\newcommand{\deldelbar}{\overline{\deldel}}

\newcommand{\lamda}{\lambda}

\newcommand{\distribution}{\gbigd\gminib}

\newcommand{\vecnbign}{{\pmb{\mathcal N}}}

\newcommand{\vecnbigt}{{\pmb{\mathcal T}}}

\newcommand{\vecnbigv}{{\pmb{\mathcal V}}}

\newcommand{\rtriplecat}{\nbigr\textrm{-Tri}}
\newcommand{\Dtriplecat}{\nbigd\textrm{-Tri}}
\newcommand{\Dcomplextriplecat}{\nbigc(\nbigd)\textrm{-Tri}}
\newcommand{\Ddoublecomplextriplecat}{\nbigc^{(2)}(\nbigd)\textrm{-Tri}}
\newcommand{\cnumcomplextriplecat}{\nbigc(\cnum)\textrm{-Tri}}
\newcommand{\cnumtriplecat}{\cnum\textrm{-Tri}}

\newcommand{\Etilde}{\widetilde{E}}

\newcommand{\Vhat}{\widehat{V}}

\newcommand{\gbar}{\overline{g}}

\newcommand{\vtilde}{\widetilde{v}}

\newcommand{\ftilde}{\widetilde{f}}
\newcommand{\utilde}{\widetilde{u}}
\newcommand{\Ftilde}{\widetilde{F}}

\newcommand{\gminiabar}{\overline{\gminia}}

\newcommand{\sigmatilde}{\widetilde{\sigma}}
\newcommand{\nbigltilde}{\widetilde{\nbigl}}
\newcommand{\nbigdhat}{\widehat{\nbigd}}

\newcommand{\ellsitabar}{\underline{\ell}}

\newcommand{\Xtilde}{\widetilde{X}}

\newcommand{\DDhat}{\widehat{\DD}}

\newcommand{\Ltilde}{\widetilde{L}}
\newcommand{\nbigmtilde}{\widetilde{\nbigm}}

\def\ord{\mathop{\rm ord}\nolimits}

\def\MT{\mathop{\rm MT}\nolimits}
\def\MTW{\mathop{\rm MTW}\nolimits}
\def\MTM{\mathop{\rm MTM}\nolimits}
\def\MTN{\mathop{\rm MTN}\nolimits}
\def\MTS{\mathop{\rm MTS}\nolimits}
\def\PTS{\mathop{\rm PTS}\nolimits}

\def\Gal{\mathop{\rm Gal}\nolimits}

\def\moderate{\mathop{\rm mod}\nolimits}

\def\Hol{\mathop{\rm Hol}\nolimits}

\def\exchange{\mathop{\rm ex}\nolimits}

\def\TNIL{\mathop{\rm TNIL}\nolimits}
\def\Glue{\mathop{\rm Glue}\nolimits}

\def\DR{\mathop{\rm DR}\nolimits}

\def\good{\mathop{\rm good}\nolimits}
\def\Map{\mathop{\rm Map}\nolimits}
\def\hol{\mathop{\rm hol}\nolimits}

\def\Per{\mathop{\rm Per}\nolimits}

\def\fil{\mathop{\rm fil}\nolimits}
\def\RMF{\mathop{\rm RMF}\nolimits}
\def\Vect{\mathop{\rm Vect}\nolimits}
\def\VTS{\mathop{\rm TS}\nolimits}

\def\adm{\mathop{\rm adm}\nolimits}
\def\sp{\mathop{\rm sp}\nolimits}

\def\herm{\mathop{\rm herm}\nolimits}
\def\Glu{\mathop{\rm Glu}\nolimits}

\def\Forget{\mathop{\rm Forget}\nolimits}
\def\St{\mathop{\rm St}\nolimits}
\def\res{\mathop{\rm res}\nolimits}

\def\integral{\mathop{\rm int}\nolimits}
\def\pt{\mathop{\rm pt}\nolimits}
\def\pre{\mathop{\rm pre}\nolimits}
\def\sm{\mathop{\rm sm}\nolimits}
\def\loc{\mathop{\rm loc}\nolimits}
\def\MF{\mathop{\rm MF}\nolimits}

\newcommand{\Dbar}{\overline{D}}

\newcommand{\Wtilde}{\widetilde{W}}
\newcommand{\Ptilde}{\widetilde{P}}

\newcommand{\iotabar}{\overline{\iota}}

\newcommand{\vecj}{{\pmb j}}
\newcommand{\vecD}{{\pmb D}}

\newcommand{\nbigstilde}{\widetilde{\nbigs}}

\newcommand{\nbigmhat}{\widehat{\nbigm}}

\newcommand{\Nhat}{\widehat{N}}

\newcommand{\gtilde}{\widetilde{g}}

\newcommand{\nbigttilde}{\widetilde{\nbigt}}

\newcommand{\Ctilde}{\widetilde{C}}

\newcommand{\nbigftilde}{\widetilde{\nbigf}}

\newcommand{\Mtilde}{\widetilde{M}}

\newcommand{\Ybar}{\overline{Y}}

\newcommand{\nbiglbar}{\overline{\nbigl}}

\newcommand{\betabar}{\overline{\beta}}

\newcommand{\nbigktilde}{\widetilde{\nbigk}}

\newcommand{\vecnbigi}{{\pmb \nbigi}}

\newcommand{\nbigmlambda}{\nbigm^{\lambda}}
\newcommand{\nrhom}{R{\mathcal Hom}}
\newcommand{\DDD}{\pmb D}

\newcommand{\vecH}{{\pmb H}}
\newcommand{\nbiggzero}{\nbigg^{(\lambda_0)}}

\newcommand{\Hhat}{\widehat{H}}

\newcommand{\Ntilde}{\widetilde{N}}
\newcommand{\nbigxzero}{\nbigx^{(\lambda_0)}}
\newcommand{\nbigdzero}{\nbigd^{(\lambda_0)}}
\newcommand{\nutilde}{\widetilde{\nu}}

\newcommand{\gminiatilde}{\widetilde{\gminia}}

\newcommand{\nbightilde}{\widetilde{\nbigh}}

\newcommand{\Vbar}{\overline{V}}
\newcommand{\Wbar}{\overline{W}}
\newcommand{\Nbar}{\overline{N}}

\newcommand{\Xbar}{\overline{X}}

\newcommand{\gammatilde}{\widetilde{\gamma}}

\newcommand{\veci}{\pmb i}
\newcommand{\nbigatilde}{\widetilde{\nbiga}}
\newcommand{\starbar}{\overline{\star}}

\newcommand{\nbigmbar}{\overline{\nbigm}}
\newcommand{\nbigfbar}{\overline{\nbigf}}

\newcommand{\vecgbigi}{\pmb \gbigi}

\newcommand{\vecNtilde}{\widetilde{\vecN}}
\newcommand{\nbigmzero}{\nbigm^{(\lambda_0)}}
\newcommand{\nbigqtilde}{\widetilde{\nbigq}}
\newcommand{\Lhat}{\widehat{L}}

\newcommand{\vecNhat}{\widehat{\vecN}}

\newcommand{\IItilde}{\widetilde{\II}}

\newcommand{\nbigtbar}{\overline{\nbigt}}
\newcommand{\Lbar}{\overline{L}}
\newcommand{\Mbar}{\overline{M}}

\newcommand{\vecetilde}{\widetilde{\vece}}
\newcommand{\etilde}{\widetilde{e}}

\newcommand{\nbigdhatzero}{\widehat{\nbigd}^{(\lambda_0)}}
\newcommand{\vecLtilde}{\widetilde{\vecL}}
\newcommand{\nbigptilde}{\widetilde{\nbigp}}

\newcommand{\newTate}{\pmb{T}}

\newcommand{\MTWint}{\MTW^{\integral}}
\newcommand{\MTint}{\MT^{\integral}}
\newcommand{\MTMint}{\MTM^{\integral}}
\newcommand{\MTSint}{\MTS^{\integral}}
\newcommand{\IMTMint}{\IMTM^{\integral}}
\newcommand{\VTSint}{\VTS^{\integral}}

\newcommand{\MTSpol}{\nbigp}
\newcommand{\IMTM}{\nbigm}
\newcommand{\bikkuri}{!}
\newcommand{\nbigotilde}{\widetilde{\nbigo}}
\newcommand{\Sigmabar}{\overline{\Sigma}}

\newcommand{\lambdatilde}{\widetilde{\lambda}}
\newcommand{\Upsilontilde}{\widetilde{\Upsilon}}
\newtheorem{thm}{Theorem}[section]
\newtheorem{cor}[thm]{Corollary}

\newtheorem{rem}[thm]{Remark}
\newtheorem{lem}[thm]{Lemma}
\newtheorem{prop}[thm]{Proposition}
\newtheorem{df}[thm]{Definition}
\newtheorem{example}[thm]{Example}

\makeindex
\begin{document}

\frontmatter%%%%%%%%%%%%%%%%%%%%%%%%%%%%%%%%%%%%%%%%%%%%%%%%%%%%%%

\begin{abstract}

%\preface

We introduce mixed twistor $D$-modules,
and establish the fundamental functorial property.
We also prove that they are described as
the gluing of admissible variations of mixed twistor structure.
In a sense, mixed twistor $D$-modules could be 
regarded as a twistor version of
mixed Hodge modules due to M. Saito .
In the other sense, they could be 
a mixed version of pure twistor $D$-modules
studied by C. Sabbah and the author.
A theory of mixed twistor $D$-modules
is one of the ultimate goals in the study
suggested by Simpson's Meta Theorem,
and the author hopes that 
it would play a basic role in the Hodge theory for 
holonomic $D$-modules possibly with irregular singularity.

\end{abstract}
\begin{altabstract}
\end{altabstract}

\altkeywords{}

\maketitle
\tableofcontents

\mainmatter%%%%%%%%%%%%%%%%%%%%%%%%%%%%%%%%%%%%%%%%%%%%%%%%%%%%%%%

\chapter{Introduction}
In this paper, we introduce mixed twistor $D$-modules,
and establish their fundamental functorial property.
We also prove that they are described as
the gluing of admissible variations of mixed twistor structure.
In a sense, mixed twistor $D$-modules could be 
regarded as a twistor version of
mixed Hodge modules due to M. Saito 
(\cite{saito1} and \cite{saito2}).
In the other sense, they could be regarded as 
a mixed version of pure twistor $D$-modules
which was introduced by
C. Sabbah (\cite{sabbah2} and \cite{sabbah5})
and studied by himself and the author
(\cite{mochi2} and \cite{mochi7}).
A theory of mixed twistor $D$-module
is one of the ultimate goals in the project
to generalize theories of Hodge structure
to those of twistor structures,
suggested by C. Simpson's Meta Theorem.
Actually, we closely follow the fundamental strategy
due to Saito for mixed Hodge modules,
although there are several issues to overcome,
partially because of the difference of the ingredients
in Hodge modules and twistor $D$-modules.
The author hopes that 
this study would be a part of the foundation
for the further study on the Hodge theory for 
holonomic $D$-modules possibly with irregular singularity. 
He also hopes that it would be a help
for readers to get into the deep theory due to Saito.

\section{Mixed Hodge modules}

Let us briefly review a stream in the Hodge theory,
which is related to the functoriality with respect to various operations,
which was finally accomplished with great generality
by the theory of mixed Hodge modules due to M. Saito.
The author regrets that 
the following is quite restricted by 
his personal interest, 
and that it is not exhaustive.
The readers can find a more thorough review
in \cite{peters-steenbrink}.

\vspace{.1in}

A variation of Hodge structure on a complex manifold $X$
is a pair of $\rnum$-local system of finite rank $\nbigv$
and a Hodge filtration $F$ which is 
a decreasing filtration of holomorphic subbundles
of $\nbigv\otimes_{\rnum}\nbigo_X$ indexed by integers
satisfying the Griffiths transversality.
Although we may replace $\rnum$
with other algebras such as $\seisuu$ and $\real$,
we omit such details here.
A variation of Hodge structure is called pure of weight $w$, 
if each restriction $(\nbigv,F)_{|Q}$ $(Q\in X)$
is a pure Hodge structure of weight $w$.
It is called polarizable, if moreover it has a polarization,
i.e., there exists a flat $(-1)^w$-symmetric pairing 
$S$ of $\nbigv$ such that
$S_{|Q}$ $(Q\in X)$ is a polarization of
the pure Hodge structure $(\nbigv,F)_{|Q}$.
A variation of mixed Hodge structure is
a variation of Hodge structure $(\nbigv,F)$
with a weight filtration $W$ of $\nbigv$
which is an increasing filtration indexed by integers,
such that $\Gr_w^W(\nbigv,F)$ $(w\in \seisuu)$
are pure of weight $w$.
It is called graded polarizable,
if each $\Gr_w^W(\nbigv,F)$ is polarizable.
In this paper, 
we almost always impose the polarizability
(resp. the graded polarizability)
to variations of pure (resp. mixed) Hodge structure.
So, we often omit the adjective ``graded polarizable''.

The notion of polarized variation of Hodge structure
was originally discovered by P. Griffiths
as {\em something} on 
the Gauss-Manin connections associated 
to smooth families of projective varieties.
This can already be regarded 
as one of the most basic and interesting cases of
the functoriality of Hodge structure
for projective push-forward.
The seminal work of Griffiths opened
several interesting research projects,
for example,
the study of polarized variation of Hodge structure
with singularity, which we will return later.

Inspired by the dream of motives,
P. Deligne discovered the notion of mixed Hodge structure,
and he proved a deep theorem which ensures
that the cohomology group of 
any complex algebraic variety
is naturally equipped with a mixed Hodge structure.
This can be regarded as one of the most important cases
of the functoriality of mixed Hodge structure.
He also proved the functoriality in various cases.
For example, 
if we are given a smooth family of smooth projective
varieties $f:\nbigy\lrarr S$ 
and a graded polarizable variation of mixed Hodge structure
$(\nbigv,F,W)$ on $\nbigy$,
then it was proved that 
the filtered local system 
$R^if_{\ast}(\nbigv,W)$ is equipped with
a naturally induced Hodge filtration $R^if_{\ast}F$,
and that $R^if_{\ast}(\nbigv,F,W)$ is 
graded polarizable variation of mixed Hodge structure.
He also observed crucial properties
of the induced variation of mixed Hodge structure,
including the Hard Lefschetz Theorem
in the pure case.
His insight has been quite influential on the subsequent works.

\vspace{.1in}
It is natural to ask what happens in the other more general cases.
For example,
if we are given a polarizable variation of Hodge structure 
$(\nbigv,F)$ on a quasi projective variety $Y$
which is not extendable on any projective completion of $Y$,
it is asked whether the cohomology group
$H^{p}(Y,\nbigv)$ or its variant
may have mixed or pure Hodge structure.
What happens for the singular family of singular varieties?
Finally, all of these questions were answered
by the functoriality of mixed Hodge modules.
But, for the theory,
it is fundamental to understand 
the asymptotic behaviour of polarized pure Hodge structure
and admissible variation of mixed Hodge structure.

\vspace{.1in}
As mentioned, the work of Griffiths naturally lead 
to the study of polarized variations of pure Hodge structure
with singularity.
An extremely important contribution was done
by W. Schmid \cite{sch}.
Let $X:=\{(z_1,\ldots,z_n)\,|\,|z_i|<1\}$
and
$D:=\bigcup_{i=1}^{\ell}\{z_i=0\}$.
Let $(\nbigv,F)$ be a polarizable variation of 
pure Hodge structure on $X-D$.
For simplicity, we assume that
the local monodromy endomorphisms
around any irreducible components of
$D$ are unipotent.
The nilpotent orbit theorem of Schmid
ensures that, around any $P\in D$,
the polarized variation of Hodge structure
can be approximated by an easier one called a nilpotent orbit.
It is not only interesting itself
but also the most important foundation
for the further investigation.
One of the important consequences
is that we obtain a nice object on  $X$,
not only on $X\setminus D$.
Namely, let $(V,\nabla)$ be the Deligne extension of
$\nbigv$ on $X$,
i.e., $V$ is the locally free $\nbigo_X$-module
with a logarithmic connection $\nabla$
such that 
(i) $(V,\nabla)_{|X\setminus D}=
 \nbigv\otimes_{\rnum}\nbigo_{X\setminus D}$,
(ii) the residues of $\nabla$ are nilpotent.
Then, $F$ is extended to a filtration of $V$
by holomorphic subbundles.
In the one dimensional case,
the study of singular polarized variation of pure Hodge structure
was accomplished by his $SL(2)$-orbit theorem,
which ensures that, if $n=\ell=1$,
the polarized variation of Hodge structure
can be approximated by an easier one called an $SL(2)$-orbit
around singularity.
As consequences, in the one variable case,
he obtained that the weight filtration of
the nilpotent part of the local monodromy
controls the growth order of the norms of
flat sections with respect to the hermitian metric
associated to the polarization.
He also obtained the polarized mixed Hodge structure
from the asymptotic data around the singularity,
which is called the limit mixed Hodge structure.
Note that, for the polarized variation of pure Hodge structure
associated to a degenerating family of smooth projective 
varieties,
the asymptotic behaviour was 
intensively studied by J. Steenbrink 
with a different method \cite{Steenbrink-limit-Hodge}.

The higher dimensional case was accomplished
by the definitive works
by E. Cattani, A. Kaplan, M. Kashiwara, T. Kawai and W. Schmid
(\cite{ck}, \cite{cks1}, \cite{cks2},
 \cite{Kashiwara-asymptotic-behaviour},
 \cite{k3}).
In the above situation,
for each point $P\in D$,
we obtain the limit mixed Hodge structure
from the asymptotic data around $P$,
which is a polarized mixed Hodge structure
in several variables.
It turned out that the limit mixed Hodge structure
controls the behaviour of $(\nbigv,F)$ around $P$.
They obtained a generalization of the norm estimate.
They also obtained a rather strong constraint on the nilpotent parts 
of the local monodromy along the loops around $\{z_i=0\}$
$(i=1,\ldots,\ell)$.
Moreover, they proved various interesting property
of polarized mixed Hodge structure,
which are significant for their study on $L^2$-cohomology.
Although it requires much more preparation
to describe their results precisely, which we do not intend here,
they are quite impressive.

\vspace{.1in}
As for singular graded polarizable
variations of mixed Hodge structure,
it was one of the main issues to clarify 
what condition should be imposed at the boundary.
Thanks to the studies of 
Kashiwara, Steenbrink and S. Zucker
(\cite{kashiwara-mixed-Hodge}, \cite{steenbrink-zucker}, 
\cite{Zucker-VMHS}),
it turned out that the admissibility condition is 
the most appropriate one.
Let us recall it in the case that
$X=\{(z_1,\ldots,z_n)\,|\,|z_i|<1\}$ and 
$D=\bigcup_{i=1}^{\ell}\{z_i=0\}$ as above.
Let $(\nbigv,F,W)$ be a graded polarizable variation of
mixed Hodge structure on $X\setminus D$.
For simplicity, suppose that the monodromy $g_i$
along the loops around $\{z_i=0\}$ are unipotent.
Let $N_i:=\log g_i$.
Let $(V,\nabla)$ be the Deligne extension of $\nbigv$,
which is equipped with the flat filtration $W$.
We should impose that the filtration $F$ is extended to 
a filtration of $V$ by holomorphic subbundles
such that 
$\Gr^F\Gr^W(V)$ is a locally free $\nbigo_X$-module.
We should also impose the existence of
a relative monodromy weight filtration $M(N_i;W)$
of $N_i$ with respect to 
the induced filtration $W$ on the space of the multi-valued
flat sections of $\nbigv$.
It was introduced by Steenbrink-Zucker in the case $n=1$,
and by Kashiwara in the higher dimensional case.
(The condition here is a priori stronger but equivalent
by results in \cite{kashiwara-mixed-Hodge}.)
Moreover,
Kashiwara introduced and studied
{\em infinitesimal mixed Hodge modules},
which is the ``mixed version'' of polarizable mixed Hodge structures.
By generalizing some construction of Steenbrink-Zucker,
he constructed some natural filtrations
which are crucial in the study on mixed Hodge modules.

\vspace{.1in}

As mentioned, one of the motivations
of singular variation of Hodge structure
was to establish the functoriality of Hodge structure,
as a generalization of the results of Deligne.
Let $X$ be a smooth projective variety
with a normal crossing hypersurface $D$.
Let $(\nbigv,F)$ be any polarizable variation of pure Hodge structure
on $X\setminus D$.
One of the issues in those days was to show that
there exists a natural pure Hodge structure
on the intersection cohomology group of $\nbigv$.
If $(\nbigv,F)$ has no singularity at $D$,
then it follows from the result of Deligne.
In the singular case,
the contribution of Zucker \cite{z}
is quite important.
Namely, 
he studied the issue in the case $\dim X=1$,
and he proved that the intersection cohomology group
is isomorphic to the $L^2$-cohomology.
He also developed the $L^2$-harmonic theory
for singular polarized variation of pure Hodge structure
on projective curves.
As a result, he obtained a naturally induced
pure Hodge structure
on the intersection cohomology of $\nbigv$.
In the higher dimensional case,
Cattani-Kaplan-Schmid and Kashiwara-Kawai
established it by making good use of their results
on the asymptotic behaviour of polarized variation of 
Hodge structure,
and by generalizing the method of Zucker.
As for the mixed case,
for an admissible variation of mixed Hodge structure
on curves, Steenbrink-Zucker proved that 
the various naturally associated cohomology groups
have mixed Hodge structure,
based on their results on the asymptotic behaviour.

\vspace{.1in}
This stream of research for functoriality
was eventually accomplished with much more great generality
by the theory of mixed Hodge modules due to M. Saito.
A cohomology theory can be regarded as a part of
the theory of six functors on
the derived categories of some type of sheaves.
The theory of mixed Hodge modules 
ensures that the derived functors for 
 $\rnum$-perverse sheaves on complex algebraic varieties
can be enriched by mixed Hodge structures.
(In this introduction,
we consider only polarizable pure Hodge modules
and graded polarizable mixed Hodge modules,
we omit the adjectives ``polarizable'' or ``graded polarizable''.)

\vspace{.1in}
Very roughly, a Hodge module on a complex manifold $X$
consists of a perverse sheaf $P$ 
with a Hodge filtration $F$ on
the regular holonomic $D$-module $\nbigm$
corresponding to $P$,
i.e., 
(i) $\DR_X(\nbigm)\simeq P\otimes\cnum$,
(ii) $F$ is an increasing filtration of $\nbigm$
 by coherent $\nbigo_X$-modules indexed by integers
 such that $F_j(\nbigd_X)\cdot F_i(\nbigm)\subset
 F_{i+j}(\nbigm)$,
 where $F_j(\nbigd_X)$ denotes the natural filtration
 of the sheaf of holomorphic differential operators on $X$
 by orders. 
In his highly original and genius work,
Saito invented the appropriate definitions of 
pure and mixed conditions for such filtered $D$-modules,
and he established their fundamental properties.
The most important theorems in the theory 
are the functoriality with respect to six operations,
and the description of pure and mixed Hodge modules.

For the functoriality in the pure case,
he proved the Hard Lefschetz Theorem.
Namely, let $f:X\lrarr Y$ be a projective morphism
of smooth projective varieties.
Let $(P,F)$ be any polarizable pure Hodge module
of weight $w$ on $X$.
Then, the $i$-th cohomology of the push-forward
$f^i_{\dagger}\nbigv$ is equipped with
a naturally induced Hodge filtration $f^i_{\dagger}F$,
so that $f^i_{\dagger}(P,F)$ is 
a polarizable pure Hodge module of weight $w+i$.
Moreover, for the morphism 
$L:f^i_{\dagger}P\lrarr f^{i+2}_{\dagger}P$
induced by the first Chern class of a relatively ample line bundle,
the morphism 
$L^i:f^{-i}_{\dagger}P\lrarr f^{i}_{\dagger}P$
is an isomorphism.
This theorem is a generalization of
the classical and important theorem of
Beilinson-Bernstein-Deligne-Gabber
on perverse sheaves of geometric origin.

As for the functoriality in the mixed case,
he constructed the six operations
together with the nearby and vanishing cycle functors
for the derived category of mixed Hodge modules
on algebraic varieties,
which are compatible with those for
the derived category of perverse sheaves.
We also have nearby and vanishing cycle functors.

Because the definitions of pure and mixed Hodge modules
are complicated,
it is important to know what objects are contained
in the categories.
Saito proved that 
a polarizable (resp. graded polarizable)
variation of pure (resp. mixed) Hodge structure on $X$
naturally gives a pure (mixed) Hodge module on $X$,
as expected.
Hence, the simplest variation of pure Hodge structure
$\rnum_X$ naturally gives a pure Hodge module.
Hence, if a perverse sheaf $P$ on $X$
is obtained from $\rnum_Y$ on some $Y$
by successive use of six functors,
it naturally underlies a mixed Hodge module.
If a perverse sheaf is of geometric origin,
then it naturally underlies a pure Hodge module.
In particular, 
the category of pure and mixed Hodge modules 
contain many objects.
Moreover, 
he proved the more general results for the description.
In the pure case,
he proved the following.
\begin{itemize}
\item
 Let $Z\subset X$ be a closed irreducible
 complex analytic subvariety.
 Let $U\subset Z$ be a complement of
 a closed analytic subset, such that $U$ is smooth.
 Let $\iota:U\lrarr X$ be the inclusion.
 Let $(\nbigv,F)$ be any polarizable variation of Hodge structure 
 on $U$.
 Then, the perverse sheaf $\iota_{\ast!}\nbigv$,
 which is the minimal extension of $\nbigv$,
 is naturally equipped with the Hodge filtration $F$
 so that $(\iota_{\ast!}\nbigv,F)$ 
 is a polarizable pure Hodge module.
\item
 Conversely, any polarizable pure Hodge module 
 is the direct sum of such minimal extensions.
\end{itemize}
Hence, for example,
suppose that 
we are given a polarizable variation of Hodge structure
$(\nbigv,F)$ on $X\setminus D$,
where $X$ is a complex manifold, and 
$D$ is a normal crossing hypersurface.
We obtain the pure Hodge module $(P,F)$ on $X$,
obtained as the minimal extension of $(\nbigv,F)$,
by the above description.
For the canonical map $a_X$ of $X$ to a point,
the $i$-th cohomology of 
the push-forward $a^i_{X\dagger}(P)$
is naturally equipped with the Hodge filtration
by the functoriality of the pure Hodge modules.
It means that the intersection cohomology of $\nbigv$
is equipped with a naturally induced pure Hodge structure,
which implies 
the theorem of Zucker, Cattani-Kaplan-Schmid
and Kashiwara-Kawai.

In the mixed case, Saito established the following:
\begin{itemize}
\item
Let $X$, $Z$, $U$ and $\iota$ be as above.
%Let $Z\subset X$ be any closed irreducible complex analytic
%subvariety.
%Let $U\subset Z$ be the complement of
%a closed complex analytic subset,
%such that $U$ is smooth.
%Let $\iota:U\lrarr X$ denote the inclusion.
Let $(\nbigv,F,W)$ be
an admissible variation of mixed Hodge structure.
Then, the perverse sheaves
$\iota_{\ast}\nbigv$ and $\iota_!\nbigv$
are equipped with
induced Hodge filtrations $\Ftilde$
and weight filtrations $\Wtilde$
such that
$(\iota_{\star}\nbigv,\Ftilde,\Wtilde)$
$(\star=\ast,!)$
are mixed Hodge modules.
\item
Conversely, any mixed Hodge modules on $X$
are obtained as the ``gluing'' of 
admissible variation of mixed Hodge structures.
\end{itemize}
It implies that, for example,
we have a natural mixed Hodge structure
on various cohomology groups
associated to an admissible variation of mixed Hodge structure.

\begin{rem}
The theory of pure and mixed Hodge modules
can be regarded as a counterpart of 
the theory of pure and mixed $\ell$-adic sheaves
on algebraic varieties over finite fields,
which has been influential in 
various fields of mathematics
including number theory and representation theory.
See a very nice book {\rm\cite{hotta-tanisaki}} 
for more details on 
the philosophical background of Hodge modules,
and for applications of the theory of Hodge modules
to representation theory.
\hfill\qed
\end{rem}

\section{From Hodge toward twistor}

As mentioned,
it is our purpose in this paper
to study a twistor version of mixed Hodge modules.
It is C. Simpson who introduced the notion of twistor structure
as an underlying structure of Hodge structure.
He proposed a principle called Simpson's Meta Theorem,
which says that stories of Hodge structure should be
generalized to stories of twistor structure.

\vspace{.1in}

He introduced twistor structure
to understand harmonic bundles
in a deeper way.
Let $(V,\nabla)$ be a flat bundle on a complex manifold $X$
with a hermitian metric $h$.
We have a unique decomposition
$\nabla=\nabla^u+\Phi$
into a unitary connection 
and a self-adjoint section of $\End(V)\otimes\Omega^1$.
We have the decompositions into
$\nabla^u=\delbar_V+\del_V$
and $\Phi=\theta^{\dagger}+\theta$
into the $(0,1)$-part and the $(1,0)$-part.
Then, $(V,\nabla,h)$ is called a harmonic bundle,
if $(V,\delbar_V,\theta)$ is a Higgs bundle.
In that case, the metric $h$ is called pluri-harmonic.

One of the most important classes of 
harmonic bundles is polarized variation of Hodge structure.
The hermitian metrics induced by polarizations
of polarizable variations of Hodge structure
are pluri-harmonic.
From the beginning of his study,
Simpson was motivated 
by the investigation of
polarized variations of Hodge structure.
He gave a method to construct polarized variation of Hodge structure
by using the Kobayashi-Hitchin correspondence for harmonic bundles.
He observed that various properties of polarized variation of
Hodge structure are naturally generalized
to those for harmonic bundles.
For example, he developed the harmonic theory
for harmonic bundles
as a generalization of that
for polarized variation of Hodge structure,
and he proved that the push-forward of
any harmonic bundle of 
any smooth family of projective varieties
is naturally a harmonic bundle.
To pursue this analogy in a deeper level,
he introduced twistor structure,
and observed that harmonic bundles can be regarded
as polarized variations of pure twistor structure.
Thus, he established the analogy between
polarized variations of Hodge structure
and harmonic bundles
in the level of the definitions.
This important idea enables us to consider
a twistor version of various objects
appeared in the Hodge theory.

This is quite efficient in the study of 
the asymptotic behaviour of 
singular harmonic bundles,
which was studied by Simpson himself
and the author. 
The twistor viewpoint suggests us
how to formulate generalizations
of results of Cattani, Kaplan, Kashiwara, Kawai and Schmid.
Indeed, we obtain a nice object on $X$
from a harmonic bundle on $X\setminus D$,
and we also obtain the limit mixed twistor structure,
which is quite useful to control the nilpotent part of the residues.
However, we would like to mention that
there are also some phenomena 
which do not appear for polarized variation of Hodge structure,
such as KMS-structure and Stokes structure,
and that the proofs are not necessarily 
given in parallel ways.

It is also suggested by Simpson's Meta Theorem
that we should have a twistor version of 
the theory of pure and mixed Hodge modules.
In the pure case, it was pursued by C. Sabbah and the author.
Sabbah prepared the notion of $\nbigr$-triples 
as an ingredient to define twistor $D$-modules,
which can be regarded as generalization of 
tuples of $\real$-perverse sheaf and filtered $D$-module,
and suitable to consider even in the case
that the underlying $D$-modules are irregular.
Based on Saito's strategy,
he gave the appropriate definition of pure twistor $D$-modules
and the framework to prove the Hard Lefschetz Theorem,
i.e., the functoriality for projective push-forward.
The correspondence between tame harmonic bundles
and regular pure twistor $D$-modules
was established in \cite{mochi2}.
In the wild case, the basic properties were established
in \cite{mochi7}.

It is interesting to have the correspondence between 
semisimple holonomic $D$-modules
and pure twistor $D$-modules  on projective varieties,
which does not appear in the theory of pure Hodge modules.
As a result, we obtain that semisimplicity of algebraic
holonomic $D$-modules is preserved 
by projective push-forward.

\vspace{.1in}

In this paper,
we introduce mixed twistor $D$-modules,
and prove the fundamental properties.
For the author, it is one of the ultimate goals of the research
for years, driven by Simpson's Meta Theorem.

There are various intermediate objects
between twistor structure and Hodge structure
such as integrable twistor structure and TERP structure.
So, we could have variants of mixed twistor $D$-modules
by considering additional structures.
Because the twistor structure could be most basic among them,
the author hopes that 
mixed twistor $D$-modules would play 
a basic role in the study of Hodge structure
on holonomic $D$-modules.

\section{Mixed twistor $D$-modules}

In the rest of this introduction,
let us briefly review the theory of pure twistor $D$-modules,
and explain our issues in the study of mixed twistor $D$-modules.

\subsection{Pure twistor $D$-modules}

For any complex manifold $X$,
the product 
$\cnum_{\lambda}\times X$ is denoted by
$\nbigx$.
We have the sheaf of relative differential operators
$\nbigd_{\nbigx/\cnum_{\lambda}}$,
and the relative tangent sheaf
$\Theta_{\nbigx/\cnum_{\lambda}}$.
Then, $\nbigr_X$ is the sheaf of subalgebras
generated by $\nbigo_{\nbigx}$
and $\lambda\cdot\Theta_{\nbigx/\cnum_{\lambda}}$.

We have two basic
conditions on $\nbigr_X$-modules.
One is strictness,
i.e., flat over $\nbigo_{\cnum_{\lambda}}$.
The other is holonomicity.
Namely, the characteristic variety of
any coherent $\nbigr_X$-module $\nbigm$ is defined
as in the case of $D$-modules,
denoted by $Ch(\nbigm)$.
It is a subvariety in
$\cnum_{\lambda}\times T^{\ast}X$.
If $Ch(\nbigm)$ is contained in the product of
$\cnum_{\lambda}$ and
a Lagrangian subvariety in $T^{\ast}X$,
the $\nbigr_X$-module $\nbigm$ is called holonomic.

An $\nbigr_X$-triple is a tuple of
$\nbigr_X$-modules $\nbigm_i$ $(i=1,2)$
with a sesqui-linear pairing $C$.
To explain what is sesqui-linear paring,
we need a preparation.
Let $\vecS$ denote the circle 
$\{\lambda\in\cnum_{\lambda}\,|\,|\lambda|=1\}$.
Let $\sigma:\vecS\lrarr\vecS$ be given by
$\sigma(\lambda)=-\lambda=-\lambdabar^{-1}$.
The induced involution
$\vecS\times X\lrarr\vecS\times X$ is also denoted 
by $\sigma$.

Let $\distribution_{\vecS\times X/\vecS}$
denote the sheaf of distributions on
$\vecS\times X$ which are continuous
in the $\vecS$-direction in an appropriate sense.
This sheaf is naturally  a module over 
$\nbigr_{X|\vecS\times X}\otimes
 \sigma^{\ast}\nbigr_{X|\vecS\times X}$.
Then, a sesqui-linear pairing of
$\nbigr_X$-modules $\nbigm_i$ $(i=1,2)$
is an
$\nbigr_{X|\vecS\times X}\otimes
 \sigma^{\ast}
 \nbigr_{X|\vecS\times X}$-homomorphism
$\nbigm_{1|\vecS\times X}
 \otimes
 \sigma^{\ast}\nbigm_{2|\vecS\times X}
\lrarr
 \distribution_{\vecS\times X/\vecS}$.
An $\nbigr_X$-triple is called 
strict (resp. holonomic),
if the underlying $\nbigr$-modules are 
strict (resp. holonomic).
The category of pure twistor $D$-module
is constructed as a full subcategory of
strict holonomic $\nbigr$-triples.

Let us recall how to impose some conditions
on strict holonomic $\nbigr$-triples.
In the case of variation of Hodge structures,
which is a $\rnum$-local system with a Hodge filtration,
it is defined to be pure,
if its restriction to the fiber over each point is pure.
But, for $\nbigr$-triples
or even for $D$-modules,
the restriction to a point is not so well behaved.
Instead, for holonomic $D$-modules,
there is a nice theory for restriction to
hypersurfaces.
Namely, we have the nearby and 
vanishing cycle functors,
which describes the behaviour of
the holonomic $D$-modules around 
the hypersurface in some degree.
To define some condition for holonomic $D$-modules,
it seems natural to consider
the conditions on nearby and vanishing cycle sheaves
inductively, instead of the restriction to a point.
Similarly,
to define some condition for $\nbigr$-triples,
we would like to consider the condition
on the appropriately defined
nearby and vanishing cycle functors
for $\nbigr$-triples.
This is a basic strategy due to Saito,
and it may lead us to an inductive definitions
of pure and mixed twistor $D$-modules,
as a variant of pure and mixed Hodge modules.

A strict holonomic $\nbigr$-triple
$\nbigt$ is called pure of weight $w$,
if the following holds.
First,
we impose that,
for any open subset $U\subset X$
with a holomorphic function $g$,
$\nbigt_{|U}$ is strictly $S$-decomposable along $g$.
It implies that 
 we have the decomposition
 $\nbigt=\bigoplus \nbigt_Z$ by
 strict support,
 where $Z$ runs through
 closed irreducible subsets of $X$.
% (More precisely, we impose that
% $\nbigt$ is strictly $S$-decomposable.)
Then, we impose the conditions
 on each $\nbigt_Z$.
 If $Z$ is a point,
 $\nbigt_Z$ is supposed to be
the push-forward of 
 pure twistor structure of weight $w$
 by the inclusion of $Z$ into $X$.
 In the positive dimensional case,
 for any open subset of $X$
with a holomorphic function,
 if we take the graduation of the weight filtration
 of the naturally induced nilpotent morphism
 on the nearby cycle functor along the function,
 the $m$-th graded pieces are pure of weight $w+m$.
Then, inductively,
the notion of pure twistor $D$-module is defined.
Precisely, we should consider the polarizable object.
A polarization of $\nbigt$ is defined
as a Hermitian sesqui-linear duality of weight $w$
satisfying some condition on positivity,
which is also given in an inductive way
using the nearby cycle functor.

Let $\MT(X,w)$ denote the category of
polarizable wild pure twistor $D$-module of weight $w$.
Let us recall some of their fundamental properties;
(i) The category $\MT(X,w)$ is abelian and semisimple;
(ii)
 For objects $\nbigt_i\in\MT(X,w_i)$
 with a morphism
 $f:\nbigt_1\lrarr\nbigt_2$
 as $\nbigr$-triples such that $w_1>w_2$,
we have $f=0$;
(iii)
For any projective morphism $f:X\lrarr Y$,
and for any  $\nbigt\in\MT(X,w)$,
the $i$-th cohomology of 
the push-forward $f^i_{\dagger}\nbigt$
is an object in $\MT(Y,w+i)$.
 Moreover, 
 $f_{\dagger}\nbigt\simeq
 \bigoplus f_{\dagger}^i\nbigt[-i]$
in the derived category of $\nbigr$-triples;
(iv) Let $Z\subset X$ be a closed complex analytic subset.
 Let $Z_0\subset Z$ be a closed complex analytic subset
 such that $Z\setminus Z_0$ is smooth.
 Then, a wild harmonic bundle on $(Z,Z_0)$
 is naturally extended to pure twistor $D$-module on $X$;
(v) Conversely, any pure twistor $D$-modules are
 the direct sum of such minimal extensions of
 wild harmonic bundles;
(vi)
 Any semisimple algebraic holonomic $D$-module
 naturally underlies 
a polarizable wild pure twistor $D$-module.

\subsection{Mixed twistor $D$-modules}

To define mixed twistor $D$-modules,
we first consider filtered $\nbigr$-triples
$(\nbigt,L)$ such that
$\Gr^L_w(\nbigt)$ are pure of weight $w$,
where $L$ are locally
finite increasing complete exhaustive filtrations
indexed by integers.
Such naive objects are called pre-mixed twistor $D$-module.
They have nice functoriality for projective push-forward.
However, we need to impose additional conditions
for other standard functoriality
such as push-forward for open embeddings
and pull back.
Very briefly, to define mixed twistor $D$-module,
we impose 
(i) the filtered strict compatibility of
$L$ and the $V$-filtrations,
and
(ii) the existence of relative monodromy filtrations
on the nearby and vanishing cycle sheaves,
and they give the weight of mixed twistor $D$-modules
with smaller supports.
(It will be explained in \S\ref{section;11.4.9.1}.)

It is easy to show that mixed twistor $D$-modules
have nice functorial property 
for nearby and vanishing cycle functors
and projective push-forward.
However, it is not so easy to show the functorial property,
for example,
for the open embedding $X\setminus H\lrarr X$,
where $H$ is a hypersurface.
To establish more detailed property,
we need a concrete description of mixed twistor $D$-modules
as the gluing of admissible variations of mixed twistor structure.

\subsection{Gluing procedure}

For perverse sheaves and holonomic $D$-modules,
there are well established theories
to glue objects on $\{f=0\}$ and 
objects on $\{f\neq 0\}$
(\cite{beilinson2},
\cite{MacPherson-Vilonen}, \cite{Verdier}).
We need such gluing procedure for 
$\nbigr$-triples.
Because of the difference of ingredients,
it is not easy to generalize the method of gluing 
prepared in \cite{saito2} for Hodge modules
to the case of  $\nbigr$-triples.
Instead, we adopt the method of Beilinson in \cite{beilinson2},
which reduces the issue to the construction of
canonical prolongations $\nbigt[\star t]$ $(\star=\ast,!)$.
(See \S\ref{section;11.4.9.2}--\S\ref{section;11.4.3.1}.)

\subsection{Admissible variation of mixed twistor structure}

We prepare a general theory for
admissible variation of mixed twistor structure
(\S\ref{section;11.4.3.4}),
which is a natural generalization
of admissible variation of mixed Hodge structure.
Very briefly,
it is a filtered $\nbigr$-triple $(\nbigv,L)$
on a complex manifold $X$
with poles along a normal crossing hypersurface $D$.
We impose the conditions
(i) each $\Gr^L_w(\nbigv)$ comes from 
 a good wild harmonic bundle,
(ii) $\nbigv$ has good-KMS structure along $D$
 compatible with $L$,
(iii) the residues along the divisors have
 relative monodromy filtrations.
It is important to understand the specialization of
admissible variation of mixed twistor structure
along the divisors.
For that purpose,
it is essential to study the relative monodromy filtrations
and their compatibility.
So, as in \cite{kashiwara-mixed-Hodge},
we study the infinitesimal version
of admissible variation of mixed twistor structure,
called infinitesimal mixed twistor module
in \S\ref{section;11.4.3.21}.
We can show that it has nice property
as in the Hodge case.

Then, we study the canonical prolongation of
admissible mixed twistor structure $(\nbigv,L)$ on $(X,D)$
to pre-mixed twistor $D$-modules on $X$.
Let $D=D^{(1)}\cup D^{(2)}$ be a decomposition.
Recall that a good meromorphic flat bundle $V$ on $(X,D)$
is extended to a $D$-module $V[\ast D^{(1)}!D^{(2)}]$ on $X$.
We prepare similar procedure 
to make a pre-mixed twistor $D$-module
$(\nbigv,L)[\ast D^{(1)}!D^{(2)}]$ from $(\nbigv,L)$.
First, we construct the underlying $\nbigr$-triple
$\nbigv[\ast D^{(1)}!D^{(2)}]$
in \S\ref{section;11.4.9.10}.
One of the main tasks is to construct
a correct weight filtration $\Ltilde$
on $\nbigv[\ast D^{(1)}!D^{(2)}]$.
By applying the procedure
in \S\ref{section;11.4.9.10}
to each $L_j\nbigv$,
we obtain a naively induced filtration $L$
on $\nbigv[\ast D^{(1)}!D^{(2)}]$.
But, this is not the correct filtration.
Indeed, $\Gr_w^L(\nbigv[\ast D^{(1)}!D^{(2)}])$
are not pure of weight $w$, in general.
We need much more considerations
for the construction of the correct weight filtration.
It is essentially contained in 
\cite{kashiwara-mixed-Hodge}
and \cite{saito2},
but we shall give rather details,
which is one of the main themes
in \S\ref{section;11.4.9.11}
and \S\ref{section;11.4.3.21}--\S\ref{section;11.4.3.5}.

Once we obtain canonical prolongations
of admissible variations of mixed twistor structures
across normal crossing hypersurfaces,
we can glue them to obtain pre-mixed twistor $D$-modules,
which are called good pre-mixed twistor $D$-modules.
It is one of the main theorems to show that
any good pre-mixed twistor $D$-module
is a mixed twistor $D$-module
(Theorem \ref{thm;10.11.13.11}).
Then, by a rather formal argument,
we can show that 
any mixed twistor $D$-module can be expressed
as  gluing of admissible variation of mixed twistor $D$-modules
as in \S\ref{subsection;11.4.3.40}.
We can deduce some basic functoriality
by using this description.
See \S\ref{subsection;11.4.9.12}--\S\ref{subsection;11.4.9.13}.

\begin{rem}
In this paper, we will often omit ``variation of''
just for simplification.
For example,
an admissible variation of mixed twistor structure
is often called admissible mixed twistor structure.
\hfill\qed
\end{rem}

\subsection{Dual}

In \S\ref{section;11.4.9.20},
we study the dual of mixed twistor $D$-module,
which should be compatible with the dual for
the underlying $D$-modules.
Briefly, there are two issues we should address.
While the dual of a holonomic $D$-module 
in the derived category
is also a holonomic $D$-module,
we cannot expect such a property
for a general holonomic $\nbigr$-module.
This issue already appeared in the Hodge case,
and solved by Saito.
Even in the twistor case,
we can apply Saito's method in a rather straightforward way.
The other issue is the construction of
sesqui-linear pairing for the dual,
which did not appear in the Hodge case.
It is non-trivial even for sesqui-linear pairings 
of holonomic $D$-modules.
Let $M_i$ $(i=1,2)$ be holonomic $D$-modules.
Let $C:M_1\otimes \overline{M_2}\lrarr
\distribution_X$
be a $D_X\otimes D_{\Xbar}$-homomorphism.
We need to construct an induced pairing
for $\DDD M_1\otimes \overline{\DDD M_2}
\lrarr\distribution_X$.
We obviously have such a pairing,
if $M_i$ are smooth, i.e., flat bundles.
It is not difficult to construct it
in the case of regular holonomic $D$-modules,
thanks to the Riemann-Hilbert correspondence.
But, at this moment,
some additional arguments are required
in the non-regular case.
(It will be given in \S\ref{section;11.4.6.10}.)
Once we have such a pairing in the case of $D$-modules,
it is rather formal to construct it
in the case of mixed twistor $D$-modules.
By using the dual,
we can introduce the notion of real structure
for mixed twistor $D$-modules,
which has nice functorial property.

\section{Acknowledgements} 

This study grew out of my attempt to
understand the works due to 
Masaki Kashiwara \cite{kashiwara-mixed-Hodge}, 
and Morihiko Saito \cite{saito1}--\cite{saito4}.
The readers can find most essential ideas in their papers.
Another important source of the ideas
is the influential work of Alexander Beilinson 
\cite{beilinson1}, \cite{beilinson2}.
I thank Morihiko Saito for some discussions.
I thank Claude Sabbah for numerous discussions
on many occasions.
I hope that this study would be useful
for further investigation on
the theory of holonomic $D$-modules
with irregular Hodge structure.
I deeply thank Carlos Simpson.
It is impossible for me to mention
what I owe to him.
I just mention here that
his most fundamental principle (Simpson's Meta-Theorem)
invited me to this study.
I am grateful to Christian Schnell
for some motivating discussions.
Special thanks go to 
Pierre Deligne,
Kenji Fukaya, 
William Fulton,
David Gieseker,
Akira Kono, 
Mikiya Masuda,
Masa-Hiko Saito, 
Michael Thaddeus,
and Kari Vilonen.
I would like to express my gratitude to
Yves Andr\'e,
Philip Boalch, 
Claus Hertling,
Maxim Kontsevich,
Kyoji Saito,
and Christian Sevenheck,
for some discussions.
I thank Akira Ishii and Yoshifumi Tsuchimoto 
for their constant encouragement.

I thank the colleagues and the staff of the RIMS
for excellent environment for the work.
I thank the organizer of the conference
``International Conference on
 Noncommutative Geometry and Physics''
in which I gave a talk on this topic.
This work was supported by
the Grant-in-Aid for Scientific Research (C)
(No. 22540078),
Japan Society for the Promotion of Science.

\part{Gluing and specialization of $\nbigr$-triples}

\chapter{Preliminary}

In \S\ref{subsection;13.4.12.1},
we shall review basic theory of
$\nbigr$-triples \cite{sabbah2} 
and variants.
In \S\ref{subsection;13.4.12.2},
we give some procedure
to produce an $\nbigr$-triple
from a given $\nbigr$-triple
with a commuting tuple of nilpotent morphisms.
It is a reformulation of the construction of
twistor nilpotent orbit in \cite{mochi2} and \cite{mochi8}.
In \S\ref{subsection;11.2.1.10},
we introduce Beilinson triple,
which is useful in the study of gluing.

\section{$\nbigr$-triples}
\label{subsection;13.4.12.1}

\subsection{$\nbigr$-module}

Let $X$ be a complex manifold.
We set $\nbigx:=\cnum_{\lambda}\times X$.
\index{complex manifold $\nbigx$}
The projection $\nbigx\lrarr X$ 
is denoted by $p_{\lambda}$.
Let $\Theta_{\nbigx/\cnum_{\lambda}}$
denote the relative tangent sheaf of 
$\nbigx$ over $\cnum_{\lambda}$.
Let $D_{\nbigx/\cnum_{\lambda}}$
denote the sheaf of relative differential operators
on $\nbigx$ over $\cnum_{\lambda}$.
Recall that
$\nbigr_X$ denote the sheaf of subalgebras 
in $D_{\nbigx/\cnum_{\lambda}}$
generated by $\lambda\,\Theta_{\nbigx/\cnum_{\lambda}}$.

Let $D$ be a hypersurface of $X$.
We set $\nbigd:=\cnum_{\lambda}\times D$.
We set 
$\nbigr_{X(\ast D)}:=
 \nbigr_X\otimes_{\nbigo_{\nbigx}} 
 \nbigo_{\nbigx}(\ast\nbigd)$,
where 
$\nbigo_{\nbigx}(\ast \nbigd)$ be the sheaf of
meromorphic functions on $\nbigx$
whose poles are contained in $\nbigd$.
\index{sheaf $\nbigr_{X(\ast D)}$}
The sheaf of algebras $\nbigr_{X(\ast D)}$ 
is Noetherian.
The notions of left and 
right $\nbigr_{X(\ast D)}$-modules
are naturally defined.
They are exchanged by a standard formalism
for $D$-modules or $\nbigr$-modules.
Namely, let 
$\omega_{\nbigx}:=
 \lambda^{-n}\cdot p_{\lambda}^{\ast}\omega_X$
as the subsheaf of 
$p_{\lambda}^{\ast}\omega_X\otimes
 \nbigo_{\nbigx}\bigl(\ast(\{0\}\times X)\bigr)$,
where $n=\dim X$.
Then, for any left $\nbigr_{X(\ast D)}$-module $\nbigm$,
we have the natural right 
$\nbigr_{X(\ast D)}$-module structure on
$\nbigm^r:=\omega_{\nbigx}\otimes_{\nbigo_{\nbigx}}\nbigm$.
We consider left $\nbigr_{X(\ast D)}$-modules
in this paper,
unless otherwise specified.
A left $\nbigr_{X(\ast D)}$-module $\nbigm$
is naturally regarded as a family of flat $\lambda$-connections.
Namely,
$\nbigm$ is equipped with a differential operator
$\DD:\nbigm\lrarr\nbigm\otimes\Omega^1_{\nbigx}$,
where $\Omega_{\nbigx}:=\lambda^{-1}p_{\lambda}^{\ast}\Omega_X^1$,
such that 
(i) $\DD(f\,s)=\lambda\,d(f)\,s+f\,\DD(s)$
for $f\in\nbigo_{\nbigx}(\ast\nbigd)$
and $s\in\nbigm$,
(ii) $\DD\circ\DD=0$.

An $\nbigr_{X(\ast D)}$-module is called strict,
if it is $\nbigo_{\cnum_{\lambda}}$-flat.
\index{strict}
If an $\nbigr_{X(\ast D)}$-module $M$
is (i) pseudo-coherent 
as an $\nbigo_{\nbigx(\ast \nbigd)}$-module,
(ii) locally finitely generated 
as an $\nbigr_{X(\ast D)}$-module,
then it is a coherent $\nbigr_{X(\ast D)}$-module.
The sheaf of algebras $\nbigr_{X(\ast D)}$ is 
naturally filtered by the order of differential operators,
and we have the notion of coherent filtration for 
$\nbigr_{X(\ast D)}$-modules
as in the case of $D$-modules.
Let $\nbigm$ be an $\nbigr_{X(\ast D)}$-module on
an open set  $U\subset\nbigx$.
We say that $\nbigm$ is {\em good},
if, for any $\lambda_0\in\cnum_{\lambda}$,
and for any compact subset $K$ of $X$
such that $\{\lambda_0\}\times K\subset U$,
there exist a neighbourhood $U'$ of 
$\{\lambda_0\}\times K$ in $U$,
and a finite filtration $F$ of $\nbigm_{|U'}$
such that
$\Gr^F(\nbigm_{|U'})$ has a coherent filtration.

The push forward and the pull back
of $\nbigr_{X(\ast D)}$-modules are defined 
by the formula of those for $\nbigr_X$-modules
in \cite{sabbah2}.
Let $f:X'\lrarr X$ be a morphism of complex manifolds.
Let $D$ be a hypersurface of $X$,
and we put $D':=f^{-1}(D)$.
Let $\nbigr_{X'\rarr X}:=
 \nbigo_{\nbigx'}\otimes_{f^{-1}\nbigo_{\nbigx}}
 f^{-1}\nbigr_X$,
which is naturally an
$\bigl(\nbigr_{X'},f^{-1}\nbigr_{X}\bigr)$-module.
\index{sheaf $\nbigr_{X'\rarr X}$}
We set
$\nbigr_{X\larr X'}:=
 \omega_{\nbigx'}\otimes
 \nbigr_{X'\rarr X}
 \otimes f^{-1}\omega_{\nbigx}$.
\index{sheaf $\nbigr_{X\larr X'}$}
For an $\nbigr_{X'(\ast D')}$-module $\nbigm'$,
we have $f_{\dagger}(\nbigm'):=
 Rf_{\ast}\bigl(\nbigr_{X\larr X'}
 \otimes_{\nbigr_{X'}}^L\nbigm'\bigr)$
in $D^b(\nbigr_{X(\ast D)})$.
\index{push-forward $f_{\dagger}\nbigm$}
Note that
\[
 \nbigr_{X\larr X'}
 \otimes_{\nbigr_{X'}}^L\nbigm'
\simeq
 \nbigr_{X\larr X'}(\ast D')
 \otimes_{\nbigr_{X'(\ast D')}}^L\nbigm',
\]
and that an $f^{-1}\nbigr_{X(\ast D)}$-injective 
resolution of 
$\nbigr_{X\larr X'}\otimes_{\nbigr_{X'}}^L\nbigm'$
is naturally an $f^{-1}\nbigr_X$-injective resolution.
Similarly,
for an $\nbigr_{X(\ast D)}$-module $\nbign$,
we have
$f^{\dagger}\nbign:=
 \nbigr_{X'\rarr X}\otimes^L_{f^{-1}\nbigr_{X}}
 f^{-1}\nbign$
in $D^b(\nbigr_{X'(\ast D')})$.
\index{pull back $f^{\dagger}\nbign$}
If $\nbigm'$ is good,
$f_{\dagger}\nbigm$ is cohomologically good,
which can be proved by the argument
in the case of $D$-modules.

\begin{lem}
\label{lem;11.4.1.1}
Assume that $f$ is proper and birational,
and it induces an isomorphism
$X'\setminus D'\simeq X\setminus D$.
Then, we have natural isomorphisms
$f_{\dagger}\nbigm'\simeq f_{\ast}\nbigm'$ 
and $f^{\dagger}\nbigm\simeq
 f^{\ast}\nbigm:=
 \nbigo_{\nbigx'}\otimes_{\nbigo_{\nbigx}}
 f^{-1}\nbigm$.
They give an equivalence of
the categories of good $\nbigr_{X(\ast D)}$-modules
and good $\nbigr_{X'(\ast D')}$-modules.
\end{lem}
\pf
Because
$\nbigr_{X'\rarr X}(\ast D')
\simeq\nbigr_{X'(\ast D')}$,
we obtain the first claim.
For an $\nbigr_{X(\ast D)}$-module 
$\nbigm$ of the form
$\nbigr_{X(\ast D)}
 \otimes_{\nbigo_{\nbigx}} 
 M$,
where $M$ is $\nbigo_{\nbigx}$-coherent,
we have
$f^{\dagger}\nbigm
\simeq
 \nbigr_{X'(\ast D')}
 \otimes_{\nbigo_{\nbigx'}}
 f^{\ast}M$.
For an $\nbigr_{X'(\ast D')}$-module 
$\nbigm'$ of the form
$\nbigr_{X'(\ast D')}
 \otimes_{\nbigo_{\nbigx'}} 
 M'$,
where $M'$ is $\nbigo_{\nbigx'}$-coherent,
we have
$f_{\dagger}\nbigm'\simeq
 \nbigr_{X(\ast D)}
 \otimes_{\nbigo_{\nbigx}}
 f_{\ast}(M')$.
Then, the second claim is clear.
\hfill\qed

\subsection{Strict specializability for $\nbigr$-modules}
\label{subsection;11.2.22.2}

\index{strict specializable}

The notion of strict specializability 
for $\nbigr_{X(\ast D)}$-module
is defined as in the case of $\nbigr_{X}$-modules
\cite{sabbah2}.
Let $\cnum_t$ be a complex line 
with a coordinate $t$.
Let $X_0$ be a complex manifold.
We put $X:=X_0\times \cnum_t$.
We identify $X_0$ and $X_0\times\{0\}$.
Let $D$ be a hypersurface of $X$.
(See also \S\ref{subsection;11.4.3.3}.)

Let $\Theta_{\nbigx/\cnum_{\lambda}}(\log \nbigx_0)$
denote the sheaf of logarithmic relative tangent sheaf
on $\nbigx$ over $\cnum_{\lambda}$.
We recall that $V_0\nbigr_X$ denotes the 
sheaf of subalgebras in $\nbigr_X$
generated by 
$\lambda\Theta_{\nbigx/\cnum_{\lambda}}(\log\nbigx_0)$.
Note that it depends only on $t$,
i.e.,
independent from the choice of a decomposition
into the product $X=X_0\times\cnum_t$.

\subsubsection{The case $X_0\not\subset D$}
We set $D_0:=X_0\cap D$.
We put 
$V_0\nbigr_{X(\ast D)}:=
 V_0\nbigr_X\otimes\nbigo_{\nbigx}(\ast\nbigd)$.
For any point $\lambda_0\in\cnum_{\lambda}$,
let $\nbigxzero$ denote a small neighbourhood 
of $\{\lambda_0\}\times X$.
We use the symbol $\nbigxzero_0$ in a similar meaning.
Let $\nbigm$ be a coherent $\nbigr_{X(\ast D)}$-module
on $\nbigx$.
It is called 
strictly specializable
along $t=0$ at $\lambda_0$,
if $\nbigm_{|\nbigxzero}$ is equipped with a filtration $\Vzero$
by coherent $V_0\nbigr_{X(\ast D)}$-modules
indexed by $\real$
satisfying Conditions 22.3.1 and 22.3.2
in \cite{mochi7}.
\index{filtration $\Vzero$}
Such a filtration is unique, if it exists.
Note that the condition is independent from
the choice of a decomposition into the product
$X=X_0\times\cnum_t$.
We say that $\nbigm$ is called strictly specializable
along $t$,
if it is strictly specializable along $t$
at any $\lambda_0$.
In this case,
we obtain 
an $\nbigr_{X_0(\ast D_0)}$-module
$\psi_{t,u}(\nbigm)$
as in the case of $\nbigr$-modules.
\index{nearby cycle sheaf $\psi_{t,u}(\nbigm)$}
Its push-forward by $X_0\lrarr X$
is also denoted by $\psi_{t,u}(\nbigm)$.
If $\nbigm_i$ $(i=1,2)$ are strictly specializable along $t$
with a morphism $F:\nbigm_1\lrarr\nbigm_2$,
then $F$ is compatible with the $V$-filtrations,
and we have an induced morphism
$\psi_{t,u}(F):\psi_{t,u}(\nbigm_1)\lrarr\psi_{t,u}(\nbigm_2)$.
If the cokernel of $\psi_{t,u}(F)$ is strict,
$F$ is called strictly specializable.
\index{strict specializable (morphism)}
In that case,
$F$ is strictly compatible with
the $V$-filtrations.
We have the naturally defined morphisms
$\deldel_t:\psi_{t,u}(\nbigm)\lrarr
 \psi_{t,u+\vecdelta}(\nbigm)$,
and 
$t:\psi_{t,u}(\nbigm)\lrarr
 \psi_{t,u-\vecdelta}(\nbigm)$,
where $\vecdelta=(1,0)\in\real\times\cnum$.
We set
$\psitilde_{t,u}(\nbigm):=
 \varinjlim\psi_{t,u-N\vecdelta}(\nbigm)$
for $u\in\real\times\cnum$.
\index{nearby cycle sheaf $\psitilde_{t,u}(\nbigm)$}
It is isomorphic to
$\psi_{t,u}(\nbigm(\ast t))$ below.

\subsubsection{The case $X_0\subset D$}
\label{subsection;13.5.4.31}

If $X_0\subset D$,
we decompose 
$D=X_0\cup D'$,
where $X_0\not\subset D'$.
We put $D_0':=D'\cap X_0$.
Let $\nbigm$ be a coherent
$\nbigr_{X(\ast D)}
=\nbigr_{X(\ast D')}(\ast t)$-module.
It is called strictly specializable along $t$ at $\lambda_0$
if $\nbigm_{|\nbigxzero}$ is equipped with a filtration $\Vzero$
by coherent $V_0\nbigr_{X(\ast D')}$-modules
indexed by $\real$,
satisfying the conditions in Definition 22.4.1
in \cite{mochi7}.
Similarly,
we obtain 
$\nbigr_{X_0(\ast D_0')}$-modules
$\psi_{t,u}(\nbigm)$.
The induced morphisms
$t:\psi_{t,u}(\nbigm)\lrarr\psi_{t,u-\vecdelta}(\nbigm)$
are isomorphisms.
It is also denoted by
$\psitilde_{t,u}(\nbigm)$.

\subsection{Some sheaves}

We recall some sheaves, following Sabbah
in \cite{sabbah4} and \cite{sabbah5},
to which we refer for more details and precision.
Let $X$ be an $n$-dimensional complex manifold.
Let $T$ be a real $C^{\infty}$-manifold.
Let $V$ be an open subset of $X\times T$.
Let $\nbige^{(n,n)}_{X\times T/T,c}(V)$
denote the space of $C^{\infty}$-sections
of $\Omega^{n,n}_{X\times T/T}$ on $V$
with compact supports.
\index{sheaf $\nbige^{(n,n)}_{X\times T/T,c}$}
Let $\Diff_{X\times T/T}(V)$ denote the space of
$C^{\infty}$-differential operators on $V$
relative to $T$,
i.e., we consider only differentials
in the $X$-direction on $V$.
For any compact subset $K\subset V$
and $P\in \Diff_{X\times T/T}(V)$,
we have the semi-norm
$\bigl\|\varphi\bigr\|_{P,K}=\sup_K|P\varphi|$.
For any closed subset $Z\subset X$,
let $\nbige^{<Z,(n,n)}_{X\times T/T,c}(V)$
denote the subspace of
$\nbige^{(n,n)}_{X\times T/T,c}(V)$,
which consists of the sections $\varphi$
such that $(P\varphi)_{|\widehat{(Z\times T)\cap V}}=0$ for
any $P\in\Diff_{X\times T/T}(V)$.
\index{sheaf $\nbige^{<Z,(n,n)}_{X\times T/T,c}$}
We have the induced semi-norms
$\|\cdot\|_{P,K}$ on the space 
$\nbige^{<Z,(n,n)}_{X\times T/T,c}(V)$.
By the semi-norms,
the spaces 
$\nbige^{(n,n)}_{X\times T/T,c}(V)$
and $\nbige^{<Z,(n,n)}_{X\times T/T,c}(V)$
are locally convex topological spaces.
Let $C^0_c(T)$ denote the space of 
continuous functions on $T$
with compact supports.
It is a normed vector space
with the sup norms.

The space of continuous $C^{0}(T)$-linear maps 
$\nbige^{(n,n)}_{X\times T/T,c}(V)\lrarr
C^0_c(T)$ is denoted by $\distribution_{X\times T/T}(V)$.
\index{sheaf $\distribution_{X\times T/T}$}
Let $D$ be a hypersurface of $X$.
Let $\distribution^{\moderate\,D}_{X\times T/T}(V)$
denote the space of
continuous $C^{\infty}(T)$-linear maps
from $\nbige^{<D,(n,n)}_{X\times T/T,c}(V)$
to $C^{\infty}_c(T)$.
\index{sheaf $\distribution^{\moderate\,D}_{X\times T/T}$}
Any elements of
$\distribution^{\moderate\,D}_{X\times T/T}(V)$
are called distributions with moderate growth
along $D$.
Let $\distribution_{X\times T/T,D}(V)$ be 
the space of continuous $C^{\infty}$-linear maps
$\Phi:\nbige^{(n,n)}_{X\times T/T,c}(V)\lrarr
 C^0_c(T)$
whose supports are contained in $D\times T$,
i.e., 
$\Phi(\varphi)=0$ if 
$\varphi=0$ on some neighbourhood of $D\times T$.
\index{sheaf $\distribution_{X\times T/T,D}$}
They give the sheaves
$\distribution_{X\times T/T}$,
$\distribution^{\moderate\,D}_{X\times T/T}$
and
$\distribution_{X\times T/T,D}$ on $X\times T$.
We have the following lemma
as in \cite{malgrange2} 
and \cite{mochi7}.

\begin{lem}
\label{lem;07.10.25.1}
We have the exact sequence
$0\lrarr \distribution_{X\times T/T,D}
\lrarr
 \distribution_{X\times T/T}
\lrarr
 \distribution^{\moderate\,D}_{X\times T/T}
\lrarr 0$.
In particular,
$\distribution^{\moderate\,D}_{X\times T/T}$
is isomorphic to the image of
$\distribution_{X\times T/T}
\lrarr j_{\ast}j^{-1}\distribution_{X\times T/T}$,
where $j:X\setminus D\lrarr X$.
\hfill\qed
\end{lem}

The following lemma can be proved
by using a standard argument in the theory of distributions.
Let $V$ be an open subset of $X$
with a real coordinate $(x_1,\ldots,x_{2n})$.
For a given $J=(j_1,\ldots,j_{2n})\in\seisuu_{\geq 0}^{2n}$,
we set $\del^J:=\prod_{i=1}^{2n}\del_{x_i}^{j_i}$
and $|J|=\sum j_i$.
\begin{lem}
Let $\Phi$ be an element of
$\distribution^{\moderate\,D}_{X\times T/T}(V)$.
Let $K$ be any compact region in $V$.
We have $m\in\seisuu_{>0}$
and $C>0$ such that
$ \sup_{T}\bigl|\Phi(\varphi)\bigr|
\leq
 C\,\sup_{|J|\leq m}|\del^J\varphi|_{K}$
for any 
$\varphi\in 
 \nbige^{<D\,(n,n)}_{X\times T/T,c}(V)$
with $\Supp(\varphi)\subset K$.
\hfill\qed
\end{lem}

\begin{cor}
\label{cor;08.12.24.1}
Let $g$ be a holomorphic function
such that $g^{-1}(0)\subset D$.
Let $\Phi\in
 \distribution^{\moderate\,D}_{X\times T/T}(V)$.
For any compact subset $K$ in $V$,
there exists $C_K>0$ such that
we have the well defined pairing
$\Phi\bigl(|g|^{s}\cdot \varphi\bigr)\in C^0_c(T)$
for $\varphi\in \nbige^{(n,n)}_{X\times T/T,c}(V)$
with $\Supp(\varphi)\subset K$
and for $\Re(s)>C_K$.
It depends on $s$ in a holomorphic way.
\hfill\qed
\end{cor}

Let $\nbigc^{\infty}_{X\times T}$
be the sheaf of $C^{\infty}$-functions on $X\times T$.
\index{sheaf $\nbigc^{\infty}_{X\times T}$}
Let $D$ be a hypersurface of $X$.
Let $U\subset X\times T$ be an open subset.
Let $\varphi$ be a $C^{\infty}$-function on 
an open subset $U\setminus (D\times T)$.
We say that $\varphi$ is a $C^{\infty}$-function
of moderate growth along $D$,
if the following holds
for any $(Q_1,Q_2)\in U\cap (D\times T)$.
\begin{itemize}
\item
For any $P\in\Diff_{X\times T/T}(U)$,
there exists $N>0$
such that 
$\bigl|P\varphi\bigr|
=O\bigl(|g|^{-N}\bigr)$ 
around $(Q_1,Q_2)$.
\end{itemize}
Let $\nbigc^{\infty\,\moderate D}_{X\times T}$
denote the sheaf of $C^{\infty}$-functions
on $X\times T$ with moderate growth along $D$.
It is standard to prove that
the sheaf is $c$-soft.
\index{sheaf $\nbigc^{\infty\,\moderate D}_{X\times T}$}

\vspace{.1in}

Let $f:X'\lrarr X$ be a morphism of complex manifolds.
Let $D$ be a hypersurface of $X$.
We put $D':=f^{-1}(D)$.

\begin{lem}
\label{lem;11.1.31.1}
Assume that $f$ is proper and birational.
We have a natural isomorphism
$a_1:f_{!}\distribution^{\moderate D'}_{X'\times T/T}
\simeq
 \distribution^{\moderate D}_{X\times T/T}$
and an epimorphism
$a_2:f^{-1}\nbigc^{\infty\,\moderate D}_{X\times T}
\lrarr
\nbigc^{\infty\,\moderate D'}_{X'\times T}$.
We also have 
$ f^{-1}\distribution^{\moderate D}_{X\times T/T}
\lrarr
\distribution^{\moderate D'}_{X'\times T/T}$
and 
$\nbigc^{\infty\,\moderate D}_{X\times T}
\simeq
f_{\ast}\nbigc^{\infty\,\moderate D'}_{X'\times T}$.
\end{lem}
\pf
The natural morphism
$\nbigc^{\infty\,\moderate D}(X\times T)
\lrarr
 \nbigc^{\infty\,\moderate\,D'}(X'\times T)$
is bijective.
Because 
$\nbigc^{\infty\,\moderate D'}_{X'\times T}$
is $c$-soft,
we obtain that $a_2$ is an epimorphism.
It is easy to prove that the induced morphism
$\nbigc^{\infty\,\moderate D}_{X\times T}
\lrarr
f_{\ast}\nbigc^{\infty\,\moderate D'}_{X'\times T}$
is an isomorphism.

Let $V$ be an open subset in $X\times T$,
and $V':=f^{-1}(V)$.
The natural continuous morphism
$\nbige^{(n,n)\,<D}_{X\times T/T,c}(V)
\lrarr
 \nbige^{(n,n)\,<D}_{X'\times T/T,c}(V')$
is bijective.
It is a homeomorphism,
which follows from Lemma \ref{lem;11.4.11.1} below.
Then, we obtain the second claim
for $\distribution^{\moderate D'}_{X'\times T/T}$
and 
$\distribution^{\moderate D}_{X\times T/T}$.
\hfill\qed

\begin{lem}
\label{lem;11.4.11.1}
Let $Y$ be a complex manifold
with a hypersurface $D_1$.
Let $g$ be a holomorphic function on $Y$
such that $g^{-1}(0)\subset D_1$.
Let $K\subset Y$ be a compact subset.
Let $P\in\Diff_{Y\times T/T}$
and $N\in\seisuu_{\geq 0}$.
Then, there exist $C>0$ and
$P_i\in\Diff_{Y\times T/T}$ $(i=1,\ldots,M)$
such that
$\sup_K|g^{-N}P\varphi|
\leq
 C\,\sum_i
 \sup_K|P_i\varphi|$
for any 
$\varphi\in
 \nbige^{(n,n)<D_1}_{Y\times T/T,c}(K\times T)$.
\end{lem}
\pf
We give only a sketch of the proof.
We have only to consider the case $P=1$.
We need estimates only around
any point of $K\times T$.
First, let us consider the case
$dg$ is nowhere vanishing on $X$.
We may assume that $X$ is equipped with
a holomorphic coordinate 
$(z_1,\ldots,z_n)$ and $z_1=g$.
Let $z_i=x_i+\sqrt{-1}y_i=r_i\,e^{\sqrt{-1}\theta_i}$.
By using the relation
$r_1^L\del_{r_1}^L=
 \sum_{q=0}^L (\!(L,q)\!)\,
 x_1^qy_1^{L-q}\,\del_{x_1}^q\del_{y_1}^{L-q}$,
where $(\!(L,q)\!)$ are binomial coefficients,
we can easily deduce
the desired estimate for
$\bigl|g^{-N}\varphi\bigr|$.
For example,
in the case of one variable and $N=1$,
we have
\[
 \bigl|z^{-1}f(z)\bigr|
\leq
 \int_0^1\bigl|(\del_rf)(sz)\bigr|ds
\leq
 \sup\bigl|
 \del_rf
 \bigr|
\leq 
 \sup\bigl|\del_xf\bigr|
+\sup\bigl|\del_yf\bigr|.
\]
If $g=\prod_{i=1}^{\ell}z_i^{m_i}$,
we set $g_j=\prod_{i=1}^jz_i^{m_i}$,
and we deduce 
$|g^{-N}\varphi|\leq
 C\,\sum_i \sup\bigl|g_{j-1}^{-N}Q_{j,i}\varphi\bigr|$
by an easy descending induction on $j$.
Let us consider the general case.
We take a projective birational morphism
$F:X_1\lrarr X$ such that
$F^{-1}(D)$ is normal crossing.
By the above consideration,
we have an estimate 
$\bigl|F^{-1}(g^{-N}\varphi)\bigr|
\leq
 C\,\sum\sup\bigl|P_i\,F^{-1}(\varphi)\bigr|$
on $X_1\times T$
for some $P_i\in\Diff_{X_1\times T/T}$.
It is easy to find $\Ptilde_{i,j}\in\Diff_{X\times T/T}$ 
$(j=1,\ldots,M_i)$
such that
$\bigl|P_i\,F^{-1}(\varphi)\bigr|
\leq
 \sum \sup \bigl|F^{-1}(\Ptilde_{i,j}\varphi)\bigr|$
for any $\varphi$.
Thus, we obtain the estimate on $X\times T$.
\hfill\qed

\begin{rem}
$\distribution^{\moderate D}_{X\times T/T}$
is also denoted by
$\distribution_{X(\ast D)\times T/T}$
or 
$\distribution_{X\times T/T}(\ast D)$
in the following.
\hfill\qed
\end{rem}
\index{sheaf $\distribution_{X(\ast D)\times T/T}$}
\index{sheaf $\distribution_{X\times T/T}(\ast D)$}

Suppose that $X=X_0\times\cnum_t$,
and $X_0\times\{0\}\not\subset D$.
We put $D_0:=(X_0\times\{0\})\cap D$.
For an open subset $V\subset X\times T$,
we set $V_0:=V\cap(X_0\times\{0\}\times T)$.
We have the naturally defined map
given by the restriction:
\begin{equation}
 \label{eq;13.5.8.10}
 \nbige^{(n,n)<D}_{X\times T/T,c}(V)
\lrarr
 \nbige^{(n-1,n-1)<D_0}_{X_0\times T/T,c}(V_0)
 \otimes(\cnum \,dt\,d\tbar)
\end{equation}

\begin{lem}
\label{lem;13.5.8.11}
The map
{\rm(\ref{eq;13.5.8.10})} is surjective.
\end{lem}
\pf
We set $H_1:=(X_0\times\{0\})\cup D$.
We take a projective birational map
$F:X'\lrarr X$ such that
(i) $H_1':=F^{-1}(H_1)$ is normal crossing,
(ii) $X'\setminus D'\simeq X\setminus D$,
where $D':=F^{-1}(D)$.
Let $X_1\subset X'$ be the strict transform of
$(X_0\times\{0\})$.
Let $V':=F^{-1}(V)$.
Let $V_0'$ be the closure of
$V_0\setminus D_0$ in $X'$.
Then, by a direct construction,
it is easy to see that
\[
\nbige^{(n,n)<D'}_{X'\times T/T,c}(V')
\lrarr
\nbige^{(n-1,n-1)<D_0'}_{X_1\times T/T,c}(V_0')
\otimes
(\cnum F^{\ast}(dt\,d\tbar))
\]
is surjective.
Then, we obtain that (\ref{eq;13.5.8.10}) is surjective.
\hfill\qed

\subsection{$\nbigr$-triple}
\label{subsection;13.5.4.40}

Let $X$ be a complex manifold,
and let $D$ be a hypersurface of $X$.
Let $\sigma:\cnum_{\lambda}^{\ast}\lrarr
 \cnum_{\lambda}^{\ast}$
be given by $\sigma(\lambda)=-\lambdabar^{-1}$.
We set $\vecS=
 \bigl\{\lambda\in\cnum\,\big| \,|\lambda|=1\bigr\}$.
\index{set $\vecS$}
\index{map $\sigma$}
If $\lambda\in\vecS$,
we have $\sigma(\lambda)=-\lambda$.
The induced maps
$\cnum_{\lambda}^{\ast}\times X
 \lrarr\cnum_{\lambda}^{\ast}\times X$
and 
$\vecS\times X\lrarr\vecS\times X$
are also denoted by $\sigma$.

We have the natural action of
$\nbigr_{X(\ast D)|\vecS\times X}$ on 
$\distribution^{\moderate D}_{\vecS\times X/\vecS}$.
We also have the action of
$\sigma^{\ast}\nbigr_{X(\ast D)|\vecS\times X}$
on $\distribution^{\moderate D}_{\vecS\times X/\vecS}$
by $\sigma^{\ast}(f)\bullet F=
 \overline{\sigma^{\ast}(f)}\cdot F$
and $\sigma^{\ast}(\deldel_i)\bullet F
=-\lambda^{-1}\delbar_iF$.
We set $\deldelbar_i:=-\lambda^{-1}\delbar_i$.

Let $\nbigm_i$ $(i=1,2)$ be
$\nbigr_{X(\ast D)}$-modules.
A hermitian sesqui-linear pairing
of $\nbigm_1$ and $\nbigm_2$
is an 
$\nbigr_{X(\ast D)|\vecS\times X}\otimes
 \sigma^{\ast}\nbigr_{X(\ast D)|\vecS\times X}$-morphism
$C:\nbigm_{1|\vecS\times X}
 \otimes\sigma^{\ast}\nbigm_{2|\vecS\times X}
\lrarr
 \distribution_{\vecS\times X/\vecS}^{\moderate\,D}$.
\index{hermitian sesqui-linear pairing}
Such a tuple $(\nbigm_1,\nbigm_2,C)$
is called an $\nbigr_{X(\ast D)}$-triple.
It is called good (coherent, strict, etc.),
if the underlying $\nbigr_{X(\ast D)}$-modules
are good (coherent, strict, etc.).
A morphism
$(\nbigm_1,\nbigm_2,C)\lrarr
 (\nbigm_1',\nbigm_2',C')$ is a pair of
morphisms
$\varphi_1:\nbigm_1'\lrarr\nbigm_1$
and $\varphi_2:\nbigm_2\lrarr\nbigm_2'$
such that
\[
C\bigl(
 \varphi_1(m_1'),\sigma^{\ast}(m_2)
 \bigr)
=C'\bigl(m_1',\sigma^{\ast}\varphi_2(m_2)\bigr).
\]
The category of
$\nbigr_{X(\ast D)}$-triples
is abelian,
and denoted by $\rtriplecat(X,D)$.
\index{$\nbigr_{X(\ast D)}$-triple}
\index{category $\rtriplecat(X,D)$}
For $\nbigt=(\nbigm_1,\nbigm_2,C)
 \in\rtriplecat(X,D)$,
we have the hermitian adjoint
$\nbigt^{\ast}:=(\nbigm_2,\nbigm_1,C^{\ast})$,
where $C^{\ast}(x,\sigma^{\ast}y):=
 \overline{\sigma^{\ast}C(y,\sigma^{\ast}x)}$.
\index{hermitian adjoint}

The push-forward for $\nbigr_{X(\ast D)}$-triples
is defined as in the case of $\nbigr$-triples
\cite{sabbah2}.
The following lemma follows from
Lemma \ref{lem;11.4.1.1} and
Lemma \ref{lem;11.1.31.1}.

\begin{lem}
\label{lem;11.2.21.1}
Let $f:X'\lrarr X$ be a proper birational morphism.
Let $D\subset X$ be a hypersurface,
and we put $D':=f^{-1}(D)$.
Assume that $f$ gives an isomorphism
$X'\setminus D'\simeq X\setminus D$.
Then, $f_{\dagger}$ induces
an equivalence of the categories of
good $\nbigr_{X'(\ast D')}$-triples
and good $\nbigr_{X(\ast D)}$-triples.
\hfill\qed
\end{lem}

Let us consider the specialization 
in the case $X=X_0\times\cnum_t$.
We set $D_0:=D\cap(X_0\times\{0\})$.
Let $\nbigt=(\nbigm',\nbigm'',C)$
be a coherent $\nbigr_{X(\ast D)}(\ast t)$-triple,
which is strictly specializable along $t$.
We have the induced 
$\nbigr_{X_0(\ast D_0)}$-triple
$\psi_{t,u}(\nbigt)$
for $u\in\real\times\cnum$
as in the case of $\nbigr$-triples.
(See \cite{sabbah2}. Note Lemma \ref{lem;13.5.8.11}.)
\index{nearby cycle functor $\psi_{t,u}(\nbigt)$}
The induced pairing is independent of 
the choice of a decomposition
into the product $X=X_0\times\cnum_t$.
It is also denoted by $\psitilde_{t,u}(\nbigt)$.
The pair of morphisms $(t^{-1},t)$ gives 
an isomorphism
$\psi_{t,u}(\nbigt)\simeq
 \psi_{t,u-\vecdelta}(\nbigt)$.
We have the induced morphisms
$-\del_tt:\psi_{t,u}(\nbigm')\,\lambda\lrarr
 \psi_{t,u}(\nbigm')$
and 
$-\del_tt:\psi_{t,u}(\nbigm'')\lrarr
 \psi_{t,u}(\nbigm'')\,\lambda^{-1}$.
The nilpotent part is denoted by $N'$
and $N''$.
Then, we obtain
$\nbign:=(N',N''):
 \psi_{t,u}(\nbigt)\lrarr
 \psi_{t,u}(\nbigt)\otimes\newTate(-1)$,
where
$\newTate(-1)
=\bigl(
 \nbigo_{\nbigx}\,\lambda,\,
 \nbigo_{\nbigx}\,\lambda^{-1},\,
 C_0
 \bigr)$
and $C_0(s_1,\sigma^{\ast}s_2)=s_1\cdot 
 \overline{\sigma^{\ast}s_2}$.

\begin{rem}
In {\rm\cite{sabbah2}},  {\rm\cite{mochi2}}
{\rm\cite{mochi8}},
we preferred to use $-\nbign$
as the canonically defined nilpotent map.
See Remark {\rm\ref{rem;13.8.1.1}}.
\hfill\qed
\end{rem}

Let $\nbigt$ be a coherent $\nbigr_{X(\ast D)}$-triple,
which is strictly specializable along $t$.
We naturally obtain a coherent
$\nbigr_{X(\ast D)}(\ast t)$-triple $\nbigt(\ast t)$,
which is strictly specializable along $t$.
We define
$\psitilde_{t,u}(\nbigt):=
 \psi_{t,u}(\nbigt(\ast t))$,
as in the case of $\nbigr$-triples.
\index{nearby cycle functor $\psitilde_{t,u}(\nbigt)$}

For a holomorphic function $g$,
we define $\psitilde_{g,u}$
by considering $\iota_{g\dagger}\nbigt$,
where $\iota_g:X\lrarr X\times\cnum_t$
is the graph,
as usual.

\vspace{.1in}
Let $j:\cnum_{\lambda}\lrarr\cnum_{\lambda}$
be given by $j(\lambda)=-\lambda$.
\index{map $j$}
The induced morphism
$\nbigx\lrarr\nbigx$ is also denoted by $j$.
Let $\nbigt=(\nbigm_1,\nbigm_2,C)
 \in\rtriplecat(X,D)$.
We have the naturally defined pairing
$j^{\ast}C:j^{\ast}\nbigm_{1|\vecS\times X}
\times\sigma^{\ast}j^{\ast}\nbigm_{2|\vecS\times X}
\lrarr
 \distribution_{\vecS\times X/\vecS}^{\moderate D}$
by which we have
$j^{\ast}\nbigt:=
 (j^{\ast}\nbigm_1,j^{\ast}\nbigm_2,j^{\ast}C)
\in\rtriplecat(X,D)$.
We have a natural identifications
$j^{\ast}(\nbigt^{\ast})
=(j^{\ast}\nbigt)^{\ast}$.

\subsection{Integrable $\nbigr$-triple}

\index{integrable $\nbigr$-triple}

Let us recall the concepts of integrable $\nbigr$-module
and integrable $\nbigr$-triple,
introduced in Section 7 of \cite{sabbah2}.
Let $D_{\nbigx}$ denote the sheaf of holomorphic differential 
operators on $\nbigx$,
and we set $D_{\nbigx(\ast\nbigd)}:=
D_{\nbigx}\otimes\nbigo_{\nbigx}(\ast\nbigd)$.
Let $\nbigr_{X(\ast D)}\langle\lambda^2\del_{\lambda}\rangle
\subset D_{\nbigx}$ be the sheaf of subalgebras 
generated by $\nbigr_X$ and $\lambda^2\del_{\lambda}$.
A module over $\nbigr_{X(\ast D)}\langle\lambda^2\del_{\lambda}\rangle$
is called an integrable $\nbigr_{X(\ast D)}$-module.
For integrable $\nbigm_i$ $(i=1,2)$,
an $\nbigr_X$-homomorphism
$\nbigm_1\lrarr\nbigm_2$ is called integrable,
if it gives an 
$\nbigr_X\langle\lambda^2\del_{\lambda}\rangle$-homomorphism.

Let $\nbigm_i$ $(i=1,2)$ be integrable $\nbigr_{X(\ast D)}$-modules.
A hermitian pairing $C$ of $\nbigm_i$ is called integrable,
if the following holds:
\[
\del_{\theta}
C\bigl(m_1,\overline{\sigma^{\ast}m_2}\bigr)
=
 C\bigl(
 \del_{\theta}m_1,\overline{\sigma^{\ast}m_2}
 \bigr)
+C\bigl(
 m_1,\overline{\sigma^{\ast}(\del_{\theta}m_2)}
 \bigr)
\]
Here, $\del_{\theta}m_i:=
 i\lambda\del_{\lambda}m_i
=\bigl(
 i\lambda\del_{\lambda}-i\lambdabar\del_{\lambdabar}
\bigr)m_i$.
If $C$ is extended to a pairing
$\nbigm_{1|\cnum_{\lambda}^{\ast}\times X}
 \times
 \sigma^{\ast}\nbigm_{2|\cnum_{\lambda}^{\ast}\times X}
\lrarr
\distribution_{\cnum_{\lambda}^{\ast}\times X
 /\cnum_{\lambda}^{\ast}}^{an\,\moderate D}$,
where $\distribution_{\cnum_{\lambda}^{\ast}\times
 X/\cnum_{\lambda}^{\ast}}
^{an,\,\moderate D}$
is the sheaf of $\lambda$-holomorphic distributions on
$\cnum_{\lambda}^{\ast}\times X$ with moderate growth along $D$,
then we have
\[
 \lambda\del_{\lambda}C\bigl(
 m_1,\overline{\sigma^{\ast}(m_2)}
 \bigr)
=C\bigl(
\lambda\del_{\lambda}m_1,\overline{\sigma^{\ast}(m_2)}
 \bigr)
-C\bigl(
m_1,\overline{\sigma^{\ast}(\lambda\del_{\lambda}m_2)}
 \bigr).
\]
A tuple of integrable $\nbigr_{X(\ast D)}$-modules
and an integrable hermitian sesqui-linear pairing is called
an integrable $\nbigr_{X(\ast D)}$-triple.
If $\nbigt$ is an integrable $\nbigr_{X(\ast D)}$-triple,
$\nbigt^{\ast}$
and $j^{\ast}\nbigt$ are also naturally integrable.
For integrable $\nbigt_i$ $(i=1,2)$,
a morphism $\nbigt_1\lrarr\nbigt_2$
in $\rtriplecat(X,D)$ is called integrable,
if the underlying morphisms of $\nbigr_{X(\ast D)}$-modules
are integrable.

Let $f:X'\lrarr X$ be a morphism of complex manifolds.
Let $D$ be a hypersurface of $X$,
and we put $D':=f^{-1}(D)$.
Let $\nbigt'$ be a good integrable $\nbigr_{X'(\ast D')}$-module.
We have the push-forward
$f^i_{\dagger}\nbigt'$ as $\nbigr_{X(\ast D)}$-triple.
As observed in \cite{sabbah2},
$f^i_{\dagger}\nbigt'$ is naturally an integrable
$\nbigr_{X(\ast D)}$-triple.
Similarly,
for an integrable $\nbigr_{X(\ast D)}$-triple $\nbigt$ on $X$,
the pull back
$f^{\dagger}\nbigt$ is naturally integrable.
We have the integrable variant of Lemma \ref{lem;11.2.21.1}.
\begin{lem}
\label{lem;11.2.22.1}
Assume that $f$ is proper
and induces an isomorphism
$X'\setminus D'\simeq X\setminus D$.
Then, $f_{\dagger}$ induces
an equivalence of the categories of
good integrable $\nbigr_{X'(\ast D')}$-triples
and good integrable $\nbigr_{X(\ast D)}$-triples.
\hfill\qed
\end{lem}

We recall the following in \cite{sabbah2}.
\begin{lem}
Let $X$ and $D$ be as in 
{\rm\S\ref{subsection;11.2.22.2}}.
Let $\nbigm$ be a coherent
$\nbigr_{X(\ast D)}$-module,
which is strictly specializable along $t$.
If $\nbigm$ is integrable,
the $V$-filtrations are preserved by 
$\lambda^2\del_{\lambda}$.
In particular,
the nearby cycle sheaves
$\psi_{t,u}(\nbigm)$
are naturally integrable.
\hfill\qed
\end{lem}

\begin{cor}
Let $X$ and $D$ be as in 
{\rm\S\ref{subsection;11.2.22.2}}.
Let $\nbigt$ be a coherent $\nbigr_X(\ast D)$-triple
which is strictly specializable along $t$.
Then, 
$\psi_{t,u}(\nbigt)$ are also integrable.
\hfill\qed
\end{cor}

\subsection{Smooth $\nbigr$-triple and some functoriality}
\label{subsection;10.12.27.1}

\index{smooth $\nbigr$-triple}

Let $X$ be a complex manifold
with a normal crossing hypersurface $D$.
An $\nbigr_{X(\ast D)}$-module $\nbigm$ is called smooth,
if it is a locally free $\nbigo_{\nbigx}(\ast\nbigd)$-module.
\index{smooth $\nbigr_{X(\ast D)}$-module}
A smooth hermitian sesqui-linear pairing of
smooth $\nbigr_X$-modules $\nbigm_i$ $(i=1,2)$
is an
$\nbigr_{X(\ast D)|\vecS\times X}\otimes
 \sigma^{\ast}
 \nbigr_{X(\ast D)|\vecS\times X}$-homomorphism
$C:\nbigm_{1|\vecS\times X}
 \otimes
 \sigma^{\ast}\nbigm_{2|\vecS\times X}
\lrarr
 \nbigc^{\infty\,\moderate \,D}_{\vecS\times X}$
such that
its restriction to $\vecS\times (X\setminus D)$
is extended to 
a pairing of
$\nbigm_{1|\cnum_{\lambda}^{\ast}\times (X\setminus D)}$
and 
$\sigma^{\ast}\!
 \nbigm_{1|\cnum_{\lambda}^{\ast}\times (X\setminus D)}$
to 
the sheaf of $C^{\infty}$-functions on 
$\cnum_{\lambda}^{\ast}\times (X\setminus D)$
which are $\lambda$-holomorphic.
Let $\rtriplecat_{\sm}(X,D)$ denote
the full subcategory of smooth $\nbigr_{X(\ast D)}$-triples
in $\rtriplecat(X,D)$.
\index{category $\rtriplecat_{\sm}(X,D)$}

\subsubsection{Tensor}
\index{tensor}

Let us consider
$\nbigt=(\nbigm_1,\nbigm_2,C)\in\rtriplecat_{\sm}(X,D)$
and $\nbigt'=(\nbigm_1',\nbigm_2',C')\in
 \rtriplecat(X,D)$.
We set
$\nbigm''_i:=\nbigm_i
 \otimes_{\nbigo_{\nbigx}}\nbigm_i'$.
We have the naturally defined pairing
$C''$ of $\nbigm''_i$ $(i=1,2)$,
given as follows:
\[
 C''\bigl(m_1\otimes m_1',
 \sigma^{\ast}(m_2\otimes m_2')\bigr)
=C(m_1,\sigma^{\ast}m_2)
 \cdot C(m_1',\sigma^{\ast}m_2')
\]
The pairing $C''$ is also denoted by $C\otimes C'$.
The $\nbigr_{X(\ast D)}$-triple 
$(\nbigm_1'',\nbigm_2'',C'')$
is denoted by $\nbigt\otimes\nbigt'$.
If $\nbigt'$ is smooth,
$\nbigt\otimes\nbigt'$ is also smooth.
If $\nbigt$ and $\nbigt'$ are integrable,
$\nbigt\otimes\nbigt'$ is also integrable.

\subsubsection{Pull back}
\index{pull back}

Let $f:X_1\lrarr X_2$ be a morphism of
complex manifolds.
Let $D_i$ be hypersurfaces of $X_i$
such that $D_1=f^{-1}(D_2)$.
Let $\nbigt=(\nbigm_1,\nbigm_2,C)
\in\rtriplecat_{\sm}(X_2,D_2)$.
We set
$f^{\ast}\nbigm_i:=
 \nbigo_{\nbigx_2}
 \otimes_{f^{-1}\nbigo_{\nbigx_1}}
 f^{-1}\nbigm_i$,
which are naturally smooth $\nbigr_{X_1(\ast D_1)}$-modules.
We have a pairing $f^{-1}C$ of
$f^{-1}\nbigm_1$ and $f^{-1}\nbigm_2$
obtained as the composition of the following:
\[
 f^{-1}(\nbigm_1)_{|\vecS\times X_1}
\times
 \sigma^{\ast}f^{-1}(\nbigm_2)_{|\vecS\times X_1}
\lrarr
 f^{-1}\nbigc^{\infty\,\moderate D_2}_{\vecS\times X_2}
\lrarr
 \nbigc^{\infty\,\moderate D_1}_{\vecS\times X_1}
\]
It induces a hermitian sesqui-linear
pairing $f^{\ast}C$
of $f^{\ast}\nbigm_1$ and $f^{\ast}\nbigm_2$.
Thus, we obtain a smooth 
$\nbigr_{X_1(\ast D_1)}$-triple
$f^{\ast}\nbigt:=
 (f^{\ast}\nbigm_1,f^{\ast}\nbigm_2,f^{\ast} C)$.
If $\nbigt$ is integrable,
$f^{\ast}\nbigt$ is also integrable.

\begin{lem}
Assume that $f$ is birational proper morphism,
and it gives an 
isomorphism $X_1\setminus D_1\simeq X_2\setminus D_2$.
Then, $f^{\ast}$ gives an equivalence
of the categories
$\rtriplecat_{\sm}(X_i,D_i)$.
\hfill\qed
\end{lem}

\subsubsection{A projection formula}

\index{projection formula}

Let $f:X_1\lrarr X_2$ be as above.
Let $\nbigt'$ be a good $\nbigr_{X_1(\ast D_1)}$-triple.
Assume that the support of $\nbigt'$ is proper over $X_2$.
Let $\nbigt\in\rtriplecat_{sm}(X_2,D_2)$.
\begin{lem}
\label{lem;10.12.24.3}
We have a natural isomorphism
$f_{\dagger}\bigl(
 f^{\ast}(\nbigt)\otimes\nbigt'
 \bigr)
\simeq
 \nbigt\otimes f_{\dagger}\nbigt'$.
\end{lem}
\pf
Let us construct an isomorphism
in the level of $\nbigr$-modules.
We have only to consider the case
$\nbigm'=
M'\otimes\nbigr_{X_1}\otimes\omega_{\nbigx_1}^{-1}$,
where $M'$ is a coherent $\nbigo_{\nbigx_1}$-module. 
We will freely use an isomorphism
in Lemma \ref{lem;11.1.14.20} below.
Let $\nbigm$ be 
a smooth $\nbigr_{X_2}(\ast H_2)$-module.
Then, we have the following isomorphism:
\begin{multline}
 f_{\dagger}\bigl(
 f^{\ast}\nbigm\otimes\nbigm'
 \bigr)
\simeq
 \nbigr_{X_2}\otimes\omega_{\nbigx_2}^{-1}
\otimes
 f_{\ast}\bigl(
 M'\otimes \Forget(f^{\ast}\nbigm)
 \bigr) \\
\simeq
 \bigl(
 \nbigr_{X_2}\otimes\omega_{\nbigx_2}^{-1}
 \bigr)
\otimes_{\nbigo_{\nbigx_2}}
 \bigl(
 f_{\ast}(M')\otimes\Forget(\nbigm)
 \bigr) \\
 \simeq
\nbigm\otimes_{\nbigo_{\nbigx_2}}
\bigl(
 f_{\ast}(M')\otimes
 \nbigr_{X_2}\otimes\omega_{\nbigx_2}^{-1}
\bigr)
\simeq
 \nbigm\otimes
 f_{\dagger}(\nbigm')
\end{multline}

If $\nbigm$ is isomorphic to
$\nbigo_{\nbigx_2}^{\oplus \,r}$
as an $\nbigr_{X_2}$-module,
then the isomorphism is equal to
the natural one
$f_{\dagger}\bigl(
 \nbigm^{\prime\,\oplus\,r}
 \bigr)
\simeq
 f_{\dagger}(\nbigm')^{\oplus\,r}$.
We can find such an isomorphism locally
around any point of
$\vecS\times (X_2\setminus D_2)$.
Then, we can easily compare
the pairings
$f_{\dagger}\bigl(
 f^{\ast}(C)\otimes C'
 \bigr)$
and $C\otimes f_{\dagger}C'$.
\hfill\qed

\subsubsection{Appendix}

Let $X$ be a complex manifold.
Let $N$ be an $\nbigo_{\nbigx}$-module.
Let $\nbigm$ be an $\nbigr_X$-module.
The underlying $\nbigo_{\nbigx}$-module of $\nbigm$
is denoted by $\Forget(\nbigm)$.
We obtain an $\nbigr_X$-module
on $\Bigl(
 \bigl(
 N\otimes_{\nbigo_{\nbigx}} 
\nbigr_X\bigr)
 \otimes_{\nbigo_{\nbigx}}\omega_{\nbigx}^{-1}
 \Bigr)\otimes_{\nbigo_{\nbigx}}\nbigm$,
where we use the left multiplication of $\nbigo_{\nbigx}$
on $\nbigr_X$ for the tensor product of
$N$ and $\nbigr_X$.
We also have an $\nbigr_X$-module
$\Bigl[
 \bigl(
 N\otimes_{\nbigo_{\nbigx}} \Forget(\nbigm)\bigr)
\otimes_{\nbigo_{\nbigx}} \nbigr_X
 \Bigr]
 \otimes_{\nbigo_{\nbigx}}\omega_{\nbigx}^{-1}$,
where we use the left multiplication of $\nbigo_{\nbigx}$
on $\nbigr_X$ for the tensor product
in the bracket.
Let us recall the following isomorphism
from \cite{sabbah2} and \cite{saito1}.
\begin{lem}
\label{lem;11.1.14.20}
We have a natural $\nbigr_X$-isomorphism
\[
\bigl(
 N\otimes \Forget(\nbigm)\bigr)
\otimes \nbigr_X\otimes\omega_{\nbigx}^{-1}
\simeq
 \bigl(N\otimes\nbigr_X\otimes\omega_{\nbigx}^{-1}
 \bigr)
 \otimes\nbigm.
\]
\end{lem}
\pf
A natural $\nbigo_{\nbigx}$-morphism
$N\otimes\Forget(\nbigm)\otimes\omega_{\nbigx}^{-1}
\lrarr
 \bigl(
 N\otimes\nbigr_X\otimes\omega_{\nbigx}^{-1}\bigr)
 \otimes \nbigm$
is naturally extended to an $\nbigr_X$-morphism
$\bigl(
N\otimes\Forget(\nbigm)\otimes
 \nbigr_X
\bigr)\otimes\omega_{\nbigx}^{-1}
\lrarr
 \bigl(
 N\otimes\nbigr_X\otimes\omega_{\nbigx}^{-1}\bigr)
 \otimes \nbigm$,
which is an isomorphism.
\hfill\qed

\subsection{Variation of twistor structure}
\label{subsection;13.5.3.10}

\index{variation of twistor structure}

%Let $\nbigt=(\nbigm_1,\nbigm_2,C)\in
% \rtriplecat_{sm}(X,D)$.
%A smooth $\nbigr_{X(\ast D)}$-triple
%$\nbigt=(\nbigm_1,\nbigm_2,C)$
%An object $\nbigt=(\nbigm_1,\nbigm_2,C)\in
% \rtriplecat_{sm}(X,D)$

A variation of twistor structure on $(X,D)$
(or simply a twistor structure on $(X,D)$)
is a smooth $\nbigr_{X(\ast D)}$-triple
$\nbigt=(\nbigm_1,\nbigm_2,C)$
satisfying the following conditions:
%is called a variation of twistor structure on $(X,D)$,
%(or simply a twistor structure on $(X,D)$),
%if the following holds:
\begin{itemize}
\item If $X$ is a point and $D=\emptyset$,
the extended pairing
$\nbigm'_{|\cnum_{\lambda}^{\ast}}
\otimes
 \sigma^{\ast}\nbigm''_{|\cnum_{\lambda}^{\ast}}
\lrarr
 \nbigo_{\cnum_{\lambda}^{\ast}}$
is perfect,
i.e., the induced morphism
\[
 \nbigm'_{|\cnum_{\lambda}^{\ast}}
\lrarr
 \nhom_{\nbigo_{\cnum_{\lambda}^{\ast}}}
 \bigl(
 \sigma^{\ast}\nbigm''_{|\cnum_{\lambda}^{\ast}},\,
 \nbigo_{\cnum_{\lambda}^{\ast}}
 \bigr)
\]
is an isomorphism. 
\item In the general case,
$\iota_{P}^{\ast}C$ is a twistor structure
for any point $\iota_P:\{P\}\lrarr X\setminus D$.
\end{itemize}
Let $\VTS(X,D)\subset
\rtriplecat_{\sm}(X,D)$
denote the full subcategory of twistor structures on $(X,D)$.
\index{category $\VTS(X,D)$}
For $\nbigt_i\in\VTS(X,D)$,
we have $\nbigt_1\otimes\nbigt_2\in\VTS(X,D)$.
For $f:(X_1,D_1)\lrarr (X_2,D_2)$,
we naturally have the pull back
$f^{\ast}:\VTS(X_2,D_2)\lrarr\VTS(X_1,D_1)$.

\subsubsection{Dual}

Let $\nbigt=(\nbigm_1,\nbigm_2,C)$
be a variation of twistor structure on $(X,D)$.
We set $\nbigm_i^{\lor}:=
 \nhom_{\nbigo_{\nbigx}(\ast\nbigd)}
 \bigl(\nbigm_i,\nbigo_{\nbigx}(\ast\nbigd)\bigr)$.
\index{dual $\nbigm^{\lor}$}
We shall observe that
we have a naturally induced pairing
\[
C^{\lor}:
 \nbigm^{\lor}_{1|\vecS\times X}
\times
 \sigma^{\ast}
 \nbigm^{\lor}_{2|\vecS\times X}
\lrarr
 \nbigc^{\infty\,\moderate D}_{\vecS\times X},
\]
and we shall set
$\nbigt^{\lor}:=
(\nbigm_1^{\lor},\nbigm_2^{\lor},C^{\lor})$,
which is called the dual of $\nbigt$.
\index{dual $\nbigt^{\lor}$}
For a neighbourhood  $U(\lambda_0)$ of $\lambda_0$
in $\cnum_{\lambda}$,
let $\vecI(\lambda_0):=\vecS\cap U(\lambda_0)$.

\vspace{.1in}
First, let us consider the case $D=\emptyset$.
We set
$\nbigm'_{C^{\infty},\vecS\times X}:=
 \nbigm'\otimes\nbigc^{\infty}_{\vecS\times X}$
and
$\nbigm''_{C^{\infty},\vecS\times X}:=
 \nbigm''\otimes\nbigc^{\infty}_{\vecS\times X}$
on $\vecS\times X$.
Then, $\nbigm'_{C^{\infty},\vecS\times X}$
and $\nbigm''_{C^{\infty},\vecS\times X}$
are naturally $\Diff_{\vecS\times X/\vecS}$-modules.
For a section $P\in \Diff_{\vecS\times X/\vecS}$ 
on $\vecI(\lambda_0)\times X$
and $m''\in \nbigm''_{C^{\infty},\vecS\times X}$ on
$\vecI(-\lambda_0)\times X$,
we set 
$P\cdot \sigma^{\ast}m'':=
 \sigma^{\ast}\bigl(
 \overline{\sigma^{\ast}(P)}\cdot m''
 \bigr)$.
Thus, $\sigma^{\ast}(\nbigm''_{C^{\infty},\vecS\times X})$
is a $\Diff_{\vecS\times X/\vecS}$-module.

\vspace{.1in}

We have the induced pairing of
$\nbigc_{\vecS\times X}^{\infty}$-modules
\begin{equation}
\label{eq;10.12.24.10}
 C:\nbigm'_{C^{\infty},\vecS\times X}
\times
 \sigma^{\ast}\nbigm''_{C^{\infty},\vecS\times X}
\lrarr
 \nbigc^{\infty}_{\vecS\times X}.
\end{equation}
We have induced morphisms of
$\nbigc_{\vecS\times X}^{\infty}$-modules
\[
 \Psi_C:\nbigm'_{C^{\infty},\vecS\times X}
\lrarr
 \sigma^{\ast}
 \nbigm^{\prime\prime\lor}_{C^{\infty},\vecS\times X},
\quad\quad
 \Phi_C:
 \sigma^{\ast}\nbigm''_{C^{\infty},\vecS\times X}
\lrarr
 \nbigm^{\prime\lor}_{C^{\infty},\vecS\times X}
\]
given as follows:
\begin{equation}
 \label{eq;10.12.25.1}
 \bigl\langle
 \Psi_C(m'),\sigma^{\ast}(m'')
 \bigr\rangle
:=C(m',\sigma^{\ast}m''),
\quad\quad
 \bigl\langle
 \Phi_C(\sigma^{\ast}m''),
 m'
 \bigr\rangle
:=C(m',\sigma^{\ast}m'')
\end{equation}
Note that we use the identification
$\sigma^{\ast}(\nbigm^{\prime\prime\lor})
\simeq
 \sigma^{\ast}(\nbigm'')^{\lor}$
given by the pairing
$\sigma^{\ast}
 \nbigm^{\prime\prime\lor}_{C^{\infty},\vecS\times X}
 \times
 \sigma^{\ast}\nbigm''_{C^{\infty},\vecS\times X}
\lrarr
 \nbigc^{\infty}_{\vecS\times X}$:
\[
 \langle
 \sigma^{\ast}n'',\sigma^{\ast}m''
 \rangle
=\overline{\sigma^{\ast}
 \langle n'',m''\rangle}.
\]
We can check that
$\Psi_C$ and $\Phi_C$
are $\Diff_{\vecS\times X/\vecS}$-homomorphisms
by direct computations.
Because $\nbigt\in\VTS(X,D)$,
$\Psi_C$ and $\Phi_C$ are isomorphisms.
Hence, we have the induced 
$\Diff_{\vecS\times X/\vecS}$-morphism:
\[
 C^{\lor}:
 \nbigm^{\prime\lor}_{C^{\infty},\vecS\times X}
\otimes
 \sigma^{\ast}\nbigm^{\prime\prime\,\lor}
 _{C^{\infty},\vecS\times X}
\lrarr
 \nbigc^{\infty}_{\vecS\times X},
\,\,\,\,
 C^{\lor}(n',\sigma^{\ast}n'')
:=C\bigl(
 \Psi^{-1}_C(\sigma^{\ast}n''),
 \Phi^{-1}_C(n')
 \bigr)
\]
The composition
$\nbigm^{\prime\,\lor}_{|\vecS\times X}
\times
 \sigma^{\ast}
 \nbigm^{\prime\prime\lor}_{|\vecS\times X}
\lrarr
 \nbigm^{\prime\,\lor}_{C^{\infty},\vecS\times X}
\times
 \sigma^{\ast}
 \nbigm^{\prime\prime\lor}_{C^{\infty},\vecS\times X}
\lrarr
 \nbigc^{\infty}_{\vecS\times X}$
is also denoted by $C^{\lor}$,
which is the desired smooth hermitian pairing.

\vspace{.1in}

Let us consider the case that
$D$ is not necessarily empty.
By the previous consideration,
we have the smooth pairing
$C_0^{\lor}:
\nbigm'_{|\vecS\times(X\setminus D)}
\times
\sigma^{\ast}\nbigm''_{|\vecS\times(X\setminus D)}
\lrarr
 \nbigc_{\vecS\times(X\setminus D)}^{\infty}$.
\begin{lem}
\label{lem;10.12.24.20}
$C_0^{\lor}$ is extended to
a pairing
$C^{\lor}:
 \nbigm'_{|\vecS\times X}
\times
\sigma^{\ast}\nbigm''_{|\vecS\times X}
\lrarr
 \nbigc_{\vecS\times X}^{\infty\,\moderate D}$.
\end{lem}
\pf
Let us consider the case $\rank\nbigm'=\rank\nbigm''=1$.
Let $P\in D$.
We will shrink $X$ around $P$.
We may assume that the monodromy of
$\nbigm'$ and $\nbigm''$ are trivial.
Let $\lambda_0\in\vecS$.
Let $U(\lambda_0)$ 
be a small neighbourhood of $\lambda_0$
in $\cnum$,
and let $U(-\lambda_0):=\sigma(U(\lambda_0))$.
We can find a frame $m'$
of $\nbigm'$ on $U(\lambda_0)\times X$,
and a frame $m''$ of $\nbigm''$
on $U(-\lambda_0)\times X$.
We have meromorphic functions $\gminia'$
and $\gminia''$
such that
(i) the poles are contained in $\nbigd$,
(ii) $\DD m'=m'\,d\gminia'$
and $\DD m''=m'' d\gminia''$.
\begin{lem}
\label{lem;10.12.24.21}
We have $\gminia'=\sigma^{\ast}\gminia''$
modulo holomorphic functions.
\end{lem}
\pf
We set $F=C(m',\sigma^{\ast}m'')$.
We obtain 
$F=G(\lambda)\,\exp
 \Bigl(2\sqrt{-1}\Image(\lambda^{-1}\gminia')
+\lambda
 \overline{\bigl(\gminia'-\sigma^{\ast}\gminia''\bigr)}
 \Bigr)$.
Because $F$ is of moderate growth,
we obtain the claim of Lemma \ref{lem;10.12.24.21}.
\hfill\qed

\vspace{.1in}
In particular,
$\gminia'$ and $\gminia''$
are independent of $\lambda$.
By considering the tensor product
with the rank one object,
we may assume that
$\gminia'=\gminia''=0$.
Then, $C_0^{\lor}$ is clearly extended.
Moreover, we have 
a bound of $C(m',\sigma^{\ast}m'')$
from below.
Then, we can deduce 
the claim of Lemma \ref{lem;10.12.24.20}
by using the exterior product.
\hfill\qed

\vspace{.1in}
The following lemma is easy to see.
\begin{lem}
If $\nbigt\in\VTS(X,D)$ is integrable,
$\nbigt^{\lor}$ is naturally integrable.
\hfill\qed
\end{lem}

\subsubsection{Real structure}

\index{real structure (smooth $\nbigr$-triple)}

When $C$ is non-degenerate,
we have natural identifications
$(\nbigt^{\ast})^{\lor}
=(\nbigt^{\lor})^{\ast}$
and
$j^{\ast}(\nbigt^{\lor})
=(j^{\ast}\nbigt)^{\lor}$.
We formally set
$\gammatilde_{\sm}^{\ast}\nbigt:=
 j^{\ast}\bigl(\nbigt^{\lor}\bigr)^{\ast}$.
\index{functor $\gammatilde_{sm}^{\ast}\nbigt$}
We naturally have 
$\gammatilde_{\sm}^{\ast}\circ\gammatilde_{\sm}^{\ast}
 (\nbigt)\simeq\nbigt$.
A real structure of $\nbigt$ is 
defined to be an isomorphism
$\kappa:\gammatilde_{\sm}^{\ast}\nbigt\lrarr\nbigt$
such that
$\kappa\circ\gammatilde_{\sm}^{\ast}\kappa=\id$.

\begin{rem}
$\gammatilde_{\sm}^{\ast}$ is the counterpart
of the pull back of vector bundles on $\proj^1$ by
$\gamma:\proj^1\lrarr\proj^1$,
where $\gamma([z_0:z_1])=[\zbar_1:\zbar_0]$.
\hfill\qed
\end{rem}

Let $f:X_1\lrarr X_2$ be a morphism
of complex manifolds.
Let $D_i$ $(i=1,2)$ be hypersurfaces of $X_i$
such that $D_1=f^{-1}(D_2)$.
Let $\nbigt\in\VTS(X_2,D_2)$.
We have natural isomorphisms
$f^{\ast}(\nbigt^{\ast})
=(f^{\ast}\nbigt)^{\ast}$,
$f^{\ast}(\nbigt^{\lor})
\simeq
 f^{\ast}(\nbigt)^{\lor}$,
$f^{\ast}j^{\ast}(\nbigt)\simeq j^{\ast}f^{\ast}(\nbigt)$.
A real structure is functorial
with respect to the pull back.

\begin{rem}
Let $\nbigt\in\VTS(X,D)$.
Let $\nbigm_i$ be the underlying $\nbigr_{X(\ast D)}$-modules.
If $\nbigt\in\VTS(X,D)$ is equipped with a real structure,
each $\nbigmlambda_{i|P}$ 
for $(\lambda,P)\in\vecS\times(X\setminus D)$
has an induced real structure.
If moreover $\nbigt$ is equipped with an integrable structure,
the flat bundle
$\nbigm_{i|\cnum_{\lambda}^{\ast}\times(X\setminus D)}$
has an induced real structure.
\hfill\qed
\end{rem}

\subsection{Tate object}
\label{subsection;10.12.26.1}

\index{Tate object $\newTate(w)$}

We have integrable $\nbigr$-triples
$\newTate(w)\!:=\!\bigl(
 \nbigo_{\cnum_{\lambda}}\lambda^{-w},
 \nbigo_{\cnum_{\lambda}}\lambda^w,
 C_0\bigr)$
for $w\in\seisuu$,
where $C_0$ is the natural pairing:
\[
 C_0(f,\sigma^{\ast}g)
=f\cdot\overline{\sigma^{\ast}g}.
\]
We naturally have
$\newTate(w)^{\ast}=\bigl(
 \nbigo_{\cnum_{\lambda}}\lambda^{w},\,
 \nbigo_{\cnum_{\lambda}}\lambda^{-w},\,
 C_0\bigr)$.
We shall use the twisted identification
$c_w:\newTate(w)^{\ast}
 \simeq\newTate(-w)$
given by 
$\bigl((-1)^w,(-1)^w\bigr)$.

We shall also use the natural identifications
$d_w:\newTate(w)^{\lor}
\simeq
 \newTate(-w)$
and
$j^{\ast}\newTate(w)
 \simeq\newTate(w)$.
Then, as the composition of the above morphisms,
we obtain an induced real structure
$\kappa:\gammatilde_{\sm}^{\ast}
\newTate(w)\simeq
 \newTate(w)$
given by $\bigl((-1)^w,(-1)^w\bigr)$.
\begin{rem}
The real part of 
the fiber $\nbigo_{\cnum_{\lambda}|1}$
is given by $(\sqrt{-1})^w\real$.
\hfill\qed
\end{rem}

Let $\Tate^S(w)$ be the smooth $\nbigr$-triple
given by
$\bigl(\nbigo_{\cnum_{\lambda}},\,
 \nbigo_{\cnum_{\lambda}},\,
 (\sqrt{-1}\lambda)^{-2w}\bigr)$.
\index{Tate object $\Tate^S(w)$}
Note that we have an isomorphism
$\Psi_w:\newTate(w)\simeq\Tate^S(w)$
given by the pair of the following morphisms:
\[
\begin{CD}
 \nbigo_{\cnum_{\lambda}}\lambda^{-w}
 @<{(\sqrt{-1}\lambda)^{-w}}<<
 \nbigo_{\cnum_{\lambda}},
\end{CD}
\quad\quad\quad
\begin{CD}
 \nbigo_{\cnum_{\lambda}}\lambda^{w}
 @>{(-\sqrt{-1}\lambda)^{-w}}>>
 \nbigo_{\cnum_{\lambda}}
\end{CD}
\]
Then, the following diagram is commutative:
\[
 \begin{CD}
 \newTate(w)^{\ast}
 @>{c_w}>>
 \newTate(-w) \\
 @A{\Psi_w^{\ast}}AA @V{\Psi_{-w}}VV \\
 \Tate^S(w)^{\ast}
 @>{c_w^{S}:=(\id,\id)}>>
 \Tate^S(-w)
 \end{CD}
\]
Hence, we may replace 
the family
$\bigl( \Tate^S(w),c_w^S\,\big|\,
 w\in\seisuu\bigr)$
with 
$\bigl(
\newTate(w),c_w\,\big|\,
 w\in\seisuu \bigr)$.

\begin{rem}
Let $X$ be a complex manifold.
For the canonical morphism $a_X$ of $X$
to a point, 
the pull back $a_X^{\ast}\newTate(w)$
is denoted by $\newTate(w)_X$,
or just by $\newTate(w)$,
when there is no risk of confusion.
\hfill\qed
\end{rem}
\begin{rem}
\label{rem;13.8.1.1}
Let $X=X_0\times\cnum_t$.
Let $\nbigt=(\nbigm',\nbigm'',C)$
be any $\nbigr_X$-triple
which is strictly specializable along $t$.
Under the above identification
$\Tate^S(-1)\simeq\newTate(-1)$,
the morphism 
$\nbign:\psi_{t,u}(\nbigt)
 \lrarr\psi_{t,u}(\nbigt)
 \otimes\newTate(-1)$
in {\rm\S\ref{subsection;13.5.4.40}}
is identified with
$-\nbign$ in {\rm\S3.6.1} of {\rm\cite{sabbah2}}.
Note that
$\del_tt=(\sqrt{-1}\lambda)^{-1}(\sqrt{-1}\deldel_tt)
=(-\sqrt{-1}\lambda)^{-1}(-\sqrt{-1}\deldel_tt)$.
\hfill\qed
\end{rem}

\subsubsection{Variant}
\label{subsection;10.12.27.2}

\index{$\nbigr$-triple $\nbigu(p,q)$}

For $(p,q)\in\seisuu^2$,
we set 
let $\nbigu(p,q):=
\bigl(
 \nbigo_{\nbigx}\,\lambda^p,
 \nbigo_{\nbigx}\,\lambda^q,C_0
 \bigr)$.
We have an isomorphism
$\Psi:\Tate^S\bigl((q-p)/2\bigr)
\lrarr
 \nbigu(p,q)$ given by
$\bigl(
 (\sqrt{-1}\lambda)^{-p},
(-\sqrt{-1}\lambda)^{q}
 \bigr)$.
Then, the following diagram is commutative:
\[
 \begin{CD}
 \Tate^{S}\bigl((q-p)/2\bigr)
 @>>>
 \Tate^S\bigl((q-p)/2\bigr)^{\ast}
 \otimes
 \Tate^S\bigl((q-p)\bigr) \\
 @V{\Psi}VV @AA{\Psi^{\ast}\otimes \Phi}A
 \\
 \nbigu(p,q)
 @>{((-1)^p,(-1)^p)}>>
 \nbigu(p,q)^{\ast}\otimes
 \newTate\bigl(-(p-q)\bigr)
 \end{CD}
\]
Hence, we use the polarization
of $\nbigu(p,q)$ given by
$\bigl((-1)^p,(-1)^p\bigr)$.

\subsection{Other basic examples of smooth $\nbigr$-triples
 with integrable and real structure of rank one}
\label{subsection;10.12.18.1}

Let $\gminia$ be a meromorphic function on $X$
whose pole is contained in $D$.
Let $L(\gminia)=\nbigo_{X(\ast D)}\,e$
with $h(e,e)=1$ and $\theta=d\gminia$,
which gives a wild harmonic bundle.
We have $\delbar e=0$,
$\del e=0$ and $\theta^{\dagger}=d\gminiabar$.
The associated family of $\lambda$-flat bundles
$\nbigl(\gminia)$ has a frame
$v=e\,\exp(-\lambda\gminiabar)$.
We have the pairing
$C(v,\sigma^{\ast}v)=
 h(v,\sigma^{\ast}v)
=\exp(-\lambda\gminiabar
+\lambda^{-1}\gminia)$
on $\vecS\times X$,
which takes values in
$\nbigc^{\infty\,\moderate D}_{\vecS\times X}$.
We obtain a smooth $\nbigr_{X(\ast D)}$-triple
$\nbigt_{\gminia}=\bigl(
 \nbigl(\gminia),\nbigl(\gminia),C
 \bigr)$.
\index{$\nbigr$-triple $\nbigt_{\gminia}$}
We have the natural isomorphism
$\nbigt_{\gminia}^{\ast}=\nbigt_{\gminia}$
given by $(\id,\id)$.
We also have natural isomorphisms
$j^{\ast}\nbigt_{\gminia}\simeq
 \nbigt_{-\gminia}$
and
$\nbigt_{\gminia}^{\lor}
\simeq
 \nbigt_{-\gminia}$.
Hence, we have a natural isomorphism
$\gammatilde_{\sm}^{\ast}\nbigt_{\gminia}
\simeq\nbigt_{\gminia}$,
which is a real structure
of $\nbigt_{\gminia}$.

\vspace{.1in}

Let $a\in\real$.
We consider a harmonic bundle
$L(a)=\nbigo_{\cnum_z^{\ast}}\,e$
with $h(e,e)=|z|^{-2a}$ and $\theta=0$.
We have $\delbar e=0$,
$\del e=e\,(-a)\,dz/z$
and $\theta^{\dagger}=0$.
We have 
$(\delbar+\lambda\theta^{\dagger})e=0$,
$\DD e=e(-\lambda\,a)\,dz/z$
and $\lambda^2\del_{\lambda}\,e=0$.
We have the pairing
$C(e,\sigma^{\ast}e)=|z|^{-2a}$.
We obtain an $\nbigr$-triple
$\nbigt_a=\bigl(\nbigl(a),\nbigl(a),C\bigr)$.

We have a natural isomorphism
$\nbigt_a^{\ast}\simeq\nbigt_a$
given by $(\id,\id)$
with which $\nbigt_a$ is a polarized
pure twistor $D(\ast z)$-module of weight $0$.
We have a natural isomorphism
$j^{\ast}\nbigt_a\simeq \nbigt_a$.
We have a natural isomorphism
$\nbigt_a^{\lor}
\simeq
 \nbigt_{-a}$.
Hence, for $0<a<1/2$,
we have a real structure of
$\nbigt_a\oplus\nbigt_{-a}$,
given as follows:
\[
 \gammatilde_{\sm}^{\ast}\nbigt_a
\oplus
 \gammatilde_{\sm}^{\ast}\nbigt_{-a}
\simeq 
 \nbigt_{-a}\oplus
 \nbigt_a
\simeq
 \nbigt_a\oplus\nbigt_{-a},
\]
where the second isomorphism is
the exchange of the components.
If $a=1/2$,
we have a real structure of $\nbigt_{1/2}$
given by
$\gammatilde_{\sm}^{\ast}\nbigt_{1/2}
\simeq
 \nbigt_{-1/2}
\simeq
 \nbigt_{1/2}$.

\vspace{.1in}

Let $\alpha\in\cnum$.
We have a harmonic bundle
 $L(\alpha)=\nbigo_{\cnum_z^{\ast}}\,e$
with $h(e,e)=1$ and $\theta\,e=e\,\alpha\,dz/z$.
We put $v=e\,\exp(-\lambda\alphabar\,\log|z|^2)$,
which gives a frame of $\nbigl(\alpha)$
on $\cnum_{\lambda}\times\cnum_z^{\ast}$.
We have $z\del_zv=
 v\,(\lambda^{-1}\alpha+\lambda\alphabar)$.
We have the pairing
$C(v,\sigma^{\ast}v)=
 |z|^{2(-\lambda\alphabar+\lambda^{-1}\alpha)}$.
We obtain an $\nbigr$-triple
$\nbigt_{\alpha}=\bigl(
 \nbigl(\alpha),\nbigl(\alpha),C_{\alpha}
 \bigr)$.
We have natural isomorphisms
$\nbigt_{\alpha}^{\ast}
\simeq\nbigt_{\alpha}$
and
$\nbigt_{\alpha}^{\lor}
\simeq
 \nbigt_{-\alpha}$
and 
$j^{\ast}\nbigt_{\alpha}
\simeq
 \nbigt_{-\alpha}$.
Hence, we have a real structure
$\gammatilde_{\sm}^{\ast}\nbigt_{\alpha}
\simeq
 \nbigt_{\alpha}$.
But, it is easy to observe that
$\nbigt_{\alpha}$ cannot be integrable,
i.e., $\nbigl(\alpha)$ cannot underlie
an integrable $\nbigr$-module.
Indeed, assume that it is integrable.
We have a meromorphic function $A$
whose poles is contained in
$\cnum^{\ast}_{\lambda}\times\{0\}$
such that
$\del_{\lambda} v=v\,A$
on $\cnum_{\lambda}^{\ast}\times X$.
We obtain the relation
$z\del_z A=-\alphabar-\lambda^2\alpha$,
which contradicts with the condition on $A$.

\subsubsection{Ramified exponential twist}

Let $\varphi_n:\cnum_{s}\lrarr \cnum_t$ be given by
$\varphi_n(s)=s^n$.
We set $X^{(n)}:=X_0\times\cnum_s$.
The induced map
$X^{(n)}\lrarr X$ is also denoted by $\varphi_n$.
We give a complement to \S22.4.2 and 
\S22.11.2 of \cite{mochi7}.
(See also \cite{sabbah2} and \cite{sabbah5}.)
Let $\nbigm$ be a strict coherent $\nbigr_{X^{(n)}}(\ast s)$-module.
We have $\varphi_{n\dagger}\nbigm=\varphi_{n\ast}\nbigm$.
\begin{lem}
$\nbigm$ is strictly specializable along $s$,
if and only if $\varphi_{n\dagger}\nbigm$ is 
strictly specializable along $t$.
In that case,
we have a natural isomorphism
Hence, we have a natural isomorphism
$\psitilde_{-\vecdelta+u}(\nbigm)
\simeq
 \psitilde_{-\vecdelta+u/n}(\varphi_{n\dagger}\nbigm)$.
\end{lem}
\pf
Suppose that 
$\nbigm$ is strictly specializable along $s$.
Then, $\varphi_{n\dagger}\nbigm=\varphi_{n\ast}\nbigm$
is strictly specializable along $t$.
Indeed, the $V$-filtration at $\lambda_0$ is given by
\[
 \Vzero_{-1+a}(\varphi_{n\dagger}\nbigm)
=\varphi_{n\ast}
 \bigl(\Vzero_{-1+a/n}\nbigm\bigr).
\]
Note the equality
$-\deldel_{t_n}t_n+\eigenmap(\lambda,-\vecdelta+u)
=n(-\deldel_tt+\eigenmap(\lambda,-\vecdelta+u/n))$.
In particular,
we have a natural isomorphism
$\psitilde_{-\vecdelta+u}(\nbigm)
\simeq
 \psitilde_{-\vecdelta+u/n}(\varphi_{n\dagger}\nbigm)$.
Conversely,
suppose that 
$\varphi_{n\dagger}(\nbigm)$
is strictly specializable along $t$.
As remarked in Lemma 22.4.7,
$\varphi_n^{\dagger}\varphi_{n\dagger}(\nbigm)$
is strictly specializable along $s$.
Because 
$\nbigm$ is a direct summand of
$\varphi_n^{\dagger}\varphi_{n\dagger}(\nbigm)$,
we obtain that
$\nbigm$ is strictly specializable along $s$.
\hfill\qed

\begin{cor}
A strict coherent $\nbigr_{X^{(n)}}(\ast s)$-triple
$\nbigt_1$ is strictly specializable along $s$,
if and only if
$\varphi_{n\dagger}\nbigt$ is 
strictly specializable along $t$.
In that case,
we have a natural isomorphism:
\[
 \psitilde_{-\vecdelta+u}(\nbigt_1)
\simeq
 \psitilde_{-\vecdelta+u/n}\bigl(\varphi_{\dagger}\nbigt_1\bigr).
\]
\end{cor}
\pf
The comparison of the induced Hermitian sesqui-linear pairings
can be done by an easy computation.
\hfill\qed

\vspace{.1in}

Let $\gminia\in s^{-1}\cnum[s^{-1}]$.
We have the $\nbigr_{\cnum_s}(\ast s)$-triple
$\nbigt_{-\gminia}$.
Via the pull back by the projection
$X^{(n)}\lrarr\cnum_s$,
it induces a smooth $\nbigr_{X^{(n)}}(\ast s)$-triple,
which is also denoted by $\nbigt_{-\gminia}$.
Let $\nbigt$ be a strict coherent $\nbigr_X(\ast t)$-triple.
Recall that,
if $\varphi_n^{\dagger}(\nbigt)\otimes\nbigt_{-\gminia}$
is strictly specializable along $s$,
we set
\[
 \psitilde_{t,\gminia,u}(\nbigt)
:=\psitilde_{s,u}\bigl(
 \varphi_n^{\dagger}(\nbigt)
 \otimes\nbigt_{-\gminia}\bigr).
\]
We have a smooth $\nbigr_{X}(\ast t)$-triple
$\varphi_{n\dagger}(\nbigt_{-\gminia})$.
It is easy to observe that
$\varphi_n^{\dagger}\nbigt
 \otimes\nbigt_{-\gminia}$
is strictly specializable along $s$,
if and only if
$\nbigt\otimes
 \varphi_{n\dagger}\nbigt_{-\gminia}$
is strictly specializable along $t$.
By the above consideration,
we have natural isomorphism
\begin{equation}
\label{eq;13.5.9.200}
 \psitilde_{t,\gminia,-\vecdelta+u}(\nbigt)
\simeq
 \psitilde_{t,-\vecdelta+u/n}
 \bigl(
 \nbigt\otimes\varphi_{n\dagger}\nbigt_{-\gminia}
 \bigr).
\end{equation}

\section{Deformation associated to nilpotent morphisms}
\label{subsection;13.4.12.2}

An $\nbigr$-triple with a tuple of nilpotent morphisms
has a natural deformation.
Such a deformation was used in the contexts of
variation of twistor structures in \cite{mochi2}
and \cite{mochi8},
which originated from the procedure
to construct nilpotent orbit in the Hodge theory.

\index{deformation associated to nilpotent maps}

\subsection{Twistor nilpotent orbit
in $\nbigr$-triple}
\label{subsection;10.12.25.20}

\index{Twistor nilpotent orbit}

Let $X$ be a complex manifold,
and let $D$ be a hypersurface.
Let $\Lambda$ be a finite set.
Let $\Lambda\textrm{-}\rtriplecat(X,D)$
be the category 
of $\nbigr_{X(\ast D)}$-triple $\nbigt$
with a tuple $\vecnbign=(\nbign_i\,|\,i\in\Lambda)$ of 
mutually commuting morphisms
$\nbign_i:\nbigt\lrarr\nbigt\otimes\newTate(-1)$.
\index{category $\Lambda\textrm{-}\rtriplecat(X,D)$}
A morphism
$F:(\nbigt_1,\vecnbign)\lrarr(\nbigt_2,\vecnbign)$
in $\Lambda\textrm{-}\rtriplecat(X,D)$
is a morphism $F$ in $\rtriplecat(X,D)$
such that $F\circ \nbign_i=\nbign_i\circ F$
for $i\in\Lambda$.

For $I\subset\Lambda$,
we put
$X\{I\}:=X\times\cnum^I$
and $D\{I\}:=
 \bigl(D\times\cnum^I\bigr)
 \cup\bigcup_{i\in I}\{z_i=0\}$.
We shall construct a functor
\begin{equation}
 \TNIL_I:
 \Lambda\textrm{-}\rtriplecat(X,D)
\lrarr
 \Lambda\textrm{-}\rtriplecat\bigl(X\{I\},D\{I\}\bigr)
\label{eq;10.8.5.4}
\end{equation}
\index{functor $\TNIL_I$}
Let $(\nbigt,\vecnbign)\in\Lambda\textrm{-}\rtriplecat(X,D)$.
Let $\nbigt=(\nbigm',\nbigm'',C)$.
Recall that morphisms
$\nbign_i=(\nbign_i',\nbign_i'')$
are given as
\[
\begin{CD}
 \nbigm'@<{\nbign_i'}<<\nbigm'\cdot\lambda,
\end{CD}
\quad\quad
\begin{CD}
  \nbigm''@>{\nbign_i''}>> \nbigm''\cdot\lambda^{-1}.
\end{CD}
\]
The underlying 
$\nbigr_{X\{I\}(\ast D\{I\})}$-modules of
$\TNIL_I(\nbigt,\vecN)$ are given as follows,
with the natural actions of $\nbigr_X$:
\[
 \nbigm'\otimes
 \nbigo_{\nbigx\{I\}}(\ast \nbigd\{I\}),
\quad
 z_i\deldel_im'
=-\lambda\,\nbign_i'm',
\]
\[
 \nbigm''\otimes
 \nbigo_{\nbigx\{I\}}(\ast \nbigd\{I\}),
\quad
 z_i\deldel_im''
=-\lambda\,\nbign_i''m''.
\]
Here, 
$m'$ and $m''$ denote the sections
$m'\otimes 1$ and $m''\otimes 1$.
We would like to define the pairing $\Ctilde$
of $\TNIL_I(\nbigt,\vecN)$ 
by the following formula:
\begin{equation}
\label{eq;10.9.29.50}
 \Ctilde(m',\sigma^{\ast}m'')
=C\Bigl(
 \exp\Bigl(\sum_{i\in I}
 -\nbign_i'\log|z_i|^2\Bigr)m',\,\sigma^{\ast}m''
 \Bigr)
\end{equation}
Because
$C\bigl(
 \nbign_i'm',\sigma^{\ast}m''
 \bigr)
=C\bigl(m',\sigma^{\ast}(\nbign_i''m'')\bigr)$,
the right hand side of {\rm(\ref{eq;10.9.29.50})}
is equal to the following:
\[
 C\Bigl(
 m',\,\sigma^{\ast}\Bigl(
 \exp\Bigl(\sum_{i\in I}
-\nbign_i''\log|z_i|^2
 \Bigr)m''
 \Bigr)
 \Bigr)
\]

Let us prove that $\Ctilde$ is 
an $\nbigr_{X\{I\}(\ast D\{I\})}
 \times
 \sigma^{\ast}
 \nbigr_{X\{I\}(\ast D\{I\})}$-homomorphism.
For that purpose,
we have only to prove
the following equalities:
\begin{equation}
 \label{eq;10.8.5.3}
 z_i\del_i\Ctilde(m',\sigma^{\ast}m'')
=\Ctilde(z_i\del_im',\sigma^{\ast}m''),
\quad\quad
 \zbar_i\del_{\zbar_i}
 \Ctilde(m',\sigma^{\ast}m'')
=\Ctilde\bigl(m',\sigma^{\ast}(z_i\del_im'')\bigr)
\end{equation}
We have the following:
\begin{multline}
z_i\del_i\Ctilde(m',\sigma^{\ast}m'') 
=z_i\del_iC\Bigl(
 \exp\Bigl(-\sum_{j\in I}\nbign_j'\log|z_j|^2\Bigr)
 m',\sigma^{\ast}m''
 \Bigr)\\
=C\Bigl(
 -\exp\Bigl(-\sum_{j\in I}\nbign_i'\log|z_j|^2\Bigr)
 \,\nbign_i'\,m',\,\,\sigma^{\ast}m''
 \Bigr)
=\Ctilde(z_i\del_im',\sigma^{\ast}m'')
\end{multline}
We also have the following:
\begin{multline}
 \zbar_i\del_{\zbar_i}
 \Ctilde(m',\sigma^{\ast}m'')
=\zbar_i\del_{\zbar_i}
 C\Bigl(m',
 \sigma^{\ast}\Bigl(
 \exp\Bigl(-\sum_{j\in I}\nbign_j''\log|z_i|^2
 \Bigr)m''
 \Bigr)
 \Bigr) \\
=C\Bigl(
 m',\sigma^{\ast}\Bigl(
 -\exp\Bigl(
 -\sum_{j\in I}
 \nbign_j''\log|z_j|^2
 \Bigr)\,\nbign_i''\,m''
 \Bigr)
 \Bigr)
=\Ctilde\bigl(m',\sigma^{\ast}(z_i\del_im'')\bigr)
\end{multline}
Hence, we obtain (\ref{eq;10.8.5.3}),
and $\Ctilde$ gives 
a hermitian sesqui-linear pairing
of $\nbigm'$ and $\nbigm''$.

By using the commutativity of
$\nbign_j$ and $\nbign_i$ $(i\in I)$,
we obtain
$\Ctilde\bigl(
 \nbign_j'm',\sigma^{\ast}m''
 \bigr)
=\Ctilde\bigl(
 m',\sigma^{\ast}\nbign_j''m''
 \bigr)$.
Hence, we have morphisms
\[
\TNIL_I(\nbign_j):\TNIL_I(\nbigt,\vecnbign)
\lrarr\TNIL_I(\nbigt,\vecnbign)\otimes
 \newTate(-1)
\]
for $j\in\Lambda$.
In particular,
an $\nbigr$-triple
$\TNIL_I(\nbigt,\vecnbign)$
with $\TNIL_I(\vecnbign)
=\bigl(\TNIL_I(\nbign_i)\,|\,i\in\Lambda\bigr)$
is an object of
$\Lambda\textrm{-}\rtriplecat\bigl(X\{I\},D\{I\}\bigr)$.
Thus, we obtain the functor
(\ref{eq;10.8.5.4}).

\vspace{.1in}
It is easy to see that $\TNIL_I(\nbigt,\vecnbign)$
is strictly specializable along $z_i$ $(i\in I)$,
and we have
\[
 \psitilde_{z_i,-\vecdelta}\TNIL_I(\nbigt,\vecnbign)
=\TNIL_{I\setminus i}(\nbigt,\vecnbign).
\]
The following lemma is easy to see.
\begin{lem}
If $(\nbigt,\vecnbign)$ is integrable,
$\TNIL_I(\nbigt,\vecnbign)$ is naturally integrable.
\hfill\qed
\end{lem}

\subsubsection{Hermitian adjoint}
For $(\nbigt,\vecnbign)\in
\Lambda\textrm{-}\rtriplecat(X,D)$,
we set
$(\nbigt,\vecnbign)^{\ast}:=
 (\nbigt^{\ast},-\vecnbign^{\ast})$,
which gives a contravariant functor
on $\Lambda\textrm{-}\rtriplecat(X,D)$.

\begin{lem}
$\bigl(\TNIL_I(\nbigt,\vecnbign)\bigr)^{\ast}$
is naturally identified with
$\TNIL_I\bigl((\nbigt,\vecnbign)^{\ast}\bigr)$.
\end{lem}
\pf
Let $\nbigt=(\nbigm',\nbigm'',C)$
and $\nbign_i=(\nbign_i',\nbign_i'')$.
Because the identification
$c_{-1}:\newTate(-1)^{\ast}
\simeq
 \newTate(1)$
is given by $(-1,-1)$,
we have $-\nbign_i^{\ast}=(\nbign_i'',\nbign_i')$.
Hence, the underlying $\nbigr$-modules
are naturally identified.

Let us compare the pairings.
Let $\Ctilde^{\ast}$ be the pairing
for $\bigl(\TNIL_I(\nbigt,\vecnbign)\bigr)^{\ast}$.
We have
\begin{multline}
 \Ctilde^{\ast}(m'',\sigma^{\ast}m')
=\overline{\sigma^{\ast}
 \Ctilde(m',\sigma^{\ast}m'')}
=\overline{
 \sigma^{\ast}C\Bigl(
 \exp\Bigl(-\sum_{i\in I}
 \nbign_i'\log|z_i|^2\Bigr)m',\,\sigma^{\ast}m''
 \Bigr)} \\
=C^{\ast}\Bigl(
 m'',\sigma^{\ast}\Bigl(
 \exp\Bigl(-\sum_{i\in I}\nbign_i'
 \log|z_i|^2\Bigr)m'
 \Bigr)
 \Bigr)
\end{multline}
This is exactly the pairing for
$\TNIL_I\bigl((\nbigt,\vecnbign)^{\ast}\bigr)$.
\hfill\qed

\subsubsection{Dual}
If $\nbigt\in\VTS(X,D)$,
we set 
$(\nbigt,\vecnbign)^{\lor}:=
 (\nbigt^{\lor},-\vecnbign^{\lor})$.

\begin{lem}
$\bigl(\TNIL_I(\nbigt,\vecnbign)\bigr)^{\lor}$
is naturally identified with
$\TNIL_I\bigl((\nbigt,\vecnbign)^{\lor}\bigr)$.
\end{lem}
\pf
The underlying $\nbigr_X$-modules are
naturally identified.
Let us compare the pairings.
We have only to consider the case that $\Lambda$
consists of one element.
For $\nbigt=(\nbigm',\nbigm'',C)$,
let $\Psi_C$ and $\Phi_C$ be defined
as in (\ref{eq;10.12.25.1}).
Similarly,
we obtain $\Psi_{\Ctilde}$
and $\Phi_{\Ctilde}$ from the pairing $\Ctilde$.
Note that we have
$\Psi_C\circ\nbign'_i=
 \sigma^{\ast}(\nbign''_i)^{\lor}\circ\Psi_C$.
Indeed,
\begin{multline}
 \bigl\langle
 \Psi_{C}(\nbign'_i\,m'),\sigma^{\ast}m''
 \bigr\rangle
=C(\nbign'_i\,m',\sigma^{\ast}m'')
=C\bigl(m',\sigma^{\ast}(\nbign''_i\,m'')\bigr)
 \\
=\bigl\langle
 \Psi_C(m'),\,\sigma^{\ast}(\nbign''_i\,m'')
 \bigr\rangle 
=\bigl\langle 
 \sigma^{\ast}(\nbign''_i)^{\lor}
 \Psi_C(m'),\,\sigma^{\ast}m''
 \bigr\rangle
\end{multline}
Similarly,
we have
$\Phi_C\circ\sigma^{\ast}\nbign''_i
=(\nbign'_i)^{\lor}\circ\Phi_C$.
We have
\[
 \Psi_{\Ctilde}
=\Psi_C\circ\exp\Bigl(-\sum_i\log|z_i|^2\nbign'_i\Bigr)
=\exp\Bigl(-\sum_i\log|z_i|^2\sigma^{\ast}(\nbign''_i)^{\lor}\Bigr)
\circ\Psi_C
\]
Indeed, we have the equalities:
\begin{multline}
 \langle
 \Psi_{\Ctilde}(m'),\sigma^{\ast}m''
 \rangle
=\Ctilde(m',\sigma^{\ast}m'')
=C\Bigl(
 \exp\Bigl(-\sum_i\log|z_i|^2\nbign'_i\Bigr)\,m',\,
 \sigma^{\ast}m''
 \Bigr) \\
=\Bigl\langle
 \Psi_C\Bigl(
 \exp\Bigl(-\sum_i\log|z_i|^2\nbign'_i\Bigr)\,m'
 \Bigr),\,\sigma^{\ast}m''
 \Bigr\rangle
\end{multline}
Similarly,
$\Phi_{\Ctilde}
=\Phi_C\circ
 \exp\Bigl(-\sum_i\log|z_i|^2\sigma^{\ast}\nbign''_i\Bigr)
=\exp\Bigl(
 -\sum_i\log|z_i|^2(\nbign'_i)^{\lor}
 \Bigr)\circ\Phi_C$.
Then, we obtain the following:
\begin{multline}
 (\Ctilde)^{\lor}(n',\sigma^{\ast}n'')
=C\Bigl(
 \exp\Bigl(-\sum_i\log|z_i|^2\nbign'_i\Bigr)
 \Psi_{\Ctilde}^{-1}(\sigma^{\ast}n''),
 \Phi_{\Ctilde}^{-1}(n')
\Bigr)
 \\
=C\Bigl(
 \Psi_C^{-1}(\sigma^{\ast}n''),\,
 \Phi_C^{-1}\Bigl(
 \exp\Bigl(\sum_i\log|z_i|^2(\nbign'_i)^{\lor}\Bigr)\,n'
 \Bigr)
 \Bigr) \\
=C^{\lor}\Bigl(
 \exp\Bigl(\sum_i\log|z_i|^2(\nbign'_i)^{\lor}\Bigr)n',\,
 \sigma^{\ast}n''
 \Bigr)
=\widetilde{\bigl(C^{\lor}\bigr)}
 (n',\sigma^{\ast}n'')
\end{multline}
Thus, we are done.
\hfill\qed

\subsubsection{Real structure}
We also have a natural isomorphism
$j^{\ast}\TNIL_I(\nbigt,\vecnbign)
\simeq
 \TNIL_I j^{\ast}(\nbigt,\vecnbign)$.
Hence, we have a natural isomorphism
\[
 \gammatilde^{\ast}_{\sm}\TNIL_I(\nbigt,\vecnbign)
\simeq
 \TNIL_I\gammatilde^{\ast}_{\sm}(\nbigt,\vecnbign).
\]

A real structure of an object
$(\nbigt,\vecnbign)$ in
$\Lambda\textrm{-}\rtriplecat(X,D)$
is an isomorphism
$\kappa:\gammatilde^{\ast}_{\sm}(\nbigt,\vecnbign)
\simeq
 (\nbigt,\vecnbign)$
such that $\kappa\circ\gammatilde^{\ast}_{\sm}\kappa=\id$.
A real structure of $(\nbigt,\vecnbign)$
naturally induces a real structure of
$\TNIL_I(\nbigt,\vecnbign)$.

\subsection{Variant}

Let $(\nbigt,\vecnbign)\in
 \Lambda\textrm{-}\rtriplecat(X,D)$.
Let $\vecvarphi=(\varphi_i\,|\,i\in\Lambda)$ be 
a tuple of meromorphic functions
whose poles are contained in $D$.
We shall construct an $\nbigr_{X(\ast D)}$-triple
$\Def_{\vecvarphi}(\nbigt,\vecnbign)$.
\index{$\nbigr$-triple $\Def_{\vecvarphi}(\nbigt,\vecnbign)$}
We define a new $\nbigr$-action on $\nbigm'$ by
$f\bullet m'=f\,m'$ for $f\in\nbigo_{\nbigx}$
and
$v\bullet m'=v\,m'-\sum_i\nbign_i'(m')\,v(\varphi_i)$
for $v\in\lambda\,\Theta_X$.
We define a new $\nbigr$-action on $\nbigm''$
in the same way.
The $\nbigr_{X(\ast D)}$-modules are denoted by
$\nbigm'_{(\vecvarphi,\vecnbign)}$
and $\nbigm''_{(\vecvarphi,\vecnbign)}$.
We define a hermitian sesqui-linear pairing
$C_{(\vecvarphi,\vecnbign)}$ of
$\nbigm'_{(\vecvarphi,\vecnbign)}$ and 
$\nbigm''_{(\vecvarphi,\vecnbign)}$
as follows:
\[
 C_{(\vecvarphi,\vecnbign)}:=
 C\Bigl(
 \exp\bigl(-\vecvarphi\cdot\vecnbign'\bigr)\,m',\,
 \sigma^{\ast}\bigl(
 \exp\bigl(-\vecvarphi\cdot\vecnbign''\bigr)\,m''
 \bigr)
 \Bigr)
\]
Here,
$\vecvarphi\cdot\vecnbign':=
 \sum_{i\in\Lambda}\varphi_i\,\nbign'_i$
and 
$\vecvarphi\cdot\vecnbign'':=
 \sum_{i\in\Lambda}\varphi_i\,\nbign''_i$.
The triple is denoted by
$\Def_{\vecvarphi}(\nbigt,\vecnbign)$.
The following lemma can be checked
as in \S\ref{subsection;10.12.25.20}.
\begin{lem}
We have natural isomorphisms
\[
 \Def_{\vecvarphi}\bigl(
(\nbigt,\vecnbign)\bigr)^{\ast}
\simeq
  \Def_{\vecvarphi}\bigl(
(\nbigt,\vecnbign)^{\ast}\bigr),
\]
\[
 \Def_{\vecvarphi}\bigl(
(\nbigt,\vecnbign)\bigr)^{\lor}
\simeq
  \Def_{\vecvarphi}\bigl(
(\nbigt,\vecnbign)^{\lor}\bigr),
\]
\[
  j^{\ast}\Def_{\vecvarphi}
(\nbigt,\vecnbign)
\simeq
 \Def_{\vecvarphi}\bigl(
 j^{\ast}(\nbigt,\vecnbign)\bigr).
\]
In particular,
we have
$\gammatilde^{\ast}_{\sm}
 \Def_{\vecvarphi}\bigl(\nbigt,\vecnbign\bigr)
\simeq
 \Def_{\vecvarphi}\bigl(
 \gammatilde^{\ast}_{\sm}(\nbigt,\vecnbign)\bigr)$.
\hfill\qed
\end{lem}

A real structure of $(\nbigt,\vecnbign)$
naturally induces a real structure of
$\Def_{\vecvarphi}(\nbigt,\vecnbign)$.
If $(\nbigt,\vecnbign)$ is integrable,
$\Def_{\vecvarphi}(\nbigt,\vecnbign)$
is also integrable.
The following lemma 
can be checked directly.
\begin{lem}
$\Def_{\vecvarphi'}\bigl(
 \Def_{\vecvarphi}
(\nbigt,\vecnbign)\bigr)
\simeq
\Def_{\vecvarphi'+\vecvarphi}
 (\nbigt,\vecnbign)$ naturally.
\hfill\qed
\end{lem}

Let $Y$ be a complex manifold
with a hypersurface $D_Y$.
Let 
\[
 \iota_j:(Y,D_Y)\subset (X\{I\},D\{I\})
\quad\quad (j=1,2)
\]
be open embeddings
such that
(i) the compositions
 $Y\stackrel{\iota_j}\lrarr X\{I\}\lrarr X$
 are equal, where the latter is the natural projection,
(ii) $\iota_1^{\ast}(z_i)=e^{\varphi_i}\iota_2^{\ast}(z_i)$
for $i\in I$,
where $\varphi_i$ are holomorphic on $Y$.
The following lemma can be checked directly.
\begin{lem}
$\iota_1^{\ast}\TNIL_I(\nbigt,\vecnbign)
\simeq\Def_{\vecvarphi}
 \Bigl(
 \iota_2^{\ast}\TNIL_I(\nbigt,\vecnbign)
 \Bigr)$ naturally.
\hfill\qed
\end{lem} 

\section{Beilinson triples}
\label{subsection;11.2.1.10}

We introduce a special type of smooth $\nbigr$-triples,
which we call Beilinson triples.
It will be used in the construction of
Beilinson functor in \S\ref{section;11.4.3.1}.

\index{Beilinson triple}

\subsection{Triples on a point}
\label{subsection;13.3.27.10}

We put
$A:=\nbigo_{\cnum_{\lambda}}
 [\lambda s,(\lambda s)^{-1}]$.
We set
$A^a:=(\lambda s)^a\,\nbigo[\lambda s]
 \subset A$.
For $a\leq b$, we put
\[
 \II_1^{a,b}:=A^{-b+1}/A^{-a+1}
\simeq
 \bigoplus_{a\leq i<b}
 \nbigo_{\cnum_{\lambda}}\cdot(\lambda s)^{-i},
\quad
 \II_2^{a,b}:=A^a/A^b
\simeq
 \bigoplus_{a\leq i<b}
 \nbigo_{\cnum_{\lambda}}\cdot(\lambda s)^i.
\]
\index{sheaves $\II_j^{a,b}$}
Let $C_0^{(i)}$ be the pairing of 
$\nbigo_{\cnum_{\lambda}}\cdot(\lambda s)^{-i}$
and 
$\nbigo_{\cnum_{\lambda}}\cdot(\lambda s)^i$
given by
$C_0^{(i)}(f\,s^{-i},\sigma^{\ast}g\,s^i)
=f\cdot\overline{\sigma^{\ast}g}$.
They induce a hermitian sesqui-linear
pairing $C_{\II}$ of
$\II_1^{a,b}$ and $\II_2^{a,b}$.
The integrable $\nbigr$-triple is denoted by
$\II^{a,b}$.
\index{$\nbigr$-triple $\II^{a,b}$}
We have a natural identification
$\II^{a,b}
=\oplus_{a\leq i<b}
 \II^{i,i+1}$.
We shall also use the identification
$\II^{i,i+1}\simeq\newTate(i)$
given by
$s^j\longleftrightarrow 1$ $(j=i,-i)$.
We have natural isomorphisms
$\Upsilon_j:
 \II^{a,b}\otimes\newTate(j)
\simeq
 \II^{a+j,b+j}$
given by the multiplication of $s^j$.

The multiplication of $-s$ induces
$\nbign_{\II}=(\nbign_{\II,1},\nbign_{\II,2}):
 \II^{a,b}\lrarr\II^{a,b}\otimes\newTate(-1)$:
\[
\begin{CD}
\II_1^{a,b}
 @<{\nbign_{\II,1}}<<
\II_1^{a,b}\otimes\nbigo_{\cnum_{\lambda}}\lambda,
\end{CD}
\quad\quad
\begin{CD}
 \II_2^{a,b}
 @>{\nbign_{\II,2}}>>
 \II_2^{a,b}\otimes\nbigo_{\cnum_{\lambda}}
 \lambda^{-1}.
\end{CD}
\]

\subsubsection{Hermitian adjoint}
We shall use the identifications
\[
\nbigs^{a,b}:(\II^{a,b})^{\ast}\simeq
 \II^{-b+1,-a+1}
\]
given by
$\bigoplus_{a\leq i<b}(-1)^ic_i$,
where
$c_i:(\II^{i,i+1})^{\ast}
\simeq
 \II^{-i,-i+1}$
are induced by
$c_i:\newTate(i)^{\ast}
\simeq
 \newTate(-i)$
given in \S\ref{subsection;10.12.26.1}.
They give isomorphisms
$\nbigs^{a,b}:
 (\II^{a,b},\nbign_{\II})^{\ast}
\simeq
 (\II^{-b+1,-a+1},\nbign_{\II})$.
Note that
$-\nbign_{\II}^{\ast}
=(\nbign_{\II,2},\nbign_{\II,1})$.
Under the identifications,
we have
$\Upsilon_j^{\ast}=(-1)^j\Upsilon_{-j}$.

\subsubsection{Dual}
We shall use the identifications
$(\II^{a,b})^{\lor}\simeq
 \II^{-b+1,-a+1}$
given by 
$\bigoplus_{a\leq i<b}(-1)^id_i$,
where $d_i:(\II^{i,i+1})^{\lor}
\simeq \II^{-i,-i+1}$
are natural isomorphisms.
Then, we have
\[
 (\II^{a,b},\nbign_{\II})^{\lor}
\simeq
 (\II^{-b+1,-a+1},\nbign_{\II}).
\]
Under the identifications,
we have
$\Upsilon_j^{\lor}=(-1)^j\Upsilon_{-j}$.

\subsubsection{Real structure}
We naturally have
$j^{\ast}(\II^{a,b},\nbign_{\II})
\simeq
 (\II^{a,b},\nbign_{\II})$.
By composition of the above isomorphisms
we obtain a real structure
$\kappa_{\II}:\gammatilde_{\sm}^{\ast}
 (\II^{a,b},\nbign_{\II})
\simeq
 (\II^{a,b},\nbign_{\II})$.
The real structure is the same as 
the direct sum
$\bigoplus_{a\leq i<b}\newTate(i)$,
under the natural isomorphism
$\II^{a,b}\simeq
 \bigoplus_{a\leq i<b}\newTate(i)$.

\subsection{The associated twistor nilpotent orbit}
\label{subsection;11.4.12.10}

We obtain integrable
$\nbigr_{\cnum_t}(\ast t)$-triples
$\IItilde^{a,b}:=\TNIL(\II^{a,b},\nbign_{\II})
=\bigl(
 \IItilde_1^{a,b},\IItilde_2^{a,b},\Ctilde_{\II}
 \bigr)$.
\index{$\nbigr$-triple $\IItilde^{a,b}$}
The underlying $\nbigr$-modules are as follows:
\[
 \IItilde_1^{a,b}=
 \bigoplus_{a\leq i<b}
 \nbigo_{\cnum^2_{\lambda,t}}(\ast t)
 \cdot (\lambda\,s)^{-i},
\quad
 \IItilde_2^{a,b}=
 \bigoplus_{a\leq i<b}
 \nbigo_{\cnum^2_{\lambda,t}}(\ast t)
 \cdot(\lambda\,s)^i,
\quad
 t\deldel_t(\lambda\,s)^i
=(\lambda\,s)^{i+1}
\]
\index{sheaves $\IItilde_j^{a,b}$}
The pairing $\Ctilde_{\II}$ is given as follows:
\begin{equation}
 \label{eq;10.12.25.10}
 \Ctilde_{\II}\bigl(
 (\lambda\,s)^i,\sigma^{\ast}(\lambda\,s)^j
 \bigr)
=
\frac{(\log|t|^2)^{-i-j}}{(-i-j)!}(-1)^j\lambda^{i-j}
 \cdot\chi_{i+j\leq 0}
\end{equation}
Here, 
$\chi_{i+j\leq 0}$ is $1$ if $i+j\leq 0$,
and $0$ if $i+j>0$.
Indeed, we have the following equalities:
\begin{multline}
 \Ctilde_{\II}\Bigl(
 (\lambda\,s)^i,\,\sigma^{\ast}\bigl(
 (\lambda\,s)^j
 \bigr)
 \Bigr)
=C_{\II}\Bigl(
 \exp\bigl(-\nbign_{\II,1}\log|t|^2\bigr)\,(\lambda\,s)^i,\,
\sigma^{\ast}\bigl(
 (\lambda\,s)^j
 \bigr)
 \Bigr) \\
=
\sum_{k=0}^{\infty}
 \frac{(\log|t|^2)^k}{k!}
 C_{\II}\bigl(
 (\lambda\,s)^is^k,\,\,
 \sigma^{\ast}\bigl(
 (\lambda\,s)^j
 \bigr)
 \bigr)
=\sum_{k=0}^{\infty}
 \frac{(\log|t|^2)^k}{k!}
 \lambda^i(-\lambda)^{-j}
 \delta_{i+k+j,0}
\end{multline}
Here, 
$\delta_{m,0}$ is $0$ if $m\neq 0$
and $1$ if $m=0$.
Then, we obtain (\ref{eq;10.12.25.10}).

We have the induced isomorphisms
\[
(\IItilde^{a,b})^{\ast}
\simeq
 \IItilde^{-b+1,-a+1},
\quad
j^{\ast}\IItilde^{a,b}\simeq \IItilde^{a,b},
\quad
 (\IItilde^{a,b})^{\lor}
\simeq
 \IItilde^{-b+1,-a+1}.
\]
It is equipped with
the induced real structure
$\kappa:\gammatilde_{\sm}^{\ast}
 \IItilde^{a,b}\simeq
 \IItilde^{a,b}$.

\subsubsection{Pull back}

Let $X$ be a complex manifold.
Let $f$ be a holomorphic function
on $X$.
We obtain a smooth $\nbigr(\ast f)$-triple
$\IItilde_{f}^{a,b}:=f^{\ast}\IItilde^{a,b}$.
The following lemma is clear
by construction.
\begin{lem}
\label{lem;10.12.25.21}
Let $g=e^{\varphi}\cdot f$,
where $\varphi$ is holomorphic.
We have a natural isomorphism
$\IItilde^{a,b}_g
\simeq
 \Def_{\varphi}(\IItilde_f^{a,b},\nbign)$.
\hfill\qed
\end{lem}

\subsection{Appendix}

Let $\nbigm$ be an $\nbigr_X$-module
with a nilpotent morphism
$\nbign:\nbigm\lrarr\nbigm$
such that $\nbign^L=0$
for some $L\in\seisuu_{>0}$.
We put $\nbign''_{\II}:=\lambda\,\nbign_{\II,2}$.
For integers $M_1\leq M_2$,
we consider the morphism
$\nbign\otimes\id+\id\otimes\nbign''_{\II}:
 \nbigm\otimes\II_2^{M_1,M_2}
\lrarr
 \nbigm\otimes\II_2^{M_1,M_2}$.
The inclusion
$\II_2^{M_1,M_1+1}\lrarr 
 \II_2^{M_1,M_2}$ induces the following morphism:
\begin{equation}
\label{eq;11.1.16.50}
 \nbigm\otimes\II_2^{M_1,M_1+1}
\lrarr
 \Cok\bigl(
  \nbigm\otimes\II_2^{M_1,M_2}
\lrarr
 \nbigm\otimes\II_2^{M_1,M_2}
 \bigr)
\end{equation}
The projection
$\II_2^{M_1,M_2}\lrarr 
 \II_2^{M_2-1,M_2}$ induces the following morphism:
\begin{equation}
\label{eq;11.1.16.51}
 \Ker\bigl(
  \nbigm\otimes\II_2^{M_1,M_2}
\lrarr
 \nbigm\otimes\II_2^{M_1,M_2}
 \bigr)
\lrarr
\nbigm\otimes\II_2^{M_2-1,M_2}
\end{equation}
The following lemma can be 
checked by direct computations.
\begin{lem}
\label{lem;11.1.16.53}
If $M_2-M_1>2L$,
then the morphisms
{\rm(\ref{eq;11.1.16.50})}
and 
{\rm(\ref{eq;11.1.16.51})}
are isomorphisms.
\hfill\qed
\end{lem}

Let $a\leq b$.
Let $M>\!>\max\{|a|,|b|\}$.
We have the following natural
commutative diagrams:
\begin{equation}
 \begin{CD}
 \Ker_1@>>>
 \nbigm\otimes\II_2^{b,M+1} 
 @>{\nbign\otimes\id+\id\otimes \nbign''_{\II}}>>
 \nbigm\otimes\II_2^{a,M+1}
 @>>>\Cok_1\\
 @V{\varphi_1}VV @VVV @VVV @V{\varphi_2}VV\\
 \Ker_2@>>>
 \nbigm\otimes\II_2^{b,M} 
 @>{\nbign\otimes\id+\id\otimes \nbign''_{\II}}>>
 \nbigm\otimes\II_2^{a,M}
 @>>>\Cok_2
 \end{CD}
\end{equation}
\begin{equation}
 \begin{CD}
 \Ker_3@>>>
 \nbigm\otimes\II_2^{-M,b} 
 @>{\nbign\otimes\id+\id\otimes \nbign''_{\II}}>>
 \nbigm\otimes\II_2^{-M,a}
 @>>>\Cok_3\\
 @V{\varphi_3}VV @VVV @VVV @V{\varphi_4}VV\\
 \Ker_4@>>>
 \nbigm\otimes\II_2^{-M-1,b} 
 @>{\nbign\otimes\id+\id\otimes \nbign''_{\II}}>>
 \nbigm\otimes\II_2^{-M-1,a}
 @>>>\Cok_4
 \end{CD}
\end{equation}

We easily obtain the following lemma.
\begin{lem}
\label{lem;10.8.6.2}
The morphisms $\varphi_i$ $(i=2,3)$
are isomorphisms.
\hfill\qed
\end{lem}

We mention variants of Lemma \ref{lem;11.1.16.53}.
Let $M>\!>\max\{|a|,|b|\}$.
Let $\psi_1$ be the composite of the following morphisms:
\[
 \begin{CD}
 \nbigm\otimes\II_2^{-M,b}
 @>{\rm projection}>>
 \nbigm\otimes\II_2^{-M,a} 
 @>{\nbign\otimes\id+\id\otimes\nbign''_{\II}}>>
 \nbigm\otimes\II_2^{-M,a} 
 \end{CD}
\]
Let $\psi_2$ be the composite of the following morphisms:
\[
 \begin{CD}
 \nbigm\otimes\II_2^{b,M}
@>{\nbign\otimes\id+\id\otimes\nbign''_{\II}}>>
 \nbigm\otimes\II_2^{b,M}
@>{\rm inclusion}>>
 \nbigm\otimes\II_2^{a,M}
 \end{CD}
\]
The natural inclusion
$\nbigm\otimes\II_2^{-M,-M+1}
\lrarr
 \nbigm\otimes\II_2^{-M,a}$
induces 
$\nbigm\otimes\II_2^{-M,-M+1}
\lrarr
 \Cok\psi_1$.
The natural projection
$\nbigm\otimes\II_2^{b,M}
\lrarr
\nbigm\otimes\II_2^{M-1,M}$
induces
$\Ker\psi_2\lrarr
 \nbigm\otimes\II_2^{M-1,M}$.
The following lemma can be checked easily.
\begin{lem}
\label{lem;11.1.16.101}
The morphisms
$\nbigm\otimes\II_2^{-M,-M+1}
\lrarr
 \Cok\psi_1$
and
$\Ker\psi_2\lrarr
 \nbigm\otimes\II_2^{M-1,M}$
are isomorphisms.
\hfill\qed
\end{lem}

\chapter{Canonical prolongations}
\label{section;11.4.9.2}

Let $M$ be a holonomic $D$-module on 
a complex manifold $X$ with a hypersurface $D$.
We have the canonically defined $D$-modules
$M[\ast D]:=M\otimes\nbigo_X(\ast D)$ and 
$M[!D]:=\DDD_X\bigl(
 (\DDD_XM)(\ast D)\bigr)$,
where $\DDD_X$ denotes the dual of $D_X$-modules.
We would like to define such prolongations
for $\nbigr_X$-triples in the case that
$D$ is given as
$\{t=0\}$ for a coordinate function $t$,
and that $M$ is strictly specializable along $t$.
\index{sheaf \mbox{$M[\bikkuri D]$}}
\index{sheaf $M[\ast D]$}

\section{Canonical prolongations of $\nbigr(\ast t)$-modules}
\label{subsection;11.2.15.1}

First, we study the canonical prolongations of
$\nbigr(\ast t)$-modules,
for a coordinate function $t$,
i.e., $t$ is a holomorphic function
whose derivative $dt$ is nowhere vanishing.
Then, the canonical prolongations
of $\nbigr$-modules $\nbigm$ are given
as the canonical prolongation of
$\nbigm(\ast t)$,
as in \S\ref{subsection;11.4.3.2}.

\subsection{Strictly specializable $\nbigr(\ast t)$-modules}
\label{subsection;11.4.3.3}

Let $X_0$ be a complex manifold,
and let $X$ be an open subset of 
$X_0\times\cnum_t$.
We shall often identify $X_0$ and $X_0\times\{0\}$.
Let $H$ be a hypersurface of $X$
such that $(X_0\times\{0\})\not\subset H$.
Let $\nbigm$ be a coherent
$\nbigr_{X(\ast H)}(\ast t)$-module
which is strictly specializable along $t$.
Let $\Vzero=\bigl(
\Vzero_a\,\big|\,a\in\real \bigr)$
be the $V$-filtration of $\nbigm$
at $\lambda_0$ on a neighbourhood
$\nbigxzero$ of $\{\lambda_0\}\times X$,
i.e.,
(i) each $\Vzero_a(\nbigm)$ is
 $V_0\nbigr_{X(\ast H)}$-coherent,
(ii) each $\Gr^{\Vzero}_a(\nbigm)$ is strict,
(iii) $t\,\Vzero_a(\nbigm)=\Vzero_{a-1}(\nbigm)$,
 $\deldel_t\Vzero_a(\nbigm)\subset\Vzero_{a+1}(\nbigm)$,
(iv) for each $P\in X_0$,
 there exists a discrete subset
 $\nbigs\subset\real\times\cnum$
 such that
\[
 \prod_{\substack{u\in\nbigs\\
 \paramap(\lambda_0,u)=a}}
 \bigl(-\deldel_t t+\eigenmap(\lambda,u)\bigr)
\]
is nilpotent on $\Gr^{\Vzero}_a(\nbigm)$ around $P$.
Let $\KMS(\nbigm,P)$ denote the minimum
among such $\nbigs$.
We have the decomposition
\[
 \Gr^{\Vzero}_a(\nbigm)
=\bigoplus_{\substack{u\in\real\times\cnum\\
 \paramap(\lambda_0,u)=a}}
 \psizero_{t,u}(\nbigm),
\]
such that $-\deldel_tt+\eigenmap(\lambda,u)$
is locally nilpotent on $\psizero_{t,u}(\nbigm)$.
For a section $f\in\nbigm$,
we have $\inf\bigl\{
 a\in\real\,\big|\,f\in\Vzero_a(\nbigm)
 \bigr\}$ in $\real\cup\{-\infty\}$,
which is denoted by $\deg^{\Vzero}(f)$.
It is also denoted by $\deg(f)$,
if there is no risk of confusion.

\subsection{The $\nbigr$-module $\nbigm[\ast t]$}

\index{sheaf $\nbigm[\ast t]$}
Let $\nbigm^{(\lambda_0)}[\ast t]$
denote the $\nbigr_{X(\ast H)}$-submodule of
$\nbigm_{|\nbigxzero}$
generated by
$\Vzero_{0}\nbigm$.

\begin{lem}
\mbox{{}}
\begin{itemize}
\item
$\nbigm^{(\lambda_0)}[\ast t]$
is $\nbigr_{X(\ast H)}$-coherent,
and strictly specializable.
The $V$-filtration is given by
\begin{equation}
\label{eq;10.4.19.1}
 \Vzero_a(\nbigm^{(\lambda_0)}[\ast t])
=\nbigm^{(\lambda_0)}[\ast t]\cap \Vzero_a(\nbigm).
\end{equation}
\item
$\Vzero_a$ are described as follows:
\[
 \Vzero_{a}(\nbigm^{(\lambda_0)}[\ast t])=
 \left\{
 \begin{array}{ll}
 \Vzero_a(\nbigm) & (a\leq 0)\\
\deldel_t\Vzero_{a-1}\bigl(
 \nbigm^{(\lambda_0)}[\ast t]
 \bigr)
+\Vzero_{<a}\bigl(
 \nbigm^{(\lambda_0)}[\ast t]
 \bigr) & (a>0) 
 \end{array}
 \right.
\]
In particular,
we have
\begin{equation}
\label{eq;10.4.19.3}
 \deldel_t:
 \psizero_u(\nbigmzero[\ast t])
\stackrel{\simeq}{\lrarr}
 \psizero_{u+\vecdelta}(\nbigmzero[\ast t]),
\quad\quad
 \paramap(\lambda_0,u)> -1
\end{equation}
\begin{equation}
\label{eq;10.4.19.4}
 t: \psizero_u(\nbigmzero[\ast t])
\stackrel{\simeq}{\lrarr}
 \psizero_{u-\vecdelta}(\nbigmzero[\ast t]),
\quad\quad
 \paramap(\lambda_0,u)\leq 0
\end{equation}
\end{itemize}
\end{lem}
\pf
It is clear that
$\nbigmzero[\ast t]$ is 
$\nbigr_{X(\ast H)}$-coherent.
Let $\Vzero(\nbigmzero[\ast t])$
be the filtration given by
(\ref{eq;10.4.19.1}).
We clearly have
(i) $t\Vzero_a(\nbigmzero[\ast t])\subset 
\Vzero_{a-1}(\nbigmzero[\ast t])$,
(ii) $\deldel_t\Vzero_a(\nbigmzero[\ast t])
 \subset \Vzero_{a+1}(\nbigmzero[\ast t])$,
(iii) 
$\Gr^{\Vzero}_{a}(\nbigm^{(\lambda_0)}[\ast t])
\subset
 \Gr^{\Vzero}_a(\nbigm)$
is strict.
We have the decomposition
\[
 \Gr^{\Vzero}_{a}(\nbigm^{(\lambda_0)}[\ast t])
=\bigoplus_{\substack{u\in \real\times\cnum\\
 \paramap(\lambda_0,u)=a }}
 \psizero_u(\nbigm^{(\lambda_0)}[\ast t])
\]
such that
$-\deldel_tt+\eigenmap(\lambda,u)$
are nilpotent on 
$\psizero_u(\nbigm^{(\lambda_0)}[\ast t])$.

Let $a\leq 0$.
Because $\Vzero_a\!(\nbigm)\!\subset\!
 \nbigmzero[\ast t]$,
we have
$\Vzero_a\!(\!\nbigmzero[\ast t]\!)
\!=\!\Vzero_a\!(\nbigm\!)$.
They are $V_0\nbigr_{X(\ast H)}$-coherent,
and we have
$t\cdot\Vzero_a(\nbigmzero[\ast t])
=\Vzero_{a-1}(\nbigmzero[\ast t])$
for $a\leq 0$.

Take $u\in\KMS(\nbigm)$ such that
$\eigenmap(\lambda_0,u)\neq 0$.
The induced action of $-\deldel_tt$
on $\psizero_{u}(\nbigmzero[\ast t])$
is invertible,
and hence 
$t:\psizero_u(\nbigmzero[\ast t])
\lrarr \psizero_{u-\vecdelta}(\nbigmzero[\ast t])$
is injective,
and 
$\deldel_t:
 \psizero_{u-\vecdelta}(\nbigmzero[\ast t])
\lrarr \psizero_{u}(\nbigmzero[\ast t])$
is surjective.

If $\paramap(\lambda_0,u)>-1$,
the function $\eigenmap(\lambda,u+\vecdelta)$
of $\lambda$ is not constantly $0$.
The action of
$-t\deldel_t+\eigenmap(\lambda,u+\vecdelta)$
on $\psizero_u(\nbigm)$ is locally nilpotent.
Hence,
the action of $-t\deldel_t$
on $\psizero_u(\nbigm)$ is injective.
Therefore,
$\deldel_t:\psizero_u(\nbigm)\lrarr
 \psizero_{u+\vecdelta}(\nbigm)$
is injective.
It implies that
$\deldel_t:
 \Gr^{\Vzero}_a(\nbigm)
\lrarr
 \Gr^{\Vzero}_{a+1}(\nbigm)$
is injective for any $a>-1$.
Hence, if $[g]\neq 0$ in
$\Gr^{\Vzero}_a(\nbigm)$
for some $-1<a\leq 0$,
then $[\deldel_t^Ng]\neq 0$
in $\Gr^{\Vzero}_{a+N}(\nbigm)$.

For any $f\in\nbigmzero[\ast t]$,
we have an expression
$f=\sum_{j=0}^N\deldel_t^jf_j$,
where $f_j\in\Vzero_0\nbigmzero[\ast t]$.
We can assume
$-1<\deg(f_j)\leq 0$ for $j\geq 1$.
If $f$ is not contained in $\Vzero_{0}$,
we have
$j+\deg(f_j)<N+\deg(f_N)$
for any $j<N$.

Let us prove
$\Vzero_a\nbigm^{(\lambda_0)}[\ast t]
=\Vzero_{<a}\nbigmzero[\ast t]
+\deldel_t\Vzero_{a-1}
 \nbigm^{(\lambda_0)}[\ast t]$
for $a>0$.
Take $f\in \Vzero_a\nbigmzero[\ast t]$
such that the induced section $[f]$ of 
$\Gr^{\Vzero}_a\bigl(
 \nbigmzero[\ast t]\bigr)$
is non-zero.
Then, we have $N+\deg(f_N)=a$,
and $[f]=[\deldel_t^Nf_N]$
in $\Gr^{\Vzero}_a\nbigmzero[\ast t]$.
It follows that
$\Vzero_a\!\nbigm^{(\lambda_0)}[\ast t]
\!=\!\Vzero_{<a}\!\nbigm[\ast t]
+\deldel_t\Vzero_{a-1}
 \nbigm^{(\lambda_0)}[\ast\! t]$.
In particular,
we obtain that
$\Vzero_a\!\nbigm[\ast t]$ $(a>0)$
are $V_0\nbigr_{X(\ast H)}$-coherent.
Hence the filtration $\Vzero$ 
gives a $V$-filtration.
\hfill\qed

\begin{lem}
$\nbigr_{X}\otimes_{V_0\nbigr_{X}}\Vzero_0\nbigm
\simeq
 \nbigmzero[\ast t]$ naturally
around any$\{\lambda_0\}\times X$.
\end{lem}
\pf
We put $\nbigm_1:=\nbigr_X\otimes_{V_0\nbigr_X}
 \Vzero_0\nbigm$.
We have a naturally defined surjection
$\nbigm_1\lrarr \nbigmzero[\ast t]$.
The composite of the morphisms
$\Vzero_0\nbigm
\lrarr
 \nbigm_1
\lrarr \nbigm$ is injective,
where the first morphism is given by
$m\longmapsto 1\otimes m$.
For $a\leq 0$,
let $\Vzero_a\bigl(\nbigm_1\bigr)$
denote the image of $\Vzero_a\nbigm$.
For $a>0$,
we define
$\Vzero_a\bigl(\nbigm_1\bigr):=
 \sum_{\substack{b+n\leq a\\
b\leq 0,n\in\seisuu_{\geq 0}}}
 \deldel_t^n\Vzero_b\bigl(
\nbigm_1\bigr)$.
For $a>0$, we have $-1<a_0\leq 0$
such that $n:=a-a_0\in\seisuu$.
Then, the following diagram is commutative:
\[
 \begin{CD}
 \Gr^{\Vzero}_{a_0}\nbigm_1
 @>{\simeq}>>
 \Gr^{\Vzero}_{a_0}\nbigm[\ast t]\\ 
 @V{\deldel_t^n}VV @V{\deldel_t^n}V{\simeq}V \\
 \Gr^{\Vzero}_{a}\nbigm_1
 @>>>
 \Gr^{\Vzero}_{a}\nbigm[\ast t]\\ 
 \end{CD}
\]
Hence, we obtain 
$\Gr^{\Vzero}_{a_0}\nbigm_1
\lrarr
 \Gr^{\Vzero}_{a}\nbigm_1$
is injective.
It is also surjective by construction.
Hence, 
$\Gr^{\Vzero}_{a_0}\nbigm_1
\lrarr
 \Gr^{\Vzero}_{a}\nbigm_1$
and
$\Gr^{\Vzero}_{a}\nbigm_1
\lrarr
 \Gr^{\Vzero}_{a}\nbigm[\ast t]$
are isomorphisms.
Then, we obtain that
$\nbigm_1\lrarr\nbigm[\ast t]$
is an isomorphism.
\hfill\qed

\vspace{.1in}
The following obvious lemma
will be used implicitly.

\begin{lem}
If $\lambda_0$ is generic,
i.e.,
$\eigenmap(\lambda_0,u)\neq 0$
for any $u\in \KMS(\nbigm,P)$,
then
$\nbigm^{(\lambda_0)}[\ast t]
=\nbigm_{|\nbigxzero}$.
\hfill\qed
\end{lem}

We obtain a globally defined module.

\begin{lem}
Let $\lambda_1$ be sufficiently close to $\lambda_0$,
and let $\nbigx^{(\lambda_1)}\subset\nbigxzero$
be a neighbourhood of $\{\lambda_1\}\times X$.
Then, 
we have
$\nbigm^{(\lambda_1)}[\ast t]
=\nbigm^{(\lambda_0)}[\ast t]_{|
 \nbigx^{(\lambda_1)}}$.
Therefore,
we have a globally defined $\nbigr_{X(\ast H)}$-module
$\nbigm[\ast t]$.
\end{lem}
\pf
Take a sufficiently small $\epsilon>0$.
Then, $\nbigm^{(\lambda_0)}[\ast t]$
is generated by $\Vzero_{\epsilon}\!(\nbigm)$
over $\nbigr_{X(\ast H)}$.
Because 
$\Vzero_{\epsilon}(\nbigm)_{|\nbigx^{(\lambda_1)}}
=V^{(\lambda_1)}_{\epsilon}(\nbigm)$,
we are done.
\hfill\qed

\vspace{.1in}
The module $\nbigm[\ast t]$ has the following
universal property.

\begin{lem}
Let $\nbigm_1$ be a coherent
$\nbigr_{X(\ast H)}$-module
such that 
(i) $\nbigm_1(\ast t)=\nbigm$,
(ii) $\nbigm_1$ is strictly specializable along $t$.
Then, the natural morphism
$\nbigm_1\lrarr\nbigm_1(\ast t)=\nbigm$
factors through $\nbigm[\ast t]$.
\end{lem}
\pf
We have only to argue around
$\{\lambda_0\}\times X$ for any $\lambda_0$.
By a standard argument,
we can prove that 
the image of $\Vzero_0\nbigm_1$ in $\nbigm$
is contained in $\Vzero_0\nbigm$.
Because $\nbigm_1$ is generated by
$\Vzero_0\nbigm_1$ over $\nbigr_{X(\ast H)}$,
the image of $\nbigm_1$ in $\nbigm$
is contained in $\nbigm[\ast t]$.
\hfill\qed

\begin{cor}
\label{cor;13.5.8.30}
$\nbigm[\ast t]$ is independent of
the choice of a decomposition
$X=X_0\times\cnum_t$.
Moreover,
for any nowhere vanishing function $A$,
we naturally have
$\nbigm[\ast t]\simeq\nbigm[\ast(At)]$.
\hfill\qed
\end{cor}

The following lemma is clear by construction.
\begin{lem}
If $\nbigm$ is integrable,
$\nbigm[\ast t]$ is also naturally integrable.
\hfill\qed
\end{lem}

\begin{lem}
If $\nbigm$ is strict,
then $\nbigm[\ast t]$ is also strict.
\end{lem}
\pf
Because $\nbigm[\ast t]\subset\nbigm(\ast t)$,
the claim is clear.
\hfill\qed

\subsection{The $\nbigr$-module $\nbigm[!t]$}

\index{sheaf \mbox{$\nbigm[\bikkuri t]$}} 

Let $\nbigmzero[!t]$ be the $\nbigr_{X(\ast H)}$-module
on $\nbigxzero$ defined as follows:
\[
 \nbigm^{(\lambda_0)}[!t]:=
 \nbigr_{X(\ast H)}\otimes_{V_0\nbigr_{X(\ast H)}}
 \Vzero_{<0}\nbigm
\]
Because the composite of the natural morphisms
$\Vzero_{<0}\nbigm\lrarr 
 \nbigm^{(\lambda_0)}[!t]\lrarr\nbigm$
is injective,
we can naturally regard
$\Vzero_{<0}\nbigm$ 
as a $V_0\nbigr_{X(\ast H)}$-submodule of
$\nbigm^{(\lambda_0)}[!t]$,
which will be used implicitly.

\begin{lem}
\mbox{{}}\label{lem;10.8.20.3}
\begin{itemize}
\item
$\nbigm^{(\lambda_0)}[!t]$ is
$\nbigr_{X(\ast H)}$-coherent
and strictly specializable along $t$.
\item
We have a natural isomorphism
$\Vzero_{a}\bigl(\nbigmzero[!t]\bigr)
 \simeq 
 \Vzero_{a}(\nbigm)$
for any $a<0$.
\item
The natural morphism
$\deldel_t:\psizero_{-\vecdelta}(\nbigmzero[!t])
 \lrarr\psizero_{0}(\nbigmzero[!t])$
is an isomorphism.
\item
We have a globally defined
$\nbigr_{X(\ast H)}$-module
$\nbigm[!t]$.
\end{itemize}
\end{lem}
\pf
We have only to consider the issues
locally on $X$,
which we will implicitly use in the following argument.
We consider the following $V_0\nbigr_{X(\ast H)}$-submodules:
\begin{equation}
 \label{eq;13.3.27.1}
\Vzero_a\nbigm[!t]:=\left\{
 \begin{array}{ll}
 \Vzero_a\nbigm & (a<0)\\
 \deldel_t\Vzero_{a-1}\nbigm[!t]
+\Vzero_{<a}\nbigm[!t] & (a\geq 0)
 \end{array}
 \right.
\end{equation}
We clearly have 
(i) $t\Vzero_{a}\subset \Vzero_{a-1}$,
(ii) $\deldel_t\Vzero_{a}\subset \Vzero_{a+1}$,
(iii) $\Gr^{\Vzero}_a\!\!\!\nbigm[!t]
\simeq
 \Gr^{\Vzero}_a\!\!\!\nbigm$
for $a<0$,
(iv) the natural morphism
$\deldel_t:
 \Gr^{\Vzero}_{a}\nbigm[!t]
\lrarr
 \Gr^{\Vzero}_{a+1}\nbigm[!t]$
is onto for $a\geq -1$.
Let us prove the claims
(A1) $\Gr^{\Vzero}_a(\nbigm[!t])$
is strict,
(A2) the induced action of
\[
 \prod_{\substack{u\in\KMS(\nbigm)\\
 \paramap(\lambda_0,u)=a}}
 (-\deldel t+\eigenmap(\lambda,u))
\]
is nilpotent on $\Gr^{\Vzero}_a(\nbigm[!t])$.
If $a<0$, they are clear by the construction.

Let us consider the case $a=0$.
If $f\in \Vzero_0\nbigm[!t]$,
we have
$f=f_0+\deldel_t\otimes f_1$,
where $f_0\in\Vzero_{<0}$
and $f_1\in\Vzero_{-1}$.
For $\lambda\neq 0$,
we put 
$\nbigm^{\lambda}:=
 \nbigm_{|\{\lambda\}\times X}$,
which can naturally be regarded as a $D_X$-module.
Note that 
we naturally have 
$\nbigmlambda[!t]\simeq
 D_X\otimes_{V_0D_X}V_{<0}\nbigmlambda$
for $D$-modules.
(See Appendix below.)
We have the following commutative diagram
for any generic $\lambda$:
\[
 \begin{CD}
 \Vzero_0\nbigm[!t]
 @>{|\lambda}>>
 V_0\bigl(\nbigmlambda[!t]\bigr)\\ 
 @A{\deldel_t}AA @A{\deldel_t}AA \\
 \Vzero_{-1}\nbigm
 @>>{|\lambda}>
 V_{-1}\nbigmlambda
 \end{CD}
\]
If $\deldel_t(f_{1|\lambda})
 \in V_{<0}\nbigmlambda[!t]$,
we obtain
$f_{1|\lambda}\in
 V_{<-1}\nbigmlambda$.
Hence, if
$\deldel_t\otimes f_1\in
 \Vzero_{<0}\nbigm[!t]$,
then we have
$f_{1|\lambda}\in V_{<-1}\nbigmlambda$
for any generic $\lambda\neq 0$.
We obtain
$f_{1}\in\Vzero_{<-1}\nbigm$,
because $\Gr^{\Vzero}_{-1}\nbigm$ is strict.
Hence, the map
$\deldel_t:\Gr^{\Vzero}_{-1}\nbigm[!t]
\lrarr
 \Gr^{\Vzero}_0\nbigm[!t]$
is an isomorphism.
It implies (A1) and (A2)
for $\Gr^{\Vzero}_0(\nbigm[!t])$.

Let us consider the case $a>0$.
For $0<a\leq 1$,
we obtain the injectivity of
$t:\Gr^{\Vzero}_a\nbigm[!t]
\lrarr \Gr^{\Vzero}_{a-1}\nbigm[!t]$
from the injectivity of $t\deldel_t$ on
$\Gr^{\Vzero}_{a-1}\nbigm[!t]$.
Then, we obtain (A1) and (A2) for 
$\Gr^{\Vzero}_a\nbigm[!t]$.
By an easy inductive argument,
we obtain (A1) and (A2)
for any $a\geq 0$.
We also obtain that 
$\deldel_t:
 \Gr^{\Vzero}_{a}(\nbigm)
\lrarr \Gr^{\Vzero}_{a+1}(\nbigm)$
is an isomorphism for $a>-1$.

Because $\Vzero_{<0}\nbigm$ is
$V_0\nbigr_{X(\ast H)}$-coherent,
we obtain that $\nbigm[!t]$ is
$\nbigr_{X(\ast H)}$-coherent.
Note that $V_0\nbigr_{X(\ast H)}$
and $\nbigr_{X(\ast H)}$ are Noetherian.
Because $V_a(\nbigm[!t])$ are
$\nbigo_X(\ast H)$-pseudo-coherent
and locally $V_0\nbigr_{X(\ast H)}$-finitely generated,
it is $V_0\nbigr_{X(\ast H)}$-coherent.
Thus, we obtain that
$\nbigm[!t]$ is strictly specializable
along $t$,
and the $V$-filtration is given by $\Vzero$.

If $\epsilon>0$ is sufficiently small,
a natural morphism
$\nbigr_{X(\ast H)}
 \otimes_{V_0\nbigr_{X(\ast H)}}
\Vzero_{<-\epsilon}\nbigm
 \lrarr \nbigmzero[!t]$
is an isomorphism.
Because
$\Vzero_{<-\epsilon|\nbigx^{(\lambda_1)}}=
V^{(\lambda_1)}_{<-\epsilon}$,
we obtain the global well definedness
of $\nbigm[!t]$.
\hfill\qed

\begin{cor}
For a generic $\lambda$,
the specialization
$(\nbigm[!t])^{\lambda}$
is naturally isomorphic to
$\nbigmlambda[!t]$.
\hfill\qed
\end{cor}

We have the following universality.
\begin{lem}
Let $\nbigm_1$ be a coherent strict
$\nbigr_{X(\ast H)}$-module such that
(i) $\nbigm_1$ is strictly specializable along $t$,
(ii) $\nbigm_1(\ast t)=\nbigm$.
Then, we have a uniquely defined morphism
$\nbigm[!t]\lrarr\nbigm_1$.
\end{lem}
\pf
We have only to consider the issue around
$\{\lambda_0\}\times X$.
We have 
$\Vzero_a(\nbigm_1)=\Vzero_a(\nbigm)$
for any $a<0$.
In particular, 
$\Vzero_{<0}\nbigm\subset\nbigm_1$.
Hence, we obtain the uniquely induced morphism
$\nbigm[!t]=\nbigr_{X(\ast H)}
 \otimes_{V_0\nbigr_{X(\ast H)}}
 \Vzero_{<0}\nbigm
\lrarr \nbigm_1$.
\hfill\qed

\begin{cor}
\label{cor;13.5.8.31}
$\nbigm[!t]$ is independent of
the choice of a decomposition
$X=X_0\times\cnum_t$.
Moreover,
for any nowhere vanishing function $A$,
we naturally have
$\nbigm[!t]\simeq\nbigm[!(At)]$.
\hfill\qed
\end{cor}

\vspace{.1in}
The following lemma is also clear
by construction.
\begin{lem}
If $\nbigm$ is integrable,
$\nbigm[!t]$ is also integrable.
\hfill\qed
\end{lem}

\begin{lem}
If $\nbigm$ is strict,
$\nbigm[!t]$ is also strict.
\end{lem}
\pf
We have only to argue the issue
locally around any point
$(\lambda,P)\in \cnum_{\lambda}\times \{t=0\}$.
Because $\Vzero_{<0}(\nbigm)\subset \nbigm(\ast t)$,
it is strict.
By construction,
each $\psi_{t,u}(\nbigm[!t])$ is strict.
Hence, we obtain that
$\nbigm[!t]$ is strict.
\hfill\qed

\subsubsection{Appendix: The case of $D$-modules}

Let $M$ be a coherent $D_X$-module
with a $V$-filtration along $t$.
(We fix an appropriate total order $\leq_{\cnum}$ on $\cnum$.)
For simplicity, we assume $M=M(\ast t)$.
Assume that we are given
morphisms of $D_{X_0}$-modules
$\psi_{1}M
\stackrel{u}{\lrarr}Q\stackrel{v}{\lrarr}
 \psi_{0}M$
such that $v\circ u$
is equal to the induced map
$\del_t:\psi_{-1}M\lrarr\psi_{0}M$.
Although the following lemma is known,
we give a construction by hand
for our understanding.
\begin{lem}
We have a coherent $D$-module $\Mtilde$
with a $V$-filtration along $t$
such that
(i) $\Mtilde(\ast t)=M$,
(ii) we have $\psi_{0}(\Mtilde)\simeq Q$,
and 
the following commutative diagram:
\[
 \begin{CD}
 \psi_{-1}(\Mtilde)
 @>{\del_t}>>
 \psi_0(\Mtilde)
 @>{t}>>
 \psi_{-1}(\Mtilde) \\
 @V{=}VV @V{\simeq}VV @V{=}VV \\
 \psi_{-1}(M)
 @>{u}>> Q @>{t\circ v}>>
 \psi_{-1}(M)
 \end{CD}
\]
\end{lem}
\pf
In general, let $N$ be a nilpotent endomorphism
of a $D_{X_0}$-module $Q$.
Let $\iota:X_0\times\{0\}\lrarr X$ be the inclusion.
We set $G(Q,N):=\iota_{\ast}Q\otimes_{\cnum}\cnum[\del_t]$
as a sheaf on $X$.
We define the action of $t$ on $G(Q,N)$ by
$t(m\,\del_t^j)=(-j-N)m\,\del_t^{j-1}$.
Because it is nilpotent,
we obtain an $\nbigo_X$-action on $G(Q,N)$.
We define an action of $\del_tt$ on $G(Q,N)$ by
$(\del_tt)(m\,\del_t^j)=(-j-N)m\,\del_t^j$.
For $f\in\nbigo_X$,
let $\mu_f:G(Q,N)\lrarr G(Q,N)$ be
the multiplication of $f$.
We can check that
$[(\del_tt),\mu_f]=\mu_{[\del_tt,f]}$.
Hence, with the above actions,
$G(Q,N)$ is a $V_0D_X$-module.
Let $\del_t:G(Q,N)\lrarr G(Q,N)$ be given by
$\del_t(m\,\del_t^j)=m\,\del_t^{j+1}$.
Then, we have 
$\mu_t\circ\del_t=(\del_tt)-\id$.
We also have
$\del_t\circ\mu_t=(\del_tt)+N\pi_0$,
where $\pi_0$ is the projection of
$G(Q,N)$ onto the $0$-th degree part $Q$.
Hence, we have
$[\del_t,\mu_t]=\id+N\circ\pi_0$.

Let $N$ on $Q$ be the composite of the morphisms
$Q\stackrel{v}\lrarr \psi_{0}(M)
\stackrel{-t}{\lrarr}\psi_{-1}(M)
\stackrel{u}{\lrarr}Q$.
Let $N$ on $\psi_0(M)$ be $-\del_tt$.
Because $N\circ v=v\circ N$,
we naturally have $V_0D_X$-homomorphism
$G(Q,N)\lrarr G(\psi_0(M),N)$,
denoted by $\rho_1$.
Note that $G(\psi_0(M),N)$ is naturally
a direct summand of $M/V_{<0}M$
as $V_0D_X$-module.
Let $\rho_2$ be the composite of
the $V_0D_X$-homomorphisms
$M\lrarr M/V_{<0}M
\lrarr G(\psi_0(M),N)$.
We consider the $V_0D_X$-submodule $\Mtilde$
of $G(Q,N)\oplus M$ given as follows:
\[
 \Mtilde=\bigl\{
 (m_1,m_2)\in G(Q,N)\oplus M\,\big|\,
 \rho_1(m_1)=\rho_2(m_2)
 \bigr\}
\]

We shall make a $\del_t$-action on $\Mtilde$.
Let $L>_{\cnum}1$.
We have the decomposition
$M/V_{-L}M
\simeq
 \bigoplus_{-L<_{\cnum} b}
 \psi_bM$
of the $p^{-1}D_{X_0}$-module
compatible with the action of $t\del_t$,
where $p:X\lrarr X_0$ be the projection.
Let $\pi_{-1}$ denote 
the projection 
$M\lrarr \psi_{-1}$,
which is independent of a choice of $L$.
Then, let $\del_t:\Mtilde\lrarr\Mtilde$
be the $\cnum$-linear morphism 
given as follows:
\[
 \del_t(m_1,m_2)=\bigl(
 \del_tm_1+u(\pi_{-1}(m_2)),\,
 \del_tm_2
 \bigr)
\]
By construction,
we have $[\del_t,P]=0$ for $P\in p^{-1}\nbigo_{X_0}$.
For $f\in\nbigo_X$,
let $\mu_f:\Mtilde\lrarr\Mtilde$ be the multiplication of $f$.
We can check $[\del_t,\mu_f]=\mu_{\del_tf}$
by a direct computation.
Hence we obtain the structure of
an $D_X$-module on $\Mtilde$.
It is easy to observe that
$\Mtilde$ is the desired one.
\hfill\qed

\vspace{.1in}
Applying the above construction,
in the case $Q=\psi_{-1}M$,
$u=\id$ and $v=\del_t$,
we obtain the following.
\begin{cor}
We have a coherent $D_X$-module $M[!t]$
with a $V$-filtration along $t$,
such that
$\del_t:\psi_{-1}(M[!t])\lrarr\psi_0(M[!t])$
is an isomorphism.
\hfill\qed
\end{cor}

We have a natural isomorphism
$D_X\otimes_{V_0D_X}V_{<0}M\simeq M[!t]$.
Indeed, the natural map
$D_X\otimes_{V_0D_X}V_{<0}M\lrarr M[!t]$
is surjective by the construction.
We define a filtration $V$ of 
$D_X\otimes_{V_0D_X}V_{<0}M$
as in (\ref{eq;13.3.27.1}).
The natural map preserves the filtration $V$.
We have 
$V_{a}\bigl(
D_X\otimes_{V_0D_X}V_{<0}M \bigr)
\simeq
V_{a}(M[!t])
\simeq
 V_a(M)$ for $a<0$.
We can prove that the induced map
on $\Gr^V_a$ is an isomorphism for each $a\geq 0$,
by using the argument in Lemma \ref{lem;10.8.20.3}.
Similarly,
we obtain a natural isomorphism
$D_X\otimes_{V_0D_X}V_0M\simeq M=M(\ast t)$.

\subsection{Characterization}

Let $\nbigm$ be any $\nbigr_{X(\ast H)}(\ast t)$-module
which is strictly specializable along $t$.

\begin{lem}
\mbox{{}}\label{lem;10.4.7.1}
Let $\nbigq$ be a coherent $\nbigr_{X(\ast H)}$-module
which is strictly specializable along $t$
such that 
$\nbigq(\ast t)\simeq \nbigm$.
\begin{itemize}
\item
If $t:\psi_0(\nbigq)\simeq
 \psi_{-\vecdelta}(\nbigq)$,
we naturally have $\nbigq\simeq\nbigm[\ast t]$.
\item
If $\deldel_t:\psi_{-\vecdelta}(\nbigq)\simeq
 \psi_{0}(\nbigq)$,
we naturally have $\nbigq\simeq\nbigm[!t]$.
\end{itemize}
\end{lem}
\pf
Let us prove the first claim.
Both $\nbigq$ and $\nbigm[\ast t]$
are naturally $\nbigr_{X(\ast H)}$-submodules of $\nbigm$.
It is easy to observe the coincidence of $\Vzero_0$
for any $\lambda_0$.
Because they are locally generated by $\Vzero_0$,
they are the same.

Let us consider the second claim.
By the universal property,
we have the naturally induced morphism
$\nbigm[!t]\lrarr\nbigq$.
We have only to check that it is an isomorphism
locally around each $(\lambda_0,P)\in\nbigx$.
We have
$\Vzero_{<0}\nbigm[!t]\simeq
 \Vzero_{<0}\nbigq$.
By the condition, we obtain that 
$\psizero_0(\nbigm[!t])\lrarr
\psizero_0(\nbigq)$ is an isomorphism.
Then, we can check that 
$\psizero_u(\nbigm[!t])\lrarr
\psizero_u(\nbigq)$ is an isomorphism
for each $u\in\real\times\cnum$.
Then, the second claim follows.
\hfill\qed

\subsection{Morphisms}

Let $\nbigm_i$ $(i=1,2)$ be 
$\nbigr_{X(\ast H)}(\ast t)$-modules
which are strictly specializable along $t$.
For $\star=\ast,!$,
we have a natural map
given by the localization
\begin{equation}
\label{eq;13.3.27.2}
 \Hom_{\nbigr_{X(\ast H)}}
 \bigl(\nbigm_1[\star t],
 \nbigm_2[\star t]\bigr)
\lrarr
 \Hom_{\nbigr_{X(\ast H)}(\ast t)}
 \bigl(\nbigm_1,\nbigm_2\bigr). 
\end{equation}

\begin{lem}
\label{lem;11.1.17.10}
The map {\rm(\ref{eq;13.3.27.2})} is an isomorphism.
If $\nbigm_i$ are integrable,
we have a bijection
between integrable homomorphisms.
\end{lem}
\pf
Any $\nbigr_{X(\ast H)}$-morphism
$f\!:\!\nbigm_1\!\lrarr\!\nbigm_2$
induces
$\Vzero_{<0}\nbigm_1\lrarr \Vzero_{<0}\nbigm_2$
and $t^{-1}\Vzero_{-1}\nbigm_1
 \lrarr t^{-1}\Vzero_{-1}\nbigm_2$.
We obtain
$\nbigm_1[!t]\lrarr\nbigm_2[!t]$
and $\nbigm_1[\ast t]\lrarr\nbigm_2[\ast t]$,
respectively.
It gives the converse of (\ref{eq;13.3.27.2}).
\hfill\qed

\vspace{.1in}
Let $f:\nbigm_1\lrarr\nbigm_2$ be a morphism.
Recall that it preserves $V$-filtrations
$\Vzero$ at each $\lambda_0$,
and we have the induced morphisms
$\psi_u(f):\psi_u(\nbigm_1)\lrarr
 \psi_u(\nbigm_2)$.
Recall the following lemma.
\begin{lem}
\label{lem;11.1.15.1}
Assume that $f$ is strictly specializable,
i.e., 
$\Cok\psi_u(f)$ are strict.
Then, the following holds:
\begin{itemize}
\item
$f$ is strict with respect to
the filtrations $\Vzero$ for each $\lambda_0$.
\item
$\Ker(f)$, $\Image(f)$ and $\Cok(f)$
are strictly specializable along $t$.
The $V$-filtrations are the same as the naturally
induced filtrations $\Vzero$.
\item
$\psi_u\Ker(f)\simeq\Ker\psi_u(f)$,
$\psi_u\Image(f)\simeq\Image\psi_u(f)$
and $\psi_u\Cok(f)\simeq\Cok\psi_u(f)$
naturally.
\hfill\qed
\end{itemize}
\end{lem}

As we have already remarked,
we have the induced morphism
$f[\star t]:\nbigm_1[\star t]\lrarr \nbigm_2[\star t]$.
\begin{lem}
\label{lem;11.1.17.11}
Under the assumption of Lemma 
{\rm\ref{lem;11.1.15.1}},
we have the following natural isomorphisms:
\[
 (\Ker f)[\star t]\simeq\Ker\bigl(f[\star t]\bigr),
 \quad
 (\Image f)[\star t]\simeq \Image\bigl(f[\star t]\bigr),
 \quad
 (\Cok f)[\star t]\simeq \Cok\bigl(f[\star t]\bigr)
\]
\end{lem}
\pf
We naturally have 
$(\Ker f)[\star t](\ast t)
\simeq(\Ker f)(\ast t)\simeq
 \Ker\bigl(f[\star t]\bigr)(\ast t)$.
Then, we obtain the first isomorphism
by the characterization.
The others are obtained similarly.
\hfill\qed

\vspace{.1in}

Let $\nbigm^{\bullet}$ be a bounded complex of
$\nbigr_{X(\ast H)}(\ast t)$-modules
such that each $\nbigm^{\bullet}$
is strict and strictly specializable along $t$.
\begin{cor}
\label{cor;11.1.15.30}
Assume that
$\nbigh^{\bullet}(\psi_u\nbigm^{\bullet})$
are strict for any $u\in\real\times\cnum$.
Then, the following holds:
\begin{itemize}
\item
The differential of $\nbigm^{\bullet}$
is strict with respect to the $V$-filtrations.
\item
$\nbigh^{\bullet}(\nbigm^{\bullet})$
are strictly specializable along $t$.
The $V$-filtrations are the same as
the naturally induced filtrations $\Vzero$.
\item
$\psi_u\nbigh^{\bullet}
\simeq
 \nbigh^{\bullet}\psi_u$
and
$\nbigh^{\bullet}\bigl(
 \nbigm^{\bullet}[\star t]\bigr)
\simeq
 \nbigh^{\bullet}
 \bigl(\nbigm^{\bullet}\bigr)[\star t]$
naturally.
\hfill\qed
\end{itemize}
\end{cor}

\subsection{Canonical prolongations of $\nbigr$-modules}
\label{subsection;11.4.3.2}

Let $\nbigm$ be any coherent 
$\nbigr_{X(\ast H)}$-module
which is strictly specializable along $t$.
Note that $\nbigm(\ast t)$ is a coherent
$\nbigr_{X(\ast H)}(\ast t)$-module,
and it is strictly specializable along $t$.
For $\star=\ast,!$, we define
$\nbigm[\star t]:=
 \bigl(\nbigm(\ast t)\bigr)[\star t]$.
If $\nbigm$ is integrable,
$\nbigm[\star t]$ are naturally also integrable.

\begin{lem}
We have a natural morphism
$\iota:\nbigm\lrarr\nbigm[\ast t]$.
We have the following naturally defined
isomorphisms:
\begin{equation}
 \label{eq;10.4.19.5}
\begin{array}{l}
\Ker(\iota)\simeq
\Ker\Bigl(
 \psi_0(\nbigm)
\stackrel{t}{\lrarr}
 \psi_{-\vecdelta}(\nbigm)
\Bigr)[\deldel_t]
 \\
 \Cok(\iota)\simeq
\Cok\Bigl(
 \psi_0(\nbigm)
\stackrel{t}{\lrarr}
 \psi_{-\vecdelta}(\nbigm)
\Bigr)[\deldel_t]
\end{array}
\end{equation}
For any $u\not\in\seisuu_{\geq 0}\times\{0\}$,
we have a natural isomorphism
$\psizero_u(\iota):
 \psizero_u(\nbigm)
\simeq
 \psizero_u(\nbigm[\ast t])$.
\end{lem}
\pf
We have a naturally defined morphism
$\nbigm\lrarr\nbigm(\ast t)$,
for which the image of
$\Vzero_0(\nbigm)$ is contained in
$\Vzero_0\bigl(\nbigm(\ast t)\bigr)
=\Vzero_0(\nbigm[\ast t])$.
Because $\nbigm_{|\nbigxzero}$
is generated by
$\Vzero_0(\nbigm)$ over $\nbigr_{X(\ast H)}$,
we obtain
$\nbigm\lrarr\nbigm[\ast t]$.
Let us consider
$\psizero_{u}(\iota): \psizero_u(\nbigm)\lrarr
 \psizero_u(\nbigm[\ast t])$.
It is an isomorphism
if $\paramap(\lambda_0,u)<0$.
If $u\not\in\seisuu_{\geq 0}\times\{0\}$,
we obtain $\psi_u(\iota)$
is an isomorphism
by using an easy induction and 
isomorphisms (\ref{eq;10.4.19.3}).

The kernel and the cokernel of
$\nbigm\lrarr\nbigm[\ast t]$
are naturally isomorphic to those of
\[
 \nbigm\big/\Vzero_{<0}\nbigm
\lrarr
 \nbigm[\ast t]\big/\Vzero_{<0}\nbigm[\ast t],
\]
which are naturally isomorphic
to those for
$\psizero_{0}(\nbigm)[\deldel_t]
 \lrarr\psizero_{0}(\nbigm[\ast t])[\deldel_t]$
by the above consideration.
Note that the induced morphism
$\psizero_0(\iota)$ is naturally identified with
$t:\psizero_0(\nbigm)\lrarr\psizero_{-\vecdelta}(\nbigm)$
under the natural identification
$t:\psizero_0(\nbigm[\ast t])\simeq
 \psizero_{-\vecdelta}(\nbigm[\ast t])
=\psizero_{-\vecdelta}(\nbigm)$.
Then, we obtain (\ref{eq;10.4.19.5}).
\hfill\qed

\vspace{.1in}
Similarly, we have the following lemma.

\begin{lem}
We have a naturally defined morphism
$\iota:\nbigm[!t]\lrarr\nbigm$.
The induced morphism
$\psi_u(\iota)$ is an isomorphism
unless $u\in\seisuu_{\geq 0}\times\{0\}$.
We have natural isomorphisms
\[
\begin{array}{l}
 \Ker(\iota)\simeq
 \Ker\Bigl(
 \psi_{-\vecdelta}(\nbigm)
\stackrel{\deldel_t}\lrarr
 \psi_0(\nbigm)
 \Bigr)[\deldel_t],
 \\
 \Cok(\iota)\simeq
 \Cok\Bigl(
 \psi_{-\vecdelta}(\nbigm)
\stackrel{\deldel_t}\lrarr
 \psi_0(\nbigm)
 \Bigr)[\deldel_t].
\end{array}
\]
\end{lem}
\pf
It can be proved as in the case of
$\nbigm[\ast t]$.
\hfill\qed

\vspace{.1in}
Let $\nbigm_i$ $(i=1,2)$ be strict coherent
$\nbigr_{X(\ast H)}$-modules,
which are strictly specializable along $t$.
We obtain the following lemma from
Lemma \ref{lem;11.1.17.10}.
\begin{lem}
We have a natural morphism
\[
\Hom_{\nbigr_{X(\ast H)}}(\nbigm_1,\nbigm_2)
\lrarr
 \Hom_{\nbigr_{X(\ast H)}}
 (\nbigm_1[\star t],\nbigm_2[\star t]).
\]
If $\nbigm_i$ are integrable,
we have the morphism
between spaces of integrable homomorphisms.
\hfill\qed
\end{lem}
We also obtain the following lemma.
\begin{cor}
We have natural bijections:
\[
 \Hom_{\nbigr_{X(\ast H)}}
 \bigl(\nbigm_1[\ast t],\nbigm_2[\ast t]\bigr)
\stackrel{a_1}{\simeq}
 \Hom_{\nbigr_{X(\ast H)}}
 \bigl(\nbigm_1,\nbigm_2[\ast t]\bigr)
\]
\[
 \Hom_{\nbigr_{X(\ast H)}}
 \bigl(\nbigm_1[!t],\nbigm_2[!t]\bigr)
\stackrel{a_2}{\simeq}
 \Hom_{\nbigr_{X(\ast H)}}
 \bigl(\nbigm_2[!t],\nbigm_2\bigr).
\]
If $\nbigm_i$ are integrable,
we have the bijections of integrable homomorphisms.
\end{cor}
\pf
The morphisms
$\nbigm_1\lrarr\nbigm_1[\ast t]$
and
$\nbigm_2[!t]\lrarr\nbigm_2$
induce $a_i$.
We can construct the converse
by using Lemma \ref{lem;11.1.17.10}.
\hfill\qed

\begin{lem}
\label{lem;13.5.10.1}
Let $\nbigm$ be any coherent $\nbigr_{X(\ast H)}$-module
which is strictly specializable along $t$.
Then, $\nbigm[\star t]$ is independent 
of the choice of a decomposition into the product
$X_0\times\cnum_t$.
If $A$ is nowhere vanishing holomorphic function on $X$,
we have
$\nbigm[\star (At)]\simeq
 \nbigm[\star t]$
naturally.
\end{lem}
\pf
It follows from Corollary \ref{cor;13.5.8.30}
and Corollary \ref{cor;13.5.8.31}.
\hfill\qed

\section{Canonical prolongations of $\nbigr$-triples}
\label{subsection;11.1.15.20}

\subsection{Canonical prolongations of 
$\nbigr(\ast t)$-triples}

\index{canonical prolongation}

We use the setting in \S\ref{subsection;11.4.3.3}.
Let $\nbigt=(\nbigm',\nbigm'',C)$ be 
an $\nbigr_{X(\ast H)}(\ast t)$-triple
which is strictly specializable along $t$,
i.e., $\nbigm'$ and $\nbigm''$ are
$\nbigr_{X(\ast H)}(\ast t)$-modules
which are strictly specializable along $t$,
and $C$ is a sesqui-linear pairing
of $\nbigm'$ and $\nbigm''$,
which is an 
$\nbigr_{X(\ast H)}\otimes
 \sigma^{\ast}\nbigr_{X(\ast H)}$-homomorphism
\[
 C:\nbigm'_{|\vecS\times X}
 \otimes
 \sigma^{\ast}\nbigm''_{|\vecS\times X}
\lrarr\distribution^{\moderate H}_{\vecS\times X/\vecS}(\ast t).
\]

\begin{prop}
\label{prop;10.4.19.10}
We have Hermitian sesqui-linear pairings
\[
 C[!t]:\nbigm'[\ast t]_{|\vecS\times X}\otimes
 \sigma^{\ast}\nbigm''[!t]_{|\vecS\times X}
 \lrarr \distribution^{\moderate H}_{\vecS\times X/\vecS}
\]
\[
 C[\ast t]:\nbigm'[!t]_{|\vecS\times X}\otimes
 \sigma^{\ast}\nbigm''[\ast t]_{|\vecS\times X}
 \lrarr \distribution^{\moderate H}_{\vecS\times X/\vecS}
\]
such that
$C[\star t]_{|\vecS\times(X\setminus\{t=0\})}
=C_{|\vecS\times(X\setminus\{t=0\})}$.
\index{paring $C[\ast t]$}
\index{pairing $C[\bikkuri t]$}
They are determined uniquely
by the conditions.
In particular,
we obtain the following $\nbigr$-triples:
\[
 \nbigt[!t]:=\bigl(
 \nbigm'[\ast t],\nbigm''[!t],C[!t]
 \bigr),
\quad\quad
 \nbigt[\ast t]:=\bigl(
 \nbigm'[!t],\nbigm''[\ast t],C[\ast t]
 \bigr)
\]
\index{$\nbigr$-triple $\nbigt[\ast t]$}
\index{$\nbigr$-triple $\nbigt[\bikkuri t]$}
\end{prop}
\pf
By the uniqueness,
we have only to consider it locally.
We have the given pairing
\[
 V_{<0}\nbigm'_{|\vecS\times X}
 \otimes
 \sigma^{\ast}
 V_{<0}\nbigm''_{|\vecS\times X}
\lrarr
 \distribution^{\moderate H}_{\vecS\times X/\vecS}(\ast t).
\]
Let us observe that it is extended to
the following pairing:
\begin{equation}
 \label{eq;10.4.20.1}
 t^{-1}V_{-1}\nbigm'_{|\vecS\times X}
 \otimes
 \sigma^{\ast}
 V_{<0}\nbigm''_{|\vecS\times X}
\lrarr
 \distribution^{\moderate H}_{\vecS\times X/\vecS}
\end{equation}
Let $X_1\subset X$ be open.
Let $U(\lambda_0)$ be a neighbourhood of $\lambda_0$
in $\cnum_{\lambda}$.
Let $u$ and $v$ be local sections of
$\nbigm'$ and $\nbigm''$
on $U(\lambda_0)\times X_1$
and $U(-\lambda_0)\times X_1$
respectively.
We put $\vecI(\lambda_0):=U(\lambda_0)\cap\vecS$.
Let $\phi$ be a $C^{\infty}$-section of
$p_{\lambda}^{\ast}\Omega^{n,n}_{X}$
on $\vecI(\lambda_0)\times X_1$ 
with compact support
such that $\phi_{|\vecI(\lambda_0)\times\Hhat}=0$.
Note that
$\bigl\langle
 C(u,v),\,|t|^{2s}\phi
 \bigr\rangle$ is well defined
for $s\in\nbigh:=
 \bigl\{s\in\cnum\,\big|\Re s>\!>0\bigr\}$,
and holomorphic for $s$
and continuous for $\lambda\in\vecS$.
By a standard argument,
if $u\in t^{-1}\Vzero_{-1}\nbigm'$
and $v\in V^{(-\lambda_0)}_{<0}\nbigm''$,
it is extended to a function on
$\vecS\times\bigl\{s\in\cnum\,\big|\,
 \Re(s)>-\epsilon
 \bigr\}$ for some $\epsilon>0$.
(See the proof of Lemma 20.10.9 of
\cite{mochi7}, for example.)
We can take the value at $s=0$,
denoted by $\Ctilde(u,v)$.
It gives a section of
$\distribution^{\moderate H}_{\vecS\times X/\vecS}$.
Thus, we obtain (\ref{eq;10.4.20.1}).

\begin{lem}
Let $\lambda_0\in\vecS$.
Assume $\sum_{i=0}^N\del_t^if_i=0$
in $\nbigm'[\ast t]=0$
for $f_i\in \Vzero_{0}\nbigm'[\ast t]$.
Then, for any $g\in V^{(-\lambda_0)}_{<0}\nbigm''$
we have
\begin{equation}
\label{eq;10.4.19.11}
 \sum_{i=0}^N 
 \del_t^i\Ctilde(f_i,\sigma^{\ast}g)
=0.
\end{equation}
Assume $\sum_{j=0}^N \del_t^jg_j=0$
in $\nbigm''[!t]$
for some $g_j\in V^{(-\lambda_0)}_{<0}\nbigm''$.
Then, for any $f\in \Vzero_0\nbigm'[\ast t]$,
we have
\[
 \sum_{j=0}^N
 \delbar^j_t
 \Ctilde(f,\sigma^{\ast}g_j)
=0.
\]
\end{lem}
\pf
Let us prove (\ref{eq;10.4.19.11}).
The other can be proved in a similar way.
We use an induction on $N$.
In the case $N=0$,
the claim is trivial.
If $\sum_{i=0}^N\del_t^if_i=0$,
we have
$f_N\in\Vzero_{-1}\nbigm[\ast t]$,
and hence $\del_tf_N\in \Vzero_0\nbigm[\ast t]$.
By using the hypothesis of the induction,
we obtain
\[
 \sum_{i=0}^{N-2}
 \del_t^{i}
 \Ctilde(f_i,\sigma^{\ast}g)
+\del_t^{N-1}\Ctilde\bigl(
 f_{N-1}+\del_tf_N,\sigma^{\ast}g
 \bigr)=0.
\]
So, we have only to prove that
$\Ctilde(\del_t f,\sigma^{\ast}g)
=\del_t\Ctilde(f,\sigma^{\ast}g)$
for $f\in \Vzero_{-1}(\nbigm')$.
For a test form $\phi$,
we have
\begin{multline}
 \bigl\langle
 \Ctilde(\del_tf,\sigma^{\ast}g),\, \phi
 \bigr\rangle
=\bigl\langle
 C(\del_tf,\sigma^{\ast}g),\,
 |t|^{2s}\phi
 \bigr\rangle_{|s=0}
 \\
=-\bigl\langle
 C(f,\sigma^{\ast}g),\,
 |t|^{2s}\,\del_t\phi
 \bigr\rangle_{|s=0}
-\Bigl(
 s\bigl\langle
 C(f,\sigma^{\ast}g),\,
 |t|^{2s-2}\tbar\,\phi
 \bigr\rangle\Bigr)_{|s=0} \\
=\bigl\langle
 \del_t\Ctilde(f,\sigma^{\ast}g),\,
 \phi\bigr\rangle
-\Bigl(
 s\bigl\langle
 C(f,\sigma^{\ast}g),\,
 |t|^{2s-2}\tbar\,\phi
 \bigr\rangle\Bigr)_{|s=0}
\end{multline}
By using $f\in \Vzero_{-1}\nbigm'$
and $g\in V^{(-\lambda_0)}_{<0}\nbigm''$,
we can prove that
$\bigl\langle C(f,\sigma^{\ast}g),\,
 |t|^{2s-2}\tbar\,\phi\bigr\rangle$
is holomorphic at $s=0$.
Hence, we obtain 
$\Ctilde(\del_t f,\sigma^{\ast}g)
=\del_t\Ctilde(f,\sigma^{\ast}g)$
for $f\in \Vzero_{-1}(\nbigm')$.
\hfill\qed

\vspace{.1in}

Then, we can naturally extend $\Ctilde$
to the pairing $C[!t]$ of $\nbigm'[\ast t]$ and
$\nbigm''[!t]$.
Let $C'$ be another pairing
of $\nbigm'[\ast t]$ and $\nbigm''[!t]$
whose restriction to $\vecS\times(X\setminus \{t=0\})$ 
is equal to $C$.
Then, we can prove that $C_1:=C[!t]-C'$
is $0$ on $\Vzero_0\nbigm'\otimes
 \sigma^{\ast}V^{(-\lambda_0)}_{<0}\nbigm''$
by a standard argument.
(See the proof of Lemma 22.10.8 of \cite{mochi7},
for example.)
Then, we obtain that $C_1=0$
on $\nbigm'[\ast t]\otimes\sigma^{\ast}\nbigm''[!t]$.
Thus, the proof of Proposition 
\ref{prop;10.4.19.10} is finished.
\hfill\qed

\vspace{.1in}
The following lemma is obvious.
\begin{lem}
We have
$\bigl(\nbigt\otimes\newTate(w)\bigr)
 [\star t]
=\nbigt[\star t]\otimes\newTate(w)$
and
$\nbigt^{\ast}[!t]
=\bigl(\nbigt[\ast t]\bigr)^{\ast}$.
\hfill\qed
\end{lem}

\begin{lem}
If $\nbigt$ is integrable,
then $\nbigt[\star t]$ are also
integrable.
\end{lem}
\pf
Let us consider $C[!t]$.
We set 
\[
 C_0(m',\sigma^{\ast}m''):=
 \del_{\theta}C(m',\sigma^{\ast}m'')
-C(\del_{\theta}m',\sigma^{\ast}m'')
-C(m',\sigma^{\ast}\del_{\theta}m'').
\]
It is $0$ outside $\{t=0\}$.
It is standard to prove
the vanishing of $C_0$
on $V_0\otimes\sigma^{\ast}V_{<0}$.
(See the proof of Lemma 22.10.8 of \cite{mochi7},
for example.)
Then, we obtain the vanishing
on $\nbigm'[\ast t]\otimes\sigma^{\ast}\nbigm''[!t]$.
\hfill\qed

\subsection{Morphisms}

Let $\nbigt_i$ $(i=1,2)$ be 
coherent $\nbigr_{X(\ast H)}(\ast t)$-triples,
which are strictly specializable along $t$.
\begin{lem}
\label{lem;13.5.10.2}
Let $\star=\ast,!$.
Morphisms $\nbigt_1\lrarr\nbigt_2$
of $\nbigr_{X(\ast H)}(\ast t)$-triples
bijectively correspond to 
morphisms
$\nbigt_1[\star t]\lrarr\nbigt_2[\star t]$
of $\nbigr_{X(\ast H)}$-triples.
If $\nbigt_i$ are integrable,
we have such bijections
for integrable morphisms.
\end{lem}
\pf
A morphism $\nbigt_1[\star t]\lrarr\nbigt_2[\star t]$
naturally induces
$\nbigt_1\lrarr\nbigt_2$.
Let $\nbigt_1\lrarr\nbigt_2$ be a morphism
of $\nbigr_{X(\ast H)}(\ast t)$-triples.
Let us observe that we have an induced morphism
$\nbigt_1[\star t]\lrarr\nbigt_2[\star t]$.
By Lemma \ref{lem;11.1.17.10},
we have the morphisms of the underlying
$\nbigr_{X(\ast H)}$-modules.
We have only to prove the compatibility
of Hermitian sesqui-linear pairings.
By construction,
they are compatible on
$\Vzero_{0}\otimes\sigma^{\ast}\Vzero_{<0}$
or $\Vzero_{<0}\otimes\sigma^{\ast}\Vzero_0$.
Then, the claim is easy to see.
\hfill\qed

\vspace{.1in}

Let $f:\nbigt_1\lrarr\nbigt_2$ be a morphism.
We have the induced morphism
$f[\star t]:\nbigt_1[\star t]\lrarr\nbigt_2[\star t]$.
We obtain the following lemma from
Lemma \ref{lem;11.1.17.11}.
\begin{lem}
\label{lem;13.5.10.3}
Suppose that $f$ is strictly specializable.
Then, we naturally have
$\Ker (f[\star t])\simeq
 \Ker (f)[\star t]$,
$\Image(f[\star t])\simeq
 \Image (f)[\star t]$,
and 
$\Cok(f[\star t])\simeq
 \Cok(f)[\star t]$.
\hfill\qed
\end{lem}

Let $\nbigt^{\bullet}$ be a bounded complex of
$\nbigr_{X(\ast H)}(\ast t)$-triples
such that each $\nbigt^{p}$
is strictly specializable along $t$.
We obtain the following from Corollary \ref{cor;11.1.15.30}.
\begin{lem}
\label{lem;11.1.17.12}
Assume that
$\nbigh^{\bullet}(\psi_u\nbigt^{\bullet})$
are strict for any $u\in\real\times\cnum$.
Then, the following holds:
\begin{itemize}
\item
$\nbigh^{\bullet}(\nbigt^{\bullet})$
are strictly specializable along $t$.
\item
We have natural isomorphisms
$\psi_u\nbigh^{\bullet}(\nbigt^{\bullet})
\simeq
 \nbigh^{\bullet}\psi_u(\nbigt^{\bullet})$
and 
$\nbigh^{\bullet}\bigl(
 \nbigt^{\bullet}[\star t]\bigr)
\simeq
 \nbigh^{\bullet}
 \bigl(\nbigt^{\bullet}\bigr)[\star t]$.
\hfill\qed
\end{itemize}
\end{lem}

\subsection{Canonical prolongations of $\nbigr$-triples}

Let $\nbigt$ be an $\nbigr_{X(\ast H)}$-triple
which is strictly specializable along $t$.
By applying the previous construction
to the $\nbigr_{X(\ast H)}(\ast t)$-triple
$\nbigt(\ast t)$, 
we obtain $\nbigr_{X(\ast H)}$-triples
$\nbigt[\star t]$ for $\star=\ast,!$.
If $\nbigt$ is integrable,
$\nbigt[\star t]$ are also integrable.

\begin{lem}
We have natural morphisms
$\nbigt[!t]\lrarr\nbigt\lrarr\nbigt[\ast t]$.
\end{lem}
\pf
We consider only the morphism
$\nbigt[!t]\lrarr\nbigt$.
The other can be proved similarly.
We have the morphisms
$\nbigm''[!t]\lrarr \nbigm''$
and 
$\nbigm'\lrarr\nbigm'[\ast t]$.
Let us check that they are compatible
with pairings.
Let $f\in\Vzero_{0}\nbigm'$
and $g\in\Vzero_{<0}\nbigm'$.
We have only to prove
$C[!t](f,\sigma^{\ast}g)
=C(f,\sigma^{\ast}g)$.
As in the proof of Proposition \ref{prop;10.4.19.10},
$\bigl\langle
 C(f,\sigma^{\ast}g),|t|^{2s}\phi
 \bigr\rangle$ is 
extended to a function on
$\vecS\times\bigl\{s\in\cnum\,\big|\,
 \Re(s)>-\epsilon
 \bigr\}$ for some $\epsilon>0$.
Moreover,
we have
$\bigl\langle
 C(f,\sigma^{\ast}g),|t|^{2s}\phi
 \bigr\rangle_{s=0}
=\bigl\langle
 C(f,\sigma^{\ast}g),\phi
 \bigr\rangle$.
(See the proof of Lemma 20.10.9 of
\cite{mochi7}, for example.)
Then, we obtain the desired compatibility
of the pairings.
\hfill\qed

\begin{prop}
$\nbigt[\star t]$ is independent of the choice of
a decomposition into the product $X_0\times \cnum_t$.
For any nowhere vanishing holomorphic function $A$
on $X$,
we have natural isomorphisms
$\nbigt[\star(At)]\simeq\nbigt[\star t]$.
\end{prop}
\pf
We have natural morphisms
$\nbigt[\star(At)]\lrarr\nbigt[\star t]$.
The underlying $\nbigr_{X(\ast H)}$-modules
are isomorphic by Lemma \ref{lem;13.5.10.1}
\hfill\qed

\subsection{Compatibility of canonical
prolongation with push-forward}
\label{subsection;11.1.17.20}

Let $F_0:X_0\lrarr Y_0$
be a morphism of complex manifolds.
Let $F:X\lrarr Y$ be $F_0\times\id$,
where $X=X_0\times\cnum_t$
and $Y=Y_0\times\cnum_t$.
Let $H_Y$ be a hypersurface of $Y$
such that $\{t=0\}\not\subset H_Y$.
We put $H:=F^{-1}(H_Y)$.
Let $\nbigm$ be a coherent
$\nbigr_{X(\ast H)}(\ast t)$-module,
which is strictly specializable along $t$.
Assume the following:
\begin{itemize}
\item
 $F_{\dagger}^i\psi_{u}(\nbigm)$
 are strict for any $u\in\real\times\cnum$.
\end{itemize}
According to \cite{sabbah2},
$F_{\dagger}^i\nbigm$ are strictly
specializable $\nbigr_{Y(\ast H_Y)}(\ast t)$-module,
and we naturally have
$\psi_uF_{\dagger}^i\nbigm\simeq
 F_{\dagger}^i\psi_u\nbigm$.

\begin{lem}
\label{lem;11.1.17.21}
Under the assumption,
we have 
$F_{\dagger}^i(\nbigm[\star t])
\simeq
 \bigl(F_{\dagger}^i\nbigm\bigr) [\star t]$
naturally.
\end{lem}
\pf
It easily follows from
the characterization in
Lemma \ref{lem;10.4.7.1}.
\hfill\qed

\begin{cor}
\label{cor;10.8.3.11}
Let $\nbigt$ be a coherent $\nbigr_{X(\ast H)}(\ast t)$-triple,
which is strictly specializable along $t$.
Assume that
$F_{\dagger}^i\psi_{u}\nbigt$ is strict.
Then, we naturally have
$F_{\dagger}^i\bigl(
 \nbigt[\star t]\bigr)
\simeq
\bigl(F_{\dagger}^i\nbigt\bigr)[\star t]$.
\hfill\qed
\end{cor}

\section{Canonical prolongations across hypersurfaces}
\label{subsection;10.12.27.20}

\subsection{Canonical prolongations across holomorphic functions}

\subsubsection{$\nbigr_X(\ast g)$-modules and $\nbigr_X(\ast g)$-triples}

Let $g$ be a holomorphic function on 
a complex manifold $X$.
Let $\iota_g:X\lrarr X\times\cnum_t$ 
denote the embedding
$\iota_g(x)=(x,g(x))$.

\index{canonical prolongation}

\begin{df}
Let $\nbigm$ 
be a coherent $\nbigr_{X}(\ast g)$-module
which is strictly specializable along $g$.
It is called localizable along $g$,
if there exist $\nbigr_X$-modules
$\nbigm[\star g]$ $(\star=\ast,!)$
such that 
$\iota_{g\dagger}(\nbigm[\star g])
\simeq
 \bigl(
 \iota_{g\dagger}(\nbigm)
\bigr)[\star t]$.
\hfill\qed
\end{df}

\index{localizable}

Note that such $\nbigr_X$-modules
$\nbigm[\star g]$
are uniquely determined up to canonical isomorphisms,
because they are recovered 
as the kernel of the multiplication of $t-g$
on $\iota_{g\dagger}(\nbigm)[\star t]$.
But, in general,
it is not clear whether
$\iota_{g\dagger}(\nbigm)[\star t]$
are strictly specializable along $t-g$.
The notation is verified by the following lemma.
\begin{lem}
If $g$ is a coordinate function,
then $\nbigm[\star g]$ given in {\rm\S\ref{subsection;11.2.15.1}}
satisfy $\iota_{g\dagger}(\nbigm[\star g])
\simeq
 (\iota_{g\dagger}\nbigm)[\star t]$.
\end{lem}
\pf
We may assume that
$X$ is an open subset of $X_1\times\cnum_g$.
We have the bi-holomorphic map
$F:\cnum_y\times\cnum_s\lrarr\cnum_g\times\cnum_t$
given by
$F(y,s)=(y+s,y)$.
It induces a bi-holomorphic map 
$X_1\times\cnum_y\times\cnum_s
\lrarr
 X_1\times\cnum_g\times\cnum_t$,
which is also denoted by $F$.
Let $\iota'_{g}:X\lrarr X_1\times\cnum_y\times\cnum_s$
be given by
$\iota'_g(x,g)=(x,g,0)$.
We have
$F\circ\iota'_g=\iota_g$.
We also have $t\circ F=y$.
Then,
$\iota'_{g\dagger}\nbigm$
is clearly strictly specializable along $y=0$,
and we have
$\iota'_{g\dagger}(\nbigm[\star g])
\simeq
 \iota'_{g\dagger}(\nbigm)[\star y]$.
By using the bi-holomorphic map,
we obtain the claim of the lemma.
\hfill\qed

\begin{lem}
If $\nbigm$ is integrable
and localizable along $g$,
then $\nbigm[\star g]$ is also integrable.
\end{lem}
\pf
Because $\nbigm[\star g]$ is obtained
as the kernel of the multiplication of
$t-g$ on the integrable $(\iota_{g\dagger}\nbigm)[\star t]$,
it is integrable.
\hfill\qed

\begin{lem}
Let $A$ be a nowhere vanishing holomorphic function.
We set $g_1:=A\,g$.
Then, $\nbigm$ is strictly specializable along $g$,
if and only if it is strictly specializable along $g_1$.
Moreover,
it is localizable along $g$,
if and only if
it is localizable along $g_1$,
and we have
$\nbigm[\star g]=\nbigm[\star g_1]$
in that case.
\end{lem}
\pf
Let $F:X\times\cnum_t\lrarr X\times\cnum_s$
be given by
$F(x,t)=\bigl(x,A(x)t\bigr)$.
We have $\iota_{g_1}=F\circ \iota_g$.
Then, the claim is clear.
\hfill\qed

\vspace{.1in}
Similarly,
for $\nbigr_X(\ast g)$-triples,
we use the following.
\begin{df}
Let $\nbigt$ be a coherent $\nbigr_{X}(\ast g)$-triple
which is strictly specializable along $g$.
It is called localizable along $g$,
if there exist $\nbigr_{X}$-triple $\nbigt[\star g]$
for $\star=\ast,!$
with isomorphisms
$\iota_{g\dagger}(\nbigt[\star g])
 \simeq
 (\iota_{g\dagger}\nbigt)[\star t]$.
\hfill\qed
\end{df}
Such $\nbigr_X$-triples are uniquely determined
up to canonical isomorphisms, if they exist.
The following lemma is clear.
\begin{lem}
Let $\nbigt$ be localizable along $g$.
\begin{itemize}
\item
 If $g$ is a coordinate function,
 $\nbigt[\star g]$ are equal to those in
 {\rm\S\ref{subsection;11.1.15.20}}.
\item
 If $\nbigt$ is integrable, then
 $\nbigt[\star g]$ are also integrable.
\item
 $\nbigt[\star g]$ are independent
 of the choice of a decomposition into the product
 $X_0\times\cnum_t$.
 If $A$ is a nowhere vanishing holomorphic function,
 then $\nbigt[\star(Ag)]\simeq \nbigt[\star g]$
 naturally.
\hfill\qed
\end{itemize}
\end{lem}

\subsubsection{Morphisms}

Let $\nbigt_i$ be $\nbigr_X(\ast g)$-modules
which are localizable along $g$.
We obtain the following lemma from
Lemma \ref{lem;13.5.10.2}.
\begin{lem}
\label{lem;13.5.10.10}
Let $\star=\ast,!$.
Morphisms $\nbigt_1\lrarr\nbigt_2$
of $\nbigr_{X}(\ast g)$-triples
bijectively correspond to 
morphisms
$\nbigt_1[\star g]\lrarr\nbigt_2[\star g]$
of $\nbigr_{X}$-triples.
If $\nbigt_i$ are integrable,
we have such bijections
for integrable morphisms.
Similar claims hold for 
$\nbigr_X(\ast g)$-modules.
\hfill\qed
\end{lem}

Let $f:\nbigt_1\lrarr\nbigt_2$ be a morphism.
We have the induced morphism
$f[\star g]:\nbigt_1[\star g]\lrarr\nbigt_2[\star g]$.
We obtain the following lemma from
Lemma \ref{lem;13.5.10.3}.
\begin{lem}
\label{lem;13.5.10.11}
Suppose that $f$ is strictly specializable along $g$,
i.e.,  the induced morphism
$\iota_{g\dagger}\nbigt_1\lrarr\iota_{g\dagger}\nbigt_2$
is strictly specializable along $t$.
Then, 
$\Ker f$,
$\Image f$
and $\Cok f$ are strictly specializable and localizable
along $g$,
and we naturally have
\[
\Ker (f[\star g])\simeq
 \Ker (f)[\star g],\,\,\,
\Image(f[\star g])\simeq
 \Image (f)[\star g],\,\,\,
\Cok(f[\star g])\simeq
 \Cok(f)[\star g].
\]
\hfill\qed
\end{lem}

Let $\nbigt^{\bullet}$ be a bounded complex of
$\nbigr_{X}(\ast g)$-triples
such that each $\nbigt^{p}$
is localizable along $g$.
We obtain the following from Lemma \ref{lem;11.1.17.12}.
\begin{lem}
\label{lem;13.5.10.13}
Assume that
$\nbigh^{\bullet}(\psi_{g,u}\nbigt^{\bullet})$
are strict for any $u\in\real\times\cnum$.
Then, the following holds:
\begin{itemize}
\item
$\nbigh^{\bullet}(\nbigt^{\bullet})$
are strictly specializable and localizable along $g$.
\item
We have natural isomorphisms
$\psi_u\nbigh^{\bullet}(\nbigt^{\bullet})
\simeq
 \nbigh^{\bullet}\psi_u(\nbigt^{\bullet})$
and 
$\nbigh^{\bullet}\bigl(
 \nbigt^{\bullet}[\star g]\bigr)
\simeq
 \nbigh^{\bullet}
 \bigl(\nbigt^{\bullet}\bigr)[\star g]$.
\hfill\qed
\end{itemize}
\end{lem}

\subsubsection{Compatibility with the push-forward}
Let $F:X\lrarr Y$ be a morphism of complex manifolds.
Let $g_Y$ be a holomorphic function on $Y$,
and we set $g_X:=g_Y\circ F$.
Let $\nbigm$ be a coherent
$\nbigr_{X}(\ast g_X)$-module
which is strictly specializable along $g_X$.
Assume the following:
\begin{itemize}
\item
 The support of $\nbigm$ is proper over $Y$.
\item
 $F_{\dagger}^i\psitilde_{g_X,u}(\nbigm)$
 are strict for any $u\in\real\times\cnum$.
\end{itemize}
According to \cite{sabbah2},
$F_{\dagger}^i\nbigm$ are strictly
specializable $\nbigr_{Y}(\ast g_Y)$-modules.
\begin{lem}
Assume that
$\nbigm$ is localizable along $g_X$.
Then, 
$(F^i_{\dagger}\nbigm)$ is localizable along $g_Y$,
and we have a natural isomorphism
$F^i_{\dagger}(\nbigm[\star g_X])
\simeq (F^i_{\dagger}\nbigm)[\star g_Y]$.
\end{lem}
\pf
We have a natural isomorphism
$F^i_{\dagger}(\iota_{g_X\dagger}\nbigm[\star g_X])
\simeq 
 \iota_{g_Y\dagger}\bigl(
 F^i_{\dagger}(\nbigm[\star g_X])
 \bigr)$.
We also have
$F^i_{\dagger}(\iota_{g_X\dagger}\nbigm[\star g_X])
\simeq
 F^i_{\dagger}\bigl(
 (\iota_{g_X\dagger}\nbigm)[\star t]
 \bigr)
\simeq
 F^i_{\dagger}(\iota_{g_X\dagger}\nbigm)[\star t]$.
Then, the claim follows.
\hfill\qed

\vspace{.1in}

Let $\nbigt$ be a coherent $\nbigr_{X}(\ast g_X)$-triple
which is strictly specializable along $g_X$.
We suppose the following.
\begin{itemize}
\item
 The support of $\nbigm$ is proper over $Y$.
\item
(ii) $F^i_{\dagger}\psi_{g_X,u}\nbigt$ 
 are strict for any $i\in\seisuu$ 
and $u\in\real\times\cnum$.
\end{itemize}
\begin{cor}
Assume $\nbigt$ is localizable along $g_X$.
Then,
$F^i_{\dagger}(\nbigt)$ is localizable along $g_Y$,
and we have natural isomorphisms
$F^i_{\dagger}(\nbigt)[\star g_Y]
\simeq
 F^i_{\dagger}\bigl(\nbigt[\star g_X]\bigr)$.
\hfill\qed
\end{cor}

\subsubsection{$\nbigr$-modules and $\nbigr$-triples}

\begin{df}
Let $\nbigm$ (resp. $\nbigt$)
be a coherent $\nbigr_{X}$-module
(resp. $\nbigr_X$-triple)
which is strictly specializable along $g$.
It is called localizable along $g$,
if and only if 
$\nbigm(\ast g)$ 
(resp. $\nbigt(\ast g)$)
is localizable along $g$.
\hfill\qed
\end{df}
In that case,
we set
$\nbigm[\star g]:=
 \bigl(
 \nbigm(\ast g)
 \bigr)[\star g]$
and 
$\nbigt[\star g]:=
 \bigl(
 \nbigt(\ast g)
 \bigr)[\star g]$.

\subsection{Canonical prolongations across hypersurfaces}

Let $D$ be an effective divisor of $X$.
\begin{df}
A coherent $\nbigr_{X(\ast D)}$-module
(resp. $\nbigr_{X(\ast D)}$-triple)
$\nbigm$ (resp. $\nbigt$)
is called strictly specializable along $D$,
if the following holds:
\begin{itemize}
\item 
Let $U$ be any open subset of $X$
with a generator $g_U$ of $\nbigo(-D)_{|U}$.
Then,
 $\nbigm_{|U}$ (resp. $\nbigt_{|U}$)
is strictly specializable along $g_U$.
\end{itemize}
We say that  $\nbigm$ (resp. $\nbigt$)
is localizable along $D$,
if the following holds:
\begin{itemize}
\item 
Let $U$ and $g_U$ be as above.
Then,
 $\nbigm_{|U}$ (resp. $\nbigt_{|U}$)
is localizable along $g_U$.
\hfill\qed
\end{itemize}
\end{df}

If $\nbigm$ is localizable along $D$,
let $\nbigm[\star D]$
denote an $\nbigr_X$-module
with an isomorphism
$\bigl( \nbigm[\star D]\bigr)(\ast D)
\simeq \nbigm$
determined by the following condition:
\begin{itemize}
\item
 Let $U$ and $g_U$ be as above.
 Then,
 $\nbigm[\star D]_{|U}
\simeq
 \nbigm_{|U}[\star g_U]$.
\end{itemize}
We use the symbol
$\nbigt[\star D]$ in similar meanings.

\subsubsection{Morphisms}

Let $\nbigt_i$ be $\nbigr_X(\ast D)$-modules
which are localizable along $D$.
We obtain the following lemma from
Lemma \ref{lem;13.5.10.10}.
\begin{lem}
Let $\star=\ast,!$.
Morphisms $\nbigt_1\lrarr\nbigt_2$
of $\nbigr_{X}(\ast D)$-triples
bijectively correspond to 
morphisms
$\nbigt_1[\star D]\lrarr\nbigt_2[\star D]$
of $\nbigr_{X}$-triples.
If $\nbigt_i$ are integrable,
we have such bijections
for integrable morphisms.
Similar claims hold for 
$\nbigr_X(\ast D)$-modules.
\hfill\qed
\end{lem}
\index{strictly specializable}

Let $f:\nbigt_1\lrarr\nbigt_2$ be a morphism.
We have the induced morphism
$f[\star D]:\nbigt_1[\star D]\lrarr\nbigt_2[\star D]$.
We obtain the following lemma from
Lemma \ref{lem;13.5.10.11}.
\begin{lem}
Suppose that $f$ is strictly specializable along $D$,
i.e.,  for any open subset $U\subset X$ with a generator
$g_U$ of $\nbigo(-D)_{|U}$,
$f_{|U}$ is strictly specializable along $g_U$.
Then, 
$\Ker f$,
$\Image f$
and $\Cok f$ are strictly specializable and localizable
along $D$,
and we naturally have
$\Ker (f[\star D])\simeq
 \Ker (f)[\star D]$,
$\Image(f[\star D])\simeq
 \Image (f)[\star D]$,
and 
$\Cok(f[\star D])\simeq
 \Cok(f)[\star D]$.
\hfill\qed
\end{lem}

Let $\nbigt^{\bullet}$ be a bounded complex of
$\nbigr_{X}(\ast D)$-triples
such that each $\nbigt^{p}$
is localizable along $D$.
We obtain the following from Lemma \ref{lem;13.5.10.13}.
\begin{lem}
Assume that the complex $\nbigt^{\bullet}$ is strictly specializable
in the following sense:
\begin{itemize}
\item
 For any open $U\subset X$ 
 with a generator $g_U$ of $\nbigo(-D)_{|U}$,
 $\nbigh^{\bullet}(\psi_{g_U,u}\nbigt^{\bullet})$
are strict for any $u\in\real\times\cnum$. 
\end{itemize}
Then, the following holds:
\begin{itemize}
\item
$\nbigh^{\bullet}(\nbigt^{\bullet})$
are strictly specializable and localizable along $D$.
\item
We have natural isomorphisms
$\nbigh^{\bullet}\bigl(
 \nbigt^{\bullet}[\star D]\bigr)
\simeq
 \nbigh^{\bullet}
 \bigl(\nbigt^{\bullet}\bigr)[\star D]$.
\hfill\qed
\end{itemize}
\end{lem}

\subsubsection{Compatibility with push-forward}

Let $F:X\lrarr Y$ be a morphism of complex manifolds.
Let $D_Y$ be an effective divisor of $Y$,
and we set $D_X:=F^{-1}(D_Y)$.
Let $\nbigt$ be a coherent $\nbigr_X$-triple
which is strictly specializable along $D_X$.
Assume the following:
\begin{itemize}
\item
 The support of $\nbigt$ is proper over $Y$.
\item
 Let $U_Y\subset Y$ be open with a generator $g_{U_Y}$ 
 of $\nbigo(-D_Y)_{|U}$.
 We set
 $U_X:=F^{-1}(U_Y)$
 and $g_{U_X}:=g_{U_Y}\circ F_{|U_X}$.
 We suppose that
 $F_{\dagger}^i\psitilde_{g_{U_X},u}(\nbigt_{|U_X})$
 are strict for any $u\in\real\times\cnum$.
\end{itemize}

\begin{lem}
If moreover $\nbigt$ is localizable along $D_X$,
then $(F^i_{\dagger}\nbigm)$ is localizable along $D_Y$,
and we have a natural isomorphism
$F^i_{\dagger}(\nbigm[\star D_X])
\simeq (F^i_{\dagger}\nbigm)[\star D_Y]$.
A similar claim holds for $\nbigr_X$-modules.
\hfill\qed
\end{lem}

\subsubsection{$\nbigr$-modules and $\nbigr$-triples}

\begin{df}
Let $\nbigm$ (resp. $\nbigt$)
be a coherent $\nbigr_{X}$-module
(resp. $\nbigr_X$-triple)
which is strictly specializable along $D$.
It is called localizable along $D$,
if and only if 
$\nbigm(\ast D)$ 
(resp. $\nbigt(\ast D)$)
is localizable along $D$.
\hfill\qed
\end{df}
In that case,
we set
$\nbigm[\star D]:=
 \bigl(
 \nbigm(\ast D)
 \bigr)[\star D]$
and 
$\nbigt[\star D]:=
 \bigl(
 \nbigt(\ast D)
 \bigr)[\star D]$.

\chapter{Gluing and specialization of $\nbigr$-triples}
\label{section;11.4.3.1}

We consider how to define specializations
and gluings of $\nbigr$-triples
(\S\ref{subsection;13.4.12.10}--\ref{subsection;13.4.12.11}).
We apply an idea of Beilinson in \cite{beilinson2}.
In particular, we obtain a way to define nearby cycle functor.
In \S\ref{subsection;11.1.24.10},
we shall compare it with the nearby cycle functor
introduced in \cite{sabbah2}.
We introduce a condition of admissible specializability
in \S\ref{subsection;13.4.12.12},
which will be used in the definition of
mixed twistor $D$-modules.

\section{Beilinson functors for $\nbigr$-modules}
\label{subsection;13.4.12.10}

\index{Beilinson functors}

\subsection{The functors $\Pi^{a,b}$, 
$\Pi^{a,b}_{!}$, $\Pi^{a,b}_{\ast}$
and $\Pi^{a,b}_{\ast!}$ for $\nbigr$-module}

\index{functor $\Pi^{a,b}$}
\index{functor $\Pi^{a,b}_{\ast}$}
\index{functor $\Pi^{a,b}_{\bikkuri}$}
\index{functor $\Pi^{a,b}_{\ast\bikkuri}$}

Let $X$ and $H$ be as in 
\S\ref{subsection;11.2.15.1}.
Let $a\leq b$ be integers.
For an $\nbigr_{X(\ast H)}(\ast t)$-module
$\nbigm$,
we set
$\Pi^{a,b}\nbigm:=
 \nbigm\otimes\IItilde^{a,b}_2$.
(See \S\ref{subsection;11.2.1.10}
for $\IItilde^{a,b}_2$.)
We have naturally defined morphisms
$\Pi^{a,b}\nbigm\lrarr\Pi^{c,d}\nbigm$
for $a>c$ and $b>d$.
If $\nbigm$ is coherent and strictly specializable along $t$,
so is $\Pi^{a,b}\nbigm$.
In that case,
we set
$\Pi^{a,b}_{\star}\nbigm:=
 \Pi^{a,b}\nbigm[\star t]$.
We define
\[
 \Pi^{a,b}_{\ast!}\nbigm:=
\varprojlim_{N\to\infty}
 \Cok\Bigl(
 \Pi^{b,N}_!\nbigm
\lrarr
 \Pi^{a,N}_{\ast}\nbigm
 \Bigr)
\]

\begin{lem}
\label{lem;11.1.16.102}
Let $P\in\nbigx$.
There exists $N(P)>0$ such that,
for any $N>N(P)$,
on a neighbourhood $\nbigx_P$ of $P$,
the morphism
$\Pi^{b,N}_!\nbigm
\lrarr
 \Pi^{a,N}_{\ast}\nbigm$ is 
strictly specializable along $t$,
and its cokernel is independent of $N$
in the sense that
the naturally defined morphism
\begin{equation}
 \label{eq;11.1.16.22}
 \Cok\bigl(
 \Pi_!^{b,N+1}\nbigm\lrarr\Pi_{\ast}^{a,N+1}\nbigm
 \bigr)
\lrarr
 \Cok\bigl(
 \Pi_!^{b,N}\nbigm\lrarr\Pi_{\ast}^{a,N}\nbigm
 \bigr)
\end{equation}
is an isomorphism
on $\nbigx_P$.
\end{lem}
\pf
In this proof, $N$ is a sufficiently large number,
and we omit to distinguish a neighbourhood $\nbigx_P$.
If $u\not\in\seisuu_{\geq 0}\times \{0\}$,
we have isomorphisms
$\psi_{u}\bigl(
 \Pi^{a,N}_{\star}\nbigm
 \bigr)\simeq
 \psi_{u}(\nbigm)
 \otimes\II_2^{a,N}$
and the following commutative diagram:
\begin{equation}
 \label{eq;11.1.16.21}
 \begin{CD}
\psi_u\bigl(\Pi_!^{b,N}\nbigm\bigr)
 @>>>
\psi_u\bigl(\Pi_{\ast}^{a,N}\nbigm\bigr)
 \\
@V{\simeq}VV @V{\simeq}VV \\
\psi_{u}(\nbigm)\otimes\II_2^{b,N}
 @>{\id\otimes\iota}>>
\psi_{u}(\nbigm)\otimes\II_2^{a,N}
 \end{CD}
\end{equation}
Here $\iota$ is a natural inclusion
$\II_2^{b,N}\lrarr\II_{2}^{a,N}$.
Hence, 
$\psi_u\bigl(\Pi_!^{b,N}\nbigm\bigr)
\lrarr
\psi_u\bigl(\Pi_{\ast}^{a,N}\nbigm\bigr)$
is strict in this case.

We have the following isomorphisms:
\[
\begin{CD}
 \psi_{0}\bigl(\Pi_{!}^{a,N}(\nbigm)\bigr)
@<{-\deldel_t}<{\simeq}<
 \psi_{-\vecdelta}\bigl(\Pi_!^{a,N}\nbigm\bigr)
\simeq
 \psi_{-\vecdelta}(\nbigm)
 \otimes\II^{a,N}_2
\end{CD}
\]
\[
 \begin{CD}
  \psi_{0}\bigl(\Pi_{\ast}^{a,N}(\nbigm)\bigr)
@>{t}>{\simeq}>
 \psi_{-\vecdelta}\bigl(\Pi_{\ast}^{a,N}\nbigm\bigr)
\simeq
 \psi_{-\vecdelta}(\nbigm)
 \otimes\II^{a,N}_2
 \end{CD}
\]
We have the following commutative diagram:
\begin{equation}
 \label{eq;11.1.16.20}
 \begin{CD}
\psi_0\bigl(\Pi_!^{b,N}\nbigm\bigr)
 @>>>
\psi_0\bigl(\Pi_{\ast}^{b,N}\nbigm\bigr)
 @>>>
\psi_0\bigl(\Pi_{\ast}^{a,N}\nbigm\bigr)
 \\
@A{\simeq}AA @V{\simeq}VV @V{\simeq}VV \\
\lambda\,\psi_{-\vecdelta}(\nbigm)\otimes\II_2^{b,N}
 @>{\nbign_{\nbigm}\otimes\id+\id\otimes \nbign''_{\II}}>>
\psi_{-\vecdelta}(\nbigm)\otimes\II_2^{b,N}
 @>{\id\otimes\iota}>>
\psi_{-\vecdelta}(\nbigm)\otimes\II_2^{a,N}
 \end{CD}
\end{equation}
Here, $\nbign_{\nbigm}$ is induced by $-t\del_t$.
Then, it is easy to check that
$\psi_0\bigl(\Pi_!^{b,N}\nbigm\bigr)
\lrarr
\psi_0\bigl(\Pi_{\ast}^{a,N}\nbigm\bigr)$
is strict.
The case $u\in\seisuu_{\geq 0}\times\{0\}$ also follows.
Hence, the natural morphism
$\Pi_!^{b,N}\nbigm\lrarr
\Pi_{\ast}^{a,N}\nbigm$ 
is strict.

The following morphism
of $\nbigr_{X(\ast H)}(\ast t)$-modules is an isomorphism:
\[
 \Cok\bigl(
 \Pi^{b,N+1}\nbigm
\lrarr
 \Pi^{a,N+1}\nbigm
 \bigr)
\lrarr
 \Cok\bigl(
 \Pi^{b,N}\nbigm
\lrarr
 \Pi^{a,N}\nbigm
 \bigr)
\]
We have the following identifications:
{\small
\begin{equation}
\label{eq;11.1.15.21}
\begin{CD}
  \psi_u\Cok\Bigl(
\Pi_{!}^{b,N+1}\nbigm
\lrarr
\Pi_{\ast}^{a,N+1}\nbigm
 \Bigr)
@>>>
\psi_u
 \Cok\Bigl(
\Pi_!^{b,N}\nbigm
\lrarr
\Pi_{\ast}^{a,N}\nbigm
 \Bigr) \\
 @V{\simeq}VV @V{\simeq}VV \\
  \Cok\Bigl(
\psi_u\bigl(\Pi_{!}^{b,N+1}\nbigm\bigr)
\lrarr
\psi_u\bigl(\Pi_{\ast}^{a,N+1}\nbigm\bigr)
 \Bigr)
@>>>
 \Cok\Bigl(
 \psi_u\bigl(\Pi_!^{b,N}\nbigm\bigr)
\lrarr
 \psi_u\bigl(\Pi_{\ast}^{a,N}\nbigm\bigr)
 \Bigr)
\end{CD}
\end{equation}
}
By using the identification (\ref{eq;11.1.16.21}),
we obtain that the horizontal arrows in
(\ref{eq;11.1.15.21}) are isomorphisms
in the case $u\not\in\seisuu_{\geq 0}\times\{0\}$.
By using Lemma \ref{lem;10.8.6.2},
with an easy diagram chasing,
we obtain that the horizontal arrows in
(\ref{eq;11.1.15.21}) are isomorphisms 
in the case $u=(0,0)$.
It also follows that they are isomorphisms
in the case $u\in\seisuu_{>0}\times\{0\}$.
Thus, (\ref{eq;11.1.16.22}) is an isomorphism.
\hfill\qed

\vspace{.1in}
We obtain the following lemma
by a similar argument.
\begin{lem}
\label{lem;11.1.17.1}
Let $P\in\nbigx$.
There exists $N(P)>0$ such that,
for any $N>N(P)$, on a neighbourhood $\nbigx_P$ of $P$,
the morphism 
$\Pi^{-N,b}_!\nbigm\lrarr\Pi^{-N,a}_{\ast}\nbigm$
is strictly specializable along $t$,
and its kernel is independent of the choice of $N$
in the sense that the following naturally defined
morphism 
\[
 \Ker\Bigl(
 \Pi^{-N,b}_!\nbigm\lrarr\Pi^{-N,a}_{\ast}\nbigm
 \Bigr)
\lrarr
  \Ker\Bigl(
 \Pi^{-N-1,b}_!\nbigm\lrarr\Pi^{-N-1,a}_{\ast}\nbigm
 \Bigr)
\]
is an isomorphism.
\hfill\qed
\end{lem}

The following lemma is clear by construction.
\begin{lem}
If $\nbigm$ is integrable,
$\Pi_{\star!}^{a,b}\nbigm$ is naturally integrable.
\hfill\qed
\end{lem}

\subsection{Another description}

We have another description of
the functor $\Pi^{a,b}_{\ast!}$.
\begin{lem}
\label{lem;11.1.17.2}
We have the following natural isomorphism:
\begin{equation}
\label{eq;11.1.16.105}
 \Pi^{a,b}_{\ast!}(\nbigm)
\simeq
 \varinjlim_{N\to\infty}
 \Ker\Bigl(
 \Pi^{-N,b}_!\nbigm
\lrarr
 \Pi_{\ast}^{-N,a}\nbigm
 \Bigr)
\end{equation}
\end{lem}
\pf
First, let us construct such an isomorphism
on a neighbourhood $\nbigx_P$ of $P$.
We have the following natural commutative diagram:
\[
 \begin{CD}
 \Pi^{b,N}_!\nbigm
 @>>>
 \Pi^{-N,N}_!\nbigm
 @>>>
 \Pi^{-N,b}_!\nbigm\\
 @VVV @VVV @VVV \\
 \Pi^{a,N}_{\ast}\nbigm
 @>>>
 \Pi^{-N,N}_{\ast}\nbigm
 @>>>
 \Pi^{-N,a}_{\ast}\nbigm
 \end{CD}
\]
Hence, we have only to prove the following morphisms
are isomorphisms on $\nbigx_P$,
if $N$ is sufficiently large:
\begin{equation}
 \label{eq;11.1.16.30}
 \Cok\Bigl(
 \Pi_!^{-N,N}\nbigm
\lrarr
 \Pi_{\ast}^{-N,N}\nbigm
 \Bigr)
\lrarr
 \Cok\Bigl(
 \Pi_!^{-N,b}\nbigm
\lrarr
 \Pi_{\ast}^{-N,a}\nbigm
 \Bigr)
\end{equation}
\begin{equation}
 \label{eq;11.1.16.31}
 \Ker\Bigl(
 \Pi_!^{b,N}\nbigm
\lrarr
 \Pi_{\ast}^{a,N}\nbigm
 \Bigr)
\lrarr
 \Ker\Bigl(
 \Pi_!^{-N,N}\nbigm
\lrarr
 \Pi_{\ast}^{-N,N}\nbigm
 \Bigr)
\end{equation}
We have
$ \Cok\Bigl(
 \Pi^{-N,N}\nbigm
\lrarr
 \Pi^{-N,N}\nbigm
 \Bigr)
\simeq
 \Cok\Bigl(
 \Pi^{-N,b}\nbigm
\lrarr
 \Pi^{-N,a}\nbigm
 \Bigr)\simeq 0$.
Let us consider the following morphism:
\begin{multline}
\label{eq;11.1.16.100}
 \Cok\Bigl(
 \psi_u\bigl(
 \Pi_{!}^{-N,N}\nbigm
\bigr)
\lrarr
\psi_u\bigl(
 \Pi_{\ast}^{-N,N}\nbigm
\bigr)
 \Bigr)
\lrarr 
 \\
 \Cok\Bigl(
\psi_u\bigl(
 \Pi_!^{-N,b}\nbigm
\bigr)
\lrarr
\psi_u\bigl(
 \Pi_{\ast}^{-N,a}\nbigm
\bigr)
 \Bigr)
\end{multline}
If $u\not\in\seisuu_{\geq 0}\times\{0\}$,
the both sides are $0$.
For $u=(0,0)$,
we obtain that (\ref{eq;11.1.16.100})
is an isomorphism
from  Lemma \ref{lem;11.1.16.101}.
Hence, (\ref{eq;11.1.16.100})
is an isomorphism for each $u$.
Then, we can deduce (\ref{eq;11.1.16.30})
is an isomorphism
by using an argument in the proof of 
Lemma \ref{lem;11.1.16.102}.
We obtain that (\ref{eq;11.1.16.31})
is an isomorphism by a similar argument.
Thus, we obtain the isomorphism (\ref{eq;11.1.16.105})
on a neighbourhood of $\nbigx_P$.
By varying $P$,
we obtain the isomorphism globally.
\hfill\qed

\subsection{The induced morphism}

Let $\nbigm_i$ $(i=1,2)$ be
coherent $\nbigr_{X(\ast H)}(\ast t)$-modules
which are strictly specializable along $t$.
For a morphism $f:\nbigm_1\lrarr\nbigm_2$,
we have the induced morphisms
$\Pi^{a,b}_{\star}(f):
 \Pi^{a,b}_{\star}(\nbigm_1)\lrarr
 \Pi^{a,b}_{\star}(\nbigm_2)$
and
$\Pi^{a,b}_{\ast!}(f):
 \Pi^{a,b}_{\ast!}(\nbigm_1)\lrarr
 \Pi^{a,b}_{\ast!}(\nbigm_2)$.

\begin{lem}
\label{lem;11.1.17.5}
If $f$ is strictly specializable,
we have the following natural isomorphisms:
\[
 \Ker\bigl(\Pi^{a,b}_{\ast!}f\bigr)
\simeq
 \Pi^{a,b}_{\ast!}\Ker f,
\quad
 \Image\bigl(\Pi^{a,b}_{\ast!}f\bigr)
\simeq
 \Pi^{a,b}_{\ast!}\Image f,
\quad
 \Cok\bigl(\Pi^{a,b}_{\ast!}f\bigr)
\simeq
 \Pi^{a,b}_{\ast!}\Cok f.
\]
\end{lem}
\pf
We obtain the claim for $\Ker$ from
the following commutative diagram:
\[
 \begin{CD}
 \Pi^{a,b}_{\ast!}\nbigm_1
 @>>>
 \Pi_!^{-N,b}\nbigm_1
 @>>>
 \Pi_{\ast}^{-N,a}\nbigm_1\\
 @VVV @VVV @VVV \\
 \Pi^{a,b}_{\ast!}\nbigm_2
 @>>>
 \Pi_!^{-N,b}\nbigm_2
 @>>>
 \Pi_{\ast}^{-N,a}\nbigm_2
 \end{CD}
\]
We obtain the claim for $\Cok$
from the following commutative diagram:
\[
 \begin{CD}
 \Pi_!^{b,N}\nbigm_1
 @>>>
 \Pi_{\ast}^{a,N}\nbigm_1
 @>>>
 \Pi^{a,b}_{\ast!}\nbigm_1 \\
 @VVV @VVV @VVV \\
 \Pi_!^{b,N}\nbigm_2
 @>>>
 \Pi_{\ast}^{a,N}\nbigm_2
 @>>>
 \Pi^{a,b}_{\ast!}\nbigm_2
 \end{CD}
\]
The claim for the image follows from
the claims for the kernel and the cokernel
with an easy diagram chase.
\hfill\qed

\begin{cor}
\label{cor;11.1.20.1}
Let $\nbigm^{\bullet}$ be as in 
Corollary {\rm\ref{cor;11.1.15.30}}.
 We have natural isomorphisms
 $\Pi^{a,b}_{\ast!}\nbigh^{\bullet}\nbigm^{\bullet}
\simeq 
 \nbigh^{\bullet}
 \Pi^{a,b}_{\ast!}\nbigm^{\bullet}$.
\hfill\qed
\end{cor}

\subsection{Compatibility with the push-forward}
\label{subsection;11.1.20.2}

Let us consider the compatibility of the functors
$\Pi^{a.b}_{\star}$ and $\Pi^{a,b}_{\ast!}$
with the push-forward.
We use the notation in 
\S\ref{subsection;11.1.17.20}.
Let $\nbigm$ be a good coherent $\nbigr_X(\ast t)$-module
which is strictly specializable along $t$.
Assume that the support of $\nbigm$ is
proper over $Y$ with respect to $F$.
We have a natural isomorphism
$F^i_{\dagger}\bigl(
 \Pi^{a,b}\nbigm\bigr)
\simeq
 \Pi^{a,b}F^i_{\dagger}(\nbigm)$
of coherent $\nbigr_{Y(\ast H_Y)}$-modules.
Assume that
$F_{\dagger}^i\psitilde_u\nbigm$
is strict
for any $u\in\real\times\cnum$
and $i\in\seisuu$.
Then, $F^i_{\dagger}\bigl(
 \Pi^{a,b}\nbigm\bigr)$
are strictly specializable along $t$,
and 
we have
$F_{\dagger}^i\bigl(
 \Pi_{\star}^{a,b}\nbigm\bigr)
\simeq
 \Pi_{\star}^{a,b}F^{i}_{\dagger}\nbigm$
by Lemma \ref{lem;11.1.17.21}.
We have the following naturally defined morphism:
\begin{equation}
\label{eq;11.1.17.31}
 \Pi^{a,b}_{\ast !}F_{\dagger}^i\nbigm
\lrarr
 F_{\dagger}^i\Pi^{a,b}_{\ast !}\nbigm.
\end{equation}

\begin{prop}
\label{prop;11.1.17.22}
The morphism {\rm(\ref{eq;11.1.17.31})} is 
an isomorphism.
\end{prop}
\pf
We consider locally around $P\in \nbigy$.
Let $\nbigk_N$ be the kernel of
$\Pi_!^{b,N}\nbigm\lrarr
 \Pi_{\ast}^{a,N}\nbigm$.
If $N$ is sufficiently large,
we have the following exact sequence:
\[
 0\lrarr \psi_u\nbigk_N\lrarr
 \psi_u(\Pi_!^{b,N}\nbigm)\lrarr
 \psi_u(\Pi_{\ast}^{a,N}\nbigm)\lrarr
 \psi_u\Pi^{a,b}_{\ast !}\nbigm
 \lrarr 0
\]
It can be rewritten as follows:
\begin{equation}
 \label{eq;11.1.17.20}
 0\lrarr
 \psi_u\nbigk_N
\lrarr
 \psi_u\nbigm\otimes \II^{b,N}
\lrarr
 \psi_u\nbigm\otimes \II^{a,N}
\lrarr
 \psi_u\Pi^{a,b}_{\ast !}\nbigm
\lrarr 0
\end{equation}
Note that, as an $\nbigr_{X_0(\ast H_0)}$-complex,
the exact sequence (\ref{eq;11.1.17.20})
has a splitting,
and as $\nbigr_{X_0(\ast H_0)}$-modules,
$\psi_u\nbigk_N$
and $\psi_u\Pi^{a,b}_{\ast !}\nbigm$
are direct sums of
some copies of $\psi_u\nbigm$.
Hence,
$F^i_{\dagger}\psi_u\Pi^{a,b}_{\ast !}\nbigm$
and $F^i_{\dagger}\psi_u\nbigk_N$
are strict,
and the induced sequence
\[
 0\lrarr
 F_{\dagger}^i\psi_u\nbigk_N\lrarr
 F_{\dagger}^i\psi_u\nbigm\otimes \II^{b,N}\lrarr
 F_{\dagger}^i\psi_u\nbigm\otimes \II^{a,N}\lrarr
 F_{\dagger}^i\psi_u\Pi^{a,b}_{\ast !}\nbigm
 \lrarr 0
\]
is exact.
The following complex
of $\nbigr_Y(\ast H_Y)$-modules
is clearly exact:
\[
 0\lrarr
 F_{\dagger}^i\nbigk_N(\ast t)\lrarr
 F_{\dagger}^i\nbigm\otimes
 \IItilde^{b,N}
 \lrarr
 F_{\dagger}^i\nbigm\otimes
 \IItilde^{a,N}
 \lrarr
 F_{\dagger}^i\Pi^{a,b}_{\ast !}\nbigm(\ast t)
 \lrarr 0
\]
We obtain that
$0\lrarr
 F_{\dagger}^i\nbigk_N\lrarr
 F_{\dagger}^i
 \Pi_!^{b,N}\nbigm
 \lrarr
 F_{\dagger}^i
 \Pi_{\ast}^{a,N}\nbigm
 \lrarr
 F_{\dagger}^i\Pi^{a,b}_{\ast !}\nbigm
 \lrarr 0$
is exact.
In particular,
we obtain that
(\ref{eq;11.1.17.31}) is an isomorphism.
\hfill\qed

\subsection{The functors $\psi^{(a)}$ and $\Xi^{(a)}$}

\index{functor $\psi^{(a)}$}
\index{functor $\Xi^{(a)}$}

Let $\nbigm$ be a coherent 
$\nbigr_{X(\ast H)}$-module
which is strictly specializable along $t$.
We define
\[
 \psi^{(a)}(\nbigm):=
 \Pi^{a,a}_{\ast!}(\nbigm),
\quad
 \Xi^{(a)}(\nbigm):=
 \Pi^{a,a+1}_{\ast!}(\nbigm).
\]
We have naturally defined exact sequences:
\[
 \begin{CD}
 0@>>>
 \nbigm[!t]\otimes \,s^a
 @>{\alpha_a}>>
 \Xi^{(a)}(\nbigm)
 @>{\beta_a}>>
 \psi^{(a)}(\nbigm)
 @>>> 0
 \end{CD}
\]
\[
 \begin{CD}
 0@>>>
 \psi^{(a+1)}(\nbigm)
 @>{\gamma_a}>>
 \Xi^{(a)}(\nbigm)
 @>{\delta_a}>>
 \nbigm[\ast t]\otimes\,s^a
 @>>> 0
 \end{CD}
\]
According to Corollary \ref{cor;11.1.20.1},
we have natural isomorphisms
\[
\psi^{(a)}\nbigh^{\bullet}(\nbigm^{\bullet})
\simeq
 \nbigh^{\bullet}\psi^{(a)}(\nbigm^{\bullet}),
\quad\quad
\Xi^{(a)}\nbigh^{\bullet}(\nbigm^{\bullet})
\simeq
 \nbigh^{\bullet}\Xi^{(a)}(\nbigm^{\bullet})
\]
for a complex $\nbigm^{\bullet}$ 
as in Corollary \ref{cor;11.1.15.30}.
In the situation of \S\ref{subsection;11.1.20.2},
we have natural isomorphisms
\[
 F^i_{\dagger}\psi^{(a)}(\nbigm)
\simeq
 \psi^{(a)}F^i_{\dagger}(\nbigm),
\quad\quad
 F^i_{\dagger}\Xi^{(a)}(\nbigm)
\simeq
 \Xi^{(a)}F^i_{\dagger}(\nbigm).
\]

\subsection{Beilinson functors for $\nbigr$-modules}

\index{Beilinson functors}
\index{functor $\Pi^{a,b}$}
\index{functor $\Pi^{a,b}_{\ast}$}
\index{functor $\Pi^{a,b}_{\bikkuri}$}
\index{functor $\Pi^{a,b}_{\ast\bikkuri}$}
\index{functor $\psi^{(a)}$}
\index{functor $\Xi^{(a)}$}

Let $\nbigm$ be a coherent $\nbigr_{X(\ast H)}$-module
which is strictly specializable along $t$.
We obtain a coherent $\nbigr_{X(\ast H)}(\ast t)$-module
$\nbigmtilde:=\nbigm(\ast t)$,
which is strictly specializable along $t$.
We define
$\Pi^{a,b}(\nbigm):=
\Pi^{a,b}(\nbigmtilde)$,
$\Pi^{a,b}_{\star}(\nbigm):=
 \Pi^{a,b}_{\star}(\nbigmtilde)$
and
$\Pi^{a,b}_{\ast!}(\nbigm):=
 \Pi^{a,b}_{\ast!}(\nbigmtilde)$.
We also define
$\psi^{(a)}(\nbigm):=
 \psi^{(a)}(\nbigmtilde)$
and 
$\Xi^{(a)}(\nbigm):=
 \Xi^{(a)}(\nbigmtilde)$.

Recall we have the morphisms
$\nbigm[!t]\stackrel{\iota_1}{\lrarr}\nbigm
 \stackrel{\iota_2}{\lrarr}\nbigm[\ast t]$.
Let $\phi^{(a)}(\nbigm)$ be defined
as the cohomology of the following complex:
\[
\begin{CD}
 \nbigm[!t]\otimes s^a
 @>{\alpha_a+\iota_1}>>
 \Xi^{(a)}(\nbigm)
\oplus
\bigl(\nbigm\otimes s^a\bigr)
 @>{\delta_1-\iota_2}>>
 \nbigm[\ast t]\otimes s^a
\end{CD}
\]
\begin{lem}
We have the following natural isomorphisms:
\[
\psi_0\bigl(
 \psi^{(a)}(\nbigm)
 \bigr)
\simeq
\psi_{-\vecdelta}(\nbigm)\,s^a,
\]
\[
 \psi_0\bigl(
 \Xi^{(a)}(\nbigm)\bigr)
\simeq
 \psi_{-\vecdelta}(\nbigm)\,s^a
\oplus
 \psi_{-\vecdelta}(\nbigm)\,s^{a+1},
\]
\[
 \psi_0\bigl(
 \phi^{(a)}(\nbigm)
 \bigr)
\simeq
 \psi_{0}(\nbigm)\,s^a
\]
\end{lem}
\pf
We have the following isomorphism:
\begin{multline}
 \psi_0(\Pi^{a,b}_{\ast!}(\nbigm))
\simeq
 \Cok\Bigl(
 \psi_{0}(\Pi_!^{b,N}
\nbigm)
\lrarr
 \psi_0(\Pi_{\ast}^{a,N}\nbigm)
 \Bigr)
 \\
 \simeq
 \Cok\Bigl(
\psi_{-\vecdelta}(\nbigm)
\otimes\II_2^{b,N}
\lrarr
\psi_{-\vecdelta}(\nbigm)
\otimes\II_2^{a,N}
 \Bigr)
\end{multline}
We obtain the first two isomorphisms.
Let us consider $\phi^{(a)}$ in the case $a=0$.
By construction,
$\psi_0\phi^{(0)}(\nbigm)$ is obtained
as the cohomology of the following:
\[
\begin{CD}
 \psi_{-\vecdelta}(\nbigm)
@>>>
\Bigl(
\psi_{-\vecdelta}(\nbigm)\,s
\oplus
 \psi_{-\vecdelta}(\nbigm)
\Bigr)
\oplus
 \psi_0(\nbigm)
@>>>
 \psi_{-\vecdelta}(\nbigm)
\end{CD}
\]
The morphism
$\psi_{-\vecdelta}(\nbigm)\lrarr
 \psi_{-\vecdelta}(\nbigm)\,s
\oplus
 \psi_{-\vecdelta}(\nbigm)$
is given by
$(s,N)$.
The morphism
$\psi_{-\vecdelta}(\nbigm)\,s
\oplus
 \psi_{-\vecdelta}(\nbigm)
\lrarr
\psi_{-\vecdelta}(\nbigm)$
is the projection.
Then, it is easy to check
that the cohomology is naturally isomorphic to
$\psi_0(\nbigm)$.
\hfill\qed

\vspace{.1in}
We can reconstruct $\nbigm$ as the cohomology of
the following,
as in \cite{beilinson2}:
\[
\begin{CD}
 \psi^{(1)}(\nbigm)
 @>>>
 \Xi^{(0)}(\nbigm)\oplus
 \phi^{(0)}(\nbigm)
 @>>>
 \psi^{(0)}(\nbigm)
\end{CD}
\]
\begin{rem}
$\phi^{(0)}$ will often be denoted by $\phi$.
\hfill\qed
\end{rem}

\section{Beilinson functors for $\nbigr$-triples}
\label{subsection;13.4.12.11}

\subsection{Functors $\Pi^{a,b}$, $\Pi^{a,b}_{\ast}$
 and $\Pi^{a,b}_{!}$ for $\nbigr(\ast t)$-triple}

\index{Beilinson functors}
\index{functor $\Pi^{a,b}$}
\index{functor $\Pi^{a,b}_{\ast}$}
\index{functor $\Pi^{a,b}_{\bikkuri}$}

Let $\nbigt=(\nbigm',\nbigm'',C)$
be a strictly specializable
$\nbigr_{X(\ast H)}(\ast t)$-triple.
As in \S\ref{subsection;10.12.27.1},
we obtain the following $\nbigr_{X(\ast H)}(\ast t)$-triple:
\[
 \Pi^{a,b}\nbigt:=
\nbigt\otimes\IItilde^{a,b}=
\bigl(
 \nbigm'\otimes \IItilde^{a,b}_1,\,
 \nbigm''\otimes \IItilde^{a,b}_2,\,
 C\otimes C_{\IItilde}
 \bigr)
\]
Here,
$C\otimes C_{\IItilde}$ is given as follows:
\begin{multline}
\label{eq;10.8.5.10}
 \bigl(C\otimes C_{\IItilde}\bigr)
 \bigl(u\otimes (\lambda s)^i,\,
 \sigma^{\ast}(v\otimes(\lambda s)^j)
 \bigr)
=C(u,\sigma^{\ast}v)\cdot
 C_{\IItilde}\bigl((\lambda s)^i,\,
 \sigma^{\ast}(\lambda s)^j
 \bigr) \\
=C(u,\sigma^{\ast}v)
 \frac{(\log|t|^2)^{-i-j}}{(-i-j)!}
 (-1)^j
 \lambda^{i-j}
 \chi_{i+j\leq 0}
\end{multline}
(See \S\ref{subsection;11.4.12.10}.)
The $\nbigr_{X(\ast H)}(\ast t)$-triple
$\Pi^{a,b}\nbigt$ is strictly specializable
along $t$.
Then, we obtain the following $\nbigr_{X(\ast H)}$-triples:
\[
 \Pi^{a,b}_{\ast}\nbigt:=\Pi^{a,b}\nbigt[\ast t],
\quad\quad
\Pi^{a,b}_!\nbigt:=\Pi^{a,b}\nbigt[!t].
\]
If $\nbigt$ is integrable,
$\Pi^{a,b}\nbigt$
and $\Pi^{a,b}_{\star}\nbigt$ $(\star=\ast,!)$
are also naturally integrable.

\begin{lem}
Let $u\in(\real\times\cnum)\setminus(\seisuu_{\geq 0}\times\{0\})$.
We have
\[
 \psi_{u}\bigl(
 \Pi^{a,b}\nbigt
 \bigr)
=\psi_{u}(\nbigt)
 \otimes \II^{a,b}
\simeq\Bigl(
 \psi_{u}(\nbigm')
 \otimes \II_1^{a,b},\,
 \psi_{u}(\nbigm'')
 \otimes \II_2^{a,b},\,
 \psi_{u}C\otimes 
 C_{\II}
 \Bigr)
\]
\end{lem}
\pf
We have only to check the claim for the pairing.
We have the following equality
for any positive integer $M$:
\begin{multline}
 \underset{s+\lambda^{-1}\eigenmap(\lambda,u)}{\Res}
 \int
 \bigl\langle
 C(u,\sigma^{\ast}v),\, \phi
 \bigr\rangle\,
 (\log|t|^2)^M\,
 |t|^{2s}\,dt\,d\tbar
 \\
=
 \underset{s+\lambda^{-1}\eigenmap(\lambda,u)}{\Res}
 \frac{d^M}{ds^M}
 \int
 \bigl\langle
 C(u,\sigma^{\ast}v),\, \phi
 \bigr\rangle\,
 |t|^{2s}\,dt\,d\tbar
=0
\end{multline}
Then, the claim follows.
\hfill\qed

\vspace{.1in}
We have
$\Pi^{a,a+1}\nbigt\simeq
 \nbigt\otimes\newTate(a)$,
and hence
\[
 \nbigt[\star t]\otimes \newTate(a)
\simeq
 \Cok\Bigl(
 \Pi_{\star}^{a+1,N}\nbigt
\lrarr
 \Pi_{\star}^{a,N}\nbigt
 \Bigr)
\quad
 (N>a)
\]
\[
 \nbigt[\star t]\otimes \newTate(b)
\simeq
 \Ker\Bigl(
 \Pi_{\star}^{N,b+1}\nbigt
\lrarr
 \Pi_{\star}^{N,b}\nbigt
 \Bigr)
\quad (N<b)
\]

\subsubsection{Relation with Hermitian adjoint}
Note that 
$\nbigs^{a,b}$ 
(\S\ref{subsection;13.3.27.10})
induces an isomorphism
\[
\Pi^{a,b}(\nbigt^{\ast})
\simeq
 (\Pi^{-b+1,-a+1}\nbigt)^{\ast}. 
\]
If $b=a+1$,
it is equal to 
the canonical isomorphism
$\nbigt^{\ast}\otimes\newTate(a)
\simeq
 \bigl(\nbigt\otimes\newTate(-a)\bigr)^{\ast}$,
multiplied by $(-1)^a$.
We have the induced isomorphisms
\[
 \Pi_!^{a,b}\bigl(\nbigt^{\ast}\bigr)
\simeq
 \bigl(
 \Pi_{\ast}^{-b+1,-a+1}\nbigt
 \bigr)^{\ast},
\quad\quad
\Pi_{\ast}^{a,b}\bigl(\nbigt^{\ast}\bigr)
\simeq
 \bigl(
 \Pi_!^{-b+1,-a+1}\nbigt
 \bigr)^{\ast}.
\]
For $a\geq a'$ and $b\geq b'$,
the following natural diagrams are commutative:
{\small
\[
 \begin{CD}
 \Pi_!^{a,b}(\nbigt^{\ast})
 @>{\simeq}>>
 \bigl(
 \Pi_{\ast}^{-b+1,-a+1}\nbigt
 \bigr)^{\ast} \\
 @VVV @VVV \\
 \Pi_!^{a',b'}(\nbigt^{\ast})
 @>{\simeq}>>
 \bigl(
 \Pi_{\ast}^{-b'+1,-a'+1}\nbigt
 \bigr)^{\ast} 
 \end{CD}
\quad\quad
 \begin{CD}
 \Pi_{\ast}^{a,b}(\nbigt^{\ast})
 @>{\simeq}>>
 \bigl(
 \Pi_{!}^{-b+1,-a+1}\nbigt
 \bigr)^{\ast} \\
 @VVV @VVV \\
 \Pi_{\ast}^{a',b'}(\nbigt^{\ast})
 @>{\simeq}>>
 \bigl(
 \Pi_{!}^{-b'+1,-a'+1}\nbigt
 \bigr)^{\ast}
 \end{CD}
\]
}
If $b=a+1$,
they are equal to
the canonical isomorphisms,
multiplied by $(-1)^a$.
\[
 \nbigt^{\ast}[!t]\otimes\newTate(a)
\simeq
 (\nbigt[\ast t])^{\ast}\otimes\newTate(a)
\simeq
 \bigl(
 \nbigt[\ast t]\otimes\newTate(-a)
 \bigr)^{\ast}
\]
\[
 \nbigt^{\ast}[\ast t]\otimes\newTate(a)
\simeq
 (\nbigt[!t])^{\ast}\otimes\newTate(a)
\simeq
 \bigl(
 \nbigt[!t]\otimes\newTate(-a)
 \bigr)^{\ast}
\]

Let $\nbigs:\nbigt\lrarr
 \nbigt^{\ast}\otimes\newTate(-w)$
be a Hermitian sesqui-linear duality.
We obtain the induced morphism:
\[
 \Pi^{a,b}\nbigs:
 \Pi^{a,b}\nbigt\lrarr
 \Pi^{a,b}\bigl(
 \nbigt^{\ast}\otimes\newTate(-w)
 \bigr)
\simeq
 \bigl(\Pi^{-b+1,-a+1}\nbigt\bigr)^{\ast}
 \otimes\newTate(-w)
\]
For $u\in(\real\times\cnum)\setminus(\seisuu_{\geq 0}\times\{0\})$,
we have the following commutativity:
\[
 \begin{CD}
 \psi_{u}(\Pi^{a,b}\nbigt)
 @>{\psi_{u}(\Pi^{a,b}\nbigs)}>>
 \psi_{u}\Bigl(
 (\Pi^{-b+1,-a+1}\nbigt)^{\ast}\otimes
 \newTate(-w)\Bigr)\\
 @VVV @VVV \\
 \psi_{u}(\nbigt)
 \otimes\II^{a,b}
 @>{\psi_{u}(\nbigs)\otimes
 \nbigs^{a,b}}>>
 \psi_{u}(\nbigt^{\ast})
 \otimes
 (\II^{-b+1,-a+1})^{\ast}
 \otimes\newTate(-w)
 \end{CD}
\]

\subsection{Functors $\Pi^{a,b}_{\ast !}$, $\psi^{(a)}$ and
   $\Xi^{(a)}$}

\index{functor $\Pi^{a,b}_{\ast\bikkuri}$}
\index{functor $\psi^{(a)}$}
\index{functor $\Xi^{(a)}$}

Let $\nbigt$ be any $\nbigr_{X(\ast H)}(\ast t)$-triple
which is strictly specializable along $t$.

\begin{lem}
Let $P\in\nbigx$.
There exists a large number $N(P)>0$ such that,
for any $N>N(P)$, 
on a small neighbourhood $\nbigx_P$ of $P$,
the following naturally defined morphisms
\[
 \Cok\bigl(
 \Pi_!^{b,N+1}\nbigt\lrarr
 \Pi_{\ast}^{a,N+1}\nbigt
 \bigr)
\lrarr
 \Cok\bigl(
 \Pi_!^{b,N}\nbigt\lrarr\Pi_{\ast}^{a,N}\nbigt
 \bigr)
\]
\[
 \Ker\bigl(
 \Pi_!^{-N,b}\nbigt\lrarr
 \Pi_{\ast}^{-N,a}\nbigt
 \bigr)
\lrarr
 \Cok\bigl(
 \Pi_!^{-N-1,b}\nbigt\lrarr
 \Pi_{\ast}^{-N-1,a}\nbigt
 \bigr)
\]
are isomorphisms.
In this sense,
$\Cok\bigl(
 \Pi_!^{b,N}\nbigt\lrarr
 \Pi_{\ast}^{a,N}\nbigt
 \bigr)$ 
and
$\Ker\bigl(
 \Pi_!^{-N,b}\nbigt\lrarr
 \Pi_{\ast}^{-N,a}\nbigt
\bigr)$
are independent of $N>N(P)$.
Moreover, they are naturally isomorphic.
\end{lem}
\pf
We obtain the first claim from 
Lemma \ref{lem;11.1.16.102}
and Lemma \ref{lem;11.1.17.1}.
We have the following natural commutative diagram:
\[
 \begin{CD}
 \Pi^{b,N}_!\nbigt
 @>>>
 \Pi^{-N,N}_!\nbigt
 @>>>
 \Pi^{-N,b}_!\nbigt\\
 @VVV @VVV @VVV \\
 \Pi^{a,N}_{\ast}\nbigt
 @>>>
 \Pi^{-N,N}_{\ast}\nbigt
 @>>>
 \Pi^{-N,a}_!\nbigt
 \end{CD}
\]
Then, by using Lemma \ref{lem;11.1.17.2},
we obtain the second claim.
\hfill\qed

\vspace{.1in}

Then, we define
\[
  \Pi^{a,b}_{\ast !}\nbigt:=
 \varprojlim_{N\to\infty}
 \Cok\bigl(
 \Pi_!^{b,N}\nbigt\lrarr
 \Pi_{\ast}^{a,N}\nbigt
 \bigr)
\simeq
 \varinjlim_{N\to\infty}
 \Ker\bigl(
 \Pi_!^{-N,b}\nbigt\lrarr
 \Pi_{\ast}^{-N,a}\nbigt
 \bigr)
\]
If $\nbigt$ is integrable,
$\Pi^{a,b}_{\ast!}\nbigt$
are also naturally integrable.
Let $\nbigt_i$ $(i=1,2)$ be
coherent $\nbigr_{X(\ast H)}(\ast t)$-triples,
which are strictly specializable along $t$.
For any morphism $f:\nbigt_1\lrarr\nbigt_2$,
we have the induced morphisms
$\Pi^{a,b}_{\star}(f):
 \Pi^{a,b}_{\star}(\nbigt_1)\lrarr
 \Pi^{a,b}_{\star}(\nbigt_2)$
and
$\Pi^{a,b}_{\ast!}(f):
 \Pi^{a,b}_{\ast!}(\nbigt_1)\lrarr
 \Pi^{a,b}_{\ast!}(\nbigt_2)$.
We obtain the following lemma from
Lemma \ref{lem;11.1.17.5}.
\begin{lem}
If $f$ is strictly specializable,
we have the following natural isomorphisms:
\[
 \Ker\bigl(\Pi^{a,b}_{\ast!}f\bigr)
\simeq
 \Pi^{a,b}_{\ast!}\Ker f,
\quad
 \Image\bigl(\Pi^{a,b}_{\ast!}f\bigr)
\simeq
 \Pi^{a,b}_{\ast!}\Image f,
\quad
 \Cok\bigl(\Pi^{a,b}_{\ast!}f\bigr)
\simeq
 \Pi^{a,b}_{\ast!}\Cok f.
\]
\hfill\qed
\end{lem}

\begin{cor}
Let $\nbigt^{\bullet}$ be as in 
Lemma {\rm\ref{lem;11.1.17.12}}.
 We have natural isomorphisms
 $\Pi^{a,b}_{\ast!}\nbigh^{\bullet}\nbigm^{\bullet}
\simeq 
 \nbigh^{\bullet}
 \Pi^{a,b}_{\ast!}\nbigm^{\bullet}$.
\hfill\qed
\end{cor}

\subsubsection{Functors $\psi^{(a)}$ and  $\Xi^{(a)}$}

In particular, for any $a\in\seisuu$,
we define
\[
 \psi^{(a)}\nbigt:=
 \Pi^{a,a}_{\ast !}\nbigt,
\quad
 \Xi^{(a)}\nbigt:=
 \Pi^{a,a+1}_{\ast !}\nbigt.
\]
If $\nbigt$ is integrable,
they have naturally induced integrable structure.
We naturally have the following exact sequences:
\[
0\lrarr\nbigt[!t]\otimes\newTate(a)
 \stackrel{\alpha_a}{\lrarr}
 \Xi^{(a)}\nbigt
 \stackrel{\beta_a}{\lrarr}
 \psi^{(a)}\nbigt\lrarr 0\]
\[
 0\lrarr \psi^{(a+1)}\nbigt
 \stackrel{\gamma_a}{\lrarr}
 \Xi^{(a)}\nbigt
 \stackrel{\delta_a}{\lrarr}
 \nbigt[\ast t]\otimes\newTate(a)\lrarr 0
\]
We also have natural identifications
\[
 \psi^{(a+1)}\nbigt
\simeq
 \psi^{(a)}\nbigt\otimes\newTate(1),
\quad
 \Xi^{(a+1)}\nbigt
\simeq
 \Xi^{(a)}\nbigt\otimes\newTate(1).
\]
The composite
$\beta_a\circ\gamma_a:\psi^{(a+1)}\nbigt
\lrarr
 \psi^{(a)}\nbigt
\simeq
 \psi^{(a+1)}\nbigt\otimes\newTate(-1)$
is induced by the natural morphisms
$\Pi^{a+1,N}\nbigt[\star t]\lrarr
 \Pi^{a,N}\nbigt[\star t]$.

\subsubsection{Relation with Hermitian adjoint}

\begin{lem}
\label{lem;13.3.27.11}
We have a natural isomorphism
$\Pi^{a,b}_{\ast !}(\nbigt^{\ast})
\simeq\bigl(
 \Pi^{-b+1,-a+1}_{\ast !}\nbigt
 \bigr)^{\ast}$.
In particular, we naturally have
\[
 \bigl(\psi^{(a)}\nbigt\bigr)^{\ast}
\simeq
 \psi^{(-a+1)}\bigl(\nbigt^{\ast}\bigr),
\quad\quad
 \bigl(\Xi^{(a)}\nbigt\bigr)^{\ast}
\simeq
 \Xi^{(-a)}\bigl(\nbigt^{\ast}\bigr).
\]
\end{lem}
\pf
We have the following natural isomorphisms:
\begin{multline}
\Pi_{\ast !}^{a,b}(\nbigt^{\ast})
=\Cok\Bigl(
 \Pi_{!}^{b,N}(\nbigt^{\ast})
\lrarr
 \Pi_{\ast}^{a,N}(\nbigt^{\ast})
 \Bigr) \\
\simeq
\Cok\Bigl(
 \bigl(\Pi_{\ast}^{-N+1,-b+1}\nbigt\bigr)^{\ast}
\lrarr
 \bigl(\Pi_{!}^{-N+1,-a+1}\nbigt\bigr)^{\ast}
 \Bigr)
 \\
=\Ker\Bigl(
 \Pi_{\ast}^{-N+1,-b+1}\nbigt
\llarr
\Pi_{!}^{-N+1,-a+1}\nbigt
 \Bigr)^{\ast}
\simeq
 \Pi^{-b+1,-a+1}_{\ast!}(\nbigt)^{\ast}
\end{multline}
Thus, we are done.
\hfill\qed

\vspace{.1in}

For $a\geq a'$ and $b\geq b'$,
the following diagram is commutative:
\[
 \begin{CD}
 \Pi^{a,b}_{\ast!}(\nbigt^{\ast})
 @>>>
 \Pi^{-b+1,-a+1}_{\ast!}(\nbigt)^{\ast}
 \\
 @VVV @VVV \\
 \Pi^{a',b'}_{\ast!}(\nbigt^{\ast})
 @>>>
 \Pi^{-b'+1,-a'+1}_{\ast!}(\nbigt)^{\ast}
 \end{CD}
\]
In particular, the following diagram is commutative:
\[
 \begin{CD}
 \psi^{(a+1)}(\nbigt^{\ast})
 @>{\gamma_{a}(\nbigt^{\ast})}>>
 \Xi^{(a)}(\nbigt^{\ast})
 @>{\beta_{a}(\nbigt^{\ast})}>>
 \psi^{(a)}(\nbigt^{\ast}) 
 \\
 @VVV @VVV @VVV \\
 \psi^{(-a)}(\nbigt)^{\ast}
 @>{\beta_{-a}(\nbigt)^{\ast}}>>
 \Xi^{(-a)}(\nbigt)^{\ast}
 @>{\gamma_{-a}(\nbigt)^{\ast}}>>
 \psi^{(-a+1)}(\nbigt)^{\ast}
 \end{CD}
\]

We have the isomorphisms
$\Upsilon_{a,\nbigt}:
 \psi^{(a+1)}(\nbigt)\simeq
 \psi^{(a)}(\nbigt)\otimes\newTate(1)$
and
$\Upsilon_{-a,\nbigt^{\ast}}:
 \psi^{(-a+1)}(\nbigt^{\ast})\simeq
 \psi^{(-a)}(\nbigt^{\ast})\otimes\newTate(1)$.
The following diagram is commutative:
\[
 \begin{CD}
 \psi^{(a)}(\nbigt)^{\ast}
 @>{\simeq}>{-(\Upsilon_{a,\nbigt})^{\ast}}>
 \psi^{(a+1)}(\nbigt)^{\ast}
 \otimes\newTate(1)\\
 @V{\simeq}V{g_1}V @V{\simeq}V{g_2}V \\
 \psi^{(-a+1)}(\nbigt^{\ast})
 @>{\simeq}>{\Upsilon_{-a,\nbigt^{\ast}}}>
 \psi^{(-a)}(\nbigt^{\ast})
 \otimes\newTate(1)
 \end{CD}
\]
Here, the vertical isomorphisms
are as in Lemma \ref{lem;13.3.27.11}.

\subsection{Vanishing cycle functor for $\nbigr$-triple}

\index{vanishing cycle functor}
\index{functor $\psi^{(a)}$}
\index{functor $\Xi^{(a)}$}
\index{functor $\phi^{(a)}$}

Let $\nbigt$ be an $\nbigr_{X(\ast H)}$-triple
which is strictly specializable along $t$.
Then, applying the above construction to
$\nbigttilde:=\nbigt(\ast t)$,
we define
\[
 \Xi^{(a)}\nbigt:=
 \Xi^{(a)}\nbigttilde,
\quad\quad
 \psi^{(a)}\nbigt:=
 \psi^{(a)}\nbigttilde.
\]
We shall introduce 
the vanishing cycle functor $\phi^{(a)}$.
We have the canonical morphisms
$\nbigt[!t]\stackrel{\iota_1}{\lrarr}
 \nbigt\stackrel{\iota_2}{\lrarr}\nbigt[\ast t]$.
We have
$\iota_{2}\circ\iota_{1}
=\delta_a\circ\alpha_a$,
because the restrictions to $X\setminus\{t=0\}$
are equal.
Then, we define $\phi^{(a)}(\nbigt)$
as the cohomology of the following complex:
\[
\begin{CD}
 \nbigt[!t]\otimes\newTate(a)
@>{\alpha_a+\iota_1}>>
 \Xi^{(a)}\nbigt\oplus
 \bigl(\nbigt\otimes\newTate(a)\bigr)
@>{\delta_a-\iota_2}>>
 \nbigt[\ast t]\otimes\newTate(a)
\end{CD}
\]
We naturally have
$\phi^{(a)}(\nbigt)\simeq
 \phi^{(a+1)}(\nbigt)\otimes
 \newTate(-1)$.
In particular,
we set $\phi(\nbigt):=\phi^{(0)}(\nbigt)$.
If $\nbigt$ is integrable,
$\phi^{(a)}(\nbigt)$ has a naturally induced
integrable structure.

\vspace{.1in}
The morphisms
$\beta_a$ and $\gamma_a$
induce the following morphisms:
\[
\begin{CD}
 \psi^{(a+1)}\nbigt
 @>{\can^{(a)}}>>
 \phi^{(a)}\nbigt
 @>{\var^{(a)}}>>
 \psi^{(a)}\nbigt
\end{CD}
\]
As in \cite{beilinson2},
we can reconstruct
$\nbigt\otimes\newTate(a)$
as the cohomology of
\[
\begin{CD}
 \psi^{(a+1)}\nbigt
 @>{\gamma_{a}+\can^{(a)}}>>
 \Xi^{(a)}\nbigt\oplus\phi^{(a)}\nbigt
 @>{\delta_{a}-\var^{(a)}}>>
 \psi^{(a)}\nbigt
\end{CD}
\]

We have natural isomorphisms
$\psi^{(a)}(\nbigt)^{\ast}
\simeq
 \psi^{(-a+1)}(\nbigt^{\ast})$
and 
$\Xi^{(a)}(\nbigt)^{\ast}
\simeq
 \Xi^{(-a)}(\nbigt^{\ast})$.
We also obtain an induced isomorphism
$\phi^{(a)}(\nbigt)^{\ast}
\simeq
 \phi^{(-a)}(\nbigt^{\ast})$.

\subsection{Gluing of $\nbigr$-triples}
\label{subsection;11.2.22.21}

Let $\nbigt$ be a coherent $\nbigr_{X(\ast H)}(\ast t)$-triple
which is strictly specializable along $t$.
Let $\nbigq$ be a strict coherent $\nbigr_{X(\ast H)}$-triple
with morphisms
\[
 \psi^{(1)}\nbigt\stackrel{u}{\lrarr}
 \nbigq\stackrel{v}{\lrarr}
 \psi^{(0)}\nbigt,
\]
such that (i) $v\circ u=\delta_0\circ\gamma_0$,
(ii) $\Supp\nbigq\subset\{t=0\}$.
Then, we obtain 
an $\nbigr$-triple
$\Glue(\nbigt,\nbigq,u,v)$
as the cohomology of the following complex:
\[
 \begin{CD}
 \psi^{(1)}\nbigt
 @>{\gamma_{0}+u}>>
 \Xi^{(0)}\nbigt\oplus\nbigq
 @>{\delta_{0}-v}>>
 \psi^{(0)}\nbigt
\end{CD}
\]
\index{$\nbigr$-triple $\Glue(\nbigt,\nbigq,u,v)$}
We naturally have 
$\Glue(\nbigt,\nbigq,u,v)(\ast t)\simeq
 \nbigt$
and
$\phi^{(0)}\bigl(
 \Glue(\nbigt,\nbigq,u,v)\bigr)\simeq\nbigq$.
Under the isomorphisms,
we have $\can^{(0)}=u$ and $\var^{(0)}=v$.
We naturally have 
$\Glue(\nbigt,\nbigq,u,v)^{\ast}
\simeq
 \Glue(\nbigt^{\ast},
 \nbigq^{\ast},-v^{\ast},-u^{\ast})$.
If $\nbigt$, $\nbigq$, $u$ and $v$ are integrable,
the object
$\Glue(\nbigt,\nbigq,u,v)$ is also 
naturally integrable.

\vspace{.1in}
We give a remark for Lemma \ref{lem;11.1.18.10} below.
We have the naturally induced morphisms
\[
\Xi^{(0)}(\nbigt)\lrarr
 \Xi^{(-1)}(\nbigt)
\simeq
 \Xi^{(0)}(\nbigt)\otimes\newTate(-1),
\]
\[
 \psi^{(a)}(\nbigt)
\lrarr
\psi^{(a-1)}(\nbigt)
\simeq
 \psi^{(a)}(\nbigt)\otimes\newTate(-1).
\]
We also have
$\nbigq\lrarr\psi^{(0)}(\nbigt)
\simeq
 \psi^{(1)}(\nbigt)\otimes\newTate(-1)
\lrarr
\nbigq\otimes\newTate(-1)$.
The morphisms induce 
$N:\Glue(\nbigt,\nbigq,u,v)
\lrarr
 \Glue(\nbigt,\nbigq,u,v)
\otimes
 \newTate(-1)$.
\begin{lem}
\label{lem;11.1.18.2}
The morphism $N$ is $0$.
\hfill\qed
\end{lem}
\pf
The map
$\Ker\bigl(
 \Xi^{(0)}\nbigt
\oplus\nbigq
\stackrel{\delta_0-v}{\lrarr}
 \psi^{(0)}\nbigt
 \bigr)
\lrarr
 \bigl(
 \Xi^{(0)}\nbigt
\oplus
 \nbigq
 \bigr)
 \otimes\newTate(-1)$
factors through
$\psi^{(0)}\nbigt$
by construction,
which implies the claim of the lemma.
\hfill\qed

\subsection{Dependence on the function $t$}

To distinguish the dependence of
$\Pi^{a,b}_{\star}$ on $t$,
we use the symbol $\Pi^{a,b}_{t\star}$.
We use the symbols
$\Pi^{a,b}_{t,\ast!}$,
$\Xi^{(a)}_t$, $\psi^{(a)}_t$
and $\phi^{(a)}_t$ 
in similar meanings.
We have the following morphism
denoted by $\nbign_t$:
\[
 \Pi^{a,b}_{t\star}(\nbigt)
\lrarr
 \Pi^{a-1,b-1}_{t\star}(\nbigt)
\simeq
 \Pi^{a,b}_{t\star}(\nbigt)
\otimes
 \newTate(-1)
\]
The induced morphisms for
$\Pi^{a,b}_{t,\ast!}(\nbigt)$,
$\Xi^{(a)}_t(\nbigt)$, $\psi^{(a)}_t(\nbigt)$
and $\phi^{(a)}_t(\nbigt)$ 
are also denoted by $\nbign_t$.

Let $\varphi$ be a holomorphic function.
Let $s=e^{\varphi}\, t$.
Let us compare the functors for $t$
and $s$.

\begin{lem}
\label{lem;11.2.22.10}
We naturally have
$\bigl(
 \Pi^{a,b}_{s,\star}(\nbigt),\nbign_s
 \bigr)
\simeq
 \Def_{\varphi}\bigl(
 \Pi^{a,b}_{t,\star}(\nbigt),\nbign_t
\bigr)$.
\end{lem}
\pf
We have
$\bigl(
 \Pi^{a,b}_{s}(\nbigt),\nbign_s
 \bigr)
\simeq
 \Def_{\varphi}\bigl(
 \Pi^{a,b}_{t}(\nbigt),\nbign_t
\bigr)$
from Lemma \ref{lem;10.12.25.21}.
Then, the claim of the lemma immediately follows.
\hfill\qed

\begin{cor}
\label{cor;11.1.18.1}
We naturally have
$\bigl(
 \Pi^{a,b}_{s,\ast!}(\nbigt),\nbign_s
 \bigr)
\simeq
 \Def_{\varphi}\bigl(
 \Pi^{a,b}_{t,\ast!}(\nbigt),\nbign_t
\bigr)$.
In particular, 
\[
\bigl(
 \psi^{(a)}_{s}(\nbigt),\nbign_s
 \bigr)
\simeq
 \Def_{\varphi}\bigl(
 \psi^{(a)}_{t}(\nbigt),\nbign_t
\bigr),
\quad
\bigl(
 \Xi^{(a)}_{s}(\nbigt),\nbign_s
 \bigr)
\simeq
 \Def_{\varphi}\bigl(
 \Xi^{(a)}_{t}(\nbigt),\nbign_t
\bigr).
\]
We also have
$\bigl(
 \phi^{(a)}_{s}(\nbigt),\nbign_s
 \bigr)
\simeq
 \Def_{\varphi}\bigl(
 \phi^{(a)}_{t}(\nbigt),\nbign_t
\bigr)$.
In particular,
\[
 \bigl(
 \Gr^{W(\nbign_s)}
 \psi^{(a)}_{s}(\nbigt),\nbign^{(0)}_s
 \bigr)
\simeq
\bigl(
 \Gr^{W(\nbign_t)}
 \psi^{(a)}_{t}(\nbigt),\nbign^{(0)}_t
\bigr),
\]
\[
 \bigl(
 \Gr^{W(\nbign_s)}
 \phi^{(a)}_{s}(\nbigt),\nbign^{(0)}_s
 \bigr)
\simeq
\bigl(
 \Gr^{W(\nbign_t)}
 \phi^{(a)}_{t}(\nbigt),\nbign^{(0)}_t
\bigr),
\]
where $W(\nbign_{\kappa})$ $(\kappa=s,t)$
denote the weight filtrations of $\nbign_{\kappa}$,
and $\nbign^{(0)}_{\kappa}$ are induced nilpotent morphisms.
\hfill\qed
\end{cor}

Let $\nbigt'$ be an $\nbigr_X(\ast t)$-module,
which is strictly specializable along $t$.
Let $\nbigq$ be a strict $\nbigr_{X_0}$-module
with morphisms
\[
\begin{CD}
 \psi_t^{(1)}(\nbigt')
 @>{u}>>
 \nbigq
 @>{v}>>
 \psi_t^{(0)}(\nbigt')
\end{CD}
\]
such that $v\circ u=\nbign_t$.
Then, we have an $\nbigr_X$-module
$\Glue_t(\nbigt',\nbigq,u,v)$.

We have $\nbign_t:\nbigq\lrarr \nbigq\otimes\newTate(-1)$
induced by
$\nbigq\lrarr\psi^{(0)}_t(\nbigt')
\simeq
 \psi^{(1)}(\nbigt')\otimes\newTate(-1)
\lrarr \nbigq\otimes\newTate(-1)$.
We put
$(\nbigqtilde,\nbign_s):=
\Def_{\varphi}(\nbigq,\nbign_t)$.
By Corollary \ref{cor;11.1.18.1},
we have the following naturally induced morphisms:
\[
\begin{CD}
 \psi_s^{(1)}(\nbigt')
 @>{\utilde}>>
 \nbigqtilde
 @>{\vtilde}>>
 \psi_s^{(0)}(\nbigt')
\end{CD}
\]
We have $\vtilde\circ\utilde=\nbign_s$.
We obtain an $\nbigr$-triple
$\Glue_{s}(\nbigt',\nbigqtilde,\utilde,\vtilde)$.
\begin{lem}
\label{lem;11.1.18.10}
$\Glue_t(\nbigt',\nbigq,u,v)
\simeq
 \Glue_{s}(\nbigt',\nbigqtilde,\utilde,\vtilde)$
naturally.
\end{lem}
\pf
We have the induced morphism
\[
 \nbign_t:\Glue(\nbigt',\nbigq,u,v)
\lrarr\Glue(\nbigt',\nbigq,u,v)\otimes\newTate(-1),
\]
and by construction, we have a natural isomorphism
\[
 \Def_{\varphi}\bigl(
 \Glue(\nbigt',\nbigq,u,v),\nbign_t
 \bigr)
\simeq
 \bigl(
 \Glue(\nbigt',\nbigqtilde,\utilde,\vtilde),\nbign_s
\bigr).
\]
Because $\nbign_t=0$ on
$\Glue(\nbigt',\nbigq,u,v)$
as remarked in Lemma \ref{lem;11.1.18.2},
we obtain the desired isomorphism.
\hfill\qed

\subsection{Compatibility with push-forward}

Let us consider the compatibility of the functors
$\Pi^{a.b}_{!\ast}$ with the push-forward.
We use the notation in
\S\ref{subsection;11.1.17.20}.
Let $\nbigt$ be a good
$\nbigr_{X(\ast H_X)}(\ast t)$-triple
which is strictly specializable along $t$.
Assume that the support of $\nbigt$ is
proper over $Y$
with respect to $F$.
We have a natural isomorphism
$F_{\dagger}\bigl(
 \Pi^{a,b}\nbigt\bigr)
\simeq
 \Pi^{a,b}F_{\dagger}(\nbigt)$
of good $\nbigr_{Y(\ast H_Y)}(\ast t)$-triples.
Assume that
$F_{\dagger}^i\psitilde_u\nbigt$
is strict
for any $u\in\real\times\cnum$
and $i\in\seisuu$.
Then, 
according to \cite{sabbah2},
$F^i_{\dagger}\bigl(
 \Pi^{a,b}\nbigt\bigr)$
are strictly specializable along $t$,
and 
we have
$F_{\dagger}^i\bigl(
 \Pi_{\star}^{a,b}\nbigt\bigr)
\simeq
 \Pi_{\star}^{a,b}F^{i}_{\dagger}\nbigt$
according to Corollary \ref{cor;10.8.3.11}.
Then, we obtain the following morphism:
\begin{equation}
\label{eq;11.1.17.30}
\Pi^{a,b}_{\ast !} F_{\dagger}^i\nbigt
\lrarr
 F_{\dagger}^i\Pi^{a,b}_{\ast !}\nbigt
\end{equation}

\begin{prop}
\label{prop;10.8.7.2}
The morphism {\rm(\ref{eq;11.1.17.30})}
is an isomorphism.
In particular, we have natural isomorphisms
\[
 F^i_{\dagger}\psi^{(a)}\nbigt
\simeq
 \psi^{(a)}F^i_{\dagger}\nbigt,
\quad
 F_{\dagger}^i\Xi^{(a)}\nbigt
\simeq
 \Xi^{(a)}F_{\dagger}^i\nbigt. 
\]
\end{prop}
\pf
It follows from Proposition \ref{prop;11.1.17.22}.
\hfill\qed

\vspace{.1in}
Let $\nbigt$ be an $\nbigr$-triple
which is strictly specializable along $t$.
Assume that
$F^i_{\dagger}\psitilde_u\nbigt$ $(u\in\real\times\cnum)$
and $F_{\dagger}^i\phi^{(0)}\nbigt$
are strict.
According to \cite{sabbah2},
$F^i_{\dagger}\nbigt$
are strictly specializable along $t$.
The following lemma is proved in \cite{sabbah2}
with a different method.

\begin{cor}
We have a natural isomorphism
$F_{\dagger}^i\phi^{(0)}(\nbigt)
\simeq
 \phi^{(0)} F_{\dagger}^i\nbigt$.
\end{cor}
\pf
We have the following descriptions:
\[
 \phi^{(0)}(\nbigt)=H^1\Bigl(
 \nbigt[!t]\lrarr
 \Xi^{(0)}(\nbigt)\oplus\nbigt\lrarr
 \nbigt[\ast t]
 \Bigr)
\]
\[
 \phi^{(0)}\bigl(F_{\dagger}^i\nbigt\bigr)
=H^1\bigl(
 (F_{\dagger}^i\nbigt)[!t]
\lrarr
 \Xi^{(0)} F_{\dagger}^i\nbigt
\oplus
 F_{\dagger}^i\nbigt
\lrarr
 (F_{\dagger}^i\nbigt)[\ast t]
 \bigr)
\]
Then, the claim of the lemma follows
from Proposition \ref{prop;10.8.7.2}.
\hfill\qed

\subsection{Beilinson functors along 
general holomorphic functions}

Let $X$ be a complex manifold.
Let $g$ be a holomorphic function on $X$.
Let $\iota_g:X\lrarr X\times\cnum_t$ be the graph.
Let $\nbigt$ be a coherent $\nbigr_X$-triple,
which is strictly specializable along $g$.
We define
$\psi^{(a)}_g(\nbigt):=
 \psi^{(a)}\bigl(
 \iota_{g\dagger}\nbigt
 \bigr)$
and
$\phi^{(a)}_g(\nbigt):=
 \phi^{(a)}\bigl(
 \iota_{g\dagger}\nbigt
 \bigr)$.
If $\nbigt$ is integrable,
they are also naturally integrable.

Assume that $\nbigt$ is strict for simplicity.
If there exists
an $\nbigr_X$-triples $\nbigt'$
such that 
$\iota_{g\dagger}\nbigt'=
 \Pi^{a,b}_{!}(\iota_{g\dagger}\nbigt)$,
it is uniquely determined up to isomorphisms,
and denoted by
$\Pi^{a,b}_{g!}(\nbigt)$.
We use the notation 
$\Pi^{a,b}_{g\ast}(\nbigt)$,
$\Pi^{a,b}_{\ast !}(\iota_{g\dagger}\nbigt)$
and 
$\Xi^{(a)}_g(\nbigt)$
similar meanings.
If $\nbigt$ is integrable,
they are naturally integrable.

\section{Comparison of nearby cycle functors}
\label{subsection;11.1.24.10}

\subsection{Statements}

Let $\nbigt$ be an $\nbigr_{X(\ast H)}(\ast t)$-triple
which is strictly specializable along $t$.
Let $\iota$ denote the inclusion
$X_0\times\{0\}\lrarr X$.
In this subsection,
we distinguish the $\nbigr_{X_0}$-triple
$\psitilde_{-\vecdelta}(\nbigt)$
and the $\nbigr_X$-triple
$\iota_{\dagger}\psitilde_{-\vecdelta}(\nbigt)$.
We shall compare 
the $\nbigr_X$-triples
$\iota_{\dagger}\psitilde_{-\vecdelta}(\nbigt)$
and $\psi^{(1)}(\nbigt)$.
We shall prove the following proposition
in \S\ref{subsection;13.3.28.1}--\ref{subsection;13.3.28.2}.

\begin{prop}
\label{prop;10.12.19.1}
We have a natural isomorphism
$\Psi:
 \iota_{\dagger}
 \psitilde_{-\vecdelta}(\nbigt)
\simeq
 \psi^{(1)}(\nbigt)\otimes\nbigu(1,0)$
with the following property:
\begin{itemize}
\item
 Let 
 $\iota_{\dagger}\psitilde_{-\vecdelta}(\nbigt)
\simeq \psi^{(0)}(\nbigt)\otimes\nbigu(0,1)$ be 
the induced isomorphism.
The following diagram is commutative:
\begin{equation}
 \label{eq;11.1.24.1}
 \begin{CD}
 \psi^{(1)}(\nbigt)
 @>{\delta_0\circ\gamma_0}>>
 \psi^{(0)}(\nbigt)
 \\
 @V{\simeq}VV @V{\simeq}VV \\
 \iota_{\dagger}\psitilde_{-\vecdelta}(\nbigt)
 \otimes\nbigu(-1,0)
 @>{\nbign}>>
 \iota_{\dagger}\psitilde_{-\vecdelta}(\nbigt)
 \otimes\nbigu(0,-1)
 \end{CD}
\end{equation}
Here, 
$\nbign=(N',N'')$ is as in {\rm\S\ref{subsection;13.5.4.40}},
i.e., $N'$ and $N''$ are the nilpotent part of
 $t\del_t$.
If $\nbigt$ is integrable,
the above isomorphism is also integrable.
\end{itemize}
See {\rm \S\ref{subsection;10.12.27.2}}
for  the integrable $\nbigr$-triples $\nbigu(p,q)$.
\end{prop}

Let $\nbigs:\nbigt\lrarr
 \nbigt^{\ast}\otimes\newTate(-w)$
be a hermitian adjoint.
We have the induced morphism
\[
 \Pi^{1,N}\nbigt\lrarr
 (\Pi^{-N+1,0}\nbigt)^{\ast}
 \otimes \newTate(-w). 
\]
Hence, we have the following
naturally induced morphisms:
\begin{multline}
 \psi^{(1)}(\nbigt)
\lrarr
 \psi^{(1)}(\nbigt^{\ast})
 \otimes\newTate(-w)
\simeq
 \psi^{(0)}(\nbigt^{\ast})
 \otimes\newTate(-w+1)
 \\
\simeq
 \psi^{(1)}(\nbigt)^{\ast}
 \otimes\newTate(-w+1)
\end{multline}
The composite is denoted by
$\psi^{(1)}(\nbigs)$.
We shall prove the following proposition
in \S\ref{subsection;13.3.28.3}.

\begin{prop}
\label{prop;10.12.19.5}
The following diagram is commutative:
\[
 \begin{CD}
 \psi^{(1)}(\nbigt)
 @>{\psi^{(1)}(\nbigs)}>>
 \psi^{(1)}(\nbigt)^{\ast}\otimes
 \newTate(-w+1)\\
 @V{\simeq}VV @V{\simeq}VV \\
 \iota_{\dagger}
 \psitilde_{-\vecdelta}(\nbigt)\otimes\nbigu(-1,0)
 @>>>
 \bigl(
 \iota_{\dagger}
 \psitilde_{-\vecdelta}(\nbigt)\otimes\nbigu(-1,0)
 \bigr)^{\ast}\otimes\newTate(-w+1)
 \end{CD}
\]
Here, the vertical arrows are 
induced by the isomorphism
in Proposition {\rm\ref{prop;10.12.19.1}},
and the lower horizontal arrow is
induced by $\psitilde_{-\vecdelta}(\nbigs)$
and the polarization $(-1,-1)$ of
$\nbigu(-1,0)$,
given in {\rm\S\ref{subsection;10.12.27.2}}.
\end{prop}

\begin{rem}
If $\nbigt$ is integrable,
our map in $\Psi$ in 
Proposition {\rm\ref{prop;10.12.19.1}}
is commutative with the action of
$\lambda^2\del_{\lambda}$.
\hfill\qed
\end{rem}

\subsection{Preliminary (1)}
\label{subsection;13.3.28.1}

Let $\nbigm$ be an $\nbigr(\ast t)$-module
strictly specializable along $t$.
Assume that
the cokernel of the nilpotent part of
$-\deldel_tt$ on $\psi_{-\vecdelta}(\nbigm)$
is strict.
We have the following commutative diagram:
\[
 \begin{CD}
 \psi_0(\nbigm[!t])
 @>{\psi_0(\varphi)}>>
 \psi_0(\nbigm[\ast t])\\
 @A{-\deldel_t}A{\simeq}A 
 @V{t}V{\simeq}V\\
 \psi_{-\vecdelta}(\nbigm)
 @>{-t\deldel_t}>>
 \psi_{-\vecdelta}(\nbigm)
 \end{CD}
\]
Hence, the cokernel of $\psi_u(\varphi)$ 
are strict for any $u$.

Let $\nbigk(\nbigm)$ and $\nbigq(\nbigm)$ be 
the kernel and the cokernel of
$\varphi:\nbigm[!t]\lrarr\nbigm[\ast t]$.
We have the exact sequence:
\[
 0\lrarr\psi_0\nbigk(\nbigm)\lrarr
 \psi_0(\nbigm[!t])
 \stackrel{\psi_0(\varphi)}{\lrarr}
 \psi_0(\nbigm[\ast t])\lrarr
 \psi_0\nbigq(\nbigm)\lrarr 0
\]
Let $N:
 \psitilde_{-\vecdelta}(\nbigm)\lrarr
 \psitilde_{-\vecdelta}(\nbigm)\lambda^{-1}$
be induced by $-t\del_t$.
We obtain the following commutative diagram:
\[
 \begin{CD}
  \psi_0\nbigk(\nbigm)@>>>
 \psi_0(\nbigm[!t])
 @>{\psi_0(\varphi)}>>
 \psi_0(\nbigm[\ast t])@>>>
 \psi_0\nbigq(\nbigm)
 \\
 @A{\mu'_{\nbigk(\nbigm)}}A{\simeq}A
 @A{-\lambda\del_t}A{\simeq}A 
@V{\lambda^{-1}t}V{\simeq}V 
 @V{\mu'_{\nbigq(\nbigm)}}V{\simeq}V
 \\
 \Ker N
 @>>>
 \psi_{-\vecdelta}(\nbigm)
 @>{N}>>
 \psi_{-\vecdelta}(\nbigm)\,\lambda^{-1}
 @>>>
 \Cok N
 \end{CD}
\]
We obtain the isomorphism
$\mu_{\nbigk(\nbigm)}:
 \iota_{\dagger}\Ker N
\lrarr
 \nbigk(\nbigm)$
given as follows:
\[
 \iota_{\dagger}(\Ker N)
=\bigoplus_{n=0}^{\infty}
 \Ker N\cdot (dt/\lambda)^{-1}
 \cdot \deldel_t^n
\lrarr
 \nbigk(\nbigm)
=\bigoplus_{n=0}^{\infty}
 \psi_0(\nbigk(\nbigm))\,\deldel_t^n
\]
\[
 \mu_{\nbigk(\nbigm)}\Bigl(
\sum a_n\cdot(dt/\lambda)^{-1}\,\deldel_t^n
 \Bigr)
:=\sum \mu'_{\nbigk(\nbigm)}(a_n)\,\deldel_t^n
\]
We also obtain the following isomorphism
$\nbigq(\nbigm)
\lrarr 
 \iota_{\dagger}(\Cok N)$:
\[
 \nbigq(\nbigm)
=\bigoplus_{n=0}^{\infty}
 \psi_0(\nbigq(\nbigm))\,\deldel_t^n
\lrarr
 \iota_{\dagger}(\Cok N)
=\bigoplus_{n=0}^{\infty}
 \Cok N\cdot (dt/\lambda)^{-1}
 \cdot \deldel_t^n
\]
\[
 \mu_{\nbigq(\nbigm)}\Bigl(
 \sum b_n\,\deldel_t^n \Bigr)
:=\sum \mu'_{\nbigq(\nbigm)}(b_n)
  \cdot(dt/\lambda)^{-1}
  \,\deldel_t^n
\]

\subsection{Preliminary (2)}
\label{subsection;10.12.19.2}

Let $\nbigt=(\nbigm',\nbigm'',C)$ be
a coherent $\nbigr_{X(\ast H)}(\ast t)$-triple 
which is strictly specializable along $t$.
We have the induced morphism
$\nbign:\psitilde_{-\vecdelta}(\nbigt)
\lrarr
\psitilde_{-\vecdelta}(\nbigt)
 \otimes\newTate(-1)$.
We have
$\Ker\nbign
\!=\!\bigl(
 \lambda\cdot\Cok N',\Ker N'',C_1
 \bigr)$
and
$\Cok\nbign
\!=\!\bigl(
 \lambda\cdot\Ker N',\Cok N'',C_2
 \bigr)$,
where $C_i$ are naturally induced pairings.
Let $\nbigk$ and $\nbigq$
be the kernel and the cokernel
of $\nbigt[!t]\lrarr\nbigt[\ast t]$.
\begin{lem}
The pair of morphisms
$(\mu_{\nbigk(\nbigm')},\mu_{\nbigq(\nbigm'')})$
gives an isomorphism
$\nbigq
\stackrel{\simeq}{\lrarr}
 \iota_{\dagger}\Cok\nbign\otimes\nbigu(-1,0)$.
Similarly,
the pair of morphisms
$(\mu_{\nbigq(\nbigm')},\mu_{\nbigk(\nbigm'')})$
gives an isomorphism
$\iota_{\dagger}
 \Ker\nbign\otimes\nbigu(-1,0)
\stackrel{\simeq}{\lrarr}
 \nbigk$.
\end{lem}
\pf
We have only to check the compatibility of
pairings.
Let us check the first claim.
Let $\lambda\,t^{-1}f\in V_0(\nbigq(\nbigm''))$
be lifted to $V_0\nbigm''[\ast t]$,
and $-\deldel_t\otimes g\in V_0(\nbigk(\nbigm'))$
be mapped to $V_0\nbigm'[!t]$.
Note $t\deldel_t g=0$ in $V_{<0}\nbigm'$.
Let $\phi=(\sqrt{-1}/2\pi)\cdot \chi\cdot \varphi\cdot
 dt\,d\tbar$,
where $\varphi$ is a $C^{\infty}$-top form on
$X_0$ with compact support,
and $\chi$ is a cut function on $\cnum_t$
around $t=0$.
We have
\begin{multline}
\Bigl\langle
 C\bigl(-\deldel_t\otimes g,
 \sigma^{\ast}(\lambda t^{-1}f)\bigr),\phi
\Bigr\rangle
=
\bigl\langle
 C\bigl(\del_t\otimes g,
 \sigma^{\ast}(t^{-1}f)\bigr),\phi
\bigr\rangle \\
=-\bigl\langle
 C(g,\sigma^{\ast}f),\tbar^{-1}\del_t\phi
 \bigr\rangle 
=-\bigl\langle
 C(g,\sigma^{\ast}f),\,
 \tbar^{-1}|t|^{2s}\del_t\phi
 \bigr\rangle_{|s=0}
 \\
=-\bigl\langle
 C(g,\sigma^{\ast}f),\,
 \del_t(\tbar^{-1}|t|^{2s}\phi)
 \bigr\rangle_{|s=0}
+\bigl\langle
 C(g,\sigma^{\ast}f),\,
 s\,|t|^{2(s-1)}\phi
 \bigr\rangle_{|s=0} \\
=\bigl\langle
 \psi_{-\vecdelta}C([g],\sigma^{\ast}[f]),\,
 \varphi
 \bigr\rangle
=\bigl\langle
 \iota_{\dagger}\psi_{-\vecdelta}C\bigl(
 [g]\cdot (dt/\lambda)^{-1},\,\,
 \sigma^{\ast}\bigl(
 [f]\cdot(dt/\lambda)^{-1}
 \bigr)
 \bigr),\phi
 \bigr\rangle
\end{multline}
Let us check the second claim.
Let $-\deldel_t\otimes g\in V_0(\nbigk(\nbigm''))$
be lifted to $V_0\nbigm''[!t]$,
and let 
$\lambda\,t^{-1}f\in V_0(\nbigq(\nbigm'))$
be lifted to $V_0\nbigm'[\ast t]$.
We have the following:
\begin{multline}
\bigl\langle
 C\bigl(\lambda\,t^{-1}f,
 \sigma^{\ast}(-\deldel_t\otimes g)\bigr),\,
 \phi \bigr\rangle
=\bigl\langle
 C\bigl(t^{-1}f,
 \sigma^{\ast}(\del_t\otimes g)\bigr),\,
 \phi \bigr\rangle
 \\
=-\bigl\langle
 C(f,\sigma^{\ast}g),\,
 t^{-1}\delbar_t\phi
 \bigr\rangle 
=-\bigl\langle
 C(f,\sigma^{\ast}g),\,
 t^{-1}|t|^{2s}
 \delbar_t\phi
 \bigr\rangle_{|s=0}
 \\
=-\bigl\langle
 C(f,\sigma^{\ast}g),\,
 \delbar_t\bigl(
 t^{-1}|t|^{2s}\phi
 \bigr)
 \bigr\rangle_{|s=0}
+\bigl\langle
 C(f,\sigma^{\ast}g),\,
 s\, |t|^{2(s-1)}\phi
 \bigr\rangle_{|s=0} \\
=\bigl\langle
 \psitilde_{-\vecdelta}
 C([f],\sigma^{\ast}[g]),\,\varphi
 \bigr\rangle
=\bigl\langle
 \iota_{\dagger}\psitilde_{-\vecdelta}
 C\bigl([f]\,(dt/\lambda)^{-1},\, 
 \sigma^{\ast}([g]\,(dt/\lambda)^{-1})\bigr),
 \,\,\phi
 \bigr\rangle
\end{multline}
Thus, we are done.
\hfill\qed

\subsection{Construction of isomorphisms}
\label{subsection;13.3.28.2}

We set
$\II^{a,b}_0:=\II^{a,b}_2=\II^{-b+1,-a+1}_1$.
We have the isomorphism
\[
 \psitilde_{-\vecdelta}(\nbigt)
\stackrel{\simeq}{\lrarr}
\Cok\Bigl(
 \nbign:
 \psitilde_{-\vecdelta}(\Pi^{1,N}\nbigt)
\lrarr
 \psitilde_{-\vecdelta}(\Pi^{1,N}\nbigt)
 \otimes\newTate(-1)
\Bigr)
\]
given as follows:
\begin{itemize}
\item
We use
$\psitilde_{-\vecdelta}(\nbigm')
\stackrel{\simeq}{\llarr}
 \Ker\Bigl(
 \psitilde_{-\vecdelta}(\nbigm')
\otimes \II_0^{-N+1,0}
\llarr
 \psitilde_{-\vecdelta}(\nbigm')
\otimes \II_0^{-N+1,0}\cdot \lambda
 \Bigr)$
induced by the projection
with $s\longmapsto -s$:
\[
 \II_0^{-N+1,0}\lambda
\lrarr
 \nbigo_{\cnum_{\lambda}}\,\lambda;
\quad
 \sum_{j=1}^{N-1}a_j\,(\lambda s)^{-j}\,\lambda
 \longmapsto -a_{1}.
\] 
\item
We use
$\psitilde_{-\vecdelta}(\nbigm'')
\stackrel{\simeq}{\lrarr}
 \Cok\Bigl(
 \psitilde_{-\vecdelta}(\nbigm'')
 \otimes \II_0^{1,N}
\lrarr
 \psitilde_{-\vecdelta}(\nbigm'')
 \otimes \II_0^{1,N}\,\lambda^{-1}
 \Bigr)$
induced by the inclusion
with $s\longmapsto -s$:
\[
 \nbigo_{\cnum_{\lambda}}
\lrarr
 \II_0^{1,N}\,\lambda^{-1};
 \quad
 a\longmapsto -a\,s.
\]
\end{itemize}
We have the isomorphism
\[
 \psitilde_{-\vecdelta}(\nbigt)
\stackrel{\simeq}{\llarr}
\Ker\Bigl(
 \nbign:\psitilde_{-\vecdelta}(\Pi^{-N,1}\nbigt)
\lrarr
 \psitilde_{-\vecdelta}(\Pi^{-N,1}\nbigt)
\otimes\newTate(-1)
 \Bigr)
\]
given as follows:
\begin{itemize}
\item
 We use
 $\psitilde_{-\vecdelta}(\nbigm')
 \stackrel{\simeq}{\lrarr}
 \Cok\Bigl(
 \psitilde_{-\vecdelta}(\nbigm')\otimes \II_0^{0,N+1}
\llarr
 \psitilde_{-\vecdelta}(\nbigm')\otimes \II_0^{0,N+1}\,\lambda
 \Bigr)$
given by the inclusion
$\nbigo_{\cnum_{\lambda}}
\lrarr
 \II_0^{0,N+1}$;
 $a\longmapsto a$.
\item
 We use 
 $\psitilde_{-\vecdelta}(\nbigm'')
 \stackrel{\simeq}{\llarr}
 \Ker\Bigl(
 \psitilde_{-\vecdelta}(\nbigm'')\otimes \II_0^{-N,1}
\lrarr
 \psitilde_{-\vecdelta}(\nbigm'')\otimes \II_0^{-N,1}
 \,\lambda^{-1}
 \Bigr)$ 
induced by the projection 
$\II_0^{-N,1}\lrarr \nbigo_{\cnum_{\lambda}}$;
$\sum_{j=0}^{N} a_j\,(\lambda s)^{-j}
\longmapsto a_0$.
\end{itemize}

Then, we obtain the following isomorphisms:
\[
 \Psi_1:\psi^{(1)}(\nbigt)
\simeq
 \iota_{\dagger}\Cok(\nbign)\otimes\nbigu(-1,0)
\simeq 
 \iota_{\dagger}\psitilde_{-\vecdelta}(\nbigt)
 \otimes\nbigu(-1,0)
\]
\[
 \Psi_2:
\psi^{(1)}(\nbigt)
\simeq
 \iota_{\dagger}\Ker(\nbign)\otimes\nbigu(-1,0)
\simeq
 \iota_{\dagger}
 \psitilde_{-\vecdelta}(\nbigt)\otimes\nbigu(-1,0)
\]
\begin{lem}
We have $\Psi_1=\Psi_2$.
They are denoted by $\Psi$.
\end{lem}
\pf
It can be checked by a direct computation.
We will give an indication.
The first component of
$\Psi_2\circ\Psi_1^{-1}$
is the composition of the following:
\begin{multline}
\label{eq;10.12.19.3}
\iota_{\dagger}\psitilde_{-\vecdelta}(\nbigm')
 \,\lambda^{-1}
\simeq
\iota_{\dagger}\Cok\Bigl(
 \psitilde_{-\vecdelta}(\nbigm')\otimes
 \II_0^{0,N}\,\lambda^{-1}
\llarr
 \psitilde_{-\vecdelta}(\nbigm')\otimes
 \II_0^{0,N}
 \Bigr) \\
\simeq 
\nbigq(\Pi^{0,N}\nbigm') 
\stackrel{\Phi}{\simeq}
 \nbigk(\Pi^{-N+1,0}\nbigm') \\
\simeq
 \iota_{\dagger}
 \Ker\Bigl(
 \psitilde_{-\vecdelta}(\nbigm')
 \otimes \II_0^{-N+1,0}\,\lambda^{-1}
 \llarr
\psitilde_{-\vecdelta}(\nbigm')
 \otimes \II_0^{-N+1,0}
 \Bigr)
\simeq
 \iota_{\dagger}
 \psitilde_{-\vecdelta}(\nbigm')\,\lambda^{-1}
\end{multline}
The morphism $\Phi$ is induced by
the following:
\[
 \begin{CD}
 \Ker @>{=}>>
 \Ker @>>>
 \nbigk(\Pi^{-N+1,0}\nbigm') \\
 @VVV @VVV @VVV \\
 \Pi_!^{0,N+1}\nbigm'
 @>>>
 \Pi_!^{-N+1,N+1}\nbigm'
 @>>>
 \Pi_!^{-N+1,0}\nbigm' \\
 @VVV @VVV @VVV \\
 \Pi_{\ast}^{0,N+1}\nbigm'
 @>>>
 \Pi_{\ast}^{-N+1,N+1}\nbigm'
 @>>>
 \Pi_{\ast}^{-N+1,0}\nbigm' \\
 @VVV @VVV @VVV \\
 \nbigq(\Pi^{0,N+1}\nbigm')
 @>>>
 \Cok
 @>{=}>>
 \Cok
 \end{CD}
\]
Hence, the composite of the morphisms
in (\ref{eq;10.12.19.3}) is induced by
the following diagram:
{\small
\[
 \begin{CD}
 @.
{\scriptstyle
 \nbigp'
 \otimes \II_0^{0,N}}
 @>>>
{\scriptstyle
 \nbigp'
 \otimes \II_0^{-N,N}}
 @>>>
{\scriptstyle
 \nbigp'
 \otimes \II_0^{-N,0}}
 @>>>
{\scriptstyle
\nbigp'\,\lambda^{-1}}
 \\
 @. @V{-t\del_t}VV @VVV @VVV @. \\
{\scriptstyle
\nbigp'\,\lambda^{-1}}
 @>>>
{\scriptstyle
\nbigp'
 \otimes \II_0^{0,N}\,\lambda^{-1}}
 @>>>
{\scriptstyle
\nbigp'
 \otimes \II_0^{-N,N}\,\lambda^{-1}}
 @>>>
{\scriptstyle
 \nbigp'
 \otimes \II_0^{-N,0}\,\lambda^{-1}}
 \end{CD}
\]
}
Here, $\nbigp'$ denotes $\psitilde_{-\vecdelta}(\nbigm')$.
Hence, the composite is the identity.
Let us look at the second component of
$\Psi_2^{-1}\circ\Psi_1$.
It is the composite of the following morphisms:
\begin{multline}
 \iota_{\dagger}\psitilde_{-\vecdelta}(\nbigm'')
\simeq
 \iota_{\dagger}\Cok\Bigl(
 \psitilde_{-\vecdelta}(\nbigm'')
 \otimes \II_0^{1,N}
 \lrarr
 \psitilde_{-\vecdelta}(\nbigm'')
 \otimes \II_0^{1,N}\,\lambda^{-1}
 \Bigr)
 \\
\simeq
 \nbigq(\Pi^{1,N}\nbigm'')
\simeq
 \nbigq(\Pi^{1,N}\nbigm'')
 \\
\simeq
 \iota_{\dagger}\Ker\Bigl(
 \psitilde_{-\vecdelta}(\nbigm'')
 \otimes \II_0^{-N,1}
\lrarr
 \psitilde_{-\vecdelta}(\nbigm'')
 \otimes \II_0^{-N,1}\,\lambda^{-1}
 \Bigr)
\simeq
 \iota_{\dagger}
 \psitilde_{-\vecdelta}(\nbigm'')
\end{multline}
The composite of the morphisms
is obtained from the following diagram:
\[
 \begin{CD}
 @.
{\scriptstyle
 \nbigp''
 \otimes \II_0^{1,N}}
 @>>>
{\scriptstyle
 \nbigp''
 \otimes \II_0^{-N,1}}
 @>>>
{\scriptstyle
 \nbigp''
 \otimes \II_0^{-N,0}}
 @>>>
{\scriptstyle
 \nbigp''
 }
 \\
 @. @V{-t\del_t}VV @VVV @VVV @. \\
{\scriptstyle
 \nbigp''}
 @>>>
{\scriptstyle
 \nbigp''
 \otimes \II_0^{1,N}\,\lambda^{-1}}
 @>>>
{\scriptstyle
 \nbigp''
 \otimes \II_0^{-N,N}\,\lambda^{-1}}
 @>>>
{\scriptstyle
 \nbigp''
 \otimes \II_0^{-N,1}\,\lambda^{-1}}
 \end{CD}
\]
Here, $\nbigp''$ denotes
$\psitilde_{-\vecdelta}(\nbigm'')$.
Then, we can check that the composite
is the identity.
\hfill\qed

\vspace{.1in}
It is easy to check that the diagram (\ref{eq;11.1.24.1})
is commutative,
because we have
$t\del_t=-N'\otimes\id+\id\otimes s$
on $\psitilde_{-\vecdelta}(\nbigm')\otimes \II_0^{p,q}$
and 
$t\del_t=-N''\otimes\id+\id\otimes s$
on $\psitilde_{-\vecdelta}(\nbigm'')\otimes \II_0^{p,q}$.
Thus, we obtain Proposition 
\ref{prop;10.12.19.1}.

\subsection{Hermitian adjoint}
\label{subsection;13.3.28.3}

Let $\nbigs:\nbigt\lrarr
 \nbigt^{\ast}\otimes \newTate(-w)$
be a hermitian adjoint.
Let us check Proposition \ref{prop;10.12.19.5}.
Under the isomorphism $\Psi$ 
in \S\ref{subsection;13.3.28.2},
$\psi^{(1)}(\nbigs)'':
 \psi^{(1)}(\nbigm'')
\lrarr
 \psi^{(1)}(\nbigm')\,\lambda^{-w+1}$
is induced by 
\begin{multline}
 \psitilde_{-\vecdelta}(\nbigm'')
\stackrel{b_1}{\simeq}
 \Cok\Bigl(
 \nbign'':\psitilde_{-\vecdelta}(\nbigm'')
 \otimes \II^{1,N}
\lrarr
 \psitilde_{-\vecdelta}(\nbigm'')
 \otimes \II^{1,N}\,\lambda^{-1}
 \Bigr) \\
\stackrel{f_1}{\lrarr}
 \Cok\Bigl(
 \nbign':\psitilde_{-\vecdelta}(\nbigm')
 \otimes \II^{1,N}
\lrarr
 \psitilde_{-\vecdelta}(\nbigm')
 \otimes \II^{1,N}\,\lambda^{-1}
 \Bigr) \cdot\lambda^{-w}
 \\
\stackrel{s^{-1}}{\lrarr}
 \Cok\Bigl(
 \nbign':\psitilde_{-\vecdelta}(\nbigm')
 \otimes \II^{0,N-1}
\lrarr
 \psitilde_{-\vecdelta}(\nbigm')
 \otimes \II^{0,N-1}\,\lambda^{-1}
 \Bigr) \cdot
\lambda^{-w+1}
 \\
 \stackrel{g_1}{\lrarr}
 \Cok\Bigl(
 \nbign':\psitilde_{-\vecdelta}(\nbigm')
 \otimes \II^{0,N-1}
\lrarr
 \psitilde_{-\vecdelta}(\nbigm')
 \otimes \II^{0,N-1}\,\lambda^{-1}
 \Bigr) \cdot
\lambda^{-w+1}
 \\
\stackrel{b_2}{\simeq}
\bigl(
 \psitilde_{-\vecdelta}(\nbigm')\,\lambda^{-1}
\bigr)
\cdot\lambda^{-w+1}
\end{multline}
Here, 
$b_i$ are isomorphisms
in \S\ref{subsection;13.3.28.2},
$f_1$ is induced by
$\nbigs''$,
and $g_1$ is induced by $\nbigs^{a,b}$.
The composite of the morphisms
is $-\psitilde_{-\vecdelta}(\nbigs'')$.

Let us look at the following morphisms:
\begin{multline}
 \lambda^{-1}\,\psitilde_{-\vecdelta}(\nbigm')
 \simeq
 \Ker\Bigl(
 \psitilde_{-\vecdelta}(\nbigm')
 \otimes \II^{-N+1,0}\,\lambda^{-1}
 \llarr
 \psitilde_{-\vecdelta}(\nbigm')
 \otimes \II^{-N+1,0}
 \Bigr) \\
\stackrel{f_2}{\llarr}
 \Ker\Bigl(
 \psitilde_{-\vecdelta}(\nbigm'')
 \otimes \II^{-N+1,0}\,\lambda^{-1}
 \llarr
 \psitilde_{-\vecdelta}(\nbigm'')
 \otimes \II^{-N+1,0}
 \Bigr)\,\lambda^{w} \\
\stackrel{s}{\llarr}
 \Ker\Bigl(
 \psitilde_{-\vecdelta}(\nbigm'')
 \otimes \II^{-N+2,1}\,\lambda^{-1}
 \llarr
 \psitilde_{-\vecdelta}(\nbigm'')
 \otimes \II^{-N+2,1}
 \Bigr)\,\lambda^{w-1} 
 \\
 \stackrel{g_2}{\llarr}
 \Ker\Bigl(
 \psitilde_{-\vecdelta}(\nbigm'')
 \otimes \II^{-N+2,1}\,\lambda^{-1}
 \llarr
 \psitilde_{-\vecdelta}(\nbigm'')
 \otimes \II^{-N+2,1}
 \Bigr)\,\lambda^{w-1} 
 \\
\simeq
 \psitilde_{-\vecdelta}(\nbigm'')
 \cdot\lamda^{w-1}
\end{multline}
Here, $f_2$ is induced by $\nbigs'$,
and $g_2$ is induced by $\nbigs^{a,b}$.
Hence, the composite is $-\psi_{-\vecdelta}(\nbigs')$.
Thus, we obtain Proposition \ref{prop;10.12.19.5}.
\hfill\qed

\section{Admissible specializability}
\label{subsection;13.4.12.12}

\subsection{Filtered $\nbigr$-modules}

\subsubsection{Filtered strictly specializability}

In this paper,
a filtered $\nbigr_X$-module means
an $\nbigr_X$-module with a locally finite increasing 
exhaustive $\seisuu$-indexed filtration
 in the category of $\nbigr_X$-modules.

Let $(\nbigm,L)$ be a filtered $\nbigr_X$-module.
Let $g$ be any holomorphic function on $X$.
\begin{df}
$(\nbigm,L)$ is called filtered strictly specializable
along $g$, if the following holds:
\begin{itemize}
\item
 Each $L_j\nbigm$ is strictly specializable along $g$
 with any ramified exponential twist.
\item
 The cokernel of
 $\psitilde_{g,\gminia,u}(L_j\nbigm)\lrarr
 \psitilde_{g,\gminia,u}(\nbigm)$
 is strict for any $u\in\real\times\cnum$
 and $\gminia\in\cnum[t_n^{-1}]$.
 The cokernel of 
 $\phi_{g}(L_j\nbigm)\lrarr
 \phi_{g}(\nbigm)$
 is also strict.
\hfill\qed
\end{itemize}
\end{df}
\index{filtered strictly specializable}

\begin{lem}
\label{lem;11.1.17.40}
If $(\nbigm,L)$ is filtered strictly specializable,
each $\bigl(L_j\nbigm/L_k\nbigm,L\bigr)$
is also filtered strictly specializable.
\end{lem}
\pf
In this case,
each $\nbigm/L_j\nbigm$ is also strictly specializable
along $g$ with any ramified exponential twist.
The morphisms
$L_j\nbigm\lrarr\nbigm\lrarr
 \nbigm/L_j\nbigm$ are strict with respect to
the $V$-filtrations for any ramified exponential twist.
We have the exact sequence
$0\lrarr \psitilde_{g,\gminia,u}(L_j\nbigm)
\lrarr\psitilde_{g,\gminia,u}\nbigm\lrarr
 \psitilde_{g,\gminia,u}(\nbigm/L_j\nbigm)\lrarr 0$.
We have
\[
 \Cok\Bigl(
 \psitilde_{g,\gminia,u}(L_m\nbigm)
\lrarr
 \psitilde_{g,\gminia,u}(L_j\nbigm)
 \Bigr)
\subset
  \Cok\Bigl(
 \psitilde_{g,\gminia,u}(L_m\nbigm)
\lrarr
 \psitilde_{g,\gminia,u}(\nbigm)
 \Bigr)
\]
\[
 \Cok\Bigl(
 \phi_g(L_m\nbigm)
\lrarr
 \phi_g(L_j\nbigm)
 \Bigr)
\subset
  \Cok\Bigl(
 \phi_g(L_m\nbigm)
\lrarr
 \phi_g(\nbigm)
 \Bigr)
\]
We also have the following isomorphisms:
\[
 \Cok\Bigl(
 \psitilde_{g,\gminia,u}(L_j\nbigm/L_m\nbigm)
\lrarr
 \psitilde_{g,\gminia,u}(\nbigm/L_m\nbigm)
 \Bigr)
\simeq
 \psitilde_{g,\gminia,u}(\nbigm/L_j\nbigm)
\]
\[
 \Cok\Bigl(
 \phi_{g}(L_j\nbigm/L_m\nbigm)
\lrarr
 \phi_{g}(\nbigm/L_m\nbigm)
 \Bigr)
\simeq
 \phi_{g}(\nbigm/L_j\nbigm)
\]
Then, the claim of Lemma \ref{lem;11.1.17.40} follows.
\hfill\qed

\vspace{.1in}
From the argument of the proof,
we also obtain the following.
\begin{lem}
Let $L(\psitilde_{g,\gminia,u}(\nbigm))$ 
be the naturally induced filtration
of $\psitilde_{g,\gminia,u}(\nbigm)$.
Then, 
$\Gr^L\psitilde_{g,\gminia,u}(\nbigm)$
is strict.
\hfill\qed
\end{lem}

Let $(\nbigm_i,L)$ $(i=1,2)$ be filtered $\nbigr_X$-modules
which are filtered strictly specializable along $g$.
Let $f:(\nbigm_1,L)\lrarr(\nbigm_2,L)$ be a morphism
of filtered $\nbigr_X$-modules.
Then,
we have the induced morphisms $\psitilde_{g,\gminia,u}(f):
 \psitilde_{g,\gminia,u}(\nbigm_1)
\lrarr
 \psitilde_{g,\gminia,u}(\nbigm_2)$
 and 
$\phi_g(f):
 \phi_g(\nbigm_1)\lrarr
 \phi_g(\nbigm_2)$,
 which are compatible with the filtrations $L$.

\subsubsection{Admissible specializability}

\begin{df}
\label{df;10.10.3.1}
Let $(\nbigm,L)$ be a filtered $\nbigr_X$-module,
which is filtered strictly specializable along $g$.
It is called admissibly specializable
along $g$, if the following holds:
\begin{description}
\item[\bf (P1)]
$\bigl(\psitilde_{g,\gminia,u}(\nbigm,L),\nbign\bigr)$
has a relative monodromy filtration
for each $u\in\real\times\cnum$
and $\gminia\in\cnum[t_n^{-1}]$.
\item[\bf (P2)]
$\bigl(\phi_{g}(\nbigm,L),\nbign\bigr)$
has a relative monodromy filtration.
\item[\bf (P3)]
For the natural morphisms
$\psi_{g,-\vecdelta}(\nbigm)
 \stackrel{u}{\lrarr}
 \phi_g(\nbigm)
 \stackrel{v}{\lrarr}
 \psi_{g,-\vecdelta}(\nbigm)$,
we have
$u\cdot M_{k}(N;L)\subset
 M_{k-1}(N;L)$
and 
$v\cdot M_{k}(N;L)\subset
 M_{k-1}(N;L)$,
where $u$ and $v$ are induced by
$\deldel_t$ and $t$.
\end{description}
(We shall give a review on relative monodromy filtration
in {\rm\S\ref{subsection;10.11.5.10}} below.)
\hfill\qed
\end{df}
\index{admissibly specializable}

The following lemma is  clear.
\begin{lem}
If $(\nbigm,L)$ is admissibly specializable,
any $\bigl(L_j\nbigm/L_k\nbigm,L\bigr)$
are also admissibly specializable.
\hfill\qed
\end{lem}

A filtered $\nbigr_X$-module $(\nbigm,L)$ is called integrable,
if each $L_j\nbigm$ is integrable.
The following lemma is clear
by Deligne's formula for relative monodromy filtration
(\ref{eq;10.9.27.21}, \ref{eq;10.9.27.22})  below.
\begin{lem}
If $(\nbigm,L)$ is integrable and admissibly specializable
along $g$,
the relative monodromy filtrations are also integrable.
\hfill\qed
\end{lem}
The following lemma is also clear
by Deligne's formula.
\begin{lem}
\label{lem;13.5.10.31}
Let $(\nbigm_i,L)$ $(i=1,2)$ be admissibly specializable
along $g$.
Let $f:(\nbigm_1,L)\lrarr(\nbigm_2,L)$ be any morphism.
Then, the induced morphisms
$\psitilde_{g,\gminia,u}(f)$
and $\phi_g(f)$
are compatible with $M(N;L)$.
\hfill\qed
\end{lem}

Let $\iota:Y\lrarr X$
be a closed immersion
of complex manifolds.
Let $(\nbigm,L)$ be a coherent $\nbigr_Y$-module.
Let $g$ be a holomorphic function on $X$.
\begin{lem}
\label{lem;10.10.4.2}
$(\nbigm,L)$ is filtered strictly specializable 
(resp. admissibly) along $g_Y:=g_{|Y}$,
if and only if 
$\iota_{\dagger}(\nbigm,L)$ 
is filtered strictly (resp. admissibly)
specializable along $g$.
\end{lem}
\pf
Assume that 
$(\nbigm,L)$ is filtered strictly specializable 
along $g_Y$.
Because 
$\iota_{\dagger}\psitilde_{g_Y,\gminia,u}(L_j\nbigm)$
is strict,
according to the compatibility of the push-forward
and the strict specializability,
$\iota_{\dagger}L_j\nbigm$ is 
strictly specializable along $g$
with any ramified exponential twist,
and $\psitilde_{g,\gminia,u}\bigl(
 \iota_{\dagger}L_j\nbigm\bigr)
\simeq
 \iota_{\dagger}\psitilde_{g_Y,\gminia,u}(L_j\nbigm)$.
Because the cokernel of
$\psitilde_{g_Y,\gminia,u}(L_j\nbigm)
\lrarr
 \psitilde_{g_Y,\gminia,u}(\nbigm)$ are strict,
the cokernel of
$\psitilde_{g,\gminia,u}(\iota_{\dagger}L_j\nbigm)
\lrarr
 \psitilde_{g,\gminia,u}(\iota_{\dagger}\nbigm)$ 
are strict.
A similar claim holds for $\phi_g$.
Hence, 
$\iota_{\dagger}(\nbigm,L)$ is filtered
strictly specializable along $g$.

Suppose that $\iota_{\dagger}(\nbigm,L)$ is
filtered strictly specializable along $g$.
Each $\iota_{\dagger}(L_j\nbigm)$ is 
strictly specializable along $g$
with any ramified exponential twist.
By a general lemma below,
each $L_j\nbigm$ is strictly specializable along $g_{Y}$,
and hence we have
$\psitilde_{g,\gminia,u}\bigl(\iota_{\dagger}L_j\nbigm\bigr)
\simeq
 \iota_{\dagger}
 \psitilde_{g_Y,\gminia,u}(L_j\nbigm)$.
We have a similar isomorphism for $\phi_g$.
Then, we obtain that
$(\nbigm,L)$ is filtered strictly specializable
along $g$.

If $(\nbigm,L)$ is 
filtered specializable along $g_Y$,
$\iota_{\dagger}\psitilde_{g_Y,\gminia,u}(\nbigm)$
and 
$\psitilde_{g,\gminia,u}(\iota_{\dagger}\nbigm)$
are naturally isomorphic,
which is compatible with the induced filtrations
and the nilpotent maps.
A similar claim holds for $\phi_g$.
Hence, 
$(\nbigm,L)$ is admissibly specializable along $g_Y$
if and only if
$\iota_{\dagger}(\nbigm,L)$
is admissibly specializable along $g$.
\hfill\qed

\begin{lem}
Let $\nbigm$ be a coherent strict $\nbigr_Y$-module.
If $\iota_{\dagger}\nbigm$ is strictly specializable along $g$,
then $\nbigm$ is strictly specializable along $g_Y$.
\end{lem}
\pf
We have only to consider the case $\codim Y=1$.
We have only to consider the issue locally around any point of $Y$.
We will shrink $X$ and $Y$ without mention.
By the graph construction,
we replace $X$ and $Y$
with $X\times\cnum_t$ and $Y\times\cnum_t$,
respectively,
and the functions $g$ and $g_Y$ are replaced with  $t$.
Then, we may assume that 
$g$ is a coordinate function,
and $g^{-1}(0)$ and $Y$ are transversal.
Then, after a bi-holomorphic transform,
we may assume that
$X$ is an open subset of
$Z\times \cnum_s\times\cnum_t$,
and $Y=\{s=0\}$ and $g=t$.
We have the grading
$\iota_{\dagger}\nbigm
=\bigoplus_j\iota_{\ast}\nbigm\,\deldel_s^j$.
Let $\Vzero$ be the $V$-filtration of
$\iota_{\dagger}\nbigm$ at $\lambda_0$.
We have only to prove that
it is compatible with the grading,
i.e., for a given 
$f=\sum_{j=0}^N f_j\deldel_s^j\in \Vzero_a(\iota_{\dagger}\nbigm)$,
each $f_j\,\deldel_s^j\in\Vzero_a(\iota_{\dagger}\nbigm)$.
We use an induction on $N$.
We have only to prove that
$f_N\deldel_s^N\in\Vzero_a(\iota_{\dagger}\nbigm)$.
We have 
$s^Nf\in \Vzero_a(\iota_{\dagger}(\nbigm))$.
We obtain $\lambda^Nf_N\in\Vzero_a(\iota_{\dagger}\nbigm)$.
By using the strictness of
$\psi^{(\lambda_0)}_b$ for any $b\in\real$,
we obtain that 
$f_N\in\Vzero_a$,
and hence
$f_N\deldel_s^N\in\Vzero_a$.
\hfill\qed

\subsection{Filtered $\nbigr$-triples}

In this paper,
a filtered $\nbigr_X$-triple means
an $\nbigr_X$-triple with a locally finite increasing 
exhaustive $\seisuu$-indexed filtration
 in the category of $\nbigr_X$-modules.

\begin{df}
A filtered $\nbigr_X$-triple
$(\nbigt,L)$ is called filtered strictly specializable
along $g$,
if the underlying filtered $\nbigr$-modules
are filtered strictly specializable along $g$.
\hfill\qed
\end{df}
\index{filtered specializable}

\begin{lem}
\mbox{{}}
Assume that $(\nbigt,L)$ is
filtered strictly specializable along $g$.
\begin{itemize}
\item
Each $(L_j\nbigt/L_k\nbigt,L)$
is filtered strictly specializable along $g$.
\item
We have the exact sequence
\[
 0\lrarr \psitilde_{g,\gminia,u}(L_j\nbigt)
\lrarr\psitilde_{g,\gminia,u}\nbigt\lrarr
 \psitilde_{g,\gminia,u}(\nbigt/L_j\nbigt)\lrarr 0
\]
 for $\gminia\in\cnum[t_n^{-1}]$
 and 
 $u\in \real\times\cnum$.
\item
We have the exact sequence
$0\lrarr \phi_{g}(L_j\nbigt)
\lrarr\phi_{g}\nbigt\lrarr
 \phi_{g}(\nbigt/L_j\nbigt)\lrarr 0$.
\end{itemize}
In particular,
we obtain a filtered $\nbigr$-triple
$\psitilde_{g,\gminia,u}(\nbigt,L)$
and $\phi_{g}(\nbigt,L)$
such that $\Gr^L$ are strict.
\hfill\qed
\end{lem}

\begin{df}
Let $(\nbigt,L)$ be a filtered $\nbigr$-triple,
which is filtered strictly specializable along $g$.
It is called admissibly specializable along $g$,
if the underlying $\nbigr_X$-modules
are admissibly specializable along $g$.
\hfill\qed
\end{df}
\index{admissibly specializable}

\begin{lem}
If 
$(\nbigt,L)$ is admissibly specializable along $g$,
$(\psitilde_{g,\gminia,u}(\nbigt,L),\nbign)$
and $(\phi_{g}(\nbigt,L),\nbign)$
have relative monodromy filtrations
$M(\nbign;L)$
in the category of $\nbigr$-triples.
\end{lem}
\pf
It follows from the following general 
Lemma \ref{lem;13.5.10.30}, which is a consequence of 
Deligne's inductive formula (\ref{eq;10.9.27.21}, \ref{eq;10.9.27.22}) 
below.
\hfill\qed

\begin{lem}
\label{lem;13.5.10.30}
Let $\nbigt=(\nbigm_1,\nbigm_2,C)$
be an $\nbigr_X$-triple
with a filtration $L$.
Let $\nbign:\nbigt\lrarr\nbigt\otimes\newTate(-1)$
be a morphism such that
$\nbign:L_j\nbigt\lrarr L_j\nbigt\otimes\newTate(-1)$.
Assume that $\nbigm_i$ $(i=1,2)$
have relative monodromy filtrations 
$\Ltilde=M(N;L)$.
Then, we have a filtration $\Ltilde$ of
$\nbigt$ which is a relative monodromy filtration
$M(\nbign;L)$.
\hfill\qed
\end{lem}

We obtain the following lemma 
from Lemma \ref{lem;13.5.10.31}.
\begin{lem}
\label{lem;13.5.10.32}
Let $(\nbigt_i,L)$ $(i=1,2)$ be admissibly specializable
along $g$.
Let $f:(\nbigt_1,L)\lrarr (\nbigt_2,L)$ be any morphism.
Then, 
the induced morphisms
$\psitilde_{g,\gminia,u}(f)$
and 
$\phi_g(f)$ are compatible
with the filtrations
$M(\nbign;L)$.
\hfill\qed
\end{lem}

\begin{lem}
\label{lem;13.5.9.10}
Suppose that $(\nbigt,L)$ is admissibly specializable along $g$.
Let $A$ be a nowhere vanishing holomorphic function on $X$.
Then, $(\nbigt,L)$ is admissibly specializable along $g_1:=Ag$,
and we have the following natural isomorphisms:
\[
\bigl(
 \Gr^{M(\nbign;L)}
 \psitilde_{g,\gminia,u}(\nbigt),
 \nbign^{(0)}
\bigr)
\simeq
 \bigl(
 \Gr^{M(\nbign;L)}
 \psitilde_{g_1,\gminia,u}(\nbigt),
 \nbign^{(0)}
\bigr)
\]
\[
\bigl(
 \Gr^{M(\nbign;L)}
 \phi_{g}(\nbigt),
 \nbign^{(0)}
\bigr)
\simeq
 \bigl(
 \Gr^{M(\nbign;L)}
 \phi_{g_1}(\nbigt),
 \nbign^{(0)}
\bigr)
\]
They preserve the canonical splitting 
the induced filtration $L$.
\end{lem}
\pf
It follows from
Corollary \ref{cor;11.1.18.1}.
\hfill\qed

\vspace{.1in}
We obtain the following lemma from 
Lemma \ref{lem;10.10.4.2}.

\begin{lem}
Let $\iota:Y\lrarr X$ be a closed embedding of
complex manifolds.
Let $(\nbigt,L)$ be a coherent strict filtered $\nbigr_Y$-triple.
Let $g$ be a holomorphic function on $X$
such that $Y\not\subset g^{-1}(0)$.
Then, $(\nbigm,L)$ is filtered strictly specializable along $g_{|Y}$,
if and only if
$\iota_{\dagger}(\nbigm,L)$ is filtered strictly specializable
along $g$.
\hfill\qed
\end{lem}

\begin{df}
Let $(\nbigt,L)$ be a filtered $\nbigr_X$-triple.
It is called admissibly (resp. filtered strictly) specializable,
if the following holds
\begin{itemize}
\item
 Let $U\subset X$ be any open subset
with a holomorphic function $g$.
 Then, $(\nbigt,L)_{|U}$ is 
 admissibly 
 (resp. filtered strictly) specializable along $g$.
\hfill\qed 
\end{itemize}
\end{df}
\index{filtered strictly specializable}
\index{admissibly specializable}

\subsection{Admissible specializability along hypersurfaces}

Let $H$ be an effective divisor of $X$.
\begin{df}
\label{subsection;13.5.9.12}
A filtered $\nbigr_X$-triple
$(\nbigt,L)$ is called admissibly 
(resp. filtered strictly) specializable along $H$,
if the following holds:
\begin{itemize}
\item
 Let $U$ be any open subset of $X$
 with a generator $g_U$ of $\nbigo(-H)_{|U}$.
 Then, $(\nbigt,L)_{|U}$ is 
 admissibly (resp. filtered strictly) specializable
 along $g_U$.
\end{itemize}
Admissible specializability
and filtered strict specializability along $H$ 
for $\nbigr_X$-modules
are defined in similar ways.
\hfill\qed
\end{df}
The following lemma is easy to see.
\begin{lem}
If $(\nbigt,L)$ is admissibly (resp. filtered strictly)
specializable along $H$,
each $L_j\nbigt/L_k\nbigt$ is 
admissibly (resp. filtered strictly)
specializable along $H$.
\hfill\qed
\end{lem}

\subsubsection{Induced graded $\nbigr$-triples}
\label{subsection;13.5.9.14}

Let $U\subset X$ be open with a generator $g_U$
of $\nbigo(-H)_{|U}$.
We set $W:=M(\nbign;L)[-2a]$ 
on $\phi^{(a)}_{g_U}(\nbigt_{|U})$,
and 
$W:=M(\nbign;L)[-2a+1]$
on $\psi^{(a)}_{g_U}(\nbigt_{|U})$.
We have the $\seisuu$-graded $\nbigr$-triples 
\[
\Gr^{W}
 \psi^{(a)}_{g_U}(\nbigt_{|U}),
\quad
\Gr^{W}
 \phi^{(a)}_{g_U}(\nbigt_{|U})
\]
on $U$
with the endomorphisms
\[
 \nbign:\Gr^{W}_k\psi^{(a)}_{g_U}(\nbigt_{|U})
\lrarr
\Gr^{W}_{k-2}\psi^{(a)}_{g_U}(\nbigt_{|U})
 \otimes\newTate(-1),
\]
\[
 \nbign:\Gr^{W}_k\phi^{(a)}_{g_U}(\nbigt_{|U})
\lrarr
\Gr^{W}_{k-2}\phi^{(a)}_{g_U}(\nbigt_{|U})
 \otimes\newTate(-1).
\]
According to Lemma \ref{lem;13.5.9.10},
they are independent of the choice of 
a generator $g_U$.
Hence, by gluing them for varied $(U,g_U)$,
we obtain graded $\nbigr_X$-triples
with endomorphisms,
denoted by
\[
\Bigl(
 \bigl(
\Gr^{W}\psi^{(a)}_H
 \bigr)(\nbigt),
 \nbign
\Bigr),
\quad
\Bigl(
 \bigl(
\Gr^{W}\phi^{(a)}_H
 \bigr)(\nbigt),
 \nbign
\Bigr).
\]
The filtration $L$
induces filtrations on
$\bigl(
\Gr^{W}\psi^{(a)}_H
 \bigr)(\nbigt)$
and
$\bigl(
\Gr^{W}\phi^{(a)}_H
 \bigr)(\nbigt)$,
which are preserved by $\nbign$.
By gluing the canonical decomposition of Kashiwara
(see \S\ref{subsection;10.12.28.1} below),
we obtain the following isomorphisms,
which are splittings of $L$:
\begin{equation}
 \label{eq;13.5.10.100}
 \bigl(
\Gr^{W}\psi^{(a)}_H
 \bigr)(\nbigt)
\simeq
 \bigl(
\Gr^{W}\psi^{(a)}_H
 \bigr)(\Gr^L\nbigt)
\end{equation}
\begin{equation}
 \label{eq;13.5.10.101}
 \bigl(
\Gr^{W}\phi^{(a)}_H
 \bigr)(\nbigt)
\simeq
 \bigl(
\Gr^{W}\phi^{(a)}_H
 \bigr)(\Gr^L\nbigt)
\end{equation}
\begin{rem}
$(\Gr^W\psi^{(a)}_H)(\nbigt)$
and 
$(\Gr^W\psi^{(a)}_H)(\Gr^L\nbigt)$
depend only on
$\nbigt(\ast H)$
and 
$\Gr^L(\nbigt(\ast H))$.
\hfill\qed
\end{rem}

By the condition of the admissible specializability,
we have the following induced morphisms:
\[
 \bigl(
 \Gr^W\psi^{(a+1)}_H
 \bigr)
 (\nbigt)
\stackrel{\alpha}{\lrarr}
  \bigl(
 \Gr^W\phi^{(a)}_H
 \bigr)
 (\nbigt)
\stackrel{\beta}{\lrarr}
  \bigl(
 \Gr^W\psi^{(a)}_H
 \bigr)
 (\nbigt)
\]
The composites $\beta\circ\alpha$
and $\alpha\circ\beta$
are identified with
the morphisms $\nbign$ on 
$(\Gr^W\!\psi^{(a)}_H)(\nbigt)$
and 
$(\Gr^W\!\phi^{(a)}_H)(\nbigt)$
up to constants.

\vspace{.1in}

Let $(\nbigt_i,L)$ $(i=1,2)$ be admissibly specializable
along $H$.
Let $f:(\nbigt_1,L)\lrarr (\nbigt_2,L)$ be any morphism.
According to Lemma \ref{lem;13.5.10.32},
we have the naturally induced morphisms
\[
\bigl(
\Gr^{W}\psi^{(a)}_H
\bigr)(f):
  \bigl(
\Gr^{W}\psi^{(a)}_H
 \bigr)(\nbigt_1)
\lrarr
  \bigl(
\Gr^{W}\psi^{(a)}_H
 \bigr)(\nbigt_2) 
\]
\[
\bigl(
\Gr^{W}\phi^{(a)}_H
\bigr)(f):
 \bigl(
\Gr^{W}\phi^{(a)}_H
 \bigr)(\nbigt_1)
\lrarr
  \bigl(
\Gr^{W}\phi^{(a)}_H
 \bigr)(\nbigt_2) 
\]

\begin{lem}
If $(\nbigt,L)$ is integrable,
the induced $\nbigr_X$-triples
$ \bigl(
\Gr^{W}\psi^{(a)}_H
 \bigr)(\nbigt)$
and 
$ \bigl(
\Gr^{W}\phi^{(a)}_H
 \bigr)(\nbigt)$
are also integrable.
\hfill\qed
\end{lem}

\chapter{Gluing of good-KMS $\nbigr$-triples}
\label{section;11.4.9.10}

We shall study 
a nice class of smooth $\nbigr$-triples with normal crossing
singularity.
\S\ref{subsection;11.4.3.10}--\ref{subsection;11.4.3.11}
are mainly preparation for \S\ref{section;11.4.3.4},
and \S\ref{subsection;11.4.3.12}--\ref{subsection;10.10.1.10}
are mainly preparation for \S\ref{section;11.4.3.5}.

\section{Smooth good-KMS $\nbigr$-modules}
\label{subsection;11.4.3.10}

\subsection{Good-KMS meromorphic prolongment}

Let $X$ be a complex manifold
with a simply normal crossing hypersurface $D$.
Let $\nbigm$ be a smooth $\nbigr_{X(\ast D)}$-module.
The natural family of flat $\lambda$-connections on $\nbigm$
is denoted by $\DD$.
In this paper,
we say that $\nbigm$ is unramifiedly good,
if the following holds:
\index{unramifiedly good}
\begin{itemize}
\item
For each point $P\in D$,
there exists a good set of irregular values\
$\Irr(\nbigm,P)
\subset
 \nbigo_X(\ast D)_P/\nbigo_{X,P}$,
such that the formal completion of $\nbigm$
at $(\lambda,P)$ is isomorphic to a direct sum
\begin{equation}
\label{eq;13.5.5.21}
 \bigoplus_{\gminia\in\Irr(\nbigm,P)}
 (\nbigm_{\gminia},\DD_{\gminia})
\end{equation}
such that
$(\nbigm_{\gminia},\DD_{\gminia}-d\gminiatilde)$
is regular, 
i.e., it has a logarithmic lattice.
Here, $\gminiatilde\in\nbigo_X(\ast D)_P$ is 
a lift of $\gminia$.
(See \S2.1 of \cite{mochi7} for good set of irregular values.)
The decomposition (\ref{eq;13.5.5.21})
is called the Hukuhara-Levelt-Turrittin type decomposition.
\end{itemize}
It is called unramifiedly good-KMS,
if moreover it has the KMS-structure
(see \S2.8 of \cite{mochi7}).
\index{unramifiedly good-KMS}
It is called good-KMS (resp. good),
if it is locally the descent of
an unramifiedly good-KMS (resp. good)
smooth $\nbigr_{X(\ast D)}$-module.
\index{good-KMS}\index{good smooth}
It is called regular-KMS
if moreover $\nbigm$ is a regular along $D$.
\index{regular-KMS}
For a decomposition $D=D^{(1)}\cup D^{(2)}$,
a good-KMS smooth $\nbigr_{X(\ast D)}$-module
$\nbigm$ is called regular along $D^{(1)}$,
if $\nbigm_{\nbigx\setminus \nbigd^{(2)}}$
is regular-KMS.

Let $D=\bigcup_{i\in\Lambda}D_i$
be the irreducible decomposition.
In the following,
for $I\subset \Lambda$,
we set $D_I:=\bigcap_{i\in I}D_i$,
$D(I):=\bigcup_{i\in I}D_i$
and $\del D_I:=\bigcup_{j\not\in I}(D_I\cap D_j)$.
We put $D_I^{\circ}:=D_I\setminus \del D_I$.
We set $I^c:=\Lambda\setminus I$.
If $I=\emptyset$,
we put $D_{\emptyset}:=X$.

\subsection{Induced bundles
on the intersection of divisors}
\label{subsection;13.5.4.20}

Let $X$ and $D$ be as above.
Let $\nbigm$ be a good-KMS $\nbigr_{X(\ast D)}$-module.
For each $\lambda_0\in\cnum_{\lambda}$,
we have the associated family of 
good filtered $\lambda$-flat bundles
$\nbigqzero_{\ast}\nbigm$
on $(\nbigxzero,\nbigdzero)$
indexed by $\real^{\Lambda}$.
We have the induced filtration
$\lefttop{i}\Fzero$ of
$\nbigqzero_{\veca}\nbigm_{|\nbigdzero_i}$
given by
\[
 \lefttop{i}\Fzero_b\bigl(
 \nbigqzero_{\veca}\nbigm_{|\nbigdzero_i}
 \bigr):=
 \Image\Bigl(
 \nbigqzero_{\veca'}\nbigm_{|\nbigdzero_i}
\lrarr
 \nbigqzero_{\veca}\nbigm_{|\nbigdzero_i}
 \Bigr), 
\]
where the $i$-th component of $\veca'$
is $b$, and the other components are equal
to those of $\veca$.
\index{filtration $\lefttop{i}\Fzero_b$}
The induced filtrations $\lefttop{i}\Fzero$ of
$\nbigqzero_{\veca}\nbigm_{|\nbigdzero_I}$ 
$(i\in I)$ are compatible.
For $\vecb\in\real^I$,
we put
$\lefttop{I}\Fzero_{\vecb}:=
 \bigcap_{i\in I}\lefttop{i}\Fzero_{b_i}$
on $\nbigdzero_I$,
and
\[
 \lefttop{I}\Gr^{\Fzero}_{\vecb}
 \nbigqzero_{\veca}\nbigm:=
 \lefttop{I}\Fzero_{\vecb}\big/
 \sum_{\vecc\lneq\vecb}
  \lefttop{I}\Fzero_{\vecc}.
\]
\index{filtration $\lefttop{I}\Fzero_b$}
\index{sheaf $\lefttop{I}\Gr^{\Fzero}_{\vecb}
 \nbigqzero_{\veca}\nbigm$}
On 
$\lefttop{I}\Gr^{\Fzero}_{\vecb}
 \nbigqzero_{\veca}\nbigm$,
we have the induced endomorphisms
$\Res_i(\DD)$ $(i\in I)$.
We have the decomposition
\[
 \lefttop{I}\Gr^{\Fzero}_{\vecb}
 \nbigqzero_{\veca}\nbigm
=\bigoplus_{\substack{
 \vecu\in(\real\times\cnum)^I\\
 \paramap(\lambda_0,\vecu)=\vecb
 }}
  \lefttop{I}\nbiggzero_{\vecu}
  \nbigqzero_{\veca}\nbigm,
\]
such that
(i) it is preserved by $\Res_i(\DD)$ $(i\in I)$,
(ii) the restriction of 
 $\Res_i(\DD)-\eigenmap(\lambda,u_i)$ to
 $\lefttop{I}\nbiggzero_{\vecu}
 \nbigqzero_{\veca}\nbigm$
 is nilpotent,
where $u_i$ is the $i$-th component of $\vecu$.
(Recall that,
for $u=(a,\alpha)$,
$\paramap(\lambda,u):=
 a-\Re(\lambda\alphabar)\in\real$
and 
$\eigenmap(\lambda,u):=
 \alpha-a\lambda-\alphabar\lambda^2\in\cnum$.)

For $I\subset\Lambda$
and $\vecb\in\real^I$,
we take $\veca\in\real^{\Lambda}$
which is mapped to $\vecb$
by the projection $\real^{\Lambda}\lrarr\real^I$,
and we put
$\lefttop{I}\nbigqzero_{\vecb}\nbigm:=
 \nbigqzero_{\veca}\nbigm
 \otimes\nbigo\bigl(\ast \nbigd(I^c)\bigr)$.
It is independent of the choice of $\veca$.
We obtain
$\lefttop{I}\nbigg_{\vecu}\nbigq_{\vecb}\nbigm$
as above.
\index{sheaf $\lefttop{I}\nbigqzero_{\vecb}\nbigm$}
In particular, for $\vecu\in (\real\times\cnum)^I$,
we obtain
\[
 \lefttop{I}\nbiggzero_{\vecu}\nbigm:=
 \lefttop{I}\nbiggzero_{\vecu}
 \lefttop{I}\nbigqzero_{\paramap(\lambda_0,\vecu)}\nbigm.
\]
If $\lambda_1$ is sufficiently close to $\lambda_0$
such that 
$\nbigx^{(\lambda_1)}
\subset\nbigxzero$,
we have
$\lefttop{I}\,\nbiggzero_{\vecu}
(\nbigm)_{|\nbigd_I^{(\lambda_1)}}
=\lefttop{I}\nbigg^{(\lambda_1)}_{\vecu}
 \bigl(\nbigm\bigr)$.
Hence, we can glue
$\lefttop{I}\nbiggzero_{\vecu}\nbigm$
for varied $\lambda_0\in\cnum_{\lambda}$,
and we obtain an 
$\nbigo_{\nbigd_I}(\ast \del\nbigd_I)$-module
$\lefttop{I}\nbigg_{\vecu}\nbigm$
on $\nbigd_I$.
\index{sheaf $\lefttop{I}\nbigg_{\vecu}\nbigm$}

\subsection{Hukuhara-Levelt-Turrittin type decomposition}
\label{subsection;13.5.4.30}

Let us consider the case $X=\Delta^n$
and $D=\bigcup_{i=1}^{\ell}\{z_i=0\}$.
We set $D_i=\{z_i=0\}$.
For $J\subset \ellsitabar:=\{1,\ldots,\ell\}$,
let $M(X,D(J))$ be the set of meromorphic functions on $X$
which may have poles at most along $D(J)$.
Let $H(X)$ be the space of holomorphic functions on $X$.
Let $\nbigm$ be 
an unramifiedly good-KMS $\nbigr_{X(\ast D)}$-module.
After shrinking $X$,
we have the good set of irregular values
$\Irr(\nbigm)$ in $M(X,D)/H(X)$.
For $I\subset\Lambda$,
let $\Irr(\nbigm,I)$ denote 
the image of $\Irr(\nbigm)$
by $M(X,D)/H(X)\lrarr M(X,D)/M(X,D(I^c))$.
We have the Hukuhara-Levelt-Turrittin type decomposition
(see \S2.4.2 of \cite{mochi7}):
\[
 (\nbigm,\DD)_{|\nbigdhat_I}
=\bigoplus_{\gminia\in\Irr(\nbigm,I)}
 \bigl(
 \nbigm_{\gminia,\nbigdhat_I},\DDhat_{\gminia}
 \bigr)
\]
On $\nbigxzero$,
the decomposition is compatible with the KMS-structure,
i.e., for $\veca\in\real^{\ell}$,
we have the induced decomposition
\[
 (\nbigqzero_{\veca}
 \nbigm,\DD)_{|\nbigdhatzero_I}
=\bigoplus_{\gminia\in\Irr(\nbigm,I)}
 \bigl(
 \nbigqzero_{\veca}
 \nbigm_{\gminia,\nbigdhat_I},\DDhat_{\gminia}
 \bigr),
\]
and 
$\DDhat_{\gminia}-d\gminia$ are logarithmic
along $\nbigd(I)$
in the following sense:
{\small
\[
\bigl(\DDhat_{\gminia}-d\gminia\bigr)
 \nbigqzero_{\veca}\nbigm_{\gminia,
 \nbigdhat_I}
\subset
 \nbigqzero_{\veca}\nbigm_{\gminia,
 \nbigdhat_I}
\otimes
\Bigl(
 \Omega^1_{\nbigxzero/\cnum_{\lambda}}
 \bigl(\log \nbigdhatzero(I)\bigr)
+\Omega^1_{\nbigxzero/\cnum_{\lambda}}
 \bigl(\ast\nbigdhatzero(I^c)\bigr)
\Bigr)
\]
}
We set 
$\nbigm_{\nbigdhat_I}^{(\reg)}:=
 \nbigm_{0,\nbigdhat_I}$,
and 
$\nbigm_{\nbigdhat_I}^{(\irr)}:=
 \bigoplus_{\gminia\neq 0}
 \nbigm_{\gminia,\nbigdhat_I}$.

Let us consider the case that $\nbigm$ is good-KMS,
which is not necessarily unramified.
Let $\varphi:(X',D')\lrarr (X,D)$ be a ramified covering
such that 
$\nbigm'=
 \varphi^{\ast}\nbigm$ is unramified.
Applying the above construction,
we obtain the decomposition
$\nbigm'_{|\nbigdhat'_I}=
 \nbigm^{\prime(\reg)}_{\nbigdhat'_I}
\oplus
 \nbigm^{\prime(\irr)}_{\nbigdhat'_I}$,
which is preserved by the action of
the Galois group of the ramified covering.
We obtain a decomposition
$\nbigm_{|\nbigdhat_I}
=\nbigm^{(\reg)}_{\nbigdhat_I}
\oplus
 \nbigm^{(\irr)}_{\nbigdhat_I}$.
It is independent of the choice of
a ramified covering
and a choice of the coordinate.
It is compatible with the KMS-structure,
i.e,
we have the induced decomposition
$\nbigqzero_{\ast}\nbigm_{|\nbigdhat_I}
=\nbigqzero_{\ast}\nbigm^{(\reg)}_{\nbigdhat_I}
\oplus
 \nbigqzero_{\ast}\nbigm^{(\irr)}_{\nbigdhat_I}$.

\subsection{Specialization}
\label{subsection;13.5.4.21}

Let $X$ and $D$ be as in 
\S\ref{subsection;13.5.4.30}.
Let $\nbigm$ be a good-KMS
$\nbigr_{X(\ast D)}$-module.
We set 
$\vecdelta_I:=(\vecdelta,\ldots,\vecdelta)\in
 (\real\times\cnum)^I$.
For $\vecu\in(\real\times\cnum)^I$,
we put 
$(\vecb,\vecbeta):=
 \kmsmap(\lambda_0,\vecu+\vecdelta_I)$,
and we obtain the following on $\nbigd_I$:
\[
  \lefttop{I}\,\psitilde_{\vecu}
 \bigl(\nbigm\bigr)
:=
 \lefttop{I}\nbigg_{\vecu+\vecdelta_I}
 \nbigm^{(\reg)}_{\nbigdhat_I}
\]
\index{sheaf $\lefttop{I}\psitilde_{\vecu}(\nbigm)$}
It is naturally a good-KMS smooth
$\nbigr_{D_I(\ast \del D_I)}$-module,
although it depends on the choice of
a coordinate system.
If $\nbigm$ is unramifiedly good,
then 
$\Irr(\lefttop{I}\psitilde_{\vecu}\nbigm)$
is 
$\bigl\{
 \gminia_{|D_I}\,\big|\,
 \gminia\in
 \varphi_I^{-1}(0)
 \bigr\}$,
where
$\varphi_I:\Irr(\nbigm)\lrarr
 M(X,D)/M(X,\nbigd(I^c))$
is the natural map.
The bundle 
$\lefttop{I}\,\psitilde_{\vecu}
 \bigl(\nbigm\bigr)$
is also equipped with
the induced endomorphisms
$\Res_i(\DD)$ $(i\in I)$,
which are independent of the choice of 
a coordinate system.
By construction,
for $I_1\sqcup I_2\subset\ellsitabar$,
we have a natural isomorphism
$\lefttop{I_1}\psitilde_{\vecu_1}
 \Bigl(
 \lefttop{I_2}\psitilde_{\vecu_2}
 (\nbigm)
 \Bigr)
\simeq
 \lefttop{I_1\sqcup I_2}
 \psitilde_{(\vecu_1,\vecu_2)}
 (\nbigm)$.

It is easy to see
that $\nbigm$ is strictly specializable along $z_i$
as an $\nbigr_{X(\ast D)}$-module
(\S\ref{subsection;13.5.4.31}).
We have
$\lefttop{i}\psitilde_u(\nbigm)
\simeq
 \psitilde_{z_i,u}(\nbigm)$.
So, 
$\lefttop{I}\psitilde_{\vecu}$
is obtained
as the composition 
$\psitilde_{z_{i_1},u_i}\circ\cdots
 \circ\psitilde_{z_{i_m},u_{i_m}}$.

\subsection{Reduction with respect to Stokes structure}
\label{subsection;13.5.5.30}

Let $X$, $D$ and $\nbigm$ be as in 
\S\ref{subsection;13.5.4.21}.
We recall the procedure of the reduction 
of $\nbigm$ with respect to
the Stokes structure.
We will shrink $X$ around the origin
without mention.

If $\nbigm$ is unramifiedly good-KMS,
we have the reduction with respect to 
the Stokes structure along $D(I)$,
and we obtain a graded good-KMS 
smooth $\nbigr_{X(\ast D)}$-module
\[
 \lefttop{I}\Gr^{\St}(\nbigm)
=\bigoplus_{\gminia\in\Irr(\nbigm,I)}
 \lefttop{I}\Gr^{\St}_{\gminia}(\nbigm).
\]
Indeed, we can choose an auxiliary sequence
$\vecm(0),\vecm(1),\ldots,\vecm(N)$
for $\Irr(\nbigm)$ as in \S2.1.2 of \cite{mochi7}.
By successive use of the reduction
with respect to the Stokes structure in the level $\vecm(j)$
in \S3.3 of \cite{mochi7},
we obtain a sequence of 
unramifiedly good-KMS $\nbigr_{X(\ast D)}$-modules
$\Gr^{\vecm(j)}(\nbigm)$.
Let $j_0$ be the minimal among 
$j$ such that
the $p$-th $(p\in I)$ components of $\vecm(j)$ are $0$.
We define
$\lefttop{I}\Gr^{\St}(\nbigm):=
 \Gr^{\vecm(j_0-1)}(\nbigm)$.

For each $\gminia$,
take an appropriate lift to $M(X,D)$
which is also denoted by $\gminia$,
then the tensor product
$\lefttop{I}\Gr^{\St}_{\gminia}(\nbigm)
 \otimes L(-\gminia)$ is regular along $D(I)$.
We have a natural isomorphism
$\lefttop{I}\Gr^{\St}(\nbigm)_{|\nbigdhat_I}
\simeq \nbigm_{|\nbigdhat_I}$.
We put
$\lefttop{I}\Gr^{\St,\reg}(\nbigm):=
 \lefttop{I}\Gr^{\St}_0(\nbigm)$
and 
$\lefttop{I}\Gr^{\St,\irr}(\nbigm):=
 \bigoplus_{\gminia\neq 0}
 \lefttop{I}\Gr^{\St}_{\gminia}(\nbigm)$.
We have
$\lefttop{I}\psitilde_{\vecu}\Bigl(
\lefttop{I}\Gr^{\St}(\nbigm)
\Bigr)
\simeq
 \lefttop{I}\psitilde_{\vecu}(\nbigm)$
naturally.

Suppose that $\nbigm$ is good-KMS.
We take a ramified covering
$\varphi:(X',D')\lrarr (X,D)$
such that
$\nbigm':=\varphi^{\ast}\nbigm$
is unramified.
By applying the above procedure,
we obtain the reduction
$\lefttop{I}\Gr^{\St}(\nbigm')
=\lefttop{I}\Gr^{\St,\reg}(\nbigm')
\oplus
 \lefttop{I}\Gr^{\St,\irr}(\nbigm')$,
on which the Galois group of the covering
naturally acts.
As the descent,
we obtain 
a good-KMS smooth $\nbigr_{X(\ast D)}$-module
$\lefttop{I}\Gr^{\St}(\nbigm)
=\lefttop{I}\Gr^{\St,\reg}(\nbigm)
\oplus
 \lefttop{I}\Gr^{\St,\irr}(\nbigm)$
on $(\nbigx,\nbigd)$.
We have natural isomorphisms
\[
 \lefttop{I}\Gr^{\St,\reg}(\nbigm)_{|\nbigdhat_I}
\simeq
 \nbigm^{(\reg)}_{|\nbigdhat_I},
\quad
 \lefttop{I}\Gr^{\St,\irr}(\nbigm)_{|\nbigdhat_I}
\simeq
 \nbigm^{(\irr)}_{|\nbigdhat_I}.
\]
We also have 
$\lefttop{I}\psitilde_{\vecu}\Bigl(
\lefttop{I}\Gr^{\St}(\nbigm)
\Bigr)
\simeq
 \lefttop{I}\psitilde_{\vecu}(\nbigm)$
naturally.

\section{Compatibility of filtrations}
\label{subsection;11.4.3.11}

\subsection{Compatibility with 
Hukuhara-Levelt-Turrittin type decomposition}

Let $X$ be a complex manifold,
and let $D$ be a simple normal crossing
hypersurface.
Let $\nbigm$ be 
a good smooth $\nbigr_{X(\ast D)}$-module.
Let $L$ be a filtration of $\nbigm$
in the category of smooth $\nbigr_{X(\ast D)}$-modules,
i.e.,
a filtration
in the category of $\nbigr_{X(\ast D)}$-modules
such that
$\Gr^L(\nbigm)$ is also a smooth $\nbigr_{X(\ast D)}$-module.

\begin{lem}
Let $P\in D$.
Suppose that $\nbigm$ is unramified around $P$,
Then, the filtration $L$ is compatible with
the formal decomposition of $\nbigm$
at $(\lambda,P)\in \nbigd$
as in {\rm(\ref{eq;13.5.5.21})}:
\begin{equation}
 \label{eq;13.5.5.20}
 \nbigm_{|\widehat{(\lambda,P)}}
=\bigoplus_{\gminia\in\Irr(\nbigm,P)}
 \nbigmhat_{k,\gminia,(\lambda,P)}.
\end{equation}
\end{lem}
\pf
Recall a standard and easy result for 
$\cnum(\!(t)\!)$-differential module.
Let $M$ be a $\cnum(\!(t)\!)$-differential module
of finite rank with a decomposition
$M=\bigoplus_{\gminia\in t^{-1}\cnum[t^{-1}]} 
M_{\gminia}$
such that $\del_t-\del_t(\gminia)$ are regular singular.
Then, any differential submodule
$M'\subset M$ is compatible with the decomposition.

Similarly,
let $M$ be a $\cnum(\!(t)\!)$-module of finite rank
with an $\cnum(\!(t)\!)$-endomorphism $f$,
with a decomposition
$(M,f)=\bigoplus_{\gminia\in t^{-1}\cnum[t^{-1}]}
(M_{\gminia},f_{\gminia})$
such that
$M_{\gminia}$ has a lattice $L_{\gminia}$
preserved by $f_{\gminia}-\del_t\gminia$.
If $M'\subset M$ is preserved by $f$,
then $M'$ is compatible with the decomposition,
which can be proved in a similar way.

Let $\varphi:\Delta\lrarr X$ be any morphism
such that $\varphi(0)=P$
and that $\varphi(C)\not\subset D$.
If $\varphi$ is general,
the induced map
$\Irr(\nbigm,P)\lrarr
 \nbigo_{\Delta}(\ast 0)_0/\nbigo_{\Delta,0}$
is injective.
By applying the above result,
we obtain that
$\varphi^{\ast}L$ is compatible with
the decomposition obtained
as the pull back of (\ref{eq;13.5.5.20}).
It follows that we obtain that
$L$ is compatible with 
(\ref{eq;13.5.5.20}).
\hfill\qed

\vspace{.1in}

We obtain the following lemma 
from the previous lemma.

\begin{lem}
Each $L_j\nbigm$ is
a good smooth $\nbigr_{X(\ast D)}$-module.
\hfill\qed
\end{lem}

\begin{rem}
The Stokes structure of 
$\nbigm$ is compatible with
the filtration $L$.
It follows from the characterization of
the Stokes filtration
in terms of the growth order of flat sections.
(See {\rm\S3.2.3 \cite{mochi7}},
for example.)
\hfill\qed
\end{rem}

\subsection{Extension of smooth good-KMS $\nbigr$-modules}

Let $\nbigk\subset\cnum_{\lambda}$
be a small neighbourhood of $\lambda_0$.
Let $X$ be a complex manifold
with a normal crossing hypersurface $D$.
We set $(\nbigx,\nbigd):=\nbigk\times(X,D)$.
Let $\nbigm$ be 
a good smooth $\nbigr_{X(\ast D)}$-module.
Let $L$ be a filtration of $\nbigm$
in the category of $\nbigr_{X(\ast D)}$-modules.

\begin{prop}
\label{prop;10.10.9.2}
Assume that
(i) $\Gr^L\nbigm$ is good-KMS,
(ii) for any smooth point $P$ of $D$,
 the restriction of $\nbigm$
 to a small neighbourhood of $P$ is good-KMS.
Then, $\nbigm$ is good-KMS.
\end{prop}
\pf
We may assume that
$X=\Delta^n$,
$D_i=\{z_i=0\}$
and $D=\bigcup_{i=1}^{\ell}D_i$.
We set $D_{[2]}:=\bigcup_{i\neq j} D_i\cap D_j$.
We may also assume that
$\nbigm$ is unramifiedly good.
For any $\lambda_0$,
by the assumption,
we have the lattices
$\nbigq_{\veca}^{(\lambda_0)\prime}\nbigm$ $(\veca\in\real^{\ell})$
of
$\nbigm_{|\nbigxzero\setminus \nbigdzero_{[2]}}$.

Let us consider the case
that $\lambda_0$ is generic
with respect to $\Gr^L\nbigm$.
The reduction $\Gr^{\St}(\nbigm)$
is also equipped with the induced filtration $L$,
and $\Gr^L\Gr^{\St}(\nbigm)$ is 
unramifiedly good-KMS.
Moreover, for any smooth point $P$ of $D$,
the restriction of 
$\Gr^{\St}\nbigm$ on a neighbourhood of $P$
has KMS structure.
Because $\lambda_0$ is generic,
$\Gr^{\St}\nbigm$ has KMS structure
at $\lambda_0$.
Indeed,
$\Gr^{\St}(\nbigm)$
has a decomposition
$\bigoplus \Gr^{\St}_{\gminia}(\nbigm)$,
and $\Gr^{\St}_{\gminia}(\nbigm)\otimes\nbigl(-\gminia)$
can be regarded as a family of 
regular singular meromorphic flat bundles.
Hence, it is easy to construct the family of lattices
with the desired property in a direct way.
By applying the Riemann-Hilbert-Birkhoff 
correspondence 
(\S4 of \cite{mochi7})
with the lattices 
$\nbigqzero_{\veca}\Gr^{\St}\nbigm$,
we obtain a lattice 
$\nbigqzero_{\veca}\nbigm$
such that
(i) it induces $\nbigqzero_{\veca}
 \Gr^{\St}\nbigm$,
(ii) $\nbigqzero_{\veca}
 \nbigm_{|\nbigx\setminus\nbigdzero_{[2]}}
\simeq
 \nbigq^{(\lambda_0)\prime}_{\veca}
 \nbigm$.
The second condition implies
$\Gr^L\nbigqzero_{\veca}\nbigm
\simeq
 \nbigqzero_{\veca}\Gr^L\nbigm$.
Then, it is easy to see that
$\nbigm$ is unramifiedly good-KMS
around $\lambda_0$.

Let us consider the case that
$\lambda_0\neq 0$ is not necessarily generic.
By the assumption and the consideration
in the generic case,
we obtain that
the restriction of $\nbigm$ to
$\nbigx\setminus\nbigd^{\lambda_0}_{[2]}$
is good-KMS,
i.e.,
we have vector bundles
$\nbigq^{(\lambda_0)\prime}_{\veca}\nbigm$
on $\nbigx\setminus \nbigd^{\lambda_0}_{[2]}$,
which induces
$\nbigqzero_{\veca}\Gr^L
 \nbigm_{|\nbigx\setminus
 \nbigd^{\lambda_0}_{[2]}}$.
We have only to prove that
$\nbigq^{(\lambda_0)\prime}_{\veca}\nbigm$
is naturally extended to 
vector bundles on $\nbigx$.
Because $\nbigd^{\lambda_0}_{[2]}$
is of codimension 3 in $\nbigx$,
it follows from Lemma \ref{lem;10.9.28.11} below.
The case $\lambda_0=0$ can be argued similarly.
Thus, we obtain Proposition \ref{prop;10.10.9.2}.
\hfill\qed

\begin{lem}
\label{lem;10.9.28.11}
Let $Y$ be a complex manifold.
Let $Z\subset Y$  be an analytic closed subset
of codimension $3$.
Let $E_i$ be vector bundles on $Y$.
Let $E_0'$ be a vector bundle on $Y-Z$
with an exact sequence
$0\lrarr E_{1|Y\setminus Z}
\lrarr E'_0\lrarr E_{2|Y\setminus Z}\lrarr 0$.
Then, $E'_0$ is naturally extended
to a vector bundle on $Y$
with an exact sequence
$0\lrarr E_{1}
\lrarr E_0\lrarr E_{2}\lrarr 0$.
\end{lem}
\pf
Let $j:Y\setminus Z\lrarr Y$ be 
the open immersion.
The claim of the lemma follows from
$R^1j_{\ast}\nbigo_{Y\setminus Z}=0$,
which holds because $\codim Z\geq 3$.
(See Lemma 5 of \cite{Douady-Bourbaki},
for example.)
\hfill\qed

\subsection{Compatibility with KMS structure}

\begin{df}
Let $\nbigm$ be a good-KMS $\nbigr_{X(\ast D)}$-module.
We say that a filtration $L$ of $\nbigm$
is compatible with the KMS-structure,
if the induced increasing sequence
$\Gr^L(\nbigqzero_{\ast}\nbigm)$
gives a $KMS$-structure 
of $\Gr^L(\nbigm)$.
\hfill\qed
\end{df}
We will give an example of a filtration below,
which is not compatible with a KMS-structure.

We obtain the following lemma
from Proposition \ref{prop;10.10.9.2}.
\begin{lem}
\label{lem;10.9.28.12}
Let $\nbigm$ be a good $\nbigr_{X(\ast D)}$-module.
Let $L$ be a filtration of $\nbigm$
in the category of smooth $\nbigr_{X(\ast D)}$-modules.
Assume the following:
\begin{itemize}
\item
$\Gr^L(\nbigm)$ has a $KMS$-structure.
\item
 For any smooth point $P$ of $D$,
 there exists a neighbourhood $X_P$ of $P$
 such that $\nbigm_{|\nbigx_P}$ has a KMS-structure
 with which $L_{|\nbigx_P}$ is compatible.
\end{itemize}
Then, $\nbigm$ has a unique KMS-structure
with which $L$ is compatible.
\hfill\qed
\end{lem}

We state it in a slightly different way.

\begin{cor}
\label{cor;10.11.9.10}
Let $\nbigm$ be 
a good-KMS $\nbigr_{X(\ast D)}$-module.
Let $L$ be a filtration of $\nbigm$
in the category of smooth $\nbigr_{X(\ast D)}$-modules.
Assume the following:
\begin{itemize}
\item
 $\Gr^L(\nbigm)$ has a good-KMS structure.
\item
 For any smooth point $P$ of $D$,
 there exists a small neighbourhood $X_P$ of $P$
 such that $L_{|\nbigx_P}$ is compatible with
 the $KMS$-structure of
 $\nbigm_{|\nbigx_P}$.
\end{itemize}
Then, $L$ is compatible with
the $KMS$-structure of $\nbigm$.
\hfill\qed
\end{cor}

\paragraph{Example}

Let $X=\cnum_z$.
Let $a,b\in\real$ such that $a\neq b$.
We consider the $\nbigr_X$-module
$\nbigm:=
 \bigoplus_{i=1,2}
 \nbigo_{\nbigx}(\ast z)\,e_i$
with $z\deldel_z e_1=-\lambda\,a\,e_1$
and $z\deldel_z e_2=e_1-\lambda\,b\,e_2$.
We put $v_2:=e_1+\lambda\,(a-b)\,e_2$.
Then, we have
$z\deldel_z v_2=-\lambda\,b\,v_2$.

Assume $a<b$.
For $c=a,b$,
let $\nbigm(c)$ be $\nbigo_{\nbigx}(\ast z)\,f$
with $z\deldel f=\lambda\,c\,f$.
We have the morphism
$\varphi:\nbigm(b)\lrarr \nbigm$
given by $f\longmapsto v$.
By a direct computation,
we can check that the cokernel of
$\psitilde_b(\varphi):
 \psitilde_b(\nbigm(b))\lrarr 
 \psitilde_b(\nbigm)$ is not strict.
The cokernel of $\varphi$
is generated by $[e_2]$,
and it is naturally isomorphic to
$\nbigm(a)$.
We can check that
the morphism $\nbigm\lrarr\nbigm(a)$
is not strict with respect to
the $V$-filtrations $V^{(0)}$.

\subsection{Curve test}

Let $\nbigm$ be a good smooth $\nbigr_{X(\ast D)}$-module.
Let $L$ be a filtration in the category of
$\nbigr_{X(\ast)}$-modules.
\begin{prop}
\label{prop;13.5.14.11}
We suppose the following:
\begin{itemize}
\item
$\Gr^L(\nbigm)$ is good-KMS.
\item
 Let $C$ be any curve in $X$
 which intersects with the smooth part of $D$
 transversally.
 Then, $\nbigm_{|C}$ has a KMS-structure
 with which $L_{|C}$ is compatible.
\end{itemize}
Then, $\nbigm$ has a KMS-structure
with which $L$ is compatible.
\end{prop}
\pf
By Lemma \ref{lem;10.9.28.12},
we may assume 
$X=\Delta^n$ and $D=\{z_1=0\}$.
As in the proof of Proposition \ref{prop;10.10.9.2},
$\nbigm$ has KMS-structure at generic $\lambda_0$.
(This case is easier. Indeed,
because the divisor is smooth, we have only to 
consider the lattice in the formal completion.)

Let us consider the case that $\lambda_0\neq 0$ is not
necessarily generic.
Let $U(\lambda_0)\subset\cnum_{\lambda}$
be a small neighbourhood of $\lambda_0$.
For any generic $\lambda_1\in U(\lambda_0)$,
we take a small neighbourhood $U(\lambda_1)\subset U(\lambda_0)$
such that any $\lambda\in U(\lambda_1)$ is generic.
We can construct
$\nbigqzero(\nbigm_{|U(\lambda_1)\times X})$
from $\nbigq^{(\lambda_1)}(\nbigm_{|U(\lambda_1)\times X})$
by the relation of $\nbigqzero$ and $\nbigq^{(\lambda_1)}$,
which is compatible with $L$.
(See \S2.8.4 of \cite{mochi7}, for example.)
By applying a general result below,
we obtain that it is extended to the KMS-structure at $\lambda_0$
which is compatible with $L$.
The case $\lambda_0=0$ can be argued similarly.
\hfill\qed

\begin{lem}
\label{lem;13.5.14.10}
Let $Z$ be a complex manifold.
We set $W:=\{(z_1,z_2)\,|\,|z_i|<1\}$
and $W^{\ast}:=W\setminus\{(0,0)\}$.
Let $E_i$ $(i=1,2)$ be a locally free sheaf on $Z\times W$.
We consider an extension of $\nbigo_{Z\times W^{\ast}}$-modules
\begin{equation}
\label{eq;13.5.14.1}
 0\lrarr E_{1|Z\times W^{\ast}}\lrarr
 E'\lrarr E_{2|Z\times W^{\ast}}\lrarr 0
\end{equation}
We suppose the following condition:
\begin{itemize}
\item
 For any $P\in Z$,
 the specialization of {\rm(\ref{eq;13.5.14.1})}
 to $\{P\}\times W^{\ast}$
 is prolonged to
\[
 0\lrarr E_{1|\{P\}\times W}
\lrarr E'_{P\times W}
\lrarr E_{2|\{P\}\times W}\lrarr 0.
\]
\end{itemize}
Then,  we have a  unique extension
$0\lrarr E_1\lrarr E\lrarr E_2\lrarr 0$ on $Z\times W$
whose restriction to
$Z\times W^{\ast}$ is {\rm(\ref{eq;13.5.14.1})}.
\end{lem}
\pf
By using Theorem 6.10 in \cite{Siu-extension},
we obtain a reflexive coherent $\nbigo_{Z\times W}$-module
$E_3$ such that
$E_{3|Z\times W^{\ast}}=E'$.
We naturally have the exact sequence
$0\lrarr E_1\lrarr E_3\lrarr E_2$.
Let us prove that $\kappa:E_3\lrarr E_2$ is surjective.
Fix a point $P\in Z$,
and we will shrink $Z$ around a fixed point $P$.
We may assume that $E_2=\nbigo_{Z\times W}$.
Let $\gminim$ denote the ideal sheaf of
$Z\times\{(0,0)\}$.
If we take a sufficiently large $N$,
the composite of
$\gminim^N\lrarr
 \nbigo_X\lrarr \Cok(\kappa)$ is $0$.
We have the induced extension
\begin{equation}
\label{eq;13.5.14.2}
 0\lrarr E_1\lrarr E_4\lrarr \gminim^N\lrarr 0.
\end{equation}
We have only to prove that (\ref{eq;13.5.14.2})
has a splitting.
Let $\gminic$ be the extension class of (\ref{eq;13.5.14.2}),
which is a section of
$\next^1(\gminim^N,E_1)$.
Note that $\gminim^N$ is flat over $\nbigo_Z$.
Let $\gminim_0\subset\nbigo_W$ be the ideal sheaf of $(0,0)$.
For any $Q\in Z$, 
we have
\begin{equation}
\label{eq;13.5.14.5}
 \next^1(\gminim^N,E_1)
 \otimes\nbigo_{\{Q\}\times W}
\simeq
 \next^1\bigl(
 \gminim^N_0,E_{1|\{Q\}\times W}
 \bigr)
\end{equation}
By the assumption,
the image of $\gminic$ in (\ref{eq;13.5.14.5})
are $0$ for any $Q$.
By the coherence of 
$\next^1(\gminim^N,E_1)$,
we obtain that $\gminic$ is $0$.
Thus, we obtain Lemma \ref{lem;13.5.14.10}.
\hfill\qed

\section{Canonical prolongations of smooth good-KMS $\nbigr$-modules}
\label{subsection;11.4.3.12}

We shall study the canonical prolongations
of good-KMS smooth $\nbigr$-modules
across normal crossing hypersurfaces.
\index{canonical prolongations}

\subsection{Goal}
\label{subsection;10.8.21.110}

Let $X$ be a complex manifold.
Let $D$ be a simple normal crossing hypersurface
with the irreducible decomposition
$D=\bigcup_{i\in\Lambda}D_i$.
We set $(\nbigx,\nbigd):=\cnum_{\lambda}\times(X,D)$.
Let $\nbigm$ be a good-KMS smooth
$\nbigr_{X(\ast D)}$-module.
An $\nbigr_X$-module $\nbigmtilde$
with an isomorphism
$\rho:\nbigmtilde\otimes\nbigo_{\nbigx}(\ast \nbigd)
\simeq \nbigm$
is called a prolongment of $\nbigm$.
We say prolongments
$(\nbigmtilde_i,\rho_i)$ $(i=1,2)$
are isomorphic,
if there exists an isomorphism
of $\nbigr_X$-modules
$F:\nbigmtilde_1\lrarr\nbigmtilde_2$ 
such that
the following diagram is commutative:
\[
 \begin{CD}
 \nbigmtilde_1\otimes\nbigo_{\nbigx}(\ast \nbigd)
 @>{F}>>
 \nbigmtilde_2\otimes\nbigo_{\nbigx}(\ast \nbigd)\\
 @V{\rho_1}VV @V{\rho_2}VV \\
 \nbigm @>{=}>>\nbigm
 \end{CD}
\]
A prolongment $(\nbigmtilde,\rho)$
is often denoted just by $\nbigmtilde$
in the following.
We shall prove the following proposition
in \S\ref{subsection;13.3.28.10}--\ref{subsection;13.3.28.11}.

\begin{prop}
\label{prop;10.9.30.20}
For any decomposition $\Lambda=I\sqcup J$,
there exists a prolongment
$\nbigm[\ast I!J]$ of $\nbigm$
with the following property,
which is unique up to isomorphisms:
\index{sheaf $\nbigm[\ast I\bikkuri J]$}
\begin{description}
\item[(P1)]
$\nbigm[\ast I!J]$ is 
$\nbigr_X$-coherent, holonomic and strict.
\item[(P2)]
For $P\in D$, we take a small coordinate
neighbourhood $(X_P;z_1,\ldots,z_n)$ around $P$
such that, for each $i\in \Lambda$,
we have $D_i\cap X_P=\emptyset$ or
$D_i\cap X_P=\{z_{k(i)}=0\}$
for some $k(i)$.
Then, 
$\nbigm[\ast I!J]_P:=
 \nbigm[\ast I!J]_{|X_P}$ is
strictly specializable along $z_k$ for any $k$,
and we have
\[
\nbigm[\ast I!J]_P[\ast z_k]=
 \nbigm[\ast I!J]_P\,\,\,\,\,
(k\in I),
\quad
 \nbigm[\ast I!J]_P[! z_k]=
 \nbigm[\ast I!J]_P\,\,\,\,\,
 (k\in J).
\]
\end{description}
\end{prop}
We shall also prove the following lemmas.
\begin{lem}
\label{lem;10.9.30.1}
Let $\varphi:\nbigm[\ast I!J]\simeq
 \nbigm[\ast I!J]$ be an isomorphism
 as prolongments of $\nbigm$.
Then, $\varphi$ is the identity.
\end{lem}

For any $I\subset\Lambda$
and $i\in\Lambda$,
we put $I_{\cup i}:=I\cup\{i\}$
and $I_{\setminus i}:=I\setminus\{i\}$.

\begin{lem}
\label{lem;10.9.30.21}
We have isomorphisms
of the following prolongments
\[
 \bigl(\nbigm[\ast I!J]_P\bigr)[\ast z_{k(i)}]
\simeq
 \nbigm[\ast I_{\cup i}!J_{\setminus i}]_P,
\quad\quad
 \bigl(\nbigm[\ast I!J]_P\bigr)[!z_{k(i)}]
\simeq
 \nbigm[\ast I_{\setminus i}!J_{\cup i}]_P.
\]
\end{lem}

\subsection{Uniqueness and Lemma \ref{lem;10.9.30.1}}
\label{subsection;13.3.28.10}

Let us return to the situation in 
\S\ref{subsection;10.8.21.110}.
Let $\nbigm[\ast I!J]_{\kappa}$ $(\kappa=1,2)$
be prolongments of $\nbigm$
satisfying the conditions (P1) and (P2).
We prove that
there uniquely exists an isomorphism
$\nbigm[\ast I!J]_1\simeq
 \nbigm[\ast I!J]_2$.

By the uniqueness in the claim,
we have only to consider the case
$X=\Delta^n$ and $D=\bigcup_{i=1}^{\ell}\{z_i=0\}$.
For $L\subset \ellsitabar$,
``$\otimes\nbigo(\ast\nbigd(L))$''
is denoted by ``$(\ast L)$'' for simplicity.
We also denote $\nbigr_{X(\ast D_L)}$
by $\nbigr_X(\ast L)$.
An $\nbigr_X(\ast L)$-module $\nbigmtilde$
with an isomorphism
$\rho:\nbigmtilde\otimes\nbigo_{\nbigx}(\ast \nbigd)
\simeq \nbigm$
is called a prolongment of $\nbigm$.
An isomorphism of prolongments as
$\nbigr_X(\ast L)$-modules
is defined as in the case of prolongments
as $\nbigr_X$-modules.

We consider
the $\nbigr_X(\ast L)$-modules
$\nbigm[\ast I!J]_{\kappa}(\ast L)$.
They are strictly specializable along $z_i$
as $\nbigr_X(\ast L)$-modules.
We have the following natural isomorphisms
as $\nbigr_X(\ast L)$-modules
for any $i\not\in L$:
\[
 \nbigm[\ast I!J]_{\kappa}(\ast L)[\ast z_i]
\simeq
 \nbigm[\ast I!J]_{\kappa}(\ast L)
\quad (i\in I)
\]
\[
 \nbigm[\ast I!J]_{\kappa}(\ast L)[! z_i]
\simeq
 \nbigm[\ast I!J]_{\kappa}(\ast L)
\quad (i\in J)
\]
Hence, we have the following natural isomorphisms
as $\nbigr_X(\ast L_{\setminus i})$-modules
for any $i\in L$:
\begin{equation}
\label{eq;10.9.30.10}
 \nbigm[\ast I!J]_{\kappa}(\ast L)[\ast z_i]
\simeq
 \nbigm[\ast I!J]_{\kappa}(\ast L_{\setminus i})
\quad (i\in I)
\end{equation}
\begin{equation}
\label{eq;10.9.30.11}
 \nbigm[\ast I!J]_{\kappa}(\ast L)[! z_i]
\simeq
 \nbigm[\ast I!J]_{\kappa}(\ast L_{\setminus i})
\quad (i\in J)
\end{equation}
We have the unique isomorphism
$\varphi:\nbigm[\ast I!J]_1(\ast \ellsitabar)
\simeq
 \nbigm[\ast I!J]_2(\ast \ellsitabar)$.
By using (\ref{eq;10.9.30.10})
and (\ref{eq;10.9.30.11}),
we obtain unique isomorphisms
$\nbigm[\ast I!J]_1(\ast L)
\simeq
 \nbigm[\ast I!J]_2(\ast L)$ for any $L\subset\ellsitabar$
as prolongments,
by using a descending induction on $|L|$.
\hfill\qed

\subsection{Local construction}
\label{subsection;10.8.21.50}

Let $X:=\Delta^n$ and 
$D:=\bigcup_{i=1}^{\ell}\{z_i=0\}$.
In this subsection,
let $V_0\nbigr_X\subset\nbigr_X$ be 
generated by
$z_i\deldel_i$ $(i\leq \ell)$
and $\deldel_i$ $(i>\ell)$
over $\nbigo_{\nbigx}$.
Let $\nbigk\subset\cnum$
be a small neighbourhood of $\lambda_0$.
We set $(\nbigxzero,\nbigdzero):=\nbigk\times(X,D)$.
Let $\nbigm$ be a good-KMS $\nbigr_{X(\ast D)}$-module
on $\nbigxzero$.
Let $(\nbigqzero_{\ast}\nbigm,\DD)$ be 
the associated filtered $\lambda$-flat bundles
at $\lambda_0$ on $(\nbigxzero,\nbigdzero)$.

Let $\ellsitabar=I\sqcup J$.
We put $\veca(I,J):=\vecdelta_I+(1-\epsilon)\vecdelta_J$
for some sufficiently small $\epsilon>0$.
We obtain a coherent $V_0\nbigr_X$-module
$V_0\nbigr_X\cdot
 \nbigqzero_{\veca(I,J)}\nbigm
\subset\nbigqzero\nbigm$.
Then, we obtain a coherent $\nbigr_X$-module
\[
 \nbigmzero[\ast I!J]:=
 \nbigr_X\otimes_{V_0\nbigr_X}
 \bigl(
 V_0\nbigr_X\cdot
 \nbigqzero_{\veca(I,J)}\nbigm
 \bigr).
\]
We shall prove that
$\nbigmzero[\ast I!J]$
has the desired property.

\begin{lem}
\label{lem;13.8.29.2}
The $\nbigr_X$-module
$\nbigmzero[\ast I!J]$ is holonomic,
and its characteristic variety is contained in
$\nbigs=\bigcup_{L\subset\ellsitabar}
\nbigk\times T^{\ast}_{D_L}X$.
\end{lem}
\pf
Let $F_0$ be the image of the natural morphism
$\nbigqzero_{\veca(I,J)}\nbigm
\lrarr
 \nbigmzero[\ast I!J]$.
By the construction,
$\nbigmzero[\ast I!J]$
is generated by $F_0$ over $\nbigr_X$.
For $\vecp\in\seisuu_{\geq\,0}^{n}$,
we put $\deldel^{\vecp}:=\prod\deldel_i^{p_i}$.
We set
$F_m:=\sum_{|\vecp|\leq m}
 \deldel^{\vecp}F_0$.
Then, $\{F_m\}$ is a coherent filtration
of $\nbigmzero[\ast I!J]$.
Let us prove that the support of
$\Gr^F\nbigmzero[\ast I!J]$
is contained in $\nbigs$.
Let $\pi:\nbigk\times T^{\ast}X\lrarr X$.

\vspace{.1in}
First, let us consider the case
$\nbigm$ is regular-KMS.
Let $Q\in D_K^{\circ}$.
For any $j\not\in K$,
we have $\deldel_j F_0\subset F_0$
around $Q$.
Hence, we obtain
$\deldel_j F_m\subset F_m$ for any $m\geq 0$.
Then, the action of $\deldel_j$
on $\Gr^F\nbigmzero[\ast I!J]$ is $0$
around $Q$.
Then, we obtain that 
$Ch(\nbigmzero[\ast I!J])
\cap\pi^{-1}(Q)
\subset
 \nbigk\times
 (N^{\ast}_{D_K}X)_Q$.

\vspace{.1in}
We put $\gminia:=\vecz^{\vecm}$
for some $\vecm\in\seisuu_{<0}^p$,
where $1\leq p\leq \ell$.
Let $L(\gminia)=\nbigo_{\nbigxzero}\,e$
with $\DD e=e\,d\gminia$.
Let $\pi:\nbigxzero\lrarr \nbigxzero$ 
be the ramified covering,
given by
$z_i^{q_i}$ for $1\leq i\leq p$
and $z_i$ for $p+1\leq i\leq n$.
Let us consider the case 
$\nbigm$ is the tensor product of
$\pi_{\ast}L(\gminia)$ and
regular-KMS $\nbigm'$.
Let $K\subset\ellsitabar$.
Let $Q\in D_K^{\circ}$.
If $K\cap\pbar=\emptyset$,
we have
$Ch(\nbigmzero[\ast I!J])\subset\nbigs$
around $Q$ by the consideration
in the regular singular case.
Let us consider the case $K\cap\pbar\neq\emptyset$.
Take $i\in K\cap \pbar$.
For $j\in \pbar\setminus K$,
we put 
$v_j:=m_j^{-1}\,z_j\deldel_j-m_i^{-1}\,z_i\deldel_i$.
Note we have
$v_j e=0$.
Hence, we have $v_jF_m\subset F_m$
around $Q$.
For $j\in\ellsitabar\setminus (\pbar\cup I)$,
we have 
$\deldel_j F_0\subset F_0$,
and hence $\deldel_j F_m\subset F_m$.
Thus, we obtain
$Ch(\nbigm[\ast I!J])_{|\pi^{-1}(Q)}
\subset
\nbigk\times (N^{\ast}_{D_I}X)_Q$.

The general case can be reduced
to the above cases,
by using the formal completion.
\hfill\qed

\subsection{Preliminary}

We have the formal decomposition
into the regular part and the irregular part:
\[
 \nbigqzero_{\ast}\nbigm_{|\nbigdhat_i}
=\nbigqzero_{\ast}\nbigm^{(\reg)}_{\nbigdhat_i}
\oplus
 \nbigqzero_{\ast}\nbigm^{(\irr)}_{\nbigdhat_i}
\]
(See \S\ref{subsection;13.5.4.30}.)
Then, we set
\[
 \lefttop{i}\psitildezero_u\bigl(
 \nbigqzero_{\veca}\nbigm
 \bigr):=
 \lefttop{i}\nbiggzero_{u+\vecdelta}
 \bigl(
 \nbigqzero_{\veca}\nbigm^{(\reg)}_{\nbigdhat_i}
 \bigr)
\]
Here, $\vecdelta=(1,0)\in\real\times\cnum$.
We obtain a good-KMS 
$\lefttop{i}\psitildezero_u\bigl(
 \nbigqzero_{\ast}\nbigm
 \bigr)$.

\vspace{.1in}

Let $I\sqcup J\sqcup\{i\}=\ellsitabar$
be a decomposition.
Let $q_j$ denote the projection to the $j$-th component.
Let $\veca(b)\in\real^{\ell}$ be
determined by
$q_j(\veca(b))=1$ for $j\in I$,
$q_j(\veca(b))=1-\epsilon$ for $j\in J$,
and $q_i(\veca(b))=b+1$,
where $\epsilon>0$ is any sufficiently small number.

\begin{lem}
\label{lem;10.8.21.35}
We have the following natural isomorphism:
\begin{equation}
\label{eq;10.8.21.20}
\frac{V_0\nbigr_{X}\cdot\nbigqzero_{\veca(b)}\nbigm}
{V_0\nbigr_X\cdot\nbigqzero_{\veca(b-\epsilon)}\nbigm}
\simeq
\bigoplus_{\paramap(\lambda_0,u)=b}
 V_0\nbigr_{D_i}\cdot
 \lefttop{i}\psitildezero_u(\nbigqzero_{\veca(b)}\nbigm)
\end{equation}
\end{lem}
\pf
We have
$V_0\nbigr_{X}\cdot
 \nbigqzero_{\veca(b)}\nbigm_{|\nbigdhat_i}=
V_0\nbigr_{X}\cdot
 \nbigqzero_{\veca(b)}\nbigm_{\nbigdhat_i}^{(\reg)}
\oplus
V_0\nbigr_{X}\cdot
 \nbigqzero_{\veca(b)}\nbigm_{\nbigdhat_i}^{(\irr)}$.
Let $V_0\nbigr_{X,\setminus i}\subset
 \nbigr_X$ be generated by
$z_j\deldel_j$ $(j\in\ellsitabar\setminus\{i\})$
and $\deldel_j$ $(j\in\nbar\setminus\ellsitabar)$
over $\nbigo_{\nbigx}$.
The following natural morphisms
are isomorphisms:
\begin{equation}
\label{eq;10.8.21.30}
V_0\nbigr_{X,\setminus i}\cdot
  \nbigqzero_{\veca(b)}\nbigm_{\nbigdhat_i}^{(\reg)}
\lrarr
V_0\nbigr_{X}\cdot
 \nbigqzero_{\veca(b)}\nbigm_{\nbigdhat_i}^{(\reg)}
\end{equation}
\begin{equation}
V_0\nbigr_{X}\cdot
 \nbigqzero_{\veca(b-\epsilon)}\nbigm_{\nbigdhat_i}^{(\irr)}
\lrarr
 V_0\nbigr_X\cdot
 \nbigqzero_{\veca(b)}\nbigm^{(\irr)}_{\nbigdhat_i}.
\end{equation}
Indeed, it can be reduced to the case
that $\nbigq\nbigm$
is unramified.
We have only to prove that their formal completions
along $\nbigk\times P$ are isomorphisms
for each $P\in D_i$.
Then, the claim can be checked by a direct computation.

We obtain the following:
\[
 \frac{V_0\nbigr_{X}\cdot\nbigqzero_{\veca(b)}\nbigm}
{V_0\nbigr_X\cdot\nbigqzero_{\veca(b-\epsilon)}\nbigm}
\simeq
  \frac{V_0\nbigr_{X}\cdot
 \nbigqzero_{\veca(b)}\nbigm_{|\nbigdhat_i}}
{V_0\nbigr_X\cdot
 \nbigqzero_{\veca(b-\epsilon)}\nbigm_{|\nbigdhat_i}}
\simeq
\bigoplus_{\paramap(\lambda_0,u)=b}
 V_0\nbigr_{D_i}\cdot
 \lefttop{i}\psitildezero_u\bigl(
 \nbigqzero_{\veca(b)}\nbigm\bigr)
\]
Thus, we are done.
\hfill\qed

\vspace{.1in}
Let $\lefttop{i}V_0\nbigr_X\subset\nbigr_X$
be generated by
$\deldel_j$ $(j\neq i)$
and $z_i\deldel_i$
over $\nbigo_{\nbigx}$.

\begin{cor}
We have the following natural isomorphism:
{\small
\begin{equation}
\label{eq;10.8.21.34}
\frac{\lefttop{i}V_0\nbigr_{X}
 \otimes_{V_0\nbigr_{X}}\bigl(
 V_0\nbigr_{X}\cdot\nbigqzero_{\veca(b)}\nbigm
 \bigr)}
{\lefttop{i}V_0\nbigr_{X}\otimes_{V_0\nbigr_{X}}
 \bigl(
 V_0\nbigr_X\cdot\nbigqzero_{\veca(b-\epsilon)}\nbigm
 \bigr)}
\simeq
\bigoplus_{\paramap(\lambda_0,u)=b}
 \nbigr_{D_i}\otimes_{V_0\nbigr_{D_i}}
 \Bigl(
 V_0\nbigr_{D_i}\cdot
 \lefttop{i}\psitildezero_u
 \bigl(\nbigqzero_{\veca(b)}\nbigm\bigr)
 \Bigr)
\end{equation}
}
\end{cor}
\pf
By Lemma \ref{lem;10.8.21.35},
we obtain the following:
{\small
\begin{equation}
\label{eq;10.8.21.33}
\frac{\lefttop{i}V_0\nbigr_{X}
 \otimes_{V_0\nbigr_{X}}\bigl(
 V_0\nbigr_{X}\cdot\nbigqzero_{\veca(b)}\nbigm
 \bigr)}
{\lefttop{i}V_0\nbigr_{X}\otimes_{V_0\nbigr_{X}}
 \bigl(
 V_0\nbigr_X\cdot\nbigqzero_{\veca(b-\epsilon)}\nbigm
 \bigr)}
\simeq
\bigoplus_{\paramap(\lambda_0,u)=b}
 \lefttop{i}V_0\nbigr_X\otimes_{V_0\nbigr_X} 
 \Bigl(
 V_0\nbigr_{D_i}\cdot
 \lefttop{i}\psitildezero_u(\nbigqzero_{\veca}\nbigm)
 \Bigr)
\end{equation}
}
It is easy to observe that the right hand side of
(\ref{eq;10.8.21.33}) is naturally isomorphic
to the right hand side of (\ref{eq;10.8.21.34}).
\hfill\qed

\subsection{Some filtrations}
\label{subsection;10.8.21.51}

For a subset $K\subset\ellsitabar$,
let $\lefttop{K}V_0\nbigr_{X}\subset\nbigr_X$
be generated by
$\deldel_i$ $(i\in \nbar\setminus K)$
and $z_i\deldel_i$ $(i\in K)$ over $\nbigo_{\nbigx}$.
Let $\ellsitabar=I\sqcup J\sqcup K$ be a decomposition.
We put 
$\veca(I,J):=\vecdelta_I+(1-\epsilon)\vecdelta_J
 \in\real^{I\sqcup J}$.
Let $\nbigqzero_{\veca(I,J)}\nbigm$
mean 
\[
 \nbigqzero_{\veca(I,J)+\vecc}\nbigm
 \otimes\nbigo_{\nbigxzero}(\ast \nbigdzero(K))
\]
for any $\vecc\in\real^K$.
We consider
the following coherent $\nbigr_X(\ast \nbigdzero(K))$-module:
\[
 \nbigmzero[\ast I!J]:=
 \lefttop{K}V_0\nbigr_X\otimes_{V_0\nbigr_X}
 \bigl(
 V_0\nbigr_X\cdot
 \nbigqzero_{\veca(I,J)}\nbigm
 \bigr)
\simeq
 \nbigr_X\otimes_{V_0\nbigr_X}
  \bigl(
 V_0\nbigr_X\cdot
 \nbigqzero_{\veca(I,J)}\nbigm
 \bigr)
\]
For $\vecb\in\real^K$,
we consider 
the following $\lefttop{K}V_0\nbigr_X$-module:
\[
 \lefttop{K}\Vzero_{\vecb}\nbigmzero[\ast I!J]:=
 \lefttop{K}V_0\nbigr_{X}
 \otimes_{V_0\nbigr_{X}}
 \bigl(
 V_0\nbigr_{X}\cdot
 \nbigqzero_{\veca(I,J)+\vecb+\vecdelta_K}\nbigm
\bigr)
\]

For $i\in I$ and $b\leq 0$,
or for $i\in J$ and $b<0$,
let $\lefttop{K,i}\Vzero_{\vecb,b}
 \nbigmzero[\ast I!J]$
be the image of the following morphism:
\begin{equation}
\label{eq;10.8.21.40}
 \lefttop{K\sqcup\{i\}}V_0
 \nbigr_{X}
 \otimes_{V_0\nbigr_{X}}\bigl(
 V_0\nbigr_{X}\cdot
 \nbigqzero_{\veca'+\vecb+\vecdelta_K}\nbigm\bigr)
\lrarr
 \lefttop{K}\Vzero_{\vecb}\nbigmzero[\ast I!J].
\end{equation}
Here, $\veca'\in\real^{I\sqcup J}$ is
determined such that
$q_j(\veca')=q_j(\veca(I,J))$
if $j\neq i$,
and $q_i(\veca')=b+1$.

For $i\in I$ and $b>0$,
or for $i\in J$ and $b\geq 0$,
we set 
\[
  \lefttop{K,i}\Vzero_{\vecb,b}\bigl(
 \nbigmzero[\ast I!J]
 \bigr)
=\sum_{(c,p)\in\nbigu(b,i)}
 \deldel_i^p \bigl(
 \lefttop{K,i}\Vzero_{\vecb,c}(
 \nbigmzero[\ast I!J])
 \bigr)
\]
Here, $\nbigu(b,i)$ denotes the set
$\bigl\{
 (c,p)\in\real_{\leq 0}\times\seisuu_{\geq 0}\,\big|\,
 c+p\leq b
 \bigr\}$ if $i\in I$,
or 
$\bigl\{
 (c,p)\in\real_{< 0}\times\seisuu_{\geq 0}\,\big|\,
 c+p\leq b
 \bigr\}$ if $i\in J$.

\begin{lem}
\label{lem;10.9.30.15}
Let us prove the following claims
by an induction on $m=|I\sqcup J|$
and $\dim X$:
\begin{description}
\item[$P(m)$:]
The morphisms {\rm(\ref{eq;10.8.21.40})} are injective.
\item[$Q(m)$:]
For $\vecc,\vecd\in\real^K$ with
$\vecc\leq\vecd$,
the natural morphisms
$\lefttop{K}\Vzero_{\vecc}\nbigmzero[\ast I!J]
\lrarr
 \lefttop{K}\Vzero_{\vecd}\nbigmzero[\ast I!J]$
are injective.
 In particular, we have the injectivity of
 the morphism
 $\lefttop{K}\Vzero_{\vecc}\nbigmzero[\ast I!J]
 \lrarr
 \nbigmzero[\ast I!J]$.
\item[$R(m)$:]
$\lefttop{K}\Vzero_{\vecc}\nbigmzero[\ast I!J]$
are strict.
\end{description}
\end{lem}
\pf
If $\dim X=0$, the claim is trivial.
The claim $P(0)$ is trivial,
and the claims $Q(0)$ and $R(0)$ are obvious.
We obtain $P(m)$ from $Q(m-1)$,
by considering the composition
of the morphism (\ref{eq;10.8.21.40})
and the natural one
\[
 \lefttop{K}\Vzero_{\vecb}\nbigmzero[\ast I!J]\lrarr
\lefttop{K}\Vzero_{\vecb}\nbigmzero[\ast I!J](\ast i).
\]

Assume $P(m)$ and 
the claims in the strictly lower dimensional case.
For $i\in I$ and $b\leq 0$,
or for $i\in J$ and $b<0$,
we obtain the following,
by using Lemma \ref{lem;10.8.21.35}:
\begin{multline}
\label{eq;10.8.21.42}
\lefttop{i}\Gr^{\Vzero}_b
 \lefttop{K}\Vzero_{\vecb}\nbigmzero[\ast I!J]
\simeq
 \bigoplus_{\paramap(\lambda_0,u)=b}
 \lefttop{K}V_0\nbigr_{D_i}
 \otimes_{V_0\nbigr_{D_i}}
 \Bigl(
  V_0\nbigr_{D_i}\cdot
\lefttop{i}\psitilde_u\bigl(
 \nbigqzero_{\veca'+\vecb}\nbigm
 \bigr)
 \Bigr)
 \\
=\bigoplus_{\paramap(\lambda_0,u)=b}
\lefttop{K}\Vzero_{\vecb}\bigl(
 \lefttop{i}\psitilde_u(\nbigm)^{(\lambda_0)}
 [\ast I_{\setminus i}\,!J_{\setminus i}]
\bigr)
\end{multline}
For $i\in I$ and $d>0$,
(resp. $i\in J$ and $d\geq 0$),
we take $p\in\seisuu_{>0}$
such that $-1<d-p\leq 0$
(resp. $-1\leq d-p<0$).
We consider the following surjective morphism:
\begin{equation}
\label{eq;10.8.21.41}
 \deldel_i^p:
 \lefttop{i}\Gr^{\Vzero}_{d-p}
 \lefttop{K}\Vzero_{\vecb}\nbigmzero[\ast I!J]
\lrarr
 \lefttop{i}\Gr^{\Vzero}_{d}
 \lefttop{K}
 \Vzero_{\vecb}\nbigmzero[\ast I!J]
\end{equation}
If $i\in I$,
the morphism $z_i^p\deldel_i^p$
on 
$\lefttop{i}\Gr^{\Vzero}_{d-p}
 \lefttop{K}\Vzero_{\vecb}\nbigmzero[\ast I!J]$
is injective by the assumption
of the induction on the base space.
Hence, we obtain that
(\ref{eq;10.8.21.41}) is injective.
If $i\in J$ and $d=0$,
we can prove that
the restriction of
(\ref{eq;10.8.21.41}) to 
$\nbigdzero_i\setminus \bigcup_{j\neq i}\nbigdzero_{j}$
is isomorphism,
by using Lemma \ref{lem;10.8.20.3}.
Then, by using the description (\ref{eq;10.8.21.42})
and the hypothesis of the induction
on $\dim X$,
we obtain that 
(\ref{eq;10.8.21.41}) is an isomorphism.
In the case $d>0$,
we can check that (\ref{eq;10.8.21.41})
is an isomorphism
by using the argument in the case $i\in I$.

\vspace{.1in}
Now, assume $P(m)$, $R(m-1)$, $Q(m-1)$
and the claims in the strictly lower dimensional case.
We obtain that
\[
 \lefttop{i}\Gr^{\Vzero}_{d}
 \lefttop{K}\Vzero_{\vecb}
 \nbigmzero[\ast I!J]
\lrarr
 \lefttop{i}\Gr^{\Vzero}_d
 \lefttop{K}\Vzero_{\vecc}\nbigmzero[\ast I!J]
\]
is injective for each $d$
by using the isomorphisms (\ref{eq;10.8.21.42})
and (\ref{eq;10.8.21.41}).
Then, $Q(m)$ follows.
We also obtain $R(m)$
from $R(m-1)$ and the strictness
in the lower dimensional case.
Thus, the proof of Lemma \ref{lem;10.9.30.15}
is finished.
\hfill\qed

\vspace{.1in}
\begin{cor}
\mbox{{}}\label{cor;10.9.30.22}
Let $I\sqcup J\subset\ellsitabar$,
and $K:=\ellsitabar\setminus (I\sqcup J)$.
\begin{itemize}
\item
$\nbigmzero[\ast I!J]$
is a coherent, holonomic and strict
$\nbigr_X(\ast K)$-module.
\item
It is strictly specializable along $z_i$
with the $V$-filtration $\lefttop{i}\Vzero$.
We naturally have
\begin{equation}
\label{eq;10.8.21.101}
 \lefttop{i}\psitildezero_u\bigl(
 \nbigmzero[\ast I!J]\bigr)
\simeq
\lefttop{i}\psitilde_u(\nbigm)^{(\lambda_0)}
 [\ast I_{\setminus i}!J_{\setminus i}].
\end{equation}
\item
The following morphisms are isomorphisms:
\[
 \deldel_i:
 \lefttop{i}\psitildezero_{-\vecdelta}
 \bigl(
 \nbigmzero[\ast I!J]\bigr)
\lrarr
 \lefttop{i}\psitildezero_0\bigl(
 \nbigmzero[\ast I !J]\bigr)
\quad\quad(i\in J)
\]
\[
 z_i:
 \lefttop{i}\psitildezero_0\bigl(
 \nbigmzero[\ast I!J]\bigr)
\lrarr
 \lefttop{i}\psitildezero_{-\vecdelta}\bigl(
 \nbigmzero[\ast I!J]\bigr)
\quad\quad(i\in I)
\]
In particular, we have isomorphisms
$\nbigmzero[\ast I!J][\ast z_i]
\simeq
 \nbigmzero[\ast I_{\cup i}!J_{\setminus i}]$
and
$\nbigmzero[\ast I!J][! z_i]
\simeq
 \nbigmzero[\ast I_{\setminus i}!J_{\cup i}]$.
\hfill\qed
\end{itemize}
\end{cor}

\begin{rem}
For any good meromorphic flat bundle $\nbigv$
on $(X,D)$,
we have a similar description of
$\nbigv[\ast I!J]$.
\hfill\qed
\end{rem}

\subsection{Globalization}
\label{subsection;13.3.28.11}

Let us return to the situation in 
\S\ref{subsection;10.8.21.110}.
Let $I\sqcup J=\Lambda$ be a decomposition.
Let $P\in D$.
We take a small coordinate neighbourhood
$(X_P;z_1,\ldots,z_n)$ of $X$ around $P$
such that $D_P:=D\cap X_P=
 \bigcup_{i=1}^{\ell}\{z_i=0\}$.
Let $U(\lambda_0)$ be a sufficiently small neighbourhood
of $\lambda_0$.
We set $(\nbigxzero_P,\nbigdzero_P):=
 U(\lambda_0)\times (X_P,X_P\cap D)$.
Apply the procedure 
in \S\ref{subsection;10.8.21.50}
to $(\nbigqzero_{\ast}\nbigm_P,\DD)
 :=
 (\nbigqzero_{\ast}\nbigm,\DD)_{|\nbigxzero_P}$,
and we obtain an $\nbigr_{X_P}$-module
$\nbigmzero_P[\ast I_P!J_P]$
on $\nbigxzero_P$
for any decomposition $\ellsitabar=I_P\sqcup J_P$.
According to Corollary \ref{cor;10.9.30.22},
they satisfy the conditions (P1) and (P2),
and the claim in Lemma \ref{lem;10.9.30.21}.

By using the uniqueness 
and Lemma \ref{lem;10.9.30.1},
we obtain an $\nbigr_X$-module $\nbigm[\ast I!J]$
by gluing $\nbigmzero_P[\ast I_P!J_P]$
for varied $(\lambda_0,P)\in\cnum_{\lambda}\times D$,
where $I_P\sqcup J_P=\ellsitabar$
is the induced decomposition
induced by $I\sqcup J=\Lambda$.
By construction,
$\nbigm[\ast I!J]$ $(I\sqcup J)$
satisfy the conditions (P1) and (P2),
and the claim in Lemma \ref{lem;10.9.30.21}.
Thus, the proof of Proposition \ref{prop;10.9.30.20}
and Lemma \ref{lem;10.9.30.21}
are finished.
\hfill\qed

\subsection{Ramified covering}

We give a remark on the functoriality
with respect to a ramified covering.
Let $X=\Delta^n$
and $D=\bigcup_{i=1}^{\ell}\{z_i=0\}$.
Let $\varphi:
 (X',D')\lrarr (X,D)$ be a ramified covering
along $(X,D)$.
Namely,
$X'=\Delta^n$,
$D'=\bigcup_{i=1}^{\ell}\{w_i=i\}$,
and 
$\varphi(w_1,\ldots,w_n)=
(w_1^{a_1},\ldots,w_{\ell}^{a_{\ell}},w_{\ell+1},\ldots,w_n)$.
Let $\nbigm'$ be a good-KMS smooth $\nbigr$-module
on $(X',D')$.
We naturally obtain a good-KMS smooth $\nbigr$-module
$\varphi_{\ast}\nbigm'$ on $(X,D)$.
\begin{prop}
For any decomposition
$I\sqcup J=\ellsitabar$,
we naturally have
\begin{equation}
\label{eq;13.5.4.1}
(\varphi_{\ast}\nbigm')[\ast I!J]
\simeq
 \varphi_{\dagger}\bigl(
 \nbigm'[\ast I!J]\bigr).
\end{equation}
\end{prop}
\pf
By using the induction on the dimension
and Corollary \ref{cor;10.9.30.22},
we can check that 
the assumption in Lemma \ref{lem;11.1.17.21}
is satisfied for
$\nbigm'[\ast I!J]$ and $\varphi$
along $z_i$ $(i=1,\ldots,\ell)$.
Then, by the lemma and the characterization of
$(\varphi_{\ast}\nbigm')[\ast I!J]$,
we obtain (\ref{eq;13.5.4.1}).
\hfill\qed

\vspace{.1in}
Let $\nbigm$ be a good-KMS $\nbigr$-module
on $(X,D)$.
We obtain a good-$KMS$ $\nbigr$-module
$\varphi^{\ast}$ on $(X',D')$.
\begin{cor}
$\nbigm[\ast I!J]$
is a direct summand of
$\varphi_{\dagger}
 \bigl(
 \varphi^{\ast}\nbigm[\ast I!J]\bigr)$.
It is the invariant part
with respect to the action of
the Galois group of $\varphi$.
\hfill\qed
\end{cor}

\section{Strict specializability along monomial functions}
\label{subsection;11.4.3.13}

\subsection{Statement}
\label{subsection;11.1.31.10}

Let $X$ be a complex manifold.
Let $D$ be a simple normal crossing hypersurface
with the irreducible decomposition 
$D=\bigcup_{i\in\Lambda}D_i$.
Let $\nbigm$ be a good-KMS smooth
$\nbigr_X(\ast D)$-module.
For simplicity, we assume the following:
\begin{description}
\item[\bf (A)]
$\nbigm$ is equipped with a filtration $L$
in the category of smooth $\nbigr_{X(\ast D)}$-modules,
such that, around any $P\in D$,
$\Gr^L(\nbigm)$ is the canonical prolongment of
a good wild harmonic bundle.
(See \S11.1 of \cite{mochi7} for the canonical prolongment
of wild harmonic bundles.)
\end{description}

Let $\Lambda=I\sqcup J$ 
be a decomposition.
Let $g$ be a holomorphic function on $X$
such that 
$g^{-1}(0)=\bigcup_{i\in K}D_i\subset D$.
\begin{prop}
\label{prop;11.1.21.1}
Assume $K\subset I$ or $K\subset J$.
\begin{itemize}
\item
$\nbigm[\ast I!J]$ is strictly specializable along $g$.
\item
For $\star=!,\ast$,
there exist
$\nbigm[\ast I!J][\star g]$,
and we have
\[
 \nbigm[\ast I!J][!g]
\simeq
 \nbigm\bigl[\ast (I\setminus K)!(J\cup K)\bigr],
\quad
 \nbigm[\ast I!J][!g]
\simeq
 \nbigm\bigl[\ast (I\cup K)!(J\setminus K)\bigr].
\]
\end{itemize}
\end{prop}

\subsection{Refinement}

We give a refined claim in the local case.
Let $X=\Delta^n$
and $D=\bigcup_{i=1}^{\ell}\{z_i=0\}$.
Let $g$ be a monomial function 
$g=\vecz^{\vecp}$,
where $\vecp\in\seisuu_{>0}^{K}$
and $K\subset\ellsitabar$.
Let $i_{g}:X\lrarr X\times\cnum_t$.
Let $\nbigk$ be a small neighbourhood of
$\lambda_0$ in $\cnum_{\lambda}$.
We set $(\nbigxzero,\nbigdzero):=\nbigk\times(X,D)$.
Let $\nbigqzero_{\ast}\nbigm$ be
a good-KMS family of filtered
$\lambda$-flat bundles.
For a decomposition 
$\ellsitabar=I\sqcup J\sqcup K$,
let us consider 
$\iota_{g\dagger}\nbigm[\ast I!J\star K]
=\iota_{g\ast}\nbigm[\ast I!J\star K]
 \otimes\cnum[\deldel_t]$
for $\star=\ast,!$.
Let $\nbigr_{X,K}\subset
 \nbigr_{X}$ be
generated by $\deldel_i$ $(i\in K)$
over $\nbigo_{\nbigx}$.
The following proposition implies
Proposition \ref{prop;11.1.21.1}.

\begin{prop}
\mbox{{}}\label{prop;10.10.2.3}
Assume that $\nbigm$ is a
good-KMS $\nbigr_X(\ast D)$-module
satisfying the condition (A')
in {\rm\S\ref{subsection;11.1.31.10}}.
Then, the following holds.
\begin{itemize}
\item
$\iota_{g\dagger}\nbigm[\ast I!J\star K]$
are strictly specializable along $t$,
and we have
\[
 \iota_{g\dagger}\nbigm[\ast I!J\star K]
\simeq
\bigl(
 \iota_{g\dagger}\nbigm[\ast I!J\star K]
\bigr)[\star t].
\]
\item
The $V$-filtration $\Uzero$ of
$\iota_{g\dagger}\nbigm[\ast I!J\ast K]$
is given as follows.
For $b\leq 0$,
\begin{equation}
\label{eq;13.5.4.10}
 \Uzero_b\bigl(
 \iota_{g\dagger}\nbigm[\ast I!J\ast K]
 \bigr)
=\nbigr_{X,K}\Bigl(
\lefttop{K}\Vzero_{b\vecp}
 \nbigm[\ast I!J]\otimes 1
 \Bigr).
\end{equation}
For $b>0$,
we have
$\Uzero_{b}=\sum_{c,j}\deldel_t^j\Uzero_c$,
where $(c,j)$ runs through
$\real_{\leq 0}\times\seisuu_{\geq 0}$
satisfying $c+j\leq b$.
(See {\rm\S\ref{subsection;10.8.21.51}}
for the filtration $\lefttop{K}\Vzero$.)
\item
The $V$-filtration $\Uzero$ of
$\iota_{g\dagger}\nbigmzero[\ast I!J! K]$
is given as follows.
For $b<0$,
\begin{equation}
 \label{eq;13.5.4.11}
 \Uzero_b\bigl(
 \iota_{g\dagger}\nbigm[\ast I!J! K]
 \bigr)
=\nbigr_{X,K}\Bigl(
\lefttop{K}\Vzero_{b\vecp}
 \nbigm[\ast I!J]\otimes 1
 \Bigr).
\end{equation}
For $b\geq 0$,
we have
$\Uzero_{b}=\sum_{c,j}\deldel_t^j\Uzero_c$,
where $(c,j)$ runs through
$\real_{< 0}\times\seisuu_{\geq 0}$
satisfying $c+j\leq b$.
\end{itemize}
\end{prop}

\subsection{Preliminary}

Recall that
for an $\nbigr_X$-module $\nbign$ on 
$\nbigxzero$,
the push-forward $i_{g\dagger}\nbign$
is naturally isomorphic to
$i_{g\ast}\nbign[\deldel_t]$,
where the action of $\nbigr_{X\times\cnum_t}$
is given as follows:
\begin{equation}
\label{eq;13.5.4.2}
 a\cdot \bigl(u\otimes\deldel_t^j\bigr)
=au\otimes\deldel_t^j\,\,\,(a\in \nbigo_X),
\quad
 \deldel_i\bigl(u\otimes\deldel_t^j\bigr)
=(\deldel_iu)\otimes\deldel_t^j
-(\del_ig\cdot u)\otimes \deldel_t^{j+1}
\end{equation}
\begin{equation}
\label{eq;13.5.4.3}
 t\cdot(u\otimes\deldel_t^j)
=(g\cdot u)\otimes\deldel_t^j
-j\lambda u\otimes\deldel_t^{j-1},
\quad
 \deldel_t\bigl(u\otimes\deldel_t^j\bigr)
=u\otimes\deldel_t^{j+1}
\end{equation}
In particular, we have
\begin{equation}
 \label{eq;13.5.4.4}
 (p_i\deldel_t t+\deldel_iz_i)
 (u\otimes\deldel_t^j)=
 -p_iu\otimes j\lambda\deldel_t^j
+\bigl(\deldel_i(z_iu)\bigr)
 \otimes\deldel_t^j.
\end{equation}

Let $\lefttop{t}V_0\nbigr_{X\times\cnum_t}
\subset\nbigr_{X\times\cnum_t}$
be the sheaf of subalgebras generated by
$t\deldel_t$ and $\nbigr_{X}$
over $\nbigo_{\cnum_{\lambda}\times X\times\cnum_t}$.
\begin{lem}
$\Uzero_b(i_{g\dagger}\nbigm[\ast I!J\star K])$ are 
$\lefttop{t}V_0\nbigr_{X\times\cnum_t}$-coherent
modules.
We have
$\bigcup_{b\in\real}
\Uzero_b(i_{g\dagger}\nbigm[\ast I!J\star K])
=i_{g\dagger}\nbigm[\ast I!J\star K]$.
\end{lem}
\pf
Let us prove the first claim.
We have only to consider the cases that
$\Uzero_b$ are expressed as
(\ref{eq;13.5.4.10}) or (\ref{eq;13.5.4.11}).
By using the relation (\ref{eq;13.5.4.4}),
we can check that
$\Uzero_b(i_{g\dagger}\nbigm[\ast I!J\star K])$ are 
$\lefttop{t}V_0\nbigr_{X\times\cnum_t}$-modules.
Let $\veca(I,J)$ be as in \S\ref{subsection;11.4.3.12}.
We have a naturally defined morphisms:
\begin{equation}
 \label{eq;13.5.4.5}
 \nbigqzero_{\veca(I,J)+b\vecp+\vecdelta_K}\nbigm
\lrarr
 \lefttop{K}\Vzero_{b\vecp}\nbigm[\ast I!J]
\lrarr
 \Uzero_b(i_{g\dagger}\nbigm[\ast I!J\ast K])
\end{equation}
By using (\ref{eq;13.5.4.4}),
we can check that
the image of (\ref{eq;13.5.4.5})
generates 
$\Uzero_b(i_{g\dagger}\nbigm[\ast I!J\ast K])$
over $\lefttop{t}V_0\nbigr_{X\times\cnum_t}$.
Then, we can deduce the coherence.
(See the last argument in the proof of 
Proposition 12.3.3 of \cite{mochi7}, for example.)

Let us prove the second claim
in the case $\star=\ast$.
Put $\nbigp:=\bigcup_{b\in\real}\Uzero_b$.
By construction,
it is an $\nbigr_{X\times\cnum_t}$-module.
We have
$\nbigp\supset 
 \lefttop{K}\Vzero_0\nbigm[\ast I!J]\otimes 1$
by the assumption.
Suppose $u\otimes\deldel_t^j\in\nbigp$.
By the first in (\ref{eq;13.5.4.2}),
we have
$(\del_ig) u\otimes \deldel_t^j\in\nbigp$.
By the second formulas in
(\ref{eq;13.5.4.3}) and (\ref{eq;13.5.4.4}),
we obtain 
$(\deldel_iu)\otimes\deldel_t^j$.
Because
$\lefttop{K}\Vzero\nbigm[\ast I!J]$
generates $\nbigm[\ast I!J]$
over $\nbigr_{X,K}$,
we obtain that
$\nbigm[\ast I!J\ast K]\otimes 1
 \subset\nbigp$,
which implies
$\nbigp=i_{g\dagger}\nbigm[\ast I!J\ast K]$.
The case $\star=!$ can be argued similarly.
\hfill\qed

\vspace{.1in}

We set
$\KMS(\nbigm,i):=
 \bigl\{
 u\in\real\times\cnum\,\big|\,
 \lefttop{i}\nbigg_u(\nbigm)\neq 0
 \bigr\}$.
For $b\in\real$,
we set
\[
 \nbigk(b,\lambda_0)
:=\bigcup_{i\in K}
 \bigl\{
 v\in\real\times\cnum\,\big|\,
 p_iv\in \KMS(\nbigm,i),\,\,
 \paramap(\lambda_0,v)=b
 \bigr\}
\]
\begin{lem}
The filtration $\Uzero$ is monodromic.
Namely,
the induced endomorphism
$\prod_{u\in \nbigk(b,\lambda_0)}(-\deldel_tt+\eigenmap(\lambda,u))$
is nilpotent on
$\Uzero_b\big/\Uzero_{<b}$.
\end{lem}
\pf
(See \S16.1 of \cite{mochi2}.)
By the relation (\ref{eq;13.5.4.4}),
for $s\in\Vzero_{b\vecp}\nbigm[\ast I!J]$
and for $i\in K$,
we have
\[
 \bigl(
 -\deldel_tt+\eigenmap(\lambda,u)
\bigr)(s\otimes 1)
-p_i^{-1}
 \bigl(
 (-\deldel_i z_i+\eigenmap(\lambda,p_iu))s
 \bigr)\otimes 1
=p_i^{-1}
 \deldel_i\bigl(
 (z_is)\otimes 1
 \bigr)
\in \Uzero_{<b}.
\]
Then, we can check the claim easily.
\hfill\qed

\vspace{.1in}

By construction,
we have $t\cdot \Uzero_b=\Uzero_{b-1}$
for $b<0$,
and $\deldel_t:\Gr^{\Uzero}_{c}\lrarr
\Gr^{\Uzero}_{c+1}$ is surjective
for $c>-1$.
If $\star=\ast$,
we also have $t\cdot \Uzero_0=\Uzero_{-1}$.
Hence, we have only to prove that
(i) $\Gr^{\Uzero}_{b}$ are strict
 for $b<0$,
(ii) $\deldel_t:\Gr^{\Uzero}_{-1}\lrarr
 \Gr^{\Uzero}_0$ is injective
 in the case $\star=!$,
which we shall consider in the following.

\subsection{Regular and pure case}

Let us consider the regular and pure case,
i.e., $\nbigm$ comes from a tame harmonic bundle.
Let $I\sqcup J\sqcup K\subset \ellsitabar$.
For $b\in\real$ and $\vecc\in\real^I$,
we consider
\begin{multline}
 \lefttop{I}\Vzero_{\vecc}
 \Uzero_b\bigl(
 \iota_{g\dagger}\nbigm[!J]
 \bigr)
:=\nbigr_{X,K}\Bigl(
 \lefttop{I}\Vzero_{\vecc}
 \lefttop{K}\Vzero_{b\vecp}\nbigm[!J]\otimes 1
 \Bigr) \\
=\nbigr_{X,K}\cdot\left\{
 \Bigl[
 \nbigr_{X,J}\otimes_{V_0\nbigr_{X,J}}
 \lefttop{I}\Vzero_{\vecc}
 \lefttop{J}\Vzero_{<0}
 \lefttop{K}\Vzero_{b\vecp}(\nbigm)
 \Bigr]\otimes 1
 \right\} \\
=\nbigr_{X,J}\otimes_{V_0\nbigr_{X,J}}
 \left\{
 \nbigr_{X,K}\Bigl[
 \lefttop{I}\Vzero_{\vecc}\lefttop{J}\Vzero_{<0}
 \lefttop{K}\Vzero_{b\vecp}(\nbigm)\Bigr]
\otimes 1
 \right\}
\\
=:\nbigr_{X,J}\otimes_{V_0\nbigr_{X,J}}
 \Bigl(
 \lefttop{I}\Vzero_{\vecc}
 \lefttop{J}\Vzero_{<0}
 \Uzero_b(\iota_{g\dagger}\nbigm)
 \Bigr).
\end{multline}
Here,
$\lefttop{I}\Vzero_{\vecc}
 \lefttop{K}\Vzero_{b\vecp}\nbigm[!J]$
is as in \S\ref{subsection;10.8.21.51}.
We obtain the following isomorphism:
\[
 \lefttop{I}\Vzero_{\vecc}
 \Gr^{\Uzero}_b\bigl(
 \iota_{g\dagger}\nbigm[!J]
 \bigr)
\simeq
 \nbigr_{X,J}\otimes_{V_0\nbigr_{X,J}}
 \Bigl(
 \lefttop{I}\Vzero_{\vecc}
 \lefttop{J}\Vzero_{<0}
 \Gr^{\Uzero}_{b}(\iota_{g\dagger}\nbigm)
 \Bigr)=:\nbigm_J
\]
\begin{lem}
\label{lem;10.10.2.2}
$\nbigm_J$ is strict.
\end{lem}
\pf
We use an induction on $|J|$.
The claim in the case $|J|=0$,
$\vecc\in\real_{<0}^I$
and $b<0$ follows from
Corollary 16.45 of \cite{mochi2}.
(Note that the minimal extension of $\nbigm$
is studied in \S16 of \cite{mochi2},
which is the image of
$\nbigm[!\ellsitabar]\lrarr
 \nbigm[\ast \ellsitabar]$
in the terminology of this paper.)
By the isomorphism given by the 
multiplication of $t$ and $z_i$ $(i\in I)$,
we obtain the claim in the case $|J|=0$.

Let $J=J_0\sqcup\{j\}$.
By the assumption of the induction,
we have the strictness of
$\nbigr_{X,J_0}\otimes
 _{V_0\nbigr_{X,J_0}}
 \lefttop{I}\Vzero_{\vecc}
 \lefttop{J}\Vzero_{<0}
 \Gr^{\Uzero}_{b}\iota_{g\dagger}\nbigm$.
We have
\begin{multline}
 \nbigr_{X,J_0}\otimes_{V_0\nbigr_{X,J_0}}
 \Bigl(
 \lefttop{I}\Vzero_{\vecc}
 \lefttop{J_0}\Vzero_{<0}
 \lefttop{j}\Gr^{\Vzero}_d\Gr^{\Uzero}_b
 \bigl(\iota_{g\dagger}\nbigm\bigr)
 \Bigr)
\simeq \\
 \nbigr_{X,J_0}\otimes_{V_0\nbigr_{X,J_0}}
 \Bigl(
 \lefttop{I}\Vzero_{\vecc}
 \lefttop{J_0}\Vzero_{<0}
  \Gr^{\Uzero}_b
 \bigl(\iota_{g\dagger}
\lefttop{j}\Gr^{\Vzero}_d
\nbigm\bigr)
 \Bigr),
\end{multline}
which is strict by the assumption of the induction.

Let us consider the filtration $\lefttop{j}\Vzero$
of $\nbigm_J$
given as follows.
For $d<0$,
we put
\[
\lefttop{j}\Vzero_d\nbigm_{J}:=
\nbigr_{X,J_0}\otimes_{V_0\nbigr_{X,J_0}}
\lefttop{I}\Vzero_{\vecc}
\lefttop{j}\Vzero_{d}
\lefttop{J_0}\Vzero_{<0}
\Gr^{\Uzero}_b(\iota_{g\dagger}\nbigm)
\]
For $d\geq 0$, we put
$\lefttop{j}\Vzero_d(\nbigm_J):=
 \sum_{c,n}\deldel_j^n
 \bigl(\lefttop{j}\Vzero_c(\nbigm_J)\bigr)$,
where $(c,n)$ runs through
$\real_{<0}\times\seisuu_{\geq 0}$
such that $c+n\leq b$.
Note that
$\lefttop{j}\Vzero_{<0}\nbigm_J$ is strict,
and $\lefttop{j}\Gr^{\Vzero}_c(\nbigm_J)$ is strict
for $c<0$,
by the above consideration.
Let us consider the morphism
\begin{equation}
\label{eq;10.8.27.1}
 \lefttop{j}\Gr^{\Vzero}_{-1}(\nbigm_J)
 \stackrel{\deldel_j}{\lrarr}
 \lefttop{j}\Gr^{\Vzero}_0(\nbigm_J).
\end{equation}
\begin{lem}
\label{lem;10.10.1.20}
The specializations of {\rm(\ref{eq;10.8.27.1})}
to any generic $\lambda$ are isomorphisms.
\end{lem}
\pf
The surjectivity is clear by construction.
Let $\pi:X\lrarr D_I$ be the projection.
For $P\in D_I^{\circ}$,
we can naturally regard 
$\nbigmlambda_{J,P}:=
 \nbigm_J\otimes\nbigo_{\{\lambda\}\times\pi^{-1}(P)}$
as a $D$-module.
If $\lambda$ is generic,
we have 
$\nbigmlambda_{J,P}[!z_j]=\nbigmlambda_{J,P}$,
and 
the specialization of $\lefttop{j}\Vzero$
gives a $V$-filtration along $z_j$,
which can be checked by a direct computation.
Hence the specialization of
(\ref{eq;10.8.27.1}) at a generic $\lambda$
is injective.
\hfill\qed

\vspace{.1in}
By Lemma \ref{lem;10.10.1.20}
and the strictness of $\lefttop{j}\Gr^{\Vzero}_{-1}\nbigm_J$,
we obtain that (\ref{eq;10.8.27.1}) is 
injective, and hence an isomorphism.
In particular,
$\lefttop{j}\Gr^{\Vzero}_0(\nbigm_J)$ is strict.
For $b>-1$, the morphism
$\deldel_j:\lefttop{j}\Gr^{\Vzero}_b(\nbigm_J)
 \lrarr
 \lefttop{j}\Gr^{\Vzero}_{b+1}(\nbigm_J)$ 
is surjective,
and $z_j\deldel_j$ on
$\lefttop{j}\Gr_b^{\Vzero}(\nbigm_J)$ is injective.
Hence, $\deldel_j$ is an isomorphism.
We also have the strictness of
$\lefttop{j}\Gr^{\Vzero}_b(\nbigm_J)$.
Thus, we obtain the claim in the case of $J$,
and the proof of Lemma \ref{lem;10.10.2.2}
is finished.
\hfill\qed

\begin{cor}
Suppose that $\nbigm$ comes from a tame harmonic bundle.
\begin{itemize}
\item
For any $J\subset \ellsitabar\setminus K$,
$\Gr^{\Uzero}_b(\iota_{g\dagger}\nbigm[!J])$
is strict,
i.e.,
$\nbigm[!J]$ is strictly specializable
along $g$.
The $V$-filtration is given in the standard way
as in Proposition {\rm\ref{prop;10.10.2.3}}.
\item
For $I\sqcup J\subset\ellsitabar\setminus K$,
$\nbigm[\ast I!J]$ is strictly specializable along $g$.
The $V$-filtration is given in the standard way
as in Proposition {\rm\ref{prop;10.10.2.3}}.
\hfill\qed
\end{itemize}
\end{cor}

Let $I\sqcup J\sqcup K=\ellsitabar$,
and we prove the claim of Proposition 
\ref{prop;10.10.2.3}
in the case that $\nbigm$ comes from
a tame harmonic bundle.
Let us consider 
$\iota_{g\dagger}\nbigm[\ast I!J\ast K]$.
We have already known the strictness of
$\Gr^{\Uzero}_{b}$ for $b<0$.
Because
$t:\Gr^{\Uzero}_0\lrarr\Gr^{\Uzero}_{-1}$
is an isomorphism,
$\Gr^{\Uzero}_0$ is also strict.
By the standard argument,
we obtain the strictness of
$\Gr^{\Uzero}_b$ for $b>0$.
Hence, we obtain the claim for
$\iota_{g\dagger}\nbigm[\ast I!J\ast K]$.

Let us consider
$\iota_{g\dagger}\nbigm[\ast I!J!K]$.
For $b<0$,
we have already known that
$\Gr^{\Uzero}_b$ are strict.
Let us consider the specialization 
to any generic $\lambda$.
We have 
\[
 \bigl(\nbigmlambda[\ast I!J!K]\bigr)[!g]
=\nbigmlambda[\ast I!J!K],
\]
and the specialization of $\Uzero$
gives a $V$-filtration.
Hence, we obtain that
$\deldel_t:\Gr^{\Uzero}_{-1|\lambda}
\lrarr\Gr^{\Uzero}_{0|\lambda}$
are isomorphisms for generic $\lambda$.
Because $\Gr^{\Uzero}_{-1}$ is strict,
we obtain that
$\deldel_t:\Gr^{\Uzero}_{-1}
\lrarr\Gr^{\Uzero}_{0}$
is an isomorphism.
In particular,
$\Gr^{\Uzero}_0$ is strict.
Then, by the standard argument,
we obtain that 
$\Gr^{\Uzero}_b$ are strict for any $b$.
Hence, we obtain the claim for
$\nbigm[\ast I!J!K]$.

\subsection{Regular and filtered case}
\label{subsection;11.1.31.11}

Let $L$ be the filtration as in the condition
(A) in \S\ref{subsection;11.1.31.10}.
Let us consider the case 
that $\Gr^L\nbigm$ comes from
a tame harmonic bundle.
In this case, 
we obtain the claim of Proposition \ref{prop;10.10.2.3}
from the following general lemma
with an easy induction.
\begin{lem}
\label{lem;10.10.2.4}
Let $0\lrarr\nbign_1\lrarr\nbign_2\lrarr\nbign_3\lrarr 0$
be an exact sequence.
They are equipped with monodromic 
filtrations $\Uzero$,
which are preserved by morphisms.
Assume the following.
\begin{itemize}
\item
$\nbign_i$ $(i=1,3)$ are strictly specializable
with $\Uzero(\nbign_i)$.
\item
$\nbign_2\lrarr\nbign_3$
is strict with respect to $\Uzero$.
\end{itemize}
Then, $\nbign_2$ is also strictly specializable
with $\Uzero$.
\end{lem}
\pf
We have only to prove that
$\nbign_1\lrarr\nbign_2$ is strict with respect to
$\Uzero$.
Because the restriction of
$\Uzero(\nbign_2)$ to $\nbign_1$
is also monodromic,
we obtain
$\Uzero_b(\nbign_1)\supset
 \Uzero_b(\nbign_2)\cap\nbign_1$
by a general result
(see Lemma 14.23 of \cite{mochi2}.)
Thus, we obtain Lemma \ref{lem;10.10.2.4}.
\hfill\qed

\subsection{Good irregular case
with unique irregular value}

We consider a ramified covering
$\pi:(X',D')\lrarr (X,D)$
given by
\[
 \pi(z_1,\ldots,z_n)
=(z_1^{e_1},\ldots,z_k^{e_k},z_{k+1},\ldots,z_n).
\]
Let $\gminia$ be a meromorphic function on $(X',D')$
such that $\gminia=\vecz^{\vecm}\gminia_1$
for some holomorphic function $\gminia_1$
with $\gminia_1(O)\neq 0$
and $\vecm\in\seisuu_{<0}^k$.
Let $\nbigm'$ be a 
smooth good-KMS $\nbigr_{X(\ast D)}$-module
satisfying the condition in 
\S\ref{subsection;11.1.31.11}.
Let us consider the case that
$\nbigm$ is obtained as
$\nbigm'\otimes\pi_{\ast}\nbigl(\gminia)$.
By using the characterization of
$\nbigm[\ast I!J\star K]$
in Proposition \ref{prop;10.9.30.20},
we have a natural isomorphism
$\nbigm[\ast I!J\star K]\simeq
 \nbigm'[\ast I!J\star K]
 \otimes\pi_{\ast}\nbigl(\gminia)$.
Let us observe 
\begin{equation}
\label{eq;13.5.4.12} 
\lefttop{K}\Vzero_{b\vecp}\bigl(
 \nbigm[\ast I!J\star K]\bigr)
=\lefttop{K}\Vzero_{b\vecp}(\nbigm'[\ast I!J\star K])
\otimes
 \pi_{\ast}\nbigl(\gminia).
\end{equation}
Both of them are
$\lefttop{K}V_0\nbigr_X$-submodules,
and contains 
$\nbigq_{\veca(I,J)+b\vecp+\vecdelta_K}$.
The left hand side is generated by
$\nbigq_{\veca(I,J)+b\vecp+\vecdelta_K}$.
Hence, the right hand side contains the left hand side.
Let 
$\nbign\subset\nbigm'[\ast I!J\star K]$
be a $V_0\nbigr_X$-submodule.
Suppose that $\nbign\otimes 1$ is contained in
the left hand side.
By considering the action of $z_i\deldel_i$ $(1\leq i\leq k)$,
we obtain that
$\nbign\otimes\pi_{\ast}\nbigl(\gminia)$
is also contained in the left hand side.
Let $F_{\ast}(\lefttop{K}V_0\nbigr_X)$
be the filtration given by the order of differential operators.
We set
$\nbign_m:=
 F_m(\lefttop{K}V_0\nbigr_X)
\nbigq_{\veca(I,J)+b\vecp+\vecdelta_K}$
in $\nbigm'[\ast I!J\star K]$,
which are $V_0\nbigr_X$-submodules.
If $\nbign_m$ is contained in the left hand side,
we obtain that
$\nbign_m\otimes\pi_{\ast}\nbigl(\gminia)$
is also contained in the left hand side
as remarked above.
Then, by a formal computation,
we obtain that
$\nbign_{m+1}$ is contained in the left hand side.
Hence, by using an induction,
we obtain that
$\bigl(
\bigcup_m\nbign_m
\bigr)\otimes\pi_{\ast}\nbigl(\gminia)$
is contained in the left hand side.
Because 
$\bigcup_m\nbign_m
=\lefttop{K}\Vzero_{b\vecp}(\nbigm'[\ast I!J\star K])$,
we obtain (\ref{eq;13.5.4.12}).

Let $q:X\times\cnum_t\lrarr X$ be the projection.
We have a natural isomorphism
\[
 i_{g\dagger}\nbigm[\ast I!J\star K]
\simeq
  i_{g\dagger}\nbigm'[\ast I!J\star K]
 \bigr)
\otimes
 q^{\ast}\bigl(
\pi_{\ast}
 \nbigl(\gminia)\bigr). 
\]
Let us observe
\begin{equation}
 \label{eq;10.10.2.5}
 \Uzero_b\bigl(
 i_{g\dagger}\nbigm[\ast I!J\star K]
 \bigr)
\simeq
 \Uzero_b\bigl(
 i_{g\dagger}\nbigm'[\ast I!J\star K]
 \bigr)
\otimes
 q^{\ast}\bigl(
\pi_{\ast}
 \nbigl(\gminia)
\bigr)
\end{equation}
Both of them 
are $\lefttop{t}V_0\nbigr_{X\times\cnum_t}$-modules,
and contain
$i_{g\ast}\bigl(
 \lefttop{K}\Vzero_{b\vecp}(\nbigm[\ast I!J\star K])
 \bigr)$.
Because the left hand side is  generated by
$i_{g\ast}\bigl(
 \lefttop{K}\Vzero_{b\vecp}(\nbigm[\ast I!J\star K])
\bigr)$,
it is contained in the right hand side.
Let $\nbign\subset 
 \Uzero_b\bigl(
 i_{g\dagger}\nbigm'[\ast I!J\star K]
 \bigr)$
be a $V_0\nbigr_{X}$-submodule.
Suppose that $\nbign$ is contained in the left hand side.
By considering the action of $z_i\deldel_i$ $(i=1,\ldots,k)$,
we obtain that 
$\nbign\otimes q^{\ast}\pi_{\ast}(\nbigl)$
is contained in the left hand side.
Let $F_{\ast}(\nbigr_{X,K})$ be the filtration
given by the order of differential operators.
We set
$\nbign_m:=
 F_m(\nbigr_{X,K})\cdot
 i_{g\ast}\bigl(
 \lefttop{K}\Vzero_{b\vecp}(\nbigm'[\ast I!J\star K])\bigr)$
in
$ \Uzero_b\bigl(
 i_{g\dagger}\nbigm'[\ast I!J\star K]
 \bigr)$.
We have
$\nbign_0\otimes\pi^{\ast}\nbigl(\gminia)
\simeq
i_{g\ast}\bigl(
 \lefttop{K}\Vzero_{b\vecp}(\nbigm[\ast I!J\star K])
\bigr)$,
which is contained in the left hand side.
Then, by an easy induction,
we obtain that
$\nbign_m\otimes
 q^{\ast}\pi_{\ast}\nbigl(\gminia)$
is contained in the left hand side.
Because
$\bigcup_m\nbign_m=
\Uzero_b\bigl(
 i_{g\dagger}\nbigm'[\ast I!J\star K]
 \bigr)$,
we obtain (\ref{eq;10.10.2.5}).

\vspace{.1in}
The claim of Proposition \ref{prop;10.10.2.3} in this case
immediately follows from (\ref{eq;10.10.2.5}).

\subsection{End of the proof of
Proposition \ref{prop;10.10.2.3}}

By using the formal completion
as in \S12.4 of \cite{mochi7},
we obtain the claims of Proposition \ref{prop;10.10.2.3}
in the general case.
\hfill\qed

\section{Good-KMS smooth $\nbigr_{X(\ast D)}$-triple}

\label{subsection;11.4.3.14}

Let $X$ be a complex manifold
with a simple normal crossing hypersurface $D$.
Let $\nbigt$ be a smooth $\nbigr_{X(\ast D)}$-triple.
It is called good(-$KMS$)
if the underlying smooth 
$\nbigr_{X(\ast D)}$-modules 
$\nbigm_i$ are good(-KMS).
\index{good-$KMS$}
We use the other adjectives 
``unramifiedly good-$KMS$'',
``regular-$KMS$''
in similar meanings.
\index{unramifiedly good-KMS}
\index{regular-$KMS$}
For simplicity, we assume the following:
\begin{itemize}
\item
Let $X_P$ be any small neighbourhood of
$P\in X$,
with a ramified covering
$\varphi_P:(X'_P,D'_P)\lrarr (X_P,D\cap X_P)$
such that
$\varphi_P^{\ast}(\nbigt)$ is unramified.
Then,
$\Irr(\nbigt,P):=
 \Irr(\nbigm_1,P)\cup
 \Irr(\nbigm_2,P)$ is a good set of irregular values.
\end{itemize}

\subsection{Reduction with respect to Stokes structure}
\label{subsection;11.2.20.3}

Let $X=\Delta^n$ and 
$D=\bigcup_{i=1}^{\ell}\{z_i=0\}$.
Let $\nbigt$ be a good-KMS smooth $\nbigr_{X(\ast D)}$-triple.
Let $I\subset\ellsitabar$.
We have the induced 
unramifiedly good-KMS smooth $\nbigr_{X(\ast D)}$-modules
$\lefttop{I}\Gr^{\St}(\nbigm_i)$
as in \S\ref{subsection;13.5.5.30}.
Let us observe that
we have the induced pairing
$\lefttop{I}\Gr^{\St}(C)$
of $\lefttop{I}\Gr^{\St}(\nbigm_i)$.
We will shrink $X$ around the origin
in the following argument.

Suppose that $\nbigt$ is unramifiedly good-KMS.
Let $\pi:\Xtilde(D)\lrarr X$ be the real blow up.
The induced morphism
$\cnum_{\lambda}\times\Xtilde(D)
\lrarr\cnum_{\lambda}\times X$
is also denoted by $\pi$.
We have the induced pairing
\[
 \Ctilde:
\bigl(
 \pi^{\ast}\nbigm_1
\bigr)_{|\vecS\times \Xtilde(D)}
\times
\bigl(
 \sigma^{\ast}\pi^{\ast}\nbigm_2
\bigr)_{|\vecS\times\Xtilde(D)}
\lrarr
 \nbigc^{\infty\,\moderate D}_{\vecS\times\Xtilde(D)}.
\]
Let $\lambda_0\in\vecS$
and $Q\in \pi^{-1}(D)$.
We have the full Stokes filtration
$\nbigftilde^{(\lambda_0,Q)}$
(resp.
$\nbigftilde^{(-\lambda_0,Q)}$)
of $\pi^{\ast}\nbigm_1$
(resp. $\pi^{\ast}\nbigm_2$)
at $(\lambda_0,Q)$
(resp. $(-\lambda_0,Q)$).
We obtain the following lemma,
by considering the growth order.
\begin{lem}
The restriction of $\Ctilde$ to
$\nbigftilde^{(\lambda_0,Q)}_{\gminia}
\times
 \sigma^{\ast}\nbigftilde^{(-\lambda_0,Q)}_{\gminib}$
vanishes,
unless $\Re(\gminia/\lambda_0)-\Re(\gminib/\lambda_0)\leq 0$.
\hfill\qed
\end{lem}
By shrinking $X$  around the origin $O$,
we obtain the pairing
\[
 \lefttop{I}\Gr^{\St}(C):
 \pi^{\ast}\bigl(
 \lefttop{I}\Gr^{\St}(\nbigm_1)
 \bigr)_{|\vecS\times \Xtilde(D)}
\times
\sigma^{\ast}\pi^{\ast}\bigl(
 \lefttop{I}\Gr^{\St}(\nbigm_2)
 \bigr)_{|\vecS\times \Xtilde(D)}
\lrarr
\nbigc^{\infty\,\moderate D}_{\vecS\times \Xtilde(D)}
\]
We obtain the induced pairing
$\lefttop{I}\Gr^{\St}(C)$
of 
$\lefttop{I}\Gr^{\St}(\nbigm_i)$ $(i=1,2)$.
It is easy to observe that
$\lefttop{I}\Gr^{\St}(C)_{|\vecS\times(X\setminus D)}$
is extended to a pairing on
$\cnum_{\lambda}^{\ast}\times (X-D)$
which is holomorphic with respect to $\lambda$.
The tuple of
$\lefttop{I}\Gr^{\St}(\nbigm_i)$ $(i=1,2)$
with $\lefttop{I}\Gr^{\St}(C)$
is denoted by $\lefttop{I}\Gr^{\St}(\nbigt)$,
which is an unramifiedly good-KMS
smooth $\nbigr_{X(\ast D)}$-triple.
We have the natural grading
$\lefttop{I}\Gr^{\St}(\nbigt)
=\bigoplus_{\gminia\in\Irr(\nbigt,I)}
 \lefttop{I}\Gr^{\St}_{\gminia}(\nbigt)$,
and we put
$\lefttop{I}\Gr^{\St,\reg}(\nbigt)
=\lefttop{I}\Gr^{\St}_0(\nbigt)$
and
$\lefttop{I}\Gr^{\St,\irr}(\nbigt)
=\bigoplus_{\gminia\neq 0}
 \lefttop{I}\Gr^{\St}_{\gminia}(\nbigt)$.

\vspace{.1in}

Let us consider the case that $\nbigt$
is not necessarily unramified.
We take a ramified covering
$\varphi:(X',D')\lrarr (X,D)$
such that
$\varphi^{\ast}\nbigt$ is unramified.
By applying the above procedure
to each good summand,
we obtain the decomposition
$\lefttop{I}\Gr^{\St}(\nbigt')
=\lefttop{I}\Gr^{\St,\reg}(\nbigt')
\oplus
\lefttop{I}\Gr^{\St,\irr}(\nbigt')$,
on which the Galois group of the ramified covering
naturally acts.
As the descent,
we obtain a good-KMS smooth
$\nbigr_{X(\ast D)}$-triple
with the decomposition:
\[
 \lefttop{I}\Gr^{\St}(\nbigt)
=\lefttop{I}\Gr^{\St,\reg}(\nbigt)
\oplus
\lefttop{I}\Gr^{\St,\irr}(\nbigt)
\]

\subsection{Specialization}

Let $X$ and $D$ be as in \S\ref{subsection;11.2.20.3}.
Let $\nbigt$ be a good-KMS smooth $\nbigr_{X(\ast D)}$-triple.
Because $\nbigt$ is strictly specializable along $z_i$
as an $\nbigr_{X(\ast D)}$-triple,
we obtain an $\nbigr_{D_i(\ast \del D_i)}$-triple
$\lefttop{i}\psitilde_{u}(\nbigt):=
 \psitilde_{z_i,u}(\nbigt)$
(\S\ref{subsection;13.5.4.40}).

\begin{lem}
We have a natural isomorphism
$\lefttop{i}\psitilde_u\lefttop{I}\Gr^{\St}(\nbigt)
\simeq
 \lefttop{I}\Gr^{\St}\lefttop{i}\psitilde_u(\nbigt)$.
It is a good-KMS smooth $\nbigr_{D_i(\ast \del D_i)}$-triple.
\end{lem}
\pf
By using Lemma 22.11.2 of \cite{mochi7},
we can reduce the issue to
the unramified case.
We have natural isomorphisms
$\lefttop{i}\psitilde_u(\nbigm_j)
\simeq
 \lefttop{i}\psitilde_u\lefttop{I}\Gr^{\St}(\nbigm_j)$.
We have only to compare the induced pairings
$\lefttop{i}\psitilde_u(C)$
and 
$\lefttop{i}\psitilde_u\lefttop{I}\Gr^{\St}(C)$,
which can be done 
by the argument in Proposition 12.7.1 of \cite{mochi7}.
Let us prove the second claim.
We set $i^c:=\{i\}^c$.
We consider the issue 
on $X\setminus D(i^c)$ for a while.
By the regularity along $D_i$,
we obtain the pairing
$\lefttop{I}\Gr^{\St\,\reg}(C)$
of 
$\nbigm_{1|\cnum^{\ast}\times(X\setminus D(i^c))}$
and
$\sigma^{\ast}
 \nbigm_{2|\cnum^{\ast}\times(X\setminus D(i^c))}$
taking values
in the sheaf of $C^{\infty}$-functions of moderate growth
along $D_i$ on $\cnum^{\ast}\times (X\setminus D(i^c))$
which are holomorphic with respect to $\lambda$.
Hence,
we obtain that
$\lefttop{i}\psitilde_u(\nbigt)_{|D_i\setminus\del D_i}$
is a smooth $\nbigr_{D_i(\ast\del D_i)}$-triple.
The underlying $\nbigr_{D_i(\ast \del D_i)}$-modules
are good-KMS smooth.
We obtain the growth estimate around $\del D_i$
directly by construction,
or by using Proposition \ref{prop;11.2.20.4} below.
\hfill\qed

\begin{lem}
We have a natural isomorphism
$\lefttop{i}\psitilde_{u_i}
 \lefttop{j}\psitilde_{u_j}(\nbigt)
\simeq
 \lefttop{j}\psitilde_{u_j}
 \lefttop{i}\psitilde_{u_i}(\nbigt)$.
\end{lem}
\pf
By using the previous lemma,
we can reduce the issue to the case
that $\nbigt$ is regular along $D_i\cup D_j$.
We have only to compare
$\lefttop{i}\psitilde_{u_i}\lefttop{j}\psitilde_{u_j}C$
and $\lefttop{j}\psitilde_{u_j}\lefttop{i}\psitilde_{u_i}C$
around generic $\lambda_0$.
It can be done by a direct computation.
Indeed,
by using the regularity along $D_i\cup D_j$
and the genericity of $\lambda_0$,
we have an expression
\[
 C(m_1,\sigma^{\ast}m_2)
=\sum a_{u_i,u_j,k_i,k_j}(\lambda,z_1,\ldots,z_n)
 |z_i|^{\eigenmap(\lambda,u_i)} 
 |z_j|^{\eigenmap(\lambda,u_j)}
 (\log |z_i|^2)^{k_i}
 (\log |z_j|^2)^{k_j}
\]
Here, $a_{u_i,u_j,k_i,k_j}$
are $C^{\infty}$ with respect to $(z_i,z_j)$.
Then, both
$\lefttop{i}\psitilde_{u_i}\lefttop{j}\psitilde_{u_j}C([m_1],[m_2])$
and
$\lefttop{j}\psitilde_{u_j}\lefttop{i}\psitilde_{u_i}C([m_1],[m_2])$
are equal to
$\bigl(a_{-\vecdelta,-\vecdelta,0,0}\bigr)_{|(z_i,z_j)=(0,0)}$
up to a constant term.

Otherwise, we can also obtain the commutativity
by using the Beilinson construction.
(See \S\ref{subsection;13.5.6.1} below.)
\hfill\qed

\vspace{.1in}

For any $I\subset\ellsitabar$
and $\vecu\in(\real\times\cnum)^I$,
we obtain 
$\lefttop{I}\psitilde_{\vecu}(\nbigt):=
 \lefttop{i_1}\psitilde_{u_1}\circ
 \cdots
 \circ\lefttop{i_m}\psitilde_{u_m}(\nbigt)$.
\index{$\nbigr$-triple $\lefttop{I}\psitilde_{\vecu}(\nbigt)$}
It is equipped with nilpotent morphisms
$\nbign_j:\lefttop{I}\psitilde_{\vecu}(\nbigt)
\lrarr
 \lefttop{I}\psitilde_{\vecu}(\nbigt)
 \otimes\newTate(-1)$ $(j\in I)$.
We have a natural isomorphism
$\lefttop{I}\psitilde_{\vecu}(\nbigt)
\simeq
 \lefttop{I}\psitilde_{\vecu}
 \lefttop{I}\Gr^{\St}(\nbigt)$.
If $\nbigt$ is unramified,
for each $\gminia\in\Irr(\nbigt,I)$,
we also obtain 
a smooth good-KMS $\nbigr_{D_I(\ast \del D_I)}$-triple
$\lefttop{I}\psitilde_{\gminia,\vecu}(\nbigt)
:=\lefttop{I}\psitilde_{\vecu}
\bigl(\nbigt\otimes L(-\gminiatilde)\bigr)$.
It depends on the choice of
a lift $\gminiatilde\in M(X,D)$
of $\gminia\in M(X,D)/H(X)$.

\subsection{Canonical prolongations}
\index{canonical prolongations}

Let $X$ be a complex manifold
with a normal crossing hypersurface $D$.
Let $\nbigt=(\nbigm_1,\nbigm_2,C)$
 be good-$KMS$ smooth $\nbigr_{X(\ast D)}$-triple.
By applying the procedure in
\S\ref{subsection;11.1.15.20}
inductively,
we obtain a uniquely determined pairing
$C[\ast I!J]$ of
$\nbigm_1[!I\ast J]$ and
$\nbigm_2[\ast I!J]$.
Thus, we obtain an $\nbigr$-triple
\[
 \nbigt[\ast I!J]:=
 \bigl(
 \nbigm_1[!I\ast J],\,
 \nbigm_2[\ast I!J],\,
 C[\ast I!J]
 \bigr). 
\]
\index{$\nbigr$-triple $\nbigt[\ast I\bikkuri J]$}
We obtain the following lemma
by the property of 
the canonical prolongments.
\begin{lem}
\mbox{{}}
\begin{itemize}
\item
 $\nbigt[\ast I!J]$ is strictly specializable
 along $z_i$, and we have the following:
\[
 \nbigt[\ast I!J][\ast z_i]
\simeq
 \nbigt[\ast I_{\cup i}!J_{\setminus i}],
 \quad
 \nbigt[\ast I!J][! z_i]
\simeq
 \nbigt[\ast I_{\setminus i}!J_{\cup i}].
\]
\item
The following morphisms are isomorphisms:
\[
 \phi_{z_i}^{(0)}\nbigt[\ast I!J]
\lrarr
 \psi_{z_i}^{(0)}\nbigt[\ast I!J],
\quad
 (i\in I)
\]
\[
 \psi^{(1)}_{z_i}\nbigt[\ast I!J]
\lrarr
 \phi^{(0)}_{z_i}\nbigt[\ast I!J],
\quad
 (i\in J)
\]
\item
For $u\in\bigl(\real\times\cnum\bigr)
 \setminus(\seisuu_{\geq 0}\times\{0\})$,
we have a natural isomorphism
$\lefttop{i}\psitilde_u\bigl(
 \nbigt[\ast I!J]\bigr)
\simeq
 \bigl(
 \lefttop{i}\psitilde_u(\nbigt)
\bigr)
 [\ast I_{\setminus i}
 !J_{\setminus i}]$.
\hfill\qed
\end{itemize}
\end{lem}

\subsection{Variant of Beilinson functors}
\label{subsection;13.5.6.1}

Let $X=\Delta^n$,
$D_i:=\{z_i=0\}$
and $D=\bigcup_{i=1}^{\ell}D_i$.
We take $K\sqcup J\sqcup I\sqcup A\sqcup B=L
 \subset\ellsitabar$.
For any function $f:K\sqcup J\lrarr \{0,1\}$, 
we put $K_i(f):=f^{-1}(i)\cap K$ for $i=0,1$.
Let $\nbigt_I$ be good-$KMS$ smooth 
$\nbigr_{D_I(\ast \del D_I)}$-triple.
We set
$\II^{-\infty,a}_g:=\varinjlim\II^{-N,a}_g$.
For $\veca=(\veca_K,\veca_J)\in\real^{K\sqcup J}$,
we put
{\small
\begin{multline}
 \nbigc_{f,\veca}(J,K,\nbigt_I)[\ast A!B]:=
 \\
 \Bigl(
 \nbigt_I\otimes
 \bigotimes_{k\in K_0(f)}
 \II^{-\infty,a_{k}+1}_{z_k}
 \otimes
 \bigotimes_{k\in J\sqcup K_1(f)}
 \II^{-\infty,a_k}_{z_k}
 \Bigr)
 \bigl[!\bigl(f^{-1}(0)\sqcup B\bigr)
 \ast(K\sqcup J\sqcup A\setminus f^{-1}(0))\bigr]
\end{multline}
}
Let $\veczero$ and $\chi_i$
be functions $K\sqcup J\lrarr \{0,1\}$
given by
$\veczero(j)=0$ for any $j\in K\sqcup J$,
and 
$\chi_i(j)=1$ $(j=i)$
or $\chi_i(j)=0$ $(j\neq i)$.
Then, let 
$\Xi^{(\veca_K)}_K\psi^{(\veca_J)}_{J}\nbigt_I[\ast A!B]$
be the $\nbigr_X(\ast D(\ellsitabar\setminus L))$-triple
obtained as the kernel of the following morphism:
\begin{equation}
 \label{eq;13.5.6.10}
 \nbigc_{\veczero,\veca}(J,K,\nbigt_I)[\ast A!B]
\lrarr
 \bigoplus_{i\in K\sqcup J}
 \nbigc_{\chi_i,\veca}(J,K,\nbigt_I)[\ast A!B]
\end{equation}
Then,
$\Xi^{(\veca_K)}_K
 \psi^{(\veca_J)}_{J}\nbigt_I[\ast A!B]
 (\ast \del D_{I\sqcup J})$
are good-KMS smooth
$\nbigr_{D_{I\sqcup J}(\ast\del D_{I\sqcup J})}$-triple.

\begin{lem}
$\Xi^{(\veca_K)}_K
 \psi^{(\veca_J)}_{J}\nbigt_I[\ast A!B]$
is strict, and strictly specializable
along any $z_j$ $(j\in \ellsitabar)$.
Moreover, 
for $j\in \ellsitabar\setminus(I\sqcup J)$,
we have natural isomorphisms
\begin{equation}
\label{eq;10.8.24.2}
 \psitilde_{z_j,u}\Xi^{(\veca_K)}_K\psi^{(\veca_J)}_J
 \nbigt_I[\ast A!B]
 \simeq
 \Xi^{(\veca_{K\setminus j})}_{K\setminus j}
 \psi^{(\veca_J)}_J\psitilde_{z_j,u}
 \nbigt_I[\ast A_{\setminus j}!B_{\setminus j}]
\end{equation}
\begin{equation}
\label{eq;10.8.24.1}
 \Xi^{(a_j)}_j\Xi^{(\veca_K)}_{K}
 \psi_J^{(\veca_J)}\nbigt_I[\ast A!B]
\simeq
 \Xi^{(\veca_{Kj})}_{Kj}
 \psi_J^{(\veca_J)}\nbigt_I
 [\ast A_{\setminus j}!B_{\setminus j}]
\end{equation}
\end{lem}
\pf
We have only to consider the case
$I=J=\emptyset$.
We use an induction on $|K|$.
If $|K|=0$,
the claim follows from Proposition \ref{prop;11.1.21.1}.
We shall prove the claim in the case $|K|=m$,
by assuming the claim in the case $|K|<m$.

Let $j\in K$.
By construction,
we have
$\Xi^{(\veca_{K})}_{K}
 \nbigt[\ast A!B](\ast z_j)
\simeq
\Xi^{(\veca_{K\setminus j})}_{K\setminus j}
 \nbigt[\ast A!B]$.
By the assumption of the induction,
it is strictly specializable along $z_i$ $(i\in\ellsitabar)$.
We obtain (\ref{eq;10.8.24.1}) in this case by construction,
which implies 
$\Xi^{(\veca_{K})}_{K}
 \nbigt[\ast A!B]$
is strictly specializable along $z_j$.
After applying $\psitilde_{z_j,u}$ to (\ref{eq;13.5.6.10}),
the morphism is strict.
Hence, we obtain (\ref{eq;10.8.24.2}) in this case.

\vspace{.1in}

Let $j\in \ellsitabar\setminus K$.
Take any $k\in K$.
Recall that
$\Xi_K^{(\veca_K)}\nbigt[\ast A!B]$
is the kernel of the following morphism:
\[
 \varphi:
 \Pi_{k!}^{-N,a_{k}+1}
 \Xi^{(\veca_{K\setminus k})}_{K\setminus k}
 \bigl(\nbigt[\ast A!B]\bigr)
\lrarr
  \Pi_{k\ast}^{-N,a_{k}}
 \Xi^{(\veca_{K\setminus k})}_{K\setminus k}
 \bigl(\nbigt[\ast A!B]\bigr)
\]
We have the following commutative diagram:
\[
 \begin{CD}
 \lefttop{j}\psitilde_{u}\Bigl(
 \Pi_{k!}^{-N,a_{k}+1}
 \Xi^{(\veca_{K\setminus k})}_{K\setminus k}
 \bigl(\nbigt[\ast A!B]\bigr)\Bigr)
@>{\lefttop{j}\psitilde_u(\varphi)}>>
 \lefttop{j}\psitilde_{u}\Bigl(
  \Pi_{k\ast}^{-N,a_{k}}
 \Xi^{(\veca_{K\setminus k})}_{K\setminus k}
 \bigl(\nbigt[\ast A!B]\bigr)\Bigr)
 \\
 @V{\simeq}VV @V{\simeq}VV \\
 \Pi^{-N,a_{k}+1}_{k!}
 \Xi^{(\veca_{K\setminus k})}_{K\setminus k}
 \bigl(\lefttop{j}\psitilde_u\nbigt[\ast A!B]\bigr)
@>>>
 \Pi^{-N,a_{k}}_{k\ast}
 \Xi^{(\veca_{K\setminus k})}_{K\setminus k}
 \bigl(\lefttop{j}\psitilde_u\nbigt[\ast A!B]\bigr)
 \end{CD}
\]
Hence, the cokernel of $\lefttop{j}\psitilde_u(\varphi)$
is naturally isomorphic to
$\Xi^{(\veca_{K\setminus k})}_{K\setminus k}
 \lefttop{k}\psi_{-\vecdelta}^{(-N)}
 \lefttop{j}\psitilde_u\nbigt[\ast A!B]$,
which is strict by the assumption
of the induction.
We obtain that
$\Xi_K^{(\veca_K)}\nbigt[\ast A!B]$
is strictly specializable along $z_j$.
If $j\in\ellsitabar\setminus (I\sqcup J\sqcup K)$,
we have the following natural isomorphisms:
\[
 \lefttop{j}\psitilde_u
 \Xi^{(\veca_K)}_K\bigl(
 \nbigt[\ast A!B]
 \bigr)
\simeq
 \Xi_k^{(a_k)}
 \Xi_{K\setminus k}^{(\veca_{K\setminus k})}
 \lefttop{j}\psitilde_u\nbigt[\ast A!B]
\simeq
 \Xi^{(\veca_K)}_K\bigl(
 \lefttop{j}\psitilde_u
 \nbigt[\ast A!B]
 \bigr)
\]
Thus, we obtain (\ref{eq;10.8.24.2}).
We obtain (\ref{eq;10.8.24.1})
by construction.
Thus, the induction can proceed.
\hfill\qed

\subsection{Growth order and the compatibility of Stokes filtrations}

Although a condition is imposed
on sesqui-linear pairings
of smooth $\nbigr$-triples,
we do not have to impose the assumption
in the good case.
Namely, the following proposition holds.
\begin{prop}
\label{prop;11.2.20.4}
Let $\nbigt$ be an $\nbigr_{X(\ast D)}$-triple
such that the underlying $\nbigr_{X(\ast D)}$-modules
are smooth and good.
Then, $\nbigt$ is smooth.
\end{prop}
\pf
We have only to consider the case that
the underlying smooth 
$\nbigr_{X(\ast D)}$-modules $\nbigm_i$ $(i=1,2)$
are unramifiedly good.
We use the notation in
\S\ref{subsection;11.2.20.3}.
Let $C_0$ denote the restriction of $C$
to $\vecS\times(X\setminus D)$.
Let $\lambda_0\in\vecS$ and $Q\in\pi^{-1}(D)$.
Let $I(\lambda_0)\subset\vecS$
and $U_Q\subset\Xtilde(D)$ be neighbourhoods of $\lambda_0$
and $Q$, respectively.
We consider the restriction of $C_0$ to
$\bigl(I(\lambda_0)\times U_Q\bigr)
 \setminus\pi^{-1}(D)$,
denoted by $C_Q$.
We have only to prove that
the restriction of $C_Q$ to 
$\nbigftilde^{(\lambda_0,Q)}_{\gminia}
\otimes
 \sigma^{\ast}\nbigftilde^{(-\lambda_0,Q)}_{\gminib}$
is $0$, unless $\Re(\gminia/\lambda)-\Re(\gminib/\lambda)\leq 0$.
By considering the specialization to curves,
which are transversal with the smooth part of $D$,
we can reduce the issue to the one dimensional case.
Then, the claim follows from Lemma 12.6.10 in \cite{mochi7}.
\hfill\qed

\subsection{$\vecnbigi$-good-KMS smooth $\nbigr$-triples}

We have a refined notion.

\begin{df}
Let $\vecnbigi=(\nbigi_P\,|\,P\in D)$
be a good system of ramified irregular values on $(X,D)$.
(See {\rm\S\ref{section;13.5.6.100}} below.)
A good(-KMS) smooth $\nbigr_{X(\ast D)}$-module $\nbigm$
is called $\vecnbigi$-good(-KMS),
if $\Irr(\nbigm,P)\subset\nbigi_P$
for any $P\in D$.
A good(-KMS) smooth $\nbigr_{X(\ast D)}$-triple
is called $\vecnbigi$-good(-KMS),
if the underlying $\nbigr_{X(\ast D)}$-modules
are $\vecnbigi$-good(-KMS).
\index{$\vecnbigi$-good-KMS}
\index{$\vecnbigi$-good}
\hfill\qed
\end{df}
\index{$\vecnbigi$-good}
\index{$\vecnbigi$-good-KMS}
A direct sum of $\vecnbigi$-good(-KMS) smooth
$\nbigr_{X(\ast D)}$-triples
is $\vecnbigi$-good(-KMS).
Let $F:X'\lrarr X$ be a morphism of complex manifolds.
If $\nbigt$ is an $\vecnbigi$-good(-KMS) $\nbigr_{X(\ast D)}$-triple,
then $F^{\ast}\nbigt$ is an $F^{-1}(\vecnbigi)$-good(-KMS)
$\nbigr_{X(\ast D)}$-triple.
Let $P\in D_I$ be any point.
If we take a holomorphic coordinate around $P$,
we obtain the $\vecnbigi(I)_{|D_I}$-good(-KMS)
smooth $\nbigr_{D_I(\ast D_I)}$-triple
$\lefttop{I}\psitilde_{\vecu}(\nbigt)$.

\section{Gluing of 
good-KMS smooth $\nbigr$-triples on
 the intersections}
\label{subsection;10.10.1.10}

Let $X$ be an open subset of $\cnum^n$.
Let $D_i:=\{z_i=0\}\cap X$
and $D=\bigcup_{i=1}^{\ell}D_i$.
We shall introduce a procedure to glue 
good-KMS smooth $\nbigr$-triples
given on the intersections $D_I$
$(I\subset\ellsitabar)$.

\subsection{A category}
\label{subsection;10.9.30.41}

Let $\vecC(X,D)$ be the category of tuples:
\[
 \vecnbigt=\bigl(
 \nbigt_I,\,(I\subset\ellsitabar);\,
 f_{I,i}, g_{I,i}\,\,
 (I\subset \ellsitabar,\,i\in\ellsitabar\setminus I)
 \bigr) 
\]

\begin{itemize}
\item
 $\nbigt_I$ are good-KMS smooth 
 $\nbigr_{D_I(\ast \del D_I)}$-triples.
\item
 $f_{I,i}$ and $g_{I,i}$ are morphisms
of $\nbigr$-triples
\[
\begin{CD}
  \psi_i^{(1)}\nbigt_I
 @>{g_{I,i}}>>
 \nbigt_{Ii}
 @>{f_{I,i}}>>
 \psi_i^{(0)}\nbigt_I
\end{CD}
\]
such that $f_{I,i}\circ g_{I,i}$
is equal to the canonical morphism
$\psi_i^{(1)}\nbigt_I\lrarr
 \psi_i^{(0)}\nbigt_I$.
We impose the commutativity of
the following diagram
for $j,k\in\ellsitabar\setminus I$
with $j\neq k$:
\begin{equation}
 \label{eq;11.1.18.3}
 \begin{CD}
 \psi_j^{(1)}(\nbigt_{Ik})
 @>{\psi_j^{(1)}(f_{I,k})}>>
 \psi_k^{(0)}\psi_j^{(1)}
 (\nbigt_{I})
 \\ 
 @V{g_{Ik,j}}VV 
 @VV{\psi^{(0)}_k(g_{I,j})}V \\
 \nbigt_{Ikj} @>{f_{Ij,k}}>>
 \psi^{(0)}_k(\nbigt_{Ij})
 \end{CD}
\end{equation}
Here,  $Ij=I\cup\{j\}$,
$Ik=I\cup\{k\}$
and $Ikj=I\cup\{k,j\}$.
\end{itemize}
A morphism
$\vecnbigt^{(1)}\lrarr\vecnbigt^{(2)}$ 
in $\vecC(X,D)$ is a tuple of morphisms
$F_I:\nbigt^{(1)}_{I}\lrarr\nbigt^{(2)}_{I}$
such that the following diagram is commutative:
\[
 \begin{CD}
   \psi_i^{(1)}\nbigt^{(1)}_{I}
 @>{g^{(1)}_{I,i}}>>
 \nbigt^{(1)}_{Ii}
 @>{f^{(1)}_{I,i}}>>
 \psi_i^{(0)}\nbigt^{(1)}_{I} \\
 @V{\psi_i^{(1)}F_I}VV 
 @V{F_{Ii}}VV 
 @V{\psi_i^{(0)}F_I}VV \\
   \psi_i^{(1)}\nbigt^{(2)}_I
 @>{g^{(2)}_{I,i}}>>
 \nbigt^{(2)}_{Ii}
 @>{f^{(2)}_{I,i}}>>
 \psi_i^{(0)}\nbigt^{(2)}_I
 \end{CD}
\]

For $p\in\ellsitabar$
and $u\in\real\times\cnum$,
we have the functors
$\lefttop{p}\psitilde_u$ and $\phi_p$
from $\vecC(X,D)$ to $\vecC(D_p,\del D_p)$,
given as follows:
\[
 \lefttop{p}\psitilde_u(\vecnbigt):=
\Bigl(
 \lefttop{p}\psitilde_u(\nbigt_I),\,\,
 (I\subset\ellsitabar\setminus p);\,\,
 \lefttop{p}\psitilde_u(f_{I,j}),\,
 \lefttop{p}\psitilde_u(g_{I,j})\,\,
 \bigl(I\subset \ellsitabar\setminus p,\,\,
 j\in \ellsitabar\setminus Ip\bigr)
\Bigr)
\]
\[
 \phi_p(\vecnbigt):=
\Bigl(
 \nbigt_{Ip},\,\,
 (I\subset\ellsitabar\setminus p);\,\,
 f_{Ip,j},\,
 g_{Ip,j}\,\,
 \bigl(I\subset \ellsitabar\setminus p,\,\,
 j\in \ellsitabar\setminus Ip\bigr)
\Bigr)
\]

We shall construct a fully faithful functor $\Psi_X$
from $\vecC(X,D)$
to the category of $\nbigr_X$-triples.

\subsection{Construction of the functor}
\label{subsection;11.2.19.20}

For a finite subset $I$,
we take a $2|I|$-dimensional hermitian vector space 
$(V,h)$
with an orthonormal base $u_k,v_k$ $(i\in I)$.
For $\veci=(i_k|k\in I)\in \{0,1\}^I$
and $\vecj=(j_k|k\in I)\in\{0,1\}^I$,
let $E(\veci,\vecj)$ denote the one dimensional subspace of
$\bigwedge^{\bullet}V$
generated by 
$\bigwedge_{k\in I}u_k^{i_k}
\wedge
 \bigwedge_{k\in J}v_k^{j_k}$.
It is naturally equipped with the hermitian metric 
denoted by $h$.
We obtain a smooth $\nbigr$-triple
$\Etilde(\veci,\vecj):=
 \bigl(\nbigo_{\cnum_{\lambda}}\otimes
 E(\veci,\vecj),
 \nbigo_{\cnum_{\lambda}}\otimes E(\veci,\vecj),\,
 C_h
 \bigr)$,
where $C_h$ is the sesqui-linear pairing
induced by $h$.
If $i'_p=i_p+1$,
$i'_k=i_k$ $(k\neq p)$
and $j'_k=j_k$ for any $k$,
we have a map
$\veca^{(1)}_p:
 \Etilde(\veci,\vecj)
\lrarr
 \Etilde(\veci',\vecj')$
given by the inner product of $u_p$
on the first component 
and the exterior product of $u_p$
on the second component.
If $j'_p=j_p+1$,
$j'_k=j_k$ $(k\neq p)$
and $i'_k=i_k$ for any $k$,
we have a map
$\veca_p^{(2)}:
 \Etilde(\veci,\vecj)
\lrarr
 \Etilde(\veci',\vecj')$
given by the inner product of $v_p$
on the first component 
and the exterior product of $v_p$
on the second component.
We have
$\veca_{p_1}^{(\kappa_1)}\circ\veca_{p_2}^{(\kappa_2)}
=-\veca_{p_2}^{(\kappa_2)}\circ\veca_{p_1}^{(\kappa_1)}$.

Let $\Gamma$ be the category
given by the following commutative diagram:
\[
 \begin{CD}
 (0,0) @>{a}>> (0,1)\\
 @V{b}VV @V{c}VV \\
 (1,0) @>{d}>> (1,1)
 \end{CD}
\]
For a finite set $I$,
let $\Gamma^I$ denote the product of the categories
of the $I$-tuples of objects in $\Gamma$.
An object $\bigl((i_k,j_k)\,\big|\,k\in I\bigr)$
is denoted as the pair of $\veci=(i_k)$
and $\vecj=(j_k)$.
For such an object,
we set $|\veci|:=\sum i_k$ and $|\vecj|:=\sum j_k$.
Let  $F$ be a functor from $\Gamma^I$
to the category of $\nbigr_X$-triples.
For $n\in\seisuu$, we define
\[
 \bigl(
 \pi_I F
 \bigr)^{n}:=
 \bigoplus_{|\veci|+|\vecj|=n+|I|}
 F(\veci,\vecj)\otimes
 \Etilde(\veci,\vecj).
\]
The morphisms
$F(\veci,\vecj)\lrarr F(\veci',\vecj')$
and $a^{(\kappa)}_{p}$
naturally give a differential.
The complex is denoted by $\pi_IF$.

\vspace{.1in}
For $1\leq m\leq \ell$,
we set $\mbar=\{1,\ldots,m\}$.
Let $\Gamma^m:=\Gamma^{\mbar}$.
For $(\veci,\vecj)\in\Gamma^{m}$,
we put
\[
I(\veci,\vecj):=\bigl\{
 k\in\mbar\,\big|\,(i_k,j_k)=(0,1)
 \bigr\},
\quad
K(\veci,\vecj):=\bigl\{
 k\in\mbar\,\big|\,(i_k,j_k)=(1,0)
 \bigr\},
\]
\[
 J_1(\veci,\vecj):=\bigl\{
 k\in\mbar\,\big|\,(i_k,j_k)=(0,0)
 \bigr\},\quad
J_0(\veci,\vecj):=\bigl\{
 k\in\mbar\,\big|\,(i_k,j_k)=(1,1)
 \bigr\}.
\]
For $\vecnbigt\in\vecC(X,D)$
and $\veca_{\mbar}\in\seisuu^m$,
we define a functor
$\nbigq_X^{\mbar}(\vecnbigt,\veca_{\mbar})$
from $\Gamma^m$
to the category of
$\nbigr_{X(\ast D(\ellsitabar-\mbar))}$-triples
given as follows:
\[
 \nbigq_X^{\mbar}\bigl(
 \vecnbigt,(\veci,\vecj),\veca_{\mbar}\bigr)
=\Xi^{(\veca_{I(\veci,\vecj)})}_{I(\veci,\vecj)}
 \psi^{(\veca_{J_1(\veci,\vecj)}
 +\vecdelta_{J_1(\veci,\vecj)})}
 _{J_1(\veci,\vecj)}
 \psi^{(\veca_{J_0(\veci,\vecj)})}_{J_0(\veci,\vecj)}
 \bigl(
 \nbigt_{K(\veci,\vecj)}
 \bigr)
\]
Here, $a_L:=(a_i\,|\,i\in L)$.
The morphisms are naturally induced ones.
Then, we obtain a complex
of $\nbigr_{X(\ast D(\ellsitabar-\mbar))}$-triples
$\pi_{\mbar}
 \nbigq^{\mbar}_X(\vecnbigt,\veca_{\mbar})$.

\begin{lem}
\mbox{{}}\label{lem;10.9.30.40}
\begin{itemize}
\item
 $\nbigh^p\pi_{\mbar}
 \nbigq^{\mbar}_X(\vecnbigt,\veca_{\mbar})=0$
unless $p\neq 0$.
\item
$\nbigh^0\pi_{\mbar}
 \nbigq^{\mbar}_X(\vecnbigt,\veca_{\mbar})$
is strict, and strictly specializable
along any $z_i$ $(i\in\ellsitabar)$.
\item
For any $p\in\ellsitabar\setminus\mbar$
and $u\in\real\times\cnum$,
we have
\[
 \lefttop{p}\psitilde_u\nbigh^0
 \pi_{\mbar}\nbigq^{\mbar}_X(\vecnbigt,\veca_{\mbar})
\simeq
  \nbigh^0
 \pi_{\mbar}\nbigq^{\mbar}_{D_p}\bigl(
  \lefttop{p}\psitilde_u\vecnbigt,\veca_{\mbar}\bigr)
\]
For any $p\in\mbar$
and $u\in\real\times\cnum$,
we have
\[
 \lefttop{p}\psitilde_u\nbigh^0
 \pi_{\mbar}\nbigq^{\mbar}_X(\vecnbigt,\veca_{\mbar})
\simeq
  \nbigh^0
 \pi_{\underline{m}\setminus p}
 \nbigq^{\underline{m}\setminus p}_{D_p}\bigl(
  \lefttop{p}\psitilde_u\vecnbigt,\veca_{\mbar\setminus p}\bigr)
\]
\item
For any $p\in\mbar$,
we have
\[
  \phi_p^{(a_p)}\nbigh^0
 \pi_{\underline{m}}
 \nbigq^{\underline{m}}_X(\vecnbigt,\veca_{\mbar})
\simeq
  \nbigh^0
 \pi_{\underline{m}\setminus p}
 \nbigq^{\underline{m}\setminus p}_{D_p}\bigl(
  \phi_p^{(a_p)}\vecnbigt,\veca_{\mbar\setminus p}\bigr)
\]
\item
For any $p\in\ellsitabar$,
we have
\[
 \Xi^{(a_p)}_p
 \nbigh^0
 \pi_{\mbar}\nbigq^{\mbar}_X(\vecnbigt,\veca_{\mbar})
\simeq
 \nbigh^0 \Xi^{(a_p)}_p
 \pi_{\mbar}\nbigq^{\mbar}_X(\vecnbigt,\veca_{\mbar})
\]
\end{itemize}
\end{lem}
\pf
We use an induction on $m$.
If $m=0$,
the claim is clear.
We put 
$\nbigt_{m}:=
 \nbigh^0\pi_{\underline{m}}
 \nbigq^{\underline{m}}_X
 (\vecnbigt,\veca_{\underline{m}})$.
We consider the following complexes:
\[
 C^{\bullet}_{0,0}:=
 \psi_m^{(a_m+1)}
 \pi_{\underline{m-1}}
 \nbigq^{\underline{m-1}}_{X}
 (\vecnbigt,\veca_{\underline{m-1}}),
\quad
 C^{\bullet} _{0,1}:=
 \Xi^{(a_m)}_m\pi_{\underline{m-1}}
 \nbigq^{\underline{m-1}}_X
 (\vecnbigt,\veca_{\underline{m-1}}),
\]
\[
 C^{\bullet}_{1,0}:=
 \pi_{\underline{m-1}}
 \nbigq^{\underline{m-1}}_{D_m}
 \bigl(\phi^{(a_m)}_m\vecnbigt,
 \veca_{\underline{m-1}}\bigr),
\quad
 C^{\bullet} _{1,1}:=
 \psi^{(a_m)}_m
 \pi_{\underline{m-1}}
 \nbigq^{\underline{m-1}}_X
 (\vecnbigt,\veca_{\underline{m-1}}).
\]
The complex
$\pi_{\mbar}\nbigq^{\mbar}_X(\vecnbigt,\veca_{\mbar})$
is associated to the following
naturally obtained triple complex:
\begin{equation}
\label{eq;10.8.24.3}
\begin{CD}
 C^{\bullet}_{0,0}
 @>>>
 C^{\bullet}_{0,1} \\
 @VVV @VVV \\
C^{\bullet}_{1,0}
 @>>>
C^{\bullet}_{1,1}
\end{CD}
\end{equation}
Let us consider the double complex
obtained as $\nbigh^0$ of  (\ref{eq;10.8.24.3}).
\begin{lem}
\label{lem;10.9.30.30}
The complex
associated to the double complex
may have non-trivial cohomology
only in the degree $0$,
i.e.,
the first claim holds in the case $m$.
\end{lem}
\pf
By the assumption in the case $m-1$,
$\nbigt_{m-1}$
is strictly specializable along $z_m$.
We also have $\nbigh^p(C_{i,j})=0$
unless $p\neq 0$,
and we have natural isomorphisms:
\[
\nbigh^0(C_{0,0})
\simeq
 \psi_m^{(a_m+1)}
 \nbigt_{m-1},\,\,\,
\nbigh^0(C_{0,1})
\simeq
 \Xi_m^{(a_m)}\nbigt_{m-1},\,\,\,
\nbigh^0(C_{1,1})
\simeq
 \psi_m^{(a_m)}
 \nbigt_{m-1}.
\]
Hence,
$\rho_1:\nbigh^0(C^{\bullet}_{0,0})
\lrarr
 \nbigh^0(C^{\bullet}_{0,1})$ is injective,
and
$\rho_2:\nbigh^0(C^{\bullet}_{0,1})
\lrarr
 \nbigh^0(C^{\bullet}_{1,1})$ is surjective.
Then the claim of Lemma \ref{lem;10.9.30.30}
follows.
\hfill\qed

\vspace{.1in}

We have the exact sequence
{\small
\begin{equation}
\label{eq;10.9.30.31}
 0\lrarr\nbigh^0(C^{\bullet}_{1,0})
\lrarr
 \Bigl(
 \nbigh^0(C^{\bullet}_{1,0})
 \oplus
 \nbigh^0(C^{\bullet}_{0,1})
\Bigr)\Big/\nbigh^0(C_{0,0})
\lrarr
\nbigh^0(C^{\bullet}_{0,1})
 \big/\nbigh^0(C_{0,0})
\lrarr 0.
\end{equation}
}
We obtain that the middle term of
(\ref{eq;10.9.30.31}) is strict.
Because $\nbigt_m$ is contained
in the middle term,
$\nbigt_m$ is also strict.
Because $\nbigt_m(\ast z_m)=\nbigt_{m-1}$,
we have 
$\lefttop{m}\psitilde_u\nbigt_m
=\lefttop{m}\psitilde_u\nbigt_{m-1}$.
Hence, we obtain the desired formula
for $\lefttop{m}\psitilde_u\nbigt_m$
from the formula for 
$\lefttop{m}\psitilde_u\nbigt_{m-1}$.
We also have 
$\Xi_m^{(a_m)}\nbigt_m
=\Xi_m^{(a_m)}\nbigt_{m-1}$.
Because 
$\Xi_m^{(a_m)}\nbigh^0C_{i,j}=0$
unless $(i,j)\neq (0,1)$,
we obtain the desired formula for
$\Xi_m^{(a_m)}\nbigt_m$.
Moreover, we obtain
the formula for $\phi_m^{(a_m)}\nbigt_m$.
By exchanging the roles,
we obtain that, for $p\leq m$,
$\nbigt_m$ is strictly specializable
along $z_p$,
and the desired formulas for
$\Xi_p^{(a_m)}\nbigt_m$,
$\lefttop{p}\psitilde_u\nbigt_m$
and $\phi_p^{(a_m)}\nbigt_m$.

\vspace{.1in}

Let $p>m$.
We have
{\small
\[
 \lefttop{p}\psitilde_uC^{\bullet}_{0,0}
=\psi_m^{(a_m+1)}
 \pi_{\underline{m-1}}
 \nbigq^{\underline{m-1}}_{D_p}
 (\lefttop{p}\psitilde_u\vecnbigt,\veca_{\underline{m-1}}),
\quad
 \lefttop{p}\psitilde_uC^{\bullet} _{0,1}=
 \Xi^{(a_m)}_m\pi_{\underline{m-1}}
 \nbigq^{\underline{m-1}}_{D_p}
 (\lefttop{p}\psitilde_u\vecnbigt,\veca_{\underline{m-1}}),
\]
\[
 \lefttop{p}\psitilde_uC^{\bullet}_{1,0}=
 \pi_{\underline{m-1}}
 \nbigq^{\underline{m-1}}_{D_p}
 \bigl(
 \phi^{(a_m)}_m
 \lefttop{p}\psitilde_u\vecnbigt,
 \veca_{\underline{m-1}}\bigr),
\quad
 \lefttop{p}\psitilde_u
 C^{\bullet} _{1,1}=
 \psi^{(a_m)}_m
 \pi_{\underline{m-1}}
 \nbigq^{\underline{m-1}}_{D_p}
 (\lefttop{p}\psitilde_u\vecnbigt,\veca_{\underline{m-1}}).
\]
}
The complex 
$\lefttop{p}\psitilde_u
 \pi_{\mbar}\nbigq^{\mbar}_X(\vecnbigt,\veca_{\mbar})$
is associated to the triple complex:
\[
 \begin{CD}
 \lefttop{p}\psitilde_uC^{\bullet}_{0,0}
 @>>>
 \lefttop{p}\psitilde_uC^{\bullet}_{0,1}\\
 @VVV @VVV \\
  \lefttop{p}\psitilde_uC^{\bullet}_{1,0}
 @>>>
 \lefttop{p}\psitilde_uC^{\bullet}_{1,1}\\
 \end{CD}
\]
Hence, by using the previous result,
we obtain that
$\lefttop{p}\psitilde_u
 \pi_{\mbar}\nbigq^{\mbar}_X(\vecnbigt,\veca_{\mbar})$
may have non-trivial cohomology only in the degree $0$,
and the $0$-th cohomology sheaf is strict.
It implies that the morphisms in the complex
$\pi_{\mbar}\nbigq^{\mbar}_X
 (\vecnbigt,\veca_{\mbar})$
is strict with respect to the $V$-filtration
along $z_p$.
Therefore, we obtain that $\nbigt_m$
is strictly specializable along $z_p$,
and we obtain the formula for
$\Xi^{(a_p)}_p\nbigt_m$.
We also have the following natural isomorphisms
for $\lefttop{p}\psitilde\nbigt_m$:
\[
 \lefttop{p}\psitilde_u\nbigt_m
\simeq
 \nbigh^0
  \lefttop{p}\psitilde_u
 \pi_{\mbar}\nbigq^{\mbar}_X(\vecnbigt,\veca_{\mbar})
\simeq
  \nbigh^0
 \pi_{\mbar}\nbigq^{\mbar}_{D_p}(
   \lefttop{p}\psitilde_u\vecnbigt,\veca_{\mbar})
\]
Thus, the induction can proceed,
and the proof of Lemma \ref{lem;10.9.30.40}
is finished.
\hfill\qed

\vspace{.1in}
Thus, we obtain the functor
$\Psi_X:=
\nbigh^0\pi_{\ellsitabar}\nbigq^{\ellsitabar}_X$
from the category $\vecC(X,D)$
to the category of $\nbigr_X$-triples.
It is independent of the choice of
a hermitian vector space $V$ with an orthonormal basis,
up to canonical isomorphisms.
By Lemma \ref{lem;10.9.30.40},
we have the following natural isomorphisms:
\[
 \lefttop{p}\psitilde_u\Psi_X(\vecnbigt)
\simeq
 \Psi_{D_p}(\lefttop{p}\psitilde_u\vecnbigt),
\quad
 \phi_p^{(a)}\Psi_X(\vecnbigt)
\simeq
 \Psi_{D_p}(\phi_p^{(a)}\vecnbigt).
\]
We have the following natural isomorphisms
for $p\neq q$:
\[
 \lefttop{p}\psitilde_u
 \lefttop{q}\psitilde_u\Psi_X(\vecnbigt)
\simeq
 \lefttop{q}\psitilde_u
 \lefttop{p}\psitilde_u\Psi_X(\vecnbigt),
\quad
 \phi_p^{(a_p)}\phi_q^{(a_q)}\Psi_X(\vecnbigt)
\simeq
 \phi_q^{(a_q)}\phi_p^{(a_p)}\Psi_X(\vecnbigt).
\]
We denote 
$\phi_{i_1}^{(a_{i_1})}
 \circ\cdots\circ\phi_{i_m}^{(a_{i_m})}
 \Psi_X(\vecnbigt)$
by 
$\phi_I^{(\veca)}\Psi_X(\vecnbigt)$.

\begin{lem}
\label{lem;11.4.4.1}
The functor $\Psi_X$ is fully faithful.
\end{lem}
\pf
We use an induction on $\dim X$.
Let $\vecF=(F_I\,|\,I\subset \ellsitabar):
 \vecnbigt_1\lrarr \vecnbigt_2$
be a morphism in $\vecC_X$
such that $\Psi_X(\vecF)=0$.
Then, we have $F_{\emptyset}=0$.
For each $p\in \ellsitabar$,
we have $\phi_p(\vecF)=0$.
It implies $F_{Ip}=0$ for any 
$I\subset\ellsitabar\setminus p$
by the assumption of the induction.
Hence, we obtain $\vecF=0$.

Let $G:\Psi_X(\vecnbigt_1)\lrarr 
 \Psi_X(\vecnbigt_2)$
be a morphism of $\nbigr$-triples.
We have the induced morphism
$G_{\emptyset}:
 \nbigt_{1,\emptyset}
\lrarr
 \nbigt_{2,\emptyset}$.
By applying $\phi_I$,
we obtain
$G_I:\nbigt_{1,I}\lrarr
 \nbigt_{2,I}$.
By using an induction,
we can prove that
$\Psi_X\bigl((G_I)\bigr)=G$.
\hfill\qed

\vspace{.1in}

Let $X'\stackrel{\iota}{\subset} X$ be any open subset.
We put $D':=D\cap X'$.
We have a natural functor
$\iota^{\ast};\vecC(X,D)\lrarr \vecC(X',D')$
given by the restriction.
It is easy to observe 
that $\Psi_{X'}\bigl(\iota^{\ast}(\vecnbigt)\bigr)$
is naturally isomorphic to
$\iota^{\ast}\Psi_{X}(\vecnbigt)$.

\subsection{Dependence on 
the coordinate system}

The category $\vecC(X,D)$ and the functor $\Psi_{X}$
depend on the coordinate system $\vecz=(z_1,\ldots,z_n)$.
To emphasize it, we use the symbols
$\vecC(X,D,\vecz)$ and $\Psi_{X,\vecz}$.
We shall study their dependence on $\vecz$.

Let $w_i:=e^{\varphi_i}\cdot z_i$,
where $\varphi_i$ are holomorphic functions on $X$.
We obtain another coordinate system
$\vecw=(w_1,\ldots,w_n)$.
We obtain a category
$\vecC(X,D,\vecw)$
and a functor
$\Psi_{X,\vecw}$.
\begin{lem}
\label{lem;11.1.18.11}
We have an equivalence
$\Def_{\vecvarphi}:
 \vecC(X,D,\vecz)
\lrarr
 \vecC(X,D,\vecw)$ such that
$\Psi_{X,\vecz}
\simeq
 \Psi_{X,\vecw}\circ\Def_{\vecvarphi}$.
\end{lem}
\pf
Let $\vecnbigt\in\vecC(X,D,\vecz)$.
For each $I\subset\ellsitabar$,
$\nbigt_I$ is equipped with a tuple
$\vecnbign_I=(\nbign_i\,\big|\,i\in I)$
of morphisms 
$\nbign_i:\nbigt_I\lrarr\nbigt_I\otimes\Tate(-1)$
induced as follows:
\[
\begin{CD}
 \nbigt_I
@>{f_{I\setminus i,i}}>>
 \psi^{(0)}\nbigt_{I\setminus i}
 \simeq
 \psi^{(1)}\nbigt_{I\setminus i}
\otimes
 \newTate(-1)
@>{g_{I\setminus i,i}}>>
 \nbigt_I\otimes\newTate(-1)
\end{CD}
\]
Let $\vecvarphi_I:=(\varphi_i\,|\,i\in I)$.
We put 
$\nbigttilde_I:=
\Def_{\vecvarphi_I}(\nbigt_I,\vecnbign_{I})$.
We have natural isomorphisms
\[
 \psi_{w_i}^{(a)}(\nbigttilde_{I\setminus i})
\simeq
 \Def_{\vecvarphi_I}
 \psi_{z_i}^{(a)}(\nbigt_{I\setminus i}).
\]
Hence, we have the following naturally induced morphisms:
\[
\begin{CD}
 \psi_{w_i}^{(1)}(\nbigttilde_I)
@>{\ftilde_{I,i}}>>
 \nbigttilde_{Ii}
@>{\gtilde_{I,i}}>>
 \psi_{w_i}^{(0)}(\nbigttilde_I)
\end{CD}
\]
They satisfy the commutativity condition induced by 
(\ref{eq;11.1.18.3}).
We put
$\Def_{\vecvarphi}(\vecnbigt)
:=\bigl(
 \nbigttilde_{I}\,;\,\ftilde_{I,i},\gtilde_{I,i}
 \bigr)$,
and thus we obtain the functor
$\Def_{\vecvarphi}:
 \vecC(X,D,\vecz)\lrarr
 \vecC(X,D,\vecw)$.
By construction, we have
$\Psi_{X,\vecw}\circ\Def_{\vecvarphi}(\vecnbigt)
\simeq
 \Def_{\vecvarphi}\Psi_{X,\vecz}(\vecnbigt)$.
By using Lemma \ref{lem;11.1.18.10} inductively,
we obtain that it is naturally isomorphic to
$\Psi_{X,\vecz}(\vecnbigt)$.
\hfill\qed

\part{Mixed twistor $D$-modules}

\chapter{Preliminary for relative monodromy filtrations}
\label{section;11.4.9.11}

We explain basic ideas to control
the weight filtrations
used in \cite{saito2}.
In \S\ref{subsection;13.4.12.20},
we recall some basic facts 
on relative monodromy filtrations
by following \cite{kashiwara-mixed-Hodge},
\cite{steenbrink-zucker},
and especially \cite{saito2}.
In \S\ref{subsection;13.4.12.21},
we review how the weight filtration of
the nearby cycle sheaf induces
a weight filtration on the vanishing cycle sheaf.
The procedure is called the transfer of filtration
in this paper.
In \S\ref{subsection;13.4.12.22},
we consider the successive use of the procedure
in a nice situation.

\section{Relative monodromy filtrations}
\label{subsection;13.4.12.20}

\subsection{Definition and basic properties}
\label{subsection;10.11.5.10}

Let $\nbiga$ be an abelian category
with additive auto equivalences 
$\vecSigma=\bigl(
 \Sigma^{p,q}\,\big|\,p,q\in\seisuu\bigr)$ 
such that
$\Sigma^{p,q}\circ
 \Sigma^{r,s}
=\Sigma^{p+r,q+s}$.
Let $C\in\nbiga$.
Let $L$ be  a finite increasing 
exhaustive complete filtration of $C$
in $\nbiga$
indexed by $\seisuu$,
i.e.,
$L$ is a tuple of subobjects
$\bigl(L_k(C)\subset C\,\big|\,k\in\seisuu\bigr)$
such that (i) $L_k(C)\subset L_{k+1}(C)$,
(ii) $L_k(C)=C$ for any sufficiently large $k$,
(iii) $L_k(C)=0$ for any sufficiently small $k$.
Such a pair $(C,L)$ is called a filtered object 
in $\nbiga$.
A morphism of filtered objects
$F:\bigl(C,L(C)\bigl)\lrarr \bigl(C',L(C')\bigr)$
is defined to be a morphism $F:C\lrarr C'$
preserving $L$,
i.e., $F\cdot L_k(C)\subset L_k(C')$.
The category of filtered objects in $\nbiga$
is denoted by $\nbiga^{\fil}$.
\index{category $\nbiga^{\fil}$}
It is an additive category.
Let $(C,L)\in\nbiga^{\fil}$.
We consider two naturally induced filtrations
$L^{(i)}$ $(i=1,2)$ on $\Sigma^{p,q} C$
given by
$L^{(1)}_k\Sigma^{p,q}C=
 \Sigma^{p,q}(L_{k}C)$
and 
$L^{(2)}_k\Sigma^{p,q}C=
 \Sigma^{p,q}(L_{k+p+q}C)$.
The object
$(\Sigma^{p,q}C,L^{(1)})$
is denoted by $(\Sigma^{p,q}C,L)$,
and 
$(\Sigma^{p,q}C,L^{(2)})$
is denoted by $\Sigma^{p,q}(C,L)$.

For $(C,L)\in\nbiga^{\fil}$,
put $\Gr^L_k(C):=L_k(C)/L_{k-1}(C)\in\nbiga$.
We have a natural isomorphism
$\Gr^L_k\Sigma^{p,q} C\simeq
 \Sigma^{p,q}\Gr^L_kC$.
If $C$ is equipped with a filtration $W$,
$\Gr^L_k(C)$ has the induced filtration given by 
$W_m\Gr^L_k(C):=
 \Image\bigl(
 W_m\cap L_k(C)\lrarr \Gr^L_k(C)
 \bigr)$.
\index{filtered object $\Sigma^{p,q}(C,L)$}
\index{filtered object $(\Sigma^{p,q}C,L)$}

We put $\vecT:=\Sigma^{1,1}$.
Let $\nbiga^{\nil}$ be the category of 
objects $C\in\nbiga$ equipped with
a nilpotent endomorphism 
$N:C\lrarr \vecT^{-1} C$.
\index{category $\nbiga^{\nil}$}
A morphism $F:(C,N)\lrarr (C',N')$ in $\nbiga^{\nil}$
is a morphism $F:C\lrarr C'$ such that
$F\circ N=N'\circ F$.
The category $\nbiga^{\nil}$ is abelian.
We set $\nbiga^{\fil,\nil}:=(\nbiga^{\nil})^{\fil}$.
For $(C,L,N)\in\nbiga^{\fil,\nil}$,
let $\Gr^L_k(N)$ denote the induced nilpotent
endomorphism of $\Gr^L_k(C)$.
\index{category $\nbiga^{\fil,\nil}$}

Let us recall the notion of relative monodromy filtration.
\begin{df}
Let $(C,L,N)\in\nbiga^{\fil,\nil}$.
A filtration $W$ of $C$ in $\nbiga$
is called a relative monodromy filtration
of $N$ with respect to $L$,
if the following holds:
\begin{itemize}
\item
 $N\cdot W_i(C)\subset
 \vecT^{-1}W_{i-2}(C)$,
 i.e.,
 $N$ induces 
 $(C,W,L)\lrarr
 \bigl(\vecT^{-1}(C,W),L\bigr)$.
\item
The induced morphisms
 $\Gr^W\Gr^L(N)^i:
 \Gr^W_{k+i}\Gr^L_k(C)
\lrarr
 \vecT^{-i}
 \Gr^W_{k-i}\Gr^L_k(C)$
are isomorphisms.
\end{itemize}
It is often denoted by $M(N;L)$ in this paper.
In this situation,
we say that $(C,L,N)$ has a relative monodromy
filtration.
\hfill\qed
\end{df}
\index{relative monodromy filtration}

According to \S1.1 of \cite{saito2}
(see also \cite{steenbrink-zucker}),
a relative monodromy filtration is uniquely determined
by Deligne's inductive formula, if it exists:
\begin{equation}
\label{eq;10.9.27.21}
 W_{-i+k}L_kC=W_{-i+k}L_{k-1}C
+N^i\vecT^{i}\bigl(W_{i+k}L_kC\bigr)
\quad\quad (i>0)
\end{equation}
\begin{equation}
\label{eq;10.9.27.22}
 W_{i+k}L_kC=
 \Ker\Bigl(
 N^{i+1}:L_kC\lrarr 
 \vecT^{-i}\bigl(
 L_kC/W_{-i-2+k}L_kC\bigr)
 \Bigr)
\quad\quad (i\geq 0)
\end{equation}
We should recall that 
a relative monodromy filtration
does not necessarily exist.
If there exists a $k\in\seisuu$ such that
$\Gr^L_m=0$ unless $m=k$,
then a relative monodromy filtration 
is the weight filtration of the nilpotent
morphism up to a shift of the degree,
and it always exists.
By Deligne's inductive formula,
relative monodromy filtrations are functorial
in the following sense.
\begin{lem}
\label{lem;10.7.23.1}
Assume that
$(C^{(i)},L,N^{(i)})\in \nbiga^{\fil,\nil}$ $(i=1,2)$
have relative monodromy filtrations
$M(N^{(i)};L)$.
Let $F:(C^{(1)},L,N^{(1)})\lrarr (C^{(2)},L,N^{(2)})$
be a morphism in $\nbiga^{\fil,\nil}$.
Then, $F\cdot M_k(N^{(1)};L)\subset M_k(N^{(2)};L)$.
\hfill\qed
\end{lem}

Let $\nbiga^{\RMF}$ be the full subcategory of
$\nbiga^{\fil,\nil}$,
whose objects have relative monodromy filtrations.
\index{category $\nbiga^{\RMF}$}
According to Lemma \ref{lem;10.7.23.1},
the correspondence
$(C,L,N)\longmapsto (C,M(N;L))$
gives a functor
$\Phi_1:\nbiga^{\RMF}\lrarr\nbiga^{\fil}$.
We have
$\Phi_1\circ\Sigma^{p,q}(C,L,N)
=\Sigma^{p,q}(C,M(N;L))$.

\subsection{Canonical decomposition}
\label{subsection;10.12.28.1}

\index{canonical decomposition}

Let us recall the notion of canonical decomposition
due to Kashiwara \cite{kashiwara-mixed-Hodge},
with a generalization in \cite{saito2},
which is one of the most fundamental 
in the study of relative monodromy filtrations.
Let $(C,L,N)\in\nbiga^{\RMF}$.
Let $M$ denote the relative monodromy filtration
for $(C,L,N)$.
We put $C^{(0)}:=\Gr^M(C)$.
It is equipped with a filtration induced by $L$,
which is also denoted by $L$.
Then, there exists a canonical splitting
of the induced filtration $L(C^{(0)})$ in $\nbiga$:
\begin{equation}
 \label{eq;10.7.19.1}
 C^{(0)}=\bigoplus C^{(0)}_i,
\quad
\mbox{\rm such that }
 L_jC^{(0)}=\bigoplus_{i\leq j}C^{(0)}_i
\end{equation}
(See \cite{kashiwara-mixed-Hodge}
and \cite{saito2} for the explicit construction
of the splitting.)
In particular, we have
$C^{(0)}_i\simeq
 \Gr^L_i\Gr^M(C)\simeq
 \Gr^M\Gr^L_i(C)$
in $\nbiga$.
It is functorial in the following sense,
which is clear by the construction
in \cite{kashiwara-mixed-Hodge}
and \cite{saito2}.
\begin{lem}
\label{lem;10.7.23.20}
Let 
$F:(C^{(1)},L,N)
\lrarr (C^{(2)},L,N)$
be a morphism in $\nbiga^{\RMF}$.
Let $M=M(N;L)$ on $C^{(i)}$.
Then, the induced morphism
$\Gr^M(F):
 \Gr^M(C^{(1)})\lrarr\Gr^M(C^{(2)})$
preserves the canonical splittings
{\rm(\ref{eq;10.7.19.1})}.
\hfill\qed
\end{lem}

\subsection{A criterion}

Let us recall a condition for the existence
of relative monodromy filtration
in \cite{steenbrink-zucker}.
Let $\nbiga$ be an abelian category.
Let $(C,L,N)\in \nbiga^{\nil,\fil}$
such that 
(i) $L_k(C)=L_{k'}(C)$ for any $k'\geq k$,
(ii) $(L_{k-1}(C),L,N)\in\nbiga^{\RMF}$.
Let $M\bigl(L_{k-1}(C)\bigr)$
denote the relative monodromy filtration
of $(L_{k-1}(C),L,N)$.
We have a naturally defined morphism
\begin{equation}
\label{eq;10.9.27.20}
\Ker\bigl(
 N^{\ell}:\vecT^{\ell}\Gr^L_k\lrarr\Gr^L_k
 \bigr)
\lrarr 
\frac{L_{k-1}}
{N^{\ell}\vecT^{\ell}L_{k-1}
 +M_{k-\ell-1}(L_{k-1}(C))}
\end{equation}

\begin{prop}[\cite{steenbrink-zucker},
see also \cite{saito2}]
\label{prop;10.11.5.1}
We have $(C,L,N)\in\nbiga^{\RMF}$
if and only if the morphisms
{\rm(\ref{eq;10.9.27.20})}
vanish for all integers $\ell>0$.
\hfill\qed
\end{prop}

\subsection{Functoriality
for tensor product and dual}

Let us consider the case that $\nbiga$
is the category $\Vect_K$
of finite dimensional vector spaces
over a ground field $K$,
and $\Sigma^{p,q}=\id$.
Recall that we have tensor product
and inner homomorphism in the category
$\Vect_K^{\fil,\nil}$,
given in the standard manner.

Let $(V^{(i)},L,N^{(i)})\in
 \Vect_K^{\fil,\nil}$ $(i=1,2)$.
The tensor product $V^{(1)}\otimes V^{(2)}$
has the induced endomorphism
$N_{V^{(1)}\otimes V^{(2)}}:=
N^{(1)}\otimes 1+1\otimes N^{(2)}$,
and the filtration of $V^{(1)}\otimes V^{(2)}$
given by
$L_k(V^{(1)}\otimes V^{(2)})=\sum_{p+q\leq k}
 L_p(V^{(1)})\otimes L_q(V^{(2)})$.
The tuple
$(V^{(1)}\otimes V^{(2)},L,N)$ is also denoted by
$(V^{(1)},L,N^{(1)})
 \otimes (V^{(2)},L,N^{(2)})$.
The space $\Hom(V^{(1)},V^{(2)})$
has the induced endomorphism
$N_{\Hom(V^{(1)},V^{(2)})}(f):=
N^{(2)}\circ f-f\circ N^{(1)}$
and the filtration $L$ given by
\[
 L_k\Hom(V^{(1)},V^{(2)}):=
 \bigl\{
 f\in\Hom(V^{(1)},V^{(2)})\,\big|\,
 f(L_{j}V^{(1)})\subset L_{j+k}V^{(2)}
\,\,\forall j
 \bigr\}.
\]
The tuple
$\bigl(
 \Hom(V^{(1)},V^{(2)}),L,N
 \bigr)$ is denoted by
$\Hom\bigl(\!
 (V^{(1)},L,N^{(1)}),(V^{(2)},L,N^{(2)})\!
 \bigr)$.
The following proposition is due to
Deligne, Steenbrink and Zucker
\cite{steenbrink-zucker}.
\begin{prop}[\cite{steenbrink-zucker}]
\label{prop;10.7.23.21}
$\Vect_K^{\RMF}$ is closed under
tensor product and inner homomorphism
in $\Vect_K^{\fil,\nil}$.
Moreover,
the relative monodromy filtrations for 
\[
 (V^{(1)},L,N^{(1)})\otimes (V^{(2)},L,N^{(2)})
\quad
\mbox{\rm and}\quad
\Hom\bigl(
 (V^{(1)},L,N^{(1)}),
 (V^{(2)},L,N^{(2)})
 \bigr) 
\]
are naturally induced by $M(N^{(i)};L)$ $(i=1,2)$.
\hfill\qed
\end{prop}
In particular,
if $(V,L,N)\in\Vect_K^{\RMF}$,
its dual is also an object in $\Vect_K^{\RMF}$,
and the relative monodromy filtration
is naturally induced by $M(N;L)$.

\section{Transfer of filtrations}
\label{subsection;13.4.12.21}

In some cases,
a relative monodromy filtration 
or a base filtration is inherited.
A fundamental general result is due to M. Saito \cite{saito2},
which we will review in \S\ref{subsection;11.4.3.20}.

\subsection{Gluing data}

We introduce some terminology
for our argument.
Let $\nbiga$ be an abelian category
with additive auto equivalences $\vecSigma$.
\begin{df}
A tuple $(C,C';u,v)$
of morphisms in $\nbiga$
\[
 \begin{CD}
 \Sigma^{1,0} C 
 @>{u}>> C' @>{v}>>
 \Sigma^{0,-1}C
 \end{CD}
\]
is called a gluing datum in $(\nbiga,\vecSigma)$.
\index{gluing datum}
Morphism of gluing datum
$(C_1,C_1';u_1,v_1)\lrarr
(C_2,C_2';u_2,v_2)$
is a commutative diagram:
\[
 \begin{CD}
 \Sigma^{1,0} C_1 @>{u_1}>> C_1' @>{v_1}>>
 \Sigma^{0,-1} C_1\\
 @V{F}VV @V{F'}VV @V{F}VV \\
 \Sigma^{1,0} C_2 @>{u_2}>> C_2' @>{v_2}>> 
 \Sigma^{0,-1}C_2
 \end{CD}
\]
The category of gluing data in 
$(\nbiga,\vecSigma)$
is denoted by $\Glu(\nbiga,\vecSigma)$.
\index{category $\Glu(\nbiga,\vecSigma)$}
\hfill\qed
\end{df}

An object of $\Glu(\nbiga,\vecSigma)^{\fil}$
is often denoted by
$(C,C';u,v;L)$,
which means an object $(C,C';u,v)$
with a filtration $L$ in $\Glu(\nbiga,\vecSigma)$.

\begin{df}
\mbox{{}}
\begin{itemize}
\item
$(C,C';u,v)\in\Glu(\nbiga,\vecSigma)$
is called $S$-decomposable,
if $C'=\Image u\oplus \Ker v$.
\index{$S$-decomposable}
\item
$(C,C';u,v;L)\in\Glu(\nbiga,\vecSigma)^{\fil}$
is called filtered $S$-decomposable,
if $\Gr^L(C,C';u,v)$ is $S$-decomposable.
\index{filtered $S$-decomposable}
\hfill\qed
\end{itemize}
\end{df}
If we say that $(C,C';u,v;L)$ is filtered $S$-decomposable
for given $(C,L),(C',L)\in\nbiga^{\fil}$
with morphisms $u:\Sigma^{1,0} C\lrarr C'$ and 
$v:C'\lrarr \Sigma^{0,-1}C$,
we implicitly imply that
$u$ and $v$ preserve the filtration $L$,
i.e.,
$u:(\Sigma^{1,0} C,L)\lrarr (C',L)$
and $v:(C',L)\lrarr (\Sigma^{0,-1}C,L)$.
We remark that this kind of condition appeared
in the study on mixed Hodge structure
in \cite{kashiwara-mixed-Hodge}
and \cite{saito2}.
The terminology ``$S$-decomposable''
is taken from \cite{sabbah2}.

\vspace{.1in}

Let $\Glu(\nbiga^{\fil},\vecSigma)$ denote
the category of gluing data in 
$\nbiga^{\fil}$,
where the action of $\Sigma^{p,q}$ on
$\nbiga^{\fil}$
is given by $(C,L)\longmapsto \Sigma^{p,q}(C,L)$.
An object of $\Glu(\nbiga^{\fil},\vecSigma)$ 
is often denoted by
$(C,C',\Ltilde;u,v)$,
i.e.,
a pair of objects
$(C,\Ltilde)$ and $(C',\Ltilde)$ 
with morphisms in $\nbiga^{\fil}$:
\[
 \begin{CD}
 \Sigma^{1,0}(C,\Ltilde)
 @>{u}>>
 (C',\Ltilde)
 @>{v}>>
 \Sigma^{0,-1}(C,\Ltilde)
 \end{CD}
\]
A morphism 
$(C_1,C_1',\Ltilde;u_1,v_1)\lrarr (C_2,C_2',\Ltilde;u_2,v_2)$
in $\Glu(\nbiga^{\fil},\vecSigma)$
is a commutative diagram:
\[
 \begin{CD}
 \Sigma^{1,0}(C_1,\Ltilde) @>{u_1}>> (C_1',\Ltilde) 
 @>{v_1}>> \Sigma^{0,-1}(C_1,\Ltilde) \\
 @VVV @VVV @VVV \\
 \Sigma^{1,0}(C_2,\Ltilde) @>{u_2}>> (C_2',\Ltilde) 
 @>{v_2}>> \Sigma^{0,-1}(C_2,\Ltilde)
 \end{CD}
\]

\subsection{Inheritance of relative monodromy filtration}
\label{subsection;10.9.24.1}

Let us recall a lemma due to Kashiwara.
Let $(C,C';u,v;L)\in\Glu(\nbiga,\vecSigma)^{\fil}$.
We put $N:=v\circ u$ and $N':=u\circ v$.

\begin{prop}[\cite{kashiwara-mixed-Hodge},
 \cite{saito2}]
\label{prop;10.10.3.2}
 Assume that
(i) the object $(C,C';u,v;L)$ is filtered $S$-decomposable,
(ii) $N'$ is nilpotent,
and $(C',L,N')\in\nbiga^{\RMF}$.
Then, we also have $(C,L,N)\in\nbiga^{\RMF}$.
Moreover,
$u$ and $v$ give the following morphisms
in $\nbiga^{\fil}$:
\[
\begin{CD}
\Sigma^{1,0}\bigl(C,M(N;L)\bigr)
@>{u}>>
 \bigl(C',M(N';L)\bigr)
@>{v}>>
\Sigma^{0,-1}\bigl(C,M(N;L)\bigr)
\end{CD}
\]
Namely, 
$u\cdot M_k(N;LC)\subset M_{k-1}(N';LC')$
and 
$v\cdot M_k(N';LC')\subset M_{k-1}(N;LC)$.
\end{prop}
\pf
We give only a remark.
In \cite{kashiwara-mixed-Hodge},
the case $\Sigma^{p,q}=\id$ is proved.
The key is Lemma 3.32 
in \cite{kashiwara-mixed-Hodge},
which was generalized in 
Corollary 1.7 \cite{saito2}.
Then, we can deduce the claim in the general case
by using the argument in \cite{kashiwara-mixed-Hodge}.
\hfill\qed

\vspace{.1in}

We have an obvious reformulation.
Let $\nbigc_1$ be the full subcategory of
$\Glu(\nbiga,\vecSigma)^{\fil}$
whose objects $(C,C';u,v;L)$ 
are filtered $S$-decomposable
and satisfy the following:
\begin{itemize}
\item
 We put $N:=v\circ u$ and $N':=u\circ v$,
 and then $(C,L,N)$ and $(C',L,N')$ are objects 
 in $\nbiga^{\RMF}$.
\item
 We put $\Ltilde(C):=M\big(N;L(C)\bigr)$
 and $\Ltilde(C'):=M\bigl(N';L(C')\bigr)$,
 and then
 we have
 $(C,C',\Ltilde;u,v)\in \Glu(\nbiga^{\fil},\vecSigma)$.
\end{itemize}
Let $\nbigc_2$ be the full subcategory of 
$\Glu(\nbiga,\vecSigma)^{\fil}$
whose objects $(C,C';u,v;L)$ are 
filtered $S$-decomposable
and satisfy the following:
\begin{itemize}
\item
We put $N':=u\circ v$,
and then $(C',L,N')$ is an object
in $\nbiga^{\RMF}$.
\end{itemize}
According to Proposition \ref{prop;10.10.3.2},
we have $\nbigc_1=\nbigc_2$.

\subsection{Transfer of filtration}
\label{subsection;11.4.3.20}
\index{transfer of filtration}

We shall recall a fundamental result due to Saito
in \cite{saito2}.
Let $(C,C';u,v)\in\Glu(\nbiga)$.
We set $N:=v\circ u$ and $N':=u\circ v$.
Assume that $C$ and $C'$ 
are equipped with two filtrations $L$ and $\Ltilde$
such that 
$(C,C';u,v;L)\in\Glu(\nbiga,\vecSigma)^{\fil}$
and $(C,C',\Ltilde;u,v)\in\Glu(\nbiga^{\fil},\vecSigma)$,
i.e., $u$ and $v$ give 
the following morphisms in $\nbiga^{\fil}$:
\begin{equation}
\label{eq;10.7.23.10}
\begin{CD}
 (\Sigma^{1,0} C,L)@>{u}>>  
 (C',L)@>{v}>> (\Sigma^{0,-1}C,L)
\end{CD}
\end{equation}
\begin{equation}
\begin{CD}
 \Sigma^{1,0}(C,\Ltilde)@>{u}>>(C',\Ltilde)
 @>{v}>>\Sigma^{0,-1}(C,\Ltilde)
\end{CD}
\end{equation}
\begin{prop}[Corollary 1.9 of \cite{saito2}]
\label{prop;10.7.19.2}
We assume
$(C,L,N)\in\nbiga^{\RMF}$
and $\Ltilde(C)=M(N;LC)$.
Then the following conditions are equivalent.
\begin{description}
\item[(A1)]
 $(C,C';u,v;L)$ is filtered $S$-decomposable,
 and we have
 $(C',L,N')\in\nbiga^{\RMF}$
 and $M(N';LC')=\Ltilde(C')$.
\item[(A2)]
 $L_kC'=
 u\bigl(\Sigma^{1,0} L_kC\bigr)
+\Bigl(
 v^{-1}(\Sigma^{0,-1}L_kC)
 \cap
 \Ltilde_kC'
 \Bigr)$ for each $k$.
\item[(A2')]
$L_kC'=v^{-1}\bigl(
 \Sigma^{0,-1}L_kC\bigr)
\cap\Bigl(
 u(\Sigma^{1,0} L_kC)+\Ltilde_kC'
 \Bigr)$ for each $k$.
\hfill\qed
\end{description}
\end{prop}

We have some immediate consequences of
Proposition \ref{prop;10.7.19.2}.
Let $(C,C',\Ltilde;u,v)\in\Glu(\nbiga^{\fil},\vecSigma)$.
Let $N:=v\circ u$ and $N':=u\circ v$.
Note that they are nilpotent.
Let $L(C)$ be a filtration of $C$
such that $(C,L,N)\in\nbiga^{\RMF}$
and $\Ltilde(C)=M\bigl(N;L(C)\bigr)$.

\begin{cor}
\label{cor;10.7.19.3}
Under the situation,
there exists a unique filtration
$L(C')$ such that 
(i) $(C,C';u,v;L)\in\Glu(\nbiga,\vecSigma)^{\fil}$ 
 is filtered $S$-decomposable,
(ii) $M\bigl(N';L(C')\bigr)=\Ltilde(C')$.
\end{cor}
\pf
The filtration $L(C')$ is given by
Saito's formula in Proposition \ref{prop;10.7.19.2}.
\hfill\qed

\begin{df}
The filtration $L(C')$
in Corollary {\rm\ref{cor;10.7.19.3}}
is called the transfer of  $L(C)$ by $(u,v)$
in this paper.
\hfill\qed
\end{df}

By Proposition \ref{prop;10.7.19.2},
the transfer is functorial.
It is reformulated as follows.
Let $\nbigc_3$ be the category of 
objects $\bigl((C,C',\Ltilde;u,v),L(C)\bigr)$,
where 
(i) $(C,C',\Ltilde;u,v)\in\Glu(\nbiga^{\fil},\vecSigma)$,
(ii) $L(C)$ is an exhaustive complete finite increasing filtration of $C$
in $\nbiga$, satisfying the following:
\begin{itemize}
\item
 We put $N:=v\circ u$,
 and then
 $(C,L,N)\in\nbiga^{\RMF}$
and $\Ltilde(C)=M\bigl(N;L(C)\bigr)$.
\end{itemize}
A morphism
$\bigl(
 (C_1,C_1',\Ltilde;u_1,v_1),L(C_1)\bigr)
\lrarr
\bigl(
 (C_2,C_2',\Ltilde;u_2,v_2),L(C_2)\bigr)$
in $\nbigc_3$ is 
a pair $(F,F')$ of morphisms
\[
 F:C_1\lrarr C_2,\quad\quad
 F':C'_1\lrarr C'_2
\]
such that
(i) $(F,F')$ gives a morphism
$(C_1,C_1',\Ltilde;u_1,v_1)\lrarr
 (C_2,C_2',\Ltilde;u_2,v_2)$
in $\Glu(\nbiga^{\fil},\vecSigma)$,
(ii) $F$ gives a morphism
$\bigl(C_1,L(C_1)\bigr)\lrarr
 \bigl(C_2,L(C_2)\bigr)$
in $\nbiga^{\fil}$.

Let $\nbigc_1$ be the category considered
in \S\ref{subsection;10.9.24.1}.
We have the functor
$\Psi:\nbigc_1\lrarr\nbigc_3$ given by
\[
 \Psi(C,C';u,v;L\bigr)
=\bigl((C,C',\Ltilde;u,v),L(C)\bigr),
\]
where 
$\Ltilde(C):=M\bigl(N;L(C)\bigr)$
and
$\Ltilde(C'):=M\bigl(N;L(C')\bigr)$.

\begin{cor}
\label{cor;10.11.5.2}
The functor $\Psi$
is an equivalence.
\end{cor}
\pf
Essential surjectivity is already stated in
Corollary \ref{cor;10.7.19.3}.
The functor is clearly faithful.
It is full by Saito's formula
in Proposition \ref{prop;10.7.19.2}.
\hfill\qed

\subsection{Special case}
\label{subsection;10.9.25.2}

As remarked in \cite{saito2},
Proposition \ref{prop;10.7.19.2}
can be regarded as a generalization
of the filtrations introduced 
in \cite{kashiwara-mixed-Hodge}
and \cite{steenbrink-zucker},
which we recall here.

Let $(C_1,L,N)\in\nbiga^{\RMF}$,
and let $W:=M(N;L)$.
We have the induced filtration
$N_{\ast}L$ on $C_1$, \index{filtration $N_{\ast}L$}
obtained as the transfer in the following case:
\[
 (C,L,\Ltilde)=\Sigma^{0,1}(C_1,L,W),
\quad
 (C',\Ltilde)=(C_1,W),
\quad
 u=N,
\quad
 v=\id.
\]
We also have the filtration $N_!L$ on $C_1$
obtained as the transfer \index{filtration $N_{\bikkuri}L$}
in the following case:
\[
 (C,L,\Ltilde)=\Sigma^{-1,0}(C_1,L,W),
\quad
 (C',\Ltilde)=(C_1,W),
\quad
 u=\id,
\quad
 v=N.
\]
Indeed, the explicit formulas are
given in 
\cite{kashiwara-mixed-Hodge}
and
\cite{steenbrink-zucker}
(we omit $\vecSigma$):
\[
 (N_{\ast}L)_k=
 NL_{k+1}+M_k(N;L)\cap L_k
=NL_{k+1}+M_k(N;L)\cap L_{k+1}
\]
\[
 (N_!L)_k=
 L_{k-1}+M_k(N;L)\cap N^{-1}L_{k-1}
=L_{k-1}+M_k(N;L)\cap N^{-1}L_{k-2}
\]

The morphism 
$N:(C_1,W)\lrarr \vecT^{-1}(C_1,W)$
induces 
$(C_1,N_!L)\lrarr
 \vecT^{-1}(C_1,N_{\ast}L)$,
which follows from the commutativity
of the following diagram:
\[
 \begin{CD}
 (C_1,W)
 @>{\id}>>
 (C_1,W)
 @>{N}>>
 \vecT^{-1}(C_1,W)
 \\
 @V{\id}VV @V{N}VV @V{\id}VV \\
 (C_1,W)
 @>{N}>>
 \vecT^{-1}(C_1,W)
 @>{\id}>>
 \vecT^{-1}(C_1,W)
 \end{CD}
\]

We define filtrations
$\Nhat_{\ast}L$ on $\Sigma^{0,-1}C_1$, 
and 
$\Nhat_!L$ on $\Sigma^{1,0}C_1$
as follows:
\index{filtration $\Nhat_{\ast}L$}
\index{filtration $\Nhat_{\bikkuri}L$}
\[
 (\Sigma^{0,-1}C_1,\Nhat_{\ast}L)
=\Sigma^{0,-1}(C_1,N_{\ast}L),
\quad\quad
 (\Sigma^{1,0}C_1,\Nhat_{!}L)
=\Sigma^{1,0}(C_1,N_{!}L)
\]

\subsection{Dual and tensor product}

Let us consider the case
that $\nbiga$ is the category of
finite dimensional vector spaces
with $\Sigma^{p,q}=\id$.
We have the functoriality of transfer for dual.
Let $(C,C';u,v;L)\in\Glu(\nbiga,\vecSigma)^{\fil}$.
Then, we have the induced morphisms
of dual filtered vector spaces
\[
 (C,L)^{\lor}\stackrel{v^{\lor}}{\lrarr}
 (C',L)^{\lor}\stackrel{u^{\lor}}{\lrarr} 
 (C,L)^{\lor}.
\]
We set
 $(C,C';u,v;L)^{\lor}:=
 \bigl(C^{\lor},(C^{\prime})^{\lor};
 v^{\lor},u^{\lor};L\bigr)
 \in\Glu(\nbiga,\vecSigma)^{\fil}$.
If $(C,C';u,v;L)$ is a filtered $S$-decomposable,
$(C,C^{\prime};u,v;L)^{\lor}$
is also filtered $S$-decomposable.
By the functoriality of relative monodromy filtration
with respect to dual,
the correspondence induces
contravariant functors $\lor_i$ on
the category $\nbigc_i$ $(i=1,2)$.

For an object
$(C,C',W;u,v)\in\Glu(\nbiga^{\fil},\vecSigma)$,
we have its dual
$(C,C',W;u,v)^{\lor}:=
 \bigl(
 C^{\lor},(C')^{\lor},W;v^{\lor},u^{\lor}
 \bigr)\in\Glu(\nbiga^{\fil},\vecSigma)$.
If $(C,C',W;u,v)$ is equipped with a filtration $L$ of $C$,
$(C,C',W;u,v)^{\lor}$ is also equipped with
an induced filtration $L$ of $C^{\lor}$.
By the functoriality of relative monodromy filtration
with respect to dual,
the correspondence of objects
$\bigl(
 (C,C',W;u,v),L(C)
 \bigr)$
and
$\bigl(
 (C,C',W;u,v)^{\lor},L(C^{\lor})
 \bigr)$
gives a contravariant functor 
$\lor_3$ on $\nbigc_3$.
The following lemma is clear
by construction.
\begin{lem}
\label{lem;10.9.25.1}
We have 
$\Psi\circ\lor_1=\lor_3\circ\Psi$
as functors from $\nbigc_1$
to $\nbigc_3$.
\hfill\qed
\end{lem}

Let us reword the claim of the lemma.
Let $\bigl((C,C',W;u,v),L(C)\bigr)\in\nbigc_3$.
We have a filtration $L(C')$ 
obtained as the transfer of $L(C)$
by $(u,v)$.
We also have 
$\bigl((C,C',W;u,v)^{\lor},L(C^{\lor})\bigr)
\in\nbigc_3$,
and we obtain a filtration $L(C^{\prime\lor})$
obtained as the transfer of $L(C^{\lor})$.
Then, according to Lemma \ref{lem;10.9.25.1},
$L(C^{\prime\lor})$
is the same as the filtration obtained
as the dual of $L(C')$.

\begin{rem}
Let $(C,L,N)\in\nbiga^{\RMF}$.
Let $N_!L$ and $N_{\ast}L$
be the filtrations in
{\rm\S\ref{subsection;10.9.25.2}}.
Then, 
as remarked in {\rm\cite{kashiwara-mixed-Hodge}},
$N_!L(C^{\lor})$ is obtained as the dual of
$N_{\ast}L(C)$, i.e.,
we have
$(N_!L)_k(C^{\lor})
=(N_{\ast}L)_{-k-1}(C)^{\bot}$.
It also follows from
the above functoriality of transfer
with respect to dual.
\hfill\qed
\end{rem}

Let $(U,L)\in\nbiga^{\fil}$.
For any $(C,C',L;u,v)\in\Glu(\nbiga,\vecSigma)^{\fil}$,
we have the naturally induced object
$(C,C',L;u,v)\otimes U\in\Glu(\nbiga,\vecSigma)^{\fil}$,
where the filtrations of $C\otimes U$
and $C'\otimes U$ are induced by $L$.
If $(C,C',L;u,v)$ is filtered $S$-decomposable,
$(C,C',L;u,v)\otimes U$ is also filtered
$S$-decomposable.
By the functoriality of relative monodromy filtration
with respect to tensor product,
it induces functors
$\otimes U$ on $\nbigc_i$ $(i=1,2)$.

For $(C,C',W;u,v)\in\Glu(\nbiga^{\fil},\vecSigma)$,
we have the naturally induced object
\[
 (C,C',W;u,v)\otimes U\in\Glu(\nbiga^{\fil},\vecSigma),
\]
where the filtrations of $C\otimes U$
and $C'\otimes U$ are induced by
$W(C)$, $W(C')$ and $L(U)$.
If $(C,C',W;u,v)$ is equipped with a filtration $L$
of $C$,
then $(C,C',W;u,v)\otimes U$
is equipped with a filtration $L$ of $C\otimes U$
induced by the filtrations $L$ of $C$ and $U$.
By the functoriality of relative monodromy filtrations,
the correspondence of
$\bigl((C,C',W;u,v),L(C)\bigr)$
and $\bigl((C,C',W;u,v)\otimes U,L(C\otimes U)\bigr)$
induces a functor 
$\otimes U$ on $\nbigc_3$.
The following lemma is clear by construction.
\begin{lem}
\label{lem;10.9.25.3}
We have
$\Psi\circ (\otimes U)
=(\otimes U)\circ \Psi$.
\hfill\qed
\end{lem}

Let us reword the lemma.
Let $(U,L)\in\nbiga^{\fil}$
and $\bigl((C,C',W;u,v),L(C)\bigr)\in\nbigc_3$.
We have a filtration $L(C')$ 
obtained as the transfer of $L(C)$
by $(u,v)$.
We have 
$\bigl((C,C',W;u,v)\otimes U,L(C\otimes U)\bigr)
\in\nbigc_3$,
and we obtain a filtration $L(C'\otimes U)$
obtained as the transfer of $L(C\otimes U)$.
Then, according to Lemma \ref{lem;10.9.25.3},
$L(C'\otimes U)$
is the same as the filtration 
induced by the filtrations $L(C')$ and $L(U)$.

\section{Pure and mixed objects}
\label{subsection;13.4.12.22}

This subsection is a preparation for
gluing of variation of admissible variation of
mixed twistor structures
(\S\ref{section;11.4.3.21}--\S\ref{section;11.4.3.5}).
The contents of this subsection are essentially 
contained in \cite{kashiwara-mixed-Hodge}
and \cite{saito2}.

\subsection{Setting}
\label{subsection;10.10.12.21}

Let $\nbiga$ be an abelian category
with a family of additive auto equivalences $\vecSigma$
as in \S\ref{subsection;10.11.5.10}.
For any finite set $\Lambda$,
we shall consider the category
$\nbiga(\Lambda)$ of objects $C$ in $\nbiga$
with a commuting $\Lambda$-tuple of nilpotent morphisms
$N_i:C\lrarr \vecT^{-1}C$ $(i\in\Lambda)$.
\index{category $\nbiga(\Lambda)$}
A morphism $(C_1,\vecN)\lrarr (C_2,\vecN)$
in $\nbiga(\Lambda)$
is a morphism $F:C_1\lrarr C_2$
such that $N_j\circ F=F\circ N_j$.
The category $\nbiga(\Lambda)$ is abelian.
It is equipped with auto equivalences given by
$(C,\vecN)\longmapsto
 (\Sigma^{p,q} C,\vecN)$,
which is also denoted by $\Sigma^{p,q}$.
When we are given
$(C,\vecN)\in\nbiga(\Lambda)$,
for any subset $\Lambda_0\subset\Lambda$,
we set $N(\Lambda_0):=\sum_{i\in\Lambda_0}N_i$,
and the weight filtration $M\bigl(N(\Lambda_0)\bigr)$
is denoted by $M(\Lambda_0)$.

Let $\Phi:\Lambda_1\lrarr\Lambda_2$ be a map.
For $(C_1,\vecN_{\Lambda_1})\in\nbiga(\Lambda_1)$,
we have a naturally induced object
$(C_1,\vecN_{\Lambda_2})\in\nbiga(\Lambda_2)$,
where 
$N_{j}:=\sum_{i\in\Phi^{-1}(j)}N_i$ 
for $j\in\Image\Phi$,
and $N_j:=0$ for $j\not\in\Image\Phi$.
It gives a functor
$\Phi_{\ast}:\nbiga(\Lambda_1)\lrarr\nbiga(\Lambda_2)$.
For $(C_2,\vecN_{\Lambda_2})\in\nbiga(\Lambda_2)$,
we have a naturally induced object
$(C_1,\vecN_{\Lambda_1})\in\nbiga(\Lambda_1)$,
where $N_j:=N_{\Phi(j)}$ for $j\in\Lambda_1$.
It gives a functor
$\Phi^{\ast}:\nbiga(\Lambda_2)
\lrarr\nbiga(\Lambda_1)$.
If $\Phi$ is a bijection,
both $\Phi_{\ast}$ and $\Phi^{\ast}$
are equivalences.

Assume that we are given 
saturated full subcategories
$\nbigp_w(\Lambda)\subset\nbiga(\Lambda)$
$(w\in\seisuu)$
for any finite sets $\Lambda$
with the following property:
\begin{description}
\item[\bf (P0)]
 $\nbigp_w(\Lambda)$ are abelian subcategories of
 $\nbiga(\Lambda)$,
 and they are semisimple.
 Any injection $\Phi:\Lambda_1\lrarr \Lambda_2$
 induces a functor
 $\Phi_{\ast}:\nbigp_w(\Lambda_1)\lrarr
 \nbigp_w(\Lambda_2)$,
 i.e., for $(C,\vecN)\in\nbigp_w(\Lambda_1)$,
 we have
 $\Phi_{\ast}(C,\vecN)\in\nbigp_w(\Lambda_2)$.
\item[\bf (P1)]
 $\Sigma^{p,q}$ induces an equivalence
 $\nbigp_{w}(\Lambda)\simeq\nbigp_{w-p-q}(\Lambda)$.
\item[\bf (P2)]
For $(C_i,\vecN)\in\nbigp_{w_i}(\Lambda)$ with
$w_1>w_2$,
we have
$\Hom_{\nbiga(\Lambda)}\!\bigl(\!
 (C_1,\vecN),\!(C_2,\vecN)
 \!\bigr)\!=\!0$.
\item[\bf (P3)]
Let $(C,\vecN)\in\nbigp_w(\Lambda)$.
\begin{description}
\item[\bf (P3.1)]
For any $\Lambda_2\subset\Lambda_1\subset\Lambda$,
there exists a relative monodromy filtration
 $M\bigl(N(\Lambda_1);M(\Lambda_2)\bigr)$,
and it is equal to 
$M(\Lambda_1)$.
\item[\bf (P3.2)]
 $\bigl(\Gr_{w_1}^{M(\Lambda_1)}(C),
 \vecN_{\Lambda_1^c}\bigr)$
 is an object in $\nbigp_{w+w_1}(\Lambda_1^c)$,
 where $\Lambda_1^c:=\Lambda\setminus\Lambda_1$.
\item[\bf (P3.3)]
 For $\bullet\in\Lambda$,
 $\bigl(\Image N_{\bullet},\vecN\bigr)$
 is an object in $\nbigp_{w+1}(\Lambda)$.
 Moreover, for $\ast\in \Lambda\setminus\bullet$,
 $\bigl(C,\Sigma^{1,0}\Image N_{\bullet};u,v;M(\ast)\bigr)$
 is filtered $S$-decomposable,
 where $u:\Sigma^{1,0} C\lrarr
 \Sigma^{1,0}\Image N_{\bullet}$ 
 and $v:\Sigma^{1,0}\Image N_{\bullet}\lrarr 
 \Sigma^{0,-1}C$ are induced by $N_{\bullet}$
 and the canonical morphism
 $\Image N_{\bullet}\lrarr \vecT^{-1}C$,
 respectively.
 Namely, 
\[
 \Sigma^{1,0} \Gr^{M(\ast)}C
\lrarr
 \Sigma^{1,0}\Gr^{M(\ast)}
 \Image N_{\bullet}
\lrarr
 \Sigma^{0,-1} \Gr^{M(\ast)}C
\]
is $S$-decomposable.
\end{description}
\end{description}

We have the category 
$\nbiga(\Lambda)^{\fil}$
of the filtered objects in $\nbiga(\Lambda)$,
on which $\Sigma^{p,q}$ naturally acts by
$\Sigma^{p,q}(C,L,\vecN)
=\bigl(\Sigma^{p,q}(C,L),\vecN\bigr)$.
Let $\nbigm'(\Lambda)
\subset\nbiga(\Lambda)$ be the saturated full
subcategory of the objects $(C,L,\vecN)$
such that 
$\Gr^L_w(C,\vecN)\in\nbigp_w(\Lambda)$.
\begin{lem}
\label{lem;10.12.28.1}
Let $F:(C_1,L,\vecN)\lrarr (C_2,L,\vecN)$
be a morphism in $\nbigm'(\Lambda)$.
\begin{itemize}
\item
 $F$ is strict with respect to $L$.
\item
 Let $(\Ker F,\vecN)$,
 $(\Image F,\vecN)$ and 
 $(\Cok F,\vecN)$
 be the kernel, the image and the cokernel 
 in $\nbiga(\Lambda)$.
 They are equipped with the naturally induced
 filtrations $L$,
 with which they are objects
 in $\nbigm'(\Lambda)$.
\item
In particular,
$\nbigm'(\Lambda)$ is abelian,
and the forgetful functor
$\nbigm'(\Lambda)\lrarr \nbiga(\Lambda)$ is exact.
\end{itemize}
\end{lem}
\pf
The first claim follows from {\bf P2}.
Then, the second claim follows from {\bf P0}
and the first claim.
\hfill\qed

\vspace{.1in}
Assume that we are given 
a saturated full subcategory
$\nbigm(\Lambda)\subset
 \nbigm'(\Lambda)$
with the following property:
\begin{description}
\item[\bf (M0)]
 Any injection $\Phi:\Lambda\lrarr\Lambda_1$
 induces a functor
 $\Phi_{\ast}:
 \nbigm(\Lambda)\lrarr\nbigm(\Lambda_1)$.
 The categories 
 $\nbigm(\Lambda)$ are abelian subcategories 
 of $\nbigm'(\Lambda)$.
 We naturally have
 $\nbigp_w(\Lambda)\subset\nbigm(\Lambda)$
 for any $w\in\seisuu$,
 i.e.,
 if $(C,L)\in\nbigm'(\Lambda)$ satisfies
 $\Gr^L_j(C)=0$ unless $j\neq w$,
 then $(C,L)\in\nbigm(\Lambda)$.
\item[\bf (M1)]
$(C,L,\vecN)\in\nbigm(\Lambda)$
if and only if
$\Sigma^{p,q}(C,L,\vecN)\in\nbigm(\Lambda)$.
\item[\bf (M2)]
 Let $(C,L,\vecN)\in\nbigm(\Lambda)$.
 \begin{description}
 \item[\bf (M2.1)]
 $\Gr^L_w(C,\vecN)\in\nbigp_w(\Lambda)$.
 \item[\bf (M2.2)]
 For any decomposition 
 $\Lambda=\Lambda_0\sqcup\Lambda_1$,
 there exists $M\bigl(N(\Lambda_1);L\bigr)=:M(\Lambda_1;L)$,
 and 
\[
 \res_{\Lambda_0}^{\Lambda}
 (C,L,\vecN):=
 (C,M(\Lambda_1;L),\vecN_{\Lambda_0})
\]
 is an object in $\nbigm(\Lambda_0)$.
 Note that $L$ gives a filtration of
 $\res^{\Lambda}_{\Lambda_0}(C,L,\vecN)$
 in $\nbigm(\Lambda_0)$.
\index{functor $\res_{\Lambda_0}^{\Lambda}$}
 \end{description}
 \item[\bf (M3)]
  Let $\bullet\in\Lambda$, and put 
 $\Lambda_0:=\Lambda\setminus\bullet$.
 Let us consider
 $(C,L,\vecN_{\Lambda})\in\nbigm(\Lambda)$
 and $(C',\Ltilde,\vecN'_{\Lambda_0})
 \in\nbigm(\Lambda_0)$ with morphisms
\[
\begin{CD}
 \Sigma^{1,0}
 \res^{\Lambda}_{\Lambda_0}
 (C,L,\vecN_{\Lambda})
 @>{u}>>
 (C',\Ltilde,\vecN'_{\Lambda_0})
 @>{v}>>
 \Sigma^{0,-1}
 \res^{\Lambda}_{\Lambda_0}
 (C,L,\vecN_{\Lambda})
\end{CD}
\]
in $\nbigm(\Lambda_0)$ such that
$v\circ u=N_{\bullet}$.
Let $L(C')$ be the filtration 
in the category of $\nbigm(\Lambda_0)$
obtained as the transfer of $L(C)$,
and put $N'_{\bullet}:=u\circ v$.
Then, 
$(C',L,\vecN'_{\Lambda})\in \nbigm(\Lambda)$.
\end{description}

The following lemma clearly follows from
{\bf M0}.
\begin{lem}
Let $(C,L,\vecN)\in\nbigm(\Lambda)$.
Then,
 $(L_jC,L,\vecN)\in\nbigm(\Lambda)$
for each $j\in \seisuu$.
\hfill\qed
\end{lem}

\begin{lem}
\label{lem;10.12.28.2}
Let $F:(C,L,\vecN)\lrarr (C',L,\vecN')$
be a morphism in $\nbigm(\Lambda)$.
Then, $F$ is strict with respect to
the filtrations $M(\Lambda_1;L)$
for any $\Lambda_1\subset\Lambda$.
\end{lem}
\pf
It follows from the conditions {\bf P2}
and {\bf M2}.
\hfill\qed

\begin{lem}
The functor 
$\res^{\Lambda}_{\Lambda_0}:
 \nbigm(\Lambda)\lrarr\nbigm(\Lambda_0)$
is exact.
\end{lem}
\pf
It follows from
{\bf M0}, Lemma \ref{lem;10.12.28.1}
and Lemma \ref{lem;10.12.28.2}.
\hfill\qed

\begin{lem}
\label{lem;10.11.5.12}
Let $(C,L,\vecN)\in\nbigm(\Lambda)$.
For any $\Lambda_2\subset\Lambda_1\subset\Lambda$,
we have
\begin{equation}
 \label{eq;13.4.11.1}
 M\bigl(
 N(\Lambda_1);M(\Lambda_2;L)
 \bigr)
=M(\Lambda_1;L). 
\end{equation}
\end{lem}
\pf
We use an induction on the length of the filtration $L$.
If $L$ is pure,
it follows from {\bf P3.1}.
Assume that 
(i) $L_k(C)=C$,
(ii) the claim holds for 
$(L_{k-1}(C),L,\vecN)$.
We have the exact sequence in
$\nbigm(\Lambda)$:
\[
 0\lrarr (L_{k-1}C,L,\vecN)
 \lrarr (C,L,\vecN)
 \lrarr (\Gr_k^LC,\vecN)
 \lrarr 0
\]
We obtain the exact sequence
in $\nbigm(\Lambda_2^c)$:
\[
 0\lrarr
 \res^{\Lambda}_{\Lambda_2^c}(L_{k-1}C,L,\vecN)
 \lrarr 
 \res^{\Lambda}_{\Lambda_2^c}(C,L,\vecN)
 \lrarr 
 \res^{\Lambda}_{\Lambda_2^c}(\Gr_k^LC,\vecN)
 \lrarr 0
\]
Because it is strict with respect to
$M\big(N(\Lambda_1);M(\Lambda_2;L)\bigr)$,
we obtain (\ref{eq;13.4.11.1}).
\hfill\qed

\subsection{A category $L\nbiga(\Lambda)$}

\index{category $L\nbiga(\Lambda)$}

We prepare a category
$L\nbiga(\Lambda)$ as follows.
An object of $L\nbiga(\Lambda)$ is a tuple
$\bigl(
 C_I\in\nbiga\,\big|\,I\subset \Lambda
 \bigr)$,
with morphisms for $I\subset J\subset\Lambda$
\[
 \begin{CD}
 \Sigma^{|J\setminus I|,0}C_I
 @>{g_{JI}}>>
 C_J
 @>{f_{IJ}}>>
 \Sigma^{0,-|J\setminus I|}C_I
 \end{CD}
\]
satisfying the compatibility conditions
for $I\subset J\subset K\subset\Lambda$:
\[
 f_{I\,J}\circ f_{J\,K}=f_{I\,K},
\quad\quad
 g_{K\,J}\circ g_{J\,I}=g_{K\,I},
\quad\quad
 g_{J,I\cap J}\circ f_{I\cap J,I}
=f_{J\,I\cup J}\circ g_{I\cup J,I}.
\]
Such a tuple
$\bigl(C_I;f_{IJ},g_{JI}\bigr)$
is often denoted by $\nbigt$.
A morphism $\nbigt_1\lrarr\nbigt_2$
is a tuple of morphisms
$F_I:C_{1,I}\lrarr C_{2,I}$
such that the following diagrams are commutative:
\[
 \begin{CD}
 \Sigma^{|J\setminus I|,0}C_{1,I}
 @>{g_{1,JI}}>>
 C_{1,J}
 @>{f_{1,IJ}}>>
 \Sigma^{0,-|J\setminus I|}C_{1,I} \\
 @V{F_I}VV @V{F_J}VV @V{F_I}VV \\
 \Sigma^{|J\setminus I|,0}C_{2,I}
 @>{g_{2,JI}}>>
 C_{2,J}
 @>{f_{2,IJ}}>>
 \Sigma^{0,-|J\setminus I|}C_{2,I}
 \end{CD}
\]
The category 
$L\nbiga(\Lambda)$ is abelian,
and equipped with naturally induced
auto equivalences $\Sigma^{p,q}$.

For any object $\nbigt=(C_I,f_{IJ},g_{JI})$,
each $C_I$ is equipped with a tuple of
morphisms
$N_i:C_I\lrarr \vecT^{-1}C_I$ $(i\in\Lambda)$
given as follows:
\[
 N_i=f_{I,I\cup\{i\}}\circ g_{I\cup\{i\},I}
 \,\,\,\,(i\not\in I),
 \quad\quad
 N_i=g_{I,I\setminus i}\circ f_{I\setminus i,I}
 \,\,\,\,(i\in I).
\]
The tuple is denoted by $\vecN_I$.

\subsection{$S$-decomposability
and strict support}

\index{$S$-decomposable}
\index{strict support}

An object $\nbigt\in L\nbiga(\Lambda)$
is called $S$-decomposable,
if the following holds:
\begin{itemize}
\item
 For $J=I\sqcup \{i\}$,
 we have
 $C_J=\Image g_{JI}\oplus\Ker f_{IJ}$,
i.e.,
the tuple
$(C_I,C_J,g_{JI},f_{IJ})$ is $S$-decomposable.
\end{itemize}
Let $\nbigt\in L\nbiga(\Lambda)$.
We say that $\nbigt$ has strict support $I$,
if the following holds:
\begin{itemize}
\item
We have
$C_{J}=0$ unless $J\supset I$.
The morphisms $g_{JI}$ (resp. $f_{IJ}$) are epimorphisms
(resp. monomorphisms)
for any $J\supset I$.
In particular, $\nbigt$ is $S$-decomposable.
\end{itemize}
The following lemma is clear.
\begin{lem}
\label{lem;10.10.12.3}
Let $\nbigt\!=\!\nbigt_1\oplus\nbigt_2$
be an object in $L\nbiga(\Lambda)$.
Then, \!$\nbigt$ is $S$-decomposable
if and only if $\nbigt_i$ $(i=1,2)$
are $S$-decomposable.
\hfill\qed
\end{lem}

\begin{lem}
\label{lem;10.10.12.1}
If $\nbigt_i\in L\nbiga(\Lambda)$ $(i=1,2)$
have strict supports $I_i$
with $I_1\neq I_2$,
then any morphism $\nbigt_1\lrarr\nbigt_2$
is $0$.
\end{lem}
\pf
Let $F:\nbigt_1\lrarr\nbigt_2$ be a morphism.
If $I_1\supsetneq I_2$,
we have the following:
\[
 \begin{CD}
 C_{1,I_1}@>>> C_{1,I_2}=0\\
 @V{F_{I_1}}VV @V{F_{I_2}}VV\\
 C_{2,I_1}@>{\rm mono}>> C_{2,I_2}
 \end{CD}
\]
Hence, we have $F_{I_1}=0$.
For any $J\supset I_1$,
we have the commutative diagram:
\[
 \begin{CD}
 C_{1,I_1} @>{\rm epi}>> C_{1,J}\\
 @V{F_{I_1}}VV @V{F_{J}}VV \\
 C_{2,I_1} @>{g_{J,I_1}}>> C_{2,J}
 \end{CD}
\]
We obtain $F_{J}=0$ for any $J\supset I$,
and hence $F=0$.
If $I_1\not\supset I_2$,
we have $C_{2,I_1}=0$
and hence $F_{I_1}:C_{1,I_1}\lrarr C_{2,I_1}$
is $0$,
which implies $F=0$ as above.
\hfill\qed

\begin{lem}
\label{lem;10.10.12.2}
Let $\nbigt\in L\nbiga(\Lambda)$ be $S$-decomposable.
Then, we have a unique decomposition
$\nbigt=\bigoplus_{I\subset\Lambda}
 \nbigt_I$
such that
each $\nbigt_I$ has strict support $I$.
\end{lem}
\pf
The uniqueness follows from
Lemma \ref{lem;10.10.12.1}.
As for the existence,
we use an induction on $|\Lambda|$.
If $|\Lambda|=0$,
the claim is clear.
For any $I\subset\Lambda$
and $i\in\Lambda\setminus I$,
let $Ii:=I\cup\{i\}$.

Let $\nbigt\in L\nbiga(\Lambda)$ be $S$-decomposable.
Assume that we have
$\Lambda_0\subset\Lambda$
such that $g_{Ii,I}$ is epimorphism for
any $i\in\Lambda_0$
and $I\subset\Lambda\setminus\{i\}$.
If $\Lambda=\Lambda_0$,
$\nbigt$ has strict support $\emptyset$.
Let us consider the case $\Lambda\neq\Lambda_0$.
Fix $j\in\Lambda\setminus\Lambda_0$.
We put
\[
 C'_I:=\left\{
 \begin{array}{ll}
 C_I & (j\not\in I)\\
 \Image g_{I,I\setminus j} & (j\in I),
 \end{array}
 \right.
\quad\quad
 C''_{I}:=\left\{
 \begin{array}{ll}
 0 & (j\not\in I)\\
 \Ker f_{I\setminus j,I} & (j\in I).
 \end{array}
 \right.
\]
We obtain a decomposition
$C_I=C_I'\oplus C_I''$ for any $I\subset\Lambda$,
which induces a decomposition
$\nbigt=\nbigt'\oplus\nbigt''$ in $L\nbiga(\Lambda)$.
By Lemma \ref{lem;10.10.12.3},
both $\nbigt'$ and $\nbigt''$ are $S$-decomposable.
We may apply the hypothesis of the induction
to $\nbigt''$.
We have the surjectivity of
$g_{Ii,I}$ for $\nbigt'$
if $i\in \Lambda_0\sqcup\{j\}$
and $i\not\in I$.
Hence, we obtain a decomposition
by strict supports by an easy induction.
\hfill\qed

\subsection{A category $L\nbiga(\Lambda_1,\Lambda_2)$}

\index{category $L\nbiga(\Lambda_1,\Lambda_2)$}

Let $\Lambda_i$ $(i=1,2)$ be finite sets.
We consider the category
$L\nbiga(\Lambda_1,\Lambda_2):=
 \bigl(L\nbiga(\Lambda_1)\bigr)(\Lambda_2)$,
i.e.,
let
$L\nbiga(\Lambda_1,\Lambda_2)$
be the category of objects $\nbigt$ in 
$L\nbiga(\Lambda_1)$
with a commuting $\Lambda_2$-tuple of
morphisms
$\nbign_j:\nbigt\lrarr\vecT^{-1}\nbigt$
$(j\in\Lambda_2)$.
It is an abelian category
equipped with naturally induced
auto equivalences $\Sigma^{p,q}$.
An object $(\nbigt,\vecnbign)\in
 L\nbiga(\Lambda_1,\Lambda_2)$
is called $S$-decomposable,
if $\nbigt$ is $S$-decomposable.
In the case,
$\vecnbign$ preserves the decomposition
in Lemma \ref{lem;10.10.12.2}.

Let $(\nbigt,\vecnbign)\in 
 L\nbiga(\Lambda_1,\Lambda_2)$.
Each underlying $C_I$ is equipped with
the morphisms $N_i$ $(i\in\Lambda_1)$
and $\nbign_j$ $(j\in\Lambda_2)$.
We set $\vecN_I^{(\Lambda_2)}:=
 \bigl(N_i\,\big|\,i\in\Lambda_1\bigr)
\sqcup
 \bigl(\nbign_j\,\big|\,j\in\Lambda_2\bigr)$.

\begin{rem}
$(\nbigt,\vecnbign)$ will often be denoted by
$\nbigt$,
if there is no risk of confusion.
\hfill\qed
\end{rem}

\subsection{Pure objects  in $L\nbiga(\Lambda_1,\Lambda_2)$}

An object 
$(\nbigt,\vecnbign)\in L\nbiga(\Lambda_1,\Lambda_2)$
is called pure of weight $w$,
if the following holds:
\begin{description}
\item[\bf (LP1)]
$(\nbigt,\vecnbign)$ is $S$-decomposable.
\item[\bf (LP2)]
Each $\bigl(C_I,\vecN_I^{(\Lambda_2)}\bigr)$
is an object in $\nbigp_w(\Lambda_1\sqcup\Lambda_2)$.
\end{description}
Let $P_wL\nbiga(\Lambda_1,\Lambda_2)\subset
 L\nbiga(\Lambda_1,\Lambda_2)$
 denote the full subcategory
of pure objects of weight $w$.
If $\Lambda_2=\emptyset$,
it is denoted by 
$P_wL\nbiga(\Lambda_1)$.

\begin{prop}
\label{prop;10.10.12.5}
The family $\bigl\{
 P_wL\nbiga(\Lambda_1,\Lambda_2)
 \bigr\}$
satisfies the conditions {\bf P0--P3}.
\end{prop}
\pf
The condition {\bf P2} for 
$P_wL\nbiga(\Lambda_1,\Lambda_2)$
follows from that for 
$\nbigp_w(\Lambda_1\sqcup\Lambda_2)$.
The condition {\bf P1} is clearly satisfied.
Any injection $\Lambda_2\lrarr\Lambda_2'$
clearly induces a functor
$P_wL\nbiga(\Lambda_1,\Lambda_2)
\lrarr
 P_wL\nbiga(\Lambda_1,\Lambda_2')$.
Let us prove that
$P_wL\nbiga(\Lambda_1,\Lambda_2)$ is abelian
and semisimple.
The following lemma is clear.
\begin{lem}
\label{lem;10.10.12.4}
Let $\nbigt_i\in
 L\nbiga(\Lambda_1,\Lambda_2)$
$(i=1,2)$.
Then,
$\nbigt_1\oplus
 \nbigt_2\in
 P_wL\nbiga(\Lambda_1,\Lambda_2)$
if and only if
$\nbigt_i\in
 P_wL\nbiga(\Lambda_1,\Lambda_2)$
for $i=1,2$.
\hfill\qed
\end{lem}

We obtain the following lemma 
from Lemma \ref{lem;10.10.12.2}
and Lemma \ref{lem;10.10.12.4}.
\begin{lem}
For any $\nbigt\in
 P_wL\nbiga(\Lambda_1,\Lambda_2)$.
we have the decomposition
$\nbigt=\bigoplus_{I\subset\Lambda_1}
 \nbigt_I$
in $P_wL\nbiga(\Lambda_1,\Lambda_2)$,
where each $\nbigt_I$ has strict support $I$.
It is uniquely determined.
\hfill\qed
\end{lem}

\begin{lem}
Let $\nbigt,\nbigt'\in P_wL\nbiga(\Lambda_1,\Lambda_2)$.
Assume that they have strict support $I$.
Let $F:\nbigt\lrarr\nbigt'$ be a morphism
in $P_wL\nbiga(\Lambda_1,\Lambda_2)$.
Then, we have 
\[
\Ker F,\Image F,\Cok F\in 
 P_wL\nbiga(\Lambda_1,\Lambda_2), 
\]
and they have strict support $I$.
\end{lem}
\pf
Let $F_I:C_I\lrarr C_I'$ $(I\subset\Lambda_1)$
denote the underlying morphisms.
Let $I\subset J$.
We omit to denote the shift by $\Sigma^{p,q}$.
The induced morphisms
$g_{J,I}:\Image (F_I)\lrarr \Image(F_J)$
and
$g_{J,I}:\Cok(F_I)\lrarr\Cok(F_J)$
are clearly epimorphisms.
The induced morphisms
$f_{I,J}:\Image(F_J)\lrarr \Image(F_I)$
and $f_{I,J}:\Ker(F_J)\lrarr \Ker(F_I)$
are clearly monomorphisms.
Hence, we have only to prove that
(i) $g_{J,I}:\Ker(F_I)\lrarr\Ker(F_J)$
are epimorphisms,
(ii) $f_{I,J}:\Cok(F_J)\lrarr \Cok(F_I)$
are monomorphisms.
We have only to consider the case
$J=I\sqcup\{i\}$.

Because the category 
$\nbigp_w(\Lambda_1\sqcup\Lambda_2)$
is semisimple,
we can take a decomposition
$C_I=\Ker F_I\oplus C$
compatible with $N_k$ $(k\in\Lambda_1\sqcup \Lambda_2)$.
Then, we have
$\Image N_i=
 \bigl(\Image N_i\cap \Ker F_I\bigr)
\oplus
 \bigl(\Image N_i\cap C\bigr)
=N_i(\Ker F_I)\oplus N_i(C)$.
We have
$N_i\Ker (F_I)\subset
 f_{IJ}(\Ker F_J)
 \subset
 \Image N_i\cap\Ker F_I$.
Then, we obtain (i).
The claim (ii) can be proved similarly.
\hfill\qed

\vspace{.1in}

Hence, the decompositions by strict support
induce an equivalence of the categories:
\[
P_wL\nbiga(\Lambda_1,\Lambda_2)\simeq
 \bigoplus_{I\subset \Lambda_1}
 \nbigp_w(I\sqcup\Lambda_2)
\]
In particular,
the category $P_wL\nbiga(\Lambda_1,\Lambda_2)$
is abelian and semisimple,
and we obtain that
$P_wL\nbiga(\Lambda_1,\Lambda_2)$
satisfies the condition {\bf P0}.

\vspace{.1in}

Let us consider {\bf P3}
for $P_wL\nbiga(\Lambda_1,\Lambda_2)$.
Let $(\nbigt,\vecnbign)\in 
 P_wL\nbiga(\Lambda_1,\Lambda_2)$.
The condition {\bf P3.1} 
follows from that of objects
in $\nbigp_w(\Lambda_1\sqcup\Lambda_2)$.

Let $\Lambda_{2}=\Lambda_{2,0}\sqcup
\Lambda_{2,1}$ be a decomposition.
We have the induced object
\[
 \bigl(
 \Gr_{w'}^{M(\Lambda_{2,1})}\nbigt,
 \vecnbign_{\Lambda_{2,0}}
 \bigr)
 \in
 L\nbiga(\Lambda_1,\Lambda_{2,0}).
\]
To prove that it is pure of weight $w+w'$,
we have only to consider the case in which
$\Lambda_{2,1}$ consists of a unique element $\ast$,
according to {\bf P3.1} for
$P_wL\nbiga(\Lambda_1,\Lambda_2)$.
It satisfies {\bf LP2},
due to the condition {\bf P3.2} 
for objects in $\nbigp_w(\Lambda_1\sqcup\Lambda_2)$.
To prove {\bf LP1} for $\bigl(
 \Gr^{M(\ast)}\nbigt,
 \vecN_{\Lambda_{2,0}}
 \bigr)$,
we may assume that
$(\nbigt,\vecnbign)$ has strict support $I$.
Then, it follows from {\bf P3.3}
for objects in $\nbigp_w(\Lambda_1\sqcup\Lambda_2)$.

Let $\bullet\in\Lambda_2$.
We obtain an object
$\bigl(\Image\nbign_{\bullet},\vecnbign\bigr)$
in $L\nbiga(\Lambda_1,\Lambda_2)$.
We shall prove that the object
$\bigl(\Image\nbign_{\bullet},\vecnbign\bigr)$
is pure of weight $w+1$.
We have only to consider the case
in which $\nbigt$ has strict support $I$.
It is clear that {\bf LP1} is satisfied.
The condition {\bf LP2}
follows from {\bf P3.3}
for $\nbigp_w(\Lambda_1\sqcup\Lambda_2)$.

Let $\ast\in\Lambda_2\setminus\bullet$.
Let us look at the induced morphisms:
\begin{equation}
\label{eq;10.10.12.10}
 \Sigma^{1,0} \Gr^{M(\ast)}\nbigt
\lrarr
 \Sigma^{1,0}\Gr^{M(\ast)}\Image \nbign_{\bullet}
\lrarr
 \Sigma^{0,-1} \Gr^{M(\ast)}\nbigt
\end{equation}
It is a tuple of the following morphisms:
\[
 \Sigma^{1,0} \Gr^{M(\ast)}C_I
\lrarr
 \Sigma^{1,0}\Gr^{M(\ast)}\Image \nbign_{\bullet,I}
\lrarr
 \Sigma^{0,-1} \Gr^{M(\ast)}C_I
\]
It is $S$-decomposable by
{\bf P3.3} for objects in
$\nbigp_w(\Lambda)$.
Hence, (\ref{eq;10.10.12.10})
is $S$-decomposable.
Thus, the proof of Proposition \ref{prop;10.10.12.5}
is finished.
\hfill\qed

\subsection{Mixed objects in
$L\nbiga(\Lambda_1,\Lambda_2)$}
\label{subsection;10.10.13.2}

A filtered object
$(\nbigt,\vecnbign,\nbigl)\!\in\!
 L\nbiga(\Lambda_1,\Lambda_2)^{\fil}$
is called mixed,
if the following holds:
\begin{description}
\item[\bf (LM1)]
$\Gr^{\nbigl}_w(\nbigt,\vecnbign)\in 
 P_wL\nbiga(\Lambda_1,\Lambda_2)$.
\item[\bf (LM2)]
Each $(C_I,\nbigl,\vecN^{(\Lambda_2)})$ 
is an object in
$\nbigm(\Lambda_1\sqcup\Lambda_2)$
for $I\subset\Lambda_1$.
\end{description}
Let $ML\nbiga(\Lambda_1,\Lambda_2)
\subset
 \nbiga(\Lambda_1,\Lambda_2)^{\fil}$
denote the full subcategory of
mixed objects.
If $\Lambda_2=\emptyset$,
it is denoted by $ML\nbiga(\Lambda_1)$.

\begin{prop}
\label{prop;10.10.12.11}
The family $\{ML\nbiga(\Lambda_1,\Lambda_2)\}$
satisfies the conditions {\bf M0}--{\bf M3}.
\end{prop}
\pf
The condition {\bf M1} is clearly satisfied.
Let $F:(\nbigt_1,\nbigl,\vecnbign)
\lrarr(\nbigt_2,\nbigl,\vecnbign)$
be a morphism in $ML\nbiga(\Lambda_1,\Lambda_2)$.
By the condition {\bf P2} for 
$P_wL\nbiga(\Lambda_1,\Lambda_2)$,
we have
the strictness of $F$ with respect to $\nbigl$.
In particular,
$\Ker F$, $\Cok F$ and $\Image F$
are naturally equipped with filtrations $\nbigl$,
and we have natural isomorphisms
$\Gr^{\nbigl}\Ker F
\simeq
 \Ker\Gr^{\nbigl}F$,
$\Gr^{\nbigl}\Image F
\simeq
 \Image\Gr^{\nbigl} F$
and 
$\Gr^{\nbigl}\Cok F
\simeq
 \Cok\Gr^{\nbigl} F$.
Hence, {\bf LM1} holds for 
$\Ker F$, $\Cok F$ and $\Image F$
by Proposition \ref{prop;10.10.12.5}.
By the condition {\bf M0} for 
$\nbigm(\Lambda_1\sqcup\Lambda_2)$,
the condition {\bf LM2}
also holds for $\Ker F$, $\Image F$ and $\Cok F$.
Hence, $ML\nbiga(\Lambda_1,\Lambda_2)$
is abelian. 
Clearly, any injection
$\Lambda_2'\lrarr\Lambda_2$
induces a functor
$ML\nbiga(\Lambda_1,\Lambda_2')
\lrarr
 ML\nbiga(\Lambda_1,\Lambda_2)$.
Hence, {\bf M0} is satisfied for 
$ML\nbiga(\Lambda_1,\Lambda_2)$.

Let us consider {\bf M2}.
For $(\nbigt,\nbigl,\vecnbign)\in
 ML\nbiga(\Lambda_1,\Lambda_2)$,
the condition {\bf M2.1}
is satisfied by definition.
Let $\Lambda_2=\Lambda_{2,0}\sqcup\Lambda_{2,1}$
be a decomposition.
Each underlying $C_I$  has
$\nbigltilde:=
 M\bigl(\nbign(\Lambda_{2,1}),\nbigl(C_I)\bigr)$
by the condition {\bf M2.2} 
for objects in $\nbigm(\Lambda_1\sqcup\Lambda_2)$.
It induces a filtration $\nbigltilde$ of 
$(\nbigt,\vecnbign_{\Lambda_{2,0}})$
in $L\nbiga(\Lambda_1,\Lambda_{2,0})$,
which is a relative monodromy filtration
$M(\nbign(\Lambda_{2,1}),\nbigl)$:
\[
 \res_{\Lambda_{2,0}}^{\Lambda_2}
 (\nbigt,\nbigl,\vecnbign)
:=(\nbigt,\nbigltilde,\vecnbign_{\Lambda_{2,0}})
\]
The condition {\bf LM2} for 
$\res_{\Lambda_{2,0}}^{\Lambda_2}
 (\nbigt,\nbigl,\vecnbign)$
follows from
{\bf M2.2} for objects in 
$\nbigm(\Lambda_1\sqcup\Lambda_2)$.
By the canonical decomposition
(\S\ref{subsection;10.12.28.1}),
we have an isomorphism
\[
\Gr^{\nbigltilde}(\nbigt,\vecnbign_{\Lambda_{2,0}})
\simeq
\Gr^{\nbigltilde}\Gr^{\nbigl}
(\nbigt,\vecnbign_{\Lambda_{2,0}}).
\]
Hence, {\bf LM1} for
$\res_{\Lambda_{2,0}}^{\Lambda_2}
 (\nbigt,\nbigl,\vecnbign)$
follows from {\bf P3.2} for
objects in 
$P_wL\nbiga(\Lambda_1,\Lambda_2)$.
Hence, we obtain that 
$\res_{\Lambda_{2,0}}^{\Lambda_2}
 (\nbigt,\nbigl,\vecnbign)$
is an object in $ML\nbiga(\Lambda_1,\Lambda_2)$,
i.e.,
{\bf M2} holds for 
$ML\nbiga(\Lambda_1,\Lambda_2)$.

\vspace{.1in}
Let us consider {\bf M3}
for $ML\nbiga(\Lambda_1,\Lambda_2)$.
Let $\bullet\in\Lambda_2$,
and we put $\Lambda_{2,0}:=
 \Lambda_2\setminus\bullet$.
We consider
$(\nbigt,\nbigl,\vecnbign_{\Lambda_2})
 \in
 ML\nbiga(\Lambda_1,\Lambda_2)$
and 
$(\nbigt',\nbigltilde,\vecnbign'_{\Lambda_{2,0}})
 \in
 ML\nbiga(\Lambda_1,\Lambda_{2,0})$
with morphisms
\[
 \Sigma^{1,0}\res_{\Lambda_{2,0}}^{\Lambda_2}
 (\nbigt,\nbigl,\vecnbign_{\Lambda_2})
\stackrel{u}{\lrarr}
 (\nbigt',\nbigltilde,\vecnbign'_{\Lambda_{2,0}})
\stackrel{v}{\lrarr}
 \Sigma^{0,-1}\res_{\Lambda_{2,0}}^{\Lambda_2}
 (\nbigt,\nbigl,\vecnbign_{\Lambda_2})
\]
in $ML\nbiga(\Lambda_1,\Lambda_{2,0})$
such that $v\circ u=\nbign_{\bullet}$.
We have the filtration
$\nbigl(\nbigt')$ obtained
as the transfer of $\nbigl(\nbigt)$.
We put $\nbign'_{\bullet}:=u\circ v$.
Let us prove that
$(\nbigt',\nbigl,\vecnbign'_{\Lambda_2})$
is an object in 
$ML\nbiga(\Lambda_1,\Lambda_2)$.
The condition {\bf LM2} for 
$ML\nbiga(\Lambda_1,\Lambda_2)$
follows from
{\bf M3} for objects in
$\nbigm(\Lambda_1\sqcup\Lambda_2)$.
By the construction of $\nbigl(\nbigt')$,
we have the decomposition
\[
 \Gr_w^{\nbigl}(\nbigt')
=\bigoplus
 \Image \Gr_w^{\nbigl}(u)
\oplus
 \Ker\Gr_w^{\nbigl}(v).
\]
Because $\nbigltilde=M(\nbign_{\bullet},\nbigl)$,
we have a natural splitting
$\Gr^{\nbigltilde}_w\Ker\Gr_w^{\nbigl}(v)
\simeq
 \Ker\Gr_w^{\nbigl}(v)$.
Hence, 
$\Ker\Gr_w^{\nbigl}(v)$
with the induced tuple of morphisms
denoted by
$\Gr_w^{\nbigl}(\vecnbign_{\Lambda_{2,0}})$
is isomorphic to a direct summand of 
\[
\Gr^{\nbigltilde}_w\Gr^{\nbigl}_w(\nbigt',
 \vecnbign_{\Lambda_{2,0}})
\simeq 
\Gr^{\nbigltilde}_w(\nbigt',\vecnbign_{\Lambda_{2,0}})
\in P_wL\nbiga(\Lambda_1,\Lambda_{2,0})
\subset
 P_wL\nbiga(\Lambda_1,\Lambda_2)
\]
By {\bf P0} for
$P_wL\nbiga(\Lambda_1,\Lambda_2)$,
we obtain that
$\bigl(\Ker\Gr^{\nbigl}(v),
 \Gr^{\nbigl}(\vecnbign_{\Lambda_{2,0}})\bigr)$
is an object in 
$P_wL\nbiga(\Lambda_1,\Lambda_2)$.
We have the natural isomorphism
\[
 \Sigma^{1,0}\bigl(
 \Image \Gr^{\nbigl}\nbign_{\bullet},
 \Gr^{\nbigl}\vecnbign_{\Lambda_2}
 \bigr)
\simeq
 \bigl(
 \Image\Gr^{\nbigl}u,
 \Gr^{\nbigl}\vecnbign_{\Lambda_2}
 \bigr).
\]
Hence, we obtain that
$\bigl(
 \Image\Gr^{\nbigl}u,
 \Gr^{\nbigl}\vecnbign_{\Lambda_2}
 \bigr)$ is an object
in $P_wL\nbiga(\Lambda_1,\Lambda_2)$
by {\bf P3.3} in Proposition \ref{prop;10.10.12.5}.
Hence, 
$(\nbigt',\nbigl,\vecnbign'_{\Lambda_2})$
is an object in 
$ML\nbiga(\Lambda_1,\Lambda_2)$,
and the condition {\bf M3}
holds for $ML\nbiga(\Lambda_1,\Lambda_2)$.
Thus, the proof of Proposition \ref{prop;10.10.12.11}
is finished.
\hfill\qed

\subsection{Some functors}
\label{subsection;10.10.12.15}

For any $K\subset\Lambda_1$,
we have the functors
\[
\psi_K,\phi_K: L\nbiga(\Lambda_1,\Lambda_2)
\lrarr 
 L\nbiga(\Lambda_1\setminus K,\Lambda_2\sqcup K)
\]
given as follows.
Let $\nbigt\in L\nbiga(\Lambda_1,\Lambda_2)$.
For $I\subset\Lambda_1\setminus K$,
we set
$\psi_K(\nbigt)_I:=C_I$
and
$\phi_K(\nbigt)_I:=C_{I\sqcup K}$.
For $I\subset J$,
we have the naturally induced morphisms
\[
 \Sigma^{|J\setminus I|,0}
\psi_K(\nbigt)_I
\lrarr
 \psi_K(\nbigt)_J
\lrarr
 \Sigma^{0,-|J\setminus I|}
 \psi_K(\nbigt)_I
\]
\[
 \Sigma^{|J\setminus I|,0}
 \phi_K(\nbigt)_I
\lrarr
 \phi_K(\nbigt)_J
\lrarr
 \Sigma^{0,-|J\setminus I|}
 \phi_K(\nbigt)_I
\]
with which they are objects in 
$L\nbiga(\Lambda_1\setminus K)$.
They are equipped with
the naturally induced morphisms
$\vecnbign_{\Lambda_2}$.
Moreover, $\psi_K(\nbigt)$ and $\phi_K(\nbigt)$
are equipped with
naturally induced morphisms 
$\vecnbign_{K}:=(\nbign_i\,|\,i\in K)$.
We set
$\vecnbign_{\Lambda_2\sqcup K}
=\vecnbign_{\Lambda_2}\sqcup\vecnbign_K$.
Hence, we obtain
$\bigl(\phi_K(\nbigt),
 \vecnbign_{\Lambda_2\sqcup K}\bigr)$
and 
$\bigl(\psi_K(\nbigt),
 \vecnbign_{\Lambda_2\sqcup K}\bigr)$.
They naturally induce the functors
\[
  \psi_K,\phi_K:
 L\nbiga(\Lambda_1,\Lambda_2)^{\fil}
\lrarr
 L\nbiga(\Lambda_1\setminus K,
 \Lambda_2\sqcup K)^{\fil}.
\]
They naturally induce
\[
 \psi_K,\phi_K:
 ML\nbiga(\Lambda_1,\Lambda_2)
\lrarr
 ML\nbiga(\Lambda_1\setminus K,\Lambda_2\sqcup K)
\]
\[
  \psi_K,\phi_K:P_wL\nbiga(\Lambda_1,\Lambda_2)
\lrarr
 P_wL\nbiga(\Lambda_1\setminus K,
 \Lambda_2\sqcup K).
\]
The induced functor
$\res_{\Lambda_2}^{\Lambda_2\sqcup K}\phi_K:
  ML\nbiga(\Lambda_1,\Lambda_2)
\lrarr
 ML\nbiga(\Lambda_1\setminus K,\Lambda_2)$
is also denoted by $\phi_K$,
if there is no risk of confusion.

\subsection{Gluing}

Fix an element $\bullet\in\Lambda_1$.
We consider the category
$\Glue(\Lambda_1,\Lambda_2,\bullet)$
given as follows.
An object in $\Glue(\Lambda_1,\Lambda_2,\bullet)$
is a tuple of
$\bigl(\nbigt',\nbigl,
 \vecnbign'\bigr)
\in ML\nbiga(\Lambda_1\setminus\bullet,
 \Lambda_2\sqcup\bullet)$
and 
$\bigl(\nbigt'',\nbigltilde,
 \vecnbign''\bigr)
\in ML\nbiga(\Lambda_1\setminus\bullet,
 \Lambda_2)$
with morphisms
$\nbigu$ and $\nbigv$ 
in $ML\nbiga(\Lambda_1\setminus\bullet,
 \Lambda_2)$,
\[
 \Sigma^{1,0}
 \res^{\Lambda_2\sqcup\bullet}_{\Lambda_2}
 \bigl(\nbigt',\nbigl,\vecnbign'
 \bigr)
\stackrel{\nbigu}{\lrarr}
 \bigl(\nbigt'',\nbigltilde,
 \vecnbign''\bigr)
\stackrel{\nbigv}{\lrarr}
 \Sigma^{0,-1}
 \res^{\Lambda_2\sqcup\bullet}_{\Lambda_2}
 \bigl(\nbigt',\nbigl,\vecnbign'
 \bigr) 
\]
such that
$\nbign'_{\bullet}=\nbigv\circ\nbigu$.
A morphism in $\Glue(\Lambda_1,\Lambda_2,\bullet)$
is a tuple of morphisms
\[
F':\bigl(\nbigt_1',\nbigl,
 \vecnbign'\bigr)
\lrarr
\bigl(\nbigt_2',\nbigl,\vecnbign'\bigr)
\quad
 \mbox{\rm in }
ML\nbiga(\Lambda_1\setminus\bullet,
 \Lambda_2\sqcup\bullet)
\]
\[
 F'':\bigl(\nbigt_1'',\nbigltilde,
 \vecnbign''\bigr)
\lrarr
\bigl(\nbigt_2'',\nbigltilde,
 \vecnbign''\bigr)
\quad
 \mbox{\rm in }
ML\nbiga(\Lambda_1\setminus\bullet,
 \Lambda_2)
\]
such that the following diagram is commutative:
{\small
\[
 \begin{CD}
  \Sigma^{1,0}
 \res^{\Lambda_2\sqcup\bullet}_{\Lambda_2}
 \bigl(\nbigt_1',\nbigl,\vecnbign'
 \bigr) 
@>>>
 \bigl(\nbigt''_1,\nbigltilde,
 \vecnbign''\bigr)
@>>>
 \Sigma^{0,-1}
 \res^{\Lambda_2\sqcup\bullet}_{\Lambda_2}
 \bigl(\nbigt_1',\nbigl,\vecnbign'
 \bigr)  \\
 @V{F'}VV @V{F''}VV @V{F'}VV \\
  \Sigma^{1,0}
 \res^{\Lambda_2\sqcup\bullet}_{\Lambda_2}
 \bigl(\nbigt_2',\nbigl,\vecnbign'
 \bigr) 
@>>>
 \bigl(\nbigt''_2,\nbigltilde,
 \vecnbign''\bigr)
@>>>
 \Sigma^{0,-1}
 \res^{\Lambda_2\sqcup\bullet}_{\Lambda_2}
 \bigl(\nbigt_2',\nbigl,\vecnbign'
 \bigr) 
 \end{CD}
\]
}
For
$(\nbigt,\nbigl,\vecnbign)
 \in
 ML\nbiga(\Lambda_1,\Lambda_2)$,
we have 
$\Gamma(\nbigt,\nbigl,\vecnbign)$
in $\Glue(\Lambda_1,\Lambda_2,\bullet)$
which is the tuple of the objects
\[
 \psi_{\bullet}
 \bigl(\nbigt,\nbigl,\vecnbign
 \bigr)
\in ML\nbiga(\Lambda_1\setminus\bullet,
 \Lambda_2\sqcup\bullet),
\]
\[
 \res^{\Lambda_2\sqcup\bullet}_{\Lambda_2}
 \phi_{\bullet}(\nbigt,\nbigl,\vecnbign)
\in ML\nbiga(\Lambda_1\setminus\bullet,
 \Lambda_2).
\]
with naturally induced morphisms
$\nbigu$ and $\nbigv$:
{\small
\[
 \Sigma^{1,0}\res^{\Lambda_2\sqcup\bullet}_{\Lambda_2}
 \psi_{\bullet}(\nbigt,\nbigl,\vecnbign)
\stackrel{\nbigu}{\lrarr}
 \res^{\Lambda_2\sqcup\bullet}_{\Lambda_2}
 \phi_{\bullet}(\nbigt,\nbigl,\vecnbign)
\stackrel{\nbigv}{\lrarr}
 \Sigma^{0,-1}\res^{\Lambda_2\sqcup\bullet}_{\Lambda_2}
 \psi_{\bullet}(\nbigt,\nbigl,\vecnbign)
\]
}
Thus, we obtain a functor
$\Gamma:
 ML\nbiga(\Lambda_1,\Lambda_2)
\lrarr
 \Glue(\Lambda_1,\Lambda_2,\bullet)$.

\begin{cor}
\label{cor;10.10.12.11}
The functor $\Gamma$ is an equivalence.
\end{cor}
\pf
Let 
$\bigl(\nbigt',\nbigl,
 \vecnbign'\bigr)$,
$\bigl(\nbigt'',\nbigltilde,
 \vecnbign''\bigr)$,
$\nbigu$ and $\nbigv$ be as above.
We have the filtration $\nbigl$
of $(\nbigt'',\vecnbign'')$
obtained as the transfer.
By {\bf M3} in Proposition \ref{prop;10.10.12.11},
$(\nbigt'',\nbigl,\vecnbign'')$
is an object in 
$ML\nbiga(\Lambda_1\setminus\bullet,
 \Lambda_2\sqcup\bullet)$.
We put
$C_I:=C_I'$ if $\bullet\not\in I$
and $C_I:=C_{I\setminus\bullet}''$ if $\bullet\in I$.
The tuple $(C_I\,|\,I\subset\Lambda_1)$
with the induced morphisms naturally
gives an object in $L\nbiga(\Lambda_1,\Lambda_2)$.
Then, it is easy to check
{\bf LM1--2} for this object,
and we obtain the essential surjectivity
of $\Gamma$.
It is obvious to see that $\Gamma$ is faithful.
It is full thanks to the functoriality of the transfer.
\hfill\qed

\subsection{Another description of
$ML\nbiga(\Lambda_1,\Lambda_2)$}
\label{subsection;10.10.13.1}

We shall introduce a category
$M'L\nbiga(\Lambda_1,\Lambda_2)$.
An object of $M'L\nbiga(\Lambda_1,\Lambda_2)$
is $(\nbigt,\vecnbign_{\Lambda_2})
\in L\nbiga(\Lambda_1,\Lambda_2)$
with a tuple of filtrations
$\vecL=\bigl(L_I\,\big|\,I\subset\Lambda_1\bigr)$
of the underlying objects
$C_I$ such that the following holds:
\begin{itemize}
\item
Each $\nbigu_I=(C_I,L_I,
 \vecN_{\Lambda_1\setminus I}^{(\Lambda_2)})$
is an object in
$\nbigm(\Lambda_1\sqcup\Lambda_2\setminus I)$,
where
$\vecN^{(\Lambda_2)}_{\Lambda_1\setminus I}
=\bigl(N_j\,\big|\,j\in\Lambda_1\setminus I\bigr)
\sqcup \vecnbign_{\Lambda_2}$.
\item
 $g_{JI}$ and $f_{IJ}$ give the following morphisms 
 in $\nbigm\bigl(\Lambda_1\sqcup
 \Lambda_2\setminus J\bigr)$:
\[
 \begin{CD}
 \Sigma^{|J\setminus I|,0}
 \res^{\Lambda_1\sqcup\Lambda_2\setminus I}
 _{\Lambda_1\sqcup\Lambda_2\setminus J}
 \bigl(\nbigu_I\bigr)
@>{g_{JI}}>>
 \nbigu_J
@>{f_{IJ}}>>
 \Sigma^{0,-|J\setminus I|}
 \res^{\Lambda_1\sqcup\Lambda_2\setminus I}
 _{\Lambda_1\sqcup\Lambda_2\setminus J}
 \bigl(\nbigu_I\bigr)
 \end{CD}
\]
\end{itemize}
A morphism
$(\nbigt_1,\vecnbign_{\Lambda_2},\vecL)
\lrarr(\nbigt_2,\vecnbign_{\Lambda_2},\vecL)$
in $M'L\nbiga(\Lambda_1,\Lambda_2)$
is a morphism
$F:(\nbigt_1,\vecnbign_{\Lambda_2})
\lrarr(\nbigt_2,\vecnbign_{\Lambda_2})$
in $L\nbiga(\Lambda_1,\Lambda_2)$
such that
each $F_I:C_{1,I}\lrarr C_{2,I}$
is compatible with the filtrations $L_I$.

We have the functor
$\Phi_{\Lambda_1,\Lambda_2}:
 ML\nbiga(\Lambda_1,\Lambda_2)
\lrarr
 M'L\nbiga(\Lambda_1,\Lambda_2)$
defined as follows.
For $(\nbigt,\nbigl,\vecnbign_{\Lambda_2})$,
we set
$L_I(C_I):=M\Bigl(N(I);\nbigl(C_I) \Bigr)$.
Then,
$(\nbigt,\vecnbign_{\Lambda_2})$
with $\vecL=(L_I\,|\,I\subset\Lambda_1)$
is an object in $M'L\nbiga(\Lambda_1,\Lambda_2)$,
which is defined to be
$\Phi_{\Lambda_1,\Lambda_2}
 (\nbigt,\nbigl,\vecnbign_{\Lambda_2})$.

\begin{thm}
\label{thm;10.10.12.12}
The above functors 
$\Phi_{\Lambda_1,\Lambda_2}$
are equivalent.
\end{thm}
\pf
The claim is clear
in the case $|\Lambda_1|=0$.
We use an induction on
$\bigl(|\Lambda_1|+|\Lambda_2|,|\Lambda_1|\bigr)$.
Take $\bullet\in\Lambda_1$.
We assume that
$\Phi_{\Lambda_1\setminus\bullet,
 \Lambda_2\sqcup\bullet}$
and 
$\Phi_{\Lambda_1\setminus\bullet, 
 \Lambda_2}$
are equivalent,
and we will prove that
$\Phi_{\Lambda_1,\Lambda_2}$
is equivalent.

\vspace{.1in}

Let $(\nbigt,\vecnbign_{\Lambda_2},\vecL)
\in ML\nbiga'(\Lambda_1,\Lambda_2)$.
Each $C_I$ is equipped with a filtration
$\Ltilde_I:=M(\nbign_{\bullet},L_I(C_I))$.
Thus, we obtain a functor
\[
 \res_{\Lambda_2\setminus\bullet}^{\Lambda_2}:
 M'L\nbiga(\Lambda_1,\Lambda_2)
\lrarr M'L\nbiga(\Lambda_1,\Lambda_2\setminus\bullet)
\]
given by
$\res_{\Lambda_2\setminus\bullet}^{\Lambda_2}
 \bigl(
 \nbigt,\vecnbign_{\Lambda_2},\vecL
 \bigr)
=\bigl(
 \nbigt,\vecnbign_{\Lambda_2\setminus\bullet},
 \widetilde{\vecL}
 \bigr)$.
We also have the functor
\[
 \res_{\Lambda_2\setminus\bullet}^{\Lambda_2}:
 ML\nbiga(\Lambda_1,\Lambda_2)
\lrarr
 ML\nbiga(\Lambda_1,\Lambda_2\setminus\bullet).
\]

\begin{lem}
\label{lem;10.10.12.20}
We have 
$\res_{\Lambda_2\setminus\bullet}^{\Lambda_2}
 \circ\Phi_{\Lambda_1,\Lambda_2}
=\Phi_{\Lambda_1,\Lambda_2\setminus \bullet}
 \circ\res_{\Lambda_2\setminus\bullet}^{\Lambda_2}$.
\end{lem}
\pf
We have 
\[
 M\Bigl(
 N(I);
 M(\nbign_{\bullet};\nbigl)
 \Bigr)
=M\Bigl(
 \nbign_{\bullet}+N(I);\nbigl
 \Bigr)
=
M\Bigl(
 \nbign_{\bullet};
 M\Bigl(N(I);\nbigl\Bigr)
 \Bigr).
\]
Then, we obtain the claim of the lemma.
\hfill\qed

\vspace{.1in}

We consider a category
$\Glue'(\Lambda_1,\Lambda_2,\bullet)$
given as follows.
An object is a tuple of
\begin{equation}
 \label{eq;10.10.12.13}
\bigl(\nbigt',\vecnbign',\vecL
 \bigr)\in
M'L\nbiga(\Lambda_1\setminus\bullet, 
 \Lambda_2\sqcup\bullet),
\quad\quad
 \bigl(\nbigt'',\vecnbign'',\vecLtilde\bigr)
\in M'L\nbiga(\Lambda_1\setminus\bullet,
 \Lambda_2)
\end{equation}
with morphisms
in $M'L\nbiga(\Lambda_1\setminus\bullet,
 \Lambda_2)$
\begin{equation}
 \label{eq;10.10.12.14}
 \Sigma^{1,0}
 \res_{\Lambda_2}^{\Lambda_2\sqcup\bullet}
 \bigl(
 \nbigt',\vecnbign',\vecL
 \bigr)
\stackrel{\nbigu}{\lrarr}
 \bigl(\nbigt'',\vecnbign'',
 \vecLtilde\bigr)
\stackrel{\nbigv}{\lrarr}
 \Sigma^{0,-1}
 \res_{\Lambda_2}^{\Lambda_2\sqcup\bullet}
 \bigl(
 \nbigt',\vecnbign',\vecL
 \bigr)
\end{equation}
such that $\nbigv\circ\nbigu=\nbign_{\bullet}$.
Morphisms in $\Glue'(\Lambda_1,\Lambda_2,\bullet)$
are naturally defined.

By Lemma \ref{lem;10.10.12.20},
we have a naturally defined functor
$\Phi_0:\Glue(\Lambda_1,\Lambda_2,\bullet)
\lrarr\Glue'(\Lambda_1,\Lambda_2,\bullet)$.

\begin{lem}
$\Phi_0:\Glue(\Lambda_1,\Lambda_2,\bullet)
\lrarr
 \Glue'(\Lambda_1,\Lambda_2,\bullet)$
is equivalent.
\end{lem}
\pf
It is easy to observe that
the functor $\Phi_0$ is fully faithful.
Let us prove the essential surjectivity.
Take an object in $\Glue'(\Lambda_1,\Lambda_2,\bullet)$,
i.e., a tuple as in (\ref{eq;10.10.12.13})
with morphisms as in (\ref{eq;10.10.12.14}).
By using the hypothesis of the induction,
we obtain the following:
\begin{itemize}
\item
 A filtration $\nbigl$ of
 $(\nbigt,\vecnbign)$ such that
$(\nbigt,\nbigl,\vecnbign)
 \in ML\nbiga(\Lambda_1\setminus\bullet,
 \Lambda_2\sqcup\bullet)$
and
$\Phi_{\Lambda_1\setminus\bullet,
 \Lambda_2\sqcup\bullet}
 (\nbigt,\nbigl,\vecnbign)
=(\nbigt,\vecnbign,\vecL)$.
\item
 A filtration $\nbigltilde$
 of $(\nbigt',\vecnbign')$
 such that
$ (\nbigt',\nbigltilde,\vecnbign')
 \in ML\nbiga(\Lambda_1\setminus\bullet,\Lambda_2)$
and
\[
 \Phi_{\Lambda_1\setminus\bullet,\Lambda_2}
 (\nbigt',\nbigltilde,\vecnbign)
=(\nbigt',\vecnbign',\vecL).
\]
\item
 A filtration $\nbigltilde$
 of $(\nbigt,\vecnbign_{\Lambda_2})$
 such that
$(\nbigt,\nbigltilde,
 \vecnbign_{\Lambda_2})
 \in ML\nbiga(\Lambda_1\setminus\bullet,\Lambda_2)$
and
\[
 \Phi_{\Lambda_1\setminus\bullet,\Lambda_2}
 (\nbigt,\nbigltilde,
 \vecnbign_{\Lambda_2})
 =\res_{\Lambda_2}^{\Lambda_2\sqcup\bullet}
(\nbigt,\vecnbign, \vecL),
\]
where $\vecnbign_{\Lambda_2}:=\bigl(
 \nbign_j\,\big|\,j\in\Lambda_2
 \bigr)$.
\end{itemize}
Because 
$\Phi_{\Lambda_1\setminus\bullet,\Lambda_2\sqcup\bullet}$
is equivalent,
$\res^{\Lambda_2\sqcup\bullet}_{\Lambda_2}
 (\nbigt,\nbigl,\vecnbign)$ is equal to
$(\nbigt,\nbigltilde,
 \vecnbign_{\Lambda_2})$,
and 
$\nbigu$ and $\nbigv$ give morphisms
in $ML\nbiga(\Lambda_1\setminus\bullet,\Lambda_2)$.
Hence, 
$(\nbigt,\nbigl,\vecnbign)$,
$(\nbigt',\nbigltilde,
 \vecnbign')$
with $(\nbigu,\nbigv)$ naturally gives
an object in $\Glue(\Lambda_1,\Lambda_2)$,
and we can deduce that
$\Phi_0$ is essentially surjective.
\hfill\qed

\vspace{.1in}
We prepare some functors.
We have the naturally induced functors
\[
 \phi_{\bullet}:
 M'L\nbiga(\Lambda_1,\Lambda_2)
\lrarr
 M'L\nbiga(\Lambda_1\setminus\bullet,
 \Lambda_2),
\]
\[
 \psi_{\bullet}:
 M'L\nbiga(\Lambda_1,\Lambda_2)
 \lrarr
 M'L\nbiga(\Lambda_1\setminus\bullet,
 \Lambda_2\sqcup\bullet).
\]
For $(\nbigt,\nbigl,\vecnbign)\in
 M'L\nbiga(\Lambda_1,\Lambda_2)$,
we have induced morphisms
$\nbigu:\Sigma^{1,0}\psi(\nbigt)\lrarr\phi(\nbigt)$
and $\nbigv:\phi(\nbigt)\lrarr\Sigma^{0,-1}\psi(\nbigt)$,
induced  by
$g_{I,I\setminus\bullet}$ and 
$f_{I\setminus\bullet,I}$.
They give morphisms 
in $M'L\nbiga(\Lambda_1\setminus\bullet,\Lambda_2)$:
\[
 \Sigma^{1,0}
 \res_{\Lambda_2}^{\Lambda_2\sqcup\bullet}
 \psi_{\bullet}(\nbigt,\vecnbign,\vecL)
\lrarr
 \phi_{\bullet}(\nbigt,\vecnbign,\vecL)
\lrarr
 \Sigma^{0,-1}
 \res_{\Lambda_2}^{\Lambda_2\sqcup\bullet}
 \psi_{\bullet}(\nbigt,\vecnbign,\vecL)
\]
Hence, we obtain a functor
$\Gamma':
M'L\nbiga(\Lambda,\Lambda_2)
\lrarr
 \Glue'(\Lambda_1,\Lambda_2,\bullet)$.
Because
$\Gamma'\circ\Phi_{\Lambda_1,
 \Lambda_2}
=\Phi_0\circ\Gamma$
by construction,
we have only to prove that
$\Gamma'$ is an equivalence.
But, it is clear by construction.
Thus, the proof of Theorem \ref{thm;10.10.12.12}
is finished.
\hfill\qed

\vspace{.1in}
For any $K\subset\Lambda_1$,
we have naturally defined functor
$\psi_K:M'L\nbiga(\Lambda_1,\Lambda_2)
\lrarr M'L\nbiga(\Lambda_1\setminus K,
 \Lambda_2\sqcup K)$,
and the following is commutative:
\begin{equation}
 \begin{CD}
 ML\nbiga(\Lambda_1,\Lambda_2)
 @>{\psi_K}>>
 ML\nbiga(\Lambda_1\setminus K,\Lambda_2\sqcup K)\\
 @V{\simeq}VV @V{\simeq}VV \\
 M'L\nbiga(\Lambda_1,\Lambda_2)
 @>{\psi_K}>>
 M'L\nbiga(\Lambda_1\setminus K,\Lambda_2\sqcup K)\\
 \end{CD}
\end{equation}

\subsection{Commutativity of the transfer}
\label{subsection;10.11.5.13}

Let $V=(C,L,\vecN)\in\nbigm(\nibar)$,
where $\nibar:=\{1,2\}$.
We obtain the following object 
$\nbigt=:V[!1\ast 2]$ in $M'L\nbiga(\nibar)$.
We put
\[
 (\nbigt_{\emptyset},L_{\emptyset}):=(C,L),
\quad
 (\nbigt_1,L_1):=
 \Sigma^{1,0}\bigl(C,M(N_1;L)\bigr),
\]
\[
  (\nbigt_2,L_2):=
 \Sigma^{0,-1}\bigl(C,M(N_2;L)\bigr),
\quad
 (\nbigt_{\nibar},L_{\nibar}):=
 \Sigma^{1,-1}
 \bigl(C,M(N_1+N_2;L)\bigr).
\]
The morphisms $f_{IJ}$ and $g_{JI}$
are given as follows:
\[
\begin{CD}
 \Sigma^{1,0}
 \bigl(\nbigt_{\emptyset},M(N_1;L_{\emptyset})\bigr)
 @>{\id}>>
 (\nbigt_{1},L_1)
 @>{N_1}>>
 \Sigma^{0,-1}\bigl(\nbigt_{\emptyset},
 M(N_1;L_{\emptyset})\bigr)
\end{CD}
\]
\[
\begin{CD}
 \Sigma^{1,0}\bigl(\nbigt_{\emptyset},M(N_2;L_{\emptyset})\bigr)
 @>{N_2}>>
 (\nbigt_{2},L_2)
 @>{\id}>>
 \Sigma^{0,-1}\bigl(\nbigt_{\emptyset},
 M(N_2;L_{\emptyset})\bigr)
\end{CD}
\]
\[
\begin{CD}
 \Sigma^{1,0}
 \bigl(\nbigt_{1},M(N_2;L_{1})\bigr)
@>{N_2}>>
 (\nbigt_{\nibar},L_{\nibar})
@>{\id}>>
 \Sigma^{0,-1}
 \bigl(\nbigt_{1},
 M(N_2;L_{1})\bigr)
\end{CD}
\]
\[
\begin{CD}
 \Sigma^{1,0}
 \bigl(\nbigt_{2},M(N_1;L_{2})\bigr)
 @>{\id}>>
 (\nbigt_{\nibar},L_{\nibar})
 @>{N_1}>>
 \Sigma^{0,-1}
 \bigl(\nbigt_{2},
 M(N_1;L_{2})\bigr)
\end{CD}
\]
By Theorem \ref{thm;10.10.12.12},
we obtain the filtration $\nbigl$ of 
$\nbigt=V[!1\ast 2]$
such that
$(\nbigt,\nbigl)\in ML\nbiga(\nibar)$.
It is easy to check the following:
\[
 \nbigl_{\emptyset}(C)=L,
\quad
 \nbigl_1(\Sigma^{1,0} C)=\Nhat_{1!}L,
\]
\[
 \nbigl_2(\Sigma^{0,-1}C)=\Nhat_{2\ast}L,
\quad
 \nbigl_{\nibar}(\Sigma^{1,-1}C)
=\Nhat_{2\ast}\Nhat_{1!}L
=\Nhat_{1!}\Nhat_{2\ast}L.
\]
In particular, we obtain the following commutativity.
\begin{lem}
\label{lem;10.10.12.22}
$N_{2\ast}N_{1!}L=N_{1!}N_{2\ast}L$.
\hfill\qed
\end{lem}

Similarly, we obtain the following lemma.
\begin{lem}
\label{lem;10.10.12.23}
$N_{2\star}N_{1\star}L
=N_{1\star}N_{2\star}L$
for $\star=\ast,!$.
\hfill\qed
\end{lem}

Let $\Lambda$ be a finite set.
\begin{prop}
\label{prop;10.10.12.25}
Let $(V,L,\vecN)$ be an object in 
$\nbigm(\Lambda)$.
Let $i,j\in\Lambda$ be distinct elements.
Let $\star_i,\star_j\in\{\ast,!\}$.
Then, we have
$N_{i\star_i}(N_{j\star_j}L)
=N_{j\star_j}(N_{i\star_i}L)$.
\end{prop}
\pf
Let $\Lambda_1:=\Lambda\setminus\{i,j\}$
and $\Lambda_2:=\{i,j\}$.
We set
$\nbigatilde:=\nbiga(\Lambda_1)$,
$\nbigptilde_w(\Lambda'):=
 \nbigp_w(\Lambda_1\sqcup\Lambda')$
and $\nbigmtilde(\Lambda')
=\nbigm(\Lambda_1\sqcup\Lambda')$.
Then, they satisfy the conditions
in \S\ref{subsection;10.10.12.21}.
Then, the claim follows from
Lemma \ref{lem;10.10.12.22}
and Lemma \ref{lem;10.10.12.23}.
\hfill\qed

\vspace{.1in}
Let $(C,L,\vecN)$ be an object in $\nbigm(\Lambda)$.
Let $J\subset\Lambda$.
We define a filtration
$\vecNhat_{J\ast}L$ of $\Sigma^{0,-|J|}C$
inductively;
take $j\in J$, put $J_0:=J\setminus\{j\}$,
and define
$\vecNhat_{J\ast}L:=\Nhat_{j\ast}\vecNhat_{J_0\ast}L$.
By Proposition \ref{prop;10.10.12.25},
it is independent of the choice of $j$.
Similarly,
we define $\vecNhat_{J!}L$ of $\Sigma^{|J|,0}C$.
We obtain $\vecN_{J\ast}L$ and $\vecN_{J!}L$
on $C$ by
$(C,\vecN_{J\ast}L):=
 \Sigma^{0,|J|}(C,\vecNhat_{J\ast}L)$
and 
$(C,\vecN_{J!}L):=
 \Sigma^{-|J|,0}(C,\vecNhat_{J!}L)$.

For $J_i\subset \Lambda$ $(i=1,2)$
with $J_1\cap J_2=\emptyset$,
we obtain a filtration
$\vecN_{J_1!}\vecN_{J_2\ast}L
=\vecN_{J_2\ast}\vecN_{J_1!}L$
on $C$.

\subsection{Canonical prolongations}
\label{subsection;11.2.18.1}

Let $V=(C,L,\vecN)$ be an object in $\nbigm(\Lambda)$.
Let $\Lambda_1\sqcup \Lambda_2=\Lambda$ 
be a decomposition.
We obtain an object 
$\nbigt=:V[\ast K_1!K_2]\in ML\nbiga(\Lambda)$
given as follows.
For $I\subset\Lambda$,
we set $I_j:=I\cap \Lambda_j$,
and 
$\nbigt_I:=
 \Sigma^{|I_2|,-|I_1|}C$.
It is equipped with a filtration
$\vecNhat_{I_1!}\vecNhat_{I_2\ast}L$.
For $I\sqcup\{i\}\subset\Lambda$,
morphisms $g_{Ii,I}$ and $f_{I,Ii}$
are given as follows:
\[
 g_{Ii,I}:=\left\{
 \begin{array}{ll}
 \id & (i\in \Lambda_2)\\
 N_i & (i\in\Lambda_1)
 \end{array}
 \right.
\quad\quad
 f_{I,Ii}:=\left\{
 \begin{array}{ll}
 N_i & (i\in \Lambda_2)\\
 \id & (i\in\Lambda_1)
 \end{array}
 \right.
\]

\chapter{Mixed twistor $D$-modules}
\label{section;11.4.9.1}

We shall introduce mixed twistor $D$-modules
(Definition \ref{df;11.4.3.30}).
The definition is not completely parallel to 
that for mixed Hodge modules in \cite{saito2}.
It is partially because we are concerned with 
only graded polarizable objects.
Since we shall later prove that any mixed twistor $D$-module
is expressed as a gluing of some 
admissible variations of mixed twistor structures
(\S\ref{subsection;11.4.3.41}
and \S\ref{subsection;11.4.3.40}),
as expected,
mixed twistor $D$-modules 
can be regarded as a twistor version of
mixed Hodge modules.

\section{Admissibly specializable pre-mixed twistor $D$-modules}

\subsection{Pre-mixed twistor $D$-modules}

\index{pre-mixed twistor $D$-module}

Let $\MT(X,w)$ denote the category of
polarizable pure twistor $D$-modules on $X$
of weight $w$. \index{category $\MT(X,w)$}
(We omit the adjective ``wild''.)
A filtered $\nbigr_X$-triple 
means a $\nbigr_X$-triple $\nbigt$
with an increasing filtration $W$
indexed by integers
such that
$W_j=0$ for $j<\!<0$
and $W_j(\nbigt)=\nbigt$ for $j>\!>0$
locally on $X$.
A filtered $\nbigr_X$-triple
$(\nbigt,\nbigw)$ is called
a pre-mixed twistor $D$-module,
if $\Gr^{\nbigw}_w(\nbigt)\in\MT(X,w)$
for each $w$.
Let $\MTW(X)$ denote the full subcategory of
pre-mixed twistor $D$-modules
in the category of filtered $\nbigr$-triples.
It is equipped with auto equivalences
$\Sigma^{p,q}(\nbigt,\nbigw)
=(\nbigt,\nbigw)\otimes\bigl(\nbigu(-p,q),\nbigw\bigr)$.
Recall the following lemma in \cite{sabbah2}
(and \cite{mochi7}).
\index{category $\MTW(X)$}

\begin{lem}
\label{lem;11.2.22.20}
The category $\MTW(X)$ is abelian.
Any morphism
$f:(\nbigt_1,\nbigw)\lrarr(\nbigt_2,\nbigw)$
is strict with respect to the weight filtration.
\hfill\qed
\end{lem}

Let $F:X\lrarr Y$ be a projective morphism.
Let $(\nbigt,\nbigw)\in \MTW(X)$.
The $\nbigr$-triple $F_{\dagger}^i\nbigt$
is equipped with an induced filtration
$\nbigw$,
i.e.,
$\nbigw_k(F_{\dagger}^i\nbigt)$
are the image of
$F_{\dagger}^i(\nbigw_{k-i}\nbigt)
\lrarr F^i_{\dagger}(\nbigt)$.
The following lemma can be found
in \cite{saito2} essentially.
\begin{prop}
\label{prop;11.2.22.31}
$(F_{\dagger}^i\nbigt,\nbigw)
\in \MTW(Y)$.
\end{prop}
\pf
If $\nbigt\in\MT(X,w)$,
we know
that $F_{\dagger}^i\nbigt
 \in \MT(Y,w+i)$
(\cite{saito1}, \cite{sabbah2}, \cite{mochi7}).
According to \cite{sabbah2},
we have a spectral sequence
$E_1^{-i,i+j}=F_{\dagger}^j\Gr^W_i(\nbigt)
\Longrightarrow
 F_{\dagger}^j\nbigt$.
Because
$F_{\dagger}^j\Gr^W_i(\nbigt)
\in\MT(Y,w+i+j)$,
we have
$E_2^{-i,i+j}\in \MT(Y,w+i+j)$,
and the spectral sequence degenerates.
\hfill\qed

\subsubsection{Some abelian categories}

Let $\MTW(X,\bullet)$ denote the category of
objects $(\nbigt,\nbigw)\in\MTW(X)$ 
equipped with a morphism $N:(\nbigt,\nbigw)
 \lrarr(\nbigt,\nbigw)\otimes\newTate(-1)$.
\index{category $\MTW(X,\bullet)$}
Morphisms in the category are naturally defined.
It is clearly an abelian category.
For $w\in\seisuu$,
let $\MTN(X,w)\subset\MTW(X,\bullet)$ be 
the full subcategory of objects $(\nbigt,\nbigw,N)$
such that $\nbigw=M(N)[-w]$.
\index{category $\MTN(X,w)$}
\begin{lem}
\label{lem;10.11.12.1}
The category $\MTN(X,w)$ is abelian.
Let $(\nbigt_i,\nbigw,N)$ be objects 
in $\MTN(X,w_i)$ with $w_1>w_2$,
and let $F:(\nbigt_1,\nbigw,N)
\lrarr(\nbigt_2,\nbigw,N)$
be a morphism in $\MTW(X,\bullet)$.
Then, $F=0$.
\end{lem}
\pf
It is easy to deduce the first claim
from Lemma \ref{lem;11.2.22.20}.
As for the second claim,
because the induced morphism
$\Gr^{\nbigw}(F)$ is $0$,
we obtain $F=0$
by Lemma \ref{lem;11.2.22.20}.
\hfill\qed

\vspace{.1in}
Let $\MTW(X,\bullet)^{\RMF}\subset
\MTW(X,\bullet)^{\fil}$
denote the full subcategory of the objects
$(\nbigt,\nbigw,L,N)$ such that 
$\Gr^L_w(\nbigt,\nbigw,N)\in
 \MTN(X,w)$.
\index{category $\MTW(X,\bullet)^{\RMF}$}
We obtain the following lemma 
from Lemma \ref{lem;10.11.12.1}.
\begin{lem}
\label{lem;10.11.12.2}
$\MTW(X,\bullet)^{\RMF}$ is abelian.
Any morphism in $\MTW(X,\bullet)^{\RMF}$
is strict with respect to the filtration $L$.
\hfill\qed
\end{lem}

\subsection{Admissible specializability for pre-mixed twistor $D$-modules}

For a holomorphic function $g$,
let $\MTW^{\sp}(X,g)\subset\MTW(X)$ 
denote the full subcategory of $(\nbigt,L)\in\MTW(X)$
which are admissibly specializable along $g$.
\index{category $\MTW^{\sp}(X,g)$}
We obtain the following lemma 
from Proposition {\rm\ref{prop;10.10.3.2}}.
\begin{lem}
Let $(\nbigt,L)\in\MTW(X)$.
We have
$(\nbigt,L)\in\MTW^{\sp}(X,g)$,
if only if 
the underlying $\nbigr$-modules are filtered strictly
specializable, and 
satisfy the conditions (P1--2) 
in Definition {\rm\ref{df;10.10.3.1}}.
\hfill\qed
\end{lem}

For $(\nbigt,L)\in\MTW^{\sp}(X,g)$,
we have the filtrations of
$\psitilde_{g,\gminia,u}(\nbigt)$
and $\phi_g(\nbigt)$
naively induced by $L$.
They are also denoted by $L$.
Let $W$ denote the filtrations of
$\psitilde_{g,\gminia,u}(\nbigt)$
and 
$\phi_g(\nbigt)$
obtained as the relative monodromy filtration
of $\nbign$ with respect to $L$.

\begin{lem}
\label{lem;11.2.22.23}
$\bigl(
 \psitilde_{g,\gminia,u}(\nbigt),W,L,\nbign\bigr)$
and 
$\bigl(\phi_g(\nbigt),W,L,\nbign\bigr)$
are objects in the category $\MTW(X,\bullet)^{\RMF}$.
\end{lem}
\pf
We have only to prove that 
$\bigl(
 \psitilde_{g,\gminia,u}(\nbigt),W\bigr)$
and 
$\bigl(\phi_g(\nbigt),W\bigr)$
are objects in $\MTW(X)$,
which follows from the canonical splitting
$\Gr^W\simeq\Gr^W\Gr^{L}$.
\hfill\qed

\subsubsection{Sub-quotients}

Following M. Saito \cite{saito2},
we prove that the admissible specializability
is stable for sub-quotients in 
$\MTW(X)$.
We consider an exact sequence
in $\MTW(X)$:
\[
 0\lrarr
 (\nbigt_1,L)
\lrarr 
 (\nbigt_2,L)
\lrarr
 (\nbigt_3,L)
\lrarr 0
\]
Let $g$ be any holomorphic function on $X$.
\begin{prop}
\label{prop;13.7.29.20}
If $(\nbigt_2,L)$
is admissibly (resp. filtered strictly)
specializable along $g$,
$(\nbigt_i,L)$ $(i=1,3)$ are also 
admissibly (resp. filtered strictly)
specializable along $g$.
\end{prop}
\pf
Let $\nbigt_i=(\nbigm_i',\nbigm_i'',C_i)$.
We may assume that $X=X_0\times\cnum_t$ 
and $g=t$.
We shall prove 
the claim by an induction on the length of $L$
together with the following claim.
\begin{lem}
\label{lem;13.7.29.21}
For any $\lambda_0\in\cnum$,
the followings complexes are exact:
\[
 0\lrarr
 \Vzero_a(\nbigm_3')
\lrarr
 \Vzero_a(\nbigm_2')
\lrarr
\Vzero_a(\nbigm_1')\lrarr 0
\]
\[
 0\lrarr
 \Vzero_a(\nbigm_1'')
\lrarr
 \Vzero_a(\nbigm_2'')
\lrarr
\Vzero_a(\nbigm_3'')
\lrarr 0
\]
\end{lem}

If $L$ is pure, the claim 
of Proposition \ref{prop;13.7.29.20}
is clear.
Because the category of polarizable pure 
twistor $D$-modules are semisimple,
the claim of Lemma \ref{lem;13.7.29.21} also holds.

We may assume that
$L_0(\nbigt_i)=\nbigt_i$,
and that the claims hold for 
$L_{-1}(\nbigt_i)$.
Suppose that 
$(\nbigt_2,L)$ is filtered strictly specializable along $t$.
Let $\lambda_0\in\cnum$.
We define $\Vzero_a(\nbigm_3'')$ as the image of
$\Vzero_a(\nbigm''_2)\lrarr\nbigm''_3$.
By the construction,
we have the natural surjection
$\Vzero_a(\nbigm''_3)\lrarr
\Vzero_a\Gr^L_0(\nbigm''_3)$.
We set
$V^{\prime(\lambda_0)}_a(L_{-1}\nbigm''_3):=
 \Vzero_a(\nbigm''_3)\cap L_{-1}\nbigm_3''$.
By construction,
we have
$\Vzero_a(L_{-1}\nbigm''_3)
\subset
 V^{\prime(\lambda_0)}_a(L_{-1}\nbigm''_3)$.
Because 
$V^{\prime(\lambda_0)}$ is a monodromic filtration 
by $V_0\nbigr_{X_0\times\cnum_t}$-coherent subsheaves,
we have
$\Vzero_a(L_{-1}\nbigm''_3)
\supset
 V^{\prime(\lambda_0)}_a(L_{-1}\nbigm''_3)$.
Hence, 
\[
 0\lrarr
 \Vzero_a(L_{-1}\nbigm''_3)
\lrarr
 \Vzero_a(\nbigm''_3)
\lrarr
 \Vzero_a(\Gr^L_0\nbigm_3'')
\lrarr 0
\]
is exact for any $a$.
By considering any ramified exponential twist,
we obtain that 
$(\nbigm_3'',L)$ is filtered strictly specializable
along $t$.

We define
$\Vzero_a(\nbigm_1''):=
 \nbigm_1''\cap
 \Vzero_a(\nbigm_2'')$.
By an easy diagram chasing,
we obtain that
$0\lrarr \Vzero_a(L_{-1}\nbigm_1'')
\lrarr
 \Vzero_a(\nbigm_1'')
\lrarr
 \Vzero_a(\Gr^L_0\nbigm_1'')
\lrarr 0$
is exact for any $a$.
Hence, 
by considering any ramified exponential twist,
we obtain that
$(\nbigm_1'',L)$ is also 
filtered strictly specializable along $t$.
We can prove the filtered strictly specializability
along $t$ for 
$(\nbigm'_i,L)$ $(i=1,3)$ in a similar way.
The claim of Lemma \ref{lem;13.7.29.21}
is also proved.

\vspace{.1in}

Suppose that
$(\nbigt_2,L)$ is admissibly specializable along $t$.
We have already known that
$(\nbigt_i,L)$ $(i=1,3)$ are filtered
strictly specializable along $t$.
We have the following commutative diagram
(the Tate twist is denoted by $(\ell)$):
{\tiny
\[
\begin{array}{ccc}
\Ker(N^{\ell}:
 \Gr^L_0\psitilde_{t,\gminia,u}(\nbigt_2)(\ell)
 \rarr
 \Gr^L_0\psitilde_{t,\gminia,u}(\nbigt_2)
 )
 &\stackrel{\gamma_2}{\lrarr}&
 \frac{L_{-1}\psitilde_{t,\gminia,u}(\nbigt_2)}
 {N^{\ell}L_{-1}\psitilde_{t,\gminia,u}(\nbigt_2)(\ell)
+ W_{-\ell-1}L_{-1}L_{-1}\psitilde_{t,\gminia,u}(\nbigt_2)}
 \\
 \gamma_1\darr & & \darr \\
\Ker(N^{\ell}:
 \Gr^L_0\psitilde_{t,\gminia,u}(\nbigt_3)(\ell)
 \rarr
 \Gr^L_0\psitilde_{t,\gminia,u}(\nbigt_3)
 )
 &\stackrel{\gamma_3}{\lrarr}&
 \frac{L_{-1}\psitilde_{t,\gminia,u}(\nbigt_3)}
 {N^{\ell}L_{-1}\psitilde_{t,\gminia,a}(\nbigt_3)(\ell)
+ W_{-\ell-1}L_{-1}L_{-1}\psitilde_{t,\gminia,u}(\nbigt_3)}
\end{array}
\]
}
By Proposition \ref{prop;10.11.5.1},
we have $\gamma_1=0$.
Because $\gamma_2$ is an epimorphism,
we obtain that $\gamma_3=0$.
Hence, by Proposition \ref{prop;10.11.5.1},
we obtain the existence of
 a relative monodromy filtration
of $N$ on $(\psitilde_{t,\gminia,u}(\nbigt_3),L)$
Similarly,
we obtain the existence of
a relative monodromy filtration of
$N$ on $(\phi^{(0)}_t(\nbigt_3),L)$.
Thus, we obtain that
$(\nbigt_3,L)$ is admissibly specializable
along $t$.
By using the Hermitian adjoint,
we obtain that $(\nbigt_1,L)$
is also admissibly specializable along $t$.
\hfill\qed

\vspace{.1in}
As an immediate corollary,
we obtain the following proposition.
\begin{prop}
\label{prop;10.10.2.12}
$\MTW^{\sp}(X,g)$ is an abelian subcategory
of $\MTW(X)$.
\hfill\qed
\end{prop}

\subsubsection{The associated gluing data}

Let $(\nbigt,L)\in\MTW^{\sp}(X,g)$.
Recall Proposition \ref{prop;10.12.19.1}.
We have the following object
in $\Glu(\MTW(X),\vecSigma)^{\fil}$,
which is filtered $S$-decomposable:
\[
 \begin{CD}
 \bigl(
 \Sigma^{1,0}
 \psitilde_{g,-\vecdelta}(\nbigt),
 L \bigr)
\lrarr
  \bigl(\phi^{(0)}_{g}(\nbigt),L\bigr)
\lrarr
 \bigl(
  \Sigma^{0,-1}
 \psitilde_{g,-\vecdelta}(\nbigt),
 L \bigr)
\end{CD}
\]
We also have the object in 
$\Glu\bigl(\MTW(X)^{\fil},\vecSigma\bigr)$:
\[
 \begin{CD}
 \Sigma^{1,0}
 \bigl(
 \psitilde_{g,-\vecdelta}(\nbigt),W
 \bigr)
\lrarr
  \bigl(\phi^{(0)}_{g}(\nbigt),W\bigr)
\lrarr
  \Sigma^{0,-1} \bigl(
 \psitilde_{g,-\vecdelta}(\nbigt),W
 \bigr)
\end{CD}
\]
We remark that 
the filtration $L\phi_g^{(0)}(\nbigt)$
is the transfer of
$L\psitilde_{g,-\vecdelta}(\nbigt)$.

\subsubsection{Independence from functions}

\begin{lem}
Let $A$ be a nowhere vanishing holomorphic function.
We set $g_1:=Ag$.
\begin{itemize}
\item
We have
$\MTW^{\sp}(X,g_1)
=\MTW^{\sp}(X,g)$.
\item
For any $(\nbigt,L)\in\MTW^{\sp}(X,g)$,
we have natural isomorphisms
\begin{equation}
 \label{eq;13.5.9.11}
 \Gr^W\psitilde_{g,-\vecdelta}(\nbigt)
\simeq
 \Gr^W\psitilde_{g_1,-\vecdelta}(\nbigt),
\quad
 \Gr^W\phi_{g}(\nbigt)
\simeq
 \Gr^W\phi_{g_1}(\nbigt),
\end{equation}
which preserve the natural nilpotent morphisms.
\end{itemize}
\end{lem}
\pf
It follows from Lemma \ref{lem;13.5.9.10}.
\hfill\qed

\vspace{.1in}
Let $\nbigt\in\MTW^{\sp}(X,g)$.
Let $\nbigs=\bigoplus\nbigs_k$
be a polarization of $\Gr^L(\nbigt)$.
We have the natural polarizations of
$\Gr^W\psi_{g,-\vecdelta}(\nbigt)
\simeq
 \Gr^W\psi_{g,-\vecdelta}(\Gr^L\nbigt)$
and 
$\Gr^W\phi_{g}(\nbigt)
\simeq
 \Gr^W\phi_{g}(\Gr^L\nbigt)$
induced by $\nbigs$ and 
the natural nilpotent maps $\nbign$
on $\psitilde_{g,-\vecdelta}(G^L\nbigt)$
and $\phi_{g}(G^L\nbigt)$.

\begin{cor}
The induced polarizations are preserved
by the isomorphisms {\rm(\ref{eq;13.5.9.11})}.
\hfill\qed
\end{cor}

\subsubsection{Admissible specializability along hypersurfaces}

Let $H$ be an effective divisor of $X$.
Let $\MTW^{\sp}(X,H)$
denote the full subcategory of 
$(\nbigt,L)\in\MTW(X)$
which are admissibly specializable along $H$
(Definition \ref{subsection;13.5.9.12}).
We obtain the following from Proposition 
\ref{prop;10.10.2.12}.
\begin{prop}
$\MTW^{\sp}(X,H)$ is an abelian subcategory of
$\MTW(X)$.
\hfill\qed
\end{prop}

As in \S\ref{subsection;13.5.9.14},
we obtain graded $\nbigr_X$-triples:
\[
  \bigl(
\Gr^{W}\psitilde_{H,-\vecdelta}
 \bigr)(\nbigt),
\quad
  \bigl(
\Gr^{W}\phi_H
 \bigr)(\nbigt).
\]
They are equipped with the nilpotent maps
\[
 \nbign:
  \bigl(
\Gr^{W}_k\psitilde_{H,-\vecdelta}
 \bigr)(\nbigt)
\lrarr
  \bigl(
\Gr^{W}_{k-2}\psitilde_{H,-\vecdelta}
 \bigr)(\nbigt)\otimes\newTate(-1),
\]
\[
 \nbign:
  \bigl(
\Gr^{W}_k\phi_{H}
 \bigr)(\nbigt)
\lrarr
  \bigl(
\Gr^{W}_{k-2}\phi_{H}
 \bigr)(\nbigt)\otimes\newTate(-1).
\]
They are also equipped with the filtration $L$
induced by the weight filtration of $L$,
and the canonical splitting:
\[
   \bigl(
\Gr^{W}\psitilde_{H,-\vecdelta}
 \bigr)(\nbigt)
\simeq
  \bigl(
\Gr^{W}\psitilde_{H,-\vecdelta}
 \bigr)(\Gr^L\nbigt),
\quad
  \bigl(
\Gr^{W}\phi_H
 \bigr)(\nbigt)
\simeq
   \bigl(
\Gr^{W}\phi_H
 \bigr)(\Gr^L\nbigt).
\]
A polarization 
$\nbigs=\bigoplus \nbigs_k$
of $\Gr^L(\nbigt)$
with $\nbign$
canonically induces Hermitian adjoint morphisms
of $\Gr^W_m\psitilde_{H,-\vecdelta}\Gr^L_k \nbigt$
and 
$\Gr^W_m\phi_{H}\Gr^L_k \nbigt$.
\begin{prop}
\label{prop;13.5.10.102}
$\bigl(\Gr^W_m\psitilde_{H,-\vecdelta}\bigr)
 \bigl(
 \Gr^L_k \nbigt\bigr)$
and 
$\bigl(\Gr^W_m\phi_{H}\bigr)
 \bigl(
\Gr^L_k \nbigt\bigr)$
are polarizable pure twistor $D$-modules
of weight $m$.
The induced Hermitian adjoint morphisms
are polarizations.
\end{prop}
\pf
Once we construct global objects,
we have only to check the local conditions,
and it follows from
Lemma \ref{lem;11.2.22.23}.
\hfill\qed

\begin{cor}
$\bigl(\Gr^W_m\psitilde_{H,-\vecdelta}\bigr)(\nbigt)$
and 
$\bigl(\Gr^W_m\phi_{H}\bigr)(\nbigt)$
are polarizable wild pure twistor $D$-modules
of weight $m$.
\hfill\qed
\end{cor}

\subsection{Admissible specializability and push-forward}

Let $F:X\lrarr Y$ be a projective morphism.
Let $g_Y$ be a holomorphic function on $Y$.
We put $g_X:=F^{-1}(g_Y)$.
Let $(\nbigt,L)\in \MTW^{\sp}(X,g_X)$.
\begin{lem}
\label{lem;10.10.2.10}
We have 
$F_{\dagger}^{\bullet}(\nbigt,L)
 \in \MTW^{\sp}(Y,g_Y)$.
\end{lem}
\pf
Let $W$ denote the relative monodromy filtrations of
$\psitilde_{g_X,\gminia,u}(\nbigt)$
and $\phi_{g_X}(\nbigt)$.
We have
$\bigl(
 \psitilde_{g_X,\gminia,u}(\nbigt),W
\bigr)\in\MTW(X)$
and
$\bigl(\phi_{g_X}(\nbigt),W\bigr)
\in \MTW(X)$.
They induce filtrations of
$F_{\dagger}^{\bullet}\bigl(
 \psitilde_{g_X,\gminia,u}(\nbigt)\bigr)$
and $F_{\dagger}^{\bullet}
 \bigl(\phi_{g_X}(\nbigt)\bigr)$,
which are also denoted by $W$.
We have
$\bigl(
 F_{\dagger}^{\bullet}\bigl(
 \psitilde_{g_X,\gminia,u}(\nbigt)\bigr),W
\bigr)\in\MTW(Y)$
and $\bigl(
F_{\dagger}^{\bullet}
 \bigl(\phi_{g_X}(\nbigt)\bigr),W
\bigr)\in\MTW(Y)$.
In particular, they are strict.
Hence, we obtain that
$F_{\dagger}(\nbigt)$ is
strictly specializable along $g_X$
with ramified exponential twist,
and we have
$\psitilde_{g_Y,\gminia,u}
 F_{\dagger}^{\bullet}\nbigt
\simeq
 F_{\dagger}^{\bullet}\psitilde_{g_X,\gminia,u}\nbigt$
and
$\phi_{g_Y}F_{\dagger}^{\bullet}\nbigt
\simeq
 F_{\dagger}^{\bullet}\phi_{g_X}\nbigt$.
Because the complex
\[
  F_{\dagger}^{\bullet}
 \psitilde_{g_X,\gminia,u}
 \bigl(L_k\nbigt\bigr)
\lrarr
 F_{\dagger}^{\bullet}
 \psitilde_{g_X,\gminia,u}
 (L_j\nbigt)
\lrarr
  F_{\dagger}^{\bullet}
 \psitilde_{g_X,\gminia,u}\bigl(
 L_j\nbigt/L_k\nbigt\bigr)
\]
is exact,
we obtain that
$F_{\dagger}^{\bullet}(\nbigt)$
is filtered strictly specializable.

The morphism
$\bigl(
 \psitilde_{g,\gminia,u}
 L_jF^{\bullet}_{\dagger}(\nbigt),W\bigr)
\lrarr
 \bigl(
 \psitilde_{g,\gminia,u}
 \Gr^L_jF^{\bullet}_{\dagger}(\nbigt),W
 \bigr)$
is a morphism in $\MTW(Y)$,
and hence it is strict with respect to $W$.
Recall that the filtration $W$ on
$\psitilde_{g,\gminia,u}
 \Gr^L_jF^{\bullet}_{\dagger}(\nbigt)$
is equal to $M(N)[j]$.
Hence, the filtration $W$ on
$\psitilde_{g,\gminia,u}
 F^{\bullet}_{\dagger}(\nbigt)$
is equal to the relative monodromy filtration.
Similarly
$(\phi_{g}F^{\bullet}_{\dagger}(\nbigt,L),N)$
has a relative monodromy filtration.
Therefore,
$F_{\dagger}^{\bullet}(\nbigt,L)$
is admissibly specializable.
\hfill\qed

\begin{cor}
Let $H_Y$ be a hypersurface of $Y$.
We put $H_X:=F^{-1}(H_Y)$.
For any $(\nbigt,L)\in \MTW^{\sp}(X,H_X)$,
we have
$F^{\bullet}_{\dagger}(\nbigt,L)
\in \MTW^{\sp}(Y,H_Y)$.
\hfill\qed
\end{cor}

\subsection{Gluing along a coordinate function}
\label{subsection;11.2.22.22}

Let $t$ be a coordinate function,
i.e., it is a holomorphic function such that
$dt$ is nowhere vanishing.
Let $(\nbigt,L)\in\MTW^{\sp}(X,t)$.
We have
$\bigl(\psitilde_{t,-\vecdelta}(\nbigt),W,L,\nbign\bigr)$
in $\MTW(X,\bullet)^{\RMF}$.
Let $(\nbigq,W)\in \MTW(X)$ with morphisms
\[
\begin{CD}
 \Sigma^{1,0}(\psitilde_{t,-\vecdelta}(\nbigt),W)
@>{u}>>
 (\nbigq,W)
@>{v}>>
 \Sigma^{0,-1}(\psitilde_{t,-\vecdelta}(\nbigt),W)
\end{CD}
\]
such that
(i) $\nbign=v\circ u$,
(ii) $\Supp\nbigq\subset\{t=0\}$.
We obtain an $\nbigr_X$-triple
$\nbigt':=\Glue(\nbigt,\nbigq,u,v)$,
as explained in \S\ref{subsection;11.2.22.21}.
Let us observe that it is naturally equipped with
a filtration $L$ such that 
$(\nbigt',L)\in\MTW(X)$.
Recall that
$\psitilde_{t,-\vecdelta}(\nbigt)$ is equipped
with the naively induced filtration $L$.
We obtain a filtration of $\nbigq$
obtained as the transfer of 
$L\psitilde_{t,-\vecdelta}(\nbigt)$
with respect to $(u,v)$,
which is also denoted by $L$.
Note that
$\Xi^{(0)}(\nbigt)$ has the naively induced filtration $L$,
and the natural morphisms
$\bigl( \psi^{(1)}_t(\nbigt),L\bigr)
\lrarr
 \bigl(\Xi^{(0)}_t(\nbigt),L\bigr)
\lrarr
 \bigl(\psi^{(0)}_t(\nbigt),L\bigr)$
are strict with respect to $L$.
Then, we obtain a filtered $\nbigr$-triple $(\nbigt',L)$
as the cohomology of the following complex:
\[
\begin{CD}
 \bigl(
 \psi^{(1)}_t(\nbigt),L
 \bigr)
@>>>
 \bigl(
 \Xi^{(0)}_t(\nbigt),L
 \bigr)
\oplus
 (\nbigq,L)
@>>>
 \bigl(
 \psi^{(0)}_t(\nbigt),L
 \bigr)
\end{CD}
\]
By construction, $\Gr^L_w(\nbigt')$
are strictly $S$-decomposable along $t$.
The components whose supports are contained in $\{t=0\}$
are isomorphic to direct summands of $\Gr^W_w(\nbigq)$.
The components whose supports are not contained in $\{t=0\}$
are isomorphic to direct summands of $\Gr^W_w(\nbigt)$.
Hence, they are objects in $\MT(X,w)$,
i.e.,
$(\nbigt',L)\in\MTW(X)$.
Moreover,
by construction,
it is easy to observe that
$(\nbigt',L)\in\MTW^{\sp}(X,t)$.

\vspace{.1in}
In particular, $(\nbigt,L)$ 
is reconstructed as the cohomology of the following
complex in filtered $\nbigr$-triples:
\begin{equation}
\label{eq;11.2.22.30}
 \begin{CD}
 \bigl(\psi^{(1)}(\nbigt),L\bigr)
@>>>
 \bigl(\phi^{(0)}(\nbigt),L\bigr)
\oplus
 \bigl(\Xi^{(0)}(\nbigt),L\bigr)
@>>>
 \bigl(\psi^{(0)}(\nbigt),L\bigr)
 \end{CD}
\end{equation}

\subsection{Localization}
\label{subsection;11.2.22.100}

\subsubsection{The case of coordinate functions}

Let $t$ be a coordinate function.
Let $(\nbigt,L)\in\MTW(X,t)^{\sp}$.
We apply the gluing construction in 
\S\ref{subsection;11.2.22.22}
to some special cases.
We obtain a filtered $\nbigr$-triple
$(\nbigt[!t],\Ltilde)$  
as the cohomology of the following:
\index{filtered $\nbigr$-triple $(\nbigt[\bikkuri t],\Ltilde)$}
\[
\begin{CD}
\bigl(
 \psi^{(1)}(\nbigt),L
\bigr)
@>>>
 \bigl(\Xi^{(0)}(\nbigt),L\bigr)
 \oplus
 \bigl(\psi^{(1)}(\nbigt),\Nhat_!L\bigr)
@>>>
\bigl(
 \psi^{(0)}(\nbigt),L
\bigr)
\end{CD}
\]
We obtain a filtered $\nbigr$-triple
$(\nbigt[\ast t],\Ltilde)$
as the cohomology of the following complex:
\index{filtered $\nbigr$-triple $(\nbigt[\ast t],\Ltilde)$}
\[
\begin{CD}
\bigl(
 \psi^{(1)}(\nbigt),L
\bigr)
@>>>
 \bigl(\Xi^{(0)}(\nbigt),L\bigr)
 \oplus
 \bigl(\psi^{(0)}(\nbigt),\Nhat_{\ast}L\bigr)
@>>>
\bigl(
 \psi^{(0)}(\nbigt),L
\bigr)
\end{CD}
\]
By Lemma \ref{lem;11.2.22.23},
we have
$(\nbigt,L)[\star t]:=(\nbigt[\star t],\Ltilde)
\in\MTW(X,t)^{\sp}$
for $\star=\ast,!$.

\begin{lem}
\label{lem;11.2.22.24}
We have natural morphisms
$(\nbigt,L)[!t]\lrarr(\nbigt,L)
 \lrarr(\nbigt,L)[\star t]$
in $\MTW(X,t)^{\sp}$.
\end{lem}
\pf
We have the following:
\[
 \begin{CD}
 \bigl(
 \psi^{(1)}(\nbigt),\Ltilde
 \bigr)
 @>{=}>>
 \bigl(
 \psi^{(1)}(\nbigt),\Ltilde
 \bigr)
 @>>>
  \bigl(
 \psi^{(0)}(\nbigt),\Ltilde
 \bigr)
 \\
 @V{=}VV @VVV @V{=}VV \\
 \bigl(
 \psi^{(1)}(\nbigt),\Ltilde
 \bigr)
 @>>>
 \bigl(
 \phi^{(0)}(\nbigt),\Ltilde
 \bigr)
 @>>>
  \bigl(
 \psi^{(0)}(\nbigt),\Ltilde
 \bigr)
 \\
 @V{=}VV @VVV @V{=}VV \\
 \bigl(
 \psi^{(1)}(\nbigt),\Ltilde
 \bigr)
 @>>>
 \bigl(
 \psi^{(0)}(\nbigt),\Ltilde
 \bigr)
 @>{=}>>
  \bigl(
 \psi^{(0)}(\nbigt),\Ltilde
 \bigr)
 \end{CD}
\]
Hence, we have the following:
\[
 \begin{CD}
 \bigl(
 \psi^{(1)}(\nbigt),L
 \bigr)
 @>{=}>>
 \bigl(
 \psi^{(1)}(\nbigt),\Nhat_!L
 \bigr)
 @>>>
  \bigl(
 \psi^{(0)}(\nbigt),L
 \bigr)
 \\
 @V{=}VV @VVV @V{=}VV \\
 \bigl(
 \psi^{(1)}(\nbigt),L
 \bigr)
 @>>>
 \bigl(
 \phi^{(0)}(\nbigt),L
 \bigr)
 @>>>
  \bigl(
 \psi^{(0)}(\nbigt),L
 \bigr)
 \\
 @V{=}VV @VVV @V{=}VV \\
 \bigl(
 \psi^{(1)}(\nbigt),L
 \bigr)
 @>>>
 \bigl(
 \psi^{(0)}(\nbigt),\Nhat_{\ast}L
 \bigr)
 @>{=}>>
  \bigl(
 \psi^{(0)}(\nbigt),L
 \bigr)
 \end{CD}
\]
Then, we obtain the claim of 
Lemma \ref{lem;11.2.22.24}.
\hfill\qed

\vspace{.1in}
We have the following characterization.
\begin{lem}
\label{lem;10.9.2.3}
Let $\star=\ast$ or $!$.
Let 
$(\nbigt,L)\in \MTW^{\sp}(X,t)$.
Assume that
$\nbigt[\star t]=\nbigt$ as an $\nbigr_X$-triple.
Then,
we have
$(\nbigt,L)\simeq
 (\nbigt,L)[\star t]$
in $\MTW(X)$.
\end{lem}
\pf
Let us consider the case $\star=!$.
The other case can be argued similarly.
By the assumption,
we have
$\psi^{(1)}_t(\nbigt) \simeq
 \phi^{(0)}(\nbigt)$
as $\nbigr$-triples.
Hence,
$\bigl(
 \psi^{(1)}_t(\nbigt),W
 \bigr)
\lrarr
 \bigl(
 \phi^{(0)}_t(\nbigt),W
 \bigr)$
is an isomorphism in 
$\MTW(X)$.
We obtain that the naively induced filtration of
$\phi^{(0)}_t(\nbigt)$
is the same as $\Nhat_!L$
under the isomorphism.
Recall that $(\nbigt,L)$ is reconstructed as 
the cohomology of (\ref{eq;11.2.22.30}).
Hence, we obtain that
$(\nbigt,L)\simeq
 (\nbigt,L)[!t]$.
\hfill\qed

\vspace{.1in}
Similarly, 
we obtain the following lemma.
\begin{lem}
\label{lem;10.10.4.20}
Let $\star=\ast$ or $!$.
Let 
$(\nbigt_i,L)\in \MTW^{\sp}(X,t)$ $(i=1,2)$.
We have a natural bijective correspondence
between morphisms
$(\nbigt_1,L)(\ast t)\lrarr(\nbigt_2,L)(\ast t)$
as filtered $\nbigr_X(\ast t)$-triples,
and morphisms
$(\nbigt_1,L)[\star t]\lrarr(\nbigt_2,L)[\star t]$
as $\MTW(X)$.
\end{lem}
\pf
We may assume 
$(\nbigt_i,L)\simeq
 (\nbigt_i,L)[\star t]$.
Let $F:(\nbigt_1,L)(\star t)\lrarr
 (\nbigt_2,L)(\star t)$ be a morphism.
It naturally induces a morphism
$\bigl(
 \psitilde_{t,-\vecdelta}(\nbigt_1),L
 \bigr)
\lrarr
 \bigl(
 \psitilde_{t,-\vecdelta}(\nbigt_2),L
 \bigr)$.
It gives a morphism
$(\psitilde_{t,-\vecdelta}(\nbigt_1),W)
\lrarr
 (\psitilde_{t,-\vecdelta}(\nbigt_2),W)$
in $\MTW(X)$.
By the assumption,
we obtain
$(\phi_t(\nbigt_1),W)
\lrarr
 (\phi_t(\nbigt_2),W)$
in $\MTW(X)$.
Because the naively induced filtrations $L$
on $\phi_t(\nbigt_i)$ are obtained
as the transfer of $L$ of
$\psitilde_{t,-\vecdelta}(\nbigt_i)$,
we obtain
$\bigl(\phi_t(\nbigt_1),L\bigr)
\lrarr
 \bigl(\phi_t(\nbigt_2),L\bigr)$.
Hence, $F$ is compatible with $L$,
i.e., we have $F:(\nbigt_1,L)\lrarr (\nbigt_2,L)$.
\hfill\qed

\begin{cor}
Let $(\nbigt_i,L)\in\MTW^{\sp}(X,t)$
$(i=1,2)$.
We have natural bijections:
\[
 \Hom_{\MTW(X)}\bigl(
 (\nbigt_1,L)[\ast t],\,(\nbigt_2,L)[\ast t]
 \bigr)
\simeq
  \Hom_{\MTW(X)}\bigl(
 (\nbigt_1,L),\,(\nbigt_2,L)[\ast t]
 \bigr)
\]
\[
 \Hom_{\MTW(X)}\bigl(
 (\nbigt_1,L)[! t],\,(\nbigt_2,L)[! t]
 \bigr)
\simeq
  \Hom_{\MTW(X)}\bigl(
 (\nbigt_1,L)[!t],\,(\nbigt_2,L)
 \bigr)
\]
\hfill\qed
\end{cor}

Let $(\nbigt,L)\in\MTW^{\sp}(X,t)$.
We obtain the graded polarizable pure twistor $D$-module
$\Gr^L(\nbigt)
=\bigoplus\Gr^L_k(\nbigt)$.
We have the decomposition
\[
 \Gr^L_k(\nbigt)
=\nbigp_{k0}
\oplus
 \nbigp_{k1},
\]
where $\nbigp_{k0}$ (resp. $\nbigp_{k1}$)
is the direct sum of
pure twistor $D$-modules whose strict supports
are not contained (resp. contained)
in $\{t=0\}$.
\begin{lem}
\label{lem;13.5.9.21}
We have the following natural isomorphisms:
\begin{equation} 
\label{eq;13.5.9.20}
 \Gr^{\Ltilde}_k(\nbigt[\ast t])
=\nbigp_{k0}
\oplus
 \Ker\bigl(
 \Gr_k^W\Gr_k^{\Nhat_{\ast}L}
 \psi_t^{(0)}(\nbigt)
\lrarr
 \Gr_k^W\Gr_k^{L}\psi_t^{(0)}(\nbigt)
 \bigr)
\end{equation}
\begin{equation} 
\label{eq;13.5.9.25}
 \Gr^{\Ltilde}_k(\nbigt[! t])
=\nbigp_{k0}
\oplus
 \Ker\bigl(
 \Gr_k^W\Gr_k^{\Nhat_{!}L}
 \psi_t^{(1)}(\nbigt)
\lrarr
 \Gr_k^W\Gr_k^{L}\psi_t^{(0)}(\nbigt)
 \bigr)
\end{equation}
\end{lem}
\pf
By the construction,
we have the natural isomorphism:
\[
  \Gr^{\Ltilde}_k(\nbigt[\ast t])
=\nbigp_{k0}
\oplus
 \Ker\bigl(
 \Gr_k^{\Nhat_{\ast}L}
 \psi^{(0)}(\nbigt)
\lrarr
 \Gr_k^{L}\psi^{(0)}(\nbigt)
 \bigr)
\]
By the construction, 
the morphism
$\Gr_k^{\Nhat_{\ast}L}
 \psi^{(0)}(\nbigt)
\lrarr
 \Gr_k^{L}\psi^{(0)}(\nbigt)$
is a morphism in $\MTW(X)$,
and 
$\Ker\bigl(
 \Gr_k^{\Nhat_{\ast}L}
 \psi^{(0)}(\nbigt)
\lrarr
 \Gr_k^{L}\psi^{(0)}(\nbigt)
 \bigr)$
is pure of weight $k$.
Hence, we obtain (\ref{eq;13.5.9.20}).
We obtain (\ref{eq;13.5.9.25}) in a similar way.
\hfill\qed

\subsubsection{The case of general holomorphic functions}

Let $g$ be a holomorphic function.
Let $\iota_g:X\lrarr X\times\cnum_t$
be the graph of $g$.
Let $(\nbigt,L)$
in $\MTW^{\sp}(X,g)$.
We obtain
$\iota_{g\dagger}(\nbigt,L)[\star t]$
in 
$\MTW(X\times\cnum_t,t)^{\sp}$.
\begin{df}
We say that
$\nbigt$ is localizable along $g$,
if there exists an object
$\nbigt[\star g]
 \in\MTW^{\sp}(X,g)$ $(\star=\ast,!)$
satisfying
\[
 \iota_{g\dagger}(\nbigt[\star g],L)
\simeq
 \iota_{g\dagger}(\nbigt,L)[\star t]
\]
in $\MTW^{\sp}(X\times\cnum_t,t)$.
\index{localizable}
\hfill\qed
\end{df}
Such $\nbigm[\star g]$ is uniquely determined 
up to canonical isomorphisms,
if it exists.

Let $\MTW^{\loc}(X,g)\subset
 \MTW^{\sp}(X,g)$ denote the full subcategory
of $(\nbigt,L)$ which are localizable along $g$.
It is clearly an abelian subcategory.

We obtain the following from Lemma \ref{lem;10.10.4.20}.
\begin{lem}
\label{lem;11.2.20.5}
Let $\star=\ast$ or $!$.
Let 
$(\nbigt_i,L)\in \MTW^{\loc}(X,g)$ $(i=1,2)$
such that
$\nbigt_i[\star g]=\nbigt_i$.
Then, we have a natural bijective correspondence
between morphisms
$(\nbigt_1,L)(\ast g)\lrarr(\nbigt_2,L)(\ast g)$
as filtered $\nbigr_X(\ast g)$-triples,
and morphisms
$(\nbigt_1,L)\lrarr(\nbigt_2,L)$
in $\MTW(X)$.
In particular,
if moreover $\nbigt_1(\ast g)\simeq\nbigt_2(\ast g)$,
the isomorphism is extended to
$\nbigt_1\simeq\nbigt_2$.
\hfill\qed
\end{lem}

\begin{lem}
\label{lem;11.1.21.30}
For any $(\nbigt,L)\in\MTW^{\loc}(X,g)$,
we have naturally defined morphisms
$(\nbigt,L)[!g]\lrarr(\nbigt,L)
 \lrarr(\nbigt,L)[\ast g]$
in $\MTW^{\loc}(X,g)$.
\hfill\qed
\end{lem}
We also have the following.
\begin{cor}
\label{cor;11.2.22.33}
Let $(\nbigt_i,L)\in\MTW^{\loc}(X,g)$ $(i=1,2)$.
We have natural bijections:
\[
 \Hom_{\MTW(X)}\bigl(
 (\nbigt_1,L)[\ast g],\,(\nbigt_2,L)[\ast g]
 \bigr)
\simeq
  \Hom_{\MTW(X)}\bigl(
 (\nbigt_1,L),\,(\nbigt_2,L)[\ast g]
 \bigr)
\]
\[
 \Hom_{\MTW(X)}\bigl(
 (\nbigt_1,L)[! g],\,(\nbigt_2,L)[! g]
 \bigr)
\simeq
  \Hom_{\MTW(X)}\bigl(
 (\nbigt_1,L)[!g],\,(\nbigt_2,L)
 \bigr)
\]
\hfill\qed
\end{cor}

Let $\nbigt\in\MTW^{\loc}(X,g)$.
We have the decomposition
$\Gr^L_k(\nbigt)
=\nbigp_{k0}\oplus\nbigp_{k1}$,
where $\nbigp_{k0}$
(resp. $\nbigp_{k1}$)
is the direct sum of 
pure twistor $D$-modules
whose strict supports are not contained 
(resp. contained) in $\{g=0\}$.
We obtain the following isomorphisms from
Lemma \ref{lem;13.5.9.21}.
\begin{equation} 
\label{eq;13.5.9.30}
 \Gr^{\Ltilde}_k(\nbigt[\ast g])
=\nbigp_{k0}
\oplus
 \Ker\bigl(
 \Gr_k^W\Gr_k^{\Nhat_{\ast}L}
 \psi_g^{(0)}(\nbigt)
\lrarr
 \Gr_k^W\Gr_k^{L}\psi_g^{(0)}(\nbigt)
 \bigr)
\end{equation}
\begin{equation} 
\label{eq;13.5.9.31}
 \Gr^{\Ltilde}_k(\nbigt[! g])
=\nbigp_{k0}
\oplus
 \Ker\bigl(
 \Gr_k^W\Gr_k^{\Nhat_{!}L}
 \psi_g^{(1)}(\nbigt)
\lrarr
 \Gr_k^W\Gr_k^{L}\psi_g^{(0)}(\nbigt)
 \bigr)
\end{equation}

\begin{lem}
\label{lem;13.5.9.40}
Let $A$ be a nowhere vanishing function on $X$.
We set $g_1:=Ag$.
\begin{itemize}
\item
$\MTW^{\loc}(X,g)=\MTW^{\loc}(X,g_1)$.
\item
For any $\nbigt\in\MTW^{\loc}(X,g)$,
we have natural isomorphisms
$\nbigt[\star g]\simeq
\nbigt[\star g_1]$.
\item
The isomorphisms {\rm(\ref{eq;13.5.9.30})}
and {\rm(\ref{eq;13.5.9.31})}
for $g$ and $g_1$
are equal 
under the isomorphisms in Lemma {\rm\ref{lem;13.5.9.10}}.
\hfill\qed
\end{itemize}
\end{lem}

\subsubsection{The case of hypersurfaces}

Let $H$ be an effective divisor of $X$.
\begin{df}
$\nbigt\in\MTW(X)$ is called 
admissibly specializable (resp. localizable)
along $H$,
if the following holds: 
\begin{itemize}
\item
Let $U\subset X$ be open with a generator
$g_U$ of $\nbigo(-H)_{|U}$.
Then,
$\nbigt_{|U}$ is admissibly specializable
(resp. localizable)
along $g_U$.
\hfill\qed
\end{itemize}
\end{df}
Let $\MTW^{\sp}(X,H)$ 
(resp. $\MTW^{\loc}(X,H)$)
denote the full subcategory of
$(\nbigt,L)\in\MTW(X)$
which are admissibly specializable 
(resp. localizable) along $H$.

For $\nbigt\in\MTW^{\loc}(X,H)$,
we obtain a filtered $\nbigr_X$-triple
$\nbigt[\star H]$
determined by the following condition:
\begin{itemize}
\item
Let $(U,g_U)$ be as above.
Then,
$\nbigt[\star H]_{|U}\simeq
 \nbigt_{|U}[\star g_U]$.
\end{itemize}

\begin{prop}
\label{prop;13.5.9.100}
$\nbigt[\star H]\in\MTW^{\sp}(X,H)$.
\end{prop}
\pf
We have the decomposition
$\Gr^{L}_k(\nbigt)
=\nbigp_{k0}\oplus\nbigp_{k1}$,
where 
$\nbigp_{k0}$ (resp. $\nbigp_{k1}$)
is the direct sum of pure twistor $D$-modules
whose strict supports are not contained
(resp. contained) in $H$.
By Lemma \ref{lem;13.5.9.40},
we obtain the following isomorphism
\begin{equation} 
\label{eq;13.5.9.41}
 \Gr^{\Ltilde}_k(\nbigt[\ast H])
=\nbigp_{k0}
\oplus
 \Ker\Bigl(
 \bigl(
\Gr_k^{\Nhat_{\ast}L}
 \Gr_k^W
 \psi_H^{(0)}
 \bigr)
 (\nbigt)
\lrarr
 \bigl(
\Gr_k^{L}
 \Gr_k^W
\psi_H^{(0)}
\bigr)
 (\nbigt)
 \Bigr)
\end{equation}
\begin{equation} 
\label{eq;13.5.9.42}
 \Gr^{\Ltilde}_k(\nbigt[! H])
=\nbigp_{k0}
\oplus
 \Ker\Bigl(
 \bigl(
 \Gr_k^{\Nhat_{!}L}
 \Gr_k^W
 \psi_H^{(1)}
 \bigr)(\nbigt)
\lrarr
 \bigl(
 \Gr_k^{L}
 \Gr_k^W
\psi_H^{(0)}
\bigr)
 (\nbigt)
 \Bigr)
\end{equation}
Then, $\Ker$ in (\ref{eq;13.5.9.41}) and
(\ref{eq;13.5.9.42}) are direct summands of
polarizable pure twistor $D$-modules
$(\Gr^W\psi_{H}^{(a)})(\nbigt)$
of weight $k$.
Note the semisimplicity of the category of
polarizable pure twistor $D$-modules.
Then, the claim of the proposition follows
from Proposition \ref{prop;13.5.10.102}.
\hfill\qed

\begin{lem}
Let $\nbigt_i\in\MTW^{\loc}(X,H)$ $(i=1,2)$
such that $\nbigt_i[\ast H]\simeq \nbigt_i$.
We have a natural bijective correspondence
of morphisms
$(\nbigt_1,L)(\ast H)\lrarr(\nbigt_2,L)(\ast H)$
as filtered $\nbigr_X(\ast H)$-triples,
with morphisms $(\nbigt_1,L)\lrarr(\nbigt_2,L)$
in $\MTW(X)$.
In particular,
if moreover $\nbigt_1(\ast H)\simeq\nbigt_2(\ast H)$,
the isomorphism is extended to
$\nbigt_1\simeq\nbigt_2$.
\hfill\qed
\end{lem}

\begin{lem}
For any $(\nbigt,L)\in\MTW^{\loc}(X,H)$,
we have naturally defined morphisms
$(\nbigt,L)[!H]\lrarr(\nbigt,L)
 \lrarr(\nbigt,L)[\ast H]$
in $\MTW(X)$.
\hfill\qed
\end{lem}

\begin{cor}
For $(\nbigt_i,L)\in\MTW^{\loc}(X,H)$ $(i=1,2)$,
we have natural bijections:
\[
 \Hom_{\MTW(X)}\bigl(
 (\nbigt_1,L)[\ast H],\,(\nbigt_2,L)[\ast H]
 \bigr)
\simeq
  \Hom_{\MTW(X)}\bigl(
 (\nbigt_1,L),\,(\nbigt_2,L)[\ast H]
 \bigr)
\]
\[
 \Hom_{\MTW(X)}\bigl(
 (\nbigt_1,L)[! H],\,(\nbigt_2,L)[! H]
 \bigr)
\simeq
  \Hom_{\MTW(X)}\bigl(
 (\nbigt_1,L)[!H],\,(\nbigt_2,L)
 \bigr)
\]
\hfill\qed
\end{cor}

\subsection{Integrable case}

An integrable pure twistor $D$-module of weight $w$
\cite{sabbah2}
is an integrable $\nbigr$-triple $\nbigt$
such that 
(i) it is a pure twistor $D$-module of weight $w$,
(ii) it has a polarization, which is integrable
as a morphism $\nbigt\lrarr\nbigt^{\ast}\otimes\newTate(-w)$.
A morphism of integrable pure twistor $D$-module
of weight $w$ is defined to be an integrable
morphism for integrable $\nbigr$-triples.
Let $\MTint(X,w)$ denote the category of
integrable pure twistor $D$-modules of weight $w$. 
\index{category $\MTint(X,w)$}

An integrable filtered $\nbigr_X$-triple $(\nbigt,L)$
is called an integrable pre-mixed twistor $D$-module,
if $\Gr^L_w(\nbigt)\in\MTint(X,w)$ for each $w$.
A morphism is defined to be an integrable morphism
of the underlying filtered $\nbigr$-triples.
Let $\MTWint(X)$ denote the category of
integrable pre-mixed twistor $D$-modules.
\index{category $\MTWint(X)$}
The following proposition is easy to see.
\begin{prop}
$\MTint(X)$ is abelian and semisimple,
and $\MTWint(X)$ is an abelian category.
\end{prop}
\pf
The first claim can be proved
as in the ordinary case.
The second claim immediately follows
from the first one.
\hfill\qed

\vspace{.1in}
The following proposition is an analogue of 
Proposition \ref{prop;11.2.22.31}.
\begin{prop}
Let $F:X\lrarr Y$ be a projective morphism.
For any object $\nbigt\in\MTint(X,w)$,
we have
$F_{\dagger}^i\nbigt\in\MTint(Y,w+i)$.
For $(\nbigt,\nbigw)\in \MTWint(X)$,
we have
$(F_{\dagger}^i\nbigt,\nbigw)
\in \MTWint(Y)$.
\hfill\qed
\end{prop}

\subsubsection{Admissibly specializability,
gluing and localizability}

For a given holomorphic function $g$,
let $\MTW^{\integral \sp}(X,g)\subset\MTWint(X)$ 
denote the full subcategory of $(\nbigt,L)\in\MTWint(X)$
which are admissibly specializable along $g$.
\index{category $\MTW^{\integral \sp}(X,g)$}
As in Proposition \ref{prop;10.10.2.12},
$\MTW^{\integral\sp}(X,g)$ is an abelian category.
We remark that 
the naively induced filtrations on
$\psitilde_{g,\gminia,u}(\nbigt)$
and $\phi_g(\nbigt)$ are integrable.
The induced nilpotent morphisms are integrable.
Hence, the relative monodromy filtrations are also integrable.

Let $t$ be a coordinate function.
Let $(\nbigt,L)\in\MTW^{\integral\sp}(X,t)$.
Then, we have
$\bigl(\psitilde_{t,-\vecdelta}(\nbigt),W\bigr)$
in $\MTWint(X)$
with the naively induced filtration $L$
and the integrable morphism
$\nbign:\bigl(\psitilde_{t,-\vecdelta}(\nbigt),W\bigr)
\lrarr
\bigl(\psitilde_{t,-\vecdelta}(\nbigt),W\bigr)
\otimes\newTate(-1)$.
We have $W=M(\nbign;L)$.
Let $(\nbigq,W)\in \MTWint(X)$ with 
integrable morphisms
\[
\begin{CD}
 \Sigma^{1,0}(\psitilde_{t,-\vecdelta}(\nbigt),W)
@>{u}>>
 (\nbigq,W)
@>{v}>>
 \Sigma^{0,-1}(\psitilde_{t,-\vecdelta}(\nbigt),W)
\end{CD}
\]
such that
(i) $\nbign=v\circ u$,
(ii) $\Supp\nbigq\subset\{t=0\}$.
As explained in 
\S\ref{subsection;11.2.22.22},
we have constructed
$\Glue(\nbigt,\nbigq,u,v)$
with the weight filtration $L$
in $\MTW^{\sp}(X,g)$.
Note that the naively induced filtrations $L$ of
$\psi^{(a)}(\nbigt)$ and its transfer $L$ to $\nbigq$
are integrable.
The naively induced filtration $L$ of
$\Xi^{(a)}(\nbigt)$ is also integrable.
Hence, we have
$\bigl(
 \Glue(\nbigt,\nbigq,u,v),L
\bigr)\in\MTWint(X)$.
In particular, we have the following.
\begin{lem}
\label{lem;11.4.5.1}
For $(\nbigt,L)\in\MTW^{\integral\sp}(X,t)$
and for $\star=\ast,!$,
we naturally have
$(\nbigt,L)[\star t]\in\MTWint(X)$,
and the morphisms
$(\nbigt,L)[!t]\lrarr (\nbigt,L)\lrarr (\nbigt,L)[\ast t]$
in $\MTWint(X)$.
\hfill\qed
\end{lem}

Let $H$ be an effective divisor of $X$.

\begin{lem}
\label{lem;13.5.10.400}
Let $(\nbigt,L)\in\MTW^{\integral}(X)$ 
be localizable along $H$.
Then, $\nbigt[\star H]\in\MTW(X)$ are 
also integrable.
\end{lem}
\pf
By the construction,
the weight filtration $\Ltilde$ of $\nbigt[\star H]$
is integrable.
In the proof of
Proposition \ref{prop;13.5.9.100},
an integrable polarization of
$\Gr^L\nbigt$ induces 
integrable polarizations of
$(\Gr^W\psi_{H}^{(a)})(\nbigt)$.
In (\ref{eq;13.5.9.41}) and
(\ref{eq;13.5.9.42}),
$\Ker$ are integrable.
Hence, we obtain that
$\Gr_kr^{\Ltilde}(\nbigt[\star H])$
have integrable polarizations.
\hfill\qed

\begin{lem}
Let $(\nbigt_i,L)\in \MTW^{\integral}(X)$ 
be localizable along $H$ such that
$\nbigt_i\simeq\nbigt_i[\star H]$.
We have a natural bijective correspondence
of morphisms
$(\nbigt_1,L)(\ast H)\lrarr(\nbigt_2,L)(\ast H)$
as filtered integrable $\nbigr_X(\ast H)$-triples,
with morphisms
$(\nbigt_1,L)\lrarr(\nbigt_2,L)$
in $\MTWint(X)$.
In particular,
if moreover $\nbigt_1(\ast H)\simeq\nbigt_2(\ast H)$,
the isomorphism is extended to
$\nbigt_1\simeq\nbigt_2$.
\hfill\qed
\end{lem}

\begin{lem}
\label{lem;11.2.22.101}
Let $(\nbigt,L)\in\MTWint(X)$ be localizable along $H$.
Then, we have natural morphisms
$(\nbigt,L)[!H]\lrarr(\nbigt,L)
 \lrarr(\nbigt,L)[\ast H]$
in $\MTW^{\integral}(X)$.
\hfill\qed
\end{lem}
The following is an analogue of Corollary
\ref{cor;11.2.22.33}.
\begin{lem}
Let $(\nbigt_i,L)\in\MTW^{\integral}(X)$ $(i=1,2)$
be localizable along $H$.
We have natural bijections:
\[
 \Hom_{\MTWint(X)}\bigl(
 (\nbigt_1,L)[\ast H],\,(\nbigt_2,L)[\ast H]
 \bigr)
\simeq
  \Hom_{\MTWint(X)}\bigl(
 (\nbigt_1,L),\,(\nbigt_2,L)[\ast H]
 \bigr)
\]
\[
 \Hom_{\MTWint(X)}\bigl(
 (\nbigt_1,L)[! H],\,(\nbigt_2,L)[! H]
 \bigr)
\simeq
  \Hom_{\MTWint(X)}\bigl(
 (\nbigt_1,L)[!H],\,(\nbigt_2,L)
 \bigr)
\]
\hfill\qed
\end{lem}

\section{Mixed twistor $D$-modules}

\subsection{Definition}

Let $X$ be a complex manifold.
We define mixed twistor $D$-modules on $X$
by using a Noetherian induction
on the support.

\begin{df}
\label{df;11.4.3.30}
Let $(\nbigt,L)\in \MTW(X)$.
It is called a mixed twistor $D$-module,
if the following holds for any open subset
$U\subset X$
with a holomorphic function $g$:
\begin{itemize}
\item
 $(\nbigt,L)_{|U}$ is admissibly specializable along $g$. 
\item
 $\bigl(
 \psitilde_{g,\gminia,u}(\nbigt),
 M(\nbign;L)
 \bigr)$
 ($\gminia\in\cnum[t_n^{-1}],
 u\in\real\times\cnum$)
 and 
 $\bigl(
 \phi_g(\nbigt),M(\nbign;L)
 \bigr)$
 are mixed twistor $D$-modules
 with strictly smaller supports,
 if $\Supp\nbigt_{|U}\not\subset g^{-1}(0)$.
\end{itemize}
\index{mixed twistor $D$-module}
(We use a Noetherian induction on 
the dimensions and the number of the components
of the maximal dimension.)
\hfill\qed
\end{df}

\begin{rem}
In this paper,
we consider only graded polarizable mixed objects.
\hfill\qed
\end{rem}

\subsection{Some basic property}

Let $\MTM(X)\subset\MTW(X)$ denote the full subcategory of
mixed twistor $D$-modules on $X$.
\index{category $\MTM(X)$}
\begin{prop}
Let $0\lrarr(\nbigt_1,L)\lrarr
 (\nbigt_2,L)
 \lrarr
 (\nbigt_3,L)\lrarr 0$
be an exact sequence in $\MTW(X)$.
If $\nbigt_2\in\MTM(X)$,
then $\nbigt_i$ $(i=1,3)$
are also mixed twistor $D$-modules
\end{prop}
\pf
Let $\Map^{\ast}(\seisuu_{\geq 0},\seisuu_{\geq 0})$ 
be the set of maps
$\rho:\seisuu_{\geq 0}\lrarr\seisuu_{\geq 0}$
such that $\rho(n)=0$ 
except for finitely many $n$.
We use the order on 
$\Map^{\ast}(\seisuu_{\geq 0},\seisuu_{\geq 0})$
given by
$\rho_1<\rho_2$
if and only if
$\rho_1(i_0)<\rho_2(i_0)$
for $i_0=\max\bigl\{i\,\big|\,
 \rho_1(i)\neq\rho_2(i) \bigr\}$.

Let $P\in X$.
Let $X_P$ be any small neighbourhood of $P$.
For any closed analytic subset $\Gamma\subset X_P$,
let $\rho(\Gamma,i)$ be the number of
the $i$-dimensional irreducible components.
We obtain
$\rho(\Gamma)
 \in \Map^{\ast}(\seisuu_{\geq 0},\seisuu_{\geq 0})$.
We use an induction on $\rho(\Gamma)$.

Suppose that 
$\Supp(\nbigt_2)=\Gamma$.
By Proposition \ref{prop;13.7.29.20},
$\nbigt_i$ $(i=1,3)$ are admissibly specializable
along any holomorphic function $g$
on $X_P$.
Suppose $\Gamma\not\subset g^{-1}(0)$.
We obtain the exact sequence in $\MTW(X_P)$
for any
$\gminia\in t_n^{-1}\cnum[t_n^{-1}]$
and $u\in\real\times\cnum$:
\[
 0\lrarr
 \psitilde_{g,\gminia,u}(\nbigt_1,L)
\lrarr
 \psitilde_{g,\gminia,u}(\nbigt_2,L)
\lrarr
 \psitilde_{g,\gminia,u}(\nbigt_3,L)
\lrarr 0
\]
By definition,
$\psitilde_{g,\gminia,u}(\nbigt_2,L)
\in\MTM(X_P)$.
By the hypothesis of the induction,
we obtain that
$\psitilde_{g,\gminia,u}(\nbigt_i,L)
 \in \MTM(X_P)$ $(i=1,3)$.
Similarly,
we obtain that 
$\phi_{g}(\nbigt_i,L)\in\MTM(X)$ $(i=1,3)$.
Hence, we obtain that
$(\nbigt_i,L)\in
 \MTM(X)$ ($i=1,3$).
\hfill\qed

\begin{rem}
It is M. Saito who emphasized
the importance of the fact that
any sub-quotients of mixed Hodge modules 
in the category of weakly mixed Hodge modules
are mixed Hodge modules.
\hfill\qed
\end{rem}

\vspace{.1in}
As an immediate corollary,
we obtain the following proposition.

\begin{prop}
\label{prop;11.2.22.41}
$\MTM(X)$ is abelian.
\hfill\qed
\end{prop}

\begin{prop}
For $(\nbigt,L)\in\MTM(X)$,
$(\nbigt,L)^{\ast}$
and
$j^{\ast}(\nbigt,L)$
are mixed twistor $D$-modules.
\end{prop}
\pf
We have only to use a Noetherian induction.
Note
the compatibility of Hermitian adjoint
and Beilinson's formalism.
\hfill\qed

\begin{prop}
\label{prop;10.11.15.21}
Let $F:X\lrarr Y$ be projective.
For $(\nbigt,L)\in\MTM(X)$,
we have $F^i_{\dagger}(\nbigt,L)
\in\MTM(Y)$.
\end{prop}
\pf
By Lemma \ref{lem;10.10.2.10},
$F_{\dagger}^i(\nbigt,L)$ is admissibly specializable.
Let $P\in Y$.
Let $g_Y$ be any holomorphic function on $Y$,
such that $\Supp F^i_{\dagger}\nbigt
 \not\subset g_Y^{-1}(0)$.
We put $g_X:=F^{-1}(g_Y)$.
Let $(\nbigt,L)\in \MTM(X)$.
By applying the hypothesis of the induction
to $\psitilde_{g_X,\gminia,u}(\nbigt,L)$
and $\phi_{g_X}(\nbigt,L)$,
we obtain that
$\psitilde_{g_Y,\gminia,u}F_{\dagger}^i(\nbigt,L)$
and 
$\phi_{g_Y}F_{\dagger}^i(\nbigt,L)$
are objects in $\MTM(Y)$.
Hence, we obtain
$F_{\dagger}^i(\nbigt,L)
\in\MTM(Y)$.
\hfill\qed

\vspace{.1in}
Let $\iota:Y\subset X$ be a closed immersion.
Let $\MTM_Y(X)\subset \MTM(X)$ 
be the full subcategory of
$(\nbigt,L)\in\MTM(X)$
such that $\Supp\nbigt\subset Y$.
\index{category $\MTM_Y(X)$}
We obtain a functor
$\iota_{\dagger}:
\MTM(Y)\lrarr\MTM_Y(X)$.
\begin{prop}
\label{prop;10.11.12.30}
The functor 
$\iota_{\dagger}:
\MTM(Y)\lrarr\MTM_Y(X)$
is an equivalence.
\end{prop}
\pf
We may assume that
$X=\Delta^n$ and $Y=\{z_1=0\}$.
We use a Noetherian induction on the support.
Let $(\nbigt,L)\in \MTM_Y(X)$.
Because $(\nbigt,L)$ is filtered
strictly specializable along $z_1$,
we have a filtered $\nbigr_Y$-triple
$(\nbigt_0,L)$ such that
$(\nbigt,L)=\iota_{\dagger}(\nbigt_0,L)$.
By Lemma \ref{lem;10.10.4.2},
$(\nbigt_0,L)$ is admissibly specializable.
Let $g$ be a holomorphic function
on $Y$ such that 
$\Supp\nbigt_0\not\subset g^{-1}(0)$.
We extend $g$ to a holomorphic function on $X$
by the pull back via the projection $X\lrarr Y$.
Then, we have
\[
 \iota_{\dagger}\psi_{g,\gminia,u}(\nbigt,L)
\simeq
 \psi_{g,\gminia,u}\iota_{\dagger}(\nbigt,L),
\quad\quad
 \iota_{\dagger}\phi_{g}(\nbigt,L)
\simeq
 \phi_{g}\iota_{\dagger}(\nbigt,L).
\]
By definition,
$\psi_{g,\gminia,u}\iota_{\dagger}(\nbigt,L)$
and
$\phi_{g}\iota_{\dagger}(\nbigt,L)$
are objects in $\MTM(X)$.
Hence, 
by the hypothesis of the induction,
$\psi_{g,\gminia,u}(\nbigt,L)$
and 
$\phi_{g}(\nbigt,L)$ are objects
in $\MTM(Y)$.
Thus, we are done.
\hfill\qed

\begin{rem}
The forgetful functor
$\Xi_{\DR}$ from $\MTM(X)$
to the category of holonomic $D$-modules
is faithful.
Indeed,
let $F:(\nbigt_1,L)\lrarr (\nbigt_2,L)$
be a morphism in $\MTM(X)$
such that $\Xi_{DR}(F)=0$.
By a result in the pure case,
we obtain $\Gr^L(F)=0$.
By the strictness, we obtain $F=0$.
\hfill\qed
\end{rem}

\subsection{Integrable case}

\index{integrable mixed twistor $D$-module}
An integrable pre-mixed twistor $D$-module
is called integrable mixed twistor $D$-module,
if the underlying pre-mixed twistor $D$-module
is a mixed twistor $D$-module.
Let $\MTMint(X)\subset\MTWint(X)$ denote the full
subcategory of integrable mixed twistor $D$-modules.
\index{category $\MTMint(X)$}
The following proposition can be proved 
by a Noetherian induction
as in the case of Proposition \ref{prop;11.2.22.41}.
\begin{prop}
$\MTMint(X)$ is closed for sub-quotients
in $\MTWint(X)$.
In particular,
$\MTMint(X)$ is abelian.
\hfill\qed
\end{prop}

We mention some basic property.
\begin{prop}
Let $(\nbigt,L)\in\MTMint(X)$.
\begin{itemize}
\item
For any open $U\subset X$
with a holomorphic function $g$,
$\bigl(
 \psitilde_{g,\gminia,u}(\nbigt),M(\nbign;L)
 \bigr)$
and
$\bigl(
 \phi_g(\nbigt),M(\nbign;L)
 \bigr)$
are also objects in
$\MTMint(X)$.
\item
$(\nbigt,L)^{\ast}$
and $j^{\ast}(\nbigt,L)$
are objects in
$\MTMint(X)$.
\item
For a projective  morphism $F:X\lrarr Y$,
we have $F^{i}_{\dagger}(\nbigt,L)
\in\MTMint(Y)$.
\hfill\qed
\end{itemize}
\end{prop}

Let $\iota:Y\subset X$ be a closed immersion.
Let $\MTMint_Y(X)\subset \MTMint(X)$ 
be the full subcategory of
$(\nbigt,L)\in\MTMint(X)$
such that $\Supp\nbigt\subset Y$.
We obtain a functor
$\iota_{\dagger}:
\MTMint(Y)\lrarr\MTMint_Y(X)$.
\begin{prop}
The functor 
$\iota_{\dagger}:
\MTMint(Y)\lrarr\MTMint_Y(X)$
is an equivalence.
\hfill\qed
\end{prop}

\chapter{Infinitesimal mixed twistor modules}
\label{section;11.4.3.21}

It is basic to study mixed twistor $D$-modules
with normal crossing singularity.
In this section, we shall study its infinitesimal version
(or linear version),
which is infinitesimal mixed twistor module.
In the Hodge case,
infinitesimal mixed Hodge modules were
introduced and studied in \cite{kashiwara-mixed-Hodge}.
A kind of gluing procedure of infinitesimal mixed Hodge modules
was studied in 
\cite{kashiwara-mixed-Hodge} and \cite{saito2}.
Our main purpose in this chapter
is to explain the twistor version of
the statements and the ideas
in \cite{kashiwara-mixed-Hodge}
and \cite{saito2}.

\section{Preliminary}

\subsection{Mixed twistor structure}
\label{subsection;11.2.22.40}

\index{variation of mixed twistor structure}
\index{mixed twistor structure}

Let $X$ be a complex manifold.
Let $(\nbigt,W)$ be a filtered object
in the category of variations of twistor structure on $X$
(\S\ref{subsection;13.5.3.10}),
where $W$ is a finite complete exhaustive
increasing filtration indexed by $\seisuu$.
It is called a variation of mixed twistor structure on $X$
if each $\Gr^W_w(\nbigt)$ is a polarizable pure twistor
$D$-module of weight $w$.
In this paper,
we often call it 
a mixed twistor structure on $X$,
i.e.,
omit ``variation of'',
if there is no risk of confusion.
It is called pure of weight $w$,
if $\Gr^L_m=0$ unless $m=w$.
Let $\PTS(X,w)\subset\VTS(X)$
denote the full subcategory of polarizable pure 
twistor structure on $X$.
It is abelian and semisimple.
For $\nbigt_i\in\PTS(X,w_i)$ with $w_1>w_2$,
any morphism $\nbigt_1\lrarr\nbigt_2$
in the category of $\nbigr$-triples is $0$.
Let $\MTS(X)\subset\VTS(X)^{\fil}$ denote 
the category of mixed twistor structure on $X$.
\index{category $\MTS(X)$}
It is an abelian category.
For any morphism
$F:(\nbigt_1,W)\lrarr(\nbigt_2,W)$
in $\MTS(X)$,
$F$ is strictly compatible with the filtration $W$.
\begin{rem}
It might be more appropriate that 
the above object is called graded polarizable variation of 
mixed twistor structure.
We omit to distinguish it,
because we consider only graded polarizable ones.
\hfill\qed
\end{rem}

The category $\MTS(X)$ is naturally equipped with
tensor product and inner homomorphism.
It is also equipped with additive auto equivalences
$\Sigma^{p,q}$ given by the tensor product of
$\nbigu(-p,q)$.

\begin{rem}
A mixed twistor structure $(\nbigt,W)$
is often denoted just by $\nbigt$,
if there is no risk of confusion.
For a filtration $L$ of
a mixed twistor structure $\nbigt$ in $\MTS(X)$,
$\Sigma^{p,q}(\nbigt,L)$
and 
$\bigl(\Sigma^{p,q}(\nbigt),L\bigr)$
are also denoted by 
$(\nbigt,L)\otimes\nbigu(-p,q)$ and
$(\nbigt\otimes\nbigu(-p,q),L)$,
respectively.
If $p=q$,
they are also denoted by
$(\nbigt,L)\otimes\newTate(q)$
and $\bigl(\nbigt\otimes\newTate(q),L\bigr)$,
or $\vecT^q(\nbigt,L)$
and $\bigl(\vecT^q\nbigt,L\bigr)$,
respectively.
\hfill\qed
\end{rem}

\begin{lem}
\label{lem;10.11.12.31}
Let $(\nbigt,L)$ be a mixed twistor $D$-module
on $X$.
Assume that the underlying $D$-module is
a flat bundle.
Then, 
$(\nbigt,L)$ comes from 
a variation of mixed twistor structure.
\end{lem}
\pf
In the pure case,
it is easy to prove the claim
by using the correspondence between
pure twistor $D$-modules and 
wild harmonic bundles.
Let us consider the mixed case.
We may assume $X=\Delta^n$.
Let $\nbigt=(\nbigm_1,\nbigm_2,C)$.
By using the result in the pure case,
we obtain that 
the $\nbigr_X$-modules $\nbigm_i$ are smooth.
By using flat frames, 
we can check that
the pairing takes values
in the sheaf of continuous functions
on $\vecS\times X$
which are $C^{\infty}$ in the $X$-direction.
By successive use of $\psi_{z_i,-\vecdelta}$,
we obtain that 
its restriction to $\vecS\times \{P\}$ 
can be extended to holomorphic function
on $\cnum_{\lambda}^{\ast}$.
Hence, we obtain a pairing of
$\nbigm_{1|\cnum_{\lambda}^{\ast}\times X}$
and $\sigma^{\ast}
 \nbigm_{2|\cnum_{\lambda}^{\ast}\times X}$
valued in the sheaf of $C^{\infty}$-functions
on $\cnum_{\lambda}^{\ast}\times X$
which are holomorphic 
in the $\cnum_{\lambda}^{\ast}$-direction.
Hence $(\nbigt,L)$ comes from
a smooth $\nbigr$-triple.
Because the $w$-th graded piece corresponds to
the pure twistor $D$-module of weight $w$,
we obtain that 
$(\nbigt,L)$ comes from
a variation of mixed twistor structure.
\hfill\qed

\subsection{Reduction}

We shall use the following lemma
implicitly.
\begin{lem}
Let $(\nbigt,W)\in\MTS(X)$
with subobjects $(\nbigt_i,W)\subset (\nbigt,W)$
$(i=1,2)$.
If $\Gr^W(\nbigt_1)=\Gr^W(\nbigt_2)$
in $\Gr^W(\nbigt)$,
then we have $\nbigt_1=\nbigt_2$.
\end{lem}
\pf
Consider
$F:(\nbigt_1,W)\lrarr (\nbigt,W)/(\nbigt_2,W)$.
If $\Gr^W(F)=0$, we have
$F=0$,
which implies the claim of the lemma.
\hfill\qed

\vspace{.1in}

Let $\bigl((\nbigt,W),L,N\bigr)\in\MTS(X)^{\fil,\nil}$.
We put $\nbigt^{(0)}:=\Gr^W(\nbigt)$.
It is equipped with naturally induced
filtrations $W^{(0)}$ and $L^{(0)}$.
We also have an induced map
$N^{(0)}:(\nbigt^{(0)},W^{(0)})\lrarr
(\nbigt^{(0)},W^{(0)})\otimes\newTate(-1)$.
Thus, we obtain
$\bigl(
 (\nbigt^{(0)},W^{(0)}),L^{(0)},N^{(0)}
 \bigr)$
in $\MTS(X)^{\fil,\nil}$.

\begin{lem}
\label{lem;10.7.22.10}
$\bigl((\nbigt,W),L,N\bigr)$
has a relative monodromy filtration,
if and only if 
$\bigl((\nbigt^{(0)},W^{(0)}),L^{(0)},N^{(0)}\bigr)$
has a relative monodromy filtration.
\end{lem}
\pf
Let $M'$ be the filtration of $(\nbigt,W)$ in $\MTS(X)$,
given by Deligne's inductive formula for $N$.
The induced filtration $M^{\prime\,(0)}$
on $(\nbigt^{(0)},W^{(0)})$ satisfies
Deligne's inductive formula for $N^{(0)}$.
Note 
$\Gr^{L}(M)^{(0)}\simeq
 \Gr^{L^{(0)}}(M^{(0)})$.
Thus, we are done.
\hfill\qed

\vspace{.1in}

Let $\MTS(X)^{\RMF}$ denote the full subcategory
of $\MTS(X)^{\fil,\nil}$,
whose objects have relative monodromy filtrations.

\begin{lem}
The category $\MTS(X)^{\RMF}$ is 
equipped with tensor product
and inner homomorphism 
as in Proposition {\rm\ref{prop;10.7.23.21}}.
\end{lem}
\pf
By using Proposition \ref{prop;10.11.5.1},
it can be reduced to
the case $\nbiga=\Vect_{\cnum}$.
\hfill\qed

\subsection{Some conditions for the existence 
of relative monodromy filtration}

Let $(\nbigt,L)$ be a filtered smooth $\nbigr_X$-triple.
Let $N:\nbigt\lrarr \nbigt\otimes\newTate(-1)$
be a morphism such that
$N\cdot L_k(\nbigt)\subset 
 L_{k}(\nbigt)\otimes\newTate(-1)$.
Because $\Gr^L(\nbigt,L,N)$ is graded,
it has a relative monodromy filtration 
$W\bigl(\Gr^L(\nbigt)\bigr)$
in the category of $\nbigr_X$-triples.
Assume the following:
\begin{itemize}
\item
$(\Gr^L(\nbigt),W)$ is a mixed twistor structure
on $X$.
\end{itemize}
The underlying $\nbigr_X$-modules of $\nbigt$
are denoted by $\nbigm_1$ and $\nbigm_2$.

\begin{prop}
\label{prop;11.2.1.1}
Assume that, for each $P\in X$,
there exists a subset
$U_P\subset \cnum_{\lambda}$
such that (i) $|U_P|=\infty$,
(ii) $(\nbigm_1,L,N)_{|(\lambda,P)}
 \in \Vect_{\cnum}^{\RMF}$
for any $\lambda\in U_P$.
Then, 
$\nbigt$ is equipped
with a filtration $W(\nbigt)$
such that
(i) $(\nbigt,W)$ is a mixed twistor structure,
(ii) $W$ is a relative monodromy filtration of
 $(\nbigt,L,N)$ in $\MTS(X)$.
\end{prop}
\pf
We have only to consider the case
$X$ is a point $\{P\}$.
The set $U_P$ is denoted just by $U$.
Because $\Gr^L(\nbigt)$ is a mixed twistor
structure,
the pairing $C$ of $\nbigm_1$ and $\nbigm_2$
is non-degenerate.
Hence, we can regard $\nbigt$
as a vector bundle $V$ on $\proj^1$,
obtained as the gluing of $\nbigm_1^{\lor}$
and $\sigma^{\ast}\nbigm_2$.
It is equipped with a filtration $L$
and a nilpotent morphism
$N:V\lrarr V\otimes\nbigo_{\proj^1}(2)$
preserving $L$.
The relative monodromy filtration $W$ on $\Gr^L(V)$
gives a mixed twistor structure.

\begin{lem}
\label{lem;11.2.1.2}
If $(V,L,N)_{|(\lambda,P)}
 \in\Vect_{\cnum}^{\RMF}$ 
for any $\lambda\in\proj^1$,
then $V$ is equipped
with a filtration $W(V)$
such that
(i) $(V,W)$ is a mixed twistor structure,
(ii) $W$ is a relative monodromy filtration of
 $(V,L,N)$.
\end{lem}
\pf
Let $M(N_{|\lambda};L_{|\lambda})$ 
denote the relative monodromy filtration
of $(V,L,N)_{|\lambda}$.
For each $k\in\seisuu$,
the rank of $M_k(N_{|\lambda};L_{|\lambda})$
is independent of $\lambda\in \proj^1$.
Then, we obtain that
it depends on $\lambda$ continuously,
from Deligne's inductive formula
(\ref{eq;10.9.27.21})
and (\ref{eq;10.9.27.22}).
The property (i) follows from the canonical decomposition
$\Gr^W\simeq \Gr^W\Gr^L$.
Thus, we obtain Lemma \ref{lem;11.2.1.2}.
\hfill\qed

\vspace{.1in}

According to the lemma,
we have only to prove that
$(V,L,N)_{|\lambda}\in \Vect_{\cnum}^{\RMF}$
for any $\lambda\in\proj^1$.
We use an induction on
the length of the filtration $L$.
We assume that
(i) $V=L_k(V)$,
(ii) the claim holds for $L_{k-1}(V)$,
and we shall prove that
the claim holds for $V$.
We consider the morphisms (\ref{eq;10.9.27.20})
for $(V,L,N)_{|\lambda}$ $(\lambda\in \proj^1)$.
The assumption implies that
(\ref{eq;10.9.27.20}) for $\lambda\in U$ vanishes.
Note that the both hand sides for $\lambda\in \proj^1$
in (\ref{eq;10.9.27.20})
give vector bundles on $\proj^1$.
Hence, by the continuity,
we obtain the vanishing of
(\ref{eq;10.9.27.20})
for any $\lambda\in\proj^1$.
Thus, Proposition \ref{prop;11.2.1.1}
is proved.
\hfill\qed

\begin{cor}
Suppose that $X$ is connected.
Assume that, for a point $P\in X$,
there exists a subset
$U_P\subset \cnum_{\lambda}$
such that (i) $|U_P|=\infty$,
(ii) $(\nbigm_1,L,N)_{|(\lambda,P)}
 \in \Vect_{\cnum}^{\RMF}$
for any $\lambda\in U_P$.
Then, 
$\nbigt$ is equipped
with a filtration $W(\nbigt)$
such that
(i) $(\nbigt,W)$ is a mixed twistor structure,
(ii) $W$ is a relative monodromy filtration of
 $(\nbigt,L,N)$ in $\MTS(X)$.
\hfill\qed
\end{cor}

\vspace{.1in}
The following lemma can be proved
similarly and more easily.
\begin{lem}
\label{lem;11.2.1.3}
Let $\bigl((\nbigt,W),N,L\bigr)\in\MTS(X)^{\fil,\nil}$.
Assume that, for each $P\in X$,
there exists $U_P\in \cnum_{\lambda}$
such that
(i) $|U_P|=\infty$,
(ii) $(\nbigm_1,L,N)_{|(\lambda,P)}
 \in\Vect^{\RMF}_{\cnum}$.
Then, there exists a relative monodromy filtration
$M(N;L)$ in $\MTS(X)$.
\hfill\qed
\end{lem}

\section{Polarizable mixed twistor structure}

\subsection{Statements}

Let $X$ be a complex manifold.
We consider an abelian category
$\nbiga=\MTS(X)$
with additive auto equivalences
$\Sigma^{p,q}(\nbigt)=\nbigt\otimes\nbigu(-p,q)$.
Then, for any finite set $\Lambda$,
we obtain the abelian category
$\MTS(X,\Lambda):=
 \MTS(X)(\Lambda)$
as in \S\ref{subsection;10.11.5.10}.
For an object $(\nbigt,W,\vecN)\in\MTS(X,\Lambda)$,
we set $N(\Lambda):=\sum_{j\in\Lambda}N_j$.
\index{category $\MTS(X,\Lambda)$}
\index{map $N(\Lambda)$}

An object $(\nbigt,W,\vecN)\in\MTS(X,\Lambda)$
is called
a $(w,\Lambda)$-polarizable mixed twistor structure,
if (i) $W=M\bigl(N(\Lambda)\bigr)[w]$,
(ii) there exists a Hermitian sesqui-linear duality 
$\nbigs:\nbigt\lrarr\nbigt^{\ast}\otimes\newTate(-w)$
of weight $w$,
such that
$\nbigs\circ N_i=-N_i^{\ast}\circ\nbigs$
and that
$(N(\Lambda)^{\ast})^{\ell}\circ\nbigs$
induces a polarization of
$P\Gr^W_{\ell}(\nbigt)$.
(Recall Remark \ref{rem;13.8.1.1}.)
Such $\nbigs$ is called a polarization of
$(\nbigt,W,\vecN)$.
Let $\MTSpol(X,w,\Lambda)\subset
 \MTS(X,\Lambda)$ denote
the full subcategory of $(w,\Lambda)$-polarizable
mixed twistor structure on $X$.
\index{category $\MTSpol(X,w,\Lambda)$}
The following proposition is essentially proved
in \cite{mochi2},
based on the results in the Hodge case
in \cite{ck}, \cite{cks1}, \cite{cks2} and \cite{k3}.
We will give an indication
in \S\ref{subsection;11.1.29.20}.

\begin{prop}
\label{prop;10.11.5.5}
The family of the categories $\MTSpol(X,w,\Lambda)$
satisfies the property {\bf P0--3} in
{\rm\S\ref{subsection;10.10.12.21}}.
\end{prop}

We state some complementary property.
First, we give remarks on the ambiguity of
polarizations.
\begin{lem}
\label{lem;10.7.9.2}
Let $(\nbigt,W,\vecN)\in\MTSpol(X,w,\Lambda)$.
If it is simple,
i.e., there is no non-trivial subobject,
then a polarization of $(\nbigt,W,\vecN)$
is unique up to constant multiplication.
\end{lem}
\pf
Let $S_i$ $(i=1,2)$ be polarizations
of $(\nbigt,W,\vecN)$.
They induce an endomorphism of
$(\nbigt,W,\vecN)$
in $\MTSpol(X,w,\Lambda)$.
Because $(\nbigt,W,\vecN)$ is simple,
it is a scalar multiplication,
which implies 
$S_1=\alpha\,S_2$
for some $\alpha\in\cnum$.
Because they are polarizations,
we obtain $\alpha$ is a positive number.
\hfill\qed

\vspace{.1in}

According to Proposition \ref{prop;10.11.5.5},
$(\nbigt,W,\vecN)\in\MTSpol(X,w,\Lambda)$
has a canonical decomposition
\begin{equation}
 \label{eq;10.7.22.2}
 (\nbigt,W,\vecN)
\simeq
 \bigoplus_i 
 (\nbigt_i,W_i,\vecN_i)\otimes U_i,
\end{equation}
where 
(i) $(\nbigt_i,W_i,\vecN_i)\not\simeq
 (\nbigt_j,W_j,\vecN_j)$ for $i\neq j$,
(ii) each $(\nbigt_i,W_i,\vecN_i)$ is irreducible,
(iii) $U_i$ are vector spaces.
(We regard a vector space
as a constant pure twistor structure on $X$.)
We take a polarization $S_i$ of each
$(\nbigt_i,W_i,\vecN_i)$,
which is unique up to positive multiplication.
We argue the following proposition
in \S\ref{subsection;11.2.1.4}.
\begin{prop}
\label{prop;11.1.29.21}
Any polarization of $(\nbigt,W,\vecN)$
is of the form 
$\bigoplus S_i\otimes h_i$,
where $h_i$ are hermitian metrics of $U_i$.
\end{prop}

\subsubsection{Some operations}

Let $(\nbigt^{(i)},W^{(i)},\vecN^{(i)})
 \in\MTSpol(X,w_i,\Lambda)$
$(i=1,2)$.
We have the induced filtration
$\Wtilde$ on 
$\nbigttilde:=
 \nbigt^{(1)}\otimes\nbigt^{(2)}$,
and $(\nbigttilde,\Wtilde)$ is a mixed twistor structure.
We have the induced morphisms
$\Ntilde_j:=N^{(1)}_j\otimes\id
+\id\otimes N^{(2)}_j$
for $j\in\Lambda$.
It is easy to observe that
$(\nbigttilde,\Wtilde,\vecNtilde)
\in\MTSpol(X,w_1+w_2,\Lambda)$.
In this situation,
$(\nbigttilde,\Wtilde,\vecNtilde)$
is denoted by 
$(\nbigt^{(1)},W^{(1)},\vecN^{(1)})
\otimes
 (\nbigt^{(2)},W^{(2)},\vecN^{(2)})$.

Similarly,
for $(\nbigt,W,\vecN)\in\MTSpol(X,w,\Lambda)$,
we have the induced filtration $\Wbar$
on $\nbigtbar:=\nbigt^{\lor}$, and the induced
morphisms $\Nbar_j:=-N^{\lor}_j$.
Then, $(\nbigtbar,\Wbar,\overline{\vecN})
\in\MTSpol(X,-w,\Lambda)$.
It is denoted by $(\nbigt,W,\vecN)^{\lor}$.

\vspace{.1in}
Recall the operations $j^{\ast}$ and $\gammatilde_{sm}^{\ast}$
on $\VTS(X)$.
It naturally induces operations on
$\MTS(X,\Lambda)$,
and it preserves $\MTSpol(X,\Lambda,w)$.

\subsection{Proof of Proposition \ref{prop;10.11.5.5}}
\label{subsection;11.1.29.20}

The property {\bf P1} clearly holds.
Let us prove {\bf P2}.
We have only to consider the case 
$|\Lambda|=1$
and $X=\{P\}$.
Let $(\nbigt_i,W,\vecN)\in\MTSpol(w_i,\Lambda)$
with $w_1>w_2$.
A morphism 
$F:(\nbigt_1,W,N)\lrarr(\nbigt_2,W,N)$
induces
$\TNIL(F):\TNIL(\nbigt_1,N)\lrarr \TNIL(\nbigt_2,N)$
on $\Delta^{\ast}(R)$
for any $R>0$.
(See \S\ref{subsection;10.12.25.20}
for $\TNIL$.)
If $R$ is sufficiently small,
$\TNIL(V_i,N)$ are pure of weight $w_i$
(see \cite{mochi8}),
and hence we have $\TNIL(F)=0$.
It implies $F=0$.
Hence, {\bf P2} holds.

\subsubsection{The property P0}
We prepare some lemmas.

\begin{lem}
Let $(\nbigt_i,W,\vecN)\in\MTS(X,\Lambda)$ $(i=1,2)$.
Then, 
$\bigoplus_{i=1,2}(\nbigt_i,W,\vecN)
 \in\MTSpol(X,w,\Lambda)$
if and only if 
$(\nbigt_i,W,\vecN)\in\MTSpol(X,w,\Lambda)$
$(i=1,2)$.
\end{lem}
\pf
If $S_i$ are polarizations of
$(\nbigt_i,W,\vecN)$,
$S_1\oplus S_2$ is a polarization of
$\bigoplus (\nbigt_i,W,\vecN)$.
If $S$ is a polarization of 
$\bigoplus(\nbigt_i,W,\vecN)$,
we have the induced Hermitian sesqui-linear duality $S_i$
of $\nbigt_i$ $(i=1,2)$,
given by
$\nbigt_i\lrarr \bigoplus_{j=1,2}\nbigt_j
\lrarr\bigoplus\nbigt_j^{\ast}\otimes\newTate(-w)
\lrarr \nbigt_i^{\ast}\otimes\newTate(-w)$.
It is easy to check that
$S_i$ are are polarizations of
$(\nbigt_i,W,\vecN)$.
\hfill\qed

\vspace{.1in}
Let $(\nbigt,W,\vecN,S)$
be a $(w,\Lambda)$-polarized
mixed twistor structure.
Let $(\nbigt',\!W',\!\vecN')\!\subset 
 (\nbigt,W,\vecN)$ be a subobject
in the category $\MTS(X,\Lambda)$.
We assume that the monodromy weight filtration
$M\bigl(N'(\Lambda)\bigr)$ on $\nbigt'$
satisfies $W'=M\bigl(N'(\Lambda)\bigr)[w]$.
Let $\nbigt''$ be the kernel of the composite
of the following morphisms:
\[
\begin{CD}
\nbigt
@>{S}>>
\nbigt^{\ast}\otimes\newTate(-w)
@>>>
 (\nbigt')^{\ast}\otimes\newTate(-w) 
\end{CD}
\]
It induces a subobject
$(\nbigt'',W'',\vecN)\subset (\nbigt,W,\vecN)$
in the category $\MTS(X,\Lambda)$.

\begin{lem}
\label{lem;10.7.22.1}
We have $\nbigt'\cap \nbigt''=0$.
Namely, 
we have a decomposition
which is orthogonal with respect to $S$:
\[
 (\nbigt,W,\vecN)
=(\nbigt',W',\vecN')
\oplus
 (\nbigt'',W'',\vecN'')
\]
In particular,
$(\nbigt',W',\vecN')$
and $(\nbigt'',W'',\vecN'')$
are polarizable.
\end{lem}
\pf
We have only to consider 
the case $|\Lambda|=1$
and $X=\{P\}$.
We have the induced morphism
$\TNIL(\nbigt',N')\subset \TNIL(\nbigt,N)$
on $\Delta^{\ast}(R)$.
If $R$ is sufficiently small,
both $\TNIL(\nbigt,N)$ and $\TNIL(\nbigt',N')$
are pure of weight $w$,
and the pairing $S$ induces a polarization 
of $\TNIL(\nbigt,N)$.
We have the orthogonal decomposition
$\TNIL(\nbigt,N)=\TNIL(\nbigt',N')
 \oplus\TNIL(\nbigt',N')^{\bot}$.
We obtain
$\nbigt'\cap \nbigt''=0$.
\hfill\qed

\begin{lem}
\label{lem;10.7.9.1}
The category $\MTSpol(X,w,\Lambda)$
is abelian and semisimple,
i.e.,
{\bf P0} holds.
\end{lem}
\pf
Let $F:(\nbigt_1,W,\vecN)\lrarr
 (\nbigt_2,W,\vecN)$
be a morphism 
in $\MTSpol(X,w,\Lambda)$.
We have the kernel, the image and the cokernel
in $\MTS(X,\Lambda)$,
denoted by
$(\Ker F,W,\vecN)$,
$(\Image F,W,\vecN)$, and $(\Cok F,W,\vecN)$.
Let us prove that they are
$(w,\Lambda)$-polarizable mixed twistor structure.
It is easy to prove that the filtrations $W$
on $\Ker F$, $\Image F$ and $\Cok F$
are equal to 
$M\bigl(N(\Lambda)\bigr)[w]$.
By Lemma \ref{lem;10.7.22.1},
they are $(w,\Lambda)$-polarizable
mixed twistor structures.
Hence, $\MTSpol(X,w,\Lambda)$
is abelian.

Let $(\nbigt,W,\vecN)\in\MTSpol(X,w,\Lambda)$,
and let $(\nbigt',W',\vecN')\subset (\nbigt,W,\vecN)$
be a subobject in $\MTSpol(X,w,\Lambda)$.
By Lemma \ref{lem;10.7.22.1},
we have a decomposition
$(\nbigt,W,\vecN)=
 (\nbigt',W',\vecN')\oplus(\nbigt'',W'',\vecN'')$
in $\MTSpol(X,w,\Lambda)$.
Hence, $\MTSpol(X,w,\Lambda)$
is semisimple.
\hfill\qed

\subsubsection{Property P3}

Let $(\nbigt,W,\vecN)\in\MTSpol(X,w,\Lambda)$.
For $I\subset\Lambda$
and $\veca\in\real^{I}_{>0}$,
we put $N(\veca):=\sum_{i\in I} a_i\,N_i$.
\begin{lem}
\mbox{{}}\label{lem;10.11.5.3}
\begin{itemize}
\item
The filtrations $M\bigl(N(\veca)\bigr)$
are independent of
$\veca\in\real^I_{>0}$.
It is denoted by $M(I)$.
\item
Let $I\subset J\subset\Lambda$.
For $\veca\in\real_{>0}^J$,
the relative monodromy filtration
$M\bigl(N(\veca);M(I)\bigr)$ exists,
and it is equal to $M(J)$.
\end{itemize}
In particular, {\bf P3.1} holds.
\end{lem}
\pf
We have only to consider the case
that $X$ is a point.
The claims are known in the Hodge case
(\cite{ck}).
The twistor case can be easily reduced 
to the Hodge case.
(See \S3 of \cite{mochi2}).
\hfill\qed

\vspace{.1in}

We take $\bullet\in\Lambda$,
and put $\Lambda_0:=\Lambda\setminus\bullet$.
We put 
$(\nbigt_k^{(1)},W^{(1)}):=
 \Gr^{M(N_{\bullet})}_k(\nbigt,W)$.
Let $\vecN^{(1)}$ denote 
the tuple of morphisms
$N^{(1)}_j:\nbigt_k^{(1)}\lrarr 
 \nbigt_k^{(1)}\otimes\newTate(-1)$
induced by $N_j$ $(j\in\Lambda\setminus\bullet)$.

\begin{lem}
\label{lem;10.7.22.20}
$\bigl(\nbigt^{(1)}_k,W^{(1)},\vecN^{(1)}\bigr)
\in\MTSpol(X,w+k,\Lambda_0)$.
A polarization is naturally induced
by $N_{\bullet}$ and a polarization $\nbigs$ 
of $(\nbigt,W,\vecN)$.
On the primitive part,
it is induced as
$(N_{\bullet}^{\ast})^k\nbigs$.
\end{lem}
\pf
We have only to consider the case that
$X$ is a point.
By considering $\Gr^W$,
we can reduce the issue to the Hodge case,
where the claim is known
by the work due to
Cattani-Kaplan-Schmid and Kashiwara-Kawai.
(See \cite{cks1}, \cite{cks2} and \cite{k3}.)
It is also easy to apply their argument in our case.
\hfill\qed

\vspace{.1in}

Let $(\nbigt^{(2)},W^{(2)},\vecN^{(2)})$
denote the image of
$N_{\bullet}:(\nbigt,W,\vecN)
 \lrarr(\nbigt,W,\vecN)\otimes\newTate(-1)$
in the category $\MTS(X,\Lambda)$.
Let $\nbigs$ be a polarization of $(\nbigt,W,\vecN)$.
It is easy to observe that the composite
$\Image(N_{\bullet})
\lrarr
 \nbigt\otimes\newTate(-1)
\stackrel{\nbigs}{\lrarr}
 \nbigt^{\ast}\otimes\newTate(-w-1)$
factors through
$\Image(N_{\bullet})^{\ast}\otimes
 \newTate(-w-1)$,
by using 
$N_{\bullet}\circ\nbigs=-N_{\bullet}^{\ast}\circ\nbigs$.
Namely,
$\nbigs$ and $N_{\bullet}$
induces a sesqui-linear duality
$\nbigstilde$ of $\nbigt^{(2)}$.

\begin{lem}[Proposition 3.126 \cite{mochi2}]
\label{lem;10.7.22.22}
$(\nbigt^{(2)},W^{(2)},\vecN^{(2)})
\in\MTSpol(X,w+1,\Lambda)$,
and 
$\nbigstilde$ is a polarization of 
$(\nbigt^{(2)},W^{(2)},\vecN^{(2)})$.
\hfill\qed
\end{lem}

We set $\nbigt^{(3)}:=\Sigma^{1,0} \nbigt^{(2)}$.
We have naturally induced morphisms
\[
\begin{CD}
 \Sigma^{1,0}(\nbigt)
 @>{u}>>
\nbigt^{(3)}
 @>{v}>>
 \Sigma^{0,-1}(\nbigt)
\end{CD}
\]

\begin{lem}
\label{lem;10.9.25.20}
Let $\ast \in\Lambda_0$.
We have the decomposition
\[
 \Gr^{M(N_{\ast})}(\nbigt^{(3)})
=\Image \Gr^{M(N_{\ast})}(u)
\oplus
 \Ker\Gr^{M(N_{\ast})}(v).
\]
\end{lem}
\pf
It follows from Proposition 3.134 of \cite{mochi2}.
\hfill\qed

\vspace{.1in}
We obtain {\bf P3.2} 
by an inductive use of Lemma \ref{lem;10.7.22.20}
with Lemma \ref{lem;10.11.5.3}.
We obtain {\bf P3.3}
from Lemma \ref{lem;10.7.22.22}
and Lemma \ref{lem;10.9.25.20}.
Thus, the proof of Proposition \ref{prop;10.11.5.5}
is finished.
\hfill\qed

\subsection{Proof of Proposition \ref{prop;11.1.29.21}}
\label{subsection;11.2.1.4}

The restriction of $S$
to $(\nbigt_i\otimes U_i)\otimes\sigma^{\ast}
 \bigl(\nbigt_j\otimes U_j\bigr)$
induces a morphism
$(\nbigt_i,W_i,\vecN_i)\otimes U_i
\lrarr
 (\nbigt_j,W_j,\vecN_j)\otimes U_j$
in $\MTSpol(w,\Lambda)$.
It has to be $0$ if $i\neq j$.
Hence, (\ref{eq;10.7.22.2}) is
orthogonal with respect to $S$.
By using Lemma \ref{lem;10.7.22.1},
we obtain that
a polarization $S$ of $(\nbigt,W,\vecN)$
is of the form $S=\bigoplus S_i\otimes h_i$.
\hfill\qed

\section{Infinitesimal mixed twistor modules}

\subsection{Definition}
We have the category of filtered objects
in $\MTS(X,\Lambda)$,
denoted by $\MTS(X,\Lambda)^{\fil}$.
We consider a twistor version of
infinitesimal mixed Hodge modules
introduced by Kashiwara \cite{kashiwara-mixed-Hodge}.
\index{infinitesimal mixed twistor module}
\index{pre-infinitesimal mixed twistor module}
\index{$\Lambda$-IMTM}
\index{$\Lambda$-pre-IMTM}

\begin{df}
Let $(\nbigt,W,L,\vecN)\in\MTS(X,\Lambda)^{\fil}$.
\begin{itemize}
\item
It is called a variation of $\Lambda$-pre-infinitesimal 
mixed twistor module on $X$,
or simply
a $\Lambda$-pre-IMTM on $X$,
if $\Gr^L_w(\nbigt,W,\vecN)$ is 
 a $(w,\Lambda)$-polarizable mixed twistor structure
 on $X$.
\item
It is called a variation of
 $\Lambda$-infinitesimal mixed twistor module on $X$,
 or simply $\Lambda$-IMTM on $X$,
 if moreover
 there exists a relative monodromy filtration
 $M(N_j;L)$ for any $j\in\Lambda$.
\hfill\qed
\end{itemize}
\end{df}
If we do not have to distinguish $\Lambda$,
we use ``IMTM'' instead of
``$\Lambda$-IMTM''.
The full subcategory of
$\Lambda$-IMTM
(resp. $\Lambda$-pre-IMTM)
in $\MTS(X,\Lambda)^{\fil}$
is denoted by
$\IMTM(X,\Lambda)$
(resp. $\IMTM^{\pre}(X,\Lambda)$).
\index{category $\IMTM(X,\Lambda)$}
\index{category $\IMTM^{\pre}(X,\Lambda)$}
Note the following lemma,
which follows from Lemma \ref{lem;10.12.28.1}
and Proposition \ref{prop;10.11.5.5}.
\begin{lem}
\label{lem;11.1.29.23}
$\IMTM^{\pre}(X,\Lambda)$ is
an abelian category.
Any morphism
in $\IMTM^{\pre}(X,\Lambda)$
is strict with respect to the filtration $L$.
\hfill\qed
\end{lem}

\begin{rem}
The definitions of IMTM and pre-IMTM
are not given in a parallel way to those 
of infinitesimal mixed Hodge module (IMHM)
and pre-IMHM in {\rm\cite{kashiwara-mixed-Hodge}}.
For pre-IMHM,
the weight filtration $W$ is given only
for $\Gr^L(\nbigt)$.
For IMHM,
the existence of relative monodromy filtration
$M\bigl(N(J);L\bigr)$ 
is assumed for each $J\subset\Lambda$.
But, it was proved in Theorem {\rm 4.4.1}
of {\rm\cite{kashiwara-mixed-Hodge}}
that, for a given pre-IMHM,
if $M\bigl(N_j;L\bigr)$ exists for each $j\in \Lambda$,
then $M\bigl(N(J);L\bigr)$ exists 
for each $J\subset\Lambda$.
\hfill\qed
\end{rem}

The following lemma is a weaker version of
Proposition \ref{prop;10.9.27.31} below.
For simplicity, assume that $X$ is connected.
\begin{lem}
\label{lem;11.2.18.10}
Let $(\nbigt,L,\vecN)\in\VTS(X,\Lambda)^{\fil}$,
i.e.,
$(\nbigt,L,\vecN)$ is a filtered object
in $\VTS(X,\Lambda)$.
Assume the following:
\begin{itemize}
\item
$\Gr_w^L(\nbigt,\vecN)\in \MTSpol(X,\Lambda,w)$.
\item
For a point $P\in X$,
there exists $U_P\subset \cnum_{\lambda}$
such that
(i) $|U_P|=\infty$,
(ii) for any $\lambda\in U_P$,
 $(\nbigm_1,N(\Lambda),L)_{|(\lambda,P)}\in
 \Vect^{\RMF}_{\cnum}$
and 
 $(\nbigm_1,N_i,L)_{|(\lambda,P)}
\in\Vect^{\RMF}_{\cnum}$.
Here, $\nbigm_1$ is one of the underlying
$\nbigr_X$-modules.
\end{itemize}
Then, 
$(\nbigt,W,L,\vecN)\in
 \IMTM(X,\Lambda)$.
\end{lem}
\pf
Applying Proposition \ref{prop;11.2.1.1}
to $(\nbigt,L,N(\Lambda))$,
we obtain the existence of
$W=M\bigl(N(\Lambda);L\bigr)$,
and we have
$(\nbigt,W,\vecN)\in \MTS(X,\Lambda)$.
By applying Lemma \ref{lem;11.2.1.3}
to $\bigl((\nbigt,W),N_i,L\bigr)$ $(i\in\Lambda)$,
we obtain the existence of
$M(N_i;L)$.
Hence, $(\nbigt,W,L,\vecN)\in \IMTM(X,\Lambda)$.
\hfill\qed

\begin{cor}
Let $(\nbigt,W,L,\vecN)\in\MTS(X,\Lambda)^{\fil}$.
Assume the following:
\begin{itemize}
\item
For a point $P\in X$,
there exists $U_P\subset \cnum_{\lambda}$
such that
(i) $|U_P|=\infty$,
(ii) for any $\lambda\in U_P$,
$(\nbigm_1,N_i,L)_{|(\lambda,P)}
\in\Vect^{\RMF}_{\cnum}$.
\end{itemize}
Then, 
$(\nbigt,W,L,\vecN)\in\IMTM(X,\Lambda)$.
\hfill\qed
\end{cor}

\subsection{Statements}

We state some basic property of IMTM.
We will prove the following theorem 
in
\S\ref{subsection;11.1.29.32}--\S\ref{subsection;11.1.29.31}
\begin{thm}
\label{thm;11.1.29.22}
The categories $\IMTM(X,\Lambda)$
have the property {\bf M0--3}
in {\rm\S\ref{subsection;10.10.12.21}}.
\end{thm}

We state some complementary property.
We have naturally defined dual and tensor product
on the category $\MTS(X,\Lambda)^{\fil}$.

\begin{prop}
$\IMTM(X,\Lambda)$
and $\IMTM^{\pre}(X,\Lambda)$
are preserved by 
the dual and the tensor product.
\end{prop}
\pf
The claims for $\Lambda$-pre-IMTM
are clear.
We have only to care
the existence of a relative monodromy filtration.
For the dual, it is clear.
For the tensor product, 
it follows from a result due to Steenbrink-Zucker
(see Proposition \ref{prop;10.7.23.21}).
\hfill\qed

\vspace{.1in}
We have naturally defined functors
$j^{\ast}$,
$\ast$,
and $\gammatilde_{\sm}^{\ast}$
on $\MTS(X,\Lambda)^{\fil}$.
The following proposition is clear.
\begin{prop}
They preserve
$\IMTM(X,\Lambda)$
and $\IMTM^{\pre}(X,\Lambda)$.
\hfill\qed
\end{prop}

\subsection{Canonical filtrations}

We reword the construction and the results
in \S\ref{subsection;10.9.25.2}
and \S\ref{subsection;10.11.5.13}.
We consider
 $(\nbigt,W,L,\vecN_{\Lambda})
\in\IMTM(X,\Lambda)$.
Take an element $\bullet\in\Lambda$,
and put $\Lambda_0:=\Lambda\setminus\bullet$.
Let $\Ltilde:=M(N_{\bullet};L)$.
By considering the morphisms
\[
\begin{CD}
 \Sigma^{1,0}(\nbigt,W,\Ltilde,\vecN_{\Lambda_0})
@>{N_{\bullet}}>>
 \Sigma^{0,-1}
 (\nbigt,W,\Ltilde,\vecN_{\Lambda_0})
@>{\id}>>
 \Sigma^{0,-1}
  (\nbigt,W,\Ltilde,\vecN_{\Lambda_0})
\end{CD}
\]
we obtain the filtration 
$\Nhat_{\bullet\,\ast}L$ of
$\Sigma^{0,-1}(\nbigt,W,\Ltilde,\vecN_{\Lambda_0})$
in the category of $\IMTM(X,\Lambda_0)$.
Similarly, by considering the morphisms
\[
\begin{CD}
 \Sigma^{1,0}(\nbigt,W,\Ltilde,\vecN_{\Lambda_0})
 @>{\id}>>
 \Sigma^{1,0}(\nbigt,W,\Ltilde,\vecN_{\Lambda_0})
@>{N_{\bullet}}>>
 \Sigma^{0,-1}
  (\nbigt,W,\Ltilde,\vecN_{\Lambda_0})
\end{CD}
\]
we obtain the filtration
$\Nhat_{\bullet\,!}L$ of
$\Sigma^{1,0}(\nbigt,W,\Ltilde,\vecN_{\Lambda_0})$
in the category of $\IMTM(X,\Lambda_0)$.
By {\bf M3} of Theorem 
\ref{thm;11.1.29.22},
$\bigl(
 \Sigma^{0,-1}(\nbigt,W),
 \Nhat_{\bullet\ast}L,\vecN_{\Lambda}
 \bigr)$
and 
$\bigl(
 \Sigma^{1,0}(\nbigt,W),
 \Nhat_{\bullet!}L,\vecN_{\Lambda}
 \bigr)$
are $\Lambda$-IMTM.
In particular,
\[
 (\nbigt,W,N_{\bullet\ast}L,\vecN_{\Lambda}):=
 \Sigma^{0,1}\bigl(
 \Sigma^{0,-1}(\nbigt,W),
 \Nhat_{\bullet\ast}L,\vecN_{\Lambda}
 \bigr)
\]
\[
 (\nbigt,W,N_{\bullet !}L,\vecN_{\Lambda}):=
 \Sigma^{-1,0}\bigl(
 \Sigma^{1,0}(\nbigt,W),
 \Nhat_{\bullet!}L,\vecN_{\Lambda}
 \bigr)
\]
are $\Lambda$-IMTM.

\begin{rem}
We can also deduce that they are $\Lambda$-IMTM,
by the reduction to the Hodge case
using $\Gr^{W}$.
The Hodge case was proved
in {\rm\cite{kashiwara-mixed-Hodge}}.
\hfill\qed
\end{rem}

We obtain the following 
as a special case of Proposition 
\ref{prop;10.10.12.25}.
\begin{prop}
Let $(\nbigt,L,\vecN)\in\IMTM(X,\Lambda)$.
Let $i,j\in\Lambda$ be distinct elements.
Let $\star(i),\star(j)\in\{\ast,!\}$.
Then, 
we have 
$N_{i\,\star(i)}(N_{j\,\star(j)}L)
=N_{j\,\star(j)}(N_{i\,\star(i)}L)$.
\hfill\qed
\end{prop}

We use the notation
$\vecN_{J\ast}\vecN_{I!}L$
($I\cap J=\emptyset$)
in the meaning 
as in \S\ref{subsection;10.11.5.13}.
\index{filtration $\vecN_{J\ast}\vecN_{I\bikkuri}L$}

\subsubsection{Canonical prolongations}

We repeat the construction 
in \S\ref{subsection;11.2.18.1}
for this situation.
We have the category
$ML\bigl(\MTS(X,\Lambda)\bigr)$
as in \S\ref{subsection;10.10.13.2}.
Let $(\nbigt,L,\vecN)\in \IMTM(X,\Lambda)$.
For a decomposition
$K_1\sqcup K_2=\Lambda$,
we have an object 
$\nbigt[\ast K_1!K_2]\in
 ML\bigl(
 \MTS(X,\Lambda)\bigr)$
given as follows.
For $I\subset\Lambda$,
we set $I_j:=I\cap K_j$.
Then, we put
$\nbigt[\ast K_1!K_2]_I:=
 \Sigma^{|I_2|,-|I_1|}\nbigt$.
It is equipped with a filtration
$\vecNhat_{I_1!}\vecNhat_{I_2\ast}L$.
\index{filtration $\vecNhat_{I_1\bikkuri}\vecNhat_{I_2\ast}L$}
For $I\sqcup\{i\}\subset\Lambda$,
morphisms $g_{Ii,I}$ and $f_{I,Ii}$
are given as follows:
\[
 g_{Ii,I}:=\left\{
 \begin{array}{ll}
 \id & (i\in \Lambda_2)\\
 N_i & (i\in\Lambda_1)
 \end{array}
 \right.
\quad\quad
 f_{I,Ii}:=\left\{
 \begin{array}{ll}
 N_i & (i\in \Lambda_2)\\
 \id & (i\in\Lambda_1)
 \end{array}
 \right.
\]

\subsection{Property M2.2}
\label{subsection;11.1.29.32}

The claims for {\bf M1} and {\bf M2.1}
are clear by definition.
Let us consider {\bf M2.2}.

\begin{prop}
\label{prop;10.7.22.21}
Let $(\nbigt,W,L,\vecN)\in
 \IMTM(X,\Lambda)$.
For any subset $J\subset\Lambda$,
we put $N(J):=\sum_{j\in J}N_j$.
\begin{itemize}
\item
There exists a relative monodromy filtration
$M\bigl(N(J);L\bigr)$.
We denote it by $M(J;L)$.
\item
Let $I\subset J\subset\Lambda$.
Then, 
$M\bigl(J;L\bigr)
=M\Bigl(N(J);M\bigl(I;L\bigr)\Bigr)$.
\end{itemize}
Note that
$M(I;L)$ is also a relative monodromy filtration of
$\sum_{i\in I}a_iN_i$ $(\veca\in\real_{>0}^I)$.
\end{prop}
\pf
In the Hodge case,
it was proved in \cite{kashiwara-mixed-Hodge}.
The twistor case can be 
reduced to the Hodge case,
by considering $\Gr^W(\nbigt)$
and using Lemma \ref{lem;10.7.22.10}.
\hfill\qed

\vspace{.1in}
We take $\bullet\in\Lambda$,
and put $\Lambda_0:=\Lambda\setminus\bullet$.

\begin{lem}
\label{lem;10.10.8.11}
Let $(\nbigt,W,L,\vecN)$ be a $\Lambda$-pre-IMTM.
Assume that there exists
a relative monodromy filtration
$M=M(N_{\bullet};L)$.
Then,
$\bigl(\nbigt,W,M,\vecN_{\Lambda_0}\bigr)$
is $\Lambda_0$-pre-IMTM,
where $\vecN_{\Lambda_0}=\bigl(
 N_i\,\big|\,i\in\Lambda_0
 \bigr)$.

If $(\nbigt,W,L,\vecN)$ is a $\Lambda$-IMTM,
$(\nbigt,W,M,\vecN_{\Lambda_0})$ is 
a $\Lambda_0$-IMTM.
\end{lem}
\pf
The first claim follows from the canonical splitting
$\Gr^M\simeq\Gr^M\Gr^L$
and Lemma \ref{lem;10.7.22.20}.
The second claim follows from
Proposition \ref{prop;10.7.22.21}.
\hfill\qed

\vspace{.1in}
Let $(\nbigt,W,L,\vecN)\in \IMTM(X,\Lambda)$.
For a decomposition
$\Lambda=\Lambda_0\sqcup\Lambda_1$,
we obtain an object in
$\MTS(X,\Lambda_0)^{\fil}$:
\[
 \res^{\Lambda}_{\Lambda_0}\bigl(
 \nbigt,W,L,\vecN\bigr)
:=
 \bigl(\nbigt,W,M(\Lambda_1;L),\vecN_{\Lambda_0}\bigr)
\]
\index{functor $\res^{\Lambda}_{\Lambda_0}$}
Here, $\vecN_{\Lambda_0}:=\bigl(
 N_j\,\big|\,j\in\Lambda_0
 \bigr)$.
We obtain the following corollary 
by an inductive use of Lemma {\rm\ref{lem;10.10.8.11}}
with Proposition {\rm\ref{prop;10.7.22.21}}.
\begin{cor}
$\res^{\Lambda}_{\Lambda_0}(\nbigt,W,L,\vecN)$
is a $\Lambda_0$-IMTM on $X$.
\hfill\qed
\end{cor}
Thus, we have proved the claim for {\bf M2.2}.
For any $I\subset \Lambda$,
we obtain a functor
$\res^{\Lambda}_{I}:
 \IMTM(X,\Lambda)\lrarr \IMTM(X,I)$.
We naturally have
$\res^{I}_J\circ\res^{\Lambda}_I
=\res^{\Lambda}_J$.

\subsection{Property M0}

We consider {\bf M0} for $\IMTM(X,\Lambda)$.
It is clear that
(i) any injection $\Phi:\Lambda\lrarr\Lambda_1$
induces
$\IMTM(X,\Lambda)\lrarr \IMTM(X,\Lambda_1)$,
(ii) we naturally have
$\MTSpol(X,w,\Lambda)\subset \IMTM(X,\Lambda)$.

Let us prove that $\IMTM(X,\Lambda)$
is abelian.
Let $F:(\nbigt,W,L,\vecN)\lrarr(\nbigt',W',L',\vecN')$
be a morphism in $\IMTM(X,\Lambda)$.
According to Lemma \ref{lem;11.1.29.23},
we have
$(\Ker F,W,L,\vecN)$,
$(\Image F,W,L,\vecN)$,
$(\Cok F,W,L,\vecN)$
in $\IMTM^{\pre}(\Lambda)$.
It remains to prove that
there exist relative monodromy filtrations
$M(N_j;L)$ on them
for any $j\in\Lambda$.
For that purpose,
we have only to prove that
$F$ is strict with respect to $M(N_j;L)$.
Fix $j\in\Lambda$,
and we put $\Lambda_0:=\Lambda\setminus\{j\}$.
Let $\vecNtilde=(N_i\,|\,i\in\Lambda_0)$
and $\vecNtilde':=(N_i'\,|\,i\in\Lambda_0)$.
Because
\[
 \bigl(\nbigt,M(N_{j};L),\vecNtilde\bigr)
\lrarr
 \bigl(\nbigt',M(N'_{j};L'),\vecNtilde'\bigr)
\]
is a morphism in $\IMTM^{\pre}(X,\Lambda_0)$,
we have the desired strictness.
Thus, we proved the claim for {\bf M0}.

\subsection{Property {\bf M3}}
\label{subsection;10.7.22.31}

We state the property {\bf M3}
in this situation.
We will prove it in 
\S\ref{subsection;11.1.29.30}--\S\ref{subsection;11.1.29.31}.
Let $\Lambda$ be a finite set.
Fix an element $\bullet\in \Lambda$,
and we put
$\Lambda_0:=\Lambda\setminus\bullet$.
Let $(\nbigt,W,L,\vecN_{\Lambda})\in
 IMTM(\Lambda)$.
We have the induced object
$(\nbigt,W,\Ltilde,\vecN_{\Lambda_0})
:=
\res_{\Lambda_0}^{\Lambda}
 (\nbigt,W,L,\vecN_{\Lambda})$
in $\IMTM(X,\Lambda_0)$.
We consider $(\nbigt',W,\Ltilde,\vecN'_{\Lambda_0})
\in \IMTM(X,\Lambda_0)$ with 
the following morphisms
\[
\begin{CD}
 \Sigma^{1,0}(\nbigt,W,\Ltilde,\vecN_{\Lambda_0})
@>{u}>>
 (\nbigt',W,\Ltilde,\vecN'_{\Lambda_0})
@>{v}>>
 \Sigma^{0,-1}(\nbigt,W,\Ltilde,\vecN_{\Lambda_0})
\end{CD}
\]
 in $\IMTM(X,\Lambda_0)$,
such that $v\circ u=N_{\bullet}$.
We set $N'_{\bullet}:=u\circ v$,
and the induced tuple
$\vecN'_{\Lambda_0}\sqcup\{N'_{\bullet}\}$
is denoted by $\vecN'_{\Lambda}$.
According to Corollary \ref{cor;10.11.5.2},
we have a filtration $L$ of $\nbigt'$
in $\IMTM(X,\Lambda_0)$,
obtained as the transfer of $L(\nbigt)$ by $(u,v)$.
We will prove the following proposition.

\begin{prop}
\label{prop;10.7.22.30}
$(\nbigt',W,L,\vecN'_{\Lambda})$
is a $\Lambda$-IMTM,
i.e.,
the categories
$\IMTM(X,\Lambda)$ have
the property {\bf M3}.
The relative monodromy filtrations
$M(N'_j;L)$ $(j\in\Lambda_0)$
are obtained as the transfer of
$M(N_j;L)$ by $(u,v)$.
\end{prop}

\subsection{Transfer for pre-IMTM}
\label{subsection;11.1.29.30}

As a preparation,
let us address a similar issue for pre-IMTM.
We consider objects
$(\nbigt,W,L,\vecN_{\Lambda})\in
 \IMTM^{\pre}(\Lambda)$
and $(\nbigt',W,\Ltilde,\vecN'_{\Lambda_0})
 \in \IMTM^{\pre}(X,\Lambda_0)$
with morphisms 
in $\IMTM^{\pre}(X.\Lambda_0)$,
\begin{equation}
\label{eq;11.2.23.1}
\begin{CD}
 \Sigma^{1,0}(\nbigt,W,\Ltilde,\vecN_{\Lambda_0})
 @>{u}>>
 (\nbigt',W,\Ltilde,\vecN'_{\Lambda_0})
 @>{v}>>
 \Sigma^{0,-1}(\nbigt,W,\Ltilde,\vecN_{\Lambda_0})
\end{CD}
\end{equation}
such that $v\circ u=N_{\bullet}$.
We have a unique filtration $L$ of
$(\nbigt',W,\Ltilde,\vecN_{\Lambda_0})$
in $\IMTM^{\pre}(X,\Lambda_0)$
obtained as the transfer of $L(\nbigt)$ by $(u,v)$.
We set $N'_{\bullet}:=u\circ v$,
and the tuple
$\vecN'_{\Lambda_0}\sqcup\{N'_{\bullet}\}$
is denoted by $\vecN'_{\Lambda}$.

\begin{lem}
\label{lem;10.7.24.1}
$(\nbigt',W,L,\vecN'_{\Lambda})
 \in \IMTM^{\pre}(X,\Lambda)$.
\end{lem}
\pf
Because $\Gr_k^L(N'_{\bullet})\!=\!0$
on $\Ker\Gr_k^L(v)$,
the induced filtration $\Ltilde$
on $\Ker\Gr_k^L(v)$ is pure of weight $k$,
i.e.,
$\Gr^{\Ltilde}_j=0$ unless $j=k$.
Hence, 
$(\Ker\Gr_k^L(v),\vecN'_{\Lambda_0})$
is naturally isomorphic to a direct summand of
$\Gr^{\Ltilde}_k\Gr^L_k(\nbigt',\vecN_{\Lambda_0})$.
By the canonical splitting 
in \S\ref{subsection;10.12.28.1},
$\Gr^{\Ltilde}_k\Gr^L_k(\nbigt',\vecN_{\Lambda_0})$
is naturally isomorphic to a direct summand of
$\Gr^{\Ltilde}_k(\nbigt',\vecN_{\Lambda_0})$.
Hence, we have
$(\Ker\Gr_k^L(v),\vecN'_{\Lambda_0})
\in\MTSpol(X,k,\Lambda_0)$.

We regard $N_{\bullet}:\nbigt\lrarr \vecT^{-1}\nbigt$.
We have
$\bigl(
\Image \Gr^L_k(N_{\bullet}),
 \Gr^L_k\vecN_{\Lambda}
 \bigr)\in
\MTSpol(X,k+1,\Lambda)$,
according to
{\bf P3.3} in Proposition \ref{prop;10.11.5.5}.
Because $\Gr^L_k(v)$
induces the following isomorphism
\[
 \bigl(
 \Image\Gr^L_k(u),\,
 \Gr^L_k(\vecN'_{\Lambda})
 \bigr)
\simeq
 \Sigma^{1,0}
 \bigl(\Image \Gr^L_k(N_{\bullet}),
 \Gr^L_k\vecN_{\Lambda}
 \bigr),
\]
we have
$\bigl(
 \Image\Gr^L_k(u),\,
 \Gr^L_k(\vecN'_{\Lambda})
 \bigr)
\in\MTSpol(X,k,\Lambda)$.
Thus, we obtain that
$(\nbigt',W,L,\vecN'_{\Lambda})
\in \IMTM^{\pre}(X,\Lambda)$.
\hfill\qed

\vspace{.1in}
For the proof of Proposition \ref{prop;10.7.22.30},
it remains to prove the existence of
relative monodromy filtrations
$M\bigl(N_j';L(\nbigt')\bigr)$ for $j\in\Lambda_0$.

\subsection{Existence of
relative monodromy filtration
 in a special case}

Let $\nibar:=\{1,2\}$.
Let $(\nbigt',W,L,\vecN')\in \IMTM^{\pre}(X,\nibar)$
and $(\nbigt,W,L,\vecN)\in \IMTM(X,\nibar)$.
Assume that we are given morphisms
in $\MTS(X)^{\fil,\nil}$
\[
\begin{CD}
 (\Sigma^{1,0}(\nbigt,W),L,N_1)
 @>{u}>>
 \bigl(\nbigt',W,L,N_1'\bigr)
 @>{v}>>
 \bigl(\Sigma^{0,-1}(\nbigt,W),L,N_1\bigr)
\end{CD}
\]
such that
(i) $v\circ u=N_2$ and $u\circ v=N_2'$,
(ii) $\bigl(\nbigt,\nbigt';u,v;L\bigr)$
 is filtered $S$-decomposable,
(iii) $(\nbigt',L,N_2')\in\MTS(X)^{\RMF}$.
By Proposition \ref{prop;10.10.3.2},
$u$ and $v$ give
\[
 \Sigma^{1,0}\bigl(\nbigt,W,M(N_2;L)\bigr)
\lrarr
 \bigl(\nbigt',W,M(N'_2;L)\bigr)
\lrarr
 \Sigma^{0,-1}\bigl(\nbigt,W,M(N_2;L)\bigr).
\]
The following lemma 
is based on an argument in \cite{saito2}.
\begin{lem}
\label{lem;10.7.20.11}
Under the assumption,
$(\nbigt',W,L,\vecN')$ is a $\nibar$-IMTM,
namely,
there exists
a relative monodromy filtration
$M(N_1';L)$.
Moreover,
$M(N_1';L)$ is obtained 
as the transfer of $M(N_1;L)$
by $(u,v)$.
\end{lem}
\pf
We put $W^{(1)}:=M(N_1;L)$.
Then, we have $(\nbigt,W^{(1)},N_2)\in \IMTM(X,1)$
and $W=M(N_2;W^{(1)})$.
We have a unique filtration $W^{(1)}$ of $\nbigt'$
obtained as the transfer of $W^{(1)}(\nbigt)$
with respect to $(u,v)$.
It satisfies the conditions
(A1) $W(\nbigt')=M(N_2';W^{(1)})$,
(A2) $\bigl(\nbigt,\nbigt';u,v;W^{(1)}\bigr)$
is filtered $S$-decomposable.
We shall prove that $W^{(1)}(\nbigt')$
gives $M(N_1';L)$.

Let us prove that
$(\nbigt',W^{(1)},N_2')$ is a $1$-IMTM,
i.e.,
$\bigl(\Gr^{W^{(1)}}_w\nbigt',
 N_2'\bigr)
\in\MTSpol(X,w,1)$.
On $\Ker\Gr^{W^{(1)}}_wv$,
the filtration $W$ is pure of weight $w$
by construction.
According to {\bf P3.3} 
in Proposition \ref{prop;10.11.5.5},
$\Image\Gr^{W^{(1)}}_w u$
is a polarizable mixed twistor structure of weight $w$.
Hence, $(\nbigt',W^{(1)},N'_2)$
is $1$-IMTM.

\vspace{.1in}

Put $(\nbigt_{k,m},W):=L_k(\nbigt,W)/L_{m}(\nbigt,W)$
and $(\nbigt'_{k,m},W):=L_k(\nbigt',W)/L_m(\nbigt',W)$
in the category $\MTS(X)$.
They are equipped with the induced 
filtrations $L$
and the induced tuple of morphisms
$\vecN$ and $\vecN'$.
Then, we have
$(\nbigt_{k,m},W,L,\vecN)
 \in\IMTM(X,\nibar)$
and 
$(\nbigt'_{k,m},W,L,\vecN')
 \in\IMTM^{\pre}(X,\nibar)$.
The induced morphisms
$u:\Sigma^{1,0}\nbigt_{k,m}\lrarr\nbigt'_{k,m}$
and $v:\nbigt_{k,m}'\lrarr \Sigma^{0,-1}\nbigt_{k,m}$
satisfy the assumptions (i) and (ii).
Hence, by the above argument,
we obtain a filtration $W^{(1)}(\nbigt_{k,m}')$
with which 
$(\nbigt'_{k,m},W,W^{(1)},N_2')\in \IMTM(X,1)$.

\vspace{.1in}

Let us consider the exact sequence
\[
 0\lrarr L_{k-1}(\nbigt',W)\lrarr
 L_k(\nbigt',W)\lrarr \Gr^{L}_k(\nbigt',W)\lrarr 0.
\]
The arrows are morphisms in $\IMTM(1)$.
In particular, they are strict with respect to $W^{(1)}$.

Note that $M(N_1')[k]$ on $\Gr^L_k(\nbigt')$
satisfies the condition (A1) for $W^{(1)}$.
It also satisfies (A2) according to 
{\bf P3.3} in Proposition \ref{prop;10.11.5.5}.
Hence, the induced filtration
$W^{(1)}$ on $\Gr^{L}_k(\nbigt')$
is the same as $M(N_1')[k]$.
Then, we can conclude that
$W^{(1)}$ gives $M(N_1;L)$ on $\nbigt$.
\hfill\qed

\begin{rem}
For any $2$-IMTM $(\nbigt,W,L,\vecN)$,
we have 
\[
 N_{2\star}M(N_1;L)=
 M(N_1;N_{2\star}L)
\]
for $\star=\ast,!$.
Indeed, we can deduce it in our situation
from Lemma {\rm\ref{lem;10.7.20.11}}.
Alternatively, by considering $\Gr^W$,
we can also reduce it to the Hodge case
proved in {\rm\cite{kashiwara-mixed-Hodge}}.
\hfill\qed
\end{rem}

\subsection{End of the proof of Proposition
\ref{prop;10.7.22.30}}
\label{subsection;11.1.29.31}

Let us return to the situation in
\S\ref{subsection;10.7.22.31}.
According to Lemma \ref{lem;10.7.24.1},
we have only to prove the existence of
$M(N_j';L\nbigt')$ for $j\in\Lambda_0$.
We have only to consider the case that
$X$ is a point.
Put $\Lambda_1:=\Lambda_0\setminus\{j\}$.
We set
$\vecN_{\Lambda_1}:=\bigl(
 N_i\,|\,i\in\Lambda_1\bigr)$
and $\vecN_{\Lambda_1}':=
 \bigl(
 N_i'\,\big|\,i\in\Lambda_1
 \bigr)$.
We obtain smooth $\nbigr$-triples
$\TNIL_{\Lambda_1}(\nbigt,\vecN_{\Lambda_1})$
and 
$\TNIL_{\Lambda_1}(\nbigt',\vecN'_{\Lambda_1})$
on $(\cnum^{\ast})^{\Lambda_1}$
with filtrations $L$.
We take a point
 $\iota_P:\{P\}\lrarr(\cnum^{\ast})^{\Lambda_1}$,
which is sufficiently close to the origin $(0,\ldots,0)$.
We set
$\nbigt_P:=
 \iota_{P}^{\ast}
 \TNIL_{\Lambda_1}(\nbigt,\vecN_{\Lambda_1})$,
which is equipped with the induced filtration $L$
and the morphisms
$N_j,N_{\bullet}:
 \nbigt_P\lrarr \vecT^{-1}\nbigt_P$.
Similarly, we put
$\nbigt'_P:=
 \iota_P^{\ast}
 \TNIL_{\Lambda_1}(\nbigt',\vecN'_{\Lambda_1})$
equipped with the induced filtration $L$
and the morphisms
$N_j',N'_{\bullet}:
 \nbigt'_P\lrarr \vecT^{-1}\nbigt'_P$.
We can apply Lemma \ref{lem;10.7.20.11}
to $(\nbigt_P,L_P,N_j,N_{\bullet})$
and $(\nbigt_P',L_P,N_j',N'_{\bullet})$,
and we obtain a relative monodromy filtration
$M\bigl(N_j';L(\nbigt'_P)\bigr)$.
By construction of $\TNIL$,
we obtain a relative monodromy filtration
$M\bigl(N_j';L(\nbigt')\bigr)$.
Thus, the proof of Proposition \ref{prop;10.7.22.30}
is finished.
\hfill\qed

\section{Nearby cycle functor
along a monomial function}

\subsection{Beilinson IMTM
and its deformation}

Recall the Beilinson triple in 
\S\ref{subsection;11.2.1.10}.
We use the same symbol to denote
the pull back via a morphism $X$ to a point.
It is naturally equipped
with the weight filtration $W$
given by
$W_k(\II^{a,b})=
 \bigoplus_{-2i\leq k}\newTate(i)$,
and $(\II^{a,b},W)$ 
is a mixed twistor structure on $X$.
The tuple
$(\II^{a,b},W,L,N_{\II})$
is a $1$-IMTM on $X$,
where $L=W$.
For any $\vecc\in\real^{\Lambda}$,
let $\vecc\,N_{\II}$ denote
the tuple $\bigl(c_i\,N_{\II}\,\big|\,i\in\Lambda\bigr)$.
We obtain 
$\bigl(
 \II^{a,b},W,L,\vecc\,N_{\II}
 \bigr)\in\IMTM(X,\Lambda)$.
Let $\vecvarphi=(\varphi_i\,|\,i\in\Lambda)$
be a tuple of holomorphic functions on $X$.
We set
$\II^{a,b}_{\vecc,\vecvarphi}:=
 \Def_{\vecvarphi}\bigl(
 \II^{a,b},W,L,\vecc\,N_{\II}
 \bigr)$
in $\IMTM(X,\Lambda)$.

\subsection{Statement}

We will omit to denote the weight filtration
of mixed twistor structure.
We consider $(\nbigt,L,\vecN)\in\IMTM(X,\Lambda)$
and $\vecm\in\seisuu_{>0}^{\Lambda}$.
We obtain the following $\Lambda$-IMTM on $X$:
\[
 \bigl(\Pi_{\vecm,\vecvarphi}^{a,b}\nbigt,L,\vecNtilde\bigr)
:=
 \bigl(\nbigt,L,\vecN\bigr)\otimes
 \II^{a,b}_{\vecm,\vecvarphi} 
\]
For any subset $I\subset\Lambda$,
we put $\Ntilde_I:=\prod_{i\in I}\Ntilde_i$.
For $\star=\ast,!$,
we have the filtrations
$\Ntilde_{I\star}L$ on 
$\Pi_{\vecm,\vecvarphi}^{a,b}\nbigt$.
Let $M$ be a sufficiently large integer.
We have an induced morphism
of $\Lambda$-IMTM:
\[
 \Ntilde_I:
 \Sigma^{|I|,0}
 \bigl(
 \Pi_{\vecm,\vecvarphi}^{0,M}\nbigt,\,
 \Ntilde_{I!}L,\,\vecNtilde
 \bigr)
\lrarr
 \Sigma^{0,-|I|}
 \bigl(\Pi_{\vecm,\vecvarphi}^{0,M}\nbigt,\,
 \Ntilde_{I\ast}L,\,\vecNtilde\bigr)
\]
The cokernel in $\IMTM(X,\Lambda)$
is denoted by 
$\bigl(
 \psi^{(0)}_{\vecm,\vecvarphi}(\nbigt)_I,\Lhat,\vecNtilde
 \bigr)$.

On the other hand,
the filtration $L$ of $(\nbigt,\vecN)$
naively induces a filtration on 
$\Pi^{a,b}_{\vecm,\vecvarphi}\nbigt$
given by
$L_k\bigl(
 \Pi^{a,b}_{\vecm,\vecvarphi}\nbigt
 \bigr)=\Pi^{a,b}_{\vecm,\vecvarphi}L_k(\nbigt)$
in $\IMTM(X,\Lambda)$.
It induces a filtration of 
$\psi^{(0)}_{\vecm,\vecvarphi}(\nbigt)_I$
in $\IMTM(X,\Lambda)$.
They are also denoted by $L$.

The morphism 
$N_{\II}:\II^{a,b}\lrarr\vecT^{-1}\II^{a,b}$
naturally induces a morphism
$\psi^{(0)}_{\vecm,\vecvarphi}(\nbigt)_I
\lrarr
 \vecT^{-1}\psi^{(0)}_{\vecm,\vecvarphi}(\nbigt)_I$,
which is also denoted by $N_{\II}$.
We shall prove the following theorem.
\begin{thm}
\label{thm;10.10.8.40}
We have $\Lhat=M(N_{\II};L)[1]$
on $\psi^{(0)}_{\vecm,\vecvarphi}(\nbigt)_I$.
\end{thm}
Here, for a filtration $F$,
we set $F[a]_{a+j}:=F_{j}$.

The Hodge version of this theorem 
appeared in \cite{saito2}.
We obtain the following as a special case.
\begin{cor}
\label{cor;11.1.29.50}
We have $\Lhat=M(N_{\II};L)[1]$ on
$\psi^{(0)}_{\vecm,\vecvarphi}(\nbigt)_{\Lambda}$.
\hfill\qed
\end{cor}

\subsection{Variant}

Let $K\subset\Lambda$ be any non-empty subset.
For any $\vecm\in\seisuu_{>0}^K$,
we obtain 
$\bigl(\Pi_{\vecm,\vecvarphi}^{a,b}\nbigt,
 L,\vecNtilde\bigr)
:=
 \bigl(\nbigt,L,\vecN\bigr)
 \otimes\II^{a,b}_{\vecm,\vecvarphi}$
in $\IMTM(X,\Lambda)$.
We have the induced morphism 
in $\IMTM(X,\Lambda)$:
\[
 \Ntilde_K:
 \Sigma^{|K|,0}
 \bigl(
 \Pi_{\vecm,\vecvarphi}^{0,M}\nbigt,\,
 \Ntilde_{K!}L,\,\vecNtilde
 \bigr)
\lrarr
 \Sigma^{0,-|K|}
 \bigl(\Pi_{\vecm,\vecvarphi}^{0,M}\nbigt,\,
 \Ntilde_{K\ast}L,\,\vecNtilde\bigr)
\]
On the cokernel 
$\psi^{(0)}_{\vecm,\vecvarphi}(\nbigt)_{K}$,
we have the induced filtration $\Lhat$
and the tuple $\vecNtilde$
so that $\bigl(
 \psi^{(0)}_{\vecm,\vecvarphi}(\nbigt)_K,\Lhat,\vecNtilde
 \bigr)$
is a $\Lambda$-IMTM.
It is also equipped with the naively induced filtration $L$.
We have a naturally induced morphism
$N_{\II}:
 \psi_{\vecm,\vecvarphi}^{(0)}(\nbigt)
 \lrarr
 \psi_{\vecm,\vecvarphi}^{(0)}(\nbigt)
 \otimes\newTate(-1)$
in $\IMTM(X,\Lambda)$.
\begin{cor}
\label{cor;10.8.30.10}
We have $\Lhat=M(N_{\II};L)[1]$.
\end{cor}
\pf
By using $\TNIL_{K^c}$,
we can reduce it to Corollary 
\ref{cor;11.1.29.50}.
\hfill\qed

\subsection{Reformulation}

Let $I':=I\sqcup\{i\}$.
We have the following commutative diagram:
\[
 \begin{CD}
 \Sigma^{|I|+1,0}
 \Pi^{a,b}_{\vecm,\vecvarphi}\nbigt
 @>{=}>>
 \Sigma^{|I'|,0}
 \Pi^{a,b}_{\vecm,\vecvarphi}\nbigt
 @>{\Ntilde_i}>>
 \Sigma^{|I|,-1}
 \Pi^{a,b}_{\vecm,\vecvarphi}\nbigt
 \\
 @V{\Ntilde_I}VV @V{\Ntilde_{I'}}VV 
 @V{\Ntilde_I}VV \\
  \Sigma^{1,-|I|}
 \Pi^{a,b}_{\vecm,\vecvarphi}\nbigt
 @>{\Ntilde_i}>>
 \Sigma^{0,-|I'|}
 \Pi^{a,b}_{\vecm,\vecvarphi}\nbigt
 @>{=}>>
 \Sigma^{0,-|I|-1}
 \Pi^{a,b}_{\vecm,\vecvarphi}\nbigt
 \end{CD}
\]
Hence, we obtain naturally induced morphisms:
\[
 \can:
 \Sigma^{1,0}\psi^{(0)}_{\vecm,\vecvarphi}(\nbigt)_I
 \lrarr
 \psi^{(0)}_{\vecm,\vecvarphi}(\nbigt)_{I'}
\quad\quad
 \var:
 \psi^{(0)}_{\vecm,\vecvarphi}(\nbigt)_{I'}
 \lrarr
 \Sigma^{0,-1}
 \psi^{(0)}_{\vecm,\vecvarphi}(\nbigt)_I
\]
The tuple
$\psi^{(0)}_{\vecm,\vecvarphi}(\nbigt)
 :=\bigl(\psi^{(0)}_{\vecm,\vecvarphi}(\nbigt)_I
 \,\big|\,
 I\subset\ellsitabar
 \bigr)$
with the filtration $\Lhat$ is
an object in
$ML\bigl(
 \MTS(X,\Lambda)\bigr)$.
(See \S\ref{subsection;10.10.13.2}
for the category $ML\bigl(\MTS(X,\Lambda)\bigr)$.)
Indeed,
it is obtained as the cokernel of
$(\Pi_{\vecm,\vecvarphi}^{0,N}\nbigt)_!\lrarr
 (\Pi_{\vecm,\vecvarphi}^{0,N}\nbigt)_{\ast}$
in $ML\bigl(\MTS(X,\Lambda)\bigr)$.

On the other hand,
the filtration $L$ naively induces
a filtration of 
$\psi^{(0)}_{\vecm,\vecvarphi}(\nbigt)$
in $ML\bigl(\MTS(X,\Lambda)\bigr)$,
which is also denoted by $L$.
The previous theorem can be reformulated
as follows.
\begin{thm}
\label{thm;10.8.30.4}
We have $\Lhat=M(N;L)[1]$
on $\psi^{(0)}_{\vecm,\vecvarphi}(\nbigt)$.
\end{thm}

\subsection{Proof}

We have only to consider the case that
$L$ is pure of weight $w$,
i.e.,
$(\nbigt,W,\vecN)\in\MTSpol(X,w,\Lambda)$.
Moreover, we may assume $w=0$.
We have the weight filtration $M(N_{\II})$
on each $\psi^{(0)}_{\vecm,\vecvarphi}(\nbigt)_I$,
which is preserved by $\can$ and $\var$.
Hence, we obtain a filtration $M(N_{\II})$
on $\psi^{(0)}_{\vecm,\vecvarphi}(\nbigt)$.
We have only to prove that $M(N_{\II})[1]=\Lhat$.
We may assume that $X$ is a point.
We have only to consider the case
that $\varphi_i=0$ $(i\in\Lambda)$.

According to Corollary 3.132 of \cite{mochi2},
$\bigl(\psi^{(0)}_{\vecm}(\nbigt)_I,
 (\vecN,N_{\II})\bigr)$
is a $(1,\Lambda\sqcup\{\bullet\})$-polarizable 
mixed twistor structure.
Hence, $\bigl(\psi^{(0)}_{\vecm}(\nbigt)_I,
 M(N_{\II})[1],\vecNtilde\bigr)$
is a $\Lambda$-IMTM.
For $I'=I\sqcup\{i\}$,
according to Lemma \ref{lem;10.9.25.20},
the following is $S$-decomposable:
\[
 \Sigma^{1,0}
\Gr^{M(N_{\II})}\psi^{(0)}_{\vecm}(\nbigt)_I
\stackrel{u}{\lrarr}
 \Gr^{M(N_{\II})}\psi^{(0)}_{\vecm}(\nbigt)_{I'}
\stackrel{v}{\lrarr}
 \Sigma^{0,-1}
  \Gr^{M(N_{\II})}\psi^{(0)}_{\vecm}(\nbigt)_{I}
\]
Namely,
$\bigl(\psi^{(0)}_{\vecm}(\nbigt),
 M(N_{\II})[1]\bigr)$ is 
an object in 
$ML\bigl(\MTS(\Lambda)\bigr)$.

\vspace{.1in}

Let $I\subset\Lambda$.\!
We have the filtrations
$M\bigl(\Ntilde(I);M(N_{\II})\![1]\!\bigr)$
and $M(\Ntilde(I);\Lhat)$
on $\psi^{(0)}_{\vecm}\!(\!\nbigt\!)_I$.
To prove $M(N_{\II})[1]=\Lhat$,
we have only to prove
$M\bigl(\Ntilde(I);M(N_{\II})[1]\bigr)
=M(\Ntilde(I);\Lhat)$
on $\psi^{(0)}_{\vecm}(\nbigt)_I$
for any $I\subset\Lambda$,
according to Theorem \ref{thm;10.10.12.12}.
If $I=\Lambda$, 
both the filtrations
$M\bigl(\Ntilde(\Lambda);M(N_{\II})[1]\bigr)$
and $M(\Ntilde(\Lambda);\Lhat)$
are the weight filtration of 
the mixed twistor structure.
Let us consider the case
$J:=\Lambda\setminus I\neq\emptyset$.
We take $Q\in\cnum^J$,
which is sufficiently close to the origin.
We have a $(0,I)$-polarizable 
mixed twistor structure
$\TNIL_{J,Q}(\nbigt)$
with $\vecN_I:=(N_i\in I)$.
Let $\vecm_I:=(m_i\,|\,i\in I)$.
On $\psi^{(0)}_{\vecm_{I}}
 \bigl(\TNIL_{J,Q}(\nbigt)\bigr)_I$,
we have
$M\bigl(\Ntilde(I);W(N)[1] \bigr)
=M\big(\Ntilde(I);\Lhat)$.
The underlying $\nbigr$-modules of
$\psi^{(0)}_{\vecm}(\nbigt)_I$
and 
$\psi^{(0)}_{\vecm_{I}}
 \bigl(
 \TNIL_{J,Q}(\nbigt)_I
 \bigr)$
are naturally isomorphic.
Therefore, we have
\[
 M\bigl(\Ntilde(I);M(N_{\II})[1]\bigr)
=M(\Ntilde(I);\Lhat)
\]
on $\psi^{(0)}_{\vecm}(\nbigt)_I$,
and the proof of 
Theorem \ref{thm;10.8.30.4}
is finished.
\hfill\qed

\section{Twistor version of a theorem of Kashiwara}

Let $(\nbigt,L,\vecN)\in\VTS(X,\Lambda)^{\fil}$.
Let $\nbigm_i$ $(i=1,2)$ denote the underlying
$\nbigr_X$-modules.
We assume the following:
\begin{itemize}
\item
 For each $w$,
 we have
 $\Gr^L_w(\nbigt,L,\vecN)
 \in \MTSpol(X,w,\Lambda)$.
\item
 For each $P\in X$,
 there exists a subset $U_P\subset \cnum_{\lambda}$
 such that
 (i) $|U_P|=\infty$,
 (ii) for any $\lambda\in U_P$ and for $i\in\Lambda$,
 we have
 $(\nbigm_1,L,N_i)_{|P}\in\Vect_{\cnum}^{\RMF}$.
\end{itemize}
The following proposition is
a twistor version of Theorem {\rm 4.4.1}
of \cite{kashiwara-mixed-Hodge}.
\begin{prop}
\label{prop;10.9.27.31}
$(\nbigt,L,\vecN)$ is a $\Lambda$-IMTM,
i.e.,
there exists
a relative monodromy filtration of
$\bigl(\nbigt,L,N(\Lambda)\bigr)$.
\end{prop}

\subsection{A purity theorem (Special case)}

Let $(\nbigt,W,L,\vecN)$ be a $\nibar$-IMTM.
For $i=1,2$,
we have the  morphisms
\[
 N_i:(\Sigma^{1,0}\nbigt,L)\lrarr 
 \Sigma^{0,-1}\bigl(\nbigt,N_{i\ast}L\bigr).
\]
They give morphisms
$N_i:\Sigma^{1,0}L_k(\nbigt)
\lrarr
 \Sigma^{0,-1}\bigl(
 N_{i\ast}L_{k-1}\nbigt\bigr)$,
i.e.,
$N_i:L_k(\nbigt)
\lrarr
 \vecT^{-1}N_{i\ast}L_{k-1}(\nbigt)$.
Similarly,
we have
$N_i:N_{j\ast}L_k(\nbigt)
\lrarr
 \vecT^{-1}N_{j\ast}N_{i\ast}L_{k-1}(\nbigt)$.
Then, we obtain the following
filtered complex $\Pi$ in $\MTS(X)$:
\[
\begin{CD}
 L_{-1}\nbigt
 @>{N_{1}\oplus N_2}>>
 \vecT^{-1}(N_{1\ast}L)_{-2}\nbigt
\oplus
 \vecT^{-1}
 (N_{2\ast}L)_{-2}\nbigt
 @>{N_2-N_1}>>
 \vecT^{-2}
 (N_{1\ast}N_{2\ast}L)_{-3}\nbigt
\end{CD}
\]
The cohomology group
$H^i(\Pi)$ with the induced weight filtration $W$
is a mixed twistor structure.

\begin{lem}
\label{lem;10.9.27.32}
$\Gr^W_jH^i(\Pi)=0$
unless $j\leq i-1$.
\end{lem}
\pf
We have
$\Gr^W_jH^i(\Pi)
\simeq
 H^i\Gr^W_j\Pi$.
Then, the claim can be reduced to the Hodge case
in \cite{kashiwara-mixed-Hodge}.
\hfill\qed

\begin{rem}
The purity theorem can be proved
in a more general situation
as in the Hodge case.
\hfill\qed
\end{rem}

\subsection{Proof of Proposition
\ref{prop;10.9.27.31}}

We have only to consider the case
that $X$  is a point $\{P\}$.
By using $\TNIL$,
we can reduce the issue
to the case $\Lambda=\{1,2\}$.
We regard $(\nbigt,L)$ as a filtered
vector bundle on $\proj^1$,
denoted by $(V,L)$.
We may assume $V=L_0(V)$.
By an inductive argument,
we may assume that the claim holds for
$(L_{-1}V,L,\vecN)$.
According to Lemma \ref{lem;11.2.18.10},
we have only to prove the existence of
relative monodromy filtration
for $N_{1|\lambda}+N_{2|\lambda}$ on
$(V,L)_{|\lambda}$ for any $\lambda\in U_P$.
Then, we have only to apply
the argument in Section 6 of \cite{kashiwara-mixed-Hodge},
with Lemma \ref{lem;10.9.27.32} above.
\hfill\qed

\section{Integrable case}

\subsection{Integrable mixed twistor structure}

\index{integrable mixed twistor structure}

Let $X$ be a complex manifold.
Let $(\nbigt,W)$ be a filtered object
in the category of  integrable smooth $\nbigr_X$-triples.
It is called an integrable variation of mixed twistor structure on $X$,
(or simply an integrable mixed twistor structure on $X$),
if each $\Gr^W_w(\nbigt)$ is integrable polarizable pure twistor
$D$-module of weight $w$.
It is called pure of weight $w$,
if $\Gr^L_m=0$ unless $m=w$.
The category of integrable mixed twistor structure on $X$
is denoted by $\MTSint(X)$.
It is an abelian category.
It is equipped with tensor product
and inner homomorphism.
It is also equipped with additive auto equivalences
$\Sigma^{p,q}$ given by the tensor product with
$\nbigu(-p,q)$.

\subsection{Integrable polarizable mixed twistor structure}

\index{integrable polarizable mixed twistor structure} 

We consider an abelian category
$\nbiga=\MTSint(X)$
with additive auto equivalences
$\Sigma^{p,q}(\nbigt)=\nbigt\otimes\nbigu(-p,q)$.
Then, for any finite set $\Lambda$,
we obtain the abelian category
$\MTSint(X,\Lambda):=
 \MTSint(X)(\Lambda)$
as in \S\ref{subsection;10.11.5.10}.
An object $(\nbigt,W,\vecN)\in\MTSint(X,\Lambda)$
is called
an integrable $(w,\Lambda)$-polarizable mixed twistor structure,
if (i) $W=M\bigl(N(\Lambda)\bigr)[w]$,
(ii) there exists an integrable Hermitian sesqui-linear duality
$\nbigs:\nbigt\lrarr\nbigt^{\ast}\otimes\newTate(-w)$
of weight $w$, which gives a polarization
of the underlying $(\nbigt,W,\vecN)\in\MTS(X,\Lambda)$.
The full subcategory of integrable $(w,\Lambda)$-polarizable
mixed twistor structure
is denoted by $\MTSpol^{\integral}(X,w,\Lambda)$.
The following is an analogue of Proposition \ref{prop;10.11.5.1}.
\begin{prop}
The categories $\MTSpol^{\integral}(X,w,\Lambda)$
satisfy the property {\bf P0--3} in
{\rm\S\ref{subsection;10.10.12.21}}.
\end{prop}
\pf
As for {\bf P0},
we have only to repeat the argument in the ordinary case.
The other property in the integrable case
can be reduced to those in the ordinary case.
\hfill\qed

\vspace{.1in}
As in the case of $\MTSpol(X,\Lambda,w)$,
any object in $\MTSpol^{\integral}(X,\Lambda,w)$
has the canonical decomposition,
and its polarization is unique up to automorphisms.
The category $\MTSpol^{\integral}(X,\Lambda,w)$
is equipped with tensor product, dual,
$j^{\ast}$ and $\gammatilde_{\sm}^{\ast}$.

\subsection{Infinitesimal mixed twistor module}

\index{infinitesimal mixed twistor module}

We consider the category of filtered objects
in $\MTSint(X,\Lambda)$,
denoted by $\MTSint(X,\Lambda)^{\fil}$.
Let $(\nbigt,W,L,\vecN)\in\MTSint(X,\Lambda)^{\fil}$.
\begin{itemize}
\item
It is called 
an integrable $\Lambda$-pre-IMTM on $X$,
if $\Gr^L_w(\nbigt,W,\vecN)$ is 
an integrable $(w,\Lambda)$-polarizable mixed twistor structure
on $X$.
\item
It is called an integrable $\Lambda$-IMTM on $X$,
if moreover there exists a relative monodromy filtration
$M(N_j;L)$ for any $j\in\Lambda$.
In other words,
$(\nbigt,W,L,\vecN)$ is an integrable $\Lambda$-IMTM,
if (i) it is an integrable pre-$\Lambda$-IMTM,
(ii) it is a $\Lambda$-IMTM.
\end{itemize}
Let 
$\IMTMint(X,\Lambda)$
(resp. $\IMTM^{\integral\pre}(X,\Lambda)$)
denote the full subcategory of integrable $\Lambda$-IMTM
(resp. integrable $\Lambda$-pre-IMTM)
in $\MTSint(X,\Lambda)^{\fil}$.
The following lemma is obvious.
\begin{lem}
$\IMTM^{\integral\pre}(\Lambda)$ is
an abelian category.
Any morphism
in $\IMTM^{\integral\pre}(\Lambda)$
is strict with respect to the filtration $L$.
\hfill\qed
\end{lem}

\begin{prop}
The categories $\IMTM(X,\Lambda)$
have the property {\bf M0--3}
in {\rm\S\ref{subsection;10.10.12.21}}.
\end{prop}
\pf
We give only an indication.
Let us consider {\bf M0}.
It is clear that
(i) any injection $\Phi:\Lambda\lrarr\Lambda_1$
induces
$\IMTMint(X,\Lambda)\lrarr \IMTMint(X,\Lambda_1)$,
(ii) we naturally have
$\MTSpol^{\integral}(X,w,\Lambda)\subset \IMTMint(X,\Lambda)$.
Let us prove that $\IMTMint(X,\Lambda)$ is abelian.
Let $F:(\nbigt,W,L,\vecN)\lrarr(\nbigt',W',L',\vecN')$
be a morphism in $\IMTMint(X,\Lambda)$.
We have the objects $(\Ker F,W,L,\vecN)$,
$(\Image F,W,L,\vecN)$,
$(\Cok F,W,L,\vecN)$
in $\IMTM(X,\Lambda)$.
They are naturally integrable smooth 
filtered $\nbigr_X$-triples,
and objects in $\IMTM^{\integral\pre}(X,\Lambda)$.
Hence, they are naturally objects 
in $\IMTMint(X,\Lambda)$.
The claims for {\bf M1} and {\bf M2.1}
are clear by definition.
Let us consider {\bf M2.2}.
Let $(\nbigt,W,L,\vecN)\in\IMTMint(X,\Lambda)$.
By {\bf M2.2} for $\IMTM(X,\Lambda)$,
we have the filtration $M(\Lambda_1;L)$
and the object
$\res^{\Lambda}_{\Lambda_0}(\nbigt,W,L,\vecN)$
in $\IMTM(X,\Lambda_0)$.
We can prove that it is naturally
an object in $\IMTMint(X,\Lambda_0)$,
by using Deligne's formula for relative monodromy filtration
and Kashiwara's canonical decomposition,
and we obtain {\bf M2.2}.
To argue {\bf M3},
let us consider the situation in 
\S\ref{subsection;10.7.22.31}
with integrability.
Let $(\nbigt,W,L,\vecN_{\Lambda})\in
 \IMTMint(\Lambda)$.
We have the induced object
$\res_{\Lambda_0}^{\Lambda}
 (\nbigt,W,L,\vecN_{\Lambda})
=:(\nbigt,W,\Ltilde,\vecN_{\Lambda_0})$
in $\IMTMint(X,\Lambda_0)$.
Let $(\nbigt',W,\Ltilde,\vecN'_{\Lambda_0})
\in \IMTMint(X,\Lambda_0)$ with 
morphisms as in (\ref{eq;11.2.23.1})  in $\IMTMint(X,\Lambda_0)$,
such that $v\circ u=N_{\bullet}$.
We set $N'_{\bullet}:=u\circ v$,
and the induced tuple
$\vecN'_{\Lambda_0}\sqcup\{N'_{\bullet}\}$
is denoted by $\vecN'_{\Lambda}$.
We have a filtration $L$ of $\nbigt'$
in $\IMTMint(X,\Lambda_0)$,
obtained as the transfer of $L(\nbigt)$ by $(u,v)$.
By {\bf M3} for $\IMTM(X,\Lambda)$,
$(\nbigt',W,L,\vecN'_{\Lambda})$
is a $\Lambda$-IMTM.
By repeating the argument for Lemma \ref{lem;10.7.24.1},
we obtain that it is an integrable $\Lambda$-pre-IMTM,
and hence it is an integrable $\Lambda$-IMTM.
\hfill\qed

\begin{rem}
Although we do not give the statements of the integral version of 
Theorem {\rm\ref{thm;10.10.8.40}}
and Proposition {\rm\ref{prop;10.9.27.31}},
they can easily be reduced to the ordinary case.
\hfill\qed
\end{rem}

\chapter{Admissible mixed twistor structure and variants}
\label{section;11.4.3.4}

In this section, we study 
admissible variation of mixed twistor structure,
which is a rather straightforward generalization of
admissible variation of mixed Hodge structure,
studied in \cite{kashiwara-mixed-Hodge}
and \cite{steenbrink-zucker}.

\section{Admissible mixed twistor structure}
\label{subsection;11.2.18.11}

\subsection{Mixed twistor structure on $(X,D)$}

Let $X$ be a complex manifold,
and $D$ be a simply 
normal crossing hypersurface of $X$
with the irreducible decomposition
$D=\bigcup_{i\in\Lambda}D_i$.
For $I\subset\Lambda$,
we set $D_I:=\bigcap_{i\in I}D_i$,
$\del D_I:=\bigcup_{j\not\in I}(D_I\cap D_j)$
and $D_I^{\circ}:=D_I\setminus \del D_I$.
Let $\vecnbigi$ be a good system of ramified irregular
values on $(X,D)$ (see \S\ref{section;13.5.6.100}).
The tuple $(X,D,\vecnbigi)$ is denoted by $\vecX$.

\vspace{.1in}

An object
$(\nbigt,L)\in\VTS(X,D)^{\fil}$
is called
a variation of mixed twistor structure on $(X,D)$,
if its restriction to $X\setminus D$
is a mixed twistor structure on $X\setminus D$.
The full subcategory
$\MTS(X,D)\subset \VTS(X,D)^{\fil}$
is abelian.
An object $(\nbigt,L)\in\MTS(X,D)$
is called a variation of mixed twistor structure on $\vecX$,
if $\nbigt$ is $\vecnbigi$-good.
Let $\MTS(\vecX)\subset\VTS(X,D)^{\fil}$ 
denote the full subcategory 
of variation of mixed twistor structure on $\vecX$.
It is an abelian subcategory.
We shall often omit ``variation of''.

\subsection{Pre-admissibility}

\index{pre-admissible}
\index{$\vecnbigi$-pre-admissible}

A mixed twistor structure $(\nbigt,L)$ on $\vecX$
is called pre-admissible,
if the following holds:
\begin{description}
\item[\bf (Adm0)]
 For each $w$, $\Gr^L_w(\nbigt)$
 is the canonical prolongment of
 an $\vecnbigi$-good wild polarizable variation of
 pure twistor structure of weight $w$.
(See \S11.1 of \cite{mochi7}.)
\item[\bf (Adm1)]
 $\nbigt$ is a smooth $\vecnbigi$-good-KMS
 $\nbigr_{X(\ast D)}$-triple,
 and the filtration $L$ is compatible with
 the KMS-structure.
\end{description}
We shall impose additional conditions
for admissibility.

We give a remark on a specialization
in the case
$X=\Delta^n$ and 
$D=\bigcup_{i=1}^{\ell}\{z_i=0\}$.
Suppose that $(\nbigt,L)\in\MTS(\vecX)$ is pre-admissible.
For $\vecu\in(\real\times\cnum)^{\ell}$,
we obtain
a smooth $\nbigr_{D_{\ellsitabar}}$-triple
$\lefttop{\ellsitabar}\psitilde_{\vecu}(\nbigt)$
which is equipped with a tuple of 
morphisms $\vecnbign=(\nbign_i\,|\,i=1,\ldots,\ell)$.
(See \S\ref{subsection;13.5.4.21}
for $\lefttop{\ellsitabar}\psitilde_{\vecu}(\nbigt)$.)
It is also equipped with the filtration $L$
in the category of $\VTS(D_{\ellsitabar})$
by {\bf (Adm1)},
which is preserved by $\nbign_i$.
Because 
$\Gr^L_w\lefttop{\ellsitabar}\psitilde_{\vecu}(\nbigt)
\simeq
 \lefttop{\ellsitabar}\psitilde_{\vecu}\Gr^L_w(\nbigy)$
are variation of twistor structure,
according to \cite{mochi7},
we have
$\lefttop{\ellsitabar}\psitilde_{\vecu}(\nbigt)
\in \VTS(D_{\ellsitabar})$.
Thus, we obtain
$\bigl(
 \lefttop{\ellsitabar}\psitilde_{\vecu}(\nbigt,L),\vecnbign\bigr)$
in $\VTS(D_{\ellsitabar},\ellsitabar)^{\fil}$.
For any point $\iota_P:\{P\}\lrarr D_{\ellsitabar}$,
we set
$\bigl(
 \lefttop{\ellsitabar}\psitilde_{\vecu}(\nbigt,L),\vecnbign
\bigr)_{|P}:=
 \iota_{P}^{\ast}
 \bigl(\lefttop{\ellsitabar}\psitilde_{\vecu}(\nbigt,L),\vecnbign\bigr)$.
If $\nbigt$ is unramified,
we also have the induced object
$\bigl(
 \lefttop{\ellsitabar}
 \psitilde_{\vecu,\gminia}(\nbigt,L),\vecnbign\bigr)$
in 
$\VTS(D_{\ellsitabar},\ellsitabar)^{\fil}$
for $\vecu\in(\real\times\cnum)^{\ell}$
and $\gminia\in\Irr(\nbigt,\ellsitabar)$.
We set
$\bigl(
 \lefttop{\ellsitabar}\psitilde_{\vecu,\gminia}(\nbigt,L),\vecnbign\bigr)_{|P}:=
 \iota_P^{\ast}\bigl(
 \lefttop{\ellsitabar}
 \psitilde_{\vecu,\gminia}(\nbigt,L),\vecnbign\bigr)$.

\subsection{Admissibility in the smooth divisor case}
\label{subsection;10.9.29.1}

\index{admissible}

Let us consider the case that $D$ is smooth. 
Let $\bigl(\nbigt,L\bigr)\in\MTS(\vecX)$
be pre-admissible.
Let $\nbigm_i$ $(i=1,2)$ be the underlying
$\nbigr_{X(\ast D)}$-modules of $\nbigt$.
We have the induced bundles
$\nbigg_u\nbigm_i$ 
$(u\in\real\times\cnum)$ 
on $\nbigd$
by taking Gr with respect to the KMS-structure.
It is equipped with the endomorphism
$\Res(\DD)$.
Let $N$ denote the nilpotent part.
It is also equipped with the induced filtration $L$
by {\bf (Adm1)}.
Then, $\nbigt$ is called admissible,
if moreover the following holds:
\begin{description}
\item[\bf (Adm2)]
 $N$ on $\bigl(\nbigg_u(\nbigm_i),L\bigr)$
 has a relative monodromy filtration.
\end{description}
(See \S\ref{subsection;13.5.4.20}
for $\nbigg_u(\nbigm_i)$.)
The following lemma is a special case of
Proposition \ref{prop;10.9.29.30} below.
\begin{lem}
\label{lem;11.2.19.1}
Let $(\nbigt,L)\in\MTS(\vecX)$ be admissible.
Assume 
$X=\Delta^n$ and $D=\{z_1=0\}$.
For each $u\in\real\times\cnum$,
$(\psitilde_u(\nbigt),L,\nbign)$ is a $1$-IMTM on $D$.
If $\nbigt$ is unramified,
for each $u\in\real\times\cnum$
and $\gminia\in\Irr(\DD)$,
$\bigl(\psitilde_{\gminia,u}(\nbigt),L,\nbign\bigr)$
is a $1$-IMTM on $D$.
\end{lem}
\pf
We have only to consider the unramified case.
By the compatibility of 
the filtration $L$ and the KMS-structure,
$\bigl(
 \Gr^L_w\psitilde_{\gminia,u}(\nbigt),\nbign
 \bigr)$
comes from $\Gr^L_w(\nbigt)$.
Hence, it is polarizable.
\hfill\qed

\begin{lem}
\label{lem;10.9.29.20}
We can replace the condition
{\bf (Adm2)}
with the following condition:
\begin{description}
\item[\bf (Adm2')]
There exists $U\subset \cnum_{\lambda}$
with $|U|=\infty$
such that
$\bigl(\nbigg_{u}(\nbigm_1),
 L,\nbign\bigr)_{(\lambda,P)}
 \in\Vect_{\cnum}^{\RMF}$
 for any $(\lambda,P)\in U\times D$.
\end{description}
\end{lem}
\pf
It follows from Proposition \ref{prop;11.2.1.1}.
\hfill\qed

\subsection{Admissibility 
in the normal crossing case}

\index{admissible}

Let $X$ be a complex manifold,
and let $D$ be a simple normal crossing 
hypersurface of $X$.
\begin{df}
Let $(\nbigt,L)\in\MTS(\vecX)$ be pre-admissible.
It is called admissible, if
for any smooth point $P\in D$,
there exists a small neighbourhood
$X_P$ of $P$ such that
$\nbigt_{|X_P}$ is admissible
in the sense of the smooth divisor case
({\S\ref{subsection;10.9.29.1}}).
\hfill\qed
\end{df}
The following lemma is clear.
\begin{lem}
\label{lem;11.2.23.21}
If $\nbigt\in\VTS(X,D)$ comes from
a good wild polarizable variation of pure twistor structure
of weight $w$,
it is naturally an admissible mixed twistor structure on $(X,D)$.
\hfill\qed
\end{lem}

\begin{prop}
\label{prop;10.9.29.30}
Let $(\nbigt,L)\in\MTS(\vecX)$ be admissible.
Let $P\in D_I^{\circ}$.
We take a coordinate around $P$.
\begin{itemize}
\item
For $\vecu\in(\real\times\cnum)^I$,
the induced 
$\bigl(
 \lefttop{I}\psitilde_{\vecu}(\nbigt,L),\vecN_I
 \bigr)_{|P}$
is an $I$-IMTM.
\item
If $\nbigt$ is unramified around $P$,
for any $\gminia\in\Irr(\nbigt,P)$
and for any $\vecu\in(\real\times\cnum)^I$,
$\bigl(
 \lefttop{I}\psitilde_{\gminia,\vecu}(\nbigt,L),\vecN
 \bigr)_{|P}$
is an $I$-IMTM.
\end{itemize}
\end{prop}
\pf
We have only to consider the unramified case.
On $\nbigd_i^{(\lambda)}$,
we have the bundle
$\lefttop{i}\nbigg_u\nbigq_{\veca}^{(\lambda)}(\nbigm)$
with the flat $\lambda$-connection 
$\lefttop{i}\DD$.
It is equipped with the $\lefttop{i}\DD$-flat
endomorphism $N_i$
and the $\lefttop{i}\DD$-flat filtration $L$.
By taking the specialization to
$\{\lambda\}\times D_i$,
we obtain a vector bundle
$\lefttop{i}\nbigg_u\nbigq_{\veca}^{\lambda}(\nbigm^{\lambda})$
with a $\lambda$-flat connection $\lefttop{i}\DD^{\lambda}$.
Its formal completion at $P$ is denoted by
$(\Vhat_P,\widehat{\DD}^{\lambda}_P)$.
It is equipped with the induced endomorphism $\widehat{N}_i$
and filtration $\widehat{L}$.
By the assumption,
$(\Vhat_P,\Lhat,\Nhat_i)$
has a relative monodromy filtration.
If $\lambda$ is generic,
we obtain that its specialization at $P$
also has a relative monodromy filtration.

By the consideration in the previous paragraph,
for any generic $\lambda$,
$N_{i|\lambda}$ $(i\in I)$ have
relative monodromy filtrations.
Then, the claim follows from 
Proposition \ref{prop;10.9.27.31}.
\hfill\qed

\vspace{.1in}

Let $f:X'\lrarr X$ be a morphism of complex manifolds,
such that $D':=f^{-1}(D)$ is normal crossing.
Let $f^{-1}\vecnbigi$ denote a good system of
ramified irregular values on $(X',D')$
obtained as the pull back of $\vecnbigi$.
Let $\vecX':=(X',D',f^{-1}\vecnbigi)$.

\begin{cor}
For any admissible variation of mixed twistor structure
$(\nbigt,L)$ on $\vecX$,
the pull back
$f^{\ast}(\nbigt,L)$ is an admissible
variation of mixed twistor structure
on $\vecX'$.
\end{cor}
\pf
It follows from Proposition \ref{prop;10.9.29.30}
and Proposition \ref{prop;10.7.22.21}.
\hfill\qed

\subsection{Category of admissible MTS}

\index{category $\MTS^{\adm}(\vecX)$}

Let $\MTS^{\adm}(\vecX)\subset \MTS(\vecX)$
be the full subcategory of 
admissible mixed twistor structure 
on $\vecX$.
It is equipped with additive auto equivalences 
$\Sigma^{p,q}$
given by the tensor product with $\nbigu(-p,q)$.

\begin{prop}
\label{prop;11.2.23.20}
$\MTS^{\adm}(\vecX)$
is an abelian subcategory of
$\MTS(\vecX)$.
\end{prop}
\pf
Let $F:(\nbigt_1,L)\lrarr(\nbigt_2,L)$
be any morphism in $\MTS^{\adm}(\vecX)$.
Note that $F$ is strict with respect to $L$.
We obtain the filtered $\nbigr_{X(\ast D)}$-triples
$(\Ker F,L)$, $(\Image F,L)$ and $(\Cok F,L)$.
Let us prove that they are objects
in $\MTS^{\adm}(\vecX)$.
If $L$ is pure of weight $w$,
it follows from the theory of
polarizable wild pure twistor $D$-modules.
Indeed,
we have the corresponding
polarizable pure twistor $D$-modules 
$\gbigt_i$ of weight $w$
with morphisms
$\Ftilde:\gbigt_1\lrarr\gbigt_2$.
We have the polarizable pure twistor $D$-modules
$\Ker\Ftilde$,
$\Image\Ftilde$
and $\Cok\Ftilde$,
which are smooth on $X\setminus D$.
By the semisimplicity of the category of
polarizable pure twistor $D$-modules,
$\Ker \Ftilde$ and $\Cok\Ftilde$
are direct summands of
$\gbigt_1$ and $\gbigt_2$,
respectively,
and $\Image\Ftilde$ is a direct summand of
both of $\gbigt_i$ $(i=1,2)$.
Because 
$\Ker F=\Ker\Ftilde(\ast D)$,
$\Cok F=\Cok\Ftilde(\ast D)$,
and
$\Image F=\Image\Ftilde(\ast D)$,
they satisfy {\bf(Adm0)} and {\bf(Adm1)}.
The condition {\bf(Adm2)} is trivial.

Let us consider the mixed case.
Because $F$ is strict with respect to $L$,
we have 
$\Gr^L\Ker F\simeq
 \Ker\Gr^LF$,
$\Gr^L\Image F\simeq
 \Image\Gr^LF$
and $\Gr^L\Cok F\simeq\Cok\Gr^LF$.
Hence, {\bf (Adm0)} is satisfied
for $\Ker F$, $\Image F$ and $\Cok F$.
Let us check {\bf (Adm1--2)}.
By {\bf (Adm0)},
we have already known that
$\Ker F$, $\Image F$ and $\Cok F$
are smooth good $\nbigr_{X(\ast D)}$-triples.

Let us check that 
they are $\vecnbigi$-good $\nbigr_{X(\ast D)}$-triples.
We have only to consider the unramified case.
Let $\nbigm_{i,c}$ $(c=1,2)$
be the underlying $\nbigr_{X(\ast D)}$-modules
of $\nbigt_i$.
We have the irregular decomposition
\[
 \nbigm_{i,c|\widehat{(\lambda,P)}}
=\bigoplus_{\gminia\in\Irr(\nbigt_i)}
 \nbigmhat_{i,c,\gminia}.
\]
Let $F_{1}:\nbigm_{2,1}\lrarr\nbigm_{1,1}$
and 
$F_{2}:\nbigm_{1,2}\lrarr\nbigm_{2,2}$
be the underlying morphisms of $F$.
It is easy to see that $F_{c|\nbigdhat}$
is compatible with the decompositions.
We have natural isomorphisms
$\Ker (F_c)_{|\widehat{(\lambda,P)}}
=\Ker(F_{c|\widehat{(\lambda,P)}})$,
$\Image (F_c)_{|\widehat{(\lambda,P)}}
=\Image(F_{c|\widehat{(\lambda,P)}})$
and $\Cok(F_c)_{|\widehat{(\lambda,P)}}
=\Cok(F_{c|\widehat{(\lambda,P)}})$.
Hence, $\Ker (F_c)$,
$\Image (F_c)$
and $\Cok(F_c)$
are good smooth $\nbigr_{X(\ast D)}$-modules,
and the set of irregular values
are contained in $\nbigi_P$.

To check the remaining claims,
we may and will assume
that $D$ is smooth.
(See Corollary \ref{cor;10.11.9.10}.)
Moreover, we may assume that
$X=\Delta^n$
and $D=\{z_1=0\}$.

Let us consider the regular singular case.
For any $u\in\real\times\cnum$,
\[
 \psitilde_u(F):\psitilde_u(\nbigt_1,L)
\lrarr
 \psitilde_u(\nbigt_2,L)
\]
is a morphism in $\MTS(D)$.
Hence, the cokernel is strict with respect to $\lambda$.
Then, we obtain that 
$F$ is strict with respect to
the KMS-structure,
and $\Ker F$, $\Image F$
and $\Cok F$ are equipped with
the induced KMS-structure.
We obtain the existence of relative monodromy filtration,
because the category of $1$-IMTM is abelian.
Thus, we are done in the regular singular case.

Let us consider the good irregular case.
By the reduction with respect to the Stokes structure
in \S\ref{subsection;11.2.20.3},
we obtain 
$\Gr^{\St}(F):
 \Gr^{\St}(\nbigt_1)\lrarr
 \Gr^{\St}(\nbigt_2)$.
By using the result in the regular singular case,
we obtain that the cokernel of
$\psitilde_{\gminia,u}(F):
 \psitilde_{\gminia,u}(\nbigt_1)\lrarr
 \psitilde_{\gminia,u}(\nbigt_2)$
are strict.
Hence, we obtain that
$F$ is strict with respect to the KMS structure.
We also obtain the existence of
the relative monodromy filtration
of the nilpotent part of the residues.
\hfill\qed

\begin{prop}
\label{prop;10.9.29.2}
\label{prop;10.9.29.40}
Let $F:\nbigt_1\lrarr\nbigt_2$
be a morphism 
in $\MTS^{\adm}(\vecX)$.
Then, 
$F$ is strictly compatible with the $KMS$-structure,
i.e.,
for $\nbigt_i=(\nbigm_{i,1},\nbigm_{i,2},C_i)$,
we have
\begin{equation}
 \label{eq;13.5.6.40}
 F(\nbigqzero_{\veca}\nbigm_{1,2})
=\nbigqzero_{\veca}(\nbigm_{2,2})\cap\Image F,
\end{equation}
\begin{equation}
 \label{eq;13.5.6.41}
 F(\nbigqzero_{\veca}\nbigm_{2,1})
=\nbigqzero_{\veca}(\nbigm_{1,1})\cap\Image F,
\end{equation}
for any $\lambda_0\in\cnum$
and $\veca\in\real^{\Lambda}$,
where $\Lambda$ is the set of irreducible components
of $D$.
\end{prop}
\pf
If $D$ is smooth, the claim has already been proved 
in the proof of Proposition \ref{prop;11.2.23.20}.
If $L$ is pure, the claim follows from that
$\Image F$ is a direct summand of
both of $\nbigt_i$,
as remarked in the proof of Proposition \ref{prop;11.2.23.20}.

We consider the general case.
By Proposition \ref{prop;11.2.23.20},
and by considering the hermitian adjoint,
we have only to consider the case
that $F$ is an epimorphism.
We obtain (\ref{eq;13.5.6.41}) 
from the smooth case and the Hartogs property.
We obtain (\ref{eq;13.5.6.40})
from the pure case by using an easy induction
on the length of $L$.
\hfill\qed

\vspace{.1in}
Let $X=\Delta^n$
and $D=\bigcup_{i=1}^{\ell}\{z_i=0\}$.
As noted in Proposition \ref{prop;10.9.29.30},
we have a functor 
\[
 \lefttop{\ellsitabar}\psitilde_{\vecu}:
 \MTS^{\adm}(\vecX)
 \lrarr \IMTM(D_{\ellsitabar},\ellsitabar). 
\]
We obtain the following corollary
from Proposition \ref{prop;10.9.29.2}.
\begin{cor}
The functor $\lefttop{\ellsitabar}\psitilde_{\vecu}$
is exact.
\hfill\qed
\end{cor}

Let $\nbign:(\nbigt,W)\lrarr(\nbigt,W)\otimes\newTate(-1)$
be a morphism in $\MTS^{\adm}(\vecX)$.
The monodromy weight filtration of $\nbign$
is a filtration in the category
$\MTS^{\adm}(\vecX)$.
It induces a filtration
of $\lefttop{\ellsitabar}\psitilde_{\vecu}(\nbigt,W)$
in $\IMTM(D_{\ellsitabar},\ellsitabar)$,
denoted by $\lefttop{\ellsitabar}\psitilde_{\vecu}M(\nbign)$.
We also have the induced morphism
$\lefttop{\ellsitabar}\psitilde_{\vecu}(\nbign):
 \lefttop{\ellsitabar}\psitilde_{\vecu}(\nbigt)
\lrarr
 \lefttop{\ellsitabar}\psitilde_{\vecu}(\nbigt)\otimes\newTate(-1)$.
Because 
$\lefttop{\ellsitabar}\psitilde_{\vecu}$ is exact,
we obtain the following corollary
by Deligne's inductive formula
for monodromy weight filtration.
\begin{cor}
\label{cor;11.2.19.11}
We have
$\lefttop{\ellsitabar}\psitilde_{\vecu}M(\nbign)=
 M\bigl(
 \lefttop{\ellsitabar}\psitilde_{\vecu}(\nbign)
\bigr)$.
\hfill\qed
\end{cor}

\subsection{Curve test}

Let $(\nbigt,L)\in\MTS(\vecX)$
satisfying {\bf (Adm0)}.
\begin{prop}
Suppose the following condition:
\begin{itemize}
\item
 Let $C$ be any curve contained in $X$
 which transversally intersects with the smooth part of $D$.
 Then, 
 $(\nbigt,L)_{|C}$ is an admissible
 mixed twistor structure on $(C,C\cap D)$.
\end{itemize}
Then, $(\nbigt,L)$ is admissible.
\end{prop}
\pf
By Proposition \ref{prop;13.5.14.11},
$(\nbigt,L)$ satisfies {\bf(Adm1)}.
It is easy to check {\bf(Adm2)}.
\hfill\qed

\section{Admissible polarizable mixed twistor structure}

\subsection{Definition}
\label{subsection;10.11.10.2}

For any finite set $\Lambda$,
we set
$\MTS^{\adm}(\vecX,\Lambda)
 \!:=\!
 \MTS^{\adm}(\vecX)(\Lambda)$.
\index{category $\MTS^{\adm}(\vecX,\Lambda)$}

\subsubsection{Unramified case}
Let $(\nbigt,W,\vecN_{\Lambda})$ be an object 
in $\MTS^{\adm}(\vecX,\Lambda)$,
which is unramified.
It is called $(w,\Lambda)$-polarizable, 
if the following holds.
\begin{description}
\item[\bf (R1)]
There exists a hermitian adjoint
$S:(\nbigt,W,\vecN_{\Lambda})
\lrarr
 (\nbigt,W,\vecN_{\Lambda})^{\ast}\otimes
 \newTate(-w)$ of weight $w$,
such that
$(\nbigt,W,\vecN_{\Lambda},S)_{|X\setminus D}
\in\MTSpol(X\setminus D,w,\Lambda)$.
\item[\bf (R2)]
For any $P\in D_I^{\circ}$,
$\gminia\in\Irr(\nbigt)_{|P}$
and $\vecu\in(\real\times\cnum)^I$,
$\bigl(
 \lefttop{I}\psitilde_{\gminia,\vecu}(\nbigt),
 \vecN_{I\sqcup\Lambda}\bigr)_{|P}$
with $\lefttop{I}\psitilde_{\gminia,\vecu}(S)$
is a $(w,\Lambda\sqcup I)$-polarized
mixed twistor structure.
\end{description}

\begin{lem}
The second condition is independent
of the choice of $S$
in the first condition.
\end{lem}
\pf
Let $S_i$ $(i=1,2)$
be polarizations of $(\nbigt,W,\vecN)$
as in the first condition.
Assume that the second condition is satisfied
for $S_1$.
Let us prove that it is also satisfied for $S_2$.
There exists a decomposition
$(\nbigt,\vecN)_{|X\setminus D}
=\bigoplus (\nbigt_j,\vecN)$
which are orthogonal with respect to
both $S_i$,
and $S_{1|\nbigt_j}=a_j\cdot S_{2|\nbigt_j}$
for some $a_j>0$.
We have the endomorphism $F$ of $\nbigt$
induced by $S$ $(i=1,2)$.
We have $\nbigt_j=\Ker(F_{|X\setminus D}-a_j)$.
Let $P\in D_I^{\circ}$.
We have the induced endomorphisms
$\lefttop{I}\psitilde_{\gminia,\vecu}(F)$
on $\lefttop{I}\psitilde_{\gminia,\vecu}(\nbigt)_{|P}$.
It is a morphism
in the category of polarizable mixed twistor
structures.
Hence, we obtain that
$F-\alpha$ are strict for any $\alpha\in\cnum$.
Then, the decomposition 
$\nbigt_{|X\setminus D}=\bigoplus \nbigt_{j}$
is extended to 
$\nbigt=\bigoplus\nbigt_j'$,
compatible with the KMS-structures.
Then, the claim is clear.
\hfill\qed

\subsubsection{Ramified case}

Let $(\nbigt,W,\vecN_{\Lambda})$
be an object in $\MTS^{\adm}(\vecX,\Lambda)$,
which is not necessarily unramified.
It is called $(w,\Lambda)$-polarizable,
if the following holds:
\begin{itemize}
\item
It is locally the descent of
an unramified admissible
$(w,\Lambda)$-polarizable mixed twistor structure.
\item
There exists a morphism
$S:(\nbigt,W,\vecN_{\Lambda})
\lrarr(\nbigt,W,\vecN_{\Lambda})^{\ast}
 \otimes\newTate(-w)$
such that
$(\nbigt,W,\vecN_{\Lambda},S)_{|X\setminus D}$
is a polarized mixed twistor structure on $X\setminus D$.
\end{itemize}

\subsection{Category of admissible
$(w,\Lambda)$-polarizable
mixed twistor structure}

Let $\nbigp(\vecX,\Lambda,w)\subset
 \MTS^{\adm}(\vecX,\Lambda)$ denote 
the full subcategory of admissible 
$(w,\Lambda)$-polarizable mixed twistor structure
on $\vecX$.

\index{category $\nbigp(X,D,\Lambda,w)$}

\begin{prop}
\label{prop;10.11.9.11}
The family of the categories 
$\nbigp(\vecX,\Lambda,w)$
has the property {\bf P0--3}.
\end{prop}
\pf
The claims for {\bf P1--2} are clear.
Let us prove {\bf P0}.
Let $F:(\nbigt_1,W,\vecN_{\Lambda})
\lrarr(\nbigt_2,W,\vecN_{\Lambda})$
be a morphism in 
$\nbigp(\vecX,\Lambda,w)$.
We have the objects
$(\Ker F,W,\vecN_{\Lambda})$,
$(\Image F,W,\vecN_{\Lambda})$
and $(\Cok F,W,\vecN_{\Lambda})$
in $\MTS^{\adm}(\vecX,\Lambda)$.
We choose polarizations $S_i$
for $(\nbigt_i,W,\vecN_{\Lambda})$.
We have the adjoint
$F^{\lor}:=S_1^{-1}\circ F\circ S_2
 :(\nbigt_2,W,\vecN_{\Lambda})
\lrarr(\nbigt_1,W,\vecN_{\Lambda})$,
giving splittings
$\nbigt_1=\Ker F\oplus\Image F^{\lor}$
and 
$\nbigt_2=\Image F\oplus \Ker F^{\lor}$
in $\MTS^{\adm}(\vecX,\Lambda)$.
The decompositions are compatible with
the filtrations $W$ and 
the tuple of morphisms $\vecN_{\Lambda}$.
Hence, we obtain that
$\nbigp(\vecX,\Lambda,w)$
is abelian and semisimple.

Let us prove {\bf P3.1}.
For $\Lambda_1\subset\Lambda$,
we have the filtration
$M(\Lambda_1):=
 M\bigl(N(\Lambda_1)\bigr)$ of $(\nbigt,W)$
in the category $\MTS^{\adm}(\vecX)$.
Let $\Lambda_2\subset\Lambda_1$.
On $X\setminus D$,
$M(\Lambda_1)_{|X\setminus D}$
is a relative monodromy filtration
of $N(\Lambda_1)_{|X\setminus D}$
with respect to $M(\Lambda_2)_{|X\setminus D}$.
Hence, it is easy to observe that
$M(\Lambda_1)$
is a relative monodromy filtration
of $N(\Lambda_1)$ with respect to
$M(\Lambda_2)$.

Let us consider {\bf P3.2}.
Let $\Lambda=\Lambda_0\sqcup\Lambda_1$
be a decomposition.
We put $\Ltilde:=M(\Lambda_1)$,
which is a filtration of
$(\nbigt,W,\vecN_{\Lambda_0})$
in $\MTS^{\adm}(\vecX,\Lambda_0)$.
\begin{lem}
\label{lem;10.10.10.12}
$\Gr^{\Ltilde}_k(\nbigt,W,\vecN_{\Lambda_0})$
is an object in 
$\nbigp(\vecX,\Lambda_0,w+k)$.
\end{lem}
\pf
We have only to consider the primitive part.
Its restriction to $X\setminus D$
is an object in 
$\MTSpol(X\setminus D,\Lambda_0,w+k)$.
Let $S$ be a polarization of
$(\nbigt,W,\vecN_{\Lambda})$.
We obtain a polarization $S_{k}$
of the primitive part
$P\Gr^{\Ltilde}_k(\nbigt,W,\vecN_{\Lambda_0})$
induced by
$S$ and $N(\Lambda_1)$.
If $\nbigt$ is unramified,
we have natural isomorphisms
$\lefttop{I}\psitilde_{\gminia,\vecu}
 P\Gr^{\Ltilde}_k(\nbigt)_{|P}
\simeq
 P\Gr^{\Ltilde}_k
 \lefttop{I}\psitilde_{\gminia,\vecu}(\nbigt)_{|P}$
for each $P$.
It is equipped with a polarization
induced by
$S$ and $N(\Lambda_1)$
(Lemma \ref{lem;10.7.22.20}),
which is also induced by $S_k$.
Hence, we obtain that
$P\Gr^{\Ltilde}_k(\nbigt,L,
 \vecN_{\Lambda_0})$
is an object in 
$\nbigp(\vecX,\Lambda_0,w+k)$.
\hfill\qed

\vspace{.1in}
Let us prove {\bf P3.3}.
Let $(\nbigt,W,\vecN)$ be an object
in $\nbigp(\vecX,\Lambda,w)$.
Let $\bullet\in\Lambda$.
We have the induced object
$(\Image N_{\bullet},W,\vecN_{\Lambda})$
in $\MTS^{\adm}(\vecX,\Lambda)$.

\begin{lem}
It is an object
in $\nbigp(\vecX,\Lambda,w+1)$.
\end{lem}
\pf
Its restriction to $X\setminus D$ is an object in
$\MTSpol(X\setminus D,\Lambda,w+1)$.
A polarization $S$ for $(\nbigt,W,\vecN_{\Lambda})$
and $N_{\bullet}$ induces
a polarization for
$(\Image N_{\bullet},W,\vecN_{\Lambda})$.
If $\nbigt$ is unramified,
we have
\[
 \lefttop{I}
 \psitilde_{\gminia,u}(\Image N_{\bullet})_{|P}
\simeq
 \Image\Bigl(
\lefttop{I}
 \psitilde_{\gminia,u}(N_{\bullet})_{|P}
 \Bigr).
\]
Hence, the claim follows from
Lemma \ref{lem;10.7.22.22}.
\hfill\qed

\vspace{.1in}

Let $\ast\in\Lambda\setminus\bullet$.
We have the following induced morphisms
in $\MTS^{\adm}(\vecX,\Lambda\setminus\ast)$:
\begin{equation}
 \label{eq;10.11.9.10}
 \Sigma^{1,0}
 \Gr^{W(N_{\ast})}\nbigt
\lrarr
 \Sigma^{1,0}\Gr^{W(N_{\ast})}\Image N_{\bullet}
\lrarr
 \Sigma^{0,-1}
 \Gr^{W(N_{\ast})}\nbigt
\end{equation}
Its restriction to $X\setminus D$ is $S$-decomposable
by Lemma \ref{lem;10.9.25.20}.
Then, (\ref{eq;10.11.9.10})
is $S$-decomposable.
Thus, the proof of Proposition 
\ref{prop;10.11.9.11} is finished.
\hfill\qed

\vspace{.1in}
We have a complement,
which is easy to see.
\begin{lem}
\mbox{{}}\label{lem;11.2.23.23}
\begin{itemize}
\item
Let $(\nbigt,W,\vecN_{\Lambda})
 \in\nbigp(\vecX,\Lambda,w)$.
Then, its dual 
$(\nbigt,W,\vecN_{\Lambda})^{\lor}$
is an object in
 $\nbigp((X,D,-\vecnbigi),\Lambda,-w)$.
\item
Let $(\nbigt_1,W,\vecN_{\Lambda})$
be an object in $\nbigp(\vecX,\Lambda,w)$.
Let $(\nbigt_2,W,\vecN_{\Lambda})$
be an admissible $(w,\Lambda)$-polarizable
mixed twistor structure on $(X,D,\veczero)$.
Then, 
$(\nbigt_1,W,\vecN_{\Lambda})
\otimes
 (\nbigt_2,W,\vecN_{\Lambda})$
is an object in 
$\nbigp(\vecX,\Lambda,w)$.
\item
The functors $j^{\ast}$,
$\ast$
and 
$\gammatilde_{\sm}^{\ast}$
on $\MTS^{\adm}(\vecX,\Lambda)$
preserve 
$\nbigp(\vecX,w,\Lambda)$.
\hfill\qed
\end{itemize}
\end{lem}
(See \S\ref{section;13.5.6.100}
for the trivial good system of ramified irregular values 
$\veczero$.)

\subsection{An equivalent condition}

Let $X=\Delta^n$, $D_i=\{z_i=0\}$
and $D=\bigcup_{i=1}^{\ell}D_i$.
Let $(\nbigt,\vecN)\in\VTS(X,D)(\Lambda)$.
Let $W:=M\bigl(N(\Lambda)\bigr)[w]$.
We assume the following:
\begin{itemize}
\item
$\nbigt$ is 
an unramifiedly $\vecnbigi$-good-KMS 
smooth $\nbigr_{X(\ast D)}$-triple.
\item
{\bf Adm0} holds for $(\nbigt,W)$.
\item
$(\nbigt,\vecN)$ satisfies
{\bf R1--2} with $w$.
\end{itemize}
Note that $W$ is not assumed to 
be compatible with the KMS-structure.
The following proposition implies that
the compatibility is automatically satisfied
under the above assumption.

\begin{prop}
\label{prop;11.2.19.10}
We have 
$(\nbigt,W,\vecN)\in\nbigp(\vecX,w,\Lambda)$.
\end{prop}
\pf
For $\veca\in\cnum^{\Lambda}$,
we put $N(\veca):=\sum a_i\,N_i$,
which gives an endomorphism
of $\nbigqzero_{\vecb}\nbigm_i$.
\begin{lem}
The conjugacy classes of
$N(\veca)_{|(\lambda,P)}$ are independent of 
$(\lambda,P)\in \cnum_{\lambda}\times X$.
\end{lem}
\pf
If we fix a point $P\in X\setminus D$,
they are independent of $\lambda$,
because of the mixed twistor property.
If we fix $\lambda\neq 0$,
they are independent of $P\in X\setminus D$
because of the flatness.
Hence, we obtain that
they are independent 
of $(\lambda,P)\in\cnum_{\lambda}\times (X\setminus D)$.

Let us fix $P\in D_I^{\circ}$.
We use the notation in \S\ref{subsection;13.5.4.20}.
The conjugacy classes of
$N(\veca)_{|(\lambda,P)}$
on
$\lefttop{I}\nbiggzero
 \nbigqzero_{\vecb}(\nbigm_i)_{|(\lambda,P)}$
are independent of $\lambda$,
because of {\bf(R2)} and the mixed twistor property.
If $\lambda$ is generic,
the filtrations $\lefttop{i}\Fzero$
have the splitting around $\lambda$
given by the generalized eigen decomposition
of $\Res_i(\DD)$.
Then, it is easy to observe that 
the conjugacy classes of $N(\veca)_{|(\lambda,P)}$
are independent around $\lambda_0$.

Let us consider the regular singular case.
If we fix a generic $\lambda$,
they are independent of $P\in X$.
Hence, we are done.

Let us consider the unramifiedly good irregular case.
We have
$\Gr^{\St}(\nbigq_{\vecb}\nbigm_i)$
with $\Gr^{St}(N(\veca))$.
By using the previous consideration,
for generic $\lambda$,
we obtain that the conjugacy classes of
$\Gr^{\St}(N(\veca))_{|(\lambda,P)}$
are independent of $P\in X$.
By considering the completion,
we obtain that
the conjugacy classes of
$\Gr^{\St}(N(\veca))_{|(\lambda,P)}$
and $N(\veca)_{|(\lambda,P)}$
are the same.
Thus, we are done.
\hfill\qed

\vspace{.1in}
We also obtain the following lemma.
\begin{lem}
\label{lem;10.10.10.1}
The weight filtration $M(\veca)$
of $N(\veca)$ is a filtration by subbundles
of $\nbigqzero_{\vecb}\nbigm_i$.
Moreover,
$M(\veca)$ is compatible
with the $KMS$-structure.
\end{lem}
\pf
The first claim follows from
the previous lemma.
To prove the compatibility,
according to Proposition \ref{prop;13.5.14.11},
we have only to consider 
the case $n=\ell=1$.
On $\nbigqzero_{b}(\nbigm_i)_{|\nbigdzero}$,
we have the parabolic filtration $\Fzero$.
By using Lemma 5.2 of \cite{mochi2},
we obtain that the weight filtration
$M(\veca)$ of $N(\veca)$ on 
$\nbigqzero_{b}(\nbigm_i)_{|\nbigdzero}$
induces the weight filtration of
$M(\veca)$ of $N(\veca)$
on $\Gr^{\Fzero}(\nbigqzero_{\vecb}\nbigm_i)$.
In particular,
$\Gr^{M(\veca)}\Gr^{\Fzero}
 (\nbigqzero_b\nbigm)$
is a locally free $\nbigo_{\nbigdzero}$-module.
Then, we obtain that
$\Gr^{M(\veca)}(\nbigqzero_{\ast}\nbigm_i)$
is a $KMS$-structure of
$\Gr^{M(\veca)}\nbigm_i$
around $\lambda_0$,
i.e.,
$M(\veca)$ is compatible with
the $KMS$-structure.
\hfill\qed

\vspace{.1in}
In particular,
$W=M\bigl(N(\Lambda)\bigr)[w]$ 
is compatible with the KMS-structure.
The existence of relative monodromy filtration
follows from the property of
polarizable mixed twistor structures.
Thus, Proposition \ref{prop;11.2.19.10} is proved.
\hfill\qed

\subsection{Specialization}
\label{subsection;11.2.23.30}

Let $X=\Delta^n$, $D_i=\{z_i=0\}$
and $D=\bigcup_{i=1}^{\ell}D_i$.
Let $(\nbigt,W,\vecN_{\Lambda})\in\nbigp(\vecX,\Lambda,w)$.
We obtain
$\lefttop{I}\psitilde_{\vecu}(\nbigt)\in \VTS(D_I,\del D_I)$
with the induced morphisms $\vecN_{\Lambda}$.
We also have the naturally induced morphisms $\vecN_I$.
We set $\vecN_{\Lambda\sqcup I}:=\vecN_{\Lambda}\sqcup\vecN_I$.
Let $\Ltilde(\lefttop{I}\psitilde_{\vecu}(\nbigt)):=M\bigl(
  N(\Lambda\sqcup I) \bigr)[w]$.
Recall that we have
the specialization of
the good system of irregular values $\vecnbigi$
to $D_I$, denoted by $\vecnbigi(I)_{|D_I}$.
(See \S\ref{subsection;13.5.6.30}.)
The induced tuple
$(D_I,\del D_I,\vecnbigi(I)_{|D_I})$
is denote by
$\vecD_I$.

\begin{prop}
\label{prop;11.2.23.31}
We have
$\bigl(\lefttop{I}\psitilde_{\vecu}(\nbigt),W,\vecN_{\Lambda\sqcup I}
 \bigr)
\in\nbigp(\vecD_I,\Lambda\sqcup I,w)$.
\end{prop}
\pf
By an inductive argument,
we have only to consider the case $I=\{i\}$.
We have two natural filtrations of
$\lefttop{i}\psitilde_u(\nbigt)$.
One is induced by $W$.
The other is obtained as the monodromy weight filtration of
$N(\Lambda)$.
By Corollary \ref{cor;11.2.19.11},
they are the same.
Because $\Ltilde$ is the relative monodromy filtration of
$N(\Lambda\sqcup \{i\})$ with respect to $W$,
we have natural isomorphisms
$\Gr^{\Ltilde}\bigl(
 \lefttop{i}\psitilde_u(\nbigt)
 \bigr)\simeq
 \Gr^{\Ltilde}\Gr^W
 \bigl(\lefttop{i}\psitilde_u(\nbigt)\bigr)
\simeq
  \Gr^{\Ltilde}
 \bigl(\lefttop{i}\psitilde_u(\Gr^W\nbigt)\bigr)$.
Hence, according to a special case of 
\S12.7 of \cite{mochi7},
we obtain that 
$\bigl(
 \lefttop{i}\psitilde_u(\nbigt),\Ltilde\bigr)$
satisfies {\bf Adm0}.
Then, it is easy to check 
$\bigl(
 \lefttop{i}\psitilde_u(\nbigt),\Ltilde,\vecN_{\Lambda\sqcup \{i\}}\bigr)$
satisfies the assumptions for Proposition \ref{prop;11.2.19.10}
for $\vecD_i$.
Hence, 
$\bigl(
 \lefttop{i}\psitilde_u(\nbigt),\Ltilde,\vecN_{\Lambda\sqcup
 \{i\}}\bigr)$ 
is an object in
 $\nbigp(\vecD_i,\Lambda\sqcup\{i\},w)$.
\hfill\qed

\vspace{.1in}
We obtain an exact functor
$\lefttop{I}\psitilde_{\vecu}:
 \nbigp(\vecX,\Lambda,w)
\lrarr
 \nbigp(\vecD_I,\Lambda\sqcup I,w)$.

\vspace{.1in}
We give a variant.
Suppose that $\vecnbigi$ is unramified.
We obtain
$\vecnbigi(-\gminia,I)_{|D_I}$
as in \S\ref{subsection;13.5.6.30}.
We set
$\vecD_I(-\gminia):=
 (D_I,\del D_I,\vecnbigi(-\gminia,I)_{|D_I})$.

\begin{cor}
\label{cor;11.2.19.12}
Let $(\nbigt,W,\vecN)\in\nbigp(\vecX,\Lambda,w)$.
For any $\vecu\in(\real\times\cnum)^I$
and $\gminia\in\Irr(\nbigt,I)$,
we have
$\bigl(
 \lefttop{I}\psitilde_{\gminia,\vecu}(\nbigt),
 \vecN_{\Lambda\sqcup I}
 \bigr)$ are objects in
 $\nbigp(\vecD_I(-\gminia),
\Lambda\sqcup I,w)$.
\hfill\qed
\end{cor}

\section{Admissible IMTM}
\label{subsection;10.11.10.10}

\subsection{Definitions}

An object $(\nbigt,W,L,\vecN)\in
\MTS^{\adm}(\vecX,\Lambda)^{\fil}$
is called pre-admissible
$\Lambda$-IMTM on $\vecX$
if the following holds:
\index{pre-admissible}
\begin{description}
\item[\bf (R3)]
 $(\nbigt,W,L,\vecN)_{|X\setminus D}
 \in \IMTM(X\setminus D,\Lambda)$,
 and
 $\Gr^L_w(\nbigt,W,\vecN)\in
 \nbigp(\vecX,\Lambda,w)$.
\end{description}
We shall impose an additional condition
for admissibility.

\subsubsection{The smooth divisor case}

Let us consider the case that $D$ is smooth.
Let $(\nbigt,W,L,\vecN_{\Lambda})$
be a pre-admissible $\Lambda$-IMTM on $\vecX$.
Let $\nbigm_i$ be the underlying smooth
$\nbigr_{X(\ast D)}$-modules,
which is equipped with the KMS-structure.
It is called admissible,
if the following holds:
\index{admissible}
\begin{description}
\item[\bf (R4)]
 The nilpotent endomorphisms
 $N$ on $\bigl(\nbigg_u(\nbigm_i),L\bigr)$
 $(i=1,2)$
 have relative monodromy filtrations.
\end{description}

\begin{lem}
Let $X=\Delta^n$ and $D=\{z_1=0\}$.
\begin{itemize}
\item
For $u\in\real\times\cnum$,
the tuple
$\bigl(\psitilde_u(\nbigt,L),
 \vecN_{\Lambda\sqcup\bullet}\bigr)$
is $(\Lambda\sqcup\bullet)$-IMTM
on $D$.
\item
Suppose that $\nbigt$ is unramified.
Then, for each $\gminia\in\Irr(\nbigt)$ and
$u\in\real\times\cnum$,
the tuple
$\bigl(\psitilde_{\gminia,u}(\nbigt,L),
 \vecN_{\Lambda\sqcup\bullet}\bigr)$
is $(\Lambda\sqcup\bullet)$-IMTM
on $D$.
\end{itemize}
\end{lem}
\pf
We have only to consider the unramified case.
By the compatibility of $L$ with the KMS-structure,
$\Gr^L_w\bigl(\psitilde_{\gminia,u}(\nbigt),
 \vecN_{\Lambda\sqcup\bullet}\bigr)$
comes from
$\Gr^L_w(\nbigt,\vecN)$.
Hence, it is a  $(w,\Lambda\sqcup\bullet)$-polarizable
mixed twistor structure.
By the assumption,
$N_{\bullet}$ has a relative monodromy filtration.
For generic $\lambda$,
$N_{i|\lambda}$ $(i\in\Lambda)$
have relative monodromy filtrations.
Then, the claim follows from 
Proposition \ref{prop;10.9.27.31}.
(We can also deduce it directly from
Proposition \ref{prop;11.2.1.1}
and Lemma \ref{lem;11.2.1.3}.)
\hfill\qed

\vspace{.1in}
We can also deduce the following lemma
from Proposition \ref{prop;10.9.27.31}
(or easier
Proposition \ref{prop;11.2.1.1}
and Lemma \ref{lem;11.2.1.3}).
\begin{lem}
We can replace the above admissibility condition
with the following:
\begin{itemize}
\item
There exists $U\subset \cnum_{\lambda}$
with $|U|=\infty$
such that
$\bigl(\nbigg_{u}(\nbigm_1),
 L,N\bigr)_{(\lambda,P)}
 \in\Vect_{\cnum}^{\RMF}$
for any $(\lambda,P)\in U\times D$.
\hfill\qed
\end{itemize}
\end{lem}

\subsubsection{The normal crossing case}

Let us consider the case that $D$ is normal crossing.
A pre-admissible 
$\Lambda$-IMTM $(\nbigt,W,L,\vecN_{\Lambda})$
on $\vecX$
is called admissible, if the following holds:
\begin{itemize}
\item
For any smooth point $P\in D$,
there exists a small neighbourhood
$X_P$ of $P$ such that
$(\nbigt,W,L,\vecN_{\Lambda})_{|X_P}$
is admissible in the sense of
the smooth divisor case.
\end{itemize}

The following proposition is an analogue of
Proposition \ref{prop;10.9.29.30}.
\begin{prop}
\label{prop;11.2.1.20}
Let $P\in D_I^{\circ}$.
We take a coordinate around $P$.
\begin{itemize}
\item
For $\vecu\in(\real\times\cnum)^I$,
the induced 
$\bigl(
 \lefttop{I}\psitilde_{\vecu}(\nbigt),
 L,\vecN_{\Lambda\sqcup I}
 \bigr)$
is a $(\Lambda\sqcup I)$-IMTM
on $D_{P,I}$.
\item
If $\nbigt$ is unramified around $P$,
for any $\gminia\in\Irr(\DD,P)$
and for any $\vecu\in(\real\times\cnum)^I$,
$\bigl(
 \lefttop{I}\psitilde_{\gminia,\vecu}(\nbigt),L,
 \vecN_{\Lambda\sqcup I}\bigr)$
is a $(\Lambda\sqcup I)$-IMTM
on $D_{P,I}$.
\hfill\qed
\end{itemize}
\end{prop}

The following lemma is obvious
by definition.
\begin{lem}
An admissible mixed twistor structure
is equivalent to
an admissible $\Lambda$-IMTM
with trivial morphisms.
An admissible
 $(w,\Lambda)$-polarizable mixed twistor structure
is equivalent to
an admissible $\Lambda$-IMTM.
\hfill\qed
\end{lem}

\subsection{Category of
 admissible variation of IMTM}

Let $\nbigm(\vecX,\Lambda)$
denote the full subcategory of admissible
$\Lambda$-IMTM on $\vecX$
in $\MTS^{\adm}(\vecX,\Lambda)^{\fil}$ 

\index{category $\nbigm(X,D,\Lambda)$}

\begin{prop}
\label{prop;10.11.9.20}
The family of the categories
$\{\nbigm(\vecX,\Lambda)\,|\,
 \Lambda\}$ have the property
{\bf M0--3}.
\end{prop}
\pf
The claim for {\bf M1} is clear.
Let us consider {\bf M0}.
Let $F:(\nbigt_1,W,L,\vecN_{\Lambda})
\lrarr(\nbigt_2,W,L,\vecN_{\Lambda})$
be a morphism in $\nbigm(\vecX,\Lambda)$.
We have the objects
$(\Ker F,W,\vecN_{\Lambda})$,
$(\Image F,W,\vecN_{\Lambda})$
and $(\Cok F,W,\vecN_{\Lambda})$
in $\MTS^{\adm}(\vecX,\Lambda)$.
Because $F$ is strict with respect to $L$,
they are equipped with naturally induced 
filtrations $L$ in the category
$\MTS^{\adm}(\vecX,\Lambda)$.
We also obtain
$\Ker \Gr^LF=\Gr^L\Ker F$,
$\Image \Gr^LF=\Gr^L\Image F$ and
$\Cok \Gr^LF=\Gr^L\Cok F$.
Hence, we obtain that the objects
$(\Ker F,W,L,\vecN_{\Lambda})$,
$(\Image F,W,L,\vecN_{\Lambda})$
and $(\Cok F,W,L,\vecN_{\Lambda})$
are pre-admissible.
To check the admissibility,
we may assume that 
$X=\Delta^n$,
$D=\{z_1=0\}$
and that $\nbigt_i$ are unramified.
Then, we have
the natural isomorphisms
\[
 \psitilde_{\gminia,u}(\Ker F)
\simeq
 \Ker \psitilde_{\gminia,u}(F),
\]
\[
  \psitilde_{\gminia,u}(\Image F)
\simeq
 \Image \psitilde_{\gminia,u}(F),
\]
\[
 \psitilde_{\gminia,u}(\Cok F)
\simeq
 \Cok \psitilde_{\gminia,u}(F).
\]
The isomorphisms are compatible
with the filtrations $W$, $L$ and a tuple of
morphisms $\vecN_{\Lambda}$,
and the nilpotent part $N_{\bullet}$
of $\Res(\DD)$.
Hence, we obtain the existence of
the relative monodromy filtration of
$N_{\bullet}$ with respect to $L$.
Thus, we obtain that
$\nbigm(\vecX,\Lambda)$ is abelian.

\vspace{.1in}

Let us consider {\bf M2}.
The claim for {\bf M2.1} is clear.
Let $\Lambda=\Lambda_0\sqcup\Lambda_1$
be a decomposition.
Let $(\nbigt,W,L,\vecN_{\Lambda})
\in\nbigm(\vecX,\Lambda)$.
Let us prove that the existence of
a relative monodromy filtration of
$N(\Lambda_1)$ with respect to $L$.
We assume the induction on 
the length of $L$.
If $L$ is pure, there is nothing to prove.
Assume that $\nbigt=L_0\nbigt$,
and that the claim holds for
$L_{-1}\nbigt$.
We construct a filtration $\Ltilde$ of $\nbigt$
in the category $\MTS^{\adm}(\vecX)$
by Deligne's inductive formula.
Its restriction to $X-D$ is 
the relative monodromy filtration
of $N(\Lambda_1)_{|X-D}$
with respect to $L_{|X-D}$.
Then, we obtain that 
$\Ltilde$ is a relative monodromy filtration
of $N(\Lambda_1)$ with respect to $L$.
Because
$\Gr^{\Ltilde}(\nbigt,W,\vecN_{\Lambda_0})
\simeq
 \Gr^{\Ltilde}\Gr^L(\nbigt,W,\vecN_{\Lambda_0})$,
$(\nbigt,W,\Ltilde,\vecN_{\Lambda_0})$
is pre-admissible.
To check the admissibility of
$(\nbigt,W,\Ltilde,\vecN_{\Lambda_0})$,
we may assume that
$D$ is smooth
and that $\nbigt$ is unramified.
We have the filtration of
$\psitilde_{\gminia,u}(\nbigt)$
induced by $\Ltilde$,
which is also denoted by $\Ltilde$.
Then, we can observe that
$\Ltilde=M(N(\Lambda_1);L)$
on $\psitilde_{\gminia,u}(\nbigt)$.
Hence, the nilpotent part of $\Res(\DD)$
has the relative monodromy filtration
with respect to $\Ltilde$,
and
$(\nbigt,W,\Ltilde,\vecN_{\Lambda_0})$
is admissible.
It is denoted by
$\res^{\Lambda}_{\Lambda_0}
 (\nbigt,W,L,\vecN_{\Lambda})$.

\vspace{.1in}
Let us consider {\bf M3}.
Let $\bullet\in\Lambda$,
and put $\Lambda_0:=\Lambda\setminus\bullet$.
Let $(\nbigt,W,L,\vecN_{\Lambda})$
be an object in $\nbigm(\vecX,\Lambda)$.
Let $(\nbigt',W,\Ltilde,\vecN_{\Lambda_0})$
be an object in $\nbigm(\vecX,\Lambda_0)$
with morphisms 
\[
 \Sigma^{1,0}
 \res^{\Lambda}_{\Lambda_0}
 (\nbigt,W,L,\vecN_{\Lambda})
\stackrel{u}{\lrarr}
 (\nbigt',W,\Ltilde,\vecN_{\Lambda_0})
\stackrel{v}{\lrarr}
 \Sigma^{0,-1}
 \res^{\Lambda}_{\Lambda_0}
 (\nbigt,W,L,\vecN_{\Lambda})
\]
in $\nbigm(\vecX,\Lambda_0)$
such that $v\circ u=N_{\bullet}$.
We obtain the filtration $L(\nbigt')$
of the object $(\nbigt',W,\Ltilde,\vecN_{\Lambda_0})$
in the category
$\nbigm(\vecX,\Lambda_0)$,
as the transfer of $L(\nbigt)$.
We set $N_{\bullet}:=u\circ v$ on $\nbigt'$.

\begin{lem}
\label{lem;10.8.25.11}
$(\nbigt',W,L,\vecN_{\Lambda})$
is an object in
$\nbigm(\vecX,\Lambda)$.
\end{lem}
\pf
We have the decomposition
$\Gr^L_w\nbigt'
=\Image\Gr_w^Lu\oplus\Ker\Gr_w^Lv$.
Because 
\[
 \Image\Gr_w^Lu\simeq
 \Sigma^{1,0}\Image \Gr^L_wN_{\bullet},
\]
it is an object in 
$\nbigp(\vecX,\Lambda,w)$.
By using the canonical decomposition,
we obtain that
$\Ker\Gr_w^Lv$  is a direct summand of
$\Gr_w^{\Ltilde}\nbigt'$.
Hence, it is an object in
$\nbigp(\vecX,\Lambda_0,w)$.
Thus,  we obtain that $\Gr^L_w\nbigt'$ is 
an object in $\nbigp(\vecX,\Lambda,w)$.

The existence of 
relative monodromy filtration of
$N_i$ $(i\in\Lambda)$ with respect to $L$
follows from Proposition \ref{prop;10.7.22.30}.
Namely,
we construct filtrations for $N_i$
by using Deligne's formula.
They give relative monodromy filtrations for $N_i$
on $X\setminus D$ by Proposition \ref{prop;10.7.22.30}.
Hence, they are relative monodromy filtration on $X$.
For the admissibility condition,
we may assume that $D$ is smooth
and that $\nbigt$ and $\nbigt'$ are unramified.
The existence of 
relative monodromy filtration of
the nilpotent part of the residue with respect to $L$
also follows from Proposition \ref{prop;10.7.22.30}.
\hfill\qed

\vspace{.1in}
Thus, the proof of Proposition 
\ref{prop;10.11.9.20} is finished.
\hfill\qed

\subsection{A remark on nearby cycle functor}
\label{subsection;11.4.13.1}

Let us give a remark related with nearby cycle functor.
Let $\II^{a,b}$ be the Beilinson IMTM.
For $\vecp\in\seisuu_{\geq 0}^{\Lambda}$,
we consider the induced 
$(\Lambda\sqcup\bullet)$-IMTM
$\II^{a,b}_{\bullet,\vecp}:=\bigl(
 \II^{a,b},\vecp N_{\II},N_{\II}\bigr)$.
We obtain a twistor nilpotent orbit
$\IItilde^{a,b}_{\bullet,\vecp}:=
 \TNIL_{\bullet}\bigl(
 \II^{a,b}_{\bullet,\vecp}\bigr)$.
Let $K\subset\ellsitabar$
and $\vecm\in\seisuu_{>0}^{K}$.
We put $g:=\vecz^{\vecm}$,
which gives $X\setminus \bigcup_{i\in K}D_i\lrarr \cnum^{\ast}$.
We put 
$\IItilde_{\vecm,\vecp}^{a,b}:=
 g^{\ast}\IItilde^{a,b}_{\bullet,\vecp}$.

Let $(\nbigt,L,\vecN)\in\nbigm(\vecX,\Lambda)$.
We obtain 
$\bigl(
 \Pi^{a,b}_{\vecm,\vecp}\nbigt,L,\vecNtilde
 \bigr):=
 (\nbigt,L,\vecN)\otimes
 \IItilde^{a,b}_{\vecm,\vecp}
 \in\nbigm(\vecX,\Lambda)$.
Take $I\subset \{i\in\Lambda\,|\,p_i>0\}$,
and we consider the following morphism
in $\nbigm(\vecX,\Lambda)$:
\[
 \Sigma^{|I|,0}\bigl(
 \Pi^{0,N}_{\vecm,\vecp}\nbigt,
 \Ntilde_{I!}L\bigr)
\lrarr
 \Sigma^{0,-|I|}\bigl(
 \Pi^{0,N}_{\vecm,\vecp}\nbigt,
 \Ntilde_{I\ast}L\bigr)
\]
The cokernel in $\nbigm(\vecX,\Lambda)$
is denoted by
$\bigl(
 \psi^{(0)}_{\vecm,\vecp}(\nbigt)_I,\Ltilde
 \bigr)$.
It is equipped with the filtration $\Ltilde$
as the object in $\nbigm(\vecX,\Lambda)$.
It is also equipped with the filtration $L$
naively induced by the filtration of $\nbigt$.
We immediately obtain the following lemma 
from Corollary \ref{cor;10.8.30.10}.
\begin{lem}
\label{lem;10.10.13.20}
We have $M(N;L)[1]=\Ltilde$.
\hfill\qed
\end{lem}

\section{Specialization of admissible mixed twistor structure}

\subsection{Statement}
\label{subsection;11.2.23.40}

Let $X:=\Delta^n$,
$D_i:=\{z_i=0\}$
and $D:=\bigcup_{i=1}^{\ell}D_i$.
Let $(\nbigt,L)\in\MTS^{\adm}(\vecX)$.
Assume that $\nbigt$ is unramified.
For $I\subset\ellsitabar$,
$\gminia\in\Irr(\nbigt,I)$
and $\vecu\in(\real\times\cnum)^I$,
we obtain an object
$\bigl(
 \lefttop{I}\psitilde_{\gminia,\vecu}(\nbigt),L,
 \vecN_I
 \bigr)$ in
$\VTS(X,D,I)^{\fil}$.
Its restriction to $D_I^{\circ}$ is an object in
$\nbigm(D_I^{\circ},I)$.
The relative monodromy filtration $W=M(N(I);L)$ of 
$\lefttop{I}\psitilde_{\gminia,\vecu}(\nbigt)_{|D_I^{\circ}}$
is extended to a filtration of
$\lefttop{I}\psitilde_{\gminia,\vecu}(\nbigt)$,
by Deligne's inductive formula.
Thus, we obtain an object
$\bigl(
 \lefttop{I}\psitilde_{\gminia,\vecu}(\nbigt),
 W
 \bigr)$
in $\MTS(\vecD_I)$.
We will prove the following proposition in 
\S\ref{subsection;11.1.14.2}.

\begin{prop}
\label{prop;10.11.10.11}
$\bigl(
 \lefttop{I}\psitilde_{\gminia,\vecu}(\nbigt),
 W
 \bigr)$
is admissible.
\end{prop}

For $(\nbigt,L)\in \MTS^{\adm}(\vecX)$
which is not necessarily unramified,
we have
$\lefttop{I}\psitilde_{\vecu}(\nbigt)
\in\VTS(D_I,\del D_I)$.
It is equipped with the filtration
$W=M\bigl(N(I);L\bigr)$.
By Proposition \ref{prop;10.11.10.11},
$\bigl(
 \lefttop{I}\psitilde_{\vecu}(\nbigt),W
 \bigr)$
is an object in
$\MTS^{\adm}(\vecD_I)$.
Thus, we obtain an exact functor
\[
 \lefttop{I}\psitilde_{\vecu}:
 \MTS^{\adm}(\vecX)
\lrarr
 \MTS^{\adm}(\vecD_I).
\]

The following corollary is clear by Proposition 
\ref{prop;10.11.10.11}
and the definition of $\nbigm(\vecX,\Lambda)$.

\begin{cor}
\label{cor;11.2.23.42}
Let $(\nbigt,W,L,\vecN_{\Lambda})$
be an object in
$\nbigm(\vecX,\Lambda)$.
\begin{itemize}
\item
For any $\vecu\in(\real\times\cnum)^I$,
$(\lefttop{I}\psitilde_{\vecu}(\nbigt),
 L,\vecN_{\Lambda\sqcup I})$
are objects in 
$\nbigm(\vecD_I,\Lambda\sqcup I)$.
Thus, we obtain an exact functor
$\lefttop{I}\psitilde_{\vecu}:
 \nbigm(\vecX,\Lambda)
\lrarr
\nbigm(\vecD_I,\Lambda\sqcup I)$.
\item
Assume that $\nbigt$ is unramified.
For any $I\subset\ellsitabar$,
$\vecu\in(\real\times\cnum)^I$
and $\gminia\in\Irr(\nbigt,I)$,
$(\lefttop{I}\psitilde_{\gminia,\vecu}(\nbigt),
 L,\vecN_{\Lambda\sqcup I})$
are objects in 
$\nbigm(\vecD_I(-\gminia),\Lambda\sqcup I)$.
(See {\rm\S\ref{subsection;11.2.23.30}}
for $\vecD_I(-\gminia)$.)
\hfill\qed
\end{itemize}
\end{cor}

\subsection{Some notation}

We introduce some notation
which will be used in \S\ref{section;11.4.3.5}.
We have the abelian category
$\nbiga_{D_I}:=
 \MTS^{\adm}(\vecD_I)$.
\index{category $\nbiga_{D_I}$}
For any finite set $\Lambda$,
we put
\[
 \nbigp_{w,D_I}(\Lambda):=
 \nbigp(\vecD_I,\Lambda,w),
\quad
 \nbigm_{D_I}(\Lambda):=
 \nbigm(\vecD_I,\Lambda).
\]
\index{category $\nbigp_{w,D_I}$}
\index{category $\nbigm_{w,D_I}$}
The family
$\bigl\{\nbigp_{w,D_I}(\Lambda)\bigr\}$
satisfies the conditions {\bf P0}--{\bf P3},
and the family $\bigl\{\nbigm_{D_I}(\Lambda)\bigr\}$
satisfies the conditions {\bf M0}--{\bf M3}
in {\rm\S\ref{subsection;10.10.12.21}}.

We obtain the categories
$ML\nbiga_{D_I}(\Lambda_1,\Lambda_2)$
and $M'L\nbiga_{D_I}(\Lambda_1,\Lambda_2)$
as in \S\ref{subsection;10.10.13.2}
and \S\ref{subsection;10.10.13.1}.
We have a natural equivalence
$ML\nbiga_{D_I}(\Lambda_1,\Lambda_2)\lrarr
M'L\nbiga_{D_I}(\Lambda_1,\Lambda_2)$
as in Theorem \ref{thm;10.10.12.12}.

For $I\subset J$,
we have a naturally defined exact functor
$\psi_{J,I}:\nbigm_{D_I}(\Lambda)
\lrarr\nbigm_{D_J}\bigl(\Lambda\sqcup (J\setminus I)\bigr)$
given by
$\lefttop{J\setminus I}\psitilde_{\vecu_0}$,
where
$\vecu_0=(-1,0)^{J\setminus I}\in
 (\real\times\cnum)^{J\setminus I}$.
\index{functor $\psi_{J,I}$}
We have the following commutative diagram
of the functors:
{\small
\begin{equation}
\label{eq;10.10.13.10}
\begin{CD}
 ML\nbiga_{D_I}(\Lambda_1,\Lambda_2)
 @>{\psi_{JI}}>>
  ML\nbiga_{D_J}\bigl(\Lambda_1,\Lambda_2
 \sqcup(J\setminus I)\bigr)
 @<{\psi_{J\setminus I}}<<
 ML\nbiga_{D_J}\bigl(
 \Lambda_1\sqcup(J\setminus I),
 \Lambda_2
 \bigr)
 \\
 @V{\simeq}VV @V{\simeq}VV 
 @V{\simeq}VV\\
 M'L\nbiga_{D_I}(\Lambda_1,\Lambda_2)
 @>{\psi_{JI}}>>
 M'L\nbiga_{D_J}\bigl(\Lambda_1,\Lambda_2
 \sqcup (J\setminus I)\bigr)
 @<{\psi_{J\setminus I}}<<
 M'L\nbiga_{D_J}\bigl(
 \Lambda_1\sqcup(J\setminus I),
 \Lambda_2
 \bigr)
\end{CD}
\end{equation}
}
The horizontal arrows are exact functors.

\subsection{Proof of Proposition \ref{prop;10.11.10.11}}
\label{subsection;11.1.14.2}

\begin{prop}
\label{prop;10.11.10.5}
Suppose that $D$ is smooth.
Let $(\nbigt,W,L,\vecN_{\Lambda})$ be an object in 
$\MTS(\vecX)(\Lambda)^{\fil}$ such that
(i) the filtration $L$ is compatible with the KMS-structure,
(ii) $(\nbigt,W,L,\vecN_{\Lambda})$
satisfies {\bf R3} and {\bf R4} 
in {\rm\S\ref{subsection;10.11.10.10}}.
Then, it is an object in $\nbigm(\vecX,\Lambda)$,
and hence
in $\MTS^{\adm}(\vecX,\Lambda)^{\fil}$.
\end{prop}
\pf
We may assume that $\nbigt$ is unramified.
By {\bf(R4)},
we have $M(N;L)$ for
$\psitilde_{\gminia,u}(\nbigt,L)$.
For generic $\lambda$,
we have $M(N_i;L)$ for 
$\psitilde_{\gminia,u}(\nbigt,L)_{|\lambda\times D}$.
Hence, 
by Proposition \ref{prop;10.9.27.31},
we obtain that
$\bigl(
\psitilde_{\gminia,u}(\nbigt,L),
 \vecN_{\Lambda\sqcup\bullet}\bigr)$
is a $(\Lambda\sqcup \bullet)$-IMTM.
Let $\nbigm_c$ $(c=1,2)$ be the underlying $\nbigr$-modules.
We have the irregular decomposition
$\nbigm_{c|\nbigdhat}
=\bigoplus_{\gminia\in\Irr(\nbigt)}
 \nbigm_{c,\gminia}$.
We have the corresponding decomposition
\[
  W_k(\nbigm_c)_{|\nbigdhat}
=\bigoplus_{\gminia\in\Irr(\nbigt)}
W_k(\nbigm_{c,\gminia}).
\] 
We have an $\nbigo_{\nbigxzero}$-coherent sheaf
$\nbigqzero_{b}W_k(\nbigm_c):=
  \nbigqzero_b\nbigm_c\cap W_k(\nbigm_c)$.

We have the induced filtration
$L$ on $\Gr^{\nbigqzero}_a\nbigm_c$.
We have the nilpotent endomorphism
$N(\Lambda)$ on 
$\Gr^{\nbigqzero}_a\nbigm_c$.
Because 
$\Gr^{\nbigqzero}_a\nbigm_c
=\bigoplus_{\gminia}
 \bigoplus_{\paramap(\lambda_0,u)=a}
 \psizero_{u,\gminia}\nbigm_c$,
there exists a relative monodromy filtration
$W=M\bigl(N(\Lambda);L\bigr)$ 
on $\Gr^{\nbigqzero}_a\nbigm_c$.
\begin{lem}
$\nbigqzero_{\ast}W_m\nbigm_c$ 
is a good-KMS family of 
meromorphic $\lambda$-flat bundles,
and we have
\begin{equation}
 \label{eq;10.10.10.10}
 W_{m}\bigl(
 \Gr^{\nbigqzero}_a\nbigm_c
 \bigr)
=\Gr^{\nbigqzero}_a\bigl(
 W_m\nbigm_c\bigr).
\end{equation}
\end{lem}
\pf
If $L$ is pure, it follows from the assumption
{\bf R3}.
We assume $L_{k}\nbigt=\nbigt$,
and the claim holds for 
$(L_{k-1}\nbigt,W, L,\vecN)$.
The filtration $W(\nbigm_c)$
is constructed by Deligne's formula:
\begin{equation}
\label{eq;10.10.10.5}
 W_{-i+k}\nbigm_c
=W_{-i+k}L_{k-1}\nbigm_c
+N(\Lambda)^i\bigl(
 W_{i+k}\nbigm_c
 \bigr)
\end{equation}
\begin{equation}
\label{eq;10.10.10.6}
 W_{i+k}\nbigm_c
=\Ker\Bigl(
 N(\Lambda)^{i+1}:
 \nbigm_c\lrarr
 \nbigm_c\big/
 W_{-i-2+k}\nbigm_c
 \Bigr)
\end{equation}

Assume that we know the claims for
$W_{-i-2+k}\nbigm_c$.
Then, 
$\nbigm_c/W_{-i-2+k}\nbigm_c$
has the induced KMS-structure,
and we have
\[
 \Gr^{\nbigqzero}_{\ast}\bigl(
 \nbigm_c/W_{-i-2+k}\nbigm_c
 \bigr)
\simeq
 \Gr^{\nbigqzero}_{\ast}\nbigm_c\big/
 \Gr^{\nbigqzero}_{\ast}W_{-i-2+k}\nbigm_c.
\]
Because
$(\psitilde_{\gminia,u}(\nbigt,L),\vecN_{\Lambda\sqcup\bullet})$
is a $(\Lambda\sqcup\bullet)$-IMTM,
we obtain that the cokernel of 
the induced morphism
\[
N(\Lambda)^{i+1}:\Gr_{\ast}^{\nbigqzero}\nbigm_c
\lrarr
  \Gr^{\nbigqzero}_{\ast}\bigl(
 \nbigm_c/W_{-i-2+k}\nbigm_c
 \bigr)
\]
is strict. Hence, we obtain that the claim holds for
$W_{i+k}\nbigm_c$.

Assume that the claim holds for
$W_{i+k}\nbigm_c$.
Because
$(\psitilde_{\gminia,u}(\nbigt,L),\vecN_{\Lambda\sqcup\bullet})$
is a $(\Lambda\sqcup\bullet)$-IMTM,
we obtain that the cokernel of the following morphism
is strict:
\[
\begin{CD}
 \Gr^{\nbigqzero}_{\ast}\bigl(
 W_{-i+k}L_{k-1}\nbigm_c
 \bigr)
\oplus
 \Gr^{\nbigqzero}_{\ast}\bigl(
 W_{i+k}\nbigm_c
 \bigr)
 @>{{\rm inclusion}+N(\Lambda)^i}>>
 \Gr^{\nbigqzero}_{\ast}\nbigm_c
 \end{CD}
\]
Hence, we obtain the claims for
$W_{-i+k}\nbigm_c$.
Thus, the induction can proceed.
\hfill\qed

\vspace{.1in}
Hence, $W$ is compatible with the KMS-structure,
and the nilpotent part of the residue
on $\Gr^{\nbigqzero}_a\nbigm_c$
has a relative monodromy filtration.
Thus, we obtain Proposition
\ref{prop;10.11.10.5}.
\hfill\qed

\vspace{.1in}
We have used the following standard lemma.

\begin{lem}
Let $\nbigm_i$ $(i=1,2)$ be 
unramifiedly $\vecnbigi$-good-KMS
$\nbigr_{X(\ast D)}$-modules.
Let $F:\nbigm_1\lrarr\nbigm_2$ be 
a morphism such that 
(i) $\Cok\psitilde_{\gminia,u}(F)$ are strict
for any $\gminia$ and $u$.
Then, $F$ is strictly compatible 
with the KMS-structure,
and $\Ker F$, $\Image F$
and $\Cok F$ have naturally induced
KMS-structure.
Moreover, we have natural isomorphisms
$\psitilde_{\gminia,u}(\Ker F)\simeq
 \Ker\psitilde_{\gminia,u}(F)$,
$\psitilde_{\gminia,u}(\Image F)\simeq
 \Image\psitilde_{\gminia,u}(F)$,
and
$\psitilde_{\gminia,u}(\Cok F)\simeq
 \Cok\psitilde_{\gminia,u}(F)$.
\hfill\qed
\end{lem}

Let us prove
Proposition \ref{prop;10.11.10.11}.
Let $(\nbigt,L)\in \MTS^{\adm}(\vecX)$.
We obtain {\bf R3} for 
$\bigl(\lefttop{I}\psitilde_{\vecu}(\nbigt),W,L,\vecN_I\bigr)$
from Proposition \ref{prop;10.9.29.30} and 
Corollary \ref{cor;11.2.19.12}.
For the remaining conditions,
we may assume that $\del D_I$ is smooth.
We obtain {\bf R4}
from Proposition \ref{prop;11.2.1.20}.
Then, applying Proposition \ref{prop;10.11.10.5},
we obtain the claim of
Proposition \ref{prop;10.11.10.11}.
\hfill\qed

\section{Integrable case}

\subsection{Admissible mixed twistor structure}
An integrable mixed twistor structure on $\vecX$
is an object $(\nbigt,\!L)$ in $\VTSint(X,\!D)^{\fil}$ 
satisfying
(i) the underlying filtered smooth $\nbigr_{X(\ast D)}$-triple
is a mixed twistor structure on $\vecX$,
(ii) each $\Gr^{L}_w(\nbigt)$ has an integrable polarization.
Let $\MTSint(\vecX)\subset\VTSint(X,D)^{\fil}$ 
denote the corresponding full subcategory.
It is an abelian category.
\index{category $\MTSint(X,D)$}

An object $(\nbigt,L)\in\MTSint(\vecX)$ is called 
admissible,
if (i) the underlying object in $\MTS(\vecX)$ is 
admissible,
(ii) for each $w$, $\Gr^L_w(\nbigt)$
 is obtained as the canonical prolongation of
 a good wild integrable polarizable variation of
 pure twistor structure of weight $w$.
Let $\MTS^{\integral\adm}(\vecX)
 \subset\MTS^{\integral}(\vecX)$
be the corresponding full subcategory.
We immediately obtain the following 
from Proposition \ref{prop;11.2.23.20}.
\index{category $\MTS^{\integral\adm}(X,D)$}
\index{category $\MTS^{\integral\adm}(\vecX)$}
\begin{prop}
$\MTS^{\integral\adm}(\vecX)$ is 
an abelian subcategory.
\hfill\qed
\end{prop}

The following lemma is an integrable analogue of
Lemma \ref{lem;11.2.23.21}.
\begin{lem}
\label{lem;11.2.23.22}
Let $\nbigt\in\VTS(X,D)$ come from
an integrable good wild polarizable variation of pure twistor structure
of weight $w$.
Then, it is naturally admissible
integrable mixed twistor structure on $(X,D)$.
\hfill\qed
\end{lem}

The following proposition can be reduced to
Proposition \ref{prop;10.9.29.30}.
\begin{prop}
Let $(\nbigt,L)\in\MTS^{\integral\adm}(\vecX)$.
Let $P\in D_I^{\circ}$.
We take a coordinate around $P$.
\begin{itemize}
\item
For $\vecu\in(\real\times\cnum)^I$,
we naturally have
$\bigl(
 \lefttop{I}\psitilde_{\vecu}(\nbigt,L),\vecN_I
 \bigr)_{|P}$
is an integrable $I$-IMTM.
\item
If $\nbigt$ is unramified around $P$,
for any $\gminia\in\Irr(\nbigt,P)$
and for any $\vecu\in(\real\times\cnum)^I$,
$\bigl(
 \lefttop{I}\psitilde_{\gminia,\vecu}(\nbigt,L),\vecN
 \bigr)_{|P}$
is naturally an integrable $I$-IMTM.
\hfill\qed
\end{itemize}
\end{prop}

Let us consider the specialization
as in \S\ref{subsection;11.2.23.40}.
Let $X=\Delta^n$
and $D=\bigcup_{i=1}^{\ell}\{z_i=0\}$.
Let $(\nbigt,L)\in\MTS^{\integral\adm}(\vecX)$.
Assume that $\nbigt$ is unramified.
For $I\subset\ellsitabar$,
$\gminia\in\Irr(\nbigt,I)$
and $\vecu\in(\real\times\cnum)^I$,
we have
$\bigl(
 \lefttop{I}\psitilde_{\gminia,\vecu}(\nbigt),
 W
 \bigr)$
in $\MTS^{\adm}(\vecD_I)$.
It is naturally a smooth filtered 
integrable  $\nbigr_{D_I(\ast \del D_I)}$-triple.
We can check that
$\Gr^W\bigl(\lefttop{I}\psitilde_{\gminia,\vecu}(\nbigt)\bigr)$
has an integrable polarization.
Hence, we obtain the following.
\begin{prop}
We naturally have
$\bigl(
 \lefttop{I}\psitilde_{\gminia,\vecu}(\nbigt),
 W
 \bigr)$
in $\MTS^{\integral\adm}(\vecD_I)$.
\hfill\qed
\end{prop}

\subsection{Admissible polarizable mixed twistor structure}

For any finite set $\Lambda$,
we consider the category
$\MTS^{\integral\adm}(\vecX,\Lambda)\!:=\!
 \MTS^{\integral\adm}(\vecX)(\Lambda)$.
\index{category $\MTS^{\integral\adm}(\vecX,\Lambda)$}
An object 
$(\nbigt,W,\vecN)$
in $\MTS^{\integral\adm}(\vecX,\Lambda)$
is called $(w,\Lambda)$-polarizable,
if 
(i) the underlying object in $\MTS^{\adm}(\vecX,\Lambda)$
is $(w,\Lambda)$-polarizable,
(ii) it has an integrable $S:(\nbigt,W,\vecN_{\Lambda})
\lrarr(\nbigt,W,\vecN_{\Lambda})^{\ast}\otimes\newTate(-w)$
such that
$(\nbigt,W,\vecN_{\Lambda},S)_{|X\setminus D}$
is a polarized mixed twistor structure on $X\setminus D$.

Let $\nbigp^{\integral}(\vecX,\Lambda,w)
\subset
 \MTS^{\integral\adm}(\vecX,\Lambda)$ denote 
the full subcategory of admissible 
$(w,\Lambda)$-polarizable mixed twistor structure
on $\vecX$.
\index{category $\nbigp^{\integral}(\vecX,\Lambda,w)$}
The following is an integrable analogue
of Proposition \ref{prop;10.11.9.11}.
\begin{prop}
The family of the categories 
$\bigl\{
 \nbigp^{\integral}(\vecX,\Lambda,w)\,\big|\,
 \Lambda,w
 \bigr\}$
has the property {\bf P0--3}.
\hfill\qed
\end{prop}
The following is an analogue of Lemma \ref{lem;11.2.23.23},
and easy to see.
\begin{lem}
\mbox{{}}
\begin{itemize}
\item
Let $(\nbigt,W,\vecN_{\Lambda})
 \in\nbigp^{\integral}(\vecX,\Lambda,w)$.
Then, its dual 
$(\nbigt,W,\vecN_{\Lambda})^{\lor}$
is an object in
 $\nbigp^{\integral}(X,D,-\vecnbigi,\Lambda,-w)$.
\item
Let $(\nbigt_1,W,\vecN_{\Lambda})$ 
be an object in $\nbigp^{\integral}(\vecX,\Lambda,w)$.
Let $(\nbigt_2,W,\vecN_{\Lambda})$
be an admissible $(w,\Lambda)$-polarizable
integrable pure twistor structure on $(X,D,\veczero)$.
Then, 
$(\nbigt_1,W,\vecN_{\Lambda})
\otimes
 (\nbigt_2,W,\vecN_{\Lambda})$
is an object in 
$\nbigp^{\integral}(\vecX,\Lambda,w)$.
\item
The functors $j^{\ast}$,
$\ast$ and
$\gammatilde_{\sm}^{\ast}$
on $\MTS^{\integral\adm}(\vecX,\Lambda)$
preserve 
$\nbigp^{\integral}(\vecX,w,\Lambda)$.
\hfill\qed
\end{itemize}
\end{lem}

Let us consider the specialization
as in \S\ref{subsection;11.2.23.30}.
Let $X=\Delta^n$, $D_i=\{z_i=0\}$
and $D=\bigcup_{i=1}^{\ell}D_i$.
Let $(\nbigt,W,\vecN_{\Lambda})\in
 \nbigp^{\integral}(\vecX,\Lambda,w)$.
We obtain
$\lefttop{I}\psitilde_{\vecu}(\nbigt)\in \VTS^{\integral}(D_I,\del D_I)$
with the induced morphisms $\vecN_{\Lambda}$.
We also have the naturally induced morphisms $\vecN_I$.
We set $\vecN_{\Lambda\sqcup I}:=\vecN_{\Lambda}\sqcup\vecN_I$.
Let $\Ltilde(\lefttop{I}\psitilde_{\vecu}(\nbigt)):=M\bigl(
  N(\Lambda\sqcup I) \bigr)[w]$.
The following proposition is an integrable analogue of
Proposition \ref{prop;11.2.23.31}.
\begin{prop}
We have
$\bigl(\lefttop{I}\psitilde_{\vecu}(\nbigt),W,\vecN_{\Lambda\sqcup I}
 \bigr)
\in\nbigp^{\integral}(\vecD_I,\Lambda\sqcup I,w)$.
\hfill\qed
\end{prop}

\subsection{Admissible IMTM}

Let $(\nbigt,W,L,\vecN)\in
\MTS^{\integral\adm}(\vecX,\Lambda)^{\fil}$.
It is called an integrable admissible $\Lambda$-IMTM,
if (i) it is an admissible $\Lambda$-IMTM,
(ii) $\Gr^L_w(\nbigt,W,\vecN)$
has an integrable polarization.
Let $\nbigm^{\integral}(\vecX,\Lambda)\subset
\MTS^{\integral\adm}(\vecX,\Lambda)$ 
denote
the full subcategory of integrable admissible
$\Lambda$-IMTM on $\vecX$.
\index{category $\nbigm^{\integral}(\vecX,\Lambda)$}
We can easily deduce the following proposition
from Proposition \ref{prop;10.11.9.20}.
\begin{prop}
The categories
$\nbigm^{\integral}(\vecX,\Lambda)$ have the property
{\bf M0--3}.
\hfill\qed
\end{prop}

The following is an integrable analogue of
Corollary \ref{cor;11.2.23.42}.
\begin{cor}
Let $(\nbigt,W,L,\vecN_{\Lambda})$
be an unramified object in
$\nbigm^{\integral}(\vecX,\Lambda)$.
\begin{itemize}
\item
For any $\vecu\in(\real\times\cnum)^I$,
$(\lefttop{I}\psitilde_{\vecu}(\nbigt),
 L,\vecN_{\Lambda\sqcup I})$
are objects in 
$\nbigm^{\integral}(\vecD_I,\Lambda\sqcup I)$.
\item
Assume that $\nbigt$ is unramified.
For any $I\subset\ellsitabar$,
$\vecu\in(\real\times\cnum)^I$
and $\gminia\in\Irr(\nbigt,I)$,
$(\lefttop{I}\psitilde_{\gminia,\vecu}(\nbigt),
 L,\vecN_{\Lambda\sqcup I})$
are objects in 
$\nbigm^{\integral}(\vecD_I(-\gminia),\Lambda\sqcup I)$.
\end{itemize}
Thus, we obtain an exact functor
$\lefttop{I}\psitilde_{\vecu}:
 \nbigm^{\integral}(\vecX,\Lambda)
\lrarr
\nbigm^{\integral}(\vecD_I,\Lambda\sqcup I)$.
\hfill\qed
\end{cor}

\chapter{Good mixed twistor $D$-modules}
\label{section;11.4.3.5}

We consider filtered $\nbigr$-triples 
locally expressed 
as the gluing of admissible mixed twistor structures
given on the intersections
of the normally crossing hypersurfaces.
We prove that they are mixed twistor $D$-modules
in \S\ref{subsection;11.4.3.41}.

\section{Good gluing data}

\subsection{An equivalence}
\label{subsection;10.8.25.12}

Let $X:=\Delta^n$
and $D:=\bigcup_{i=1}^{\ell}\{z_i=0\}$.
Let $\vecnbigi$ be a good system of 
ramified irregular values,
induced by a good set of ramified irregular values
in $\nbigotilde_{X}(\ast D)_O\big/\nbigotilde_{D,O}$,
where $O$ is the origin of $X$.
We set $\vecX:=(X,D,\vecnbigi)$.
The induced tuples
$(D_I,\del D_I,\vecnbigi(I)_{|D_I})$
are denoted by $\vecD_I$.

Let $H=\bigcup_{i\in K}\{z_i=0\}$
for some $K\subset\ellsitabar$.
We introduce categories
$\gbigg_{i}(\vecX,\ast H)$
$(i=0,1)$ given as follows.
(If $H=\emptyset$,
they are denoted by
$\gbigg_i(\vecX)$.)
We omit to denote the weight filtration
for objects in $\MTS^{\adm}(\vecD_I)$.
The argument in this subsection can work
also in the integrable case.
In the following,
$\lefttop{i}\psitilde_{-\vecdelta}$
is denoted by $\psi_i$.
For $I=(i_1,\ldots,i_m)$,
the composition
$\psi_{i_1}\circ\cdots\circ\psi_{i_m}$
is denoted by $\psi_I$.

\subsubsection{Category $\gbigg_0(\vecX,\ast H)$}
\index{category $\gbigg_0(\vecX,\ast H)$}
Objects of $\gbigg_{0}(\vecX,\!\ast\! H)$ are tuples 
$\nbigt_I\in \MTS^{\adm}(\vecD_I)$
$(I\subset\ellsitabar\setminus K)$
with morphisms in $\MTS^{\adm}(\vecD_{Ii})$
\[
\begin{CD}
 \Sigma^{1,0}
 \res_{\emptyset}^i\bigl(
 \psi_i(\nbigt_I)
 \bigr) 
@>{g_{I,i}}>>
 \nbigt_{Ii}
@>{f_{I,i}}>>
 \Sigma^{0,-1}
 \res_{\emptyset}^i\bigl(
 \psi_i(\nbigt_I)
 \bigr)
\end{CD}
\]
for $i\in\ellsitabar\setminus (I\cup K)$
such that $f_{I,i}\circ g_{I,i}=N_i$.
For $j,k\in\ellsitabar\setminus (I\cup K)$,
we impose the commutativity of 
the following diagrams:
{\small
\begin{equation}
\label{eq;10.10.1.1}
 \begin{CD}
 \Sigma^{1,0}\Bigl(
 \psi_j(\nbigt_{Ik})\Bigr)
 @>{\psi_j(f_{I,k})}>>
 \Sigma^{1,-1}\Bigl(
 \psi_{jk}(\nbigt_{I})
 \Bigr)
 \\
 @V{g_{Ik,j}}VV @V{\psi_k(g_{I,j})}VV \\
 \nbigt_{Ijk}
 @>{f_{Ij,k}}>>
 \Sigma^{0,-1}\Bigl(
 \psi_k(\nbigt_{Ij})\Bigr)
 \end{CD}
\end{equation}
}
{\small
\begin{equation}
\label{eq;10.10.1.3}
  \begin{CD}
 \Sigma^{2,0}\Bigl(
 \psi_{jk}(\nbigt_{I})
 \Bigr)
 @>{\psi_j(g_{I,k})}>>
 \Sigma^{1,0}\Bigl(
  \psi_j(\nbigt_{Ik})
 \Bigr)
 \\
 @V{\psi_k(g_{I,j})}VV
 @V{g_{Ik,j}}VV \\
 \Sigma^{1,0}\Bigl(
 \psi_k(\nbigt_{Ij})
 \Bigr)
 @>{g_{Ij,k}}>>
 \nbigt_{Ijk}
 \end{CD}
\quad\quad
  \begin{CD}
 \nbigt_{Ijk}
 @>{f_{Ik,j}}>>
 \Sigma^{0,-1}\Bigl(
 \psi_j(\nbigt_{Ik})
 \Bigr)
 \\
 @V{f_{Ij,k}}VV
 @V{\psi_{j}(f_{I,k})}VV \\
 \Sigma^{0,-1}\Bigl(
 \psi_k(\nbigt_{Ij})
 \Bigr)
 @>{\psi_{k}(f_{I,j})}>>
\Sigma^{0,-2}\Bigl(
 \psi_{jk}(\nbigt_{I})
 \Bigr)
 \end{CD}
\end{equation}
}
For $\vecnbigt^{(i)}=(\nbigt_I^{(i)})
 \in\gbigg_0(\vecX,\ast H)$ $(i=1,2)$,
a morphism
$\vecF:\vecnbigt^{(1)}\lrarr\vecnbigt^{(2)}$
in $\gbigg_0(\vecX,\ast H)$ is a tuple of morphisms
$F_I:\nbigt^{(1)}_I\lrarr\nbigt^{(2)}_I$
in $\MTS^{\adm}(\vecD_I)$
such that the following diagram is commutative:
\begin{equation}
 \label{eq;10.10.1.2}
 \begin{CD}
 \Sigma^{1,0}\Bigl(
 \psi_i(\nbigt^{(1)}_I)
 \Bigr)
@>{g^{(1)}_{I,i}}>>
 \nbigt^{(1)}_{Ii}
@>{f^{(1)}_{I,i}}>>
 \Sigma^{0,-1}\Bigl(
 \psi_i(\nbigt^{(1)}_I)
 \Bigr)
 \\
 @V{\psi_i(F_I)}VV 
 @V{F_{Ii}}VV 
 @V{\psi_i(F_I)}VV \\
\Sigma^{1,0}\Bigl(
 \psi_i(\nbigt^{(2)}_I)
 \Bigr)
@>{g^{(2)}_{I,i}}>>
 \nbigt^{(2)}_{Ii}
@>{f^{(2)}_{I,i}}>>
\Sigma^{0,-1}\Bigl(
 \psi_i(\nbigt^{(1)}_I)
\Bigr)
\end{CD}
\end{equation}

For an object  
$\vecnbigt=(\nbigt_I)$
in $\gbigg_0(\vecX,\ast H)$,
each $\nbigt_I$ is equipped with a tuple of
morphisms $\vecN_I=(N_i\,|\,i\in I)$,
given by 
$N_i:=g_{I\setminus i,i}\circ f_{I\setminus i,i}$.

\subsubsection{Category $\gbigg_1(\vecX,\ast H)$}
\index{category $\gbigg_1(\vecX,\ast H)$}

We consider objects $(\vecnbigt,\nbigl)$ in
$\gbigg_0(\vecX,\ast H)^{\fil}$
such that 
(i) $(\nbigt_I,\nbigl,\vecN_{I})
 \in \nbigm_{D_I}(I)$
 for each $I\subset\ellsitabar\setminus K$,
(ii) the tuple
 $\bigl(
 \psi_i(\nbigt_I),
 \nbigt_{Ii};
 g_{I,i},f_{I,i};\nbigl
 \bigr)$ is filtered $S$-decomposable.
Let $\gbigg_1(\vecX,\ast H)$ be the full subcategory
of such objects in 
$\gbigg_0(\vecX,\ast H)^{\fil}$.

\begin{prop}
The forgetful functor 
$\Psi:\gbigg_1(\vecX,\ast H)\lrarr
 \gbigg_0(\vecX,\ast H)$ is an equivalence.
\end{prop}
\pf
Let us prove the essential surjectivity.
Let $\vecnbigt\in\gbigg_0(\vecX,\ast H)$.
For each $I\subset J$,
we have 
\[
 \nbigu_I^{(J)}:=\bigl(
\psi_{J\setminus I}(\nbigt_I),\vecN_{J\setminus I}
\bigr)\in
 \nbigm_{D_J}(J\setminus I). 
\]
The tuple
$\bigl(\nbigu_I^{(J)}\,\big|\,
 I\subset J\bigr)$
with the following induced morphisms for 
$I_2\subset I_1$
in $\nbigm_{D_J}(J\setminus I_1)$
give an object in $ML'\nbiga_{D_J}(J)$:
\[
 \Sigma^{|I_1\setminus I_2|,0}
 \Bigl(
 \res_{J\setminus I_1}^{J\setminus I_2}
 \bigl(\nbigu_{I_2}^{(J)}\bigr)
 \Bigr)
 \lrarr
 \nbigu_{I_1}^{(J)}
 \lrarr
 \Sigma^{0,-|I_1\setminus I_2|}
 \Bigl(
 \res_{J\setminus I_1}^{J\setminus I_2}
 \bigl(\nbigu_{I_2}^{(J)}\bigr)
 \Bigr)
\]
By Theorem \ref{thm;10.10.12.12},
we have the corresponding object
in $ML\nbiga_{D_J}(J)$.
Hence, we obtain the filtrations $L^{(J)}$ of 
$\psi_{J\setminus I}(\nbigt_I)$
in $\MTS^{\adm}(\vecD_J)$
for any $I\subset J$
such that
$\bigl(
 \psi_{J\setminus I}(\nbigt_I),L^{(J)},\vecN_J
 \bigr)$
is an object in $\nbigm_{D_J}(J)$.
In particular,
we obtain a filtration $L^{(J)}$ of $\nbigt_J$,
and $(\nbigt_J,L^{(J)},\vecN_J)
\in \nbigm_{D_J}(J)$.
By using the commutative diagram
(\ref{eq;10.10.13.10}),
we obtain that 
$L^{(I)}$ on $\nbigt_I$
induces $L^{(J)}$ on 
$\psi_{J\setminus I}(\nbigt_I)$.
Hence, the tuple 
$\bigl(L^{(J)}(\nbigt_J)\,\big|\,J\subset\ellsitabar\bigr)$
gives a filtration $\nbigl$ of $\vecnbigt$
in the category $\gbigc_0(\vecX,\ast H)^{\fil}$,
and 
$(\vecnbigt,\nbigl)$
is an object in $\gbigg_1$.
Thus, we obtain the essential surjectivity.

The fully faithfulness of $\Psi$
follows from those for
$ML\nbiga_{D_I}(I)\lrarr
 ML'\nbiga_{D_I}(I)$.
\hfill\qed

\vspace{.1in}
By definition,
we have the forgetful functor
$\gbigg_i(\vecX)
\lrarr
 \gbigg_i(\vecX,\ast H)$
denoted by
$\vecnbigt\longmapsto
 \vecnbigt(\ast H)$.

\subsection{Canonical prolongments}

Let $K=K_1\sqcup K_2$ be a decomposition.
We have a functor
$[\ast K_1!K_2]:
 \gbigg_0(\vecX,\ast H)\lrarr
 \gbigg_0(\vecX)$ given as follows.
\index{functor $[\ast K_1\bikkuri K_2]$}
Let $\vecnbigt=(\nbigt_I,g_{I,i},f_{I,i})
 \in\gbigg_0(\vecX,\ast H)$.
For $I\subset \ellsitabar$, we put
\[
 \nbigttilde_I:=
 \Sigma^{|I\cap K_2|,-|I\cap K_1|}
 \res_{\emptyset}^{I\cap K}\Bigl(
 \psi_{I\cap K}\bigl(\nbigt_{I\setminus K}\bigr)
 \Bigr)
\]
The tuple is equipped with the following morphisms
in $\nbiga_{D_{Ii}}$
\[
\begin{CD}
 \Sigma^{1,0}
 \res^i_{\emptyset}\bigl(
 \psi_i\nbigttilde_I\bigr)
 @>{\psi_{K\cap I}(g_{I\setminus K,i})}>>
 \nbigttilde_{Ii}
 @>{\psi_{K\cap I}(f_{I\setminus K,i})}>>
 \Sigma^{0,-1}
 \res^i_{\emptyset}\bigl(
 \psi_i\nbigttilde_{I}\bigr)
\end{CD}
\quad\quad\quad
(i\not\in K)
\]
\[
 \begin{CD}
 \Sigma^{1,0}
 \res^i_{\emptyset}\bigl(
 \psi_i\nbigttilde_I\bigr)
 @>{\psi_{I\cap K}(N_i)}>> 
 \nbigttilde_{Ii}
 @>{\id}>>
 \Sigma^{0,-1}
 \res^i_{\emptyset}\bigl(
 \psi_i\nbigttilde_I\bigr)
 \end{CD}
\quad\quad\quad
(i\in K_1)
\]
\[
 \begin{CD}
 \Sigma^{1,0}
 \res^i_{\emptyset}\bigl(
 \psi_i\nbigttilde_I\bigr)
 @>{\id}>>
 \nbigttilde_{Ii}
 @>{\psi_{K\cap I}(N_i)}>>
\Sigma^{0,-1}
 \res^i_{\emptyset}\bigl(
 \psi_i\nbigttilde_I\bigr)
 \end{CD}
\quad\quad\quad
(i\in K_2)
\]
We obtain an object
$\vecnbigt[\ast K_1!K_2]$
in $\gbigg_0(\vecX)$,
and the functor
$[\ast K_1!K_2]:
\gbigg_0(\vecX,\ast H)\lrarr
\gbigg_0(\vecX)$.
We have the induced functor
\[
 [\ast K_1!K_2]^{\fil}:
 \gbigg_0(\vecX,\ast H)^{\fil}
\lrarr
 \gbigg_0(\vecX)^{\fil}. 
\]

We have the corresponding functor
$[\ast K_1!K_2]:\gbigg_1(\vecX,\ast H)\lrarr
 \gbigg_1(\vecX)$
given by
\[
\begin{CD}
 \gbigg_1(\vecX,\ast H)
 @>{\Psi}>{\simeq}>
 \gbigg_0(\vecX,\ast H)
 @>{[\ast K_1!K_2]}>>
 \gbigg_0(\vecX)
 @<{\Psi}<{\simeq}<
 \gbigg_1(\vecX)
\end{CD}
\]

Let $(\vecnbigt,\nbigl)\in
\gbigg_1(\vecX,\ast H)$.
We have $\vecnbigt[\ast K_1!K_2]
\in\gbigg_0(\vecX)$.
It is equipped with two naturally 
induced filtrations.
One is the filtration $\nbigltilde$ of
$(\vecnbigt,\nbigl)[\ast K_1!K_2]
\in\gbigg_1(\vecX)$.
The other is the filtration $\nbigl$ of
$\vecnbigt[\ast K_1!K_2]^{\fil}
\in\gbigg_2(\vecX)^{\fil}$.
Note that they are not the same in general.
The latter is called the naively induced filtration.

\vspace{.1in}

We obtain the induced functor
$[\ast K_1!K_2]:
 \gbigg_0(\vecX)\lrarr\gbigg_0(\vecX)$
given by
$\vecnbigt\longmapsto
 \bigl(\vecnbigt(\ast H)\bigr)[\ast K_1!K_2]$,
and the corresponding functor
$[\ast K_1!K_2]:
 \gbigg_1(\vecX)\lrarr\gbigg_1(\vecX)$.

\begin{lem}
We have the natural transformations
$[!K]\lrarr \id \lrarr [\ast K]$
as functors on 
$\gbigg_i(\vecX)$.
\end{lem}
\pf
The claim is clear in the case $i=0$,
which implies the claim in the case $i=1$.
\hfill\qed

\subsection{Nearby cycle, vanishing cycle and maximal functors}

\label{subsection;10.8.30.21}

\index{nearby cycle functor}
\index{vanishing cycle functor}
\index{maximal functor}

Let $g=\vecz^{\vecm}$,
where $\vecm\in\seisuu_{>0}^K$.
Let $\vecnbigt=(\nbigt_I,f_{I,i},g_{I,i})
 \in \gbigg_0(\vecX)$.
We have
$\vecnbigt(\ast H)\in \gbigg_0(\vecX,\ast H)$.
Let 
$\Pi^{a,b}_{\vecm}\vecnbigt(\ast H)
=(\nbigttilde_I,\ftilde_{I,i},\gtilde_{I,i})$
be an object in
$\gbigg_0(\vecX,\ast H)$
given as follows,
for $I\subset \ellsitabar\setminus K$
and $i\in\ellsitabar\setminus (K\cup I)$:
\[
 \nbigttilde_I:=
 \nbigt_I\otimes
 \IItilde^{a,b}_{\vecm},
\quad\quad
\begin{CD}
 \Sigma^{1,0}\psi_i\bigl(\nbigttilde_I\bigr) 
  @>{g_{I,i}\otimes\id}>>
 \nbigttilde_{Ii}
 @>{f_{I,i}\otimes\id}>>
 \Sigma^{0,-1}
 \psi_i\bigl(\nbigttilde_I\bigr)   
\end{CD}
\]
(See \S\ref{subsection;11.4.13.1}
for $\IItilde^{a,b}_{\vecm}$
with $\Lambda=\emptyset$.)
We set
$\Pi^{a,b}_{g\star}\vecnbigt:=
 \bigl(\Pi^{a,b}_g\vecnbigt\bigr)[\star K]
 \in \gbigg_0(\vecX)$ for $\star=\ast,!$.
Then, we define
$\Pi^{a,b}_{g\ast!}\vecnbigt
 \in\gbigg_0(\vecX)$
as the kernel of
$\Pi^{-N,a}_{g!}\vecnbigt
\lrarr
 \Pi^{-N,b}_{g\ast}\vecnbigt$
for a sufficiently large $N$.
In particular, we obtain
the following objects in $\gbigg_0(\vecX)$:
\[
\psi^{(a)}_{g}\vecnbigt:=\Pi^{a,a}_{g\ast!}\vecnbigt,
\quad\quad
\Xi^{(a)}_g\vecnbigt:=
 \Pi^{a,a+1}_{g\ast!}\vecnbigt.
\]
We define $\phi^{(0)}_g\vecnbigt\in\gbigg_0(\vecX)$
as the cohomology of the complex
in $\gbigg_0(\vecX)$:
\[
 \vecnbigt[!K]\lrarr 
 \vecnbigt\oplus\Xi^{(0)}_g\vecnbigt
\lrarr\vecnbigt[\ast K]
\]
We can reconstruct $\vecnbigt\in\gbigg_0(\vecX)$
as the cohomology of the complex
in $\gbigg_0(\vecX)$:
\begin{equation}
\label{eq;10.8.30.11}
 \psi_g^{(1)}\vecnbigt\lrarr
 \phi_g^{(0)}\vecnbigt\oplus \Xi^{(0)}_g\vecnbigt
 \lrarr
 \psi_g^{(0)}\vecnbigt
\end{equation}
The following lemma can be checked
by direct computations.
\begin{lem}
\mbox{{}}\label{lem;10.11.10.1}
\begin{enumerate}
\item
The functors $\psi_g^{(a)}$,
$\Xi_g^{(a)}$ and $\phi_g^{(a)}$
are exact.
\item
If $\vecnbigt[\ast K]=\vecnbigt$,
then 
$\phi_g^{(0)}\vecnbigt
 \simeq\psi_g^{(0)}\vecnbigt$.
If $\vecnbigt[! K]=\vecnbigt$,
then 
$\phi_g^{(0)}\vecnbigt
 \simeq\psi_g^{(1)}\vecnbigt$.
\hfill\qed
\end{enumerate}
\end{lem}

We have the corresponding functors
for $\gbigg_1(\vecX)$,
denoted by the same notation.

\begin{lem}
If $(\vecnbigt,L)\in\gbigg_1(\vecX)$ is pure,
$\psi_g^{(1)}\vecnbigt\lrarr
 \phi_g^{(0)}\vecnbigt
\lrarr\psi^{(0)}_g\vecnbigt$
is $S$-decomposable.
\end{lem}
\pf
We have only to consider the case that
$\vecnbigt$ is obtained as the image of
$\nbigt[!\ellsitabar]\lrarr \nbigt[\ast \ellsitabar]$,
where $\nbigt\in
 \gbigg_0(\vecX,\ast D)$ comes from
a wild variation of polarizable pure 
twistor structure on $X$.
By the claim $1$ in the previous lemma,
the induced morphism
$\phi_g^{(0)}\nbigt[!\ellsitabar]
\lrarr
 \phi_g^{(0)}\vecnbigt$ is surjective,
and 
$\phi_g^{(0)}\vecnbigt\lrarr
 \phi_g^{(0)}\nbigt[\ast \ellsitabar]$
is injective.
Then, the claim follows from
the claim $2$ of Lemma \ref{lem;10.11.10.1}.
\hfill\qed

\vspace{.1in}

Let $(\vecnbigt,L)\in\gbigg_1(\vecX)$.
We have two kinds of filtrations on 
$\psitilde_g(\vecnbigt)$
and $\phi^{(0)}_g(\vecnbigt)$.
One is the naively induced filtration $L$,
when we consider 
as functors on $\gbigg_0(\vecX)^{\fil}$.
The other is the weight filtration $\Ltilde$
as the objects in $\gbigg_1(\vecX)$.
We have the naturally induced morphism
$N:\phi^{(0)}_g(\vecnbigt)\lrarr
 \phi^{(0)}_g(\vecnbigt)\otimes\newTate(-1)$
and 
$N:\psitilde_g(\vecnbigt)\lrarr
 \psitilde_g(\vecnbigt)\otimes\newTate(-1)$.

\begin{prop}
\label{prop;10.8.30.20}
We have $\Ltilde=M(N;L)$
on $\psitilde_g\vecnbigt$
and $\phi_g^{(0)}\vecnbigt$.
\end{prop}
\pf
The claim for $\psitilde_g\vecnbigt$
follows from Lemma \ref{lem;10.10.13.20}.
Let us consider $\phi_g^{(0)}\vecnbigt$.
We have only to consider the case $L$ is pure.
Then, we have the decomposition
$\phi_g^{(0)}\vecnbigt
=\Image \can_g\oplus\Ker \var_g$.
We obtain the relation
$\Ltilde=M(N;L)$ on $\Image\can_g$
from the claim for $\psitilde_g\vecnbigt$.
We have the decomposition
$\vecnbigt=\vecnbigt_1\oplus\Ker\var_g$,
where $\vecnbigt_1$ has no subobject
whose support is contained in $g^{-1}(0)$.
Hence, we have $L=\Ltilde$ on $\Ker\var_g$,
and we obtain the claim for 
$\phi_g^{(0)}\vecnbigt$.
\hfill\qed

\vspace{.1in}

Note that we can reconstruct
the filtration $L$ of $\vecnbigt$
from (\ref{eq;10.8.30.11})
with the naively induced filtrations $L$ of
$\psi^{(a)}_g\vecnbigt$
$(a=0,1)$,
$\Xi^{(0)}_g\vecnbigt$
and $\phi^{(0)}_g\vecnbigt$.
We can also reconstruct $L$ of $\vecnbigt$
from (\ref{eq;10.8.30.11})
with the filtrations $\Ltilde$ of
$\psi^{(a)}_g\vecnbigt$ $(a=0,1)$,
$\Xi^{(0)}_g\vecnbigt$
and $\phi^{(0)}_g\vecnbigt$
as objects in $\gbigg_1(\vecX)$.

\subsection{Gluing along a monomial function}

Let $g$ be as above.
Let $\vecnbigt\in \gbigg_1(\vecX,\ast H)$.
We have 
$\Xi_g^{(0)}\vecnbigt$
and $\psi^{(a)}_g\vecnbigt$
in $\gbigg_1(\vecX)$.
Let $\vecnbigt'\in \gbigg_1(\vecX)$ such that
$\nbigt'_I=0$ if $I\cap K=\emptyset$.
If we are given morphisms
\[
 \psi^{(1)}_g\vecnbigt
 \stackrel{u}{\lrarr}
 \vecnbigt'
 \stackrel{v}{\lrarr}
 \psi^{(0)}_g\vecnbigt
\]
such that $v\circ u$ is equal to 
the natural morphism
$\psi^{(1)}_g\vecnbigt
 \lrarr
 \psi^{(0)}_g\vecnbigt$.
Then, we obtain 
$\Glue(\vecnbigt,\vecnbigt';u,v)
\in \gbigg_1(\vecX)$
as the cohomology of
$\psi^{(1)}_g\vecnbigt\lrarr
 \Xi^{(0)}_g\vecnbigt\oplus \vecnbigt'
\lrarr
 \psi^{(0)}_g\vecnbigt$.
Let $L^{(1)}$ denote the filtration
of $\Glue(\vecnbigt,\vecnbigt';u,v)$
as an object of $\gbigg_1(\vecX)$.

We have the naively induced filtrations $L$
on $\Xi^{(0)}_g\vecnbigt$ 
and $\psi^{(a)}_g\vecnbigt$.
We have the filtration $L$ of
$\vecnbigt'$
obtained as the transfer of 
$L\bigl(\psitilde_g(\vecnbigt)\bigr)$.
Then, we obtain a filtration $L^{(2)}$
of $\Glue(\vecnbigt,\vecnbigt';u,v)$
with which $\Glue(\vecnbigt,\vecnbigt';u,v)\in
 \gbigg_1(\vecX)$.

\begin{lem}
\label{lem;10.8.30.22}
We have $L^{(1)}=L^{(2)}$.
\end{lem}
\pf
We apply the procedure in
\S\ref{subsection;10.8.30.21}
to $\bigl(\Glue(\vecnbigt,\vecnbigt';u,v),L^{(1)}\bigr)$.
The naively induced filtration $L$
on $\phi_g^{(0)}\Glue(\vecnbigt,\vecnbigt';u,v)\simeq\vecnbigt'$
is the same the filtration
obtained as the transfer,
by Proposition \ref{prop;10.8.30.20}.
Then, the claim of the lemma follows.
\hfill\qed

\section{Good pre-mixed twistor $D$-module}

\subsection{Weak admissible specializability}

We introduce an auxiliary notion.

\begin{df}
Let $g$ be any holomorphic function on 
a complex manifold $X$.
An object $(\nbigt,L)\in\MTW(X)$
is called weakly admissibly specializable,
if the following holds:
\begin{itemize}
\item
 $L_j\nbigt$ are strictly specializable along $g$.
\item
 For each $u\in\real\times\cnum$,
 the cokernel of
 $\psitilde_{g,u}(L_j\nbigt)
 \lrarr
 \psitilde_{g,u}(\nbigt)$ is strict,
 and
 $\nbign$ on $\psitilde_{g,u}(\nbigt,L)$
 has a relative monodromy filtration.
\item
 The cokernel of
 $\phi_g(L_j\nbigt)\lrarr \phi_g(\nbigt)$
 is also strict,
and 
 $\nbign$ on 
 $\phi_g(\nbigt,L)$
 has a relative monodromy filtration.
\hfill\qed
\end{itemize}
\end{df}

\subsection{Local case}

Let $X=\Delta^n$ and 
$D=\bigcup_{i=1}^{\ell}\{z_i=0\}$.
We have the category $\gbigg_1(\vecX)$ 
in \S\ref{subsection;10.8.25.12}.
We use a symbol $\vecnbigv$
to denote an object in $\gbigg_1(\vecX)$.
We may regard it as a filtered object in 
the category $\vecC(X,D)$
in \S\ref{subsection;10.9.30.41}.
Then, we obtain a filtered $\nbigr_X$-triple
$\Psi_X(\vecnbigv)$,
as in \S\ref{subsection;11.2.19.20}.

\begin{lem}
$\Psi_X(\vecnbigv)$ is
a pre-mixed twistor $D$-module.
\end{lem}
\pf
We have only to prove 
that, if $\vecnbigv$ is pure of weight $w$,
$\Psi_X(\vecnbigv)$ is a polarizable
wild pure twistor $D$-module of weight $w$.
We may also assume that
$g_{I,i}$ are surjective,
and $f_{I,i}$ are injective.
There exists $J\subset\ellsitabar$
such that 
(i) $\nbigv_I= 0$ unless $I\supset J$,
(ii) $\nbigv_J\neq 0$.
By using the result in
\S\ref{subsection;10.10.1.10},
we obtain that
$\Psi_X(\vecnbigv)$ is strictly $S$-decomposable
along any $z_i$ $(i=1,\ldots,\ell)$.
On the other hand,
we have a wild pure twistor $D$-module
$\gbigt$ of weight $w$
such that
$\gbigt(\ast D(J^c))=\nbigv_{J}$.
Because both 
$\gbigt$ and $\Psi_X(\vecnbigv)$
are strictly $S$-decomposable along $z_i$
$(i=1,\ldots,\ell)$,
we obtain that 
$\gbigt=\Psi_X(\vecnbigv)$.
\hfill\qed

\begin{df}
Let $\MTW^{\good}(\vecX)\subset\MTW(X)$
be the essential image of 
$\Psi_X$.
Any object in $\MTW^{\good}(\vecX)$
is called a good pre-mixed twistor $D$-module
on $\vecX$.
\hfill\qed
\end{df}
\index{good pre-mixed twistor $D$-module}
\index{category $\MTW^{\good}(\vecX)$}
It is independent of the choice of
the coordinate $(z_1,\ldots,z_n)$,
by using 
the isomorphism in Lemma \ref{lem;10.12.25.21}.
(See Lemma \ref{lem;11.1.18.11}.)

By definition,
$\Psi_X$ is essentially surjective.
It is also fully faithful according to 
Lemma \ref{lem;11.4.4.1}.
Hence,
$\Psi_X$ gives an equivalence
$\gbigg_1(\vecX)\lrarr
 \MTW^{\good}(\vecX)$.

\begin{prop}
\label{prop;11.1.21.2}
Let $g=\vecz^{\vecm}$
for some $\vecm\in\seisuu_{>0}^K$
where $K\subset\ellsitabar$.
The following holds:
\begin{itemize}
\item
For $\vecnbigv\in\gbigg_1(\vecX)$,
$\nbigt:=\Psi_X(\vecnbigv)$ is 
weakly admissibly specializable along $g$.
Moreover, there exists
$\nbigt[\star g]$
in $\MTW^{\good}(\vecX)$,
and we have
$\nbigt[\star g]\simeq
 \Psi_X\bigl(\vecnbigv[\star g]\bigr)$
as good prolongment of
$\nbigt(\ast g)$.
\item
We have natural isomorphisms
in $\MTW^{\good}(\vecX)$:
\begin{equation}
 \label{eq;11.4.4.2}
\psitilde_{g,u}(\nbigt)\simeq
 \Psi_X\circ\psitilde_{g,u}(\vecnbigv),
\quad
\Xi_g(\nbigt)
\simeq
 \Psi_X\circ\Xi_g(\vecnbigv),
\quad
\phi_g(\nbigt)
\simeq
 \Psi_X\circ\phi_g(\vecnbigv)
\end{equation}
\item
In particular,
we have $\Ltilde=M(N;L)$
on $\phi_g\nbigt$ and $\psitilde_g\nbigt$,
where $\Ltilde$ denote the filtrations
as objects in $\MTW^{\good}(\vecX)$,
and $L$ denote the naively induced filtrations.
\end{itemize}
\end{prop}
\pf
Let us begin with the following lemma,
forgetting the weight filtrations.
\begin{lem}
\label{lem;11.4.4.4}
(i) $\nbigt=\Psi_X(\vecnbigv)$
 is strictly specializable along $g$ as an $\nbigr_X$-triple,
(ii) $\nbigt[\star g]$ ($\star=\ast,!$) 
exist as $\nbigr_X$-triples,
(iii) we have the isomorphisms {\rm(\ref{eq;11.4.4.2})}
as $\nbigr_X$-triples.
\end{lem}
\pf
Let us consider the case
$\vecnbigv=\nbigv[\ast I!J]$
for $\nbigv\in\MTS^{\adm}(\vecX)$.
We set 
$I_0:=I\setminus K$
and $J_0:=J\setminus K$.
By Proposition \ref{prop;11.1.21.1},
the $\nbigr_X(\ast g)$-triple
$\Psi_X(\nbigv[\ast I!J])(\ast g)$
is strictly specializable along $g$,
and
$\Psi_X(\nbigv[\ast I!J])[\star g]
\simeq
 \Psi_X(\nbigv[\ast I_0!J_0\star K])$
as $\nbigr_X$-triple.
Moreover,
we have
$\Psi_X\bigl(
 \nbigv\otimes\IItilde^{a,b}_g[\ast I!J]\bigr)
 [\star g]
=\Psi_X\bigl(
 \nbigv\otimes\IItilde^{a,b}_g[\ast I_0!J_0\star K]
 \bigr)$
($\star=\ast,!$)
as $\nbigr_X$-triples.
Hence, 
$\Xi^{(0)}_{g}\bigl(
 \Psi_X(\nbigv[\ast I!J])\bigr)$
exists
as an $\nbigr_X$-triple,
and it is given as follows:
\begin{multline}
 \Xi^{(0)}_{g}\bigl(
 \Psi_X(\nbigv[\ast I!J])\bigr)
=  \\
 \Ker\Bigl(
\Psi_X\bigl(
 \nbigv\otimes\IItilde^{-N,1}_g[\ast I_0!J_0!K]
 \bigr)
\lrarr
 \Psi_X\bigl(
\nbigv\otimes\IItilde^{-N,0}_g[\ast I_0!J_0\ast K]
 \bigr)
 \Bigr) \\
\simeq
 \Psi_X\bigl(
 \Xi^{(0)}_g(\nbigv[\ast I!J])
 \bigr).
\end{multline}
Similarly, we have
$\psi^{(a)}_g\bigl(
 \Psi_X(\nbigv[\ast I!J])\bigr)
\simeq
 \Psi_X\psi^{(a)}_g(\nbigv[\ast I!J])$.
Let $\nbigp$ be obtained as the cohomology
of the following:
\[
 \Psi_X\bigl(\nbigv[\ast I_0!J_0!K]\bigr)
\lrarr
 \Psi_X\bigl(
 \Xi^{(0)}_g\nbigv[\ast I!J]\bigr)
 \oplus
 \Psi_X\bigl(
 \nbigv[\ast I!J]\bigr)
\lrarr
 \Psi_X\bigl(
 \nbigv[\ast I_0!J_0\ast K]
 \bigr)
\]
It is naturally isomorphic to
$\Psi_X(\phi_g^{(0)}\nbigv[\ast I!J])$,
which is strict.
Hence, we obtain that
$\Psi_X(\nbigv[\ast I!J])$
is strictly specializable along $g$,
and $\nbigp\simeq
 \phi_g^{(0)}\Psi_X(\nbigv[\ast I!J])
\simeq
 \Psi_X\bigl(\phi^{(0)}_g(\nbigv[\ast I!J])\bigr)$.
We can compute 
$\psitilde_{g,u}\bigl(
 \Psi_X(\nbigv[\ast I!J])\bigr)$
as above,
and we have
$\psitilde_{g,u}\Psi_X(\nbigv[\ast I!J])
\simeq
 \Psi_X\psitilde_{g,u}(\nbigv[\ast I!J])$.

\vspace{.1in}
A general $\vecnbigv\in\gbigg_1(\vecX)$ is 
expressed as the cohomology of a complex
$\vecnbigv^0\lrarr
 \vecnbigv^1\lrarr
 \vecnbigv^2$,
where $\vecnbigv^p=
 \bigoplus_{k_p} \nbigv_{k_p}[\ast I_{k_p}!J_{k_p}]$.
We have
$\Psi_X(\vecnbigv)
=H^1\bigl(
 \Psi_X(\vecnbigv^{\bullet})
 \bigr)$.
We have 
$\psitilde_g\Psi_X(\vecnbigv^i)
\simeq
 \Psi_X\psitilde_g(\vecnbigv^i)$,
and the cohomology of the complex
$\Psi_X\psitilde_g(\vecnbigv^0)\lrarr
 \Psi_X\psitilde_g(\vecnbigv^1)\lrarr
 \Psi_X\psitilde_g(\vecnbigv^2)$
is strict.
Similarly,
we have 
$\phi^{(0)}_g\Psi_X(\vecnbigv^i)
\simeq
 \Psi_X\phi^{(0)}_g(\vecnbigv^i)$,
and the cohomology of the complex
$\Psi_X\phi^{(0)}_g(\vecnbigv^0)\rarr
 \Psi_X\phi^{(0)}_g(\vecnbigv^1)\rarr
 \Psi_X\phi^{(0)}_g(\vecnbigv^2)$
is strict.
Hence, we obtain that
$\Psi_X(\vecnbigv)
=H^1\bigl(
 \Psi_X(\vecnbigv^{\bullet})
 \bigr)$ is strictly specializable along $g$,
and we have natural isomorphisms
$\psitilde_{g,u}\Psi_X(\vecnbigv)
\simeq
 \Psi_X\psitilde_{g,u}(\vecnbigv)$
and
$\phi^{(0)}_{g}\Psi_X(\vecnbigv)
\simeq
 \Psi_X\phi^{(0)}_{g}(\vecnbigv)$.
In particular,
the following morphisms are isomorphisms:
\[
 \can:
 \psi_g^{(1)}\Psi_X(\vecnbigv[! g]),
 \stackrel{\simeq}{\lrarr}
 \phi_g^{(0)}\Psi_X(\vecnbigv[! g])
\]
\[
 \var:
 \phi_g^{(0)}\Psi_X(\vecnbigv[\ast g])
 \stackrel{\simeq}{\lrarr}
 \psi_g^{(0)}\Psi_X(\vecnbigv[\ast g]).
\]
Hence, we obtain that
$\Psi_X(\vecnbigv)[\star g]\simeq
 \Psi_X(\vecnbigv[\star g])$.
Thus, we obtain Lemma \ref{lem;11.4.4.4}.
\hfill\qed

\vspace{.1in}
Let $\vecnbigv_1\lrarr\vecnbigv_2$ be a monomorphism
in $\gbigg_1(\vecX)$.
Then, the induced morphism
$\psitilde_{g,u}\Psi_X(\vecnbigv_1)
\lrarr
 \psitilde_{g,u}\Psi_X(\vecnbigv_2)$
is identified with
$\Psi_X\psitilde_{g,u}(\vecnbigv_1)
\lrarr
 \Psi_X\psitilde_{g,u}(\vecnbigv_2)$,
and hence the cokernel is strict.
Applying this strictness to
the weight filtration of $\vecnbigv$,
we obtain that
$\Psi_X(\vecnbigv)$ is filtered specializable
along $g$.
Then, the isomorphisms
$\psitilde_{g,u}\Psi_X(\vecnbigv)
\simeq
 \Psi_X\psitilde_{g,u}(\vecnbigv)$
and
$\phi_{g}\Psi_X(\vecnbigv)
\simeq
 \Psi_X\phi_{g}(\vecnbigv)$
are compatible with the naively induced filtrations.
Hence, we obtain that 
$\Psi_X(\vecnbigv)$ is admissibly specializable
by Proposition \ref{prop;10.8.30.20}.
We obtain
$\Psi_X(\vecnbigv)[\star g]
\simeq
 \Psi_X(\vecnbigv[\star g])$
in $\MTW^{\good}(X,D)$
from the above results.
\hfill\qed

\subsection{Global case}

Let $X$ be a complex manifold
with a normal crossing hypersurface 
$D=\bigcup_{i\in\Lambda}D_i$.
Let $\vecnbigi$ be a good set of
ramified irregular values on $(X,D)$.
We set $\vecX:=(X,D,\vecnbigi)$.

\begin{df}
Let $\MTW^{\good}(\vecX)$ be the full subcategory
of $\MTW(X)$,
whose objects $\nbigt$ satisfy the following:
\begin{itemize}
\item
For any point $P\in X$,
we take a small coordinate neighbourhood
$(X_P;z_1,\ldots,z_n)$
such that $D_P:=X_P\cap D=
\bigcup_{i=1}^{\ell}\{z_i=0\}$.
We set $\vecX_P:=(X_P,D_P,\vecnbigi_{|X_P})$.
Then,
$\nbigt_{|X_P}\in
 \MTW^{\good}(\vecX_P)$.
\end{itemize}
An object in $\MTW^{\good}(\vecX)$
is called a good pre-mixed twistor $D$-module
on $\vecX$.
\hfill\qed
\end{df}
\index{category $\MTW^{\good}(\vecX)$}
\index{good pre-mixed twistor $D$-module}

\subsubsection{Weakly admissible specializability and localizability}

We obtain the following proposition
from Proposition \ref{prop;11.1.21.2}.
\begin{prop}
\label{prop;10.10.1.30}
Let $\nbigt\in\MTW^{\good}(\vecX)$.
Let $g$ be a holomorphic function on $X$
such that $g^{-1}(0)\subset D$.
\begin{itemize}
\item
$\nbigt$ is weakly admissibly specializable
along $g$.
There exist $\nbigt[\star g]$
in $\MTW^{\good}(\vecX)$
for $\star=\ast,!$.
\item
$\psitilde_{g,u}(\nbigt)$,
$\Xi_g(\nbigt)$ and 
$\phi^{(a)}_g(\nbigt)$
are objects in $\MTW^{\good}(\vecX)$.
\item
Let $\Ltilde$ be the weight filtrations
of $\psitilde_{g,u}\nbigt$
and $\phi^{(0)}_g\nbigt$ 
as objects in $\MTW^{\good}(\vecX)$.
Let $L$ be the filtrations
of $\psitilde_{g,u}\nbigt$
and $\phi^{(0)}_g\nbigt$
naively induced by the filtration of $\nbigt$.
Then, we have $\Ltilde=M(N;L)$.
\hfill\qed
\end{itemize}
\end{prop}

Let $I\subset\Lambda$.
Let $H$ be an effective divisor 
$\sum_{i\in I} m_iD_i$.
\index{localization}
\begin{prop}
\label{prop;13.5.9.101}
Let $\star=\ast$ or $!$.
Let $\nbigt\in \MTW^{\good}(\vecX)$.
\begin{itemize}
\item
 $\nbigt\in\MTW^{\loc}(X,H)$.
\item
$\nbigt[\star H]\in \MTW^{\loc}(X,H)$
satisfies the following:
\begin{itemize}
\item
Let $P$ be any point of $D$.
Let $(X_P,z_1,\ldots,z_{n})$ be a small
 coordinate neighbourhood around $P$
such that
$D=\bigcup_{i=1}^{\ell}\{z_i=0\}$
and 
 $H=\bigcup_{j\in I_P}\{z_i=0\}$.
Then, we have an isomorphism
$\nbigt[\star H]_{|X_P}\simeq
 (\nbigt_{|X_P})[\star I_P]$.
Such $\nbigt[\star H]$ is unique up to canonical isomorphisms.
\end{itemize}
\item
In particular,
$\nbigt[\star H]\in\MTW^{\good}(\vecX)$,
and it depends only on $I$.
\end{itemize}
\end{prop}
\pf
The first claim follows from Proposition \ref{prop;10.10.1.30}.
By Proposition \ref{prop;13.5.9.100},
$\nbigt[\star H]\in\MTW(X,H)$.
By Proposition \ref{prop;10.10.1.30},
$\nbigt[\star H]$ is determined by
the filtered $\nbigr$-triple satisfying the condition
in the second claim.
The third claim immediately follows from the second.
\hfill\qed

\subsubsection{Canonical prolongments
of admissible mixed twistor structure} 
\label{subsection;13.5.10.200}

Let $\nbigv\in\MTS^{\adm}(\vecX)$.
Let $\Lambda=I\sqcup J$ be a decomposition.
We have a filtered $\nbigr_X$-triple 
$\nbigv[\ast I!J]$
satisfying the following condition.
\begin{itemize}
\item
 Let $P$ be any point of $D$.
 Let $(X_P,z_1,\ldots,z_{n})$ be a small
 coordinate neighbourhood around $P$
 such that $D_P:=D\cap X_P=
 \bigcup_{i=1}^{\ell}\{z_i=0\}$.
 We set $\vecX_P:=(X_P,D_P,\vecnbigi_{|X_P})$.
 We have the decomposition
 $\ellsitabar=I_P\sqcup J_P$
 induced by $\Lambda=I\sqcup J$.
 We put $\nbigv_P:=\nbigv_{|X_P}$,
 which induces
 $\nbigv_P[\ast I_P!J_P]$
 in $\gbigg_1(\vecX_P)$.
 Then, we have
 $\nbigv[\ast I!J]_{|X_P}\simeq 
 \Psi_{X_P}(\nbigv_P[\ast I_P!J_P])$.
\end{itemize}
Such $\nbigv[\ast I!J]$ is unique
up to canonical isomorphisms.

\begin{prop}
\label{prop;13.5.10.201}
$\nbigv[\ast I!J]$ is an object 
in $\MTW^{\good}(\vecX)$.
\end{prop}
\pf
By Proposition \ref{prop;13.5.9.101},
we have only to consider the case
$J=\emptyset$.
We have only to prove that
$\Gr^{\Ltilde}_k(\nbigv)$
are polarizable pure twistor $D$-module
of weight $k$.
We regard $D$ as a reduced effective divisor.
Because admissible specializability
and localizability are checked locally,
we obtain that
$\nbigv[\ast D]$ is admissible and 
localizable along $D$
by Proposition \ref{prop;11.1.21.2}.
As explained in \S\ref{subsection;13.5.9.14},
we have the $\nbigr$-triples
$(\Gr^W\psi^{(a)}_D)(\nbigv[\ast D])$
and
$(\Gr^W\phi^{(a)}_D)(\nbigv[\ast D])$.
The natural morphisms
\[
  (\Gr^W\phi^{(a)}_D)(\nbigv[\ast D])
\lrarr
 (\Gr^W\psi^{(a)}_D)(\nbigv[\ast D])
\]
are isomorphisms.
They are equipped with
the naturally induced nilpotent maps $\nbign$.
\begin{lem}
$\bigl(\Gr^W_k\psi^{(a)}_D\bigr)(\nbigv[\ast D])
 \in \MT(X,k)$.
\end{lem}
\pf
We have the canonical splitting
(\ref{eq;13.5.10.100}).
Note that 
$(\Gr^W_k\psi^{(a)}_D)(\nbigv[\ast D])$
and 
$(\Gr^W_k\psi^{(a)}_D)(\Gr^{\Ltilde}\nbigv[\ast D])$
depend only on $\nbigv$
and $\Gr^L(\nbigv)$.
Hence, we have only to prove the claim of the lemma
in the case that $\nbigv$ is pure of weight $0$.
We may replace $\nbigv[\ast D]$
with the polarized pure twistor $D$-module
associated to $\nbigv$.
Then, the claim follows from 
Proposition \ref{prop;13.5.10.102}.
\hfill\qed

\vspace{.1in}

For any $P\in D$,
we take a small coordinate neighbourhood
$(X_P,z_1,\ldots,z_n)$ around $P$
such that 
$D_P:=X_P\cap D=\bigcup_{i=1}^{\ell}\{z_i=0\}$.
Then, the formula (\ref{eq;13.5.9.41}) holds for
$\Gr_k^{\Ltilde}(\nbigv[\ast D])_{|X_P}$.
By varying $P$
and by gluing the decompositions,
we obtain that 
$\Gr_k^{\Ltilde}(\nbigv[\ast D])$
is isomorphic to the direct sum of
the polarizable pure twistor $D$-module
associated to $\Gr^L_k(\nbigv)$,
and a direct summand of 
$\bigl(\Gr^W_k\psi^{(0)}_D\bigr)(\nbigv[\ast D])$.
Hence, we obtain that
$\Gr_k^{\Ltilde}(\nbigv[\ast D])
 \in \MT(X,k)$.
\hfill\qed

\vspace{.1in}

We obtain the following corollary 
from Proposition \ref{prop;11.1.21.2}.
\begin{cor}
\mbox{{}}\label{cor;10.10.2.1}
Let $K\subset\Lambda$
and $H=\sum_{i\in K}m_iD_i$.
Let $\nbigv\in\MTS^{\adm}(\vecX)$.
\begin{itemize}
\item
$\nbigv[\ast I!J]$ is weakly admissibly specializable along $H$.
\item
There exists 
$\bigl(
 \nbigv[\ast I!J]
 \bigr)
 [\star H]$ in $\MTW^{\good}(\vecX)$,
and we have the following isomorphisms
as good prolongment of $\nbigt$:
\[
 \bigl(
 \nbigv[\ast I!J]
 \bigr)
 [\ast H]
\simeq
 \nbigv[\ast(I\cup K)!(J\setminus K)],
\quad\quad
\bigl(\nbigv[\ast I!J]\bigr)[!H]
=\nbigv[\ast (I\setminus K)!(J\cup K)]
\]
\hfill\qed
\end{itemize}
\end{cor}

\subsubsection{Admissible specializability
along any monomial functions}

We consider $\nbigv\in\MTS^{\adm}(\vecX)$.
Let $\Lambda=I\sqcup J$ be any decomposition.
Let $g$ be any holomorphic function on $X$
such that $g^{-1}(0)\subset D$.
\begin{lem}
\label{lem;13.7.31.10}
$\nbigv[\ast I!J]$ is admissibly specializable 
along any $g$.
\end{lem}
\pf
We need to consider 
the ramified exponential twist
by $\gminia\in t^{-1/m}\cnum[t^{-1/m}]$.
Let $\Gal(\varphi_m)$ be the Galois group of
the ramified covering 
$\varphi_m:\cnum_{t_m}\lrarr \cnum_t$.
We set
$\nbigi_{-\gminia}:=\bigl\{
 -\kappa^{\ast}\gminia\,\big|\,
 \kappa\in\Gal(\varphi_m)
 \bigr\}$,
which gives a good set of ramified irregular values on 
$(\cnum_t,0)$.
(See \S\ref{section;13.5.6.100}.)
We obtain a good system of irregular values
$\vecnbigi_{X,-\gminia}:=
 (g)^{-1}\vecnbigi_{-\gminia}$
on $(X,D)$.
For each $P\in D$,
we set
$\nbigi^{(1)}_P:=\bigl\{
 \gminib+\gminic\,\big|\,
 \gminib\in\nbigi_P,\,
 \gminic\in\nbigi_{X,-\gminia,P}
 \bigr\}$.
Then, by applying Proposition \ref{prop;13.5.7.110}
below 
to the system
$\vecnbigi^{(1)}:=\bigl(
 \nbigi^{(1)}_P\,\big|\,P\in D
 \bigr)$,
we obtain a projective morphism of complex manifolds
$F:X'\lrarr X$ such that
(i) $D':=F^{-1}(D)$ is simply normal crossing,
(ii) $X'\setminus (g\circ F)^{-1}(0)
 \simeq X\setminus g^{-1}(0)$,
(iii) $\vecnbigi':=F^{-1}\vecnbigi^{(1)}$ 
 is a good system of irregular values on $(X',D')$.
Let $\vecX':=(X',D',\vecnbigi')$.
Then, we have
$\nbigv'_{-\gminia}:=
 F^{\ast}\bigl(\nbigv
 \otimes g^{\ast}\varphi_m{\dagger}\nbigt_{-\gminia}\bigr)
\in \MTS^{\adm}(\vecX')$.

Let $D':=\bigcup_{i\in\Lambda'}D_i'$
be the irreducible decomposition.
Let $\Lambda'=I'\sqcup J'$ be any decomposition.
Then, $\nbigv'_{-\gminia}[\ast I'!J']$
is weakly admissibly specializable along $g\circ F$,
according to Proposition \ref{prop;10.10.1.30}.
By an argument in Lemma \ref{lem;10.10.2.10},
we obtain that
$F_{\dagger}(\nbigv'_{-\gminia}[\ast I'!J'])$
is admissibly specializable along $g$.
Then, we can easily deduce that
$\nbigv[\ast I!J]
 \otimes
 g^{\ast}\varphi_{m\dagger}\nbigt_{-\gminia}$
is admissibly specializable along $g$.
By applying it to any $\gminia$,
we obtain the admissible specializability of
$\nbigv[\ast I!J]$ along $g$.
\hfill\qed

\subsection{Gluing}

Let $g$ be a holomorphic function on $X$
such that $D_1:=g^{-1}(0)\subset D$.
Let $\nbigt\in \MTW^{\good}(\vecX)$.
We have
$\psi_g^{(a)}(\nbigt)$ and $\Xi^{(a)}_g(\nbigt)$
in $\MTW^{\good}(\vecX)$.
They are equipped with the filtrations $\Ltilde$
as objects in $\MTW^{\good}(\vecX)$.
They are also equipped with
the naively induced filtrations $L$.

Let $\nbigt'\in\MTW^{\good}(\vecX)$
whose support is contained in $D_1$.
Assume that we are given morphisms
\[
 \psi^{(1)}_g\nbigt
 \stackrel{u}{\lrarr}
 \nbigt'
 \stackrel{v}{\lrarr}
 \psi^{(0)}_g\nbigt
\]
in $\MTW^{\good}(\vecX)$ such that $v\circ u$
is equal to the natural morphism
$\psi^{(1)}_g\nbigt\lrarr
 \psi^{(0)}_g\nbigt$.
Then, we obtain 
$\Glue(\nbigt,\nbigt',u,v)$
in $\MTW^{\good}(\vecX)$
as the cohomology of the complex
$\psi^{(1)}_g\nbigt\lrarr
 \Xi^{(0)}_g\nbigt\oplus
 \nbigt'\lrarr\psi^{(0)}_g\nbigt$.
It has the filtration $L^{(1)}$
with which  we have
$\Glue(\nbigt,\nbigt',u,v)
 \in \MTW^{\good}(\vecX)$.

We have the filtration $L$ of
$\nbigt'$ obtained as the transfer
of $L$ of $\psitilde_{g,-\vecdelta}(\nbigt)$.
Then, we obtain a filtration $L^{(2)}$
of $\Glue(\nbigt,\nbigt',u,v)$.
We obtain the following lemma
from Lemma \ref{lem;10.8.30.22}.
\begin{lem}
\label{lem;10.8.31.2}
We have $L^{(1)}=L^{(2)}$.
\hfill\qed
\end{lem}

\section{Good mixed twistor $D$-modules}
\label{subsection;11.4.3.41}

\subsection{Statement}

Let $X$ be a complex manifold.
Let $D$ be a normal crossing hypersurface.
Let $\vecnbigi$ be a good system of
ramified irregular values on $(X,D)$.
We set $\vecX:=(X,D,\vecnbigi)$.
We shall prove the following theorem.
\begin{thm}
\label{thm;10.11.13.11}
A good pre-mixed twistor $D$-module
is a mixed twistor $D$-module.
\end{thm}
After the theorem is proved,
good pre-mixed twistor $D$-module
is also called good mixed twistor $D$-module,
and $\MTW^{\good}(\vecX)$ is also denoted by
$\MTM^{\good}(\vecX)$.
\index{category $\MTM^{\good}(\vecX)$}
\index{good mixed twistor $D$-module}

In particular,
for a given admissible mixed twistor structure on $\vecX$,
we have its canonical prolongation.
Namely, let $\nbigv\in\MTS^{\adm}(\vecX)$,
and let $D=\bigcup_{i\in\Lambda}D_i$ be the irreducible
decomposition.
For a decomposition $\Lambda=I\sqcup J$,
we have the good mixed twistor $D$-module
$\nbigv[\ast I!J]$.

We also prove the following,
which is a special case of Proposition \ref{prop;10.11.15.10}
below.
\begin{prop}
\label{prop;11.1.21.10}
Let $(\nbigt,L)\in\MTM^{\good}(\vecX)$.
For any holomorphic function $f$ on $X$,
we have
$(\nbigt,L)[\star f]$ $(\star=\ast,!)$
in $\MTM(X)$.
Moreover, if $f^{-1}(0)=g^{-1}(0)$
for some holomorphic function $g$,
we have a natural isomorphism
$(\nbigt,L)[\star f]\simeq
 (\nbigt,L)[\star g]$.
\end{prop}

\begin{cor}
Let $H$ be an effective divisor of $X$.
For any object $(\nbigt,L)\in \MTM^{\good}(\vecX)$,
we have
$(\nbigt,L)[\star H]\in \MTM(X)$.
It depends only on the support of $H$.
\end{cor}
\pf
It follows from Proposition \ref{prop;11.1.21.10}
and Proposition \ref{prop;13.5.9.100}.
\hfill\qed

\subsection{Preliminary}
\label{subsection;10.10.2.22}

Let $\nbigv\in\MTS^{\adm}(\vecX)$.
Let $D_0$ be a hypersurface of $X$.
Let $\varphi:X'\lrarr X$ be a proper birational morphism
such that
$D':=\varphi^{-1}\bigl(
 D\cup D_0\bigr)$ is normal crossing.
We set $\vecX':=(X',D',\varphi^{-1}\vecnbigi)$.
We obtain an object
$\nbigv':=\varphi^{\ast}(\nbigv)(\ast D')
\in \MTS^{\adm}(\vecX')$.

\vspace{.1in}

Let $D=D_1\cup D_2$ be a decomposition
into normal crossing hypersurfaces
with $\codim D_1\cap D_2\geq 2$.
We put $D_2':=\varphi^{-1}(D_2)$,
and let 
$D'=D_1'\cup D_2'$ be the decomposition
with $\codim D_1'\cap D_2'\geq 2$.
We put $D_1'':=\varphi^{-1}(D_1)$,
and let $D'=D_1''\cup D_2''$
be the decomposition
with $\codim D_1''\cap D_2''\geq 2$.
Note that $D_1''\subset D_1'$
and $D_2'\subset D_2''$.
We shall prove the following lemma
in \S\ref{subsection;11.1.21.11}.
\begin{lem}
\label{lem;10.10.2.11}
We have natural morphisms in $\MTW(X)$
\[
\varphi_{\dagger}\nbigv'[\ast D_1''!D_2'']
 \lrarr
\nbigv[\ast D_1!D_2] 
 \lrarr
\varphi_{\dagger}\nbigv'[\ast D_1'!D_2']
\]
which induces natural isomorphisms
$\varphi_{\dagger}\nbigv'
\simeq
\nbigv(\ast D_0)
\simeq
\varphi_{\dagger}\nbigv'$.

If one of $D_i$ is empty,
$\nbigv[\ast D_1!D_2]$ is the image of
the morphism
$\varphi_{\dagger}\nbigv'[\ast D_1''!D_2'']
 \lrarr
\varphi_{\dagger}\nbigv'[\ast D_1'!D_2']$
in $\MTW(X)$.
\end{lem}

\subsection{Admissible specializability
of good pre-mixed twistor $D$-modules}

Let $\nbigv\in\MTS^{\adm}(\vecX)$.
We have 
$\nbigv[\star D]\in \MTW^{\good}(\vecX)$
for $\star=\ast,!$.

\begin{lem}
\label{lem;10.10.2.13}
Let $f$ be any holomorphic function on $X$.
Let $\star$ be $\ast$ or $!$.
\begin{itemize}
\item
$\nbigv[\star D]$ 
is admissibly specializable along $f$.
\item
There exist $(\nbigv[\star D])[\ast f]$
and $(\nbigv[\star D])[! f]$
in $\MTW^{\sp}(X,f)$.
\end{itemize}
\end{lem}
\pf
We set $D_0:=f^{-1}(0)$,
and we take $\varphi:X'\lrarr X$
as in \S\ref{subsection;10.10.2.22}.
We set $D_1:=D$ and $D_2:=\emptyset$
if $\star=\ast$,
and $D_1=\emptyset$ and $D_2:=D$
if $\star=!$.
We use the notation in 
\S\ref{subsection;10.10.2.22}.
According to Proposition \ref{prop;10.10.1.30}
and Lemma \ref{lem;13.7.31.10},
$\nbigv'[\ast D'_1!D_2']$
and $\nbigv'[\ast D_1''!D_2'']$
are admissibly specializable
along $\varphi^{\ast}f$.
According to Lemma \ref{lem;10.10.2.10},
we have
$\varphi_{\dagger}\nbigv'[\ast D'_1!D_2']$
and
$\varphi_{\dagger}\nbigv'[\ast D_1''!D_2'']$
are admissibly specializable
along $f$.
By Lemma \ref{lem;10.10.2.11}
and Proposition \ref{prop;10.10.2.12},
we obtain that $\nbigv[\star D]$
is also admissibly specializable along $f$.
It is easy to see that
$\varphi_{\dagger}\bigl(
 \nbigv'[\ast D_1'!D_2'][\ast \varphi^{\ast}f]\bigr)$
gives $(\nbigv[\star D])[\ast f]$,
and that
$\varphi_{\dagger}\bigl(
 \nbigv'[\ast D_1''!D_2'']
 [!\varphi^{\ast}f]\bigr)$
gives $(\nbigv[\star D])[!f]$.
Thus, we obtain Lemma
\ref{lem;10.10.2.13}.
\hfill\qed

\begin{cor}
\label{cor;10.8.31.5}
Let $g$ be a holomorphic function
such that $g^{-1}(0)=D$.
Then, $\psi_g^{(a)}\nbigv[\star D]$
and $\Xi_g^{(a)}\nbigv[\star D]$
are admissibly specializable along 
any holomorphic function $f$ on $X$.
\hfill\qed
\end{cor}

\begin{lem}
\label{lem;10.10.1.40}
$\nbigt\in\MTW^{\good}(\vecX)$
is admissibly specializable along 
any holomorphic function $f$ on $X$.
Moreover,
there exist 
$\nbigt[\star f]$ $(\star=\ast,!)$
in $\MTW^{\sp}(X,f)$.
\end{lem}
\pf
We have only to consider the case
$X=\Delta^n$ and 
$D=\bigcup_{i=1}^{\ell}\{z_i=0\}$.
For $\nbigt\in\MTW^{\good}(\vecX)$,
we put $\rho_1(\nbigt):=\dim\Supp\nbigt$,
and let $\rho_2(\nbigt)$ denote the number
of $I\subset\ellsitabar$ such that
$|I|+\dim\Supp\nbigt=n$
and $\phi_I(\nbigt)\neq 0$.
We set $\rho(\nbigt):=\bigl(\rho_1(\nbigt),
 \rho_2(\nbigt)\bigr)\in
 \seisuu_{\geq 0}\times\seisuu_{\geq 0}$.
We use the lexicographic order on
$\seisuu_{\geq 0}\times\seisuu_{\geq 0}$.

We take $I\subset\ellsitabar$ such that
$|I|+\dim\Supp\nbigt=n$
and $\phi_I(\nbigt)\neq 0$.
We put $h:=\prod_{i\in I^c}z_i$.
Then, $\nbigt$ can be reconstructed
as the cohomology of the complex:
\[
 \psi_h^{(1)}\nbigt\lrarr
 \phi_h^{(0)}\nbigt\oplus\Xi_h^{(0)}\nbigt
 \lrarr \psi_h^{(0)}\nbigt
\]
Because $\rho(\psi_h^{(a)}\nbigt)<\rho(\nbigt)$
and $\rho(\phi_h^{(a)}\nbigt)<\rho(\nbigt)$,
we can apply the hypothesis of the induction.
Applying Corollary \ref{cor;10.8.31.5}
to $\Xi_h^{(0)}\nbigt$,
we obtain
$\nbigt\in\MTW^{\sp}(X,g)$.
Thus, the proof of 
Lemma \ref{lem;10.10.1.40}
is finished.
\hfill\qed

\subsection{Proof of Theorem \ref{thm;10.11.13.11}
and Proposition \ref{prop;11.1.21.10}}
\label{subsection;10.11.13.10}

Let $X$ be a complex manifold.
Let $H$ be a normal crossing hypersurface of $X$.
Let $\vecnbigi$ be a good system of ramified
irregular values on $(X,D)$.
We set $\vecX=(X,H,\vecnbigi)$.
We consider the following.
\begin{description}
\item[${\pmb P(n)}$:]
 The claim of the theorem holds for
 $\nbigt\in\MTW^{\good}(\vecX)$,
 if $\dim\Supp\nbigt\leq n$.
\end{description}
We prove $P(n)$ by an induction on $n$.
Assume that $P(n-1)$ holds.

Let us consider the case that
$X:=\Delta^n$
and $H:=\bigcup_{i=1}^{\ell}\{z_i=0\}$.
We take $\nbigv\in\MTS^{\adm}(\vecX)$,
and let us prove that
$\nbigt=\nbigv[\ast H]$ 
is a mixed twistor $D$-module.
We have already known that
$\nbigt$ is admissibly specializable.
Let $g$ be any holomorphic function on $X$.
We take a projective birational morphism
$\varphi:X'\lrarr X$ such that
$H':=\varphi^{-1}(H\cup g^{-1}(0))$
is normal crossing.
We obtain $\nbigv':=\varphi^{\ast}(\nbigv)$
in $\MTS^{\adm}(\vecX')$.
We put $H_1':=\varphi^{-1}(H)$,
and let $H_2'$ be the complement of $H_1'$
in $H'$.
We obtain 
$\nbigt_1:=\nbigv'[!H_2'\ast H_1']$
and 
$\nbigt_2:=\nbigv'[\ast H']$
in $\MTW^{\good}(\vecX')$.
We put $g':=\varphi^{\ast}(g)$.
By the hypothesis of the induction,
\[
 \phi_{g'}(\nbigt_1),\quad
 \phi_{g'}(\nbigt_2),\quad
  \psitilde_{g',u}(\nbigt_1),\quad
 \psitilde_{g',u}(\nbigt_2)
\]
are mixed twistor $D$-modules.
By Lemma \ref{lem;10.10.2.11},
we have
$\nbigt
\simeq
 \Image\bigl(
 \varphi_{\dagger}(\nbigt_1)\lrarr\varphi_{\dagger}(\nbigt_2)
 \bigr)$.
We obtain 
\[
 \phi_g(\nbigt)
\simeq
 \Image\Bigl(
 \varphi_{\dagger}\phi_{g'}(\nbigt_1)\lrarr
 \varphi_{\dagger}\phi_{g'}(\nbigt_2)
\Bigr),
\]
\[
 \psitilde_{g,u}(\nbigt)
\simeq
 \Image\Bigl(
 \varphi_{\dagger}\psitilde_{g',u}(\nbigt_1)
 \lrarr
 \varphi_{\dagger}\psitilde_{g',u}(\nbigt_2)
\Bigr).
\]
Hence, we obtain
$\phi_g(\nbigt)$ and
$\psitilde_{g,u}(\nbigt)$
are mixed twistor $D$-modules.

For any $\gminia\in t^{-1/m}\cnum[t^{-1/m}]$,
by applying the argument 
in Lemma \ref{lem;13.7.31.10},
we obtain that
$\psitilde_{g',\gminia,u}(\nbigt_i)$
are mixed twistor $D$-modules.
We obtain that
$\psitilde_{g,\gminia,u}(\nbigt)$
are mixed twistor $D$-modules
as above.
Hence, we obtain
$\nbigt=\nbigv[\ast H]$ is a mixed twistor $D$-module.
Similarly,
we obtain that 
$\nbigv[!H]$ is a mixed twistor $D$-module.

\vspace{.1in}

Let us consider the case that
$X=\Delta^N$ for $N\geq n$,
and $H=\bigcup_{i=1}^{\ell}\{z_i=0\}$.
We set $H_{[n]}:=\bigcup_{|I|=N-n} H_I$,
where $H_I=\bigcap_{i\in I}H_i$.
Let $\nbigt\in \MTW^{\good}(\vecX)$
with $\Supp\nbigt\subset H_{[n]}$.
Let us prove that $\nbigt\in\MTM(X)$.
We use an induction on the number $k(\nbigt)$
of the $n$-dimensional irreducible components
of the support of $\nbigt$.
The case $k(\nbigt)=0$ follows from
the hypothesis of the induction.
We take $I\subset\ellsitabar$
with $|I|=N-n$ such that 
$\phi_I\nbigt\neq 0$.
Let $g:=\prod_{\ellsitabar\setminus I}z_i$.
By the hypothesis of the induction,
we have
$\phi^{(0)}_{g}(\nbigt),
\psi^{(a)}_g(\nbigt)\in\MTM(X)$.
By the result in the previous paragraph,
we obtain that $\Xi^{(0)}_g(\nbigt) \in \MTM(X)$.
Because we can reconstruct $\nbigt$ from
$\phi^{(0)}_{g}(\nbigt),
\psi^{(a)}_g(\nbigt)$ $(a=0,1)$
and $\Xi^{(0)}_g(\nbigt)$,
we obtain that $\nbigt$ is 
a mixed twistor $D$-module.
Thus, the proof of Theorem \ref{thm;10.11.13.11}
is finished.

\vspace{.1in}
We can observe that 
$(\nbigt[\star D])[\ast f]$
and $(\nbigt[\star D])[!f]$ in Lemma \ref{lem;10.10.2.13}
are mixed twistor $D$-modules,
by their construction and Theorem \ref{thm;10.11.13.11}.
It follows that $\nbigt[\star f]$ $(\star =\ast,!)$
in Lemma \ref{lem;10.10.1.40}
are mixed twistor $D$-modules.
Thus, the proof of Proposition \ref{prop;11.1.21.10}
is finished.

\subsection{Proof of Lemma \ref{lem;10.10.2.11}}
\label{subsection;11.1.21.11}

Let us return to the situation
in \S\ref{subsection;10.10.2.22}.
Let $\nbigm$ be one of the smooth $\nbigr_X(\ast D)$-modules
underlying $\nbigv$.
Let $\nbigm':=
 \varphi^{\ast}\nbigm\otimes
 \nbigo_{\nbigx}(\ast\nbigd')$.
We put $D_3':=\varphi^{-1}(D_0)$,
and $D_4':=\varphi^{-1}(D)$.
Let $\nbigm'_1:=\bigl(
 \nbigm'[\ast D_3']\bigr)(\ast D_4')$
and $\nbigm'_2:=\bigl(
 \nbigm'[!D_3']\bigr)(\ast D_4')$.

\begin{lem}
\label{lem;10.10.2.20}
$\nbigm$ is naturally isomorphic
to the image of
$\varphi_{\dagger}\nbigm'_2
\lrarr\varphi_{\dagger}\nbigm'_1$.
\end{lem}
\pf
By the construction of $\nbigm'_1$,
we have a naturally defined morphism
$\nbigm\lrarr\varphi_{\dagger}\nbigm'_1$.
For $\lambda\neq 0$,
we have a naturally defined morphisms
$\varphi_{\dagger}\nbigm^{\prime\lambda}_2
\lrarr
 \nbigmlambda$.
Hence, we obtain
$\varphi_{\dagger}\nbigm_2'\lrarr
 \nbigm$.

Let $\lambda_0\neq 0$
and $P\in D_0\setminus D$.
If we take a small neighbourhood
$\nbigu$ around $(\lambda_0,P)$,
we have an isomorphism
of $\nbigm$
with a natural $\nbigr_X$-module $\nbigo_{\nbigx}$
on $\nbigu$.
Hence, it is easy to observe that 
$\nbigm_{|\nbigu}$
is the image of 
$\varphi_{\dagger}\nbigm_2'\lrarr
\varphi_{\dagger}\nbigm_1'$
on $\nbigu$.
We obtain that
$\nbigm_{|\cnum_{\lambda}^{\ast}\times X}$
is the image of
$\varphi_{\dagger}\nbigm_2'\lrarr
\varphi_{\dagger}\nbigm_1'$
on $\cnum_{\lambda}^{\ast}\times X$.
Because the cokernel of
$\varphi_{\dagger}\nbigm_2'\lrarr
\varphi_{\dagger}\nbigm_1'$
is strict,
we obtain that $\nbigm$
is the image of
$\varphi_{\dagger}\nbigm_2'\lrarr
\varphi_{\dagger}\nbigm_1'$.
\hfill\qed

\vspace{.1in}
Let $L$ denote the weight filtrations of $\nbigm$
and $\nbigm'$,
induced by the weight filtration of
the admissible variation of mixed twistor structure.
The naively induced filtrations of $\nbigm_i'$ $(i=1,2)$
are denoted by $L$.
Let $\Ltilde$ denote the weight filtrations
of $\nbigm_i'$ induced by
the weight filtration 
as objects of $\MTW^{\good}(\vecX')$.
They induce the weight filtration of
$\varphi_{\dagger}\nbigm_i'$,
which are also denoted by $\Ltilde$.

We have the naturally induced morphisms
$\varphi_{\dagger}(\nbigm_2',L)
\lrarr
 (\nbigm,L)
\lrarr
 \varphi_{\dagger}(\nbigm_1',L)$.
They are strictly compatible
with the filtrations
by Lemma \ref{lem;10.10.2.20}.

\begin{lem}
The morphisms
$\varphi_{\dagger}\nbigm'_2\lrarr\nbigm
\lrarr\varphi_{\dagger}\nbigm_1'$
give morphisms
\[
 (\varphi_{\dagger}\nbigm'_2,\Ltilde)
\lrarr(\nbigm,L)
\lrarr
 (\varphi_{\dagger}\nbigm'_1,\Ltilde),
\]
which are strictly compatible with the filtrations.
\end{lem}
\pf
We have the induced morphisms
$L_k\varphi_{\dagger}\nbigm'_2
 \lrarr
 L_k\nbigm
 \lrarr
 L_k\varphi_{\dagger}\nbigm'_1$
for each $k$.
We use an induction on $k$
to prove that 
we have natural morphisms
\[
 (L_k\varphi_{\dagger}\nbigm'_2,\Ltilde)
\lrarr(L_k\nbigm,L)
\lrarr
 (L_k\varphi_{\dagger}\nbigm'_1,\Ltilde),
\]
and that they are strictly compatible with 
the filtrations.
If $k$ is sufficiently negative,
the claim is trivial.
Assume the claim in the case $k-1$.
Let us look at 
the induced morphisms
$\Gr^L_k\varphi_{\dagger}\nbigm'_2\lrarr
 \Gr^L_k\nbigm$.
Note that 
$\Gr^L_k\varphi_{\dagger}\nbigm'_2
=\Ltilde_k\Gr^L_k\varphi_{\dagger}\nbigm'_2$,
and the support of
$\Ltilde_{k-1}
\Gr^L_k\varphi_{\dagger}\nbigm'_2$
is contained in $D_0$.
Hence, the induced morphism
$\Ltilde_{k-1}\Gr^L_k\varphi_{\dagger}\nbigm'_2
\lrarr \Gr^L_k\nbigm$ is $0$,
it implies that the image of
$\Ltilde_{k-1}L_k\varphi_{\dagger}\nbigm'_2$
is contained in $L_{k-1}\nbigm$.
Let us look at the following:
\[
\begin{CD}
 \Ltilde_{k-1}L_k\varphi_{\dagger}\nbigm'_2
@>{a}>>
 L_{k-1}\nbigm
@>{b}>>
 L_{k-1}\varphi_{\dagger}\nbigm'_1
\end{CD}
\]
The image of $\Ltilde_{m}L_k\varphi_{\dagger}\nbigm_2'$
$(m<k)$ via $b\circ a$
is contained in
$\Ltilde_mL_k\varphi_{\dagger}\nbigm_1'$.
Because $b$ is injective and
strict by the assumption of the induction,
we obtain that
the image of $\Ltilde_{m}L_k\varphi_{\dagger}\nbigm_2'$
$(m<k)$ via $a$ is contained in
$L_{m}\nbigm$.
Thus, we obtain 
$(L_k\varphi_{\dagger}\nbigm'_2,\Ltilde)
 \lrarr
 (L_k\nbigm,L)$.
Moreover, it is strictly compatible with
the filtrations.

Because 
$\varphi_{\dagger}\nbigm'_2
\lrarr\varphi_{\dagger}\nbigm'_1$
is strictly compatible with $\Ltilde$,
we obtain that
$(\nbigv,L)\lrarr
 (\varphi_{\dagger}\nbigm'_1,\Ltilde)$
is strictly compatible.
Hence, the induction can go on.
\hfill\qed

\vspace{.1in}

Let $\nbigm_3:=
 \varphi_{\dagger}\nbigm'[\ast D_1'!D_2']$
and $\nbigm_4:=
 \varphi_{\dagger}\nbigm'[\ast D_1''!D_2'']$.
Suppose $D_i=g_i^{-1}(0)$
for holomorphic functions $g_i$.
We know that
$\nbigm[\ast D_1]$
and $\nbigm_3(\ast g_2)$
are admissible specializable along $g_1$,
and we have isomorphisms
$\nbigm[\ast D_1][\ast g_1]
\simeq
 \nbigm[\ast g_1]$
and 
$\nbigm_3(\ast g_2)
\simeq
 \nbigm_3(\ast \nbigd)[\ast g_1](\ast g_2)$.
Hence, we have a uniquely determined morphism
$\nbigm[\ast D_1]\lrarr\nbigm_3(\ast g_2)$
induced by $\nbigm\lrarr
 \varphi_{\dagger}\nbigm'_1=
 \nbigm_3(\ast\nbigd)$,
which is compatible with the weight filtrations.
We have that
$\nbigm[\ast D_1!D_2]$
and $\nbigm_3$
are admissible specializable along $g_2$,
and that
$\nbigm[\ast D_1!D_2][!g_2]
\simeq
 \nbigm[\ast D_1!D_2]$
and 
$\nbigm_3[!g_2]
\simeq
 \nbigm_3$.
Hence, we have a uniquely determined morphism
$\nbigm[\ast D_1!D_2]
\lrarr\nbigm_3$.
Similarly, we obtain
$\nbigm_4\lrarr\nbigm[\ast D_1!D_2]$.
By the uniqueness,
we obtain the global case.
Thus, we obtain the first claim of
Lemma \ref{lem;10.10.2.11}.

\vspace{.1in}

Suppose $D_1=\emptyset$.
For any generic $\lambda$,
the specialization $\nbigmlambda[\ast D_1!D_2]$
is the image of
the induced morphism
$(\nbigm_4)^{\lambda}
 \lrarr
(\nbigm_3)^{\lambda}$.
Then, by using the strictness,
we obtain that 
$\nbigm[\star D]\lrarr\nbigm_3$
is injective,
and that 
$\nbigm_4\lrarr\nbigm[\star D]$
is surjective,
i.e.,
$\nbigm[\star D]$ is the image of
$\nbigm_3\lrarr\nbigm[\star D]\lrarr\nbigm_4$.
Thus, the proof of Lemma \ref{lem;10.10.2.11}
is finished.
\hfill\qed

\section{Integrable case}

Let $X$, $D$, $\vecnbigi$ and $\vecX$
be as in \S\ref{subsection;11.4.3.41}.
Let $D=\bigcup_{i\in \Lambda}D$ be the irreducible decomposition.
An integrable mixed twistor $D$-module $(\nbigt,L)$
is called good on $\vecX$,
if the underlying mixed twistor $D$-module is good on $\vecX$.
Let $\MTM^{\integral\good}(\vecX)\subset
 \MTM^{\integral}(\vecX)$
be the full subcategory of good integrable
mixed twistor $D$-modules on $\vecX$.

Let $\nbigv\in\MTS^{\integral\adm}(\vecX)$.
We have the good mixed twistor $D$-module 
$\nbigv[\ast I!J]$ constructed in 
\S\ref{subsection;13.5.10.200}.
\begin{prop}
\label{prop;11.2.23.50}
We naturally have
$\nbigv[\ast I!J]\in\MTM^{\integral\good}(\vecX)$.
\end{prop}
\pf
By Lemma \ref{lem;11.4.5.1} with $z_i$ $(i\in J)$,
we have only to consider the case
$\nbigv[\ast D]$.
By the construction,
$\nbigv[\ast D]$ and its weight filtration
are integrable.
To check that
$\Gr^{\Ltilde}_k(\nbigv[\ast D])
\in \MTint(X,k)$,
we have only to enhance the argument
in the proof of Proposition
\ref{prop;13.5.10.201}
with integrability.
\hfill\qed

\begin{prop}
Let $\nbigt\in\MTM^{\integral\good}(\vecX)$.
Let $g$ be a holomorphic function on $X$
such that $g^{-1}(0)\subset D$.
\begin{itemize}
\item
There exist $\nbigt[\star g]$
in $\MTM^{\integral\good}(\vecX)$
for $\star=\ast,!$.
\item
$\psitilde_{g,\gminia,u}(\nbigt)$,
$\Xi_g(\nbigt)$ and 
$\phi^{(a)}_g(\nbigt)$
are objects in $\MTM^{\integral\good}(\vecX)$.
\hfill\qed
\end{itemize}
\end{prop}

\begin{prop}
Let $H$ be an effective divisor of $X$.
For any object $\nbigt\in\MTM^{\integral\good}(\vecX)$,
we have
$\nbigt[\star H]\in \MTM^{\integral}(X)$.
It depends only on the support of $H$.
\end{prop}
\pf
We only remark that the integrability of the polarization
in Lemma \ref{lem;13.5.10.400}.
\hfill\qed

\chapter{Some basic property}
\label{section;11.1.26.1}

We prove that 
mixed twistor $D$-modules
can be regarded as the gluing of
admissible mixed twistor structures.
It provides us with a nice understanding
on the structure of mixed twistor $D$-modules.
We also consider the localizations
and the twist by smooth mixed twistor structure.

\section{Expression as gluing of admissible mixed twistor structure}
\label{subsection;11.4.3.40}

\subsection{Cell}

Let $X$ be a complex manifold.
Let $P$ be a point of $X$.
An $n$-dimensional cell at $P$ is a tuple
$(Z,U,\varphi,\nbigv)$ as follows:
\index{cell}
\begin{description}
\item[(Cell 1)]
 $\varphi:Z\lrarr X$ is a morphism of
 complex manifolds
 such that $P\in\varphi(Z)$ and $\dim Z=n$.
 We assume that  there exists a neighbourhood of $X_P$
 of $P$ in $X$ such that
 $\varphi:\varphi^{-1}(X_P)\lrarr X_P$ is projective.
\item[(Cell 2)]
 $U\subset Z$ is the complement 
 of a simply normal crossing hypersurface $D_Z$.
 The restriction $\varphi_{|U}$ is an immersion.
 Moreover, there exists a hypersurface $H$ of $X_P$
 such that $\varphi^{-1}(H)=D_Z\cap\varphi^{-1}(X_P)$.
\item[(Cell 3)]
 $\nbigv\in\MTS(Z,D_Z)$
 satisfying {\bf Adm0}.
\end{description}
A function $g$ on $X$ is called a cell function
of a cell $\nbigc=(Z,U,\varphi,\nbigv)$,
if $(\varphi^{\ast}g)^{-1}(0)=D_Z$.
\index{cell function}
If $\nbigv$ is admissible,
$\nbigc$ is called an admissible cell.
If $\nbigv$ is integrable,
$\nbigc$ is called an integrable cell.
\index{admissible cell}

\subsection{Expression as a gluing}

Let $\nbigt\in\MTM(X)$.
Let $P\in\Supp\nbigt$.
We will shrink $X$ around $P$.
We may have a cell $\nbigc=(Z,U,\varphi,\nbigv)$
and a cell function $g$ such that
$\nbigt(\ast g)=\varphi_{\dagger}(\nbigv)$.
We shall prove the following proposition
in \S\ref{subsection;11.2.20.1}.
\begin{prop}
\label{prop;11.2.23.100}
The cell $\nbigc$ is admissible.
If $\nbigt\in\MTMint(X)$,
$\nbigc$ is also integrable.
\end{prop}

By Lemma \ref{lem;10.10.2.10} and
Proposition \ref{prop;11.1.21.10},
we have 
$\varphi_{\dagger}(\nbigv)[\star g]$
in $\MTM(X)$.
According to Lemma \ref{lem;11.1.21.30},
we have the morphisms 
$\varphi_{\dagger}(\nbigv)[!g]
\lrarr 
 \nbigt
\lrarr
 \varphi_{\dagger}(\nbigv)[\ast g]$
in $\MTM(X)$,
and we obtain
$\bigl(\phi_g^{(0)}(\nbigt),\Ltilde^{(1)}\bigr)$
as the cohomology of the following complex in $\MTM(X)$:
\[
 \varphi_{\dagger}(\nbigv)[!g]
\lrarr 
 \nbigt\oplus
 \Xi_g\varphi_{\dagger}(\nbigv)
\lrarr
 \varphi_{\dagger}(\nbigv)[\ast g]
\]
We have the filtration $L\phi_g^{(0)}(\nbigt)$,
naively induced by the weight filtration of $\nbigt$.
It is a filtration of $\bigl(\phi_g^{(0)}(\nbigt),\Ltilde^{(1)}\bigr)$
in the category $\MTM(X)$.

By definition of mixed twistor $D$-module,
we have the relative monodromy filtration
$\Ltilde:=M(\nbign;L\phi^{(0)}_g\nbigt)$,
with which
$\bigl(
 \phi_g^{(0)}(\nbigt),\Ltilde
 \bigr)\in\MTM(X)$.

\begin{lem}
\label{lem;10.11.15.3}
We have $\Ltilde=\Ltilde^{(1)}$.
\end{lem}
\pf
We have only to consider the case
that $\nbigt$ is pure of weight $0$.
Let us consider the morphisms
$\psi^{(1)}_g(\nbigt)
\stackrel{\can}{\lrarr}
 \phi^{(0)}_g(\nbigt)
\stackrel{\var}{\lrarr}
 \psi^{(0)}_g(\nbigt)$
of $\nbigr_X$-triples.
We have the following mixed twistor $D$-modules:
\[
\bigl(\psi^{(a)}_g(\nbigt),\Ltilde\bigr),
\quad
\bigl(\phi^{(0)}_g(\nbigt),\Ltilde\bigr),
\quad
\bigl(\phi^{(0)}_g(\nbigt),\Ltilde^{(1)}\bigr)
\]
The morphisms $\can$ and $\var$
are morphisms in $\MTM(X)$,
for both the filtrations $\Ltilde$ and $\Ltilde^{(1)}$
for $\phi^{(0)}_g(\nbigt)$.
We have the decomposition
$\phi^{(0)}_g(\nbigt)
=\Image\can\oplus\Ker\var$,
which is compatible with
both $\Ltilde$ and $\Ltilde^{(1)}$.
The restriction of $\Ltilde$ to $\Ker\var$
is pure of weight $0$ by definition.
Because $\nbigt=\nbigt_0\oplus\Ker\var$
as pure twistor $D$-module,
$\Ltilde^{(1)}$ is also pure of weight $0$
on $\Ker\var$.
Because $\can$ and $\var$
are strict with respect to 
both $\Ltilde$ and $\Ltilde^{(1)}$,
we have
$\Ltilde=\Ltilde^{(1)}$ on 
$\Image\can$.
Thus, we obtain Lemma \ref{lem;10.11.15.3}.
\hfill\qed

\vspace{.1in}

We obtain the local expression of $\nbigt$
as the cohomology of the following complex
in $\MTM(X)$:
\begin{equation}
\label{eq;10.11.15.10}
 \psi^{(1)}_g(\varphi_{\dagger}\nbigv)
\lrarr
 \Xi_g^{(0)}(\varphi_{\dagger}\nbigv)
\oplus
 \phi_g^{(0)}(\nbigt)
\lrarr
 \psi^{(0)}_g(\varphi_{\dagger}\nbigv)
\end{equation}

As a consequence,
we obtain the following proposition.
\begin{prop}
\label{prop;13.7.29.30}
Let $(\nbigt,L)\in \MTM(X)$.
For any $P\in X$,
there exist a neighbourhood $X_P$,
an admissible cell
$\nbigc=(Z,U,\varphi,\nbigv)$
with a cell function $g$
for $(\nbigt,L)_{|X_P}$,
and a mixed twistor $D$-module
$(\nbigt_0,L)$
with
\[
 \Supp(\nbigt_0)\subset 
 g^{-1}(0)\cap\Supp(\nbigt),
\]
such that 
$(\nbigt,L)$
is a sub-quotient of
$(\nbigt_0,L)
\oplus
 \varphi_{\dagger}(\nbigv)[!g]$
in $\MTM(X)$.
\hfill\qed
\end{prop}

\begin{rem}
It is M. Saito who emphasized the significance
and the convenience of a description 
of mixed Hodge modules
as in Proposition {\rm\ref{prop;13.7.29.30}}.
\hfill\qed
\end{rem}

\subsection{Gluing}

Let $\nbigc=(Z,U,\varphi,\nbigv)$ be an admissible cell
with a cell function $g$ at $P$.
We will shrink $X$ around $P$.
Let $\nbigt'\in\MTM(X)$ such that
$\Supp\nbigt'\subset g^{-1}(0)$.
Assume that we are given morphisms
\begin{equation}
\label{eq;10.11.15.11}
 \psi_g^{(1)}\varphi_{\dagger}(\nbigv)
 \stackrel{u}{\lrarr}
 \nbigt'
 \stackrel{v}{\lrarr}
 \psi_g^{(0)}\varphi_{\dagger}\nbigv
\end{equation}
in $\MTM(X)$, 
such that $v\circ u$ is equal to
the canonical morphism
$\psi_g^{(1)}\varphi_{\dagger}(\nbigv)\lrarr
 \psi_g^{(0)}\varphi_{\dagger}\nbigv$.
Then, we obtain the $\nbigr$-triple
$(\Glue(\nbigc,\nbigt',u,v),L^{(1)})\in \MTM(X)$
as the cohomology of the following complex
in $\MTM(X)$:
\begin{equation}
\label{eq;10.11.15.12}
 \psi_g^{(1)}\varphi_{\dagger}\nbigv
\lrarr
 \Xi_g^{(0)}\varphi_{\dagger}\nbigv
\oplus
 \nbigt'
\lrarr
 \psi_g^{(0)}\varphi_{\dagger}\nbigv
\end{equation}
We can reconstruct $\nbigt'$
as the cohomology of the complex
in $\MTM(X)$:
\begin{equation}
 \varphi_{\dagger}\nbigv[!g]
\lrarr
 \Xi_g\varphi_{\dagger}\nbigv
 \oplus \Glue(\nbigc,\nbigt',u,v)
\lrarr
 \varphi_{\dagger}\nbigv[\ast g]
\end{equation}
We have a similar expression for
an integrable mixed twistor $D$-module
as the gluing of an integrable admissible cell
and the vanishing cycle integrable mixed twistor $D$-module.

\vspace{.1in}

Let $L\psi_g^{(a)}\varphi_{\dagger}(\nbigv)$
and $L\Xi_g^{(0)}\varphi_{\dagger}(\nbigv)$
be the filtrations naively induced by 
the weight filtration of $\nbigv$.
Let $L$ be the filtration of $\nbigt'$
obtained as the transfer of
$L$ of $\psi_{g}^{(a)}\varphi_{\dagger}(\nbigv)$
with respect to (\ref{eq;10.11.15.11}).
We obtain a filtration $L^{(2)}$ of
$\Glue(\nbigc,\nbigt',u,v)$
from (\ref{eq;10.11.15.12}) with 
the filtrations $L$.

\begin{lem}
We have $L^{(1)}=L^{(2)}$
on $\Glue(\nbigc,\nbigt',u,v)$.
\end{lem}
\pf
We have the filtration $L^{(1)}$
of $\nbigt'$
naively induced by the filtrations 
$L^{(1)}$ on $\Glue(\nbigc,\nbigt',u,v)$.
By using Lemma \ref{lem;10.11.15.3},
we obtain $L^{(1)}=L$ on $\nbigt'$.
Then, we obtain $L^{(1)}=L^{(2)}$ on
$\Glue(\nbigc,\nbigt',u,v)$.
\hfill\qed

\subsection{Admissibility of cells}
\label{subsection;11.2.20.1}

Let $\nbigt$ be a mixed twistor $D$-module
on $X$.
The weight filtration is denoted by $L$.
Let $Z$ be the support of $\nbigt$.
We have the graded pure twistor $D$-module
$\Gr^L(\nbigt)=\bigoplus\Gr^L_w(\nbigt)$,
whose support is $Z$.
Let $Z_0$ be an irreducible component of $Z$.
Let $P\in Z_0$.
We study the local property of $\nbigt$
around $P$.
We shall shrink $X$ around $P$.
Then, we may have a hypersurface $H$ of $X$
such that 
(i) $Z\setminus H=Z_0\setminus H\neq\emptyset$
 is smooth,
(ii) $\nbigt_{|X\setminus H}$ comes from
 a mixed twistor structure
 on $Z_0\setminus H$.
(We use Proposition \ref{prop;10.11.12.30}
and Lemma \ref{lem;10.11.12.31}.)
If $\nbigt\in\MTM^{\integral}(X)$,
the mixed twistor structure is integrable.
We shall prove the following proposition,
which implies Proposition \ref{prop;11.2.23.100}.

\begin{prop}
\label{prop;10.11.13.2}
There exists a projective birational morphism
$\varphi:Z_1\lrarr Z_0$ with the following property:
\begin{itemize}
\item
 $H_1:=\varphi^{-1}(H)$ is simply normal crossing.
\item
 There exists an admissible mixed twistor structure
 $\nbigv$ on $(Z_1,H_1)$
 such that $\nbigt(\ast H)\simeq\varphi_{\dagger}\nbigv$.
\end{itemize}
If $\nbigt\in\MTMint(X)$,
the admissible mixed twistor structure is integrable.
\end{prop}
\pf
Let us consider the case $X=\Delta$
and $H:=\{O\}$.
Let $\nbigt$ be a mixed twistor $D$-module
such that $\nbigt_{|X\setminus H}$
is a mixed twistor structure.
Let $\nbigm_i$ $(i=1,2)$ be the underlying 
$\nbigr_X$-modules.
Because $\Gr^L_w\nbigt(\ast H)$
is obtained as the canonical meromorphic prolongation of
a wild variation of polarizable pure twistor structure of weight $w$,
$\Gr^L_w(\nbigm_i)(\ast H)$ are good $\nbigr_{X(\ast H)}$-modules.
Hence, $\nbigm_i$ are also good $\nbigr_{X(\ast H)}$-modules.
If we take an appropriate ramified covering
$\varphi:(X_2,H_2)\lrarr (X,H)$,
we have the irregular decomposition
$\varphi^{\ast}(\nbigm_i,L)_{|\Hhat_2}
=\bigoplus (\nbigmhat_{i,\gminia},L)$.
Because $(\nbigm_i,L)$ are
admissibly specializable,
each $\nbigm_i$ has $KMS$-structure
compatible with the filtration $L$,
and the condition {\bf(Adm2)}
in \S\ref{subsection;10.9.29.1}
is satisfied.
Hence, $\nbigt(\ast H)$ comes from
an admissible mixed twistor structure.
If $\nbigt$ is integrable,
then the admissible mixed twistor structure 
is integrable.

\vspace{.1in}

Let us consider the general case.
We take a projective birational morphism
$\varphi:X'\lrarr X$ such that 
(i) the proper transform $Z_0'$ of $Z_0$ is smooth,
(ii) $Z_0'$ intersects with $H':=\varphi^{-1}(H)$
in a normal crossing way.
We set $H'_0:=Z_0'\cap H'$.
Let $L$ be the weight filtration of $\nbigt$.
We have the filtered $\nbigr_{X'(\ast H')}$-triple
$(\nbigt',L)$ obtained as the lift of
$(\nbigt,L)(\ast H)$,
as in Lemma \ref{lem;11.2.21.1}.
(We use Lemma \ref{lem;11.2.22.1}
in the integrable case.)
We have the graded $\nbigr_{X'(\ast H')}$-triple
$\Gr^L(\nbigt')$.

\begin{lem}
\label{lem;13.5.13.40}
There exists a projective birational morphism
$\varphi_1:Z_1\lrarr Z_0'$ such that
 $H_1:=\varphi_1^{-1}(H'_0)$ is normal crossing,
and 
there exist
filtered $\nbigr_{Z_1}(\ast H_{1})$-triples
 $\nbigt^{(1)}_w$ $(w\in\seisuu)$ satisfying the following:
\begin{itemize}
\item
 $\nbigt^{(1)}_w$ is obtained as 
 the canonical prolongment of a graded
 good wild variation of pure twistor structure of weight $w$
 on  $(Z_1,H_{1})$.
\item
 $\bigoplus\nbigt^{(1)}_w$ is a good smooth 
 $\nbigr_{Z_1}(\ast H_1)$-triple.
\item
 $\Gr_w^L(\nbigt')\simeq
 (\iota_{Z_1'}\circ\varphi_1)_{\dagger}\nbigt^{(1)}_w$,
where $\iota_{Z_1'}:Z_1'\lrarr X'$ denotes
the inclusion.
\end{itemize}
\end{lem}
\pf
Let $\Gr^L(\nbigt)_0$ be the sum of 
the direct summands of $\Gr^L(\nbigt)$
whose strict supports are $Z_0$.
According to \cite{mochi7},
there exists a projective birational morphism
$\varphi_2:Z_2\lrarr Z_0$ such that
the following holds:
\begin{itemize}
\item 
 There exist a normal crossing hypersurface
 $H_{2}\subset Z_2$ satisfying
 $H_2\supset\varphi_2^{-1}(H_0)$,
 and good wild polarizable variations of 
 pure twistor structure $\nbigv^{(2)}_w$ of weight $w$
 on $(Z_2,H_{2})$.
 Let $\gbigt^{(2)}_{w}$ be 
 the polarizable pure twistor $D$-module of weight $w$
 on $Z_2$ obtained as the minimal extension of 
 $\nbigv^{(2)}_w$.
\item
 $\Gr^L_w(\nbigt)(\ast H)=
 (\iota_{Z_0}\circ\varphi_2)^0_{\dagger}\gbigt^{(2)}_{w})(\ast H)$,
 where $\iota_{Z_0}:Z_0\lrarr X$ denotes the inclusion.
\end{itemize}
We may assume that
$\varphi_2$ factors through $Z_0'$.

We have the polarizable variation of 
pure twistor structure $\nbigv_w^{(0)\prime}$
of weight $w$ on $Z_0'\setminus H_0'$
corresponding to $\Gr^L_w(\nbigt')_{|X'\setminus H'}$.
Applying Corollary \ref{cor;13.5.13.31} below
with the existence of
$\varphi_2$ and $\nbigv_w^{(2)}$ as above,
we can take a projective birational morphism
$\varphi_1:Z_1\lrarr Z_0'$
such that 
(i) $H_{1}:=\varphi_1^{-1}(H_0')$ is normal crossing,
(ii) each 
 $\nbigv^{(1)}_w:=
 \varphi_1^{-1}\nbigv_w^{(0)\prime}$
comes from a good wild harmonic bundle 
$(E^{(1)}_w,\delbar_{E^{(1)}_w},\theta^{(1)}_w,h^{(1)}_w)$
on $(Z_1,H_{1})$ up to shift of the weight,
(iii) $\bigoplus 
 (E^{(1)}_w,\delbar_{E^{(1)}_w},\theta^{(1)}_w,h^{(1)}_w)$
is good wild.
We may assume that
$\varphi_2$ factors through $Z_1$,
i.e.,
we may have $\psi:Z_2\lrarr Z_1$ 
such that $\varphi_2=\varphi_1\circ\psi$.
Let $\nbigt_w^{(1)}$ is the $\nbigr_{Z_1}(\ast H_1)$-triple
obtained as the canonical prolongment of $\nbigv_w^{(1)}$.
Let $\nbigt^{(2)}_w$ be the $\nbigr_{Z_2}(\ast H_{2})$-triple
obtained as the canonical prolongment of
$\nbigv^{(2)}_w$.
We have a natural isomorphism
$\psi^{\ast}\nbigt_w^{(1)}(\ast H_2)
\simeq
 \nbigt_w^{(2)}$ by the construction.
We obtain
$\psi^{\ast}\nbigt_w^{(1)}
\simeq
 \gbigt_w^{(2)}(\ast \psi^{-1}H_1)$.
Hence, we obtain a morphism
$(\iota_{Z_0'}\varphi_1)_{\dagger}\nbigt_w^{(1)}
\lrarr
 \Gr^L_w(\nbigt')(\ast H_0')$.
Because its restriction to
$X'\setminus H_0'$ is an isomorphism
by construction,
it is an isomorphism on $X'$.
Thus, we obtain Lemma \ref{lem;13.5.13.40}.
\hfill\qed

\begin{lem}
\label{lem;13.5.13.50}
There exists a filtered 
$\nbigr_{Z_0'(\ast H_0')}$-triple
$(\nbigt_0',L)$ such that 
$(\nbigt',L)$
is the push-forward of
$(\nbigt_0',L)$
by $(Z_0',H_0')\lrarr (X',H')$.
\end{lem}
\pf
We have only to prove the underlying
filtered $\nbigr_{X'(\ast H')}$-modules $\nbigm_i$
are the push-forward of
$\nbigr_{Z_0'(\ast H_0')}$-modules.
By Lemma \ref{lem;13.5.13.40},
$\Gr^L(\nbigm_i)$ are the push-forward of
graded $\nbigr_{Z_0'(\ast H_0')}$-modules 
$\bigoplus_i\nbigk_{iw}$.
By the admissible specializability of mixed twistor $D$-modules,
the restriction $(\nbigm_i,L)_{|X'\setminus H'}$
are the push-forward of smooth filtered
$\nbigr_{Z_0'\setminus H_0'}$-modules
$(\nbigk_i,L)$.
We have
$\Gr^L(\nbigk_i)=
\bigoplus_i\nbigk_{iw|Z_0'\setminus H_0'}$.
We have only to prove that
$(\nbigk_i,L)$ are uniquely extended to
filtered $\nbigr_{Z_0'(\ast H_0')}$-modules
$(\nbigktilde_i,L)$
such that
$\Gr^L(\nbigktilde_i)\simeq
 \bigoplus\nbigk_{iw}$.
The uniqueness is clear.
So, we have only to construct it locally
around any point $P$ of $H_0'$.
Then, by applying
Corollary \ref{cor;10.11.12.41} below,
we obtain Lemma \ref{lem;13.5.13.50}.
\hfill\qed

\vspace{.1in}

By using Lemma \ref{lem;11.2.21.1},
we obtain a filtered $\nbigr_{Z_1(\ast H_{1})}$-triple
$(\nbigt^{(1)},L)$ such that
(i) $\Gr^L(\nbigt^{(1)})=\bigoplus\nbigt^{(1)}_w$,
(ii) $\varphi_{1\dagger}(\nbigt^{(1)},L)
 \simeq(\nbigt'_0,L)$.
Because the underlying $\nbigr_{Z_1(\ast H_{1})}$-modules
are good and smooth,
$\nbigt^{(1)}$ is 
a good smooth $\nbigr_{Z_1(\ast H_{1})}$-triple
by Proposition \ref{prop;11.2.20.4}.
Then, we obtain that
$(\nbigt^{(1)},L)$ is an object in $\MTS(Z_1,H_{1})$.
The condition {\bf (Adm0)} is satisfied.
Proposition \ref{prop;10.11.13.2}
is reduced to the following.
\begin{prop}
$(\nbigt^{(1)},L)$ is admissible.
\end{prop}
\pf
First, let us consider the case that $\dim Z_0=1$.
We may assume $X=\Delta^n$
and $H=\{z_1=0\}$.
Let $q_1:X\lrarr\Delta^1$ be the projection
onto the first component.
We may assume that the induced morphism
$F:Z_1\lrarr \Delta^1$ is a covering
with a ramification along $0$.
We obtain $\nbigt_1:=
 q_{1\dagger}(\nbigt)\in \MTM(\Delta^1)$.
By the result in 
the one dimensional case,
$\nbigt_2:=\nbigt_1(\ast 0)$
is an object in $\MTS^{\adm}(\Delta^1,0)$.
We obtain 
$F^{\ast}\nbigt_2\in\MTS^{\adm}(Z_1,H_{1})$.
Because
$\nbigt_2=F_{\dagger}(\nbigt^{(1)})$,
$\nbigt^{(1)}$ is a direct summand of
$F^{\ast}\nbigt_2$
as a filtered smooth $\nbigr_{Z_1(\ast H_{1})}$-triple.
Hence, $\nbigv_0\in\MTS^{\adm}(Z_1,H_{1})$.

\vspace{.1in}

Let us consider the general case.
Let $C$ be a smooth curve in
$Z_1$, which intersects with
the smooth part of $H_1$ transversally.
By using the nearby cycle functor 
and the vanishing cycle functor to $\nbigt$ 
successively,
we obtain a mixed twistor $D$-module
$\nbigt_3$ such that 
(i) the support of $\nbigt_3(\ast H)$
is $\varphi_1(C)$,
(ii) the lift of $\nbigt_3(\ast H)$ to $C$ is $\nbigt^{(1)}_{|C}$.
Then, by using the result in the case $\dim Z_0=1$,
we obtain that $\nbigt^{(1)}_{|C}$ is admissible.
By Proposition \ref{prop;13.5.14.11},
$\nbigt^{(1)}$ satisfies {\bf(Adm1)}.
It also satisfies {\bf(Adm2)}.
Then, we obtain the admissibility of $(\nbigv,L)$.
\hfill\qed

\subsubsection{Strict specializability}

We give a proof of a claim
used in the proof of Proposition \ref{prop;10.11.13.2}.
Let $X:=\Delta^n$ 
and $D=\bigcup_{i=1}^{\ell}\{z_i=0\}$.
Let $Z=\bigcap_{i=\ell+1}^m\{z_i=0\}$
and $D_Z:=D\cap Z$.
Let $(\nbigm,L)$ be
a filtered $\nbigr_X(\ast D)$-module
satisfying the following conditions:
\begin{itemize}
\item
 $(\nbigm,L)_{|
 \cnum_{\lambda}\times (X\setminus D)}$
 is the push-forward of
 a smooth $\nbigr_{Z\setminus D_Z}$-module.
\item
 Each $\Gr^L_w(\nbigm)$ is the push-forward of
 a strict coherent $\nbigr_{Z(\ast D_Z)}$-module.
\end{itemize}
Let $X_1:=\{z_{\ell+1}=0\}$
and $D_1:=X_1\cap D$.
Let $i:X_1\lrarr X$ be the inclusion.
Let $f:\nbigm\lrarr\nbigm$
be given by the multiplication of $z_{\ell+1}$.
We obtain
an $\nbigr_{X_1(\ast D_1)}$-module $\Ker f$,
which is naturally equipped with a filtration $L$.

\begin{lem}
We have a natural isomorphism
$\iota_{\dagger}(\Ker f,L)
\simeq (\nbigm,L)$.
We also have
$\Gr^L(\Ker f)\simeq \Ker\Gr^L(f)$.
\end{lem}
\pf
If $L$ is pure, the claim is trivial.
We use an induction on the length of $L$.
Assume that
$L_{-1}\nbigm=0$
and $L_0\nbigm\neq 0$.
By the assumption,
$L_0\nbigm$ comes from 
a strict coherent $\nbigr_{Z(\ast D_Z)}$-module.
By a direct computation,
we can check that
the cokernel $\Cok (f_{|L_0\nbigm})$
is naturally isomorphic to
$L_0\nbigm\bigl/\lambda L_0\nbigm$.
Hence, a section of $g\in\Cok\bigl(f_{|L_0\nbigm}\bigr)$
is $0$, if its restriction to $X\setminus D$ is $0$.
We have the following commutative diagram:
\[
 \begin{CD}
 L_{0}\nbigm @>>>
 \nbigm @>>>
 \nbigm/L_0\nbigm \\
 @V{f}VV @V{f}VV @V{f}VV \\
 L_{0}\nbigm @>>>
 \nbigm @>>>
 \nbigm/L_0\nbigm
 \end{CD}
\]
We may apply the assumption of the induction
to $\nbigm/L_0\nbigm$.
Then, we have only to show that
the induced morphism
\begin{equation}
 \label{eq;10.11.12.20}
 \Ker\bigl(
 f: \nbigm/L_0\nbigm\lrarr
 \nbigm/L_0\nbigm
\bigr)
\lrarr
 \Cok\bigl(
 f:L_0\nbigm\lrarr L_0\nbigm
 \bigr)
\end{equation}
is $0$.
By the assumption,
the restriction of (\ref{eq;10.11.12.20})
to $X\setminus D$ vanish.
Hence, (\ref{eq;10.11.12.20})
vanishes because of the previous consideration.
\hfill\qed

\begin{cor}
\label{cor;10.11.12.41}
Under the assumption,
$(\nbigm,L)$ comes from
a filtered strict coherent $\nbigr_{Z(\ast D_Z)}$-module.
\hfill\qed
\end{cor}

\section{Localization}
\label{subsection;11.4.9.12}

\subsection{Localization along functions}

Let $\nbigt\in\MTM(X)$.
Let $f$ be a holomorphic function on $X$.
Let $\iota_f:X\lrarr X\times\cnum_t$ be the graph.
Recall that we have constructed 
$(\iota_{f\dagger}\nbigt)[\star t]
\in\MTW(X\times\cnum_t)$
in \S\ref{subsection;11.2.22.100}.

\index{localization}

\begin{prop}
\mbox{{}}\label{prop;10.11.15.10}
\begin{itemize}
\item
For $\star=\ast,!$,
we have $\nbigt[\star f]\in\MTM(X)$
such that
$\iota_{f\dagger}(\nbigt[\star f])
\simeq
 (\iota_{f\dagger}\nbigt)[\star t]$.
If $\nbigt\in\MTMint(X)$,
we naturally have
$\nbigt[\star f]\in\MTMint(X)$.
\item
If $f_1^{-1}(0)=f_2^{-1}(0)$,
we have 
$\nbigt[\star f_1]\simeq\nbigt[\star f_2]$
naturally.
\end{itemize}
\end{prop}
\pf
By the uniqueness,
we have only to consider the issues locally.
Let $P\in \Supp\nbigt$.
We use the Noetherian induction on the support
around $P$.
We will shrink $X$ around $P$
in the following argument.
Let $(Z,U,\varphi,\nbigv)$ be a cell of $\nbigt$
with a cell function $g$ at $P$.
We have the expression of $\nbigt$
as the cohomology of the following
complex in $\MTM(X)$:
\[
 \psi_g^{(1)}\varphi_{\dagger}(\nbigv)
\lrarr
 \Xi_g^{(0)}\varphi_{\dagger}(\nbigv)
\oplus
 \phi_g^{(0)}(\nbigt)
\lrarr
 \psi_g^{(0)}\varphi_{\dagger}(\nbigv)
\]
We have already known the claims for
good mixed twistor $D$-modules.
Hence, we obtain the claims for
$\psi_g^{(a)}\varphi_{\dagger}(\nbigv)$
and 
$\Xi_g^{(0)}\varphi_{\dagger}(\nbigv)$.
We can apply the hypothesis of the induction
to $\phi_g^{(0)}(\nbigt)$.
Then,  we obtain the claims for $\nbigt$.
The claim in the integrable case
follows from 
Lemma \ref{lem;11.2.22.101}.
\hfill\qed

\vspace{.1in}
We have naturally defined morphisms
$\nbigt[!f]\lrarr\nbigt\lrarr\nbigt[\ast f]$,
as remarked in Lemma \ref{lem;11.1.21.30}.

\subsection{Localization along hypersurfaces}

Let $H$ be a hypersurface of $X$.
Let $\nbigt\in\MTM(X)$.
By Proposition \ref{prop;10.11.15.10}
and Proposition \ref{prop;13.5.9.100},
we have $\nbigt[\star H]$ in $\MTM(X)$
with the following property.
\begin{itemize}
\item
Let $P$ be any point of $X$.
Let $X_P$ be a small neighbourhood around $P$.
We have an expression $H=\{f=0\}$ on $X_P$.
Then, $\nbigt[\star H]_{|X_P}=\nbigt_{|X_P}[\star f]$.
\end{itemize}
We obtain the full subcategory $\MTM(X,[\star H])$
of $(\nbigt,L)\in \MTM(X)$
such that $\nbigt=\nbigt[\star H]$
in $\MTM(X)$.
\index{category $\MTM(X,[\star H])$}
By Proposition \ref{prop;10.11.15.10},
we have a naturally defined functor:
\[
 \MTM(X)\lrarr \MTM(X,[\star H])
\]
We have naturally defined morphisms
$\nbigt[!H]\lrarr\nbigt\lrarr\nbigt[\ast H]$.
If $\nbigt$ is integrable,
$\nbigt[\star H]$ are also integrable.

We obtain the following from Lemma \ref{lem;11.2.20.5}.
\begin{lem}
\label{lem;11.2.21.2}
Let $\star=\ast$ or $!$.
Let 
$\nbigt_i\in \MTM(X)$ $(i=1,2)$
such that
$\nbigt_i[\star H]=\nbigt_i$.
We have a natural bijective correspondence
between morphisms
$\nbigt_1(\ast H)\lrarr\nbigt_2(\ast H)$
as filtered $\nbigr_X(\ast g)$-triples,
and morphisms $\nbigt_1\lrarr\nbigt_2$ 
in $\MTM(X)$.
If $\nbigt_i$ are integrable,
we have the bijection of integrable morphisms.
\hfill\qed
\end{lem}

Let $\MTM_H(X)\subset \MTM(X)$
be the full subcategory of
mixed twistor $D$-modules whose
supports are contained in $H$.
\begin{lem}
For $\nbigt\in\MTM(X)$,
we have $\nbigt[\star H]=0$ if and only if
$\nbigt\in\MTM_H(X)$.
\end{lem}
\pf
The if part is clear.
If $\nbigt[\star H]=0$,
we have $\nbigt(\ast H)=0$,
and hence $\nbigt\in\MTM_H(X)$.
\hfill\qed

\begin{lem}
The functor $[\star H]$ is exact.
\end{lem}
\pf
Let $\varphi:\nbigt_1\lrarr\nbigt_2$
be a morphism in $\MTM(X)$.
Let us consider the induced morphism
$\varphi[\star H]:
\nbigt_1[\star H]\lrarr\nbigt_2[\star H]$.
By using Lemma \ref{lem;11.1.17.11}
and Lemma \ref{lem;10.9.2.3},
we obtain 
$\Ker(\varphi[\star H])
\simeq
 \Ker(\varphi)[\star H]$,
$\Image(\varphi[\star H])
\simeq
 \Image(\varphi)[\star H]$
and
$\Cok(\varphi[\star H])
\simeq
 \Cok(\varphi)[\star H]$
in $\MTM(X)$.
Then, the claim is clear.
\hfill\qed

\vspace{.1in}

We obtain the following
from Corollary \ref{cor;11.2.22.33}
\begin{lem}
Let $\nbigt_i\in\MTM(X)$ $(i=1,2)$.
We have natural bijections:
\[
 \Hom_{\MTM(X)}\bigl(
 \nbigt_1[\ast H],\,\nbigt_2[\ast H]
 \bigr)
\simeq
  \Hom_{\MTM(X)}\bigl(
  \nbigt_1,\,\nbigt_2[\ast H]
 \bigr)
\]
\[
 \Hom_{\MTM(X)}\bigl(
 \nbigt_1[!H],\,\nbigt_2[!H]
 \bigr)
\simeq
  \Hom_{\MTM(X)}\bigl(
 \nbigt_1[!H],\,\nbigt_2
 \bigr)
\]
Similarly,
we have the bijections for integrable morphisms,
if $\nbigt_i$ are integrable.
\hfill\qed
\end{lem}

\begin{prop}
\mbox{{}}
\begin{itemize}
\item
Let $F:X\lrarr Y$ be a projective morphism.
Let $H_Y$ be a hypersurface of $Y$.
We put $H_X:=F^{-1}(H_Y)$.
For $\nbigt\in \MTM(X)$,
we have a natural isomorphism
$\bigl(
 F_{\dagger}(\nbigt)
\bigr)
 [\star H_Y]
\simeq
 F_{\dagger}\bigl(\nbigt[\star H_X]\bigr)$.
\item
For $(\nbigt,L)\in\MTM(X)$,
we have a natural isomorphism
$\bigl(\nbigt[\ast H]\bigr)^{\ast}
\simeq
 \bigl(\nbigt^{\ast}\bigr)[!H]$.
\end{itemize}
\end{prop}
\pf
The second claim is trivial.
As for the first claim,
we have a natural isomorphism
$\bigl(
 F_{\dagger}(\nbigt)
\bigr)
 [\star H_Y](\ast H_Y)
\simeq
 F_{\dagger}\bigl(\nbigt[\star H_X]\bigr)(\ast H_Y)$.
Because 
$F_{\dagger}\bigl(\nbigt[\star H_X]\bigr)$
satisfies the characterization of the localization,
we obtain the desired isomorphism.
\hfill\qed

\subsection{Independence from compactification}

Let $X$ be a complex manifold
with a hypersurface $H$.
Let $F:X'\lrarr X$ be a projective birational
morphism such that
$X'\setminus H'\simeq X\setminus H$,
where $H':=\varphi^{-1}(H)$.
\begin{prop}
\label{prop;11.3.30.10}
The push-forward induces
equivalences of the categories
\[
F_{\dagger}:
 \MTM(X',[\ast H'])
\lrarr
 \MTM(X,[\ast H]),
\]
\[
 F_{\dagger}:
 \MTMint(X',[\ast H'])
\lrarr
 \MTMint(X,[\ast H]).
\]
\end{prop}
\pf
We prove only the ordinary case.
The integrable case can be argued similarly.
We obtain the fully faithfulness
from Lemma \ref{lem;11.2.21.2}.
Let us prove the essential surjectivity.
Let $\nbigt\in\MTM(X,[\ast H])$.
We have the filtered $\nbigr_X(\ast H)$-triple
$\nbigt(\ast H)$.
We have the corresponding
$\nbigr_{X'}(\ast H')$-triple $\nbigt_1$,
as remarked in Lemma \ref{lem;11.2.21.1}.
We have only to show that
there exists $\nbigt'\in\MTM(X')$
such that 
$\nbigt'(\ast H')=\nbigt_1$
as filtered $\nbigr_{X'(\ast H')}$-triples.

First, let us consider the local problem.
Let $P$ be a point of $X$.
We take a small neighbourhood $X_P$ of $P$.
We set $H_P:=H\cap X_P$,
$X_P':=F^{-1}(X_P)$
and $H_P':=F^{-1}(H_P)$.
We set $\nbigt_P:=\nbigt_{|X_P}$.
Let $F_P:=F_{|X_P'}$.
\begin{lem}
\label{lem;13.5.14.30}
If we shrink $X_P$,
there exists $\nbigt'_P\in \MTM(X'_P,[\ast H_P'])$
such that
$F_{P'\dagger}\nbigt'_P\simeq\nbigt_P$.
\end{lem}
\pf
In the proof, we shrink $X$ instead of considering $X_P$.
So, we omit the subscript $P$
to simplify the notation.
We use the Noetherian induction.
We take a cell 
$\nbigc=(Z,U,\varphi,\nbigv)$ of $(\nbigt,L)$
at $P$.
We may assume that $\varphi:Z\lrarr X$
factors through $X'$,
i.e., $\varphi$ is the composition of
$\varphi':Z\lrarr X'$ and $F:X'\lrarr X$.
Let $g$ be a cell function for $\nbigc$
with $H\subset g^{-1}(0)$.
We have the expression of $\nbigt$
as the cohomology of
$\psi^{(1)}_g\varphi_{\dagger}\nbigv
\lrarr
 \Xi^{(0)}_g\varphi_{\dagger}\nbigv\oplus\phi^{(0)}_g(\nbigt)
\lrarr
 \psi^{(0)}_g\varphi_{\dagger}\nbigv$.
Let $g':=g\circ F$.
We have the mixed twistor $D$-modules
$\psi^{(a)}_{g'}\varphi'_{\dagger}\nbigv[\ast H']$
and
$\Xi^{(a)}_{g'}\varphi'_{\dagger}\nbigv[\ast H']$
on $X'$.
By the assumption of the induction,
we have $\nbigq\in\MTM(X')$
such that $\nbigq[\ast H']=\nbigq$
and $F_{\dagger}\nbigq=\phi^{(0)}_{g}\nbigt[\ast H]$.
By the fully faithfulness,
we have the morphisms
\[
 \psi_{g'}^{(1)}\varphi'_{\dagger}\nbigv[\ast H']
\lrarr\nbigq\lrarr
 \psi_{g'}^{(0)}\varphi'_{\dagger}\nbigv[\ast H']
\]
corresponding to
$\psi_{g}^{(1)}\varphi_{\dagger}\nbigv[\ast H]
\lrarr
 \phi_{g}^{(0)}\nbigt[\ast H]
\lrarr
 \psi_{g}^{(0)}\varphi_{\dagger}\nbigv[\ast H]$.
We obtain $\nbigt'\in\MTM(X')$ 
as the cohomology of
\[
 \psi_{g'}^{(1)}\varphi_{\dagger}\nbigv[\ast H']
\lrarr
\Xi_{g'}^{(0)}\varphi_{\dagger}\nbigv[\ast H']\oplus
 \nbigq\lrarr
 \psi_{g'}^{(0)}\varphi_{\dagger}\nbigv[\ast H'].
\]
It satisfies
$\nbigt'(\ast H')=\nbigt_1$.
Thus, we obtain Lemma \ref{lem;13.5.14.30}.
\hfill\qed

\vspace{.1in}
Let $\Ltilde$ denote the weight filtration of
$\nbigt_P'$ as an object in $\MTM(X'_P)$.
We have the decomposition of 
the polarizable pure twistor $D$-modules:
\begin{equation}
\label{eq;13.5.14.40}
 \Gr^{\Ltilde}_w(\nbigt_P')
=\nbigp'_{P,w}
\oplus
 \Ker\Bigl(
 \Gr^{\Nhat_{\ast}L}_w\Gr^W_w\psi_{H'_P}^{(0)}(\nbigt'_P)
\lrarr
 \Gr^{L}_w\Gr^W_w\psi_{H_P'}^{(0)}(\nbigt'_P)
 \Bigr)
\end{equation}
Here, $\nbigp'_{P,w}$ denotes the sum of the direct summand of
$\Gr^{\Ltilde}_w(\nbigt_P')$
whose strict supports are not contained in $H_P'$.

\vspace{.1in}
Let us consider the global problem.
By gluing $\nbigt'_P$ for varied $P$,
we obtain a filtered $\nbigr_{X'}$-triple
$(\nbigt',\Ltilde)$ on $X'$.
It is localizable along $H'$,
and we have
$\nbigt'[\ast H']\simeq\nbigt'$.
To establish
$\nbigt'[\ast H']\in \MTM(X',[\ast H'])$,
we have only to prove that
$\Gr^{\Ltilde}_w(\nbigt')$
is polarizable pure twistor $D$-module of weight $w$.
By gluing (\ref{eq;13.5.14.40}),
we have the following decomposition:
\begin{equation}
 \Gr^{\Ltilde}_w(\nbigt')
=\nbigp'_{w}
\oplus
 \Ker\Bigl(
 \Gr^{\Nhat_{\ast}L}_w\Gr^W_w\psi_{H'}^{(0)}(\nbigt')
\lrarr
 \Gr^{L}_w\Gr^W_w\psi_{H'}^{(0)}(\nbigt')
 \Bigr)
\end{equation}

Let $\Ltilde$ denote the weight filtration of
$\nbigt$ as an object in $\MTM(X)$.
Let $\Gr^{\Ltilde}_w(\nbigt)
=\bigoplus \nbigp_{Z,w}$
be the decomposition by the strict supports,
where $Z$ runs through closed irreducible subvariety
of $X$.
By using the correspondence between
wild harmonic bundles and wild pure twistor $D$-modules,
for $Z\not\subset H$,
we take a projective birational morphism
$\kappa_Z:Z_1\lrarr Z$
satisfying the following conditions:
\begin{itemize}
\item
 $Z_1$ is smooth with a normal crossing hypersurface
 $H_{1}$.
\item
 We have a polarizable pure twistor $D$-module 
 $\gbigp_{Z,w}$ on $Z_1$
 whose strict support is $Z_1$,
 such that
 $\nbigp_{Z,w}$ is the component
 $\kappa^0_{Z\dagger}(\gbigp_{Z,w})$
 whose strict support is $Z$. 
\item
$\kappa_Z$ factors through $X'$,
 i.e.,
 $\iota_Z\circ\kappa_Z=F\circ\kappa_Z'$
 for some $\kappa_Z':Z_1\lrarr X'$,
 where $\iota_Z:Z\lrarr X$ denote the inclusion.
\end{itemize}
Let $\nbigp'_{Z,w}$ be the polarizable pure twistor $D$-module
on $X'$ obtained as the component of
$\kappa'_{Z\dagger}\gbigp_{Z,w}$
whose strict support is $\kappa_Z'(Z_1)$.
We naturally have
$F_{\dagger}\bigl(
 \nbigp'_{Z,w}(\ast H')
 \bigr)
\simeq
 \nbigp_{Z,w}(\ast H)$.
By the uniqueness in the local construction,
we have
$\bigoplus_{Z}\nbigp'_{Z,w|X_P'}\simeq
 \nbigp'_{P,w}$.
Hence, we obtain that
$\nbigp'_{w}=
 \bigoplus_Z\nbigp'_{Z,w}$,
i.e.,
$\nbigp'_w$ is a polarizable pure twistor $D$-module
of weight $w$.

We have the canonical decomposition:
\[
 \Gr^W_w\psi^{(0)}_{H'}(\nbigt')
\simeq
 \Gr^W_w\psi^{(0)}_{H'}(\Gr^{\Ltilde}(\nbigt')) 
\simeq
 \Gr^W_w\psi^{(0)}_{H'}(\nbigp'_w)
\]
Because $\nbigp'_w$ is a polarizable pure twistor $D$-module
of weight $w$,
$\Gr^W_w\psi^{(0)}_{H'}(\nbigt')$
is a polarizable pure twistor $D$-module of weight $w$
by Proposition \ref{prop;13.5.10.102}.
Thus, the proof of Proposition \ref{prop;11.3.30.10}
is finished.
\hfill\qed

\section{Twist by admissible twistor structure}
\label{subsection;11.4.9.13}

\subsection{Smooth case}

Let $\nbigt\in\MTM(X)$
and $\nbigv\in\MTS(X)$.
We have a naturally defined
filtered $\nbigr_X$-triple
$\nbigt\otimes\nbigv$.
\begin{lem}
We have
$\nbigt\otimes\nbigv
\in\MTM(X)$.
If $\nbigt$ and $\nbigv$ are integrable,
then we naturally have
$\nbigt\otimes\nbigv
\in\MTMint(X)$.
\end{lem}
\pf
We have only to use the Noetherian induction.
\hfill\qed

\begin{lem}
Let $F:X\lrarr Y$ be a projective morphism.
For $\nbigt\in\MTM(X)$
and $\nbigv\in\MTS(Y)$,
we have a natural isomorphism
\[
 F_{\dagger}\bigl(
 (\nbigt,W)\otimes F^{\ast}(V,W)
 \bigr)
\simeq
 F_{\dagger}(\nbigt,W)\otimes(V,W).
\]
We have a similar isomorphism
in the integrable case.
\end{lem}
\pf
It follows from
Lemma \ref{lem;10.12.24.3}.
\hfill\qed

\subsection{Admissible case}

Let $\nbigt\in\MTM(X)$.
Let $\nbigt_0\in\MTS^{\adm}(X,H)$.
We naturally obtain
a filtered $\nbigr_{X(\ast H)}$-triple
$\nbigt\otimes\nbigt_0$.
\begin{lem}
\label{lem;10.12.27.30}
Assume that 
$H=f^{-1}(0)$ for some function.
Then, 
$\nbigt\otimes\nbigt_0$
is admissibly specializable along $f$,
and
$(\nbigt\otimes\nbigt_0)[\star f]$
exists in $\MTM(X)$.
If $\nbigt\in\MTMint(X)$
and $\nbigt_0\in\MTS^{\integral\adm}(X,H)$,
then
$\bigl(\nbigt\otimes\nbigt_0\bigr)[\star f]\in\MTMint(X)$.
\end{lem}
\pf
Let $P\in H$.
We may shrink $X$ around $P$.
Let us consider the case that
$\nbigt=
 \varphi_{\dagger}(\nbigv)[\starbar g]$
$(\starbar=\ast,!)$
for an admissible cell $(\varphi,Z,U,\nbigv)$ with
a cell function $g$.
We may assume that $\varphi^{\ast}(fg)^{-1}(0)$
is normal crossing.
Because
$\bigl(
 \nbigt\otimes\nbigt_0[\starbar g]
\bigr)[\star f]
=\Bigl(
 \varphi_{\dagger}\bigl(
 \nbigv\otimes\varphi^{\ast}\nbigt_0
 \bigr)[\starbar gf]
\Bigr)[\star f]$,
it is a mixed twistor $D$-module.
We obtain the claims in the cases that
$\nbigt$ is
$\Xi_g^{(a)}(\varphi_{\dagger}\nbigv)$
or
$\psi_g^{(a)}(\varphi_{\dagger}\nbigv)$.
Then, we obtain the claim in the general case
by a Noetherian induction.
The integrable case can be argued similarly.
\hfill\qed

\begin{prop}
\label{prop;13.7.29.1}
The filtered $\nbigr_X(\ast H)$-triple
$\nbigt\otimes\nbigt_0$ is naturally extended to 
a mixed twistor $D_X$-module
$(\nbigt\otimes\nbigt_0)[\star H]$.
If $\nbigt$ and $\nbigt_0$
are integrable,
then
we naturally have
$(\nbigt\otimes\nbigt_0)[\star H]
\in\MTMint(X)$.
\end{prop}
\pf
Locally it is given in Lemma \ref{lem;10.12.27.30}.
We can construct a filtered $\nbigr$-triple
$\nbigt\otimes\nbigt_0[\star H]$ as the gluing.
By the construction it is admissibly specializable
and localizable along $H$.
We can check that
it is an object in $\MTW(X)$
by the argument in the proof of 
Proposition \ref{prop;11.3.30.10}.
Then, we obtain that 
it is an object in $\MTM(X)$
by the construction.
\hfill\qed

\subsection{Beilinson functors}

As a consequence of Proposition \ref{prop;13.7.29.1},
for any holomorphic function $g$ on $X$,
we obtain Beilinson functors
$\Pi^{a,b}_{g\ast !}$ on $\MTM(X)$,
in particular,
$\Xi^{(a)}_g$ and $\psi_g^{(a)}$.
We obtain $\phi^{(a)}_g$ as the cohomology
of the natural complex
$\nbigt[!g]\lrarr
 \Xi^{(a)}_g(\nbigt)
 \oplus
 \nbigt\lrarr
 \nbigt[\ast g]$.
We can reconstruct $\nbigt$ as
the cohomology of the natural complex
$\psi^{(1)}_g(\nbigt)
\lrarr
 \Xi^{(0)}_g(\nbigt)
\oplus
 \phi^{(0)}_g(\nbigt)
\lrarr
 \psi^{(0)}(\nbigt)$.

\section{Exterior tensor product}

Let $X_i$ $(i=1,2)$ be complex manifolds.
We set $X:=X_1\times X_2$.
For any  $\nbigt_i\in\MTM(X)$,
we obtain a naturally defined
filtered $\nbigr_X$-triple
$\nbigt=\nbigt_1\boxtimes\nbigt_2$.
The weight filtration $L$ on $\nbigt$
is given by
$L_k(\nbigt)=\sum_{i+j\leq k}
 L_{i}(\nbigt_1)\boxtimes L_j(\nbigt_2)$.

\begin{prop}
\label{prop;13.7.29.5}
$\nbigt$ is a mixed twistor $D$-module on $X$.
\end{prop}

\subsection{Pure case}

Let us consider the pure case.

\begin{lem}
\label{lem;13.7.29.3}
Suppose that $\nbigt_i$ are 
polarizable pure twistor $D$-module 
of weight $w_i$ on $X_i$.
Then, 
$\nbigt_1\boxtimes\nbigt_2$
is a polarizable pure twistor $D$-module 
of weight $w_1+w_2$ on $X_1\times X_2$.
\end{lem}
\pf
We can take complex manifolds $Z_i$,
projective morphisms
$\varphi_i:Z_i\lrarr X_i$,
and polarizable pure twistor $D$-modules $\gbigt_i$
on $Z_i$,
such that
$\nbigt_i$ are direct summands of
$\varphi^0_{i\dagger}(\gbigt_i)$.
Then,
$\nbigt_1\boxtimes\nbigt_2$
is a direct summand of
$\varphi^0_{1\dagger}(\gbigt_1)
 \boxtimes
 \varphi^0_{2\dagger}(\gbigt_2)
\simeq
 (\varphi_1\times\varphi_2)^0_{\dagger}
 (\gbigt_1\boxtimes\gbigt_2)$.
We have only to prove
that
$\gbigt_1\boxtimes\gbigt_2$
is a polarizable pure twistor $D$-module
of weight $w_1+w_2$.
Hence, we have only to consider the case
$\Supp(\nbigt_i)=X_i$.
We may also assume that
there exist normal crossing hypersurface
$H_i\subset X_i$ 
and good wild variation of pure twistor structure
$\nbigv_i$ on $(X_i,H_i)$
such that $\nbigt_i$ 
are obtained as the image of
$\nbigv_i[!H_i]\lrarr
\nbigv_i[\ast H_i]$
as $\nbigr_{X_i}$-triples
Then, the $\nbigr_{X_1\times X_2}$-triple
$\nbigv_1\boxtimes\nbigv_2$
is isomorphic to the image of
$\nbigv_1[!H_1]
\boxtimes
 \nbigv_2[!H_2]
\lrarr
\nbigv_1[\ast H_1]
\boxtimes
 \nbigv_2[\ast H_2]$.
Hence, it is a polarizable pure twistor
$D$-module of weight $w_1+w_2$.
\hfill\qed

\begin{cor}
For $\nbigt_i\in\MTM(X_i)$,
the filtered $\nbigr_{X_1\times X_2}$-triple
$\nbigt_1\boxtimes\nbigt_2$
is an object in $\MTW(X_1\times X_2)$.
\hfill\qed
\end{cor}

\subsection{Admissible variation of mixed twistor structure}

Let $Y:=\Delta^n$
and $H:=\bigcup_{i=1}^{\ell}\{z_i=0\}$.
Let $\nbigv\in\MTS^{\adm}(Y,H)$.
Let $I\sqcup J=\ellsitabar$ be a decomposition.
We have a mixed twistor $D$-module
$\nbigv[\ast I!J]$ on $Y$
with the weight filtration $L$.
We have the following characterization.
\begin{lem}
\label{lem;13.7.29.2}
Suppose that $(\nbigt,L)\in\MTW(Y)$
satisfies the following conditions:
\begin{itemize}
\item
 $(\nbigt,L)(\ast H)=\nbigv$.
\item
 $(\nbigt,L)$ is admissibly specializable
 along any $z_i$ $(i=1,\ldots,\ell)$.
\item
 $(\nbigt,L)[\star_iz_i]\simeq (\nbigt,L)$
 as filtered $\nbigr$-triples,
 where
 $\star_i=\ast$ $(i\in I)$
 or $\star_i=!$ $(i\in J)$.
\end{itemize}
\end{lem}
\pf
Put $H(j)\!:=\!\bigcup_{i=j}^{\ell}\{z_i=0\}$.
We can prove
$(\nbigt,L)(\ast H(j))\!\simeq\!
 (\nbigv[\ast I!J])\bigl(\ast H(j)\bigr)$
as filtered $\nbigr_Y(\ast H(j))$-triples
by an easy induction.
\hfill\qed

\vspace{.1in}
Let $K_i$ $(i=1,2)$ be finite sets
with a finite set $L_i\subset K_i$.
Let $Y_i:=\Delta^{K_i}$
and $H_i:=\bigcup_{j\in K_i}\{z_j=0\}$.
We set $Y:=Y_1\times Y_2$
and $H:=(Y_1\times H_2)\cup(H_1\times Y_2)$.
Let $\nbigv_i\in\MTS^{\adm}(Y_i,H_i)$.
We set
$\nbigv_1\boxtimes\nbigv_2
 \in\MTS^{\adm}(Y,H)$.
Let $I_i\sqcup J_i=L_i$ be decompositions.
We put $I:=I_1\sqcup I_2$
and $J:=J_1\sqcup J_2$.
By using Lemma \ref{lem;13.7.29.2}
and Lemma \ref{lem;13.7.29.3},
we obtain the following lemma.
\begin{lem}
\label{lem;13.7.29.10}
We have a natural isomorphism
$\nbigv_1[\ast I_1!J_1]
\boxtimes
 \nbigv_2[\ast I_2!J_2]
\!\simeq\!
 \nbigv[\ast I!J]$.
In particular,
$\nbigv_1[\ast I_1!J_1]
\boxtimes
 \nbigv_2[\ast I_2!J_2]
\in
 \MTM(Y)$.
\hfill\qed
\end{lem}

\subsection{Proof of Proposition \ref{prop;13.7.29.5}}

We have already known
$\nbigt_1\boxtimes\nbigt_2
\in\MTW(X)$.
Because the conditions in Definition
\ref{df;11.4.3.30} are local,
we may and will shrink $X_i$.

Let $\nbigc=(Z,U,\varphi,\nbigv)$ 
be any admissible cell
with a cell function $g$ on $X_1$.
\begin{lem}
\label{lem;13.7.29.12}
$\Pi^{a,b}_{g\ast !}(\varphi_{\dagger}\nbigv)
\boxtimes
 \nbigt_2\in\MTM(X)$.
In particular,
$\Xi_g^{(a)}(\varphi_{\dagger}\nbigv)
\boxtimes
 \nbigt_2\in\MTM(X)$
and 
$\psi_g^{(a)}(\varphi_{\dagger}\nbigv)
\boxtimes
 \nbigt_2\in\MTM(X)$
\end{lem}
\pf
We use a Noetherian induction
on the support of $\nbigt_2$.
Shrinking $X_2$,
we can take an admissible cell
$\nbigc_2=(Z_2,U_2,\varphi_2,\nbigv_2)$
with a cell function $g_2$
for $\nbigt_2$.
Then, 
$\Pi^{a,b}_{g\ast !}(\varphi_{\dagger}\nbigv)
\boxtimes
 \nbigt_2$
is isomorphic to the cohomology
of the following complex
of filtered $\nbigr_{X_1\times X_2}$-triples:
\[
 \Pi^{a,b}_{g\ast!}(\varphi_{\dagger}\nbigv)
\boxtimes
\Bigl(
\psi^{(1)}_{g_2}(\varphi_{2\dagger}\nbigv_2)
\lrarr
 \Xi^{(0)}_{g_2}(\varphi_{2\dagger}\nbigv_2)
\oplus
 \phi^{(0)}_{g_2}(\varphi_{2\dagger}\nbigt_2)
\lrarr
 \psi^{(1)}_{g_2}(\varphi_{2\dagger}\nbigv_2)
\Bigr)
\]
By the hypothesis of the induction,
we know that
$\Pi^{a,b}_{g\ast !}(\varphi_{\dagger}\nbigv)
\boxtimes
 \psi^{(a)}_{g_2}(\varphi_{2\dagger}\nbigv_2)$
and 
$\Pi^{a,b}_{g\ast !}(\varphi_{\dagger}\nbigv)
\boxtimes
 \phi^{(a)}_{g_2}(\varphi_{2\dagger}\nbigv_2)$
are mixed twistor $D$-modules.
By using Lemma \ref{lem;13.7.29.10},
we obtain that
$\Pi^{a,b}_{g\ast !}(\varphi_{\dagger}\nbigv)
\boxtimes
 \Xi^{(0)}_{g_2}(\varphi_{2\dagger}\nbigv_2)
\in\MTM(X)$.
Then, the claim of Lemma \ref{lem;13.7.29.12}
follows.
\hfill\qed

\vspace{.1in}
Let us finish the proof of 
Proposition \ref{prop;13.7.29.5}.
We use a Noetherian induction
on the support of $\nbigt_1$.
We take an admissible cell
$\nbigc_1=(Z_1,U_1,\varphi_1,\nbigv_1)$
with a cell function $g_1$ for $\nbigt_1$.
Then, the filtered $\nbigr_X$-triple
$\nbigt_1\boxtimes\nbigt_2$
is isomorphic to the cohomology of
the following complex:
\[
 \Bigl(
 \psi^{(1)}_{g_1}(\varphi_{1\dagger}\nbigv_1)
\lrarr
 \Xi^{(0)}_{g_1}(\varphi_{1\dagger}\nbigv_1)
\oplus
 \phi^{(0)}_{g_1}(\nbigt_1)
\lrarr
 \psi^{(0)}_{g_1}(\varphi_{1\dagger}\nbigv_1)
 \Bigr)
\boxtimes
 \nbigt_2
\]
By Lemma \ref{lem;13.7.29.12},
$\psi^{(a)}_{g_1}(\varphi_{1\dagger}\nbigv_1)
 \boxtimes\nbigt_2$
and 
$\Xi^{(0)}_{g_1}(\varphi_{1\dagger}\nbigv_1)
\boxtimes\nbigt_2$
are mixed twistor $D$-modules on $X$.
By the hypothesis of the induction,
$\phi^{(0)}_{g_1}(\varphi_{1\dagger}\nbigv_1)
 \boxtimes\nbigt_2$
is a mixed twistor $D$-module.
Hence, we obtain that
$\nbigt_1\boxtimes\nbigt_2$
is a mixed twistor $D$-module on $X_1\times X_2$.
Thus, the proof of Proposition \ref{prop;13.7.29.5}
is finished.
\hfill\qed

\chapter{Dual and real structure of mixed twistor $D$-modules}
\label{section;11.4.9.20}

We introduce the dual of mixed twistor $D$-modules.
A preparation is given in \S\ref{section;11.4.6.10}.
We also introduce the notion of real structure
of mixed twistor $D$-modules.

\section{Dual of $\nbigr_X$-modules}

\subsection{Dual}
\index{dual}

Let $X$ be any complex manifold,
and let $H$ be any hypersurface.
Let $N$ be any left-$\nbigr_{X(\ast H)}$-bi-module,
i.e., it is equipped with mutually commuting 
two $\nbigr_{X(\ast H)}$-actions $\rho_i$ $(i=1,2)$.
Let $(N,\rho_i)$ denote the left $\nbigr_{X(\ast H)}$-modules
by $\rho_i$.
For any $\nbigr_{X(\ast H)}$-module $L$,
let $\nhom_{\nbigr_{X(\ast H)}}
 (L,N^{\rho_1,\rho_2})$
denote  the sheaf of $\nbigr_{X(\ast H)}$-homomorphisms
$L\lrarr (N,\rho_1)$.
The sheaf is equipped with an $\nbigr_{X(\ast H)}$-action
induced by $\rho_2$.
For any $\nbigr_{X(\ast H)}$-complex $L^{\bullet}$,
we obtain an $\nbigr_{X(\ast H)}$-complex
$\nhom_{\nbigr_{X(\ast H)}}
 (L^{\bullet},N^{\rho_1,\rho_2})$.
We have the naturally defined
$\nbigr_{X(\ast H)}$-homomorphism
\begin{equation}
 \label{eq;11.1.22.20}
 L^{\bullet}\lrarr \nhom_{\nbigr_{X(\ast H)}}\Bigl(
 \nhom_{\nbigr_{X(\ast H)}}\bigl(L^{\bullet},N^{\rho_2,\rho_1}
 \bigr),
 N^{\rho_1,\rho_2}
 \Bigr)
\end{equation}
given by
$x\longmapsto
 \bigl(
 F\longmapsto
 (-1)^{|x|\,|F|}F(x)
 \bigr)$.

\vspace{.1in}

Let $\Theta_X$ denote the tangent sheaf of $X$,
and let $\Omega^1_X:=\nhom_{\nbigo_X}\bigl(\Theta_X,\nbigo_X\bigr)$.
\index{sheaf $\Theta_X$}
\index{sheaf $\Omega^1_X$}
Let $p_{\lambda}:\nbigx\lrarr X$
denote the projection.
Let $\Theta_{\nbigx}:=\lambda\cdot
 p_{\lambda}^{\ast}\Theta_X$,
and let $\Omega^1_{\nbigx}:=\nhom_{\nbigo_{\nbigx}}\bigl(
 \Theta_{\nbigx},\nbigo_{\nbigx} \bigr)$.
\index{sheaf $\Theta_{\nbigx}$}
\index{sheaf $\Omega^1_{\nbigx}$}
We have 
$\Omega^1_{\nbigx}=
 \lambda^{-1}p_{\lambda}^{\ast}\Omega_X^1$.
Its $j$-th exterior product is denoted by 
$\Omega^j_{\nbigx}$.
In particular, we set
$\omega_{\nbigx}:=\Omega_{\nbigx}^{\dim X}$.
\index{sheaf $\omega_{\nbigx}$}

\vspace{.1in}
Recall that we have the natural two $\nbigr_{X(\ast H)}$-action
on $\nbigr_{X(\ast H)}\otimes\omega_{\nbigx}^{-1}$.
The left multiplication is denoted by $\ell$.
Let $r$ denote the action induced by the right
multiplication.
More generally, for a left $\nbigr_{X(\ast H)}$-module $N$,
we have two induced left $\nbigr_{X(\ast H)}$-module
structures on 
$N\otimes \nbigr_{X(\ast H)}\otimes\omega_{\nbigx}^{-1}$.
One is given by $\ell$ and
the left $\nbigr$-action on $N$,
which is denoted by $\ell$.
The other is induced by the right multiplication,
denoted by $r$.
We have the $\cnum$-linear isomorphism
$\Phi_N:N\otimes 
 \nbigr_{X(\ast H)}\otimes\omega_{\nbigx}^{-1}
\simeq
 \nbigr_{X(\ast H)}\otimes\omega_{\nbigx}^{-1}$,
such that
$\Phi_N\circ r=\ell\circ\Phi_N$
and $\Phi_N\circ\ell=r\circ\Phi_N$,
as in \cite{sabbah2} and \cite{saito1}.

Let $\nbigg^{\bullet}$ be 
a $\nbigr_{X(\ast H)}$-bimodule resolution of
$\nbigr_{X(\ast H)}\otimes\omega_{\nbigx}^{-1}$.
It also gives an $\nbigr_{X(\ast H)}$-module resolution
with respect to both of $\ell$ and $r$.
Let $\nbigm$ be a coherent $\nbigr_{X(\ast H)}$-module.
We define
\[
 \DDD'_{X(\ast H)} \nbigm:=
 \nhom_{\nbigr_{X(\ast H)}}\bigl(\nbigm,\nbigg^{\bullet\,\ell,r}
 \bigr),
 \quad\quad
 \DDD_{X(\ast H)}\nbigm:=
 \lambda^{d_X}\cdot\DDD'_{X(\ast H)}\nbigm.
\]
\index{sheaf $\DDD'_{X(\ast H)}\nbigm$}
\index{sheaf $\DDD_{X(\ast H)}\nbigm$}
The morphism as in (\ref{eq;11.1.22.2})
with $\Phi_{\nbigg^{\bullet}}$ induces
an isomorphism
$\nbigm\lrarr \DDD_{X(\ast H)}\circ\DDD_{X(\ast H)}(\nbigm)$.

\subsection{Compatibility with push-forward}
\label{subsection;10.9.17.20}

Let $F:X\lrarr Y$ be a morphism of complex manifolds.
We shall construct a trace morphism
\begin{equation}
 \label{eq;10.9.16.1}
 \lambda^{d_X}\cdot
 F_{\dagger}(\nbigo_{\nbigx})[d_X]
\lrarr
 \lambda^{d_Y}\cdot \nbigo_{\nbigy}[d_Y]
\end{equation}
in $D^b(\nbigr_Y)$ by a standard method.
Let $\distribution_{\nbigx/\cnum_{\lambda}}$
denote the sheaf of distributions on $\nbigx$
which are holomorphic in $\lambda$.
We set 
\[
 \distribution_{\nbigx/\cnum_{\lambda}}^{\bullet}
:=
\Tot\Bigl(
 p_{\lambda}^{-1}\Omega_X^{0,\bullet}
 \otimes_{p_{\lambda}^{-1}C^{\infty}(X)}
 \Omega_{\nbigx}^{\bullet}
 \otimes_{\nbigo_{\nbigx}}
 \distribution_{\nbigx/\cnum_{\lambda}}
\Bigr)
\]
We have a natural quasi-isomorphism
$\Omega^{\bullet}_{\nbigx}\lrarr
 \distribution_{\nbigx/\cnum_{\lambda}}^{\bullet}$,
and hence an isomorphism
$\nbigr_X\otimes_{\nbigo_{\nbigx}}
 \bigl(
 \distribution_{\nbigx/\cnum_{\lambda}}^{\bullet}
 \bigr)[d_X]
\simeq
 \omega_{\nbigx}$
in $D^b(\nbigr_X)$.
We obtain the following morphism of complexes,
by the integration multiplied with
$(2\pi\sqrt{-1})^{-d_X+d_Y}$:
\[
 \lambda^{d_X}\cdot
 F_{!}\distribution_{\nbigx/\cnum_{\lambda}}^{\bullet}[2d_X]
\lrarr
 \lambda^{d_Y}\cdot
 \distribution_{\nbigy/\cnum_{\lambda}}^{\bullet}[2d_Y]
\]
Hence, we obtain 
the trace morphism (\ref{eq;10.9.16.1})
as follows:
\begin{multline}
 F_{\dagger}\nbigo_{\nbigx}[d_X]
\lrarr
 F_!\bigl(
 \nbigr_{Y\larr X}\otimes^L_{\nbigr_X}
 \nbigo_{\nbigx}[d_X]
 \bigr) \\
=F_{!}\bigl(
 \distribution^{\bullet}_{\nbigx/\cnum_{\lambda}}
 \otimes_{F^{-1}\nbigo_{\nbigy}}
 F^{-1}(\nbigr_Y\otimes\omega_{\nbigy}^{-1})[2d_X]
 \bigr) \\
\lrarr
 \lambda^{d_Y-d_X}\cdot
 \distribution_{\nbigy/\cnum_{\lambda}}^{\bullet}
 \otimes\nbigr_Y\otimes\omega_{\nbigy}^{-1}[2d_Y]
\simeq
 \lambda^{d_Y-d_X}\cdot \nbigo_{\nbigy}[d_Y]
\end{multline}
As in the case of $D$-modules
(see (\ref{eq;11.3.21.21}) below),
we obtain the following morphism:
\begin{multline}
 F_{\dagger}\bigl(
 \nbigo_{\nbigx}[d_X]
 \otimes_{F^{-1}\nbigo_{\nbigy}}
 \bigl(\nbigr_Y\otimes\omega_{\nbigy}^{-1}\bigr)
 \bigr)
\lrarr
 \lambda^{d_Y-d_X}\nbigo_{\nbigy}[d_Y]
 \otimes_{\nbigo_{\nbigy}}
 \bigl(\nbigr_Y\otimes\omega_{\nbigy}^{-1}\bigr)
 \\
\simeq
 \lambda^{d_Y-d_X}\nbigr_Y\otimes\omega_{\nbigy}^{-1}[d_Y]
\end{multline}

Let $\nbigm$ be a coherent
$\nbigr_X$-module.
We have a natural morphism
$\varphi:
 F_{\dagger}\DDD_X\nbigm\lrarr
 \DDD_YF_{\dagger}\nbigm$
given as follows:
{\small
\begin{multline}
 F_{\dagger}\nrhom_{\nbigr_X}\bigl(
 \nbigm,\nbigr_X\otimes\omega_{\nbigx}^{-1}
\bigr)[d_X]
\lrarr
 F_{!}\nrhom_{\nbigr_X}\Bigl(
 \nbigm,\nbigo_{\nbigx}\otimes_{F^{-1}\nbigo_{\nbigy}}
 F^{-1}\bigl(\nbigr_Y\otimes\omega_{\nbigy}^{-1}\bigr)
 \Bigr)[d_X] \\
\lrarr
 F_{!}\nrhom_{F^{-1}\nbigr_Y}\Bigl(
 \nbigr_{Y\larr X}\otimes^L_{\nbigr_Y}\nbigm,\,
 \nbigr_{Y\larr X}\otimes^L_{\nbigr_Y}
 \nbigo_{\nbigx}\otimes_{F^{-1}\nbigo_{\nbigy}}
 F^{-1}\bigl(\nbigr_Y\otimes\omega_{\nbigy}^{-1}\bigr)
 \Bigr)[d_X] \\
\lrarr
 \nrhom_{\nbigr_Y}\bigl(
 F_{\dagger}\nbigm,\,
 F_{\dagger}(\nbigo_{\nbigx}\otimes_{F^{-1}\nbigo_{\nbigy}}
 F^{-1}(\nbigr_Y\otimes\omega_{\nbigy}^{-1}))
 \bigr)[d_X]
 \\
\lrarr
 \nrhom_{\nbigr_Y}\bigl(
 F_{\dagger}\nbigm,\,
 (\nbigr_Y\otimes\omega_{\nbigy}^{-1})
 \bigr)[d_Y]\cdot \lambda^{d_Y-d_X}
\end{multline}
}

The following lemma can be shown
by an argument in the proof
of Proposition 4.39 of \cite{kashiwara_text}.
\begin{lem}
\label{lem;10.9.17.10}
Assume that (i) the restriction of $F$ to $\Supp\nbigm$
is proper, (ii)  $\nbigm$ is good over $Y$.
Then, $\varphi$ is an isomorphism.
\hfill\qed
\end{lem}

\subsection{Twist by smooth $\nbigr$-modules}
\label{subsection;11.1.29.2}

Let $X$ and $H$ be as above.
Let $\nbigm$ be a holonomic $\nbigr_{X(\ast H)}$-module.
Let $V$ be a smooth $\nbigr_{X(\ast H)}$-module.
We have the natural isomorphism
$\DDD_{X(\ast H)}(\nbigm\otimes V)\simeq
 \DDD_{X(\ast H)}(\nbigm)\otimes V^{\lor}$,
given as follows.
We have two naturally induced
left $\nbigr$-module structures
$\ell$ and $r$
on $V^{\lor}\otimes_{\nbigo_{\nbigx}}
 (\nbigr_{X(\ast H)}\otimes\omega_{\nbigx}^{-1})$.
Then, we obtain
\begin{multline}
 \DDD_{X(\ast H)}(\nbigm\otimes V)
\simeq
 \nrhom_{\nbigr_X(\ast H)}\bigl(\nbigm\otimes V,
 (\nbigr_{X(\ast H)}\otimes\omega_{\nbigx}^{-1})^{\ell,r}
 \bigr)[d_X] \\
\simeq 
 \nrhom_{\nbigr_X(\ast H)}\bigl(\nbigm,
 (V^{\lor}\otimes 
 \nbigr_{X(\ast H)}\otimes\omega_{\nbigx}^{-1})^{\ell,r}
 \bigr)[d_X].
\end{multline}
By $\Phi_{V^{\lor}}$, it is isomorphic to
\begin{multline}
 \nrhom_{\nbigr_X(\ast H)}\bigl(
 \nbigm, (V^{\lor}\otimes 
 \nbigr_{X(\ast H)}\otimes\omega_{\nbigx}^{-1})^{r,\ell}
 \bigr)[d_X]
\simeq \\
 \nrhom_{\nbigr_X(\ast H)}\bigl(
 \nbigm,(\nbigr_{X(\ast H)}
 \otimes\omega_{\nbigx}^{-1})^{r,\ell}
 \bigr)
\otimes V^{\lor}[d_X]. 
\end{multline}
By $\Phi_{\nbigo_{\nbigx}(\ast \nbigh)}$,
it is isomorphic to
$\nrhom_{\nbigr_X(\ast H)}\bigl(
 \nbigm,(\nbigr_{X(\ast H)}\otimes
 \omega_{\nbigx}^{-1})^{\ell,r}
 \bigr)
\otimes V^{\lor}[d_X]$.
We obtain the desired isomorphism
$\DDD_{X(\ast H)}(\nbigm\otimes V)\simeq
 \DDD_{X(\ast H)}(\nbigm)\otimes V^{\lor}$.

\subsection{Dual of smooth triples and the double dual
(Appendix)}
\label{subsection;13.4.1.10}

Let $\nbigm$ be a smooth 
$\nbigr_{X(\ast H)}$-module.
We set $\nbigm^{\lor}:=
 \nhom_{\nbigo_{\nbigx(\ast \nbigh)}}
 (\nbigm,\nbigo_{\nbigx}(\ast \nbigh))$,
which is also naturally 
a smooth $\nbigr_{X(\ast H)}$-module.
Although it is easy to see that
$\nbigm^{\lor}$ is isomorphic to 
$\DDD'_{X(\ast H)} \nbigm$,
let us look at the isomorphism more closely,
and check the signature
of an isomorphism
$\nbigm\simeq
 \DDD'\circ \DDD'(\nbigm)
\simeq
 \DDD'(\nbigm^{\lor})
\simeq
 \nbigm$.
For simplicity of the description,
we omit to denote $(\ast \nbigh)$.

\vspace{.1in}

We have two natural $\nbigr_X$-actions $\ell$ and $r$
on $\nbigr_X\otimes\omega_{\nbigx}^{-1}$.
We have an $\nbigr_X\otimes \nbigr_X$-homomorphism
\[
 C:\bigl(\nbigr_X\otimes\Theta^{-\bullet}_{\nbigx}\bigr)
 \otimes
 \bigl(\nbigr_X\otimes\Theta^{-\bullet}_{\nbigx}\bigr)
\lrarr
 \nbigr_X\otimes\omega_{\nbigx}^{-1}[d_X]
\]
given as follows:
\[
 C(P_1\otimes\tau_1,P_2\otimes\tau_2)
=\left\{
 \begin{array}{ll}
 (-1)^{|\tau_2|}\ell(P_1)\,r(P_2)
 \tau_1\wedge\tau_2
 & (|\tau_1|+|\tau_2|=d_X)\\
 0& \mbox{\rm otherwise}
 \end{array}
 \right.
\]
If there exists a holomorphic coordinate
$(z_1,\ldots,z_n)$,
and if $\tau_j=\bigwedge_{j\in I_i}\deldel_{j}$
such that $I_1\sqcup I_2=\{1,\ldots,n\}$,
we have 
$C(P_1\otimes\tau_1,P_2\otimes \tau_2)
=(-1)^{|\tau_2|}
 P_1\,\lefttop{t}P_2\,\tau_1\wedge\tau_2$.
\begin{lem}
$C$ is a morphism of complexes.
\end{lem}
\pf
We have only to show that
$C\circ\del(P_1\otimes\tau_1\otimes
 P_2\otimes\tau_2)=0$
if $|\tau_1|+|\tau_2|=d+1$.
We may assume that there exists
a holomorphic coordinate $(z_1,\ldots,z_n)$,
and we have only to consider the case
$\tau_i=\bigwedge_{j\in I_i}\deldel_j$.
We have
\begin{multline}
\del(P_1\otimes\tau_1\otimes P_2\otimes\tau_2)
=\sum_{i=1}^{n}
 P_1\deldel_i\otimes\bigl(\iota(\lambda^{-1}dz_i)\tau_1\bigr)
\otimes (P_2\otimes\tau_2)
 \\
+\sum_{i=1}^{n}
 (-1)^{|\tau_1|}
 P_1\otimes\tau_1\otimes
 (P_2\deldel_i\otimes\iota(\lambda^{-1}dz_i)\tau_2)
\end{multline}
Hence, we have the following:
\begin{multline}
C\circ\del\bigl(
 P_1\otimes\tau_1\otimes P_2\otimes\tau_2
 \bigr)
=\sum_{i=1}^nP_1\deldel_i\lefttop{t}P_2
\otimes
 \iota(\lambda^{-1}dz_i)\tau_1\wedge\tau_2(-1)^{|\tau_2|}
 \\
+\sum_{i=1}^nP_1\lefttop{t}(P_2\deldel_i)
\otimes
 \tau_1\wedge\iota(\lambda^{-1}dz_i)\tau_2
 (-1)^{|\tau_2|+|\tau_1|-1}
\\
=\sum_{i=1}^n
 P_1\deldel_i\lefttop{t}P_2\otimes
\Bigl(
 \iota(\lambda^{-1}dz_i)\tau_1\wedge\tau_2
+(-1)^{|\tau_1|}
 \tau_1\wedge\iota(\lambda^{-1}dz_i)\tau_2
\Bigr)\,(-1)^{|\tau_2|}=0
\end{multline}
We have used
$\tau_1\wedge\tau_2=0$.
Thus, we are done.
\hfill\qed

\vspace{.1in}
We obtain an induced isomorphism
of $\nbigr_X$-complexes
\[
 \Psi:
 \nbigr_X\otimes\Theta^{-\bullet}
\lrarr
 \nhom_{\nbigr_X}\bigl(
 \nbigr_X\otimes\Theta^{-\bullet},
 \bigl(
 \nbigr_X\otimes\omega_{\nbigx}^{-1}[d_X]
 \bigr)^{\ell,r}
 \bigr),
\]
given by
$\Psi(P\otimes\tau)(Q\otimes\omega)
=(-1)^{|\tau|\,|\omega|}
 C(Q\otimes\omega,P\otimes\tau)$.

Let $C_{\nbigm}:
 \nbigm^{\lor}\otimes \nbigm\lrarr \nbigo_{\nbigx}$
be the natural perfect pairing
of the $\nbigo_{\nbigx}$-modules.
We obtain an $\nbigr_X\otimes\nbigr_{X}$-homomorphism
\[
 C:\bigl(\nbigr_X\otimes\Theta^{-\bullet}\otimes
 \nbigm^{\lor}\bigr)
 \otimes
 \bigl(\nbigr_X\otimes\Theta^{-\bullet}\otimes\nbigm\bigr)
\lrarr
 \nbigr_X\otimes\omega_{\nbigx}^{-1}[d_X]
\]
given as follows:
{\small
\[
 C(P_1\otimes\tau_1\otimes x_1,P_2\otimes\tau_2\otimes x_2)
\!=\!\left\{
 \begin{array}{ll}
 (-1)^{|\tau_2|}\ell(P_1)\,r(P_2)\,
 \bigl(
 \tau_1\wedge\tau_2\,C_{\nbigm}(x_1,x_2)
 \bigr)
 & (|\tau_1|+|\tau_2|=d_X)\\
 0& \mbox{\rm (otherwise)}
 \end{array}
 \right.
\]
}
It is a morphism of complexes.
Hence, we obtain an induced isomorphism
of $\nbigr_X$-complexes:
\[
 \Psi_{\nbigm}:
 \nbigr_X\otimes\Theta^{-\bullet}\otimes \nbigm^{\lor}
\lrarr
 \nhom_{\nbigr_X}\bigl(
 \nbigr_X\otimes\Theta^{-\bullet}\otimes \nbigm,
 \bigl(
 \nbigr_X\otimes\omega_{\nbigx}^{-1}[d_X]
 \bigr)^{\ell,r}
 \bigr)
\]
By the natural quasi isomorphism
$\nbigr_X\otimes\Theta^{-\bullet}\otimes \nbigm
\simeq \nbigm$,
we obtain an isomorphism
$\nbigm^{\lor}\simeq \DDD' \nbigm$
in the derived category.
We obtain an isomorphism
$(\Psi_{\nbigm}^{\ast})^{-1}\circ\Psi_{\nbigm^{\lor}}:
 \nbigm\simeq \DDD'\circ \DDD'(\nbigm)$
in the derived category:
\begin{multline}
 \nbigr_X\otimes\Theta^{-\bullet}\otimes \nbigm
\stackrel{\Psi_{\nbigm^{\lor}}}{\lrarr}
 \nhom_{\nbigr_X}\bigl(
 \nbigr_X\otimes\Theta^{-\bullet}\otimes\nbigm^{\lor},
 (\nbigr_X\otimes\omega_{\nbigx}^{-1})^{\ell,r}[d_X]
 \bigr)\\
\stackrel{\Psi_{\nbigm}^{\ast}}{\llarr}
 \nhom_{\nbigr_X}\Bigl(
 \nhom_{\nbigr_X}\bigl(
 \nbigr_X\otimes\Theta^{-\bullet}\otimes\nbigm,
 (\nbigr_X\otimes\omega_{\nbigx}^{-1})^{\ell,r}
 \bigr),
 (\nbigr_X\otimes\omega_{\nbigx}^{-1})^{\ell,r}[d_X]
\Bigr)
\end{multline}
Here, $\Psi_{\nbigm}^{\ast}(F):=F\circ\Psi_{\nbigm}$.

\vspace{.1in}

Let $\Phi:\nbigr_X\otimes\omega_{\nbigx}^{-1}[d_X]
\simeq \nbigr_X\otimes\omega_{\nbigx}^{-1}[d_X]$
be the automorphism,
which exchanges the actions $\ell$ and $r$.
It induces
\begin{multline}
 \Phi_{\ast}:
 \nhom_{\nbigr_X}\bigl(
 \nbigr_X\otimes\Theta^{-\bullet}\otimes \nbigm,
 (\nbigr_X\otimes\omega_{\nbigx}^{-1})^{\ell,r}[d_X]
 \bigr)
\simeq \\
  \nhom_{\nbigr_X}\bigl(
 \nbigr_X\otimes\Theta^{-\bullet}\otimes \nbigm,
 (\nbigr_X\otimes\omega_{\nbigx}^{-1})^{r,\ell}[d_X]
 \bigr).
\end{multline}
It induces the following morphism:
\begin{multline}
 \label{eq;11.1.22.3}
 \Phi^{\ast}:
 \nhom_{\nbigr_X}\Bigl(
  \nhom_{\nbigr_X}\bigl(
 \nbigr_X\otimes\Theta^{-\bullet}\otimes \nbigm,
 (\nbigr_X\otimes\omega_{\nbigx}^{-1})^{r,\ell}
 \bigr),
  (\nbigr_X\otimes\omega_{\nbigx}^{-1})^{\ell,r}[d_X]
 \Bigr)
\simeq \\
 \nhom_{\nbigr_X}\Bigl(
  \nhom_{\nbigr_X}\bigl(
 \nbigr_X\otimes\Theta^{-\bullet}\otimes \nbigm,
 (\nbigr_X\otimes\omega_{\nbigx}^{-1})^{\ell,r}
 \bigr),
 (\nbigr_X\otimes\omega_{\nbigx}^{-1})^{\ell,r}[d_X]
\Bigr) \\
=\DDD'\circ \DDD'(\nbigm)
\end{multline}
The natural morphism (\ref{eq;11.1.22.20})
and (\ref{eq;11.1.22.3})
give a quasi-isomorphism
$\Lambda:
 \nbigr_X\otimes\Theta^{-\bullet}\otimes \nbigm
\lrarr
 \DDD'\circ \DDD'(\nbigm)$.

\begin{lem}
\label{lem;11.1.22.1}
We have $\Lambda=
 (-1)^{d_X}(\Psi_{\nbigm}^{\ast})^{-1}\circ
 \Psi_{\nbigm^{\lor}}$
in the derived category.
\end{lem}
\pf
Let $P\otimes\tau\otimes x,
\nbigr_X\otimes\omega_{\nbigx}\otimes\nbigm$
and 
$Q\otimes\omega\otimes y\in 
\nbigr_X\otimes\omega_{\nbigx}\otimes\nbigm^{\lor}$.
We have the following:
\begin{multline}
\bigl(
 \Psi_{\nbigm}^{\ast}\Lambda(P\otimes\tau\otimes x)
\bigr)(Q\otimes\omega\otimes y)
=\Lambda(P\otimes\tau\otimes x)
 \bigl(\Psi_{\nbigm}(Q\otimes\omega\otimes y)\bigr)
 \\
=\Phi\bigl(
 (-1)^{|\tau|\,|\omega|}
 \Psi_{\nbigm}
 (Q\otimes\omega\otimes y)(P\otimes\tau\otimes x)
 \bigr) 
=\Phi\bigl(
 C(P\otimes\tau\otimes x,Q\otimes\omega\otimes y)
 \bigr) 
\\
=(-1)^{d_X+|\tau|\,|\omega|}
 C(Q\otimes\omega\otimes y,P\otimes\tau\otimes x)
=\bigl(
(-1)^{d_X}\Psi_{\nbigm^{\lor}}(P\otimes\tau\otimes x)
\bigr)
 (Q\otimes\omega\otimes y)
\end{multline}
Thus, we are done.
\hfill\qed

\section{Dual and strict specializability
 of $\nbigr_X$-modules}

\subsection{Statement}
\label{subsection;10.9.16.10}

Let $X=X_0\times\cnum_t$.
Let $\nbigm$ be a coherent $\nbigr_X$-module
which is strictly specializable along $t$.
Assume the following:
\begin{description}
\item[(P0)]
Either 
(i) $\nbigm=\nbigm[\ast t]$ or 
(ii) $\nbigm=\nbigm[!t]$ holds.
\item[(P1)]
 $\DDD_{X_0}\Gr_a^{\Vzero}\nbigm
 \simeq
 \nbigh^0\DDD_{X_0}\Gr_a^{\Vzero}\nbigm$,
and they are strict.
\item[(P2)]
$\DDD\nbigm(\ast t)\simeq
 \nbigh^0(\DDD\nbigm(\ast t))$,
and it is strict.
\end{description}
We shall prove the following proposition
in the rest of this subsection.
\begin{prop}
\label{prop;10.9.16.11}
Under the assumptions,
the following holds:
\begin{itemize}
\item
 We have $\DDD\nbigm=\nbigh^0\DDD\nbigm$.
\item
 It is strict, and strictly specializable 
 along $t$.
\item
 We have $\DDD\nbigm=(\DDD\nbigm)[!t]$ 
 in the case (i),
 and $\DDD\nbigm=(\DDD\nbigm)[\ast t]$
 in the case (ii).
\end{itemize}
\end{prop}

For the proof of this proposition,
we do not have to care the twist with
$\lambda^{d_X}$,
and we will omit it.

\subsection{Preliminary}

We prepare some general lemmas.

\begin{lem}
\label{lem;10.9.13.4}
Let $\nbigm$ be a coherent $\nbigr_X$-module.
Assume the following.
\begin{itemize}
\item
 $\nbigm(\ast t)$ is strict
\item
$\nbigm$ is equipped with a 
 $V_0\nbigr_X$-coherent filtration
 $\Vzero(\nbigm)$
 such that 
 (i) monodromic,
 (ii) $\Gr^{\Vzero}$ is strict.
\end{itemize}
Then, $\nbigm$ is strict.
\end{lem}
\pf
Let us show that
$\Vzero_{<0}\nbigm\lrarr\nbigm(\ast t)$
is injective.
Let $f\in\Vzero_{<0}\nbigm$
be mapped to $0$ in $\nbigm(\ast t)$.
We may assume that $tf=0$.
Assume that there exists $a\in\real$ such that
$f\in \Vzero_a\setminus\Vzero_{<a}$.
Then, we obtain a non-zero element
$[f]$ in $\Gr^{\Vzero}_a$.
By using the monodromic property,
we can easily deduce that $[f]=0$,
which is a contradiction.
Hence, we obtain
$f\in \Vzero_a$ for any $a\in\real$.
We obtain
$\nbigr_X\cdot f\in \Vzero_{<0}$.
We have the decomposition
$\nbigr_X\cdot f
=\bigoplus_{j\geq 0}
 \deldel_t^j\nbigr_{X_0}\,f$
as an $\nbigr_{X_0}$-module.
Then, it is easy to check that
$\nbigr_X\cdot f$ is not finitely generated
as $V_0\nbigr_X$-module.
It contradicts with the coherence of
$\Vzero_{<0}\nbigm$.
\hfill\qed

\vspace{.1in}
We can show the following lemma
by a similar argument.

\begin{lem}
\label{lem;10.9.13.3}
Let $\nbigm$ be a coherent $\nbigr_X$-module.
Assume that (i) it is also $V_0\nbigr_X$-coherent,
(ii) $\Supp\nbigm\subset X_0$.
Then, we have $\nbigm=0$.
\hfill\qed
\end{lem}

\subsection{$\nbigr_{X_0[t]}$-modules}
\label{subsection;10.9.13.2}

Let $X_0$ be a complex manifold.
Let $t$ be a formal variable.
We set 
$\omega_{\nbigx_0[t]}:=
 \omega_{\nbigx_0}\cdot dt/\lambda$
and
$\nbigr_{X_0[t]}:=
 \nbigr_{X_0}[t]\langle\deldel_t\rangle$.
For a coherent
$\nbigr_{X_0[t]}$-module $\nbigm$,
we set
\[
 \DDD_{X_0[t]}\nbigm:=\nrhom_{\nbigr_{X_0[t]}}
 \bigl(\nbigm,\nbigr_{X_0[t]}
 \otimes\omega_{X_0[t]}^{-1}
 \bigr)[d_{X_0}+1].
\]

Let us compute some specific examples.
Let $\nbigm_0$ be a coherent
strict $\nbigr_{X_0}$-module.
We set $\nbigm_0\langle t,\deldel_t\rangle:=
 \nbigr_{X_0[t]}\otimes_{\nbigr_{X_0}}\nbigm_0$.
For $P(t,\deldel_t)\in
 \nbigo_{\cnum_{\lambda}}\langle
 t,\deldel_t
 \rangle$,
let $R\bigl(P(t,\deldel_t)\bigr)$
denote the right multiplication of
$P(t,\deldel_t)$ on $\nbigm_0\langle t,\deldel_t\rangle$.
Let $\nbign$ be a nilpotent map on $\nbigm_0$.
Let $u\in\real\times\cnum$ with
$-1\leq \paramap(\lambda_0,u)\leq 0$.
We define
$\nbigb(\nbigm_0,u,\nbign)$
as the cokernel of the injective morphism:
\[
 R\bigl(-\deldel_tt+\eigenmap(\lambda,u)\bigr)
+\nbign:
 \nbigm_0\langle t,\deldel_t\rangle
\lrarr
 \nbigm_0\langle t,\deldel_t\rangle.
\]
For $n\in\seisuu$,
we put
\[
 \Vzero_{\paramap(\lambda_0,u)+n}
 \nbigb(\nbigm_0,u,\nbign):=
 \left\{
 \begin{array}{ll}
 \Cok\bigl(t^{-n}\nbigm_0[t]
 \lrarr
 t^{-n}\nbigm_0[t]\bigr)
 & (n\leq 0) \\
 \sum_{\substack{i+j\leq n\\ i>0}}
 \deldel_t^i\Vzero_{\paramap(\lambda_0,u)+j}
 & (n>0) 
 \end{array}
 \right.
\]
For any $b\in\real$, we have
$b_0:=\max\bigl\{
 \paramap(\lambda_0,u)+n\leq b\,\big|\,
 n\in\seisuu
 \bigr\}$,
and we put
\[
 \Vzero_{b}\nbigb(\nbigm_0,u,\nbign):=
 \Vzero_{b_0}\nbigb(\nbigm_0,u,\nbign).
\]
By construction, 
we have 
$t\cdot\Vzero_{a}\subset\Vzero_{a-1}$
and $\deldel_t\cdot\Vzero_a\subset\Vzero_{a+1}$.
If $a\leq \paramap(\lambda_0,u)$,
we have
$t:\Vzero_a\simeq \Vzero_{a-1}$.
If $a\geq \paramap(\lambda_0,u)$,
we have
$\deldel_t:
 \Gr^{\Vzero}_{a}
\simeq
 \Gr^{\Vzero}_{a+1}$.

\vspace{.1in}
Let us compute $\DDD \nbigb(\nbigm_0,u,\nbign)$.
We assume that
$\DDD_{X_0}\nbigm_0\simeq
 \nbigh^0\DDD_{X_0}\nbigm_0$,
and that it is strict,
for simplicity.
We consider it locally,
and we take an $\nbigr_{X_0}$-free resolution
$\nbigp^{\bullet}$ of $\nbigm_0$.
We have 
$\nbign^{\bullet}:\nbigp^{\bullet}\lrarr\nbigp^{\bullet}$
which induces $\nbign$.
Then, $\nbigb(\nbigm_0,u,\nbign)$
is naturally quasi-isomorphic
to the complex associated to the following
double complex:
\[
 R\bigl(-\deldel_tt+\eigenmap(\lambda,u)\bigr)
+\nbign^{\bullet}:
 \nbigp^{\bullet}\langle t,\deldel_t\rangle
\lrarr
 \nbigp^{\bullet}\langle t,\deldel_t\rangle
\]
As a right $\nbigr_{X_0[t]}$-complex,
$\DDD\nbigb(\nbigm_0,u,\nbign)$ is quasi-isomorphic to
the complex associated to the following double complex:
{\small
\[
\begin{CD}
 \nhom_{\nbigr_{X_0}}
 \bigl(\nbigp^{\bullet}(\nbigm_0),\nbigr_{X_0}\bigr)
 \langle t,\deldel_t\rangle
@>{-\deldel_tt+\eigenmap(\lambda,u)+\nbign^{\lor}}>>
 \nhom_{\nbigr_{X_0}}
 \bigl(\nbigp^{\bullet}(\nbigm_0),\nbigr_{X_0}\bigr)
 \langle t,\deldel_t\rangle
\end{CD}
\]
}
It is quasi isomorphic to the complex 
$-\deldel_tt+\eigenmap(\lambda,u)+\nbign^{\lor}:
 \DDD_{X_0}\nbigm_0\langle t,\deldel_t
 \rangle
\lrarr
 \DDD_{X_0}\nbigm_0\langle t,\deldel_t
 \rangle$.
As a left $\nbigr_{X_0}\langle t,\deldel_T\rangle$-module,
it is quasi isomorphic to
$t\deldel_t+\eigenmap(\lambda,u)+\nbign^{\lor}:
 \DDD_{X_0}\nbigm_0\langle t,\deldel_t \rangle
\lrarr
 \DDD_{X_0}\nbigm_0\langle t,\deldel_t \rangle$.
Note 
$t\deldel_t+\eigenmap(\lambda,u)
=\deldel_tt+\eigenmap(\lambda,u+\vecdelta)$.
Hence, we obtain
\[
 \DDD\nbigb\bigl(\nbigm_0,u,\nbign\bigr)
\simeq
 \nbigb\bigl(\DDD_{X_0}\nbigm_0,
 -u-\vecdelta,-\nbign^{\lor}\bigr)
\]

\subsection{Filtered free module}

Let $X=X_0\times\cnum_t$.
Recall that $V_0\nbigr_X\subset\nbigr_X$ 
is the sheaf of subalgebras
generated by $\nbigr_{X_0}$
and $t\deldel_t$ over $\nbigo_{\nbigx}$.
Let $\nbigp$ be a free $\nbigr_X$-module
with a generator $e$.
Let $a$ be a real number.
We have a filtration $V$ given as follows.
For $n\leq 0$,
we put $V_{a+n}\nbigp:=
 t^{-n}V_0\nbigr_X\cdot e$.
For $n>0$,
we put
$V_{a+n}\nbigp:=
 \sum_{i+j\leq n} \deldel_t^iV_{j}\nbigp$.
For $b\not\in a+\seisuu$,
we put 
$V_{b}\nbigp:=V_{c}\nbigp$,
where $c:=\max\{a+n\,|\,a+n\leq b,n\in\seisuu\}$.
Such a filtered $\nbigr$-module
is denoted by $\nbigr_X(e,a)$.
A filtered $\nbigr$-module is called
$(\nbigr,V)$-free,
if it is isomorphic to a direct sum
$\bigoplus_i \nbigr_X(e_i,a_i)$.
By the identification
$\nbigr_X=\nbigr_X(1,0)$,
$\nbigr_X$ is a filtered ring.

\vspace{.1in}

We have similar notions  for $\nbigr_{X_0[t]}$-modules.
Let $\nbigq$ be a free $\nbigr_{X_0[t]}$-module
with a generator $e$.
Let $a$ be a real number.
We have a filtration $V$ given as follows.
For $n\leq 0$,
we put $V_{a+n}\nbigq:=
 t^{-n}V_0\nbigr_{X_0[t]}\cdot e$.
For $n>0$,
we put
$V_{a+n}\nbigq:=
 \sum_{i+j\leq n} \deldel_t^iV_{j}\nbigq$.
For $b\not\in a+\seisuu$,
we put 
$V_{b}\nbigq:=V_{c}\nbigq$,
where $c:=\max\{a+n\,|\,a+n\leq b,n\in\seisuu\}$.
Such a filtered $\nbigr_{X_0[t]}$-module
is denoted by $\nbigr_{X_0[t]}(e,a)$.
A filtered $\nbigr_{X_0[t]}$-module
is called $(\nbigr_{X_0[t]},V)$-free,
if it is isomorphic to a direct sum
$\bigoplus_i \nbigr_{X_0[t]}(e_i,a_i)$.

\vspace{.1in}

Note that $\Gr^V\nbigr_X$ is naturally isomorphic to
$\nbigr_{X_0[t]}$,
and we have a natural isomorphism
$\Gr^V\nbigr_{X}(e,a)\simeq \nbigr_{X_0[t]}(e,a)$.
If a filtered $\nbigr_X$-module
$(\nbigl,V)$ is $(\nbigr_X,V)$-free,
$\Gr^V(\nbigl,V)$ is 
$(\nbigr_{X_0[t]},V)$-free.

\subsubsection{Strictness}

Let $(\nbigp_1,V)\stackrel{\varphi_1}{\lrarr}
 (\nbigp_2,V)\stackrel{\varphi_2}{\lrarr}(\nbigp_3,V)$
be a complex of coherent $(\nbigr_X,V)$-free modules.

\begin{lem}
\label{lem;10.9.16.20}
Assume that
$\Gr^V\nbigp_1\lrarr\Gr^V\nbigp_2\lrarr
\Gr^V\nbigp_3$ is exact.
Then, $\varphi_i$ are strict
with respect to $V$.
\end{lem}
\pf
The strictness of $\varphi_2$ is easy.
Let us argue the strictness of $\varphi_1$.
We fix a sufficiently small $b<0$.
Let $a\in\real$.
Let us observe that there exists
$N(a)>0$ such that
$\Image\varphi_1\cap 
 V_{b-N(a)}\nbigp_2
\subset
 \varphi_1\bigl(V_{a-1}\nbigp_1\bigr)$.
We obtain a $V_0\nbigr_X$-coherent submodule
$\varphi_1^{-1}(V_{b}\nbigp_2)$.
If we take sufficiently large $N(a)$,
we have 
$\varphi_1^{-1}(V_{b}\nbigp_2)
\subset V_{a+N(a)-1}\nbigp_1$.
Then, if $b$ is sufficiently small,
we have 
$t^{N(a)}V_{b}\nbigp_2=V_{b-N(a)}\nbigp_2$.
Hence, we have 
$\varphi_1^{-1}(V_{b-N(a)}\nbigp_2)
\subset V_{a-1}\nbigp_1$.

Take $f\in V_{a}\nbigp_1$
such that $\varphi_1(f)\in V_{b}\nbigp_2$
for some $b<a$.
By using the exactness of
$\Gr(\nbigp_{\bullet})$,
we can find $g\in V_{<a}\nbigp_1$
such that
$\varphi_1(f-g)\in V_{b-N(a)-1}\nbigp_2$.
Then, we can find $h\in V_{a-1}\nbigp_1$
such that $\varphi_1(f-g-h)=0$,
i.e., $\varphi_1(f)=\varphi_1(g+h)$.
By an easy inductive argument,
we obtain that $\varphi_1(f)\in
 \varphi_1\bigl(V_b(\nbigp_1)\bigr)$.
\hfill\qed

\subsubsection{Dual}
We consider the dual of $\nbigr_X(e,a)$,
i.e.,
\[
 \DDD\nbigr_X(e,a):=
 \nhom_{\nbigr_X}(\nbigr_X(e,a),\nbigr_X
 \otimes\omega_{\nbigx}).
\]
Let $e^{\lor}\in\DDD\nbigr_X(e,a)$
be given by $e^{\lor}(e)=1$.
Then, $\DDD\nbigr_X(e,a)\simeq
 \nbigr_X\cdot e^{\lor}$.
A filtration of $\DDD\nbigr_X(e,a)$ is induced by
the filtration of $\nbigr_X(e^{\lor},-a-1)$.
Similarly, we have the filtration
of $\DDD\nbigr_{X_0[t]}(e,a)$
induced by $\DDD\nbigr_{X_0[t]}(e,a)\simeq
 \nbigr_{X_0[t]}(e^{\lor},-a-1)$.
Hence, we have naturally induced filtrations
on the dual of free $(\nbigr_X,V)$-modules
or $(\nbigr_{X_0[t]},V)$.

\subsection{A filtered free resolution}
\label{subsection;10.9.13.1}

Let $X=X_0\times\cnum_t$.
Let $\nbigm$ be a strictly specializable
$\nbigr_X$-module along $t$.
For simplicity, we assume either 
(i) $\nbigm=\nbigm[\ast t]$,
or (ii) $\nbigm=\nbigm[!t]$.
Locally, we shall construct a complex of
coherent $(\nbigr_X,V)$-free modules
{\small
\[
\cdots\lrarr
 (\nbigp_{\ell+1}(\nbigm),V)
\lrarr(\nbigp_{\ell}(\nbigm),V)
\lrarr \cdots\cdots\lrarr(\nbigp_1(\nbigm),V)
 \lrarr(\nbigp_0(\nbigm),V) 
\]
}
with the following property:
\begin{itemize}
\item
$\nbigp_{\bullet}(\nbigm)$
is a free resolution of $\nbigm$.
\item
The morphisms
$\nbigp_{\ell+1}(\nbigm)\lrarr
 \nbigp_{\ell}(\nbigm)$
and $\nbigp_{0}(\nbigm)\lrarr\nbigm$
are strict with respect to the filtrations.
\item
In particular,
$\Gr^V(\nbigp_{\bullet}(\nbigm))$
gives a free $\nbigr_{X_0[t]}$-resolution
of $\Gr^V(\nbigm)$.
\end{itemize}

First, we construct
a $(\nbigr_X,V)$-free module $\nbigp_0(\nbigm)$
with a surjection $\nbigp_0(\nbigm)\lrarr \nbigm$.

\subsubsection{The case (i)}

Let us consider the case (i)
in the condition (P0).
We put
$a_0:=\min\bigl\{
 -1<a\leq 0\,\big|\,
 \Gr^{\Vzero}_a\nbigm\neq 0
  \bigr\}$.
We take a generator 
$\vecetilde_{a_0}=(\etilde_{a_0,i})$
of $\Vzero_{a_0}\nbigm$.
We take generators $\vece_{a}$ of
$\Gr^{\Vzero}_a(\nbigm)$ for $a_0<a\leq 0$,
and we take lifts $\vecetilde_a=(\etilde_{a,i})$
to $\Vzero_a(\nbigm)$.
Let 
$\nbigp_0\nbigm$ be the $(\nbigr_{X},V)$-free module
generated by $\vecetilde_a$ $(a_0\leq a\leq 0)$,
where $a$ is associated to each $\etilde_{a,i}$.
We have a naturally defined filtered morphism
$(\nbigp_0,V)\lrarr(\nbigm,\Vzero)$.

\subsubsection{The case (ii)}

We take a generator 
$\vecetilde_{-1}=(\etilde_{-1,i})$
of $\Vzero_{-1}\nbigm$.
We take generators $\vece_{a}$ of
$\Gr^{\Vzero}_a(\nbigm)$ for $-1<a<0$,
and we take lifts $\vecetilde_a=(\etilde_{a,i})$
to $\Vzero_a(\nbigm)$.
Let 
$\nbigp_0\nbigm$ be the $(\nbigr_{X},V)$-free module
generated by $\vecetilde_a$ $(-1\leq a< 0)$,
where $a$ is associated to each $\etilde_{a,i}$.
We have the naturally defined filtered morphism
$(\nbigp_0,V)\lrarr(\nbigm,\Vzero)$.

\vspace{.1in}

By construction,
for each $a\in\real$,
the morphism
$V_a\nbigp_0\lrarr\Vzero_a\nbigm$
is surjective.
By construction,
the induced morphism
$\Gr^V\nbigp_0\lrarr
 \Gr^{\Vzero}\nbigm$ is surjective.
Let $\nbigk_0(\nbigm)$
denote the kernel of 
$\nbigp_0(\nbigm)\lrarr\nbigm$.
It is equipped with a naturally induced filtration
$V$.
By construction, we have
$t\cdot V_a\nbigk_0(\nbigm)\subset
 V_{a-1}\nbigk_0(\nbigm)$
and
$\deldel_t\cdot V_a\nbigk_0(\nbigm)\subset
 V_{a+1}\nbigk_0(\nbigm)$.
\begin{itemize}
\item
We have
$t:V_{b}(\nbigk_0(\nbigm))
 \simeq V_{b-1}(\nbigk_0(\nbigm))$
for $b\leq 0$ in the case (i),
or for $b<0$ in the case (ii).
\item
We have the exact sequence
$0\lrarr \Gr^{V}\nbigk_0\lrarr
 \Gr^V\nbigp_0\lrarr \Gr^{\Vzero}\nbigm\lrarr 0$.
In particular, we have
$\deldel_t:\Gr^{V}_b\nbigk_0\simeq
 \Gr^{V}_{b+1}\nbigk_0$
for $b>-1$ in the case (i),
or for $b\geq -1$ in the case (ii).
\end{itemize}

Inductively, for any $\ell\geq 0$,
we can construct
filtered $\nbigr_X$-modules
$(\nbigp_{\ell}(\nbigm),V)$ and
$(\nbigk_{\ell}(\nbigm),V)$
with exact sequences
$0\lrarr\nbigk_{\ell+1}\lrarr 
 \nbigp_{\ell+1}\lrarr\nbigk_{\ell}\lrarr 0$,
such that
\begin{itemize}
\item
$(\nbigp_{\ell}(\nbigm),V)$
are free $(\nbigr_X,V)$-modules.
\item
For both $\nbigp_{\ell}$ and $\nbigk_{\ell}$,
the morphisms
$t:V_{b}\simeq V_{b-1}$ 
$(b\leq 0)$ and
$\deldel_t:\Gr^{V}_b\simeq\Gr^{V}_{b+1}$
$(b>-1)$ are isomorphisms
in the case (i).
The morphisms
$t:V_{b}\simeq V_{b-1}$ 
$(b< 0)$ and
$\deldel_t:\Gr^{V}_b\simeq\Gr^{V}_{b+1}$
$(b\geq -1)$ are isomorphisms
in the case (ii).
\item 
 The morphisms 
$\nbigk_{\ell+1}\lrarr 
 \nbigp_{\ell+1}$
and
$\nbigp_{\ell+1}\lrarr\nbigk_{\ell}$
are strict with respect to 
the filtrations $V$.
In particular,
the induced sequence
$0\lrarr\Gr^V\nbigk_{\ell+1}\lrarr 
 \Gr^V\nbigp_{\ell+1}\lrarr
 \Gr^V\nbigk_{\ell}\lrarr 0$
is exact.
\end{itemize}

We obtain a $(\nbigr_X,V)$-free resolution
$\nbigp_{\bullet}(\nbigm)$ of $\nbigm$.
We have a natural isomorphism
\[
\Gr^V\nhom_{\nbigr_X}\bigl(\nbigp_{\bullet},\nbigr_X\bigr)
\simeq
 \nhom_{\Gr^V\nbigr_X}\bigl(\Gr^V\nbigp_{\bullet},
 \Gr^V\nbigr_X\bigr).
\]
It naturally gives a 
$(\nbigr_{X_0[t]},V)$-free resolution
of $\Gr^V\nbigm$.

\subsection{Proof of
Proposition \ref{prop;10.9.16.11}}

Let $\nbigm$ be as in
\S\ref{subsection;10.9.16.10}.
Let $\nbigp_{\bullet}(\nbigm)$ be 
a $(\nbigr_X,V)$-free resolution of $\nbigm$
as in \S\ref{subsection;10.9.13.1}.
Let us study the dual of $\nbigm$
by using $\nbigp_{\bullet}(\nbigm)$.
The complexes
\[
 \nbigc:=\nhom_{\nbigr_X}\bigl(
 \nbigp_{\bullet}(\nbigm),\nbigr_X\otimes\omega_{\nbigx}
 \bigr)
\]
and
\[
\Gr^V(\nbigc)=\nhom_{\Gr^V\nbigr_X}\bigl(
 \Gr^V\nbigp_{\bullet}(\nbigm),\Gr^V\nbigr_X
 \otimes\omega_{\nbigx[t]}
 \bigr)
\]
express $\DDD\nbigm$
and $\DDD \Gr^V\nbigm$,
respectively.

\begin{lem}
\mbox{{}}\label{lem;10.9.16.21}
Assume $(P0)$ and $(P1)$.
Then, the following holds.
\begin{itemize}
\item
 $\nbigh^0\DDD\Gr^V(\nbigm)\simeq
 \DDD\Gr^V(\nbigm)$,
 and it is strict.
\item
 $\Vzero$ 
 on $\nbigh^0\DDD\Gr^V(\nbigm)$
 is monodromic,
 and we have 
 $t\cdot \Vzero_a=\Vzero_{a-1}$ for $a<0$
 and
 $\deldel_t:\Gr^{\Vzero}_a\simeq\Gr^{\Vzero}_{a+1}$
 for $a>-1$.
\item
 In the case (i), 
 the morphism
 $\deldel_t:
 \Gr^{\Vzero}_{-1}\DDD\Gr^V(\nbigm)
 \lrarr
 \Gr^{\Vzero}_0\DDD\Gr^V(\nbigm)$
 is an isomorphism.
 In the case (ii),
 the morphism 
 $t: \Gr^{\Vzero}_{0}\DDD\Gr^V(\nbigm)
 \lrarr
 \Gr^{\Vzero}_{-1}\DDD\Gr^V(\nbigm)$
 is an isomorphism.
\end{itemize}
\end{lem}
\pf
It follows from the computation in
\S\ref{subsection;10.9.13.2}.
\hfill\qed

\vspace{.1in}

Let us finish the proof of Proposition \ref{prop;10.9.16.11}.
By Lemma \ref{lem;10.9.16.20} and 
Lemma \ref{lem;10.9.16.21},
the morphisms in the complex $\nbigc$
is strict with respect to the filtration $V$.
Hence, we have the commutativity
$\Gr^V\nbigh^i\nbigc\simeq\nbigh^i\Gr^V\nbigc$.
According to Lemma \ref{lem;10.9.16.21},
we have 
$\Gr^{V}\nbigh^i\nbigc=0$
unless $i=0$.
We also have 
$\bigl(\nbigh^i\nbigc\bigr)(\ast t)=0$
unless $i=0$.
Hence, 
we obtain $\nbigh^i\DDD\nbigm=\nbigh^i\nbigc=0$
unless $i=0$,
by Lemma \ref{lem;10.9.13.3}.

By Lemma \ref{lem;10.9.16.21},
$\Gr^V\nbigh^0\DDD\nbigm$ is strict
and monodromic.
Hence, 
by using Lemma \ref{lem;10.9.13.4},
we obtain that 
$\nbigh^0\DDD\nbigm$ is strict
and strictly specializable along $t$.
Moreover the induced filtration $\Vzero$
gives a $V$-filtration.
By the second claim of Lemma \ref{lem;10.9.16.21},
we obtain $(\DDD\nbigm)[!t]=\DDD\nbigm$ 
in the case (i),
or $(\DDD\nbigm)[\ast t]=\DDD\nbigm$
in the case (ii).
\hfill\qed

\section{Dual of mixed twistor $D$-module}

\subsection{Statements}

Let us consider the dual 
for mixed twistor $D$-modules.
We will prove the following theorem in
\S\ref{subsection;10.9.17.12}--\S\ref{subsection;11.1.29.3}.

\begin{thm}
\label{thm;10.9.13.11}
Let $\nbigt\in \MTM(X)$,
and $\nbigt=(\nbigm_1,\nbigm_2,C)$
as an $\nbigr_X$-triple.
\begin{itemize}
\item
 $\nbigh^0(\DDD\nbigm_i)\simeq\DDD\nbigm_i$,
 and it is strict.
 In particular, it is equipped with a filtration
 induced by the weight filtration of $\nbigt$.
\item
 We have a unique sesqui-linear pairing $\DDD C$ of
 $\DDD\nbigm_1$ and
 $\DDD\nbigm_2$ such that 
\[
  (\DDD C)_{|\nbigxlambda}=
 \DDD \bigl(C_{|\nbigxlambda}\bigr)
\]
 for generic $\lambda\in\vecS$,
where the latter is defined for 
non-degenerate
hermitian pairings of holonomic $D$-modules
in Theorem {\rm\ref{thm;10.9.14.2}}.
\item
The $\nbigr$-triple
$\DDD\nbigt:=
 \bigl(\DDD\nbigm_1,\DDD\nbigm_2,\DDD C\bigr)$
with the naturally induced filtration is a mixed twistor
$D$-module.
\item
The dual $\DDD$ gives a contravariant functor
on $\MTM(X)$,
and $\DDD\circ\DDD\simeq\id$.
\item
If $\nbigt\in\MTMint(X)$,
we naturally have
$\DDD\nbigt\in\MTMint(X)$.
\end{itemize}
\end{thm}

Before going to the proof of Theorem \ref{thm;10.9.13.11},
we give some consequences.

\vspace{.1in}

For 
$\nbigt\!=\!(\nbigm_1,\nbigm_2,C)\in\MTM(X)$
and $\ell\in\seisuu$,
we set 
$\mu_{\ell}(\!\nbigt\!)\!:=\!(\nbigm_1,\nbigm_2,(-1)^{\ell}C)$.
\index{$\nbigr$-triple $\mu_{\ell}(\nbigt)$}
Once we know Theorem \ref{thm;10.9.13.11},
we obtain the following compatibility of the dual
and the push-forward,
from Theorem \ref{thm;10.9.14.3} below.
\begin{thm}
\label{thm;10.9.17.1}
Let $\nbigt\in\MTM(X)$.
Let $F:X\lrarr Y$ be a projective morphism.
Then, 
the natural isomorphisms of the underlying
$\nbigr_Y$-modules 
(Lemma {\rm\ref{lem;10.9.17.10}})
give a natural isomorphism
$\DDD F^j_{\dagger}\nbigt
 \simeq
 \mu_jF^{-j}_{\dagger}(\DDD\nbigt)$
in $\MTM(Y)$.
\hfill\qed
\end{thm}
In other words,
$\bigl(
 \epsilon(-j)\varphi_1^{-1},\epsilon(j)\varphi_2
 \bigr)$
gives an isomorphism
$ F^{j}_{\dagger}(\DDD\nbigt)
\simeq
\DDD F^{-j}_{\dagger}\nbigt$,
where 
$\varphi_1:F_{\dagger}^{-j}\DDD_X\nbigm_1
 \simeq
 \DDD_YF_{\dagger}^{j}\nbigm_1$
and 
$\varphi_2:F_{\dagger}^{j}\DDD_X\nbigm_2
 \simeq
 \DDD_YF_{\dagger}^{-j}\nbigm_2$
are the isomorphisms in Lemma \ref{lem;10.9.17.10},
and $\epsilon(m):=(-1)^{m(m-1)/2}$ for integers $m$.

\begin{prop}
We have natural isomorphisms
\[
\DDD(\nbigt^{\ast})
\simeq
 (\DDD\nbigt)^{\ast},
\quad
 \DDD j^{\ast}\nbigt
\simeq
 j^{\ast}\DDD\nbigt
\]
in $\MTM(X)$.
\end{prop}
\pf
We have natural isomorphisms
for the underlying filtered $\nbigr$-modules.
We have only to show the compatibility
of the pairings
$\DDD (C^{\ast})
=(\DDD C)^{\ast}$
and
$j^{\ast}\DDD C
=\DDD j^{\ast}C$.
Both of them can be reduced to 
the claims for the $D$-modules,
and easy to check.
(See \S\ref{subsection;11.4.6.1}.)
\hfill\qed

\begin{prop}
\label{prop;11.1.29.1}
Let $H$ be a hypersurface of $X$.
We have natural isomorphisms
$\DDD(\nbigt[\ast H])
\simeq
 (\DDD\nbigt)[! H]$
and 
$\DDD(\nbigt[! H])
\simeq
 (\DDD\nbigt)[\ast H]$
in $\MTM(X)$.
\end{prop}
\pf
Let us observe that
$\DDD(\nbigt[\ast H])
\simeq
 \DDD(\nbigt)[!H]$.
We have the isomorphisms 
of the underlying $\nbigr$-modules
(Proposition \ref{prop;10.9.16.11}).
We have the compatibility of the pairings
(see Proposition \ref{prop;10.11.17.10} below).
Hence, 
we obtain the isomorphism
$\DDD(\nbigt[\ast H])
\simeq
 \DDD(\nbigt)[!H]$
as $\nbigr$-triples.
It is compatible with the naively induced filtrations.
Then, the claim follows from Lemma \ref{lem;10.9.2.3}.
\hfill\qed

\begin{prop}
Let $\nbigv\in\MTS^{\adm}(X,H)$.
Let $\nbigt\in\MTM(X)$.
Then, we have natural isomorphisms
$\DDD\bigl((\nbigt\otimes\nbigv)[!H]\bigr)
\simeq
\bigl(\DDD\nbigt\otimes \nbigv\bigr)[\ast H]$
and
$\DDD\bigl((\nbigt\otimes\nbigv)[\ast H]\bigr)
\simeq
\bigl(\DDD\nbigt\otimes \nbigv\bigr)[!H]$.
\end{prop}
\pf
We have only to consider the case $H=\emptyset$.
(See Proposition \ref{prop;11.1.29.1}.)
We have natural isomorphisms
as in \S\ref{subsection;11.1.29.2}.
We have the compatibility of the pairings
in Lemma \ref{lem;10.12.27.23} below.
\hfill\qed

\begin{cor}
We have natural isomorphisms
$\DDD(\Pi^{a,b}_{\ast!}\nbigt)
\!\simeq\!
 \Pi^{-b+1,-a+1}_{\ast!}(\DDD\nbigt)$.
In particular, we have
\[
\DDD\psi_g^{(a)}(\nbigt)
\simeq
 \psi_g^{(-a+1)}\DDD\nbigt,
\quad\quad
\DDD\Xi_g^{(a)}(\nbigt)
\simeq
 \Xi_g^{(-a)}\DDD\nbigt.
\]
We also obtain
$\DDD\phi_g^{(a)}(\nbigt)
\simeq
 \phi_g^{(-a)}\DDD\nbigt$.
\hfill\qed
\end{cor}

\subsection{Relative monodromy filtration}
\label{subsection;10.9.17.12}

For the proof of the theorems,
we shall use an induction on the support of
the dimensions.
Let $A(n)$ denote the claim of Theorem
\ref{thm;10.9.13.11}
in the case $\dim\Supp\nbigt\leq n$.

Assume $A(n)$.
Let $\nbigt\in\MTM(X)$
with $\dim\Supp\nbigt\leq n$.
Let $W$ be the weight filtration of $\nbigt$,
and let $L$ be a filtration of $\nbigt$ in $\MTM(X)$.
Let 
 $\nbign:(\nbigt,W,L)\lrarr
 \bigl((\nbigt,W)\otimes\newTate(-1),L\bigr)$ 
be a morphism such that $W=M(\nbign;L)$.
We have $\DDD(\nbigt,W)=(\DDD\nbigt,\DDD W)$
in $\MTM(X)$.
It is equipped with the induced filtration $\DDD L$
and the induced morphism 
$\DDD\nbign:\DDD(\nbigt,W,L)
\lrarr
 \bigl(\DDD(\nbigt,W)\otimes\newTate(-1),
 \DDD L\bigr)$.
We can show the following lemma
by a standard argument.
\begin{lem}
We have $\DDD W=M\bigl(
 \DDD\nbign;\DDD L
 \bigr)$.
\end{lem}
\pf
Let us consider the case $\Gr^L_j(\nbigt)=0$
unless $j=w$.
We may assume that $w=0$.
We naturally have
$\Gr^{\DDD W}_j\DDD\nbigt\simeq
 \DDD\Gr^{W}_{-j}\nbigt$,
and the morphism
$\DDD\nbign^j:\Gr^{\DDD W}_{j}\DDD\nbigt\lrarr
 \Gr^{\DDD W}_{-j}\DDD\nbigt$ is the dual of
$\nbign^j:\Gr^{W}_{j}\nbigt
 \lrarr\Gr^{W}_{-j}\nbigt$.
Hence, it is an isomorphism,
and $\DDD W$ on $\DDD\nbigt$
is the monodromy filtration of $\DDD\nbign$.

The induced morphism
$L_j\DDD \nbigt\lrarr \Gr^L_j\DDD \nbigt$
is strict with respect to $\DDD W$.
Hence, we can deduce that
$\DDD W$ is the relative monodromy filtration
with respect to $L$.
\hfill\qed

\subsection{Dual of smooth $\nbigr$-triple}

We give a remark on the signature.
Let $\nbigv=(\nbigv_1,\nbigv_2,C)\in\VTS(X)$.
Let 
$\nbigv^{\lor}=(\nbigv_1^{\lor},\nbigv_2^{\lor},C^{\lor})$
be the dual in $\VTS(X)$.
We have the isomorphisms 
$\Phi_i:
 \nbigv_i^{\lor}\simeq
 \DDD\nbigv_i$ $(i=1,2)$
given by 
$\Phi_i(v_i)=\lambda^{d_X}\,v_i$.
(See \S\ref{subsection;13.4.1.10}.)
We have the induced pairing $C_1$ of 
$\DDD\nbigv_1$ and $\DDD\nbigv_2$.
\begin{lem}
\label{lem;10.11.18.2}
We have
$C_1^{\lambda}=
\DDD(C^{\lambda})$
for each $\lambda\in\vecS$,
where $\DDD(C^{\lambda})$ denote the dual functor
for non-degenerate $D$-triples.
(See {\rm\S\ref{section;11.4.6.10}}.)
\end{lem}
\pf
Let $\lambda_0\in\vecS$.
We have
$C_1\bigl(\lambda^{d_X}f,
 \overline{\lambda^{d_X}g}\bigr)
=C^{\lor}(f,\overline{g})$
by construction.
On the other side,
we have
\begin{multline}
 \DDD(C^{\lambda_0})\bigl(
 (\lambda^{d_X}f)_{\lambda_0},
 (\overline{\lambda^{d_X}g})_{-\lambda_0}
 \bigr) 
= \\
 (-1)^{d_X}
 \lambda_0^{d_X}
 (-\lambda_0)^{-d_X}
 (C^{\lambda_0})^{\lor}\bigl(f_{\lambda_0},
 \overline{g}_{-\lambda_0}\bigr)
=C^{\lor}(f,\overline{g})_{\lambda_0}
\end{multline}
(See Example \ref{example;13.4.1.11}.)
Then, the claim of the lemma follows.
\hfill\qed

\subsection{Dual of canonical prolongation as $\nbigr$-triples}
\label{subsection;10.9.16.2}

Let $X$ be any $n$-dimensional complex manifold
with a normal crossing hypersurface $D$.
Let $\nbigv\in\MTS^{\adm}(X,D)$.
Let $\nbigv^{\lor}\in\MTS^{\adm}(X,D)$ be its dual.
We have the description
$\nbigv=(\nbigv_1,\nbigv_2,C)$ 
and $\nbigv^{\lor}=
 (\nbigv_1^{\lor},\nbigv_2^{\lor},C^{\lor})$
as smooth $\nbigr_{X(\ast D)}$-triples.
We put $\nbigv_{i\star}:=\nbigv_{i}[\star D]$
for $\star=\ast,!$.
We use the symbol $\nbigv^{\lor}_{i\star}$
in similar meanings.
We have the $\nbigr_X$-triples
$\nbigv_{\ast}
=(\nbigv_{1!},\nbigv_{2\ast},C_{\ast})$
and $\nbigv_{!}
=(\nbigv_{1\ast},\nbigv_{2!},C_{!})$.
We also have the $\nbigr_X$-triples
$\nbigv^{\lor}_{\ast}
=\bigl(\nbigv^{\lor}_{1!},
 \nbigv^{\lor}_{2\ast},C^{\lor}_{\ast}\bigr)$
and $\nbigv^{\lor}_{!}
=\bigl(\nbigv^{\lor}_{1\ast},
 \nbigv^{\lor}_{2!},C^{\lor}_{!}\bigr)$.

\begin{lem}
\label{lem;10.9.17.2}
Assume $A(n-1)$.
Then, we have natural isomorphisms
of $\nbigr_X$-modules
$\DDD\nbigv_{i!}\simeq
 \lambda^{d_X}\nbigv^{\lor}_{i\ast}$
and 
$\DDD\nbigv_{i\ast}\simeq
 \lambda^{d_X}\nbigv^{\lor}_{i!}$.
\end{lem}
\pf
We may assume that there exists a holomorphic 
function $g$ on $X$ such that $g^{-1}(0)=D$.
Let $\iota_g:X\lrarr X\times\cnum_t$ be the graph.
By the assumption $A(n-1)$,
the condition $(P1)$
in \S\ref{subsection;10.9.16.10}
is satisfied.
The other conditions $(P0)$ and $(P2)$
are also satisfied.
Hence, we obtain that
$\DDD\nbigv_{i\ast}
=\nbigh^0\DDD\nbigv_{i\ast}$
by Proposition \ref{prop;10.9.16.11}.
Moreover,
it is strictly specializable along $g$
and we have
$\DDD(\nbigv_{i\ast})[!g]=\DDD\nbigv_{i\ast}$.
Because 
$\DDD(\nbigv_{i\ast})(\ast g)\simeq
 \lambda^{d_X}\nbigv_i^{\lor}$ naturally,
we obtain
$\DDD(\nbigv_{i\ast})\simeq
 \lambda^{d_X}\nbigv_{i!}^{\lor}$.
We obtain the other isomorphism
in a similar way.
\hfill\qed

\begin{lem}
$\nbigv_{\star}$ $(\star=\ast,!)$ have their duals
as $\nbigr$-triples,
and the isomorphisms in Lemma
{\rm\ref{lem;10.9.17.2}}
induce
$\DDD\nbigv_!\simeq\nbigv^{\lor}_{\ast}$
and $\DDD\nbigv_{\ast}\simeq\nbigv^{\lor}_!$.
\end{lem}
\pf
By Lemma \ref{lem;10.11.18.2},
we have the coincidence of pairings
on $\{\lambda\}\times (X\setminus D)$
$(\lambda\in\vecS)$
under the isomorphisms of 
the underlying $\nbigr$-modules.
Then, it is extended to
$\DDD C_{!|\nbigxlambda}\simeq
 C^{\lor}_{\ast|\nbigxlambda}$
for generic $\lambda$.
It means 
$\DDD\nbigv_{!}$ exists and
$\DDD\nbigv_!\simeq\nbigv^{\lor}_{\ast}$.
We can show the claim for the other
in a similar way.
\hfill\qed

\vspace{.1in}
We shall argue the comparison
of the weight filtrations later.

\subsection{Dual of minimal extension
in the pure case}

Let $(X,D)$ and $\nbigv$ be as in 
\S\ref{subsection;10.9.16.2}.
Let us assume that $\nbigv$ is pure.
We set $\nbigv_{!\ast}:=\Image
 \bigl(\nbigv_!\lrarr\nbigv_{\ast}\bigr)$,
which is a polarizable wild pure twistor $D$-module.
The underlying $\nbigr_X$-modules
are denoted by $\nbigv_{i!\ast}$ $(i=1,2)$.
Similarly, 
we obtain a polarizable wild pure
twistor $D$-module 
$\nbigv^{\lor}_{!\ast}$
with the underlying $\nbigr_X$-modules
$\nbigv^{\lor}_{i!\ast}$.

\begin{lem}
We have natural isomorphisms 
$\DDD\nbigv_{i!\ast}
\simeq \nbigv^{\lor}_{i!\ast}$.
The dual of $\nbigv_{!\ast}$ as $\nbigr_X$-triples
exists, and it is isomorphic to
the $\nbigr_X$-triple 
$\nbigv^{\lor}_{!\ast}$.
\end{lem}
\pf
Let $K_i$ denote the kernel of
$\nbigv_{i!}\lrarr \nbigv_{i!\ast}$,
and let $C_i$  denote the cokernel of
$\nbigv_{i!\ast}\lrarr \nbigv_{i\ast}$.
By the assumption $A(n-1)$,
we have
$\nbigh^0\DDD K_i=\DDD K_i$
and $\nbigh^0\DDD C_i=\DDD C_i$.

From the exact sequence
$0\lrarr K_i\lrarr \nbigv_{i!}\lrarr
 \nbigv_{i!\ast}\lrarr 0$,
we obtain 
$\nbigh^j\DDD\nbigv_{i!\ast}=0$
unless $j=0,-1$.
From the exact sequence
$0\lrarr \nbigv_{i!}\lrarr
 \nbigv_{i!\ast}\lrarr C_{i}\lrarr 0$,
we obtain 
$\nbigh^j\DDD\nbigv_{i!\ast}=0$
unless $j=0,1$.
Hence, we obtain
$\nbigh^j\DDD\nbigv_{i!\ast}=0$
unless $j=0$.
We also obtain that
$\nbigh^0\DDD\nbigv_{i!\ast}$
is the image of
$\nbigh^0\DDD\nbigv_{i\ast}
\lrarr
 \nbigh^0\DDD\nbigv_{i!}$.
Therefore,
$\DDD\nbigv_{i!\ast}
\simeq
 \nbigv^{\lor}_{i!\ast}$.
By using the uniqueness of pairing,
we obtain
$\DDD\nbigv_{!\ast}\simeq
 \nbigv^{\lor}_{!\ast}$.
\hfill\qed

\vspace{.1in}
We have immediate consequences
on the dual of the filtered $\nbigr$-modules
underlying good mixed twistor $D$-modules.

\begin{cor}
Let $\nbigt\in\MTM^{\good}(X,D)$.
Let $\nbigm_i$ $(i=1,2)$
be the underlying $\nbigr_X$-modules
of $\nbigt$.
Then, we have 
$\nbigh^0\DDD\nbigm_i\simeq \DDD\nbigm_i$,
and they are strict.
Moreover, they are equipped with
induced filtrations.
\hfill\qed
\end{cor}

Let $\varphi:\nbigt_1\lrarr\nbigt_2$ be
a morphism in $\MTM^{\good}(X,D)$.
Let $\nbigm_i'$ and $\nbigm_i''$
be the $\nbigr_X$-modules
underlying $\nbigt_i$.
We have the underlying morphisms
$\varphi':\nbigm_2'\lrarr\nbigm_1'$
and $\varphi'':\nbigm_1''\lrarr\nbigm_2''$.

\begin{cor}
\label{cor;10.9.17.4}
We have natural isomorphisms
\[
\Ker\DDD\varphi'=\DDD\Cok\varphi',\quad
\Image\DDD\varphi'=\DDD\Image\varphi',\quad
\Cok\DDD\varphi'=\DDD\Ker\varphi'.
\]
We have similar isomorphism for $\varphi''$.
\hfill\qed
\end{cor}

We shall use it in the special case.
Let $g$ be a holomorphic function
such that $g^{-1}(0)=D$.
Let $\nbigv\in\MTS^{\adm}(X,D)$
with the underlying $\nbigr_{X(\ast D)}$-modules.
For any $a,b$,
we have 
$\bigl(\Pi^{a,b}\nbigv\bigr)^{\lor}
\simeq
 \Pi^{-b+1,-a+1}\nbigv^{\lor}$
in $\MTS^{\adm}(X,D)$.

\begin{cor}
We have natural isomorphisms
of $\nbigr_X$-modules
\[
\DDD\Pi^{a,b}_{g!\ast}\nbigv_{i}
\simeq \Pi^{-b+1,-a+1}_{g!\ast}\nbigv^{\lor}_i.
\]
In particular, we have
$\DDD\psi^{(a)}_g\nbigv_{i}
\simeq \psi^{(-a+1)}_g\nbigv^{\lor}_i$
and
$\DDD\Xi^{(a)}_g\nbigv_{i}
\simeq \Xi^{(-a)}_g\nbigv^{\lor}_i$.
\end{cor}
\pf
Let us consider the natural morphism
$\Pi^{a,N}_{g!}\nbigv_i
\lrarr
 \Pi^{b,N}_{g\ast}\nbigv_i$
for a sufficiently large $N$.
According to Lemma \ref{lem;10.9.17.2},
its dual is naturally identified with
\[
 \Pi^{-N+1,-b+1}_{g!}\nbigv_i^{\lor}
\lrarr
 \Pi^{-N+1,-a+1}_{g\ast}\nbigv_i^{\lor}.
\]
Hence, the claim follows from
Corollary \ref{cor;10.9.17.4}.
\hfill\qed

\begin{cor}
$\Pi^{-b+1,-a+1}_{g!\ast}\nbigv^{\lor}$
is a dual of $\Pi^{a,b}_{g!\ast}\nbigv$
as an $\nbigr_X$-triple.
In particular,
$\psi^{(a)}_g(\nbigv)$
and $\Xi^{(a)}_g(\nbigv)$
have duals as $\nbigr_X$-triples,
and naturally
$\DDD\psi^{(a)}_g(\nbigv)\simeq\psi^{(-a+1)}_g(\nbigv^{\lor})$
and 
$\DDD\Xi^{(a)}_g(\nbigv)\simeq\Xi^{(-a)}_g(\nbigv^{\lor})$.
\hfill\qed
\end{cor}

\subsection{Dual of the canonical prolongation in MTM}

Let $X$ and $D$ be as above.
Let $\nbigv\in\MTS^{\adm}(X,D)$.
We have already obtained an isomorphism
$\DDD(\nbigv_!)\simeq
 \nbigv^{\lor}_{\ast}$
as $\nbigr_X$-triples.
We have the filtration $\DDD W$
obtained as the dual of $W$ on $\nbigv_!$,
and the filtration $W$ of 
$\nbigv^{\lor}_{\ast}$.

\begin{lem}
They are the same.
\end{lem}
\pf
Recall that $\nbigv_!$
is obtained as the cohomology 
of the complex of $\nbigr_X$-triples
\[
 \psi_g^{(1)}\nbigv
\lrarr\Xi_g^{(0)}\nbigv\oplus
 \phi_g^{(0)}(\nbigv_!)
 \lrarr\psi_g^{(0)}\nbigv,
\]
where $\phi_g^{(0)}\bigl(\nbigv_!\bigr)
\simeq\psi_g^{(1)}\bigl(\nbigv\bigr)$.
The weight filtration of $\nbigv_!$ is obtained from
the naively induced filtrations $L$ on 
$\psi_g^{(a)}(\nbigv)$ and 
$\Xi_g^{(0)}(\nbigv)$,
and the filtration $L$ of 
$\phi_g^{(0)}\bigl(\nbigv_!\bigr)$
obtained as the transfer of $L$ on
$\psi_g^{(a)}(\nbigv)$.

Hence, $\DDD(\nbigv_!)$ is obtained
as the cohomology of
the complex of $\nbigr$-triples
\[
 \DDD\psi_g^{(0)}(\nbigv)
\lrarr 
 \DDD\Xi_g^{(0)}(\nbigv)\oplus
 \DDD\phi_g^{(0)}(\nbigv_!)\lrarr
 \DDD\psi_g^{(1)}(\nbigv), 
\]
and the filtration 
$\DDD W$ is induced by
$\DDD L$ on 
$\DDD\psi_g^{(0)}(\nbigv)$,
$\DDD\Xi_g^{(0)}(\nbigv)$,
and 
$\DDD\phi_g^{(1)}(\nbigv_!)$.
It is easy to check that
$\DDD L$ on 
$\DDD\psi_g^{(a)}(\nbigv)$ and
$\DDD\Xi_g^{(0)}(\nbigv)$
are the same as the naively induced
filtrations under the isomorphisms
$\DDD\psi_g^{(a)}(\nbigv)
\simeq
\psi_g^{(-a+1)}(\nbigv^{\lor})$
and
$\DDD\Xi_g^{(0)}(\nbigv)
\simeq
 \Xi_g^{(0)}(\nbigv^{\lor})$.
Because $\DDD W$ 
on $\DDD\psi_g^{(a)}(\nbigv)$
is $M(\DDD\nbign;\DDD L)[1-2a]$,
we obtain $\DDD W$ is the canonical
weight filtration of 
$\psi_g^{(-a+1)}(\nbigv^{\lor})$.
Hence, $\DDD W$ on 
$\DDD\phi^{(0)}_g(\nbigv_!)$
is the same as
the canonical weight filtration
$W$ of $\DDD\phi^{(0)}_g(\nbigv_!)\simeq
 \phi^{(0)}_g(\nbigv^{\lor}_{\ast})$.
We obtain that
$\DDD L$ on $\DDD\phi^{(0)}_g(\nbigv_!)$
is the transfer of $L$ on 
$\psi_g^{(a)}(\nbigv^{\lor})$.
and hence,
$\DDD W$ on $\DDD\nbigv_!$
is the same as the canonical weight filtration
of $\nbigv^{\lor}_{\ast}$.
\hfill\qed

\begin{cor}
\label{cor;10.9.17.11}
The induced filtrations of
$\DDD\psi_g^{(a)}\nbigv$
and 
 $\DDD\Xi_g^{(a)}\nbigv$
are the canonical weight filtrations
under the identifications
$\DDD\psi_g^{(a)}\nbigv
 \simeq \psi_g^{(-a+1)}\nbigv^{\lor}$
and $\DDD\Xi_g^{(a)}\nbigv
\simeq \Xi_g^{(-a)}\nbigv^{\lor}$.
\hfill\qed
\end{cor}

\begin{rem}
In the case of $\DDD\psi_g^{(a)}\nbigv$,
we can also deduce it by using the characterization
as the relative monodromy filtration.
\hfill\qed
\end{rem}

\subsection{Proof of Theorem \ref{thm;10.9.13.11}}
\label{subsection;11.1.29.3}

To show Theorem \ref{thm;10.9.13.11},
we use an induction,
i.e., we shall show $A(n)$ by assuming $A(n-1)$.
Let $\nbigt\in\MTM(X)$
such that $\dim\Supp\nbigt=n$.
Let $P\in\Supp\nbigt$.
Because the claims are local,
we have only to show them
for the restriction of $\nbigt$
to a small neighbourhood of $P$.
We will shrink $X$ without mention
in the following argument.

We take an $n$-dimensional admissible cell
$\nbigc=(Z,U,\varphi,\nbigv)$ of $\nbigt$ at $P$.
Let $g$ be a cell function.
We have the expression of $\nbigt$ around $P$
as the cohomology of the complex in $\MTM(X)$:
\[
 \psi_g^{(1)}\varphi_{\dagger}(\nbigv)
\lrarr
 \Xi_g^{(0)}\varphi_{\dagger}(\nbigv)
\oplus
 \phi_g^{(0)}\nbigt
\lrarr
 \psi_g^{(0)}\varphi_{\dagger}(\nbigv)
\]
By the assumption $A(n-1)$,
Corollary \ref{cor;10.9.17.11} 
and the result in \S\ref{subsection;10.9.17.12},
we obtain the following complex
in $\MTM(X)$:
\begin{equation}
 \DDD\psi_g^{(0)}\varphi_{\dagger}(\nbigv)
\lrarr
 \DDD\Xi_g^{(0)}\varphi_{\dagger}(\nbigv)
\oplus
 \DDD\phi_g^{(0)}\nbigt
\lrarr
 \DDD\psi_g^{(1)}\varphi_{\dagger}(\nbigv).
\end{equation}
We obtain a mixed twistor $D$-module
$\DDD\nbigt$ as the cohomology.
The underlying filtered $\nbigr_X$-modules are 
the dual of the underlying filtered $\nbigr_X$-modules of
$\nbigt$.
Let $C_1$ and $C_2$
denote the pairings of $\nbigt$ and $\DDD\nbigt$,
respectively.
By construction, 
$C_{2|\nbigxlambda}
=\DDD C_{1|\nbigxlambda}$
for generic $\lambda$.
Hence, we obtain $C_2=\DDD C_1$.

If $\nbigt\in\MTMint(X)$,
then $\nbigv$ and $\nbigv^{\lor}$ are integrable.
Hence, by construction,
$\DDD\psi_g^{(a)}\varphi_{\dagger}(\nbigv)$
and 
$\DDD\Xi_g^{(0)}\varphi_{\dagger}(\nbigv)$
are also naturally integrable.
By the hypothesis of the induction,
$\DDD\phi_g^{(0)}\nbigt$ is integrable.
Hence, $\DDD\nbigt$ is integrable.
Thus, the induction can go further.
\hfill\qed

\section{Real structure of mixed twistor $D$-modules}

\subsection{Some functors}
\label{subsection;11.3.30.12}

Let $\DDD^{\herm}$ denote the hermitian dual,
i.e.,
$\DDD^{\herm}(\nbigt)=\nbigt^{\ast}$.
\index{functor $\DDD^{\herm}$}
Formally,
we set $\sigmatilde^{\ast}:=
 \DDD\circ\DDD^{\herm}
=\DDD^{\herm}\circ\DDD$.
\index{functor $\sigmatilde^{\ast}$}
It gives an involution on
$\MTM(X)$.
We have a natural transformation
$\sigmatilde^{\ast}\circ\sigmatilde^{\ast}
 \simeq \id$.
We define
$\gammatilde^{\ast}:=
 j^{\ast}\circ\sigmatilde^{\ast}
=\sigmatilde^{\ast}\circ j^{\ast}$.
\index{functor $\gammatilde^{\ast}$}
We naturally have
$\gammatilde^{\ast}\circ\gammatilde^{\ast}\simeq \id$.
For $\nbigt=(\nbigm',\nbigm'',C)\in\MTM(X)$,
we have
\[
 \sigmatilde^{\ast}\nbigt
=\DDD(\nbigt^{\ast})=
 \bigl(
 \DDD\nbigm'',\DDD\nbigm',\DDD C^{\ast}
 \bigr),
\quad\quad
 \gammatilde^{\ast}\nbigt=
 \Bigl(
 j^{\ast}\DDD\nbigm'', j^{\ast}\DDD\nbigm', j^{\ast}\DDD C^{\ast}
 \Bigr)
\]
with the induced weight filtrations.
The following lemma is clear by construction.
\begin{lem}\mbox{{}}
\label{lem;11.3.29.2}
\begin{itemize}
\item
We have natural isomorphism
$\sigmatilde^{\ast}\circ\gammatilde^{\ast}\simeq
 \gammatilde^{\ast}\circ\sigmatilde^{\ast}\simeq j^{\ast}$.
We also have
$\sigmatilde^{\ast}\circ j^{\ast}=
 j^{\ast}\circ\sigmatilde^{\ast}=
 \gammatilde^{\ast}$
and
$j^{\ast}\circ\gammatilde^{\ast}=
 \gammatilde^{\ast}\circ j^{\ast}=\sigmatilde^{\ast}$.
\item
Let $\nutilde^{\ast}$ be 
one of $\sigmatilde^{\ast}$, 
$\gammatilde^{\ast}$ or $j^{\ast}$.
We have natural commutativity
$\nutilde^{\ast}\circ\DDD
\simeq
 \DDD\circ\nutilde^{\ast}$
and
$\nutilde^{\ast}\circ\DDD^{\herm}
\simeq
 \DDD^{\herm}\circ\nutilde^{\ast}$
by the natural isomorphisms of
the underlying filtered $\nbigr$-modules.
\item
Let $F$ be a projective morphism.
The natural transform
$\DDD F_{\dagger}\simeq F_{\dagger}\DDD$
for $\nbigr$-modules induces
$F^p_{\dagger}\gammatilde^{\ast}(\nbigt)
\simeq
 \mu_p\circ \gammatilde^{\ast}F^{p}_{\dagger}(\nbigt)$
and
$F^p_{\dagger}\sigmatilde^{\ast}(\nbigt)
\simeq
 \mu_p\circ \sigmatilde^{\ast} F^{p}_{\dagger}(\nbigt)$.
We also have 
$j^{\ast}\circ F_{\dagger}=F_{\dagger}\circ j^{\ast}$
naturally.
\hfill\qed
\end{itemize}
\end{lem}
We consider $\nbigt=(\nbigm_1,\nbigm_2,C)\in\MTM(X)$.
Let $\varphi_1:F_{\dagger}^{-p}\DDD\nbigm_1
 \simeq \DDD F_{\dagger}^p\nbigm_1$
and
$\varphi_2:F_{\dagger}^p\DDD\nbigm_2\simeq
 \DDD F^{-p}_{\dagger}\nbigm_2$
be the natural isomorphisms.
Then,
$\bigl(\epsilon(p)\varphi_2,\epsilon(-p)\varphi_1^{-1}\bigr)$
gives an isomorphism
$F_{\dagger}^p\sigmatilde^{\ast}(\nbigt)
\simeq
 \sigmatilde^{\ast}F_{\dagger}^{p}(\nbigt)$,
and 
$\bigl(\epsilon(p)j^{\ast}\varphi_2,\epsilon(-p)j^{\ast}\varphi_1^{-1}\bigr)$
gives an isomorphism
$F_{\dagger}^p\gammatilde^{\ast}(\nbigt)
\simeq
 \gammatilde^{\ast}F_{\dagger}^{p}(\nbigt)$.

\vspace{.1in}

The following lemma is obvious by construction.
\begin{lem}
Let $\nbigt\in\MTM(X)$
and $\nbigv\in\MTS^{\adm}(X,H)$.
Then, we have the relation
$\gammatilde^{\ast}(\nbigt\otimes\nbigv[\star H])
\simeq
 \bigl(
 \gammatilde^{\ast}(\nbigt)
\otimes
 \gammatilde_{\sm}^{\ast}(\nbigv)
 \bigr)[\star H]$.
\end{lem}
\pf
We have naturally
$j^{\ast}\bigl(
 (\nbigt\otimes\nbigv)[\ast H]
 \bigr)
\simeq
 \bigl(
 j^{\ast}\nbigt\otimes
 j^{\ast}\nbigv\bigr)[\ast H]$,
$\Bigl(
 (\nbigt\otimes\nbigv)[\ast H]
\Bigr)^{\ast}
\simeq
 (\nbigt^{\ast}\otimes \nbigv^{\ast})[!H]$
and
$\DDD\Bigl(
 (\nbigt\otimes\nbigv)[!H]
\Bigr)
\simeq
 \bigl(\DDD(\nbigt)\otimes \nbigv^{\lor}\bigr)[\ast H]$,
which induce the desired isomorphism.
\hfill\qed

\subsubsection{Complex of mixed twistor $D$-modules}

Let $\nbigc(\MTM(X))$ denote the category of
bounded complexes of mixed twistor $D$-modules.
Let $(\nbigt^{\bullet},\delta^{\bullet})\in\nbigc(\MTM(X))$.
\begin{itemize}
\item
Let $\DDD^{\herm}(\nbigt^{\bullet})
\in\nbigc(\MTM(X))$
be given by
$\DDD^{\herm}(\nbigt^{\bullet})^p
=\DDD^{\herm}(\nbigt^{-p})$
with the differentials
$\DDD^{\herm}(\delta)^{p}=\DDD^{\herm}(\delta^{-p-1})$.
\index{functor $\DDD^{\herm}$}
\item
Let $\DDD(\nbigt^{\bullet})\in\nbigc(\MTM(X))$
be given by
$\DDD(\nbigt^{\bullet})^p=
 \mu_p\DDD(\nbigt^{-p})$
with the differentials
\[
 \DDD(\delta)^p=\bigl(
 (-1)^{p}\DDD\delta_1^{p},
 (-1)^{p+1}\DDD\delta_2^{-p-1}
 \bigr). 
\]
\index{functor $\DDD$}
\item
Let $j^{\ast}(\nbigt^{\bullet})\in\nbigc(\MTM(X))$
be given by
$j^{\ast}(\nbigt^{\bullet})^p
=j^{\ast}(\nbigt^p)$
with the differential
$j^{\ast}(\delta)^p=j^{\ast}(\delta^p)$.
\index{functor $j^{\ast}$}
\item
We define $\nbigs_{\ell}(\nbigt^{\bullet})^p
\in\nbigc(\MTM(X))$
by
$\nbigs_{\ell}(\nbigt^{\bullet})^p
=\nbigt^{p+\ell}$
with 
$\nbigs_{\ell}(\delta)^p=(-1)^{\ell}\delta^{p+\ell}$.
\index{functor $\nbigs_{\ell}$}
\end{itemize}
The natural quasi-isomorphisms of the underlying complexes of
$\nbigr$-modules give the commutativity of the functors
$\DDD^{\herm}$, $\DDD$ and $j^{\ast}$.
We also have
$\nbigs_{\ell}\circ\DDD^{\herm}
\simeq
 \DDD^{\herm}\circ \nbigs_{-\ell}$,
$\nbigs_{\ell}\circ\DDD
\simeq
 \DDD\circ\nbigs_{-\ell}$
and
$j^{\ast}\circ\nbigs_{\ell}
=\nbigs_{\ell}\circ j^{\ast}$.
(See \S\ref{subsection;11.3.29.1} below
for the signature.)

We define the functors
$\sigmatilde^{\ast}$
and $\gammatilde^{\ast}$
on $\nbigc(\MTM(X))$
by 
$\sigmatilde^{\ast}:=
 \DDD\circ\DDD^{\herm}$
and $\gammatilde^{\ast}:=j^{\ast}\circ\sigmatilde^{\ast}$.
The natural quasi-isomorphisms
of the underlying $\nbigr$-modules
give 
$\nbigs_{\ell}\circ\nutilde^{\ast}
\simeq
 \nutilde^{\ast}\circ\nbigs_{\ell}$,
where $\nutilde^{\ast}$ is one of
$j^{\ast}$, $\gammatilde^{\ast}$
or $\sigmatilde^{\ast}$.
The compatibility with the push-forward
is reformulated as follows.
\begin{lem}
Let $F$ be a projective morphism.
We have natural isomorphisms
$\nbigs_{-p}\circ F^p_{\dagger}(\DDD\nbigt)
\simeq
 \DDD \nbigs_p\circ F^{-p}_{\dagger}(\nbigt)$,
$\nbigs_{-p}\circ F^p_{\dagger}(\sigmatilde^{\ast}\nbigt)
\simeq
 \sigmatilde^{\ast}\nbigs_{-p}\circ F_{\dagger}^{p}(\nbigt)$
and 
$\nbigs_{-p}\circ F^p_{\dagger}(\gammatilde^{\ast}\nbigt)
\simeq
 \gammatilde^{\ast}\nbigs_{-p}\circ F_{\dagger}^{p}(\nbigt)$.
\hfill\qed
\end{lem}

\subsection{Real structure of mixed twistor $D$-modules}
\index{real structure}
\index{category $\MTM(X,\real)$}
\index{category $\MTM^{\integral}(X,\real)$}

Let $\nbigt$ be 
an (integrable) mixed twistor $D$-module.
A real structure of $\nbigt$ is an (integrable) isomorphism
$\iota:\gammatilde^{\ast}\nbigt\simeq\nbigt$
such that
$\gammatilde^{\ast}\iota\circ\iota=\id$.
A morphism of
(integrable) mixed twistor $D$-modules with real structure
is an (integrable) morphism compatible with the real structure.
Let $\MTM(X,\real)$ (resp. $\MTMint(X,\real)$)
denote the category of mixed twistor $D$-modules 
(resp. integrable mixed twistor $D$-modules)
with real structure
(resp. integrable real structure).
We obtain the following proposition.
\begin{prop}
\label{prop;11.3.24.11}
Let $\nbigt$ be an (integrable) mixed twistor $D$-module
with (integrable) real structure on $X$.
\begin{itemize}
\item
$\DDD(\nbigt)$,
$\DDD^{\herm}(\nbigt)$
and $j^{\ast}(\nbigt)$
have natural real structures.
\item
Let $F:X\lrarr Y$ be a projective morphism.
Then, $F^i_{\dagger}(\nbigt)$
are equipped with induced real structures.
\item
Let $H$ be a hypersurface of $X$.
Then,
$\nbigt[\star H]$ is naturally equipped with real structure.
\item
Let $\nbigv$ be an (integrable) admissible mixed twistor
structure on $(X,D)$.
If $\nbigv$ has a real structure,
$(\nbigt\otimes\nbigv)[\star H]$ is also
equipped with an induced real structure
for $\star=\ast,!$.
\hfill\qed
\end{itemize}
\end{prop}

\begin{cor}
Let $g$ be any holomorphic function.
Then,
$\psi^{(a)}_g(\nbigt)$,
$\Xi^{(a)}_g\nbigt$
and $\phi_g^{(a)}\nbigt$
are naturally equipped with induced real structures.
\hfill\qed
\end{cor}

Let $X$ be a complex manifold
with a hypersurface $H$.
Let $\MTM(X,[\ast H],\real)\subset\MTM(X,\real)$ be
the full subcategory of objects $(\nbigt,\iota)\in\MTM(X,\real)$
such that $\nbigt[\ast H]=\nbigt$.
We use the notation
$\MTMint(X,[\ast H],\real)$
in a similar meaning.

\begin{prop}
Let $F:(X',H')\lrarr (X,H)$ be a projective
birational morphism such that
$X'\setminus H'\simeq X\setminus H$.
Then, $F_{\dagger}$ induces
equivalences
$\MTM(X',[\ast H'],\real)\simeq
 \MTM(X,[\ast H],\real)$
and
$\MTMint(X',[\ast H'],\real)\simeq
 \MTMint(X,[\ast H],\real)$.
\hfill\qed
\end{prop}

\subsubsection{Example}
Let us look at the smooth $\nbigr$-triple
$\nbigu_X(d_X,0)$.
It is a mixed twistor $D$-module
according to Theorem \ref{thm;10.11.13.11}.
It is naturally equipped with an integrable structure.
We have the following isomorphism:
\begin{multline}
 \gammatilde^{\ast}\nbigu(d_X,0)
=\DDD j^{\ast}
 \bigl(
 \nbigo_{\nbigx},\,\nbigo_{\nbigx}\,\lambda^{d_X},\,
 C_0
 \bigr)
=\DDD\bigl(
 \nbigo_{\nbigx},\,\nbigo_{\nbigx}\,\lambda^{d_X},
 C_0\bigr) \\
=
 \bigl(
 \nbigo_{\nbigx}\lambda^{d_X},\,
 \nbigo_{\nbigx},\,(-1)^{d_X}C_0
 \bigr)
\end{multline}
Here, we fix the isomorphism
$\DDD\nbigo_{\nbigx}
\simeq
 \lambda^{d_X}\,\nbigo_{\nbigx}$
determined by the condition that
its restriction to $\lambda=1$ 
induces the isomorphism $\nu$
in \S\ref{subsection;11.3.30.2}.
Hence, for example,
$\bigl(\epsilon(-d_X),\epsilon(d_X)\bigr)$
and
$\bigl(\epsilon(d_X),\epsilon(-d_X)\bigr)$
give real structures.
(We remark the signature in Lemma \ref{lem;11.1.22.1}.)

\subsection{$\real$-Betti structure
of the underlying $D$-modules}

Let us consider $(\nbigt,W)\in\MTM(X,\real)$.
Let $\nbigm_i$ $(i=1,2)$ be the underlying
$\nbigr_X$-modules.
Take any $\lambda_0\in\vecS$.
The real structure of $\nbigt$ induces
a real structure of the non-degenerate $D$-triple
$(\nbigm_1^{-\lambda_0},\nbigm_2^{\lambda_0},C^{\lambda_0})$.
(See \S\ref{subsection;11.3.29.1}.)
As explained in \S\ref{subsection;13.8.29.20},
\ref{subsection;13.4.13.10}
and \ref{subsection;13.8.29.40},
$\nbigm_2^{\lambda_0}$
is equipped with a naturally induced 
$\real$-Betti structure
in the sense of \S7.2 of \cite{mochi9}.
Thus, 
we obtain a functor
$\Upsilon^{\lambda_0}:
 \MTM(X,\real)
\lrarr
 \Hol(X,\real)$,
where $\Hol(X,\real)$ denotes
the category of holonomic $D$-modules
on $X$ with $\real$-Betti structure.
We obtain the following from
Proposition \ref{prop;13.8.30.12}.

\begin{prop}
\label{prop;13.8.30.30}
The functor $\Upsilon^{\lambda_0}$
is compatible with
the following functors:
\begin{itemize}
\item
 The push-forward by any projective morphisms.
\item
 The dual $\DDD$.
\hfill\qed
\end{itemize}
\end{prop}

We obtain the following from 
Proposition \ref{prop;13.8.30.13}.
\begin{prop}
For any hypersurface $H$,
the localizations
$[\ast H]$ and $[!H]$
are compatible with 
$\Upsilon^{\lambda_0}$.
\hfill\qed
\end{prop}

Let $\nbigv=(\nbigv_1,\nbigv_2,C) 
 \in\MTS^{\adm}(X,H)$
with a real structure
as a smooth $\nbigr_{X(\ast H)}$-triple.
Then, the meromorphic flat bundle
$\nbigv_2^{\lambda_0}$
has a good $\real$-structure.
We obtain the following from Proposition
\ref{prop;13.8.30.14}.
\begin{prop}
Let $\nbigt\in\MTM(X,\real)$.
Let $\nbigv\in\MTS^{\adm}(X,H)$
with a real structure.
Then, we naturally have
$\Upsilon^{\lambda_0}\bigl(
 (\nbigt\otimes\nbigv)[\star H]
 \bigr)
\simeq
 \bigl(
 \Upsilon^{\lambda_0}(\nbigt)
 \otimes
 \nbigv_2^{\lambda_0}
 \bigr)[\star H]$
in $\Hol(X,\real)$.
\hfill\qed
\end{prop}

\begin{cor}
We naturally have
$\Upsilon^{\lambda_0}\circ\Xi^{(a)}
\simeq \Xi^{(a)}\circ\Upsilon^{\lambda_0}$,
$\Upsilon^{\lambda_0}\circ\psi^{(a)}
\simeq
 \psi^{(a)}\circ\Upsilon^{\lambda_0}$
and 
$\Upsilon^{\lambda_0}\circ\phi^{(a)}
\simeq
 \phi^{(a)}\circ\Upsilon^{\lambda_0}$.
\hfill\qed
\end{cor}

\subsection{Real structure of the underlying
integrable $\nbigr$-modules}

\subsubsection{Preliminary}

Let $(\nbigt,W)$ be an integrable mixed twistor $D$-module on $X$.
Let $\nbigm_i$ $(i=1,2)$ be the underlying
integrable $\nbigr_X$-modules.
Let $\pi_0:\nbigh\lrarr\cnum_{\lambda}^{\ast}$ be 
a universal covering.
Let $p_{\nbigh}:\nbigh\times X\lrarr X$
be the projection.
Fix any $\lambda_0\neq 0$,
and $\lambdatilde_0\in\nbigh$
such that $\pi_0(\lambdatilde_0)=\lambda_0$.
We set 
$\nbigm^{\lambda_0}_i:=
 \nbigo_{\{\lambda_0\}\times X}
 \otimes_{
 \nbigo_{\cnum_{\lambda}\times X}}
 \nbigm_i$.
Let 
$i_{\lambda_0}:\{\lambdatilde_0\}\times X
\lrarr \nbigh\times X$ 
denote the inclusion.
We have a natural morphism
\begin{equation}
 \label{eq;13.8.29.1}
 \Omega^{\bullet}_{\nbigh\times X}
\otimes_{\nbigo_{\nbigh\times X}}
 \pi_0^{-1}(\nbigm_{i|\cnum_{\lambda}^{\ast}\times X})
\lrarr
 i_{\lambda_0\ast}
 \bigl(
 \Omega_X^{\bullet}\otimes_{\nbigo_X}
 \nbigm^{\lambda_0}_i
\bigr)
\end{equation}

\begin{prop}
\label{prop;13.8.29.11}
We have a unique isomorphism
in $D^b_c(\nbigh\times X)$
\[
 \Phi_i^{\lambda_0}:
 \Omega^{\bullet}_{\nbigh\times X}
 \otimes_{\nbigo_{\nbigh\times X}}
 \pi_0^{-1}\nbigm_{i|\cnum_{\lambda}^{\ast}\times X}
\simeq
 p_{\nbigh}^{-1}\bigl(
 \Omega_X^{\bullet}\otimes_{\nbigo_X}
 \nbigm^{\lambda_0}_{i}
 \bigr)
\]
such that 
the composite of
$\Phi_i^{\lambda_0}$
and 
the natural morphism
$p_{\nbigh}^{-1}\bigl(
 \Omega_X^{\bullet}\otimes_{\nbigo_X}
 \nbigm^{\lambda_0}_{i}
 \bigr)
\lrarr
  i_{\lambda_0\ast}
 \bigl(
 \Omega_X^{\bullet}\otimes_{\nbigo_X}
 \nbigm^{\lambda_0}_i
\bigr)
$
is equal to {\rm(\ref{eq;13.8.29.1})}.
\end{prop}
\pf
We begin with a special case.
Let $H$ be a simply normal crossing hypersurface of $X$.
Let $(\nbigv_1,\nbigv_2,C)$ be a good admissible
mixed twistor structure on $(X,H)$.
Note that
$\nbigv_{i|\cnum^{\ast}\times X}$ 
are good meromorphic flat bundle 
on $\cnum^{\ast}\times (X,H)$.
Recall Lemma \ref{lem;13.8.29.2}.
By Theorem 11.3.3 of 
\cite{kashiwara-schapira},
the micro-supports of
$\Omega^{\bullet}_{\nbigh\times X}\otimes
 \pi_0^{-1}\nbigv_i[\star H]$
is contained in 
$\nbigh\times \nbigs$
for a Lagrangian cone in 
the cotangent bundle of $X$.
By Proposition 5.4.5 of \cite{kashiwara-schapira},
there exists an object
$G_{i\star}\in D^b_c(X)$
such that 
$\Omega^{\bullet}_{\nbigh\times X}\otimes
 \pi_0^{-1}\nbigv_i[\star H]$
is isomorphic to
$p_{\nbigh}^{-1}G_{i\star}$
in $D^b_c(\nbigh\times X)$.
By using the description of
the de Rham complex of $\nbigv_{i}[\star H]$
in \S5.1.1 of \cite{mochi9},
we can observe that
\[
 G_{i\star}\simeq
 i_{\lambda_0}^{-1}
 \Omega^{\bullet}_{\nbigh\times X}\otimes
 \pi_0^{-1}\nbigv_i[\star H]
\simeq
 \Omega^{\bullet}_X\otimes
 \nbigv_i^{\lambda_0}[\star H].
\]
Let $g$ be any holomorphic function such that
$g^{-1}(0)=H$.
We have unique isomorphisms
in $D^b_c(\nbigh\times X)$
\[
\Omega_{\nbigh\times X}^{\bullet}\otimes
\pi_0^{-1}\bigl(
 \Pi^{a,b}_{g\star}(\nbigv_i)
\bigr)
\simeq
 p_{\nbigh}^{-1}\bigl(
 \Omega_X^{\bullet}\otimes
 \Pi^{a,b}_{g\star}(\nbigv_i)^{\lambda_0}
 \bigr)
\]
with the desired property.
By Proposition 2.7.8 of \cite{kashiwara-schapira},
the natural morphism
$\Omega_{\nbigh\times X}^{\bullet}\otimes
\pi_0^{-1}\bigl(
 \Pi^{a,b}_{g!}(\nbigv_i)
\bigr)
\lrarr
 \Omega_{\nbigh\times X}^{\bullet}\otimes
\pi_0^{-1}\bigl(
 \Pi^{a,b}_{g\ast}(\nbigv_i)
\bigr)$
is obtained as the pull back of the natural morphism
$\Omega_{X}^{\bullet}\otimes
 \Pi^{a,b}_{g!}(\nbigv_i)^{\lambda_0}
\lrarr
 \Omega_{X}^{\bullet}\otimes
 \Pi^{a,b}_{g\ast}(\nbigv_i)^{\lambda_0}$.
Hence, we have a unique isomorphism
\[
\Omega_{\nbigh\times X}^{\bullet}\otimes
\pi_0^{-1}\bigl(
 \Pi^{a,b}_{g!\ast}(\nbigv_i)
\bigr)
\simeq
 p_{\nbigh}^{-1}\bigl(
 \Omega_X^{\bullet}\otimes
 \Pi^{a,b}_{g!\ast}(\nbigv_i)^{\lambda_0}
 \bigr)
\]
with the desired property.
In particular,
we have a desired isomorphism
for $\Xi^{(a)}_g(\nbigv_i)$.
Then, we obtain the desired isomorphisms
for the integrable $\nbigr_X$-modules
underlying integrable mixed twistor $D$-modules,
by using a standard Noetherian induction
on the supports.
\hfill\qed

\vspace{.1in}
Let $\nbigh_{\real}:=\pi_0^{-1}(\vecS)$.
Let $\pi_1:\nbigh_{\real}\lrarr\vecS$ be the restriction
of $\pi_0$.
For any $C^{\infty}$-manifold $Z$,
let $\nbige_{Z}^{\bullet}$
denote the $C^{\infty}$-de Rham complex of $Z$.
Let $p_{\nbigh_{\real}}:\nbigh_{\real}\times X\lrarr X$
be the projection.
Fix any $\lambda_0\in\vecS$.
Let $\iota:\nbigh_{\real}\times X\lrarr \nbigh\times X$
denote the inclusion.
The inclusion $i_{1,\lambda_0}:\{\lambda_0\}\times X\lrarr 
 \nbigh_{\real}\times X$
induces 
\begin{equation}
 \label{eq;13.8.29.10}
 \nbige^{\bullet}_{\nbigh_{\real}\times X}
 \otimes_{\iota^{-1}\nbigo_{\nbigh\times X}}
 \iota^{-1}
 \pi_0^{-1}\bigl(
 \nbigm_{i|\cnum^{\ast}\times X}
 \bigr)
\lrarr 
 i_{1,\lambda_0\ast}
\bigl(
 \nbige^{\bullet}_{X}\otimes_{\nbigo_X}
 \nbigm^{\lambda_0}_i
 \bigr).
\end{equation}

\begin{prop}
\label{prop;13.8.29.12}
We have unique isomorphisms
in $D^b_c(\nbigh_{\real}\times X)$
\[
 \Phi^{\lambda_0}_{1,i}:
  \nbige^{\bullet}_{\nbigh_{\real}\times X}
 \otimes_{\iota^{-1}\nbigo_{\nbigh\times X}}
 \iota^{-1}
 \pi_0^{-1}\bigl(
 \nbigm_{i|\cnum^{\ast}\times X}
 \bigr)
\simeq
 p_{\nbigh_{\real}}^{-1}
\bigl(
 \nbige^{\bullet}_{X}\otimes_{\nbigo_X}
 \nbigm^{\lambda_0}_i
 \bigr)
\]
such that the composite of 
$\Phi^{\lambda_0}_{1,i}$
and the natural
$p_{\nbigh_{\real}}^{-1}
\bigl(
 \nbige^{\bullet}_{X}\otimes_{\nbigo_X}
 \nbigm^{\lambda_0}_i
 \bigr)
\lrarr
 i_{1,\lambda_0\ast}
\bigl(
 \nbige^{\bullet}_{X}\otimes_{\nbigo_X}
 \nbigm^{\lambda_0}_i
 \bigr)$
is equal to 
{\rm(\ref{eq;13.8.29.10})}.
\end{prop}
\pf
Let $q:\nbigh\lrarr \nbigh_{\real}$ be the morphism
corresponding to the projection
$\cnum^{\ast}\simeq \vecS\times\real_{>0}
\lrarr \vecS$.
By using the argument in the proof of 
Proposition \ref{prop;13.8.29.11},
we obtain isomorphisms
in $D^b_c(\nbigh\times X)$
\[
  \nbige^{\bullet}_{\nbigh\times X}
\otimes_{\nbigo_{\nbigh\times X}}
 \pi_0^{-1}(\nbigm_{i|\cnum_{\lambda}^{\ast}\times X})
\simeq
 q^{-1}\Bigl(
 \nbige^{\bullet}_{\nbigh_{\real}\times X}
 \otimes_{\iota^{-1}\nbigo_{\nbigh\times X}}
 \iota^{-1}
 \pi_0^{-1}\bigl(
 \nbigm_{i|\cnum_{\lambda}^{\ast}\times X}
 \bigr)
\Bigr)
\]
such that its restriction to
$\nbigh_{\real}\times X$ 
is the natural one.
Then, the claim of Proposition \ref{prop;13.8.29.12}
is obtained as the specialization of
Proposition \ref{prop;13.8.29.11}.
\hfill\qed

\subsubsection{Real structure of the underlying
integrable $\nbigr_X$-modules}

Let $(\nbigt,W)\in\MTMint(X,\real)$.
Let $\nbigm_i$ $(i=1,2)$ be the underlying
integrable $\nbigr_X$-modules.
Take any $\lambda_0\in\vecS$.
We have the induced real structure of 
$\nbigm_2^{\lambda_0}$.
The isomorphism in Proposition \ref{prop;13.8.29.11}
induces a real structure of 
$\DR_{\nbigh\times X}\bigl(
 \pi_0^{-1}\nbigm_{2|\cnum_{\lambda}^{\ast}\times X}
\bigr)$.

\begin{lem}
The real structure is independent of 
the choices of $\lambda_0$
and $\lambdatilde_0$.
\end{lem}
\pf
By Proposition \ref{prop;13.8.29.12},
the pairings
$\DR(\nbigm_1^{-\lambda})
\otimes
 \overline{\DR(\nbigm_2^{\lambda})}
\lrarr
 \omega^{top}$
are independent of $\lambda\in\vecS$.
The isomorphisms
$\DR(\nbigm_2^{\lambda})
\simeq
 \DDD\DR(\nbigm_1^{-\lambda})$
are also independent of $\lambda$.
Then, the claim follows.
\hfill\qed

\vspace{.1in}
Hence, we obtain a pre-$\real$-Betti structure
of $\nbigm_{i|\cnum^{\ast}\times X}$.

\begin{lem}
The pre-$\real$-Betti structure of
$\nbigm_{i|\cnum^{\ast}\times X}$
is a $\real$-Betti structure.
\end{lem}
\pf
We have only to check the claim
locally around any point $P$ of $X$.
We shall shrink $X$ without mention.
We take an admissible cell
$(Z,U,\varphi,\nbigv)$ for $\nbigt$
around $P$ with a cell function $g$.
We have
$\varphi_{\dagger}(\nbigv)[\star g]
\simeq\nbigt[\star g]$.
By the standard argument with a Noetherian 
induction on the support,
we have only to prove the claim for
$\varphi_{\dagger}(\nbigv)[\star g]$.
Let $\nbigv=(\nbigv_1,\nbigv_2,C_{\nbigv})$.
We set $D_Z:=Z\setminus U$
and $g_Z:=\varphi^{\ast}(g)$.
We can observe that
$\nbigv_2$ is a good meromorphic flat bundle
on $\cnum_{\lambda}^{\ast}\times (Z,D_Z)$
with a good $\real$-structure.
Hence,
$\nbigv_2[\star g_Z]$ has an $\real$-Betti structure.
The real structure of
$\pi_0^{-1}\bigl(
 \Omega^{\bullet}_{\cnum^{\ast}\times Z}
 \otimes\nbigv_2[\star g_Z]\bigr)$
is equal to the one induced by
the isomorphism with
$p_{\nbigh}^{-1}
\bigl(
 \Omega^{\bullet}_{Z}
\otimes\nbigv_2^{\lambda_0}[\star g_Z]\bigr)$.
Then, the claim for 
$\varphi_{\dagger}(\nbigv)[\star g]$
easily follows.
\hfill\qed

\vspace{.1in}
Thus, 
we obtain a functor
$\Upsilontilde:
 \MTM(X,\real)
\lrarr
 \Hol(\cnum^{\ast}\times X,\real)$.
\begin{prop}
The functor $\Upsilontilde$
is compatible with the following functors.
\begin{itemize}
\item
The push-forward by any projective morphisms.
\item
The dual $\DDD$.
\end{itemize}
\end{prop}
\pf
Let $F:X\lrarr Y$ be any projective morphism.
Let $\nbigt\in\MTM(X,\real)$.
Let $\nbigm_i$ be the underlying $\nbigr_X$-modules.
Let $S(\lambda_0):=
 \bigl\{\lambda\in\cnum^{\ast}\,\big|\,
 |\arg(\lambda)-\arg(\lambda_0)|<\epsilon
 \bigr\}$
for a small $\epsilon$.
We have only to compare the $\real$-structures
on 
$\DR_{\cnum^{\ast}\times Y} F^i_{\dagger}(\nbigm_2)
 \bigr)_{|S(\lambda_0)\times Y}$.
The specializations at $\lambda_0\in\vecS$
are the same
by Proposition \ref{prop;13.8.30.30}.
Let $p_{S(\lambda_0)}$
denote the projections
$S(\lambda_0)\times X\lrarr X$
and $S(\lambda_0)\times Y\lrarr Y$.
Because the $\real$-structures on
$S(\lambda_0)\times Y$ are obtained by
the isomorphisms
\[
 \DR_{\cnum^{\ast}\times X}(\nbigm_2)
 _{|S(\lambda_0)\times X}
 \simeq
 p_{S(\lambda_0)}^{-1}
 \DR_X(\nbigm_2^{\lambda_0}),
\]
\[
 \DR_{\cnum^{\ast}\times Y}F_{\dagger}^i(\nbigm_2)
 _{|S(\lambda_0)\times Y}
 \simeq
 p_{S(\lambda_0)}^{-1}
 \DR_Y F_{\dagger}^i(\nbigm_2^{\lambda_0}),
\]
we obtain the desired coincidence.
The second claim can be shown similarly.
\hfill\qed

\vspace{.1in}
We obtain the following 
from the characterization of
$[\star H]$ as in Theorem 8.1.4 of \cite{mochi9}.
\begin{prop}
For any hypersurface $H$,
the localizations
$[\ast H]$ and $[!H]$
are compatible with 
$\Upsilontilde$.
\hfill\qed
\end{prop}

Let $\nbigv=(\nbigv_1,\nbigv_2,C) 
 \in\MTS^{\integral\adm}(X,H)$
with an integrable real structure
as a smooth integrable $\nbigr_{X(\ast H)}$-triple.
Then,
$\nbigv_{2|\cnum^{\ast}\times X}$
has a good $\real$-structure.
\begin{prop}
Let $\nbigt\in\MTMint(X,\real)$.
Let $\nbigv\in\MTS^{\integral\adm}(X,H)$
with an integrable real structure.
Then, we naturally have
$\Upsilontilde\bigl(
 (\nbigt\otimes\nbigv)[\star H]
\bigr)
\simeq
 \bigl(
 \Upsilontilde(\nbigt)
 \otimes
 \nbigv_2
 \bigr)[\star H]$
in $\Hol(\cnum^{\ast}\times X,\real)$.
\hfill\qed
\end{prop}

\begin{cor}
We naturally have
$\Upsilontilde\circ\Xi^{(a)}
\simeq \Xi^{(a)}\circ\Upsilontilde$,
$\Upsilontilde\circ\psi^{(a)}
\simeq
 \psi^{(a)}\circ\Upsilontilde$
and 
$\Upsilontilde\circ\phi^{(a)}
\simeq
 \phi^{(a)}\circ\Upsilontilde$.
\hfill\qed
\end{cor}

\subsubsection{Good real structure of
 integrable mixed twistor $D$-modules}

Let $\nbigt\in\MTMint(X,\real)$.
Let $\nbigm_i$ $(i=1,2)$ be the underlying 
integrable $\nbigr_X$-modules.
We obtain $D_{\cnum\times X}$-modules
$\nbigm_i(\ast \nbigx^0)$,
where $\nbigx^0:=\{0\}\times X$.
We say that the real structure of $\nbigt$ is good,
if $\nbigm_2(\ast \nbigx^0)$
has a $\real$-Betti structure
whose restriction to $\nbigm_{2|\cnum^{\ast}\times X}$
is equal to the one induced by
the real structure of $\nbigt$.
Let $\MTMint(X,\real)^{\good}
\subset\MTMint(X,\real)$
denote the full subcategory of
integrable mixed twistor $D$-modules
with real good structure.
By using the results in \cite{mochi9},
we obtain the following.

\begin{prop}
Let $\nbigt\in \MTMint(X,\real)^{\good}$.
\begin{itemize}
\item
For any projective morphism
$F:X\lrarr Y$,
we have
$F^i_{\dagger}(\nbigt)
 \in \MTMint(Y,\real)^{\integral}$.
\item
 $\DDD\nbigt\in\MTMint(X,\real)^{\good}$.
\item
For any hypersurface $H$ of $X$,
we have
$\nbigt[\star H]
 \in \MTMint(X,\real)^{\good}$.
\item
For any holomorphic function $g$ on $X$,
$\Xi^{(a)}_g(\nbigt)$,
$\psi^{(a)}_g(\nbigt)$
and 
$\phi^{(a)}_g(\nbigt)$
are objects in
$\MTMint(X,\real)^{\good}$.
\hfill\qed
\end{itemize}
\end{prop}

Let $K$ be a subfield of $\real$.
Let $\nbigt\in\MTMint(X,\real)$.
A good $K$-structure of $\nbigt$
is a $K$-Betti structure of
$\nbigm_2(\ast \nbigx^0)$
whose restriction to
$\cnum^{\ast}\times X$
induces the $\real$-Betti structure
of $\nbigm_{2|\cnum^{\ast}\times X}$.
If $\nbigt$ has a good $K$-structure,
we have $\nbigt\in\MTMint(X,\real)^{\good}$.

Let $\MTMint(X,K)^{\good}$ denote
the category of integrable mixed twistor $D$-modules
with good $K$-structure.
By using the results in \cite{mochi9},
we obtain the following.
\begin{prop}
Let $\nbigt\in \MTMint(X,K)^{\good}$.
\begin{itemize}
\item
For any projective morphism
$F:X\lrarr Y$,
we have
$F^i_{\dagger}(\nbigt)
 \in \MTMint(Y,K)^{\integral}$.
\item
 $\DDD\nbigt\in\MTMint(X,K)^{\good}$.
\item
For any hypersurface $H$ of $X$,
we have
$\nbigt[\star H]
 \in \MTMint(X,K)^{\good}$.
\item
For any holomorphic function $g$ on $X$,
$\Xi^{(a)}_g(\nbigt)$,
$\psi^{(a)}_g(\nbigt)$
and 
$\phi^{(a)}_g(\nbigt)$
are objects in
$\MTMint(X,K)^{\good}$.
\hfill\qed
\end{itemize}
\end{prop}

\section{A relation with mixed Hodge modules}

Let $P$ be a $\real$-perverse sheaf on $X$.
We have a regular holonomic $D$-module $M$
with an isomorphism
$\DR(M)\simeq P\otimes\cnum$.
A Hodge filtration of $P$ is 
a filtration $F$ of $M$ by coherent $\nbigo_X$-submodules
indexed by integers such that
$F_j(D_X)\,F_i(M)\subset F_{i+j}(M)$,
where $F_{\ast}(D_X)$ be the filtration by the order of
differential operators.
The category $MF_h(D_X,\real)$ 
of $\real$-perverse sheaves with a Hodge filtration
is naturally defined.
Objects are filtered $D_X$-module $(M,F)$
and $\real$-perverse sheaf $P$
with an isomorphism
$\DR_X(M)\simeq P\otimes\cnum$.
Morphisms are naturally defined.
Recall that filtered objects in $\MF_h(D_X,\real)$
are ingredients of mixed Hodge modules.

Let $(M,F;P)$ with $W$ be an object in 
$\MF_h(D_X,\real)^{\fil}$
which gives a mixed Hodge module.
Then, we have a naturally associated
filtered integrable $\nbigr$-triple with a real structure.
Indeed,
let $\nbigm$ be the $\nbigr_X$-module
obtained as the analytification of the Rees module
of $(M,F)$.
Note that $\DDD\nbigm$ is the analytification of
the Rees module of the dual of $(M,F)$
in the category of filtered $D$-modules.
We have the naturally defined hermitian pairing
$C:\DDD\nbigm_{|\vecS\times X}\times
\nbigmbar_{|\vecS\times X}
\lrarr\distribution_{\vecS\times X/X}$
induced by the real structure of $\DR(M)$.
Thus, we obtain an $\nbigr$-triple
$\nbigt=(\DDD\nbigm,\nbigm,C)$
with a naturally induced filtration $W$.
We set $\Phi(M,F,K,W):=(\nbigt,W)$.
\begin{lem}
For any projective morphism
$f:X\lrarr Y$
and any mixed Hodge module
$(M,F,K,W)$ on $X$,
we have a natural isomorphism
$\Phi \nbigh^if_{\ast}(M,F,K,W)
\simeq
 f^i_{\dagger}\Phi(M,F,K,W)$.
\end{lem}
\pf
It follows from the construction of the functors,
and the compatibility of 
the push-forward and the dual.
\hfill\qed

\begin{lem}
Let $a:(M_1,F,K,W)\lrarr(M_2,F,K,W)$
be a morphism of mixed Hodge modules.
Then, we have
$\Phi(\Ker(a))=\Ker\Phi(a)$
$\Phi(\Image(a))=\Image\Phi(a)$
and 
$\Phi(\Cok(a))=\Cok\Phi(a)$.
\end{lem}
\pf
It follows from that
$a:(M_1,F,W)\lrarr (M_2,F,W)$ is bi-strict.
\hfill\qed

\begin{lem}
Let $(M,F,K,W)$ be a mixed Hodge module on $X$.
\begin{itemize}
\item
If $(M,F,K,W)$ is pure,
$\Phi(M,F,K,W)$ is strictly $S$-decomposable.
\item
In general, $\Phi(M,F,K,W)$ is admissible specializable.
\item
 Let $g$ be a holomorphic function on $X$.
Let $j_g:X\setminus g^{-1}(0)\lrarr X$
 denote the inclusion.
If $j_{g\star}j_g^{-1}(M,F,K,W)\simeq (M,F,K,W)$,
then
 $\Phi(M,F,K,W)[\star g]\simeq
 \Phi(M,F,K,W)$.
\end{itemize}
\end{lem}
\pf
Let $R^F(M,F)$ denote the Rees module
associated to $(M,F)$.
Let $U$ be any open subset of $X$
with a holomorphic function $g$.
Let $i_g:U\lrarr U\times\cnum_t$
denote the graph.
The $V$-filtration of $i_{g\dagger}M$
naturally induces a $V$-filtration of
$R^Fi_{g\dagger}(M,F)$.
Because $i_{g\dagger}(M,F)$ is 
quasi-unipotent and regular along $\{t=0\}$
in the sense of \S3.2 of \cite{saito1},
we obtain that 
$t:V_{a}R^Fi_{g\dagger}(M,F)
\lrarr V_{a-1}R^Fi_{g\dagger}(M,F)$
are isomorphisms for any $a<0$,
and that
$\del_t:\Gr^V_aR^Fi_{g\dagger}(M,F[1])
\lrarr \Gr^V_{a+1}R^Fi_{g\dagger}(M,F)$
are isomorphisms for any $a>-1$.
By the construction,
$\Gr^V_aR^Fi_{g\dagger}(M,F)$
are flat over $\cnum[\lambda]$.
Hence, 
we obtain that $\nbigm$ is strictly specializable
along $g$.
By the construction,
$\DDD\nbigm$ is the analytification 
of the Rees module associated to
the dual $(\DDD M,F)$
of the filtered module $(M,F)$.
Hence,
$\DDD\nbigm$ is also strictly specializable
along $g$.
Moreover, if $(M,F,K,W)$ is pure,
we have 
$\Gr^V_0(R^F(M))=
\Image\can\oplus\Ker\var$,
where 
$\can:\Gr^V_{-1}R^Fi_{g\dagger}(M,F[1])
\lrarr \Gr^V_{0}R^Fi_{g\dagger}(M,F)$
and 
$\var:\Gr^V_{0}R^Fi_{g\dagger}(M,F)\lrarr
 \Gr^V_{-1}R^Fi_{g\dagger}(M,F)$
are induced by the action of
$t$ and $\del_t$, respectively.
Hence, $\nbigm$ is strictly $S$-decomposable
along $g$.
Similarly
$\DDD\nbigm$ is also strictly $S$-decomposable
along $g$.

For any mixed Hodge module $(M,F,K,W)$,
the filtrations $F$, $V$ and $W$ on $i_{g\dagger}M$
are compatible.
It means that
$V_aW_jR^Fi_{g\dagger}(M,F)
 \lrarr V_a\Gr^W_jR^Fi_{g\dagger}(M,F)$
are surjective for any $a$ and $j$.
It implies that we obtain
the filtered strictly specializability of $\nbigm$
along $g$.
We set
$\bigl(
 \psi_g(M),F
 \bigr)
:=\bigoplus_{-1\leq a<0}
 \Gr^V_a(M,F[1])$
and 
$\bigl(
 \phi_{g,1}(M),F
 \bigr)
:=\Gr^V_0(M,F)$.
We also set
$L_j\psi_g(M):=
 \psi(W_{j+1}M)$
and 
$L_j\phi_{g,1}(M):=
 \phi_{g,1}(W_jM)$.
Let $N$ denote the nilpotent part of the action of
$-\del_tt$.
By the condition for mixed Hodge modules,
there exists a relative monodromy filtration
$W$ on $(\psi_g(M),L)$
and $(\phi_{g,1}(M),L)$,
and $(\psi_g(M,F),W)$
and $(\phi_{g,1}(M,F),W)$
are mixed Hodge modules.
It implies that
the induced actions of 
$-\lambda\del_tt$
on $\bigoplus_{-1\leq a<0}
 \Gr^V_{a}R^Fi_{g\dagger}(M,F)$
and 
$\Gr^V_0R^Fi_{g\dagger}(M,F)$
with $L$,
have relative monodromy filtrations.
Together with the first claim,
we obtain that 
$(\nbigt,W)$ is admissibly specializable along $g$.

The third claim follows from the characterization
of $\Phi(M,F,W,K)[\star g]$.
\hfill\qed

\begin{prop}
$\Phi(M,F,K,W)$ is a mixed twistor $D$-module on $X$
for any mixed Hodge module $(M,F,K,W)$ on $X$.
\end{prop}
\pf
Let us consider the case 
$(M,F,K,W)$ is a polarizable pure Hodge module of weight $w$.
We may assume that it is irreducible.
There exists a polarizable variation of pure Hodge structure
$(\nbigh,F)$ on a Zariski open subset $U\subset \Supp(M)$,
and that 
$(M,F,K,W)$ is obtained as the minimal extension of
$(\nbigh,F)$.
Let $\nbigm_0$ be the $\nbigr_U$-module
obtained as the analytification of $R^F(\nbigh,F)$.
We have a polarizable variation of pure twistor structure
$(\DDD\nbigm_0,\nbigm_0,C)$ on $U$.
It is uniquely extended to 
a polarizable pure twistor $D$-module $\nbigt_1$
whose strict support is $\Supp(M)$.
Let $Z:=\Supp(M)\setminus U$.
We have a natural isomorphism
$\Upsilon:
 \nbigt_{1|X\setminus Z}
\simeq
 \Phi(M,F,K,W)_{|X\setminus Z}$.
We shall prove that the isomorphism
$\Upsilon$ is uniquely extended
to 
$\nbigt_1\simeq
 \Phi(M,F,K,W)$.
By the uniqueness,
we have only to check the claim
locally around any point of $\Supp(M)$.
We may assume that there exists
a holomorphic function $g$
such that 
$g^{-1}(0)\cap \Supp(M)=Z$.
By the strict $S$-decomposability of
$\nbigt_1$ and $\Phi(M,F,K,W)$,
we have only to prove that
$\nbigt_1(\ast g)
=\Phi(M,F,K,W)(\ast g)$.
We may assume that there exists
a complex manifold $Y$
with a projective birational morphism
$\varphi:Y\lrarr \Supp(M)$
such that
(i) $H_Y:=\varphi^{-1}(Z)$ is normal crossing
hypersurface,
(ii) $Y\setminus H_Y=U$.
Let $\nbightilde$ denote 
the $D_{Y}(\ast H_Y)$-module
obtained as the regular extension of $\nbigh$.
It is equipped with the filtration 
by locally free $\nbigo_{Y}(\ast H_Y)$-modules,
induced by $F$,
denoted by $\Ftilde$.
As the analytification of the Rees module
$R^{F}(\nbightilde,\Ftilde)$,
we have an $\nbigr_{Y(\ast H_Y)}$-module
$\nbigmtilde$.
Then, we can observe that
the underlying $\nbigr_{X}(\ast g)$-modules of
$\nbigt_1(\ast g)$ and 
$\Phi(M,F,K,W)(\ast g)$
are naturally isomorphic,
by the construction of
the prolongations of
variation of polarizable pure Hodge (twistor)
structure to polarizable pure Hodge (twistor)
modules.
Thus, we are done in the pure case.

\vspace{.1in}
We consider the mixed case.
We use an induction on $\dim M$.
If $\dim M=0$, the claim is easy.
We have already known
$\Phi(M,F,K,W)\in\MTW(X)$.
To check $\Phi(M,F,K,W)\in\MTM(X)$,
we have only to prove it locally around 
any point of $P\in\Supp(M)$.
We may assume that there exists a holomorphic function
$g$ on $X$ such that 
(i) $\dim\Supp(M)\cap g^{-1}(0)<\dim\Supp(M)$,
(ii) $\Supp(M)\setminus g^{-1}(0)$ is
 a complex submanifold of
 $X\setminus g^{-1}(0)$,
(iii) $(M,F,K,W)_{|X\setminus g^{-1}(0)}$
comes from an admissible variation of mixed Hodge structure
 $(\nbigh,F,W)$ on $\Supp(M)\setminus g^{-1}(0)$.
We may assume to have a complex manifold $Y$
with a projective birational morphism
$\varphi:Y\lrarr \Supp(M)$
such that
(i) $H_Y:=\varphi^{-1}(g^{-1}(0))$ is normal crossing,
(ii) $Y\setminus H_Y\simeq \Supp(M)\setminus g^{-1}(0)$.
We have an admissible variation of mixed Hodge structure
on $(Y,H_Y)$ given by $\varphi^{\ast}(\nbigh,F,W)$.
It induces an admissible variation of mixed twistor structure
$(\nbigv,W)$ on $(Y,H_Y)$.
As in the pure case, we obtain that
$\Phi(M,F,K,W)(\ast g)\simeq
 \varphi_{\dagger}(\nbigv,W)$.
Hence, by using the admissible specializability along $g$
and the characterization as in Lemma \ref{lem;11.2.20.5},
we have
$j_{g\star}j_g^{-1}\Phi(M,F,K,W)
\simeq
 \varphi_{\star}(\nbigv,W)$.
By using Beilinson's functor for mixed Hodge modules
as in \cite{saito-beilinson},
we have a description of
$(M,F,K,W)$ as the cohomology of the complex
in the category of mixed Hodge modules:
\[
 \psi_g(M,F,K,W)
\lrarr
 \Xi_g(M,F,K,W)
 \oplus
 \phi_g(M,F,K,W)
\lrarr
 \psi_g(M,F,K,W)(-1)
\]
By the construction of $\Xi_g(M,F,K,W)$,
there exist admissible variations of
mixed Hodge structure
$(\nbigv_i,W)$ $(i=1,2)$ on $(Y,H_Y)$
such that
$\Xi_g(M,F,K,W)$
is isomorphic
to the kernel of a morphism
$\varphi_{!}(\nbigv_1,W)
\lrarr
 \varphi_{\ast}(\nbigv_2,W)$.
Hence, we have
$\Phi\Xi_g(M,F,K,W)
\in\MTM(X)$.
Similarly,
we obtain that
$\Phi\psi_g(M,F,K,W)$
is a mixed twistor $D$-module.
By the hypothesis of the induction,
we have
$\Phi\phi_g(M,F,K,W)\in\MTM(X)$.
Then, we obtain that
$(M,F,K,W)\in\MTM(X)$.
\hfill\qed

\vspace{.1in}
The mixed twistor $D$-module $(\nbigt,W)$
has a natural real structure $\kappa=(\id,\id)$,
where we use the natural
identification $j^{\ast}\nbigm\simeq\nbigm$.
Thus, we obtain a functor
from the category of mixed Hodge modules
to $\MTM(X,\real)$.

\chapter{Derived category of algebraic mixed twistor $D$-modules}

In this chapter,
by following \cite{beilinson1} and \cite{saito2},
we give a summary for derived functors
on the derived category of mixed twistor $D$-modules
on algebraic varieties.
It also works for mixed twistor $D$-modules
with additional structures
such as 
integrable structure,
$\real$-structure,
good $K$-structure, etc..
See \cite{saito2} for more detailed arguments.

\section{Algebraic mixed twistor $D$-modules}

Let $X$ be a smooth complex 
quasi-projective variety.
We take a smooth projective completion
$X\subset \Xbar$ such that 
$D=\Xbar-X$ is a hypersurface.
We set
$\MTM(X):=\MTM(\Xbar,[\ast D])$,
which is independent of the choice of 
a projective completion $\Xbar$
by Proposition \ref{prop;11.3.30.10}.
\index{category $\MTM(X)$}
An object in $\MTM(X)$ is called
an algebraic mixed twistor $D$-module
on $X$.
Let $D^b(\MTM(X))$ denote
the derived category of 
bounded complexes in $\MTM(X)$.
\index{category $D^b(\MTM(X))$}
Similarly,
$\MTMint(X)$
and $D^b(\MTMint(X))$
are defined.
\index{category $\MTMint(X)$}
\index{category $D^b(\MTMint(X))$}

\section{Localization}

Let $H$ be a hypersurface of $X$.
We have the localizations
$[\ast H]:
 \MTM(X)\lrarr \MTM(X)$
and 
$[!H]:\MTM(X)\lrarr \MTM(X)$
which are exact functors.
They induce the exact functors
on $D^b\MTM(X)$.
By comparison of the Yoneda extensions,
for $\nbigt_i\in\MTM(X)$,
we have the following natural isomorphisms:
\[
 \Ext^i_{\MTM(X)}\bigl(\nbigt_1,\nbigt_2[\ast H]\bigr)
\simeq
 \Ext^i_{\MTM(X)}\bigl(\nbigt_1[\ast H],\nbigt_2[\ast H]\bigr)
\]
\[
 \Ext^i_{\MTM(X)}\bigl(\nbigt_1[!H],\nbigt_2\bigr)
\simeq
 \Ext^i_{\MTM(X)}\bigl(\nbigt_1[!H],\nbigt_2[!H]\bigr)
\]
We have similar localizations
in the integrable case.

\section{Beilinson functors}
Let $g$ be an algebraic function on $X$.
We have the exact functors
$\Pi^{a,b}_{g\star}$ on $\MTM(X)$
and $D^b\MTM(X)$,
where $\star=\ast,!$ and $a,b\in\seisuu$.
We naturally obtain the exact functors
$\Xi_g^{(a)}$, $\psi_g^{(a)}$
and $\phi_g^{(a)}$ on $\MTM(X)$
and $D^b\MTM(X)$.
We have similar functors in the integrable case.

\section{Dual and hermitian dual}

Let $\nbigt\in\MTM(X)$.
It is represented by
$\nbigtbar\in\MTM(\Xbar,[\ast \Dbar])$.
We set $\DDD_X\nbigt:=
\bigl(\DDD_{\Xbar}\nbigtbar
 \bigr)[\ast \Dbar]$,
which gives an object in $\MTM(X)$.
It is well defined
by the compatibility of push-forward
and dual.
Similarly,
we obtain $\DDD^{\herm}\nbigt$
in $\MTM(X)$.
We extend $\DDD$ and $\DDD^{\herm}$
on $\MTM(X)$ to the functors on $D^b\MTM(X)$
as in  \S\ref{subsection;11.3.30.12}.

\begin{itemize}
\item
We have natural isomorphisms
$\DDD\circ\DDD\simeq\id$
and $\DDD^{\herm}\circ\DDD^{\herm}=\id$.
We also have
$\DDD\circ\DDD^{\herm}
=\DDD^{\herm}\circ\DDD$.
\item
We have natural isomorphisms
$\DDD(\nbigt^{\bullet}[\ast H])
\simeq
 \DDD(\nbigt^{\bullet})[! H]$
and
$\DDD^{\herm}(\nbigt^{\bullet}[\ast H])
\simeq
 \DDD^{\herm}(\nbigt^{\bullet})[! H]$
\item
Let $g$ be any algebraic function on $X$.
We have natural isomorphisms
\[
\Xi_{g\ast }^{a,b}(\DDD\nbigt^{\bullet})
\simeq
 \DDD\Xi_{g\ast!}^{-b+1,-a+1}(\nbigt^{\bullet}),
\quad
 \Xi_{g\ast!}^{a,b}(\DDD^{\herm}\nbigt^{\bullet})
\simeq
 \DDD^{\herm}\Xi_{g\ast!}^{-b+1,-a+1}(\nbigt^{\bullet}).
\]
In particular, we have compatibility of
the dual and hermitian dual
with the nearby cycle functor,
the vanishing cycle functor,
and the maximal functor.
\end{itemize}
By comparing Yoneda extensions,
we have the following natural isomorphisms:
\[
 \Ext^i_{\MTM(X)}(\nbigt_1,\nbigt_2)
\simeq
 \Ext^i_{\MTM(X)}\bigl(\DDD_X\nbigt_2,
 \DDD_X\nbigt_1\bigr)
\]
\[
 \Ext^i_{\MTM(X)}(\nbigt_1,\nbigt_2)
\simeq
 \Ext^i_{\MTM(X)}\bigl(\DDD^{\herm}_X\nbigt_2,
 \DDD^{\herm}_X\nbigt_1\bigr).
\]

\section{Real structure}

As in \S\ref{subsection;11.3.30.12},
we define
$\gammatilde^{\ast}:=j^{\ast}\circ
 \DDD^{\herm}\circ\DDD$
on $D^b\MTM(X)$.
We have a natural isomorphism
$\gammatilde^{\ast}\circ\gammatilde^{\ast}
=\id$.

\begin{itemize}
\item
We have natural commutativity
$\gammatilde^{\ast}\circ\DDD
\simeq
 \DDD\circ\gammatilde^{\ast}$
and
$\gammatilde^{\ast}\circ\DDD^{\herm}
\simeq
 \DDD^{\herm}\circ\gammatilde^{\ast}$.
\item
We have natural isomorphisms
$\gammatilde^{\ast}(\nbigt^{\bullet}[\ast H])
\simeq
 \gammatilde^{\ast}(\nbigt^{\bullet})[\ast H]$.
\item
Let $g$ be an algebraic function on $X$.
We have natural isomorphisms
\[
 \Xi_{g\ast }^{a,b}(\gammatilde^{\ast}\nbigt^{\bullet})
\simeq
 \gammatilde^{\ast}\Xi_{g\ast!}^{a,b}(\nbigt^{\bullet}).
\]
In particular, we have compatibility of
$\gammatilde^{\ast}$
with the nearby cycle functor,
the vanishing cycle functor,
and the maximal functor.
\end{itemize}

A real structure of $\nbigt^{\bullet}\in D^b\MTM(X)$
is defined to be an isomorphism
$\kappa:\gammatilde^{\ast}\nbigt^{\bullet}
\simeq
 \nbigt^{\bullet}$
such that $\gammatilde^{\ast}\kappa\circ\kappa=\id$.
Let $D^b(\MTM(X),\real)$ be the full subcategory of
objects in $D^b\MTM(X)$
with a real structure.
A morphism $(\nbigt^{\bullet}_1,\kappa_1)\lrarr
 (\nbigt_2^{\bullet},\kappa_2)$
in $D^b(\MTM(X),\real)$ is
a morphism $\varphi:\nbigt_1^{\bullet}\lrarr\nbigt_2^{\bullet}$
in $D^b(\MTM(X))$
such that
$\varphi\circ \kappa_1
=\kappa_2\circ\gammatilde^{\ast}\varphi$.
The category $D^b(\MTM(X),\real)$
is naturally equipped with 
localization, Beilinson's functors,
$\DDD$ and $\DDD^{\herm}$.
The integrable case can be defined similarly.

Let us look at typical examples.
For simplicity, let $Z$ be a projective variety.
We have objects $\nbigu_Z(0,d_Z)[d_Z]$
and $\nbigu_Z(d_Z,0)[-d_Z]$ 
in $D^b\MTM(X)$.
We fix the isomorphism
$\mu_Z:\DDD\bigl(\lambda^d_Z\nbigo_{\nbigz}\bigr)
\simeq\nbigo_{\nbigz}$
whose restriction to $\lambda=1$
is given as in \S\ref{subsection;11.3.30.2}.
Then, we have natural real structures
$\bigl(\mu_Z,\DDD(\mu_Z)^{-1}\bigr)$
of $\nbigu_Z[d_Z]$,
and 
$\bigl(\mu_Z,\DDD(\mu_Z)^{-1}\bigr)$
of $\nbigu_Z[d_Z]$,
$\bigl(\DDD(\mu_Z),\mu_Z^{-1}\bigr)$
of $\nbigu_Z[-d_Z]$.

\section{External product}
Let $(\nbigt_i,L)\in\MTM(X_i)$ $(i=1,2)$.
It is easy to show the following lemma
by using Proposition \ref{prop;13.7.29.5}.
\begin{prop}
$\nbigt_1\boxtimes\nbigt_2$
with the naturally induced filtration $L$
is an object in $\MTM(X_1\times X_2)$.
As a result,
we obtain a bifunctor
$\boxtimes:
 \MTM(X_1)\times\MTM(X_2)
\lrarr 
 \MTM(X_1\times X_2)$,
compatible with 
the standard external products
$\boxtimes:
 \Hol(X_1)\times\Hol(X_2)
\lrarr
 \Hol(X_1\times X_2)$.
It is naturally extended to
\[
\boxtimes:
 D^b\MTM(X_1)\times D^b\MTM(X_2)
\lrarr 
 D^b\MTM(X_1\times X_2).
\]
The functor can be enriched with
real structures and integrable structures.
\hfill\qed
\end{prop}

\section{A version of Kashiwara's equivalence}

Let $A$ be a subvariety of $Y$,
which is not necessarily smooth.
Let $\MTM_A(Y)$ be the full subcategory of
mixed twistor $D$-modules on $Y$
whose supports are contained in $A$.
Let $D^{b}_{A}\MTM(Y)$ be
the full subcategory of $D^b\MTM(Y)$
which consists of the objects $\nbigt^{\bullet}$
such that the supports of the cohomology
$\bigoplus_i\nbigh^i\nbigt^{\bullet}$
are contained in $A$.
\begin{prop}
\label{prop;10.1.12.10}
The natural functor
$D^b\MTM_A(Y)\lrarr D^b_A\MTM(Y)$
is an equivalence.
We have similar equivalences
in the integrable case
and the real case.
\end{prop}
\pf
According to \cite{bbd},
we have only to check the following effaceability:
\begin{itemize}
\item
Let $\nbigt_i\in\MTM_A(Y)$.
For any $f\in \Ext^i_{\MTM(Y)}(\nbigt_1,\nbigt_2)$,
there exists a monomorphism
$\nbigt_2\lrarr \nbigt'$
in $\MTM_A(Y)$ such that
the image of $f$
in $\Ext^i_{\MTM(Y)}(\nbigt_1,\nbigt_2')$ is $0$.
\end{itemize}
We can show it by using the arguments
in Sections 2.2 and 2.2.1 in \cite{beilinson1}.
\hfill\qed

\section{Push-forward}

Let $f:X\lrarr Y$ be a projective morphism
of quasi-projective varieties.
We take a factorization
$X\subset \Xbar\stackrel{f'}\lrarr Y$
such that (i) $f'$ is projective,
(ii) $H=\Xbar-X$ is normal crossing.
We have a natural equivalence
between $\MTM\bigl(\Xbar,[\ast H]\bigr)$
and $\MTM(X)$.
Let $(\nbigtbar,\Lbar)\in
 \MTM\bigl(\Xbar,[\ast H]\bigr)$
correspond to $(\nbigt,L)\in\MTM(X)$.
According to Proposition \ref{prop;10.11.15.21},
we have
\[
  f^i_{\ast}(\nbigt,L):=
  (f')^i_{\dagger}(\nbigtbar,\Lbar)
 \in \MTM(Y),
\quad
 f^i_{!}(\nbigt,L):=
 (f')^i_{\dagger}\bigl(
 (\nbigtbar,\Lbar)[!H]
 \bigr)
\in\MTM(Y).
\]
They are independent of the choice of $\Xbar$
up to natural isomorphisms.
We obtain cohomological functors
$f^i_{\ast},
 f^i_{!}
 :\MTM(X)\lrarr \MTM(Y)$
for $i\in\seisuu$.
We can show the following proposition
by an argument in \cite{saito2}.

\begin{prop}
\label{prop;09.11.10.10}
For each $\star=!,\ast$,
there exists a functor of
triangulated categories
\[
 f_{\star}:
 D^b\MTM(X)\lrarr D^b\MTM(Y)
\]
such that
(i) it is compatible with the standard functor
$f_{\star}:D^b_{\hol}(X)\lrarr D^b_{\hol}(Y)$,
(ii) the induced functor 
$H^i(f_{\star}):\MTM(X)\lrarr \MTM(Y)$
is isomorphic to $f_{\star}^i$.
It is characterized by the property (i) and (ii)
up to natural equivalence.
We have similar functors
in the integrable case and the real case.
\end{prop}
\pf
We give only a remark.
Let $\nbigt\in\MTM(X)$.
We take ample hypersurfaces $H_i$ $(i=1,2)$
such that $H_i$ and $H_1\cap H_2$
are non-characteristic to $\nbigt$.
Then,
we have a natural isomorphism
$\nbigt[\ast H_1!H_2]
\simeq
 \nbigt[!H_2\ast H_1]$.
Indeed, 
we have the following morphisms:
\[
 \nbigt[\ast H_1!H_2]
\stackrel{a}{\lrarr}
 \bigl(
 \nbigt[\ast H_1!H_2]
 \bigr)[\ast H_1]
\stackrel{b}{\llarr}
 \nbigt[!H_2\ast H_1]
\]
By the assumption on $H_i$,
$\Xi_{\DR}(a)$
and $\Xi_{\DR}(b)$
are isomorphisms.
Hence, we obtain that
$a$ and $b$ are isomorphisms.
We also have
$f_{\ast}^i\bigl(
 \nbigt[\ast H_1!H_2]
 \bigr)=0$
unless $i=0$.
Then, we can construct
a resolution $\nbigt^{\bullet}$ of $\nbigt$
such that 
$f_{\ast}^i\bigl(
 \nbigt^{j}
 \bigr)=0$ unless $j=0$.
(See  \cite{mochi9}, for example.)
Then, we can construct the desired functor 
$f_{\ast}$.
We obtain $f_!$ in a similar way.
\hfill\qed

\section{Pull back}

As in \cite{saito2},
the pull back is defined to be
the adjoint of the push-forward.

\begin{prop}
\label{prop;09.11.10.12}
$f_{!}$ has the right adjoint $f^{!}$,
and
$f_{\ast}$ has the left adjoint $f^{\ast}$.
Thus, we obtain the following functors:
\[
 f^{\star}:
 D^b\MTM(Y)\lrarr
 D^b\MTM(X)
\quad\quad
 (\star=!,\ast)
\]
They are compatible with
the corresponding functors of
holonomic $\nbigd$-modules
with respect to the forgetful functor.

We have similar functors
in the integrable and the real case.
\end{prop}
\pf
We have only to follow the argument
in \cite{saito2}.
We have only to construct adjoint functors
in the cases 
(i) $f$ is a closed immersion,
(ii) $f$ is a projection $X\times Y\lrarr Y$.
We give only an indication
for the construction.

\vspace{.1in}

Let $f:X\lrarr Y$ be a closed immersion.
The open immersion $Y-X\lrarr Y$
is denoted by $j$.
Let $\nbigt^{\bullet}$ 
be a complex in $\MTM(Y)$.
Let $H_i$ $(i=1,\ldots,N)$
be sufficiently general ample hypersurfaces of $Y$
such that 
(i) $\nbigt^{\bullet}
 \lrarr \nbigt^{\bullet}[\ast H_i]$ are monomorphisms,
(ii) $\bigcap_{i=1}^N H_i=X$.
For any subset $I=(i_1,\ldots,i_m)
 \subset \{1,\ldots,N\}$,
let $\cnum_I$ be the subspace of
$\bigwedge^{m}\cnum^N$
generated by
$e_{i_1}\wedge\cdots \wedge e_{i_m}$,
where $e_i\in\cnum^N$ denotes
an element whose $j$-th entry is
$1$ $(j=i)$ or $0$ $(j\neq i)$.
We use a natural inner product of
$\cnum^N$ given b
$(e_i,e_j)=0$ $(i\neq j)$
and 
$(e_i,e_i)=1$.
For $I=I_0\sqcup\{i\}$,
the canonical morphism
$\nbigt^p[\ast H(I_0)]
\lrarr \nbigt^p[\ast H(I)]$,
the multiplication and the inner product
of $e_i$ 
induce
$\nbigt^p[\ast H(I_0)]
 \otimes\cnum_{I_0}
\lrarr
 \nbigt^p[\ast H(I)]\otimes\cnum_I$.
For $m\geq 0$,
we put
$\nbigc^{m}
 (\nbigt^{p},\ast\vecH):=
 \bigoplus_{|I|=m}
 \nbigt^p[\ast H(I)]\otimes\cnum_I$,
and we obtain the double complex
$\nbigc^{\bullet}
 (\nbigt^{\bullet},\ast \vecH)$.
The total complex is denoted by
$\Tot\nbigc^{\bullet}(\nbigt^{\bullet},\ast\vecH)$.
It is easy to observe that
the support of the cohomology of
$\Tot\nbigc^{\bullet}(\nbigt^{\bullet},\ast\vecH)$
is contained in $X$.
According to Proposition \ref{prop;10.1.12.10},
we obtain
$f^{!}\nbigt^{\bullet}:=
 \Tot\nbigc^{\bullet}
 (\nbigt^{\bullet},\ast \vecH)$
in $D^b\MTM(X)$.
Thus, we obtain a functor
$f^{!}:
 D^b\MTM(Y)\lrarr
 D^b\MTM(X)$.
Note that the underlying 
$\nbigd_Y$-complex
is naturally quasi-isomorphic to
$f^{!}\nbigm^{\bullet}$,
where $f^{!}$ is the left adjoint of
$f_{\dagger}:D^b_{\hol}(X)\lrarr D^b_{\hol}(Y)$.
It is easy to check that
$f^{!}$ is the right adjoint of $f_!$.

Let $\cnum_I^{\lor}$ denote
the dual of $\cnum_I$.
For $I=I_0\sqcup\{i\}$,
we have a natural morphism
$\nbigt[!H(I)]\otimes\cnum_I^{\lor}
\lrarr
 \nbigt[!H(I_0)]\otimes\cnum_{I_0}^{\lor}$.
For $m\geq 0$,
we set
\[
\nbigc^{-m}(\nbigt^{\bullet},!\vecH):=
 \bigoplus_{|I|=m}
 \nbigt^p[!H(I)]\otimes\cnum_I, 
\]
and we obtain the double complex
$\nbigc^{\bullet}(\nbigt^{\bullet},!\vecH)$.
As in the previous case,
the support of the cohomology of
$\Tot\nbigc^{\bullet}(\nbigt^{\bullet},!\vecH)$
is contained in $X$.
We set
$f^{\ast}(\nbigt^{\bullet}):=
 \Tot\nbigc^{\bullet}(\nbigt^{\bullet},!\vecH)
 \in D^b\MTM(X)$.
It is easy to check that
$f^{\ast}$ is the left adjoint of
$f_{\ast}=f_!$.
The constructions are compatible with the integrable structure
and the real structure.

\begin{lem}
\label{lem;11.1.24.11}
We put $d:=d_Y-d_X$.
We have natural isomorphisms
\[
f^{\ast}\nbigu_Y(p,q)\simeq
 \nbigu_X(p-d,q)[d],
\quad
f^{!}\nbigu_Y(p,q)\simeq
 \nbigu_X(p,q-d)[-d].
\]
\end{lem}
\pf
We have 
$L^if^!\nbigu_Y(p,q)=0$
unless $i=d$,
and 
$L^if^{\ast}\nbigu_Y(p,q)=0$
unless $i=-d$.
We set
$\nbigt_0:=L^df^!\nbigu_Y(p,q)$.
We shall construct
$\nbigt_0\simeq
 \nbigu_X(p,q-d)$.
We obtain isomorphisms for the other
as the Hermitian adjoint.

Let us consider the case that
$X$ is a smooth hypersurface of $Y$.
In this case,
$\nbigt_0\simeq
 \Cok(\nbigu(p,q)\lrarr \nbigu(p,q)[\ast X])$.
Let $Y_1$ be an open subset of $Y$
in the classical topology
with a holomorphic coordinate
$(z_1,\ldots,z_n)$
such that $X_1:=X\cap Y_1=\{z_1=0\}$.
We have an isomorphism
$\nbigo_{\nbigy_1}[\ast \nbigx_1]\big/
 \nbigo_{\nbigy_1}
\simeq
 i_{\dagger}\lambda^{-1}\nbigo_{\nbigx_1}$
for which
$z_1^{-1}\longleftrightarrow
 \lambda^{-1}1\cdot (dz_1/\lambda)^{-1}$.
We also have an isomorphism
$\Ker\bigl(
 \nbigo_{\nbigy_1}[!\nbigx_1]
\lrarr
 \nbigo_{\nbigy_1}
 \bigr)
\simeq
 i_{\dagger}\nbigo_{\nbigx_1}$
for which
$-\deldel_{z_1}\otimes 1
\longleftrightarrow
 1\cdot(dz_1/\lambda)^{-1}$.
They give an isomorphism of
$\nbigr$-triples.
(See \S\ref{subsection;10.12.19.2}.)
It is independent of the choice of a coordinate.
Hence, we can glue the isomorphisms
for varied $(Y_1,z_1,\ldots,z_n)$,
and we obtain a global isomorphism
$\nbigt_0\simeq
 \nbigu_X(p,q-1)$.

\vspace{.1in}

Let us consider the general case.
We take a Zariski open subset $Y_2\subset Y$
with hypersurfaces $H_1,\ldots,H_{d}$
such that
$\bigcup H_i$ is a normal crossing,
and that
$\bigcap_{i=1}^{d}H_i=Y_2\cap X=:X_2$.
By using the result in the previous paragraph
successively,
we obtain an isomorphism
$\nbigt_{0|X_2}\simeq
 \nbigu_{X}(p,q-d)_{|X_2}$.
Because it is given on an Zariski open subset of $X$,
it is extended to an isomorphism
$\nbigt_{0}\simeq
 \nbigu_{X}(p,q-d)$,
denoted by
$\varphi_{Y_2,H_1,\ldots,H_{d}}$.

Let us observe that it is independent of 
the choice of $Y_2$ and $H_1,\ldots,H_{d}$.
Note that the isomorphism is
already determined up to constant multiplications.
First, we shall observe the independence
of the order of $H_1,\ldots,H_{d}$.
We set $\underline{d}:=\{1,\ldots,d\}$.
We set
\[
\nbigc^j:=
 \bigoplus_{\substack{I\subset\underline{d}\\
 |I|=j}}
 \nbigo_{\nbigy_2}(\ast H(I))
\otimes\cnum_I
\]
The inclusion
$\nbigo_{\nbigy_2}[\ast H(I)]
\lrarr 
\nbigo_{\nbigy_2}(\ast H(I\sqcup{i}))$
and the multiplication of $e_i$
induces a map
$\delta:\nbigc^j\lrarr\nbigc^{j+1}$.
Thus, we obtain a complex
$(\nbigc^{\bullet},\delta)$.
It is quasi-isomorphic to
the $d$-th cohomology sheaf
$\nbigh^d(\nbigc^{\bullet})$.

We set $\nbigm_0:=\nbigo_{\nbigy_2}$.
Inductively,
we set
$\nbigm_i:=\nbigm_{i-1}[\ast H_i]$.
We have
$\nbigo_{\nbigy_2}[\ast H(\underline{d})]
=\nbigm_{d}$.
By the construction,
$\nbigo_{\nbigy_2}[\ast H(\underline{d})]$
is the $\nbigr_{Y_2}$-submodule
in $\nbigo_{\nbigy_2}(\ast \nbigh(\underline{d}))$
generated by
$\nbigo_{\nbigy_2}(\nbigh(\underline{d}))$.
We have a unique $\nbigr_{Y_2}$-homomorphism
\[
 \nbigo_{\nbigy_2}[\ast H(\underline{d})]
 \otimes \bigwedge^d\cnum^d
\lrarr
 \iota_{\dagger}\lambda^{-d}\nbigo_{\nbigx_2}
\]
for which 
$(z_1\cdots z_{d})^{-1}
 \otimes e_1\wedge\cdots e_d
\longmapsto
 \lambda^d\cdots
 (dz_1/\lambda)^{-1}\cdot
 (dz_d/\lambda)^{-1}$.
It induces an isomorphism
$\nbigh^d(\nbigc^{\bullet})
\simeq 
 \iota_{\dagger}\lambda^{-d}\nbigo_{\nbigx_2}$.
By using this isomorphism,
we can observe that 
$\varphi_{Y_2,H_1,\ldots,H_{d}}$
is independent of the order of $H_1,\ldots,H_d$.
Next, if $H_i=H_i'$ $(i=1,\ldots,d-1)$,
then we have
$\varphi_{Y_2,H_1,\ldots,H_d}=\varphi_{Y_2,H_1',\ldots,H_d'}$
by the construction.
We can also observe that,
if $Y_3\subset Y_2$
and $H_i'=H_i\cap Y_3$,
then we have
$\varphi_{Y_2,H_1,\ldots,H_d}=\varphi_{Y_3,H_1',\ldots,H_d'}$.
Then, we can easily obtain that
$\varphi_{Y_2,H_1,\ldots,H_d}$
is independent of $Y_2$ and $H_1,\ldots,H_d$.
\hfill\qed

\vspace{.1in}

Let us consider the case that
$f$ is the projection of
$X=Y\times Z$ to $Y$.
We put $d_Z:=\dim Z$.
We set
\[
 f^{\ast}\nbigt:=
\nbigt\boxtimes
\nbigu_Z(d_Z,0)[-d_Z],
\quad
 f^!\nbigt=\nbigt\boxtimes
 \nbigu(0,d_Z)[d_Z].
\]
Let us show that
$f^{\ast}$ is the left adjoint of $f_{\ast}$.
We have only to construct
natural transformations
$\alpha:\id\lrarr  
 f_{\ast}f^{\ast}$
and
$\beta:f^{\ast}f_{\ast}\lrarr \id$
such that
{\small
\begin{equation}
\label{eq;11.1.24.12}
 \beta\circ f^{\ast}\alpha:
  f^{\ast}\nbigt^{\bullet}
\lrarr
 f^{\ast}f_{\ast}
 f^{\ast}\nbigt^{\bullet}
\lrarr
 f^{\ast}\nbigt^{\bullet},
\quad\quad
 f_{\ast}\beta\circ\alpha:
 f_{\ast}\nbign^{\bullet}
\lrarr
 f_{\ast}f^{\ast}f_{\ast}
 \nbign^{\bullet}
\lrarr
 f_{\ast}\nbign^{\bullet}
\end{equation}
}
are the identities.
We have the following natural morphisms:
{\small
\[
 (\tr,a_Z^{\ast}):
 \nbigu_{\pt}(0,0)
\lrarr
a_{Z\ast}\nbigu_Z(d_Z,0)[-d_Z],
\quad\quad
 (a_Z^{-1},\tr):
a_{Z!}\nbigu_Z(0,d_Z)[d_Z]
\lrarr
 \nbigu_{\pt}(0,0)
\]
}
Here,
the morphism
$\tr:a_{Z!}\nbigo_{\nbigz}[d_Z]\lambda^{d_Z}\lrarr
 \nbigo_{\cnum_{\lambda}}$
is given by the trace map,
and 
$a_Z^{-1}:
 \nbigo_{\cnum_{\lambda}}\lrarr
a_{Z\ast}\nbigo_{\nbigz}[-d_Z]$
is given by the pull back.
In particular,
we obtain a natural transform
$\alpha:\id\lrarr f_{\ast}f^{\ast}$.
Note that the morphisms are compatible with
natural integrable structures and real structures.
For the construction of $\beta$,
the following diagram is used:
\[
 \begin{CD}
 Z\times Y @>{i}>>
 Z\times Z\times Y @>{q_1}>> Z\times Y\\
 @. @V{q_2}VV @V{p_1}VV \\
 @. Z\times Y @>{p_2}>> Y
 \end{CD}
\]
Here, $i$ is induced by the diagonal 
$Z\lrarr Z\times Z$,
$q_j$ are induced by the projection
$Z\times Z\lrarr Z$ onto the $j$-th component,
and $p_j$ are the projections.
We have the following morphisms of complexes
of mixed twistor $D$-modules:
\begin{equation}
\label{eq;09.11.10.21}
  f^{\ast}f_{\ast}\nbigt^{\bullet}
=
 p_2^{\ast}p_{1\ast}\nbigt^{\bullet}
\simeq
 q_{2\ast}q_1^{\ast}\nbigt^{\bullet}
\lrarr
 q_{2\ast}
 \bigl(i_{\ast}i^{\ast}
 q_1^{\ast}\nbigt^{\bullet}
 \bigr)
\simeq
 i^{\ast}q_1^{\ast}
 \nbigt^{\bullet}
\end{equation}
According to Lemma \ref{lem;11.1.24.11},
we have a natural isomorphism
$i^{\ast}\nbigu_{Z\times Z}(d_Z,0)[-d_Z]
\simeq
 \nbigu_Z(0,0)$.
It induces 
$i^{\ast}q_1^{\ast}\nbigt^{\bullet}
\simeq
 \nbigt^{\bullet}$
in $D^b\MTM(Z\times Y)$.
We define $\beta$ as the composite of
(\ref{eq;09.11.10.21})
with the isomorphism.
Note that this construction is compatible 
with that in the case of $D$-modules.
Then, we obtain that 
the transformations in (\ref{eq;11.1.24.12}) 
are the identities,
because the transforms for the underlying
$D$-modules are the identity.
The construction is compatible with
real structures and integrable structures.
\hfill\qed

\section{Tensor and inner homomorphism}
Let $X$ be an algebraic variety.
Let $\delta_X:X\lrarr X\times X$
be the diagonal morphism.
We obtain the functors
$\otimes$ and $\nrhom$
on $D^b\MTM(X)$
in the standard ways:
\[
 \nbigt_1\otimes\nbigt_2:=
 \delta_X^{\ast}\bigl(
 \nbigt_1\boxtimes\nbigt_2
 \bigr),
\quad
 \nrhom(\nbigt_1,\nbigt_2):=
 \delta_X^!\bigl(
 \DDD_X\nbigt_1\boxtimes\nbigt_2
 \bigr)
\]
They are compatible with the corresponding
functors on $D^b_{\hol}(X)$.
They can be enriched with
real and integrable structures.

\chapter{$D$-triples and their functoriality}
\label{section;11.4.6.10}

We study hermitian pairing of holonomic $D$-modules.
The main purpose is to establish
the compatibility with the dual functor of 
holonomic $D$-modules
in \S\ref{subsection;13.4.12.30}.
We also argue real structure of holonomic $D$-modules.
This chapter is a preparation for \S\ref{section;11.4.9.20}.

\section{$D$-triples and their push-forward}
\label{subsection;11.3.23.1}

\subsection{$D$-triple and $D$-complex-triple}

\index{$D$-triple}
\index{$D$-complex-triple}

We introduce the notions of $D$-triples
and $D$-complex-triples,
which are variants of $\nbigr$-triples in \cite{sabbah2}.
Let $X$ be a complex manifold.
We set $D_{X,\Xbar}:=D_X\otimes_{\cnum}D_{\Xbar}$,
which is naturally a sheaf of algebras.
Let $M_i$ $(i=1,2)$ be $D_X$-modules.
A hermitian pairing of $M_1$ and $M_2$
is a $D_{X,\Xbar}$-homomorphism
$C:M_1\otimes_{\cnum}\Mbar_2
 \lrarr
 \distribution_X$.
Such a tuple $(M_1,M_2,C)$ is called a $D_X$-triple.
A morphism of $D_X$-triples
$(M_1',M_2',C')\lrarr (M_1,M_2,C)$
is a pair of morphisms
$\varphi_1:M_1\lrarr M_1'$ and
$\varphi_2:M_2'\lrarr M_2$
such that
$C'\bigl(\varphi_1(m_1),\overline{m_2'}\bigr)=
  C\bigl(m_1,\overline{\varphi_2(m_2')}\bigr)$.
Let $\Dtriplecat(X)$ denote the category of $D_X$-triples.
\index{category $\Dtriplecat(X)$}
It is an abelian category.
A $D_X$-triple $(M_1,M_2,C)$ is called
coherent (good, holonomic, etc.),
if the underlying $D_X$-modules are 
coherent (good, holonomic, etc.).

\subsubsection{$D$-complex-triple}

Let $M_i^{\bullet}$ $(i=1,2)$
be bounded complexes of $D_X$-modules.
A hermitian pairing of $M_1^{\bullet}$ and $M_2^{\bullet}$
is a morphism of 
$D_X\otimes D_{\Xbar}$-complexes
$C:
 \Tot\bigl(
 M_1^{\bullet}
\otimes
 \Mbar_2^{\bullet}
 \bigr)
\lrarr
 \distribution_X$.
Namely, a tuple of morphisms
$C^p:M_1^{-p}\otimes \Mbar_2^p\lrarr
 \distribution_X$
such that
\[
 C^p(dx^{-p-1},y^p)
+(-1)^{p+1}C_{p+1}(x^{-p-1},dy^{p})
=0.
\]
Such $(M_1^{\bullet},M_2^{\bullet},C)$
is called a $D_X$-complex-triple.
A morphism of $D_X$-complex-triples
\[
 (M_1^{\bullet},M_2^{\bullet},C_M)
\lrarr
 (N_1^{\bullet},N_2^{\bullet},C_N)
\]
is a pair of morphisms of $D_X$-complexes
$\varphi_1:N_1^{\bullet}\lrarr M_1^{\bullet}$
and $\varphi_2:M_2^{\bullet}\lrarr N_2^{\bullet}$
such that
$C_M(\varphi_1(x),y)
=C_N(x,\varphi_2(y))$.
Let $\Dcomplextriplecat(X)$ denote the category of
$D_X$-complex-triples.
\index{category $\Dcomplextriplecat(X)$}
It is an abelian category.
A morphism in $\Dcomplextriplecat(X)$
is called a quasi-isomorphism,
if the underlying morphisms of $D_X$-complexes
are quasi isomorphisms.

\vspace{.1in}
Let $\gbigt=(M_1^{\bullet},M_2^{\bullet},C)\in\Dcomplextriplecat(X)$.
Let $H^j(M_i^{\bullet})$ be the $j$-th cohomology of 
the complexes $M_i^{\bullet}$.
We have the induced $D_X$-triple
$\bigl(H^{-j}(M_1^{\bullet}),
 H^j(M_2^{\bullet}),
 H^j(C)
 \bigr)$
denoted by $H^j(\gbigt)$.
We also have the induced $D_X$-complex-triple
\[
\bigl(
 H^{-j}(M_1^{\bullet})[j],
 H^j(M_2^{\bullet})[-j],
 H^j(C)
 \bigr)
\]
denoted by $\vecH^j(\gbigt)$.
\index{$D$-triple $H^j(\gbigt)$}
\index{$D$-complex-triple $\vecH^j(\gbigt)$}

\vspace{.1in}
For any integer $\ell$,
we put $\epsilon(\ell):=(-1)^{\ell(\ell-1)/2}$.
We set
\[
 \nbigs_{\ell}(\gbigt):=
 (M^{\bullet}_1[-\ell],M^{\bullet}_2[\ell],C[\ell])
\]
where
$C[\ell]^p(x^{-p-\ell},y^{p+\ell})
:=(-1)^{\ell p}\epsilon(\ell)\,C^{\ell+p}(x^{-p-\ell},y^{p+\ell})$
for $x^{-p-\ell}\in M_1^{-p-\ell}=M_1[-\ell]^{-p}$
and $y^{p+\ell}\in M_2^{p+\ell}=M_2[\ell]^p$.
It is called the shift functor.
\index{functor $\nbigs_{\ell}$}
\index{pairing $C[\ell]$}

\vspace{.1in}

Let $(M^{\bullet}_1,M^{\bullet}_2,C_M)$ and
$(N^{\bullet}_1,N^{\bullet}_2,C_N)$
be $D$-complex-triples.
Let $\varphi_i:M_i^{\bullet}\lrarr N_i^{\bullet}$
be quasi-isomorphisms.
We say that $C_M$ and $C_N$ are the same
under the quasi-isomorphisms,
if $C_M(m_1,\overline{m}_2)
=C_N\bigl(\varphi_1(m_1),\overline{\varphi_2(m_2)}\bigr)$.
In this case,
we have the following natural quasi-isomorphisms
of $D$-complex-triples:
\[
\begin{CD}
 (M^{\bullet}_1,M^{\bullet}_2,C_M)
 @>>>
 (M^{\bullet}_1,N^{\bullet}_2,C')
 @<<<
 (N^{\bullet}_1,N^{\bullet}_2,C_N)
\end{CD}
\]
Here,
$C'(m_1,\overline{n}_2)
=C_N(\varphi_1(m_1),n_2)$.

\subsubsection{Complex of $D$-triples}

A complex of $D_X$-triples consists of
$D_X$-triples $\nbigt^{p}$ $(p\in\seisuu)$
and morphisms
$\delta^p:\nbigt^p\lrarr\nbigt^{p+1}$
such that $\delta^{p+1}\circ\delta^p=0$.
It is described by
$\nbigt^p=(M^{-p}_1,M^p_2,C^p)$
and
$\delta^p=(\delta^p_1,\delta^p_{2})$,
where
$\delta^p_{1}:M^{-p-1}_1\lrarr M^{-p}_1$
and
$\delta^p_{2}:M^p_2\lrarr M^{p+1}_2$.
They give complexes
$M^{\bullet}_i$,
and satisfy
\begin{equation}
 \label{eq;11.3.21.1}
 C^{p+1}(x^{-p-1},\delta^p_{2} y^p)
=C^p(\delta^p_{1} x^{-p-1},y^p)
\end{equation}
Let $\nbigc\bigl(\Dtriplecat(X)\bigr)$
denote the category of
bounded complexes of $D_X$-triples.
\index{category $\nbigc\bigl(\Dtriplecat(X)\bigr)$}
For $\nbigt^{\bullet}\in\nbigc(\Dtriplecat(X))$,
its $j$-th cohomology is denoted by
$H^j\bigl(\nbigt^{\bullet}\bigr)$.
For any integer $\ell$,
we define the shift $\nbigt^{\bullet}[\ell]$
in the standard way.
\index{shift $\nbigt^{\bullet}[\ell]$}
Namely,
we set $\nbigt[q]^p:=\nbigt^{p+q}$,
and the differentials 
$\deltatilde^p:\nbigt[q]^p\lrarr 
 \nbigt[q]^{p+1}$
are given by $(-1)^q\delta^{p+q}$.

\vspace{.1in}

We put $\epsilon(p):=(-1)^{p(p-1)/2}$ for integers $p$.
For $\nbigt^{\bullet}\in\nbigc\bigl(\Dtriplecat(X)\bigr)$
with the above description,
we define an object 
$\Psi_1(\nbigt^{\bullet}):=
 \bigl(M_1^{\bullet},M_2^{\bullet},
\Ctilde
\bigr)$
in $\Dcomplextriplecat(X)$,
where
$\Ctilde^p:=\epsilon(p)\,C^p$.
Thus, we obtain a functor
\[
 \Psi_1:\nbigc(\Dtriplecat(X))
\lrarr
 \Dcomplextriplecat(X).
\]
\index{functor $\Psi_1$}
It is easy to check the following.
\begin{prop}
$\Psi_1$ is an equivalence.
We have
$\Psi_1(\nbigt^{\bullet}[\ell])
\simeq
 \nbigs_{\ell}(\Psi_1(\nbigt^{\bullet}))$
naturally.
We also have
$\vecH^j\Psi_1(\nbigt^{\bullet})
\simeq
 \Psi_1\bigl(
 H^j(\nbigt^{\bullet})[-j]
 \bigr)$.
\hfill\qed
\end{prop}

\subsubsection{$D$-double-complex-triple and total complex}
\index{$D$-double-complex-triple}

A $D_X$-double-complex-triple is
a tuple
$\bigl(M_1^{\bullet,\bullet},M_2^{\bullet,\bullet},C\bigr)$:
\begin{itemize}
\item
 $M_i^{\bullet,\bullet}$ $(i=1,2)$
 are double complexes of $D_X$-modules,
 i.e.,
 they are $\seisuu^2$-graded $D$-modules
 $\{M_i^{\vecp}\,\,|\,\vecp\in\seisuu\}$
with morphisms
$d_j:M_i^{\vecp}\lrarr M_i^{\vecp+\vecdelta_j}$
$(i=1,2,j=1,2)$
such that
$d_j\circ d_j=0$
and $d_j\circ d_{m}=d_m\circ d_j$.
Here, $\vecdelta_1=(1,0)$
and $\vecdelta_2=(0,1)$.
For simplicity, we assume the boundedness.
\item
$C:M_1^{\bullet,\bullet}
 \otimes \Mbar_2^{\bullet,\bullet}\lrarr \distribution_X$
be a morphism of
$D_X\otimes D_{\Xbar}$-double-complexes,
i.e.,
a tuple of $D_X\otimes D_{\Xbar}$-morphisms
$C^{\vecp}:
 M_1^{-\vecp}\otimes M_2^{\vecp}
\lrarr
 \distribution_X$ such that
\[
 C^{\vecp}(d_1x^{-\vecp-\vecdelta_1},y^{\vecp})
+(-1)^{p_1+1}
 C^{\vecp+\vecdelta_1}
 (x^{-\vecp-\vecdelta_1},d_1y^{\vecp})=0
\]
\[
  C^{\vecp}(d_2x^{-\vecp-\vecdelta_2},y^{\vecp})
+(-1)^{p_2+1}
 C^{\vecp+\vecdelta_2}
 (x^{-\vecp-\vecdelta_2},d_2y^{\vecp})=0
\]
\end{itemize}
Let $\Ddoublecomplextriplecat(X)$ denote the category of
$D_X$-double-complex-triples.
\index{category $\Ddoublecomplextriplecat(X)$}

Let $(M_1^{\bullet,\bullet},M_2^{\bullet,\bullet},C)$ be 
an object in $\Ddoublecomplextriplecat(X)$.
We define the total complex object
$\Tot(M_1^{\bullet,\bullet},M_2^{\bullet,\bullet},C)$
in $\Dcomplextriplecat(X)$.
The underlying $D$-complexes are
the total complexes
$\Tot(M_i^{\bullet,\bullet})$,
i.e.,
\[
 \Tot(M_i^{\bullet,\bullet})^p
=\bigoplus_{p_1+p_2=p}M_i^{\vecp},
\]
with the differential
$d \vecx^{\vecp}
=d_1x^{\vecp}+(-1)^{p_1}d_2x^{\vecp}$.
The pairings are given by
\begin{equation}
 \label{eq;13.4.2.10}
  \Ctilde^{p}=\bigoplus_{p_1+p_2=p}
 (-1)^{p_1p_2} C^{\vecp}.
\end{equation}
Let us show that they give a $D_X$-complex-triple.
We have only to check (\ref{eq;11.3.21.1}).
We have
\begin{multline}
\label{eq;11.3.21.2}
 \Ctilde(dx^{-\vecp-\vecdelta_1},y^{\vecp})
+(-1)^{p_1+p_2+1}
 \Ctilde(x^{-\vecp-\vecdelta_1},dy^{\vecp})
 \\
=(-1)^{p_1p_2}
 C^{\vecp}(d_1x^{-\vecp-\vecdelta_1},y^{\vecp})
+(-1)^{p_1+p_2+1}
 (-1)^{(p_1+1)p_2}C^{\vecp+\vecdelta_1}
 (x^{-\vecp-\vecdelta_1},d_1y^{\vecp})
\\
=(-1)^{p_1p_2}\bigl(
 C^{\vecp}(d_1x^{-\vecp-\vecdelta_1},y^{\vecp})
+(-1)^{p_1+1}
 C^{\vecp+\vecdelta_1}(x^{-\vecp-\vecdelta_1},d_1y^{\vecp})
 \bigr)
=0
\end{multline}
We also have
\begin{multline}
\label{eq;11.3.21.3}
\Ctilde^{\vecp}(dx^{-\vecp-\vecdelta_2},y^{\vecp})
+(-1)^{p_1+p_2+1}
 \Ctilde^{\vecp+\vecdelta_2}(x^{-\vecp-\vecdelta_2},dy^{\vecp})
 \\
=(-1)^{p_1p_2}
 C^{\vecp}(d_2x^{-\vecp-\vecdelta_2},y^{\vecp})
+(-1)^{p_1+p_2+1}(-1)^{p_1(p_2+1)}
 C^{\vecp+\vecdelta_2}(x^{-\vecp-\vecdelta_2},d_2y^{\vecp})
 \\
=(-1)^{p_1p_2}
 \bigl(
 C^{\vecp}(d_2x^{-\vecp-\vecdelta_2},y^{\vecp})
+(-1)^{p_2+1}C^{\vecp+\vecdelta_2}(x^{-\vecp-\vecdelta_2},
 d_2y^{\vecp})
 \bigr)=0
\end{multline}
We can easily deduce (\ref{eq;11.3.21.1})
from (\ref{eq;11.3.21.2}) and (\ref{eq;11.3.21.3}).
Thus,  we obtain a natural functor
\[
 \Tot:\Ddoublecomplextriplecat(X)\lrarr
 \Dcomplextriplecat(X).
\]

\subsection{The push-forward}
\label{subsection;11.3.22.1}

\index{push-forward $f^{(0)}_{\dagger}$}
\index{push-forward $f_{\dagger}$}

Let $\nbigt=(M_1,M_2,C)$
be a $D_X$-triple.
Let $f:X\lrarr Y$ be a morphism of complex manifolds
such that 
the restriction of $f$ to the support of $\nbigt$
is proper.
As in the case of $\nbigr$-triples \cite{sabbah2},
we shall construct 
$f^{(0)}_{\dagger}\nbigt=
 (f_{\dagger}M_1,f_{\dagger}M_2,f_{\dagger}^{(0)}C)$
in $\Dcomplextriplecat(Y)$
and correspondingly
$f_{\dagger}\nbigt=
(f_{\dagger}M_1,f_{\dagger}M_2,f_{\dagger}C)$
in $\nbigc(\Dtriplecat(Y))$
i.e.,
$\Psi_1\bigl(f_{\dagger}\nbigt\bigr)
=f_{\dagger}^{(0)}\nbigt$.
We remark that
this is a specialization of 
the push-forward for $\nbigr$-triples
given in \cite{sabbah2}.
We shall compare it with a more naive construction
in \S\ref{subsection;13.4.5.10}.

\subsubsection{Closed immersion}
If $f$ is a closed immersion,
we have $f_{\dagger}M_i=
 \omega_X\otimes f^{-1}(D_Y\otimes\omega_Y^{-1})
 \otimes_{D_X}M_i$.
Let $\eta_X$ and $\eta_Y$ denote local generators
of $\omega_X$ and $\omega_Y$,
respectively.
We put $d_X:=\dim X$ and $d_Y:=\dim Y$.
We set
\begin{multline}
 f_{\dagger}^{(0)}C\Bigl(
 (\eta_X/\eta_Y)\cdot m_1,
 \overline{(\eta_X/\eta_Y)\cdot m_2}
 \Bigr)
:=  \\
\frac{1}{\eta_Y\etabar_Y}
 f_{\ast}\bigl( \eta_X\etabar_X
 \cdot C(m_1,\overline{m}_2) \bigr)
\left(
 \frac{1}{2\pi\sqrt{-1}}
\right)^{d_X-d_Y}
 \epsilon(d_X)\,\epsilon(d_Y)
\end{multline}
Namely, for a test form $\varphi=\phi\cdot\eta_Y\etabar_Y$,
we define
\begin{multline}
 \Bigl\langle
  f_{\dagger}^{(0)}C\Bigl(
 (\eta_X/\eta_Y)\cdot m_1,
 \overline{(\eta_X/\eta_Y)\cdot m_2}
 \Bigr),\,
 \varphi
 \Bigr\rangle
:= \\
 \Bigl\langle
 C(m_1,\overline{m}_2),\,
 f^{\ast}\phi\,\eta_X\etabar_X
 \Bigr\rangle
\left(
 \frac{1}{2\pi\sqrt{-1}}
\right)^{d_X-d_Y}
 \epsilon(d_X)\,\epsilon(d_Y)
\end{multline}
In this case, we have
$f_{\dagger}C=f^{(0)}_{\dagger}C$.

Let $f:X\lrarr Y$ and $g:Y\lrarr Z$
be closed immersions.
We have natural isomorphisms
$(g\circ f)_{\dagger}M_i
\simeq
 g_{\dagger}\bigl(f_{\dagger}M_i\bigr)$
$(i=1,2)$.
The following lemma is easy to see.
\begin{lem}
We have
$(g\circ f)^{(0)}_{\dagger}(C)
=g_{\dagger}^{(0)}(f_{\dagger}^{(0)}C)$
under the isomorphisms.
\hfill\qed
\end{lem}

\begin{rem}
If $Y=\Delta^n$, $X=\{z_{\ell+1}=\cdots z_n=0\}$,
$\eta_X=dz_1\cdots dz_{\ell}$
and $\eta_Y=dz_1\cdots dz_{n}$,
we have
\[
 \eta_X\,\etabar_X=
 \eta_Y\,\etabar_Y\Big/
 \prod_{i=\ell+1}^ndz_i\,d\zbar_i
\times
 \epsilon(d_X)\,\epsilon(d_Y)
\]
Hence, for a test form $\varphi$ on $Y$,
we have
\begin{multline}
\label{eq;13.8.3.1}
 \Bigl\langle
 f_{\dagger}^{(0)}C
 \bigl(
 (\eta_X/\eta_Y)\,m_1,\,
 \overline{(\eta_X/\eta_Y)\,m_2}
 \bigr),\varphi
 \Bigr\rangle
= \\
\Bigl\langle
 C(m_1,\overline{m}_2),\,
  f^{\ast}\varphi\big/
 \prod_{i=\ell+1}^{n} dz_{i}d\zbar_i
 \Bigr\rangle
 \left(\frac{1}{2\pi\sqrt{-1}}\right)^{d_X-d_Y}
\end{multline}
When we consider $\nbigr$-triples,
$\eta_X/\eta_Y$ is replaced with
$(\lambda^{d_X}\eta_X)\big/(\lambda^{d_Y}\eta_Y)$,
Then, the signature $(-1)^{d_X-d_Y}$ appears
in the right hand side of 
{\rm(\ref{eq;13.8.3.1})}.
The formula is the same as
that in {\rm\S1.6.d \cite{sabbah2}}.
\hfill\qed
\end{rem}

\subsubsection{Projection}
Let us consider the case
$f:X=Z\times Y\lrarr Y$ is the projection.
Set $d_Z=\dim Z$.
We have
$f_{\dagger}M_i=
 f_{!}\bigl(M_i\otimes\nbige_Z^{\bullet}[d_Z]\bigr)$.
We set
\begin{multline}
 f_{\dagger}^{(0)}C\bigl(
 \eta^{d_Z-p}\,m_1,\,
 \overline{\eta^{d_Z+p}m_2}
 \bigr)
:= \\
 \int\eta^{d_Z-p}\wedge
 \etabar^{d_Z+p}
 C(m_1,\overline{m}_2)
\left(
 \frac{1}{2\pi\sqrt{-1}}
 \right)^{d_Z} 
 \epsilon(d_Z)\,(-1)^{pd_Z}
\end{multline}
We have the following:
\begin{equation}
\label{eq;13.8.3.2}
 f_{\dagger}C\bigl(
 \eta^{d_Z-p}\,m_1,\,
 \overline{\eta^{d_Z+p}m_2}
 \bigr)
=
 \int\eta^{d_Z-p}\wedge
 \etabar^{d_Z+p}
 C(m_1,\overline{m}_2)
\left(
 \frac{1}{2\pi\sqrt{-1}}
 \right)^{d_Z} 
 \epsilon(p+d_Z)
\end{equation}
Note that $\epsilon(d_Z)(-1)^{pd_Z}\epsilon(p)=\epsilon(p+d_Z)$.
The formula (\ref{eq;13.8.3.2})
is the same as that in \S1.6.d \cite{sabbah2}.

\subsubsection{Some compatibility}
For our later purpose,
we give some lemmas
on compatibility of the push-forward
for closed immersion and projection
(Lemma \ref{lem;11.3.21.4} and
Lemma \ref{lem;11.3.21.5}).
Let $\ell\leq k$.
We put
$Y:=\Delta^n$.
Let $X$ and $Z$ be the submanifold of $Y$
given as
$X:=\{z_{\ell+1}=\cdots =z_n=0\}$
and
$Z:=\{z_{k+1}=\cdots =z_n=0\}$.
Let $f:X\lrarr Y$ be the natural inclusion.
Let $g:Y\lrarr Z$ be the projection
forgetting the last components.
The composite $g\circ f:X\lrarr Z$ is
the natural inclusion.

For $p=1,\ldots,n$,
we set $\eta_p:=dz_1\cdots dz_p$.
For $p_1\leq p_2$,
we set $\eta_{p_1,p_2}:=\eta_{p_2}/\eta_{p_1}$.
Let 
$(M_1,M_2,C)\in\Dtriplecat(X)$.
We put $d_{Y/Z}:=d_Y-d_Z$.
We have
\[
 g_{\dagger}\bigl(
 f_{\dagger}M_i
 \bigr)
\simeq
 g_{\ast}\Bigl(
 \Omega_{Y/Z}^{\bullet}[d_{Y/Z}]
 \otimes f_{\dagger}M_i
 \Bigr)
\]
We have the quasi-isomorphism
$g_{\dagger}\bigl(
 f_{\dagger}M_i\bigr)
\simeq
 (g\circ f)_{\dagger}M_i$
induced by 
\[
\eta_{k,n}\otimes \bigl(
 m_i/\eta_{\ell,n}
 \bigr)
\longmapsto
 m_i/\eta_{\ell,k}. 
\]

\begin{lem}
\label{lem;11.3.21.4}
We have
$g^{(0)}_{\dagger}(f^{(0)}_{\dagger}C)
=(g\circ f)^{(0)}_{\dagger}C$
under the quasi-isomorphism.
\end{lem}
\pf
We have the following:
{\small
\begin{multline}
\Bigl\langle
 g^{(0)}_{\dagger}\bigl(
 f^{(0)}_{\dagger}C\bigr)
 \Bigl(
 \eta_{k,n}\,(m_1/\eta_{\ell,n}),\,\,
 \overline{\eta_{k,n}\,(m_2/\eta_{\ell,n})}
 \Bigr),\,\varphi
\Bigr\rangle
 \\
=\left(
 \frac{1}{2\pi\sqrt{-1}}
 \right)^{n-k}\,
 \epsilon(n-k)\,
 \int \eta_{k,n}\etabar_{k,n}
 f^{(0)}_{\dagger}C\Bigl(
 m_1/\eta_{\ell,n},\,
 \overline{m_2/\eta_{\ell,n}}
 \Bigr)
 g^{\ast}\varphi
 \\
=\left(
 \frac{1}{2\pi\sqrt{-1}}
 \right)^{n-k}
 \epsilon(n-k)\,
 \Bigl\langle
 f^{(0)}_{\dagger}C\bigl(
 m_1/\eta_{\ell,n},\overline{m_2/\eta_{\ell,n}}
 \bigr),\,\,
 \eta_{k,n}
 \etabar_{k,n}
 g^{\ast}\varphi
 \Bigr\rangle
 \\
=\left(
 \frac{1}{2\pi\sqrt{-1}}
 \right)^{\ell-k}
\epsilon(n-k)\,\epsilon(\ell)\,\epsilon(n)
 \Bigl\langle
 (\eta_n\etabar_n)^{-1}
 f_{\ast}\bigl(\eta_{\ell}\etabar_{\ell}
 C(m_1,\overline{m}_2)\bigr),\,\,
 g^{\ast}\varphi\,\eta_{k,n}\etabar_{k,n}
 \Bigr\rangle
 \\
=\left(
 \frac{1}{2\pi\sqrt{-1}}
 \right)^{\ell-k}
\epsilon(n-k)\,\epsilon(\ell)\,\epsilon(n)
 \Bigl\langle
 \eta_{\ell}\etabar_{\ell}
 C(m_1,\overline{m}_2),\,\,
 f^{\ast}\Bigl(
 \frac{g^{\ast}\varphi\,\eta_{k,n}\etabar_{k,n}} {\eta_n\etabar_n}
 \Bigr)
 \Bigr\rangle
\end{multline}
}
We also have the following:
\begin{multline}
 \Bigl\langle
 (g\circ f)^{(0)}_{\dagger}C\bigl(
  m_1/\eta_{\ell,k},\,
 \overline{m_2/\eta_{\ell,k}}
 \bigr),\,\,\varphi
 \Bigr\rangle
 \\
=\left(
 \frac{1}{2\pi\sqrt{-1}}
 \right)^{\ell-k}
 \Bigl\langle
 (\eta_k\etabar_k)^{-1}
 (g\circ f)_{\ast}\bigl(
 \eta_{\ell}\etabar_{\ell}
 C(m_1,\overline{m}_2)
 \bigr),\,
 \varphi
 \Bigr\rangle
 \times
 \epsilon(\ell)\,\epsilon(k)
\end{multline}
Let $\varphi=\eta_k\etabar_k\,\phi$.
Then, we have
\[
 (g\circ f)^{\ast}\left(
 \frac{\varphi}{\eta_k\etabar_k}
 \right)
=(g\circ f)^{\ast}\phi
=(-1)^{k(n-k)}
 f^{\ast}\left(
 \frac{g^{\ast}(\eta_k\etabar_k\phi)\,\eta_{k,n}\etabar_{k,n}}
 {\eta_n\etabar_n}
 \right)
\]
Because
$(n-k)(n-k-1)/2-n(n-1)/2-k(k-1)/2+k(n-k)\equiv 0$
modulo $2$,
we obtain 
$(g\circ f)_{\dagger}^{(0)}C
=g_{\dagger}^{(0)}f_{\dagger}^{(0)}C$.
Thus, Lemma \ref{lem;11.3.21.4} is proved.
\hfill\qed

\vspace{.1in}

We consider another compatibility.
Let $Z$ be a complex manifold.
Let us consider the following diagram:
\[
 \begin{CD}
 Z\times \Delta^{k} @>{g_1}>> 
 Z\times\Delta^{k+\ell}\\
 @V{f_1}VV @V{f_2}VV \\
 \Delta^k @>{g_2}>> \Delta^{k+\ell}
 \end{CD}
\]
Here, $g_i$ are the natural inclusions,
and $f_i$ are the projections.
Let $(M_1,M_2,C)\in\Dtriplecat(Z\times\Delta^k)$
such that the support of $M_i$ are proper
with respect to $f_1$.
We have natural isomorphisms
$g_{2\dagger}\bigl(
 f_{1\dagger}
 M_i
 \bigr)
\simeq
 f_{2\dagger}g_{1\dagger}M_i$.

\begin{lem}
\label{lem;11.3.21.5}
We have
$g_{2\dagger}^{(0)}\circ f^{(0)}_{1\dagger}(C)
=f_{2\dagger}^{(0)}\circ g^{(0)}_{1\dagger}(C)$
under the isomorphisms.
\end{lem}
\pf
We put $n=\dim Z$.
We have the following:
{\small
\begin{multline}
 g_{2\dagger}^{(0)}
\bigl(
 f^{(0)}_{1\dagger}(C)
 \bigr)
 \Bigl(
 (\eta_k/\eta_{k+\ell})\cdot
 \xi^{n-p}m_1,\,
 \overline{(\eta_k/\eta_{k+\ell})\,\xi^{n+p}m_2}
 \Bigr)
 \\
=
 \left(
 \frac{1}{2\pi\sqrt{-1}}
 \right)^{-\ell}
 (\eta_{k+\ell}\etabar_{k+\ell})^{-1}
 g_{2\ast}\Bigl(
 \eta_k\etabar_k\,
 f^{(0)}_{1\dagger}(C)(\xi^{n-p}m_1,\overline{\xi^{n+p}m_2})
 \Bigr)
 \times 
 \epsilon(k+\ell)\,\epsilon(k)
 \\
=
 \left(
 \frac{1}{2\pi\sqrt{-1}}
 \right)^{n-\ell}
 (\eta_{k+\ell}\etabar_{k+\ell})^{-1}
 g_{2\ast}\Bigl(
 \eta_k\etabar_k\,
 \int_{f_1}
 \xi^{n-p}\xibar^{n+p}
 C(m_1,\overline{m}_2)
 \Bigr) 
 \epsilon(n)\,(-1)^{np}\,
 \epsilon(k+\ell)\,\epsilon(k)
\end{multline}
}
Suppose that 
the supports
$\Supp(\xi^{n-p})$ and 
$\Supp(\xi^{n+p})$ are sufficiently small
so that 
there exists a $C^{\infty}$-local generator
$\eta_Z$ of $\omega_Z$
on a neighbourhood $\nbigu$ of
$\Supp(\xi^{n-p})\cup\Supp(\xi^{n+p})$.
We set $\eta'_{m}:=\eta_Z\cdot\eta_m$ on $\nbigu$.
We also have the following:
{\small
\begin{multline}
f^{(0)}_{2\dagger}\bigl(g^{(0)}_{1\dagger}C\bigr)
\Bigl(
\xi^{n-p}(\eta'_{k}/\eta'_{k+\ell})\,m_1,\,\,
\overline{\xi^{n+p}(\eta'_{k}/\eta'_{k+\ell})m_2}
\Bigr)
\\
=\left(\frac{1}{2\pi\sqrt{-1}}\right)^{n}
 \int_{f_2}
 \xi^{n-p}\xibar^{n+p}
 g_{1\dagger}\bigl(
 (\eta'_{k}/\eta'_{k+\ell})\,m_1,\,
 \overline{(\eta'_{k}/\eta'_{k+\ell})\,m_2}
 \bigr)
\cdot
 \epsilon(n)\,(-1)^{-np}
 \\
=\!\!\left(
 \frac{1}{2\pi\sqrt{-1}}
 \right)^{n-\ell}\!\!\!\!\!\!
 \int_{f_2}\!\!\xi^{n-p}\xibar^{n+p}\!\!
 (\!\eta'_{k+\ell}\etabar'_{k+\ell})^{-1}
 g_{1\ast}\!\bigl(
 \eta'_{k}\etabar'_{k}
 C(m_1,\overline{m}_2)\!
 \bigr)
\epsilon(n)(-1)^{-np}
\epsilon(n+k+\ell)\epsilon(n+k)
 \\
 \mbox{{}}
\end{multline}
}
We have the following equality modulo $2$:
{\small
\begin{multline}
 (k+\ell)(k+\ell-1)/2
-k(k-1)/2-(n+k+\ell)(n+k+\ell-1)/2
-(n+k)(n+k-1)/2
 \\
\equiv
 k\ell+(n+k)\ell
=n\ell
\end{multline}
}
We have the following:
\begin{multline}
\label{eq;13.4.2.2}
 \Bigl\langle
 (\eta_{k+\ell}\etabar_{k+\ell})^{-1}
 g_{2\ast}\Bigl(
 \eta_k\etabar_k\,\int_{f_1}\xi^{n-p}\xibar^{n+p}
 C(m_1,\overline{m}_2)
 \Bigr),\,
 \phi\eta_{k+\ell}\etabar_{k+\ell}
 \Bigr\rangle \\
=\Bigl\langle
 \int_{f_1}\xi^{n-p}\xibar^{n+p} C(m_1,\overline{m}_2),\,
 g_2^{\ast}(\phi)\eta_k\etabar_k
 \Bigr\rangle
\end{multline}
We have the following:
\begin{multline}
 \label{eq;11.3.18.1}
 \Bigl\langle
 \int_{f_2}\xi^{n-p}\xibar^{n+p}
 \Bigl(
 (\eta'_{k+\ell}
 \etabar'_{k+\ell})^{-1}
 g_{1\ast}\bigl(\eta'_{k}\etabar'_{k}\cdot
 C(m_1,\overline{m}_2)\bigr)
 \Bigr),\,
 \phi\cdot\eta_{k+\ell}\etabar_{k+\ell}
 \Bigr\rangle \\
=\Bigl\langle
 g_{1\ast}\bigl(\eta'_{k}\etabar'_{k}C(m_1,\overline{m}_2)\bigr),\,
 \frac{f_2^{\ast}\phi\,\xi^{n-p}\xibar^{n+p}\eta_{k+\ell}\etabar_{k+\ell}}
 {\eta'_{k+\ell}\etabar'_{k+\ell}}
 \Bigr\rangle
\end{multline}
We have the description
$\xi^{n-p}\xibar^{n+p}=A\eta_Z\etabar_Z$.
Then, we have
\[
 \frac{f_2^{\ast}\phi\,\xi^{n-p}\xibar^{n+p}\eta_{k+\ell}\etabar_{k+\ell}}
 {\eta'_{k+\ell}\etabar'_{k+\ell}}
=f_2^{\ast}\phi\,A\,(-1)^{n(k+\ell)}
\]
Then, (\ref{eq;11.3.18.1}) is rewritten as follows:
\begin{multline}
\label{eq;13.4.2.1}
 \Bigl\langle
 C(m_1,\overline{m}_2),\,
 g_1^{\ast}(f_2^{\ast}\phi^{\ast}A(-1)^{n(k+\ell)})
 \eta'_{k}\etabar'_{k}
 \Bigr\rangle
= \\
 \Bigl\langle
 C(m_1,\overline{m}_2),\,
 g_1^{\ast}(f_2^{\ast}\phi^{\ast}A(-1)^{n\ell})
 \eta_{Z}\etabar_{Z}\eta_{k}\etabar_{k}
 \Bigr\rangle
 \\
=(-1)^{n\ell}
 \Bigl\langle
 C(m_1,\overline{m}_2),\,
 \xi^{n-p}\xibar^{n+p}
 (g_2\circ f_1)^{\ast}\phi\,\eta_k\etabar_k
 \Bigr\rangle
\end{multline}
By comparing (\ref{eq;13.4.2.2}) and (\ref{eq;13.4.2.1}),
we obtain the desired equality 
in the case that
$\Supp(\xi^{n-p})\cup\Supp(\xi^{n+p})$
is sufficiently small.
We obtain the general case by using
the partition of the unity,
and thus Lemma \ref{lem;11.3.21.5} is proved.
\hfill\qed

\subsubsection{Construction of the push-forward in the general case}
In the general case,
we factor $f$ into the closed immersion
$f_1:X\lrarr X\times Y$
and the projection 
$f_2:X\times Y\lrarr Y$.
We obtain 
a $D$-complex-triple
$f^{(0)}_{2\dagger}\bigl(
 f^{(0)}_{1\dagger}(M_1,M_2,C)
 \bigr)$
on $Y$.
We obtain the following lemma
from Lemma \ref{lem;11.3.21.4}
and Lemma \ref{lem;11.3.21.5}.
\begin{lem}
If $f$ is a closed immersion or a projection,
$f^{(0)}_{2\dagger}\bigl(
 f^{(0)}_{1\dagger}(M_1,M_2,C)
 \bigr)$
is naturally isomorphic to
the object $f^{(0)}_{\dagger}(M_1,M_2,C)$
given previously.
\hfill\qed
\end{lem}

\subsubsection{Push-forward of $D$-complex-triples}

Let $(M_1^{\bullet},M_2^{\bullet},C)\in
 \Dcomplextriplecat(X)$ such that
the restriction of $f$ to the support of $M_i$
 are proper.
We have the following object in
$\Ddoublecomplextriplecat(Y)$:
\[
\Bigl(
 f_{\dagger}(M_1^{\bullet})^{\bullet},
 f_{\dagger}(M_2^{\bullet})^{\bullet},
 f^{(0)}_{\dagger}C
\Bigr)
\]
By taking the total complex,
we obtain
an object in $\Dcomplextriplecat(Y)$,
which we denote by
$f^{(0)}_{\dagger}(M_1^{\bullet},M_2^{\bullet},C)$.
Correspondingly,
we have the push-forward for complexes
of $D$-triples.

\subsubsection{Composition}

Let $f:X\lrarr Y$ and $g:Y\lrarr Z$ be morphisms 
of complex manifolds.
Let $(M_i^{\bullet},M_2^{\bullet},C)$
be a $D_X$-complex-triple such that
the restriction of $f$ and $g\circ f$ to 
the supports of $M_i$ are proper.
We have natural quasi-isomorphisms
$(g\circ f)_{\dagger}M_i\simeq
 g_{\dagger}\bigl(f_{\dagger}M_i\bigr)$.
The following lemma is implied in \cite{sabbah2}.
\begin{lem}
\label{lem;11.3.21.10}
We have
$(g\circ f)^{(0)}_{\dagger}
=g^{(0)}_{\dagger}(f^{(0)}_{\dagger}C)$
under the quasi-isomorphisms.
\end{lem}
\pf
By using Lemma \ref{lem;11.3.21.4}
and Lemma \ref{lem;11.3.21.5},
we can reduce the issue to the case that
$f:W\times Z\times Y\lrarr Z\times Y$
and $g:Z\times Y\lrarr Y$ are the projections.
We set $m:=\dim W$
and $n:=\dim Z$.
We have the following:
\begin{multline}
 (g\circ f)^{(0)}_{\dagger}C\bigl(\xi^{n-q}\eta^{m-p}m_1,
 \overline{\xi^{n+q}\eta^{m+p}m_2}
 \bigr) \\
=\left(\frac{1}{2\pi\sqrt{-1}}\right)^{m+n}
 \int\xi^{n-q}\eta^{m-p}
 \overline{\xi^{n+q}\eta^{m+p}}
 C(m_1,\overline{m}_2)\,
 \epsilon(m+n)\,(-1)^{-(m+n)(p+q)}
\end{multline}
We also have the following
(recall (\ref{eq;13.4.2.10})):
{\small
\begin{multline}
 g^{(0)}_{\dagger}\bigl(f^{(0)}_{\dagger}C\bigr)
 \Bigl(
 \xi^{n-q}\eta^{m-p}m_1,\,
 \overline{\xi^{n+q}\eta^{m+p}m_2}
 \Bigr)
 \times(-1)^{pq}
 \\
=\left(
 \frac{1}{2\pi\sqrt{-1}}
 \right)^{n}
 \int\xi^{n-q}\xibar^{n+q}
 f^{(0)}_{\dagger}C(\eta^{m-p}m_1,\,\overline{\eta^{m+p}m_2})
\,\epsilon(n)\,(-1)^{-nq+pq}
 \\
=\left(\frac{1}{2\pi\sqrt{-1}}\right)^{m+n}
 \int\xi^{n-q}\eta^{m-p}
 \xibar^{n+q}\etabar^{m+p}
 C(m_1,\overline{m}_2)
\epsilon(m)\epsilon(n)
(-1)^{-mp-nq+pq+(n+q)(m-p)}
\end{multline}
}
We have the following equality modulo $2$:
\begin{multline}
 m(m-1)/2+n(n-1)/2+mp+nq+pq+(n+q)(m-p)
\equiv \\
 (m+n)(m+n-1)/2+(m+n)(p+q)
\end{multline}
Thus, we obtain Lemma \ref{lem;11.3.21.10}.
\hfill\qed

\subsubsection{Correspondence of left and right triples (Appendix)}

The correspondence between hermitian pairings
for left and right $D$-modules
are given as follows.
Let $C$ be a hermitian pairing 
$M_1\times \Mbar_2\lrarr \distribution_X$.
The corresponding pairing of
right $D$-modules
$(\omega_X\otimes M_1)\otimes_{\cnum}
 \overline{(\omega_X\otimes M_2)}\lrarr
 \distribution_X^{d_X,d_X}$
is given as follows:
\[
 C^r\bigl(
 \eta_1\otimes m_1,\,
 \overline{\eta_2\otimes m_2}
 \bigr)
:=\eta_1\etabar_2\,
 C(m_1,\overline{m}_2)\,
 \epsilon(d_X)\,
 \left(\frac{1}{2\pi\sqrt{-1}}\right)^{d_X}
\]
For a closed immersion $f:X\lrarr Y$,
we should have
$f_{\ast}C^r(\eta_X\otimes m_1,\,\overline{\eta_X\otimes m_2})
=
 (f^{(0)}_{\dagger}C)^r\bigl(
 \eta_X\otimes m_1,\,
 \overline{\eta_X\otimes m_2}
 \bigr)$.
It implies the formula in the closed immersion,
as follows.
\begin{multline}
\Bigl\langle
 C(m_1,\overline{m}_2),\,f^{\ast}\phi\,\eta_X\etabar_X
\Bigr\rangle\,
\epsilon(d_X)\,
\left(
 \frac{1}{2\pi\sqrt{-1}}
 \right)^{d_X}=
 \\
\Bigl\langle
 f_{\ast}C^r(\eta_X\otimes m_1,\overline{\eta_X\otimes m_2}),\,
 \phi
 \Bigr\rangle
=\Bigl\langle
 (f_{\dagger}^{(0)}C)^r(\eta_X\otimes m_1,\overline{\eta_X\otimes m_2}),\,
 \phi
 \Bigr\rangle \\
=\Bigl\langle
 (f_{\dagger}^{(0)}C)\bigl(
 (\eta_X/\eta_Y)m_1,
 \overline{(\eta_X/\eta_Y)m_2}\bigr),\,
 f^{\ast}\phi\,\eta_Y\etabar_Y
\Bigr\rangle\,
\epsilon(d_Y)\,
\left(\frac{1}{2\pi\sqrt{-1}}
 \right)^{d_Y}
\end{multline}

\subsection{Hermitian adjoint of $D$-complex-triples}
\label{subsection;13.8.5.1}
\index{Hermitian adjoint}
\index{functor $\DDD^{\herm}$}
\index{functor $\DDD^{(0)\herm}$}

For $\nbigt^{\bullet}\!\in\!
\nbigc(\Dtriplecat(X))$,
we define 
$\DDD^{\herm}(\nbigt^{\bullet})$
in $\nbigc(\Dtriplecat(X))$
as follows.
The $p$-th member
$\DDD^{\herm}(\nbigt^{\bullet})^p$
is $(\nbigt^{-p})^{\ast}$.
The differential
$\DDD^{\herm}(\nbigt^{\bullet})^p
\lrarr
 \DDD^{\herm}(\nbigt^{\bullet})^{p+1}$
is defined 
by 
$\DDD^{\herm}(\delta_{-p-1})$,
as in \S1.6.c \cite{sabbah2}.
Then, $\DDD^{\herm}$ gives
a contravariant auto equivalence 
of $\nbigc(\Dtriplecat(X))$.
We have a natural isomorphism
$\DDD^{\herm}\bigl(
 \nbigt^{\bullet}[\ell]
\bigr)
=\DDD^{\herm}(\nbigt^{\bullet})[-\ell]$.

By the equivalence $\Psi_1$,
we obtain a contravariant auto equivalence
$\DDD^{(0)\herm}$ 
on $\Dcomplextriplecat(X)$:
\[
 \DDD^{(0)\herm}(M_1^{\bullet},M_2^{\bullet},C)
=(M_2^{\bullet},M_1^{\bullet},C^{\ast})
\]
Here,
$(C^{\ast})^{p}=(-1)^p\bigl(C^p\bigr)^{\ast}$.
We have a natural isomorphism
$\nbigs_{-\ell}\circ\DDD^{\herm}
\simeq
 \DDD^{\herm}\circ\nbigs_{\ell}$.

\begin{lem}
\label{lem;11.3.23.20}
Let $f:X\lrarr Y$
be a morphism of complex manifolds.
Let $\gbigf\in\Dcomplextriplecat(X)$
such that the restriction of $f$ to the support of $\gbigt$
is proper.
We have 
$f^{(0)}_{\dagger}\circ \DDD^{(0)\herm}(\gbigt)
=
 \DDD^{(0)\herm}\circ f^{(0)}_{\dagger}(\gbigt)$
in $\Dcomplextriplecat(Y)$.
As a consequence, 
for $\nbigt^{\bullet}\in\nbigc(\Dtriplecat(X))$,
we have 
$f_{\dagger}\circ \DDD^{\herm}(\nbigt^{\bullet})
=
 \DDD^{\herm}\circ f_{\dagger}(\nbigt^{\bullet})$
in $\nbigc(\Dtriplecat(Y))$.
\end{lem}
\pf
We use the notation in the construction of
push-forward in \S\ref{subsection;11.3.22.1}.
Let us consider the case that $f$ is a closed immersion.
We have the following:
{\small
\begin{multline}
f^{(0)}_{\dagger}(C^{\ast})\bigl(
 (\eta_X/\eta_Y)\,m_2,\,\,
 \overline{(\eta_X/\eta_Y)\,m_1}
 \bigr)
 \\
=(\eta_Y\etabar_Y)^{-1}\,
 f_{\ast}\bigl(
 \eta_X\,\etabar_X\,
 C^{\ast}(m_2,\overline{m}_1)
 \bigr)
 \,\left(\frac{1}{2\pi\sqrt{-1}}\right)^{d_X-d_Y}
 \epsilon(d_X)\,\epsilon(d_Y)
 \\
=(\eta_Y\etabar_Y)^{-1}
 f_{\ast}\bigl(
 \eta_X\etabar_X\,
 \overline{C(m_1,\overline{m}_2)}
 \bigr)
 \left(\frac{1}{2\pi\sqrt{-1}}\right)^{d_X-d_Y}
 \epsilon(d_X)\,\epsilon(d_Y)
 \\
=\overline{
 (\etabar_Y\eta_Y)^{-1}
 f_{\ast}\bigl(\etabar_X\eta_X
 C(m_1,\overline{m}_2)
 \bigr)
 \left(\frac{1}{2\pi\sqrt{-1}}\right)^{d_X-d_Y}
 }
 (-1)^{d_X-d_Y}
 \epsilon(d_X)\,\epsilon(d_Y)
\\
=\overline{
 (\eta_Y\etabar_Y)^{-1}
 f_{\ast}\bigl(\eta_X\etabar_X
 C(m_1,\overline{m}_2)
 \bigr)
 \left(\frac{1}{2\pi\sqrt{-1}}\right)^{d_X-d_Y}
 }
 \epsilon(d_X)\,\epsilon(d_Y)
\\
=(f^{(0)}_{\dagger}C)^{\ast}\bigl(
 (\eta_X/\eta_Y)\,m_1,
 \overline{(\eta_X/\eta_Y)\,m_2}
 \bigr)
\end{multline}
}
Let us consider the case
that $f:Z\times Y\lrarr Y$.
We put $m:=\dim Z$.
We have the following:
{\small
\begin{multline}
 (f^{(0)}_{\dagger}C^{\ast})(\eta^{m+p}m_2,\,
 \overline{\eta^{m-p}m_1})
=\left(\frac{1}{2\pi\sqrt{-1}}\right)^m
 \int\eta^{m+p}\etabar^{m-p}
 C^{\ast}(m_2,\overline{m}_1)\,
 \epsilon(m)\,(-1)^{mp}
 \\
=\overline{
 \left(\frac{1}{2\pi\sqrt{-1}}\right)^m(-1)^m
 \int\eta^{m-p}\etabar^{m+p}
 (-1)^{m-p}
 C(m_1,\overline{m}_2)\,
 \epsilon(m)\,(-1)^{mp}
 }
\\
=(-1)^p(f^{(0)}_{\dagger}C)^{\ast}(\eta^{m+p}m_2,
 \overline{\eta^{m-p}m_1}).
\end{multline}
}
Thus, we are done.
\hfill\qed

\subsection{Comparison with the naive push-forward}
\label{subsection;13.4.5.10}

\subsubsection{Trace morphism}
\label{subsection;13.8.17.1}

Let us observe that we have a natural morphism,
called the trace morphism:
\begin{equation}
\label{eq;11.3.23.10}
\tr:
 f_{!}\Bigl(
 \bigl(
 D_{Y\larr X}\otimes_{\cnum}
 D_{\Ybar\larr \Xbar}
 \bigr)
 \otimes^L_{D_{X,\Xbar}}\distribution_X
 \Bigr)
\lrarr
 \distribution_Y
\end{equation}
Indeed, we have the following:
\begin{multline}
\label{eq;11.3.26.1}
\Bigl(
\omega_X\otimes f^{-1}\bigl(
 D_Y\otimes\omega_Y^{-1}
 \bigr)
\otimes_{\cnum}
 \omega_{\Xbar}\otimes f^{-1}\bigl(
 D_{\Ybar}\otimes\omega_{\Ybar}^{-1}
 \bigr)
\Bigr)
\otimes^L_{D_{X,\Xbar}}\distribution_X
\lrarr \\
 f^{-1}(D_{Y,\Ybar}\otimes\omega^{-1}_{Y,\Ybar})
\otimes
 \Omega_X^{\bullet}[d_X]\otimes
 \Omega_{\Xbar}^{\bullet}[d_X]
\otimes_{\nbigo_{X,\Xbar}}
 \distribution_X
 \\
\lrarr
 f^{-1}(D_{Y,\Ybar}\otimes\omega^{-1}_{Y,\Ybar})
\otimes
 \distribution_X^{\bullet}[2d_X]
\end{multline}
Here, we use the pairing
$\Omega_X^{\bullet}[d_X]\otimes\Omega_{\Xbar}^{\bullet}[d_X]
\lrarr
 \distribution_X^{\bullet}[2d_X]$
given as follows.
We have the natural identification
$\Omega_X^{\bullet}[d_X]\otimes
 \Omega_{\Xbar}^{\bullet}[d_X]
\simeq 
\bigl(
 \Omega_X^{\bullet}\otimes\Omega_{\Xbar}^{\bullet}
\bigr)[2d_X]$
which is given by the multiplication of
$(-1)^{pd_X}\epsilon(d_X)$
on $\Omega^{p+d_X}_X\otimes\Omega_{\Xbar}^{q+d_X}$.
We have the natural isomorphism
$\Tot\bigl(
 \Omega_X^{\bullet}\otimes\Omega_{\Xbar}^{\bullet}
 \bigr)\otimes\distribution_X
\lrarr
 \distribution_X^{\bullet}$.
The composition gives the desired pairing.
Then, we obtain
\begin{multline}
 f_{!}\Bigl(
 \bigl(
 D_{Y\larr X}\otimes_{\cnum} D_{\Ybar\larr \Xbar}
\bigr)\otimes_{D_{X,\Xbar}}^L
 \distribution_X
 \Bigr)
\lrarr
 \bigl(
 D_{Y,\Ybar}\otimes \omega_{Y,\Ybar}^{-1}
 \bigr)
 \otimes f_!(\distribution_X^{\bullet}[2d_X])
 \\
\stackrel{A}{\lrarr}
\bigl(
 D_{Y,\Ybar}\otimes\omega_{Y,\Ybar}^{-1}
\bigr)
\otimes\distribution_Y^{\bullet}[2d_Y]
 \\
\stackrel{B}{\simeq}
 \bigl(
 \omega_Y^{-1}
\otimes
 D_Y\otimes\Omega_Y^{\bullet}[d_Y]
\bigr)
 \otimes
 \bigl(
 \omega_{\Ybar}^{-1}\otimes
 D_{\Ybar}\otimes\Omega_{\Ybar}^{\bullet}[d_Y]
\bigr)
\otimes\distribution_Y
\simeq
 \distribution_Y
\end{multline}
The morphism $A$
is induced by the integration
$f_!\distribution_X^{p+d_X,q+d_X}\lrarr
 \distribution_Y^{p+d_Y,q+d_Y}$
multiplied with $(2\pi\sqrt{-1})^{d_X-d_Y}$.
For $B$, we use the identification 
as in the case of (\ref{eq;11.3.26.1}).

\subsubsection{Naive push-forward of pairings}
Let $C:M_1\times \Mbar_2\lrarr\distribution_X$.
We have the following induced pairing:
\[
\begin{CD}
\bigl(
 D_{Y\larr X}\otimes^{L}_{D_X}M_1
\bigr)
\otimes_{\cnum}
\overline{
\bigl(
 D_{Y\larr X}\otimes^{L}_{D_X}M_2
\bigr)}
 @>>>
 \Bigl(
 D_{Y\larr X}\otimes_{\cnum}
 D_{\Ybar\larr \Xbar}
 \Bigr)
 \otimes_{D_{X,\Xbar}}^L\distribution_X
\end{CD}
\]
Hence, we obtain the following pairing 
\[
\begin{CD}
 f_{\dagger}M_1\otimes
 f_{\dagger}\Mbar_2
 @>>>
 f_!\Bigl(
 \bigl(
 D_{Y\larr X}\otimes D_{\Ybar\larr \Xbar}
 \Bigr)\otimes_{D_{X,\Xbar}}^L
\distribution_X^{\bullet}[2d_X]\bigr)
 @>{\tr}>>
 \distribution_Y
\end{CD}
\]
It is denoted by $f_{\dagger}^{(1)}C$.
We have the following comparison.
\begin{lem}
We have
$f_{\dagger}^{(1)}(C)
=f_{\dagger}^{(0)}(C)$.
\end{lem}
\pf
Let us consider the case
that $f$ is a closed immersion.
We have the following:
{\small
\begin{multline}
 f_{\dagger}^{(1)}(C)\bigl(
 (\eta_X/\eta_Y)\,m_1,\,\,
 \overline{(\eta_X/\eta_Y)\,m_2}
 \bigr) \\
=(\eta_Y\,\etabar_Y)^{-1}
 f_{\ast}\bigl(
 \eta_X\etabar_X\,
 C(m_1,\overline{m}_2)
 \bigr)
 \left(\frac{1}{2\pi\sqrt{-1}}\right)^{d_Y-d_X}
 \epsilon(d_X)\,\epsilon(d_Y)
\end{multline}
}
Here,
$\epsilon(d_X)$
and $\epsilon(d_Y)$
appear by the identifications
$\Omega_X^{\bullet}[d_X]
 \otimes
 \Omega_{\Xbar}^{\bullet}[d_X]
\simeq
 \bigl(
 \Omega_{X}^{\bullet}\otimes
 \Omega_{\Xbar}^{\bullet}
 \bigr)[2d_X]$
and
$\Omega_Y^{\bullet}[d_Y]
 \otimes
 \Omega_{\Ybar}^{\bullet}[d_Y]
\simeq
 \bigl(
 \Omega_{Y}^{\bullet}\otimes
 \Omega_{\Ybar}^{\bullet}
 \bigr)[2d_Y]$.
By definition,
we have the following:
{\small
\begin{multline}
 f^{(0)}_{\dagger}(C)
 \bigl(
 (\eta_X/\eta_Y)\,m_1,\,
 \overline{(\eta_X/\eta_Y)\,m_2}
 \bigr)  \\
=(\eta_Y\etabar_Y)^{-1}
 f_{\ast}\bigl(
 \eta_X\etabar_X\,
 C(m_1,\overline{m}_2)
 \bigr)
\cdot
\left(\frac{1}{2\pi\sqrt{-1}}\right)^{d_X-d_Y}
 \epsilon(d_X)\,\epsilon(d_Y)
\end{multline}
}
Hence, we are done in this case.

Let us consider the case $f:X=Z\times Y\lrarr Y$.
We have the identification
\[
f_{\dagger} M\simeq
 f_{\dagger}(M\otimes\nbige_Z^{\bullet}[d_Z])
\simeq
 f_{\dagger}\Bigl(M\otimes\nbige_Z^{\bullet}[d_Z]
 \otimes_{f^{-1}\nbigo_Y}
 f^{-1}\bigl(
 \Omega_Y^{\bullet}[d_Y]\otimes D_Y\otimes\omega_Y^{-1}
 \bigr)\Bigr)
\]
Let $m=\dim Z$ and $n=\dim Y$.
We have
{\small
\begin{multline}
 \label{eq;11.3.26.10}
f_{\dagger}^{(1)}C\Bigl(
 \eta_Y^{-1}\otimes 
 (m_1\cdot \xi^{m-p}\,\eta_Y),\,\,
 \overline{
 \eta_Y^{-1}\otimes
 (m_2\cdot \xi^{m+p}\eta_Y)}
\Bigr)\\
=(\eta_Y\etabar_Y)^{-1}
 \int C(m_1,\overline{m}_2)\,
 \xi^{m-p}\eta_Y\,
 \overline{\xi^{m+p}\eta_Y}\,
(-1)^{pd_X}\left(
 \frac{1}{2\pi\sqrt{-1}}
 \right)^{m}
 \epsilon(d_X)\,\epsilon(d_Y)
 \\
=(\eta_Y\etabar_Y)^{-1}
 \int C(m_1,\overline{m}_2)\,
 \xi^{m-p}\xibar^{m+p}
 \eta_Y\etabar_Y(-1)^{n(m+p)+pd_X}
 \left(\frac{1}{2\pi\sqrt{-1}}\right)^m
 \epsilon(d_X)\,\epsilon(d_Y)
\end{multline}
}
We also have the following:
{\small
\begin{equation}
\label{eq;11.3.26.11}
 f_{\dagger}^{(0)}C(m_1\,\xi^{m-p},\,\overline{m_2\,\xi^{m+p}})
=\int C(m_1,\overline{m}_2)\,
 \xi^{m-p}\xibar^{m+p}
 (-1)^{mp}\left(\frac{1}{2\pi\sqrt{-1}}\right)^{m}
 \epsilon(m)
\end{equation}
}
It is easy to check that 
(\ref{eq;11.3.26.10}) and (\ref{eq;11.3.26.11})
are equal.
Hence, we obtain
$f_{\dagger}^{(1)}C=f_{\dagger}^{(0)C}$.
\hfill\qed

\begin{rem}
We will not distinguish $f^{(0)}$ and $f^{(1)}$
in the following.
\hfill\qed
\end{rem}

\subsection{Rules for signature (Appendix)}
\label{subsection;11.4.6.20}

\subsubsection{Contravariant functor and a shift of degree}
\label{subsection;13.8.2.1}

Let $A^{\bullet}$ be a $\cnum_X$-complex
with a differential $\del$.
Recall that,
for any integer $\ell$,
the shift $A^{\bullet}[\ell]$
is the complex such that 
$(A^{\bullet}[\ell])^p=A^{\ell+p}$
with the differential $(-1)^{\ell}\del$.
Recall the rule in \cite{deligne-SGA4}
on the signature for contravariant functor.
Let $F$ be a contravariant functor from
the category of $\cnum_X$-modules
to a $\cnum$-linear abelian category $\vecC$.
Then, we obtain a complex $\nbigb$ in $\vecC$
given by $\nbigb^k:=F(A^{-k})$.
The differential
$\nbigb^k\lrarr\nbigb^{k+1}$
is given by $(-1)^{k+1}F(\del):F(A^{-k})\lrarr F(A^{-k-1})$.
The complex is denoted by $F(A^{\bullet})$.
For any integer $\ell$,
we have a natural isomorphism
$F(A^{\bullet})[-\ell]
\simeq
 F(A^{\bullet}[\ell])$
given by the multiplication of
$\epsilon(\ell)(-1)^{p\ell}$
on $F(A^{\bullet})[-\ell]^p=F(A^{-p+\ell})$,
where $\epsilon(\ell)=(-1)^{\ell(\ell-1)/2}$.
The composition
$F(A^{\bullet})[-\ell-k]
\simeq
 F(A^{\bullet}[\ell])[-k]
\simeq
 F(A^{\bullet}[\ell+k])$
is equal to the direct one
$F(A^{\bullet})[-\ell-k]
\simeq
 F(A^{\bullet}[\ell+k])$.

\subsubsection{Naive $n$-complex
 and the total complex}
\label{subsection;13.8.2.2}

Let $\nbiga$ be a naive $n$-complex of $\cnum_X$-modules,
i.e., it is $\seisuu^n$-graded $\cnum_X$-module
$\nbiga=\bigoplus_{\veck\in\seisuu^n}\nbiga^{\veck}$
equipped with
differentials $\del_i:\nbiga^{\veck}\lrarr\nbiga^{\veck+\vecdelta_i}$
$(i=1,\ldots,n)$
such that $\del_i\circ\del_i=0$ and $[\del_i,\del_j]=0$.
Here, the $j$-th component of $\vecdelta_i$
is $1$ if $i=j$, and $0$ if $i\neq j$.
We put $|\veck|=\sum k_i$
for $\veck=(k_i)\in\seisuu^n$.
Then, the total complex $\Tot(\nbiga)$
is defined as a $\seisuu$-graded $\cnum_X$-module
$\Tot(\nbiga)^{k}=\bigoplus_{|\veck|=k}
 \nbiga^{\veck}$
with the differential
$\del:\Tot(\nbiga)^k\lrarr\Tot(\nbiga)^{k+1}$
given by
$\del(a^{\veck})
=\sum_{i}(-1)^{\sum_{j<i}k_i}\del_i(a^{\veck})$.

Let $\nbiga$ be a $n$-complex.
Let $\nbiga[\ell\vecdelta_i]$ be the $n$-complex
given by
$\nbiga[\ell\vecdelta_i]^{\veck}
=\nbiga^{\veck+\ell\vecdelta_i}$
with the differentials $\del_j$ $(j\neq i)$
and $(-1)^{\ell}\del_i$.
We have an isomorphism
$\Tot(\nbiga[\ell\vecdelta_i])
\simeq
 \Tot(\nbiga)[\ell]$
given by
the multiplication of
$(-1)^{\sum_{j<i}\ell k_i}$
on $\nbiga^{\veck}$.

\subsubsection{Inner homomorphism}

Let $A^{\bullet}_i$ $(i=1,2)$ be $\cnum_X$-complexes.
We have the $\cnum_X$-double complex 
$\nbigb^{\bullet,\bullet}$
given by
$\nbigb^{p,q}:=\nhom(A_1^{-q},A_2^p)$.
We often denote it just by
$\nhom(A_1^{\bullet},A_2^{\bullet})$.
By the standard rule of the signature,
we have
$\del (f)=\del\circ f-(-1)^{p-q}f\circ\del$
for $f\in \nhom(A_1^q,A_2^p)$,
i.e.,
$\del(f)(a)=\del \bigl(f(a)\bigr)-(-1)^{|f|}f(\del a)$.

We have the isomorphism
$\nhom(A_1^{\bullet}[\ell],A_2^{\bullet})
\!\simeq\!
 \nhom(A_1^{\bullet},A_2^{\bullet})[-\ell\vecdelta_2]$
in \S\ref{subsection;13.8.2.1}.
The isomorphism in \ref{subsection;13.8.2.2}
induces the following isomorphism
\[
 \Tot\nhom(A_1^{\bullet}[\ell],A_2^{\bullet})
\simeq
 \Tot\nhom(A_1^{\bullet},A_2^{\bullet})[-\ell].
\]
which is given by the multiplication of
$\epsilon(\ell)(-1)^{r\ell}$
on $\Tot\nhom(A_1^{\bullet}[\ell],A_2^{\bullet})^r$.

\subsubsection{Pairing}
\label{subsection;13.8.3.10}

Let $I^{\bullet}$ be a $\cnum_X$-complex.
Let 
$C:\Tot (A_1^{\bullet}\otimes A_2^{\bullet})\lrarr
 I^{\bullet}$
be a morphism.
For any integer $\ell$,
we have an isomorphism
$\Tot(A_1^{\bullet}\otimes A_2^{\bullet})
\simeq
 \Tot(A_1^{\bullet}[-\ell]\otimes A_2^{\bullet}[\ell])$
given by the multiplication of 
$\epsilon(\ell)(-1)^{\ell p}$
on 
$\bigl(
 A_1^{\bullet}[-\ell]
\otimes
 A_2^{\bullet}[\ell]
 \bigr)^{p,q}$.
It induces a morphism
\[
 C[\ell]:
 \Tot\bigl(
 A_1^{\bullet}[-\ell]
\otimes
 A_2^{\bullet}[\ell]
 \bigr)
\lrarr I^{\bullet}.
\]

A morphism 
$C:\Tot (A_1^{\bullet}\otimes A_2^{\bullet})\lrarr
 I^{\bullet}$
is equivalent to 
a morphism
$\Psi_C:A_1^{\bullet}\lrarr \Tot\nhom(A_2^{\bullet},I^{\bullet})$.
The correspondence is given by $\Psi_C(a)(b)=C(a,b)$.
The following lemma can be checked 
in an elementary way.
\begin{lem}
\label{lem;13.4.5.30}
Under the identification
\[
\Tot\nhom(A_2^{\bullet}[\ell],I^{\bullet})
\simeq
 \Tot\nhom(A_2^{\bullet},I^{\bullet})[-\ell],
\]
we have
$\Psi_{C[\ell]}=\Psi_C$.
\hfill\qed
\end{lem}

\section{Localization, twist and specialization
of non-degenerate $D$-triple}
\label{section;13.8.16.1}

\subsection{Category of non-degenerate $D$-triples}

For any $D_X$-module $M$, 
we have an object
$\nbigc_XM:=\nrhom_{D_X}(M,\distribution_X)$
in the derived category of $D_{\Xbar}$-complexes.
\index{functor $\nbigc_X$}
If $M$ is holonomic,
we have the natural identification
\[
 \nbigc_XM\simeq\nbigh^0\nbigc_XM
=\nhom_{\nbigd_X}(M,\distribution_X),
\]
and it is holonomic.
(See \cite{Kashiwara-distribution},
\cite{sabbah4} and \cite{sabbah_lecture_Stokes}.
See also \cite{mochi10}.)

A $D_X$-triple $(M_1,M_2,C)$ is called holonomic,
if $M_i$ are holonomic.
It is called non-degenerate,
if moreover the induced morphism
$\varphi_{C,1,2}:\Mbar_2\lrarr \nbigc_X(M_1)$
is an isomorphism.
The condition is equivalent to
that the induced morphism
$\varphi_{C,2,1}:M_1\lrarr \nbigc_{\Xbar}(\Mbar_2)$
is an isomorphism.
Let $\Dtriplecat^{nd}(X)\subset\Dtriplecat(X)$ denote 
the full subcategory of 
non-degenerate $D_X$-triples $(M_1,M_2,C)$.
\index{category $\Dtriplecat^{nd}(X)$}
It is easy to check that
$\Dtriplecat^{nd}(X)$ is an abelian
subcategory of $\Dtriplecat(X)$.

\subsection{Localization}
\label{subsection;10.12.27.22}

Let $g$ be a holomorphic function 
on a complex manifold.
Let $\nbigt=(M_1,M_2,C)\in\Dtriplecat(X)$ be holonomic.
As in the case of $\nbigr$-triples,
we have the uniquely induced pairings
$C[!g]:M_1[\ast g]\times\overline{M_2[!g]}
\lrarr\distribution_X$
and 
$C[\ast g]:M_1[!g]\times\overline{M_2[\ast g]}
\lrarr\distribution_X$.
(See \S\ref{subsection;10.12.27.20}.)
Thus, we have
$\nbigt[!g]=\bigl(M_1[\ast g],M_2[!g],C[!g]\bigr)$
and
$\nbigt[\ast g]=
 \bigl(M_1[! g],M_2[\ast g],C[\ast g]\bigr)$.
\index{pairing $C[\star g]$}
\index{$D$-triple $\nbigt[\star g]$}

\begin{lem}
If $\nbigt$ is non-degenerate,
$\nbigt[!g]$ and $\nbigt[\ast g]$
are also non-degenerate.
\end{lem}
\pf
Let $N$ be a holonomic $D_{X}$-module
such that
$N(\ast g)\simeq
 \nbigc_{\Xbar}(\Mbar_2[\ast g])(\ast g)$.
We have an isomorphism
$\nbigc_{X}(N)(\ast g)\simeq
 \Mbar_2(\ast g)$.
Hence, we have a unique morphism
$\Mbar_2[!g]\lrarr \nbigc_X(N)$.
It induces
$N\lrarr\nbigc_{\Xbar}(\Mbar_2[!g])$.
By this universal property,
we obtain that
$\nbigc_{\Xbar}(\Mbar_2[!g])
\simeq
 \nbigc_{\Xbar}(\Mbar_2[!g])[\ast g]$.
Therefore, we obtain that
the induced morphism
$M_1[\ast g]\lrarr
 \nbigc_{\Xbar}(\Mbar_2[!g])$
is an isomorphism.
Namely, $\nbigt[!g]$ is non-degenerate.
Similarly, we obtain that
$\nbigt[\ast g]$ is non-degenerate.
\hfill\qed

\begin{lem}
\label{lem;10.12.27.21}
If $g_1^{-1}(0)=g_2^{-1}(0)$,
we have 
$\nbigt[\star g_1]=\nbigt[\star g_2]$.
\end{lem}
\pf
We have $M[\ast g_1]=M[\ast g_2]$
for any holonomic $D$-module $M$,
which implies $M[!g_1]=M[!g_2]$.
Hence, we have the coincidence
of the underlying $D_X$-modules of
$\nbigt[\star g_i]$ $(i=1,2)$.
Then, it is easy to deduce the coincidence
of the pairings on $X$
from the coincidence on $X\setminus g_i^{-1}(0)$.
\hfill\qed

\vspace{.1in}
Let $H$ be a hypersurface of $X$.
By Lemma \ref{lem;10.12.27.21},
we obtain hermitian pairings:
\[
C[\ast H]:
 M_1[!H]\times \overline{M_2[\ast H]}
\lrarr
 \distribution_X,
\quad\quad
C[!H]:
 M_1[\ast H]\times \overline{M_2[!H]}
\lrarr
 \distribution_X.
\]
Thus, we obtain the induced $D$-triples
$\nbigt[\ast H]
\!=\!\bigl(
 M_1[!H],M_2[\ast H],C[\ast H]
 \bigr)$
and
$\nbigt[!H]
\!=\!\bigl(
 M_1[\ast H],M_2[!H],C[!H]
 \bigr)$
for holonomic $\nbigt\in\Dtriplecat(X)$.
If $\nbigt$ is non-degenerate,
$\nbigt[\star H]$ are also non-degenerate.
It is easy to see the following.
\begin{lem}
We have
$\DDD^{\herm}(\nbigt[\ast H])
=
 (\DDD^{\herm}\nbigt)[!H]$
and
$\DDD^{\herm}(\nbigt[! H])
=
 (\DDD^{\herm}\nbigt)[\ast H]$.
\hfill\qed
\end{lem}

\subsection{Tensor with smooth triple}
\label{subsection;11.3.21.10}

A $D_X$-module is called smooth,
if it is coherent as an $\nbigo_X$-module,
which implies that it is a locally free $\nbigo_X$-module.
A $D_X$-triple is called smooth,
if the underlying $D_X$-modules are smooth.
The values of the hermitian pairing
is contained in the sheaf of $C^{\infty}$-functions on $X$.
Let $\nbigt_i=(\nbigm_{i,1},\nbigm_{i,2},C_i)
 \in\Dtriplecat(X)$ $(i=1,2)$.
If $\nbigt_2$ is smooth,
we naturally have
an object
$\nbigt_1\otimes\nbigt_2
:=\bigl(\nbigm_{1,1}\otimes\nbigm_{2,1},
 \nbigm_{2,1}\otimes\nbigm_{2,2},
 C_1\otimes C_2\bigr)$
in $\Dtriplecat(X)$.
We clearly have
$\DDD^{\herm}(\nbigt_1\otimes\nbigt_2)
=\DDD^{\herm}(\nbigt_1)\otimes\DDD^{\herm}(\nbigt_2)$.
\begin{lem}
If $\nbigt_i$ are non-degenerate,
$\nbigt_1\otimes\nbigt_2$ is also non-degenerate.
\end{lem}
\pf
We have only to consider the case
that $\nbigt_2=(\nbigo,\nbigo,C_0)$.
Then, the claim is clear.
\hfill\qed

\vspace{.1in}

Let us consider a more general case.
Let $H$ be a hypersurface of $X$.
A $D_X(\ast H)$-module is called smooth,
if it is a locally free $\nbigo_X(\ast H)$-module.
A smooth $D_{X}(\ast H)$-triple
is a tuple of smooth $D_X(\ast H)$-modules $V_i$ 
with a pairing
$C_V:V_1\otimes \overline{V_2}\lrarr 
\nbigc_X^{\infty\moderate H}$.
It is called non-degenerate,
if $(V_1,V_2,C_V)_{|X\setminus H}$ is non-degenerate.

Let $\nbigt_V=(V_1,V_2,C_V)$ be 
a non-degenerate smooth $D_X(\ast H)$-triple.
Let $\nbigt=(M_1,M_2,C)\in\Dtriplecat^{nd}(X)$.
We have the induced pairing
$C\otimes C_V:
 (M_1\otimes V_1)
\otimes
 \overline{(M_2\otimes V_2)}
\lrarr
 \distribution_X^{\moderate H}$,
where $\distribution_X^{\moderate H}$
denotes the sheaf of distributions on $X$
with moderate growth along $H$.
We obtain hermitian pairings,
as in \S\ref{subsection;10.12.27.22}:
\[
(C\otimes C_V)[\ast H]:
 (M_1\otimes V_1)[!H]
\otimes
 \overline{(M_2\otimes V_2)[\ast H]}
\lrarr
 \distribution_X
\]
\[
(C\otimes C_V)[!H]:
 (M_1\otimes V_1)[\ast H]
\otimes
 \overline{(M_2\otimes V_2)[!H]}
\lrarr
 \distribution_X
\]
Thus, we obtain hermitian pairings:
\[
(\nbigt\otimes\nbigt_V)[\ast H]:=
 \bigl(
 M_1\otimes V_1[!H],
 M_2\otimes V_2[\ast H],
 C\otimes C_V[\ast H]
 \bigr)
\]
\[
 (\nbigt\otimes\nbigt_V)[!H]:=
 \bigl(
 M_1\otimes V_1[\ast H],
 M_2\otimes V_2[!H],
 C\otimes C_V[!H]
 \bigr)
\]
We clearly have
\[
\DDD^{\herm}\bigl((\nbigt\otimes\nbigt_V)[!H]\bigr)
= 
 \DDD^{\herm}(\nbigt)\otimes\DDD^{\herm}(\nbigt_V)[\ast H]
\]
\[
 \DDD^{\herm}\bigl((\nbigt\otimes\nbigt_V)[\ast H]\bigr)
= 
 \DDD^{\herm}(\nbigt)\otimes\DDD^{\herm}(\nbigt_V)[!H].
\]
\begin{lem}
$(\nbigt\otimes\nbigt_V)[\star H]$ $(\star=\ast,!)$
are also non-degenerate.
\end{lem}
\pf
Let us consider the case $\star=\ast$.
By the previous lemma,
the restriction of the induced morphism
$M_1\otimes V_1[!H]
\lrarr
 \nbigc_{\Xbar}\bigl(
 M_2\otimes V_2[\ast H]
 \bigr)$ to $X\setminus H$
is an isomorphism.
Hence, 
$M_1\otimes V_1[!H]
\lrarr
 \nbigc_{\Xbar}\bigl(
 M_2\otimes V_2[\ast H]
 \bigr)$ is an isomorphism.
The other case can be shown similarly.
\hfill\qed

\subsection{Beilinson functors for $D$-triples}
\label{subsection;11.3.21.11}

\index{Beilinson functors}
\index{functor $\Pi^{a,b}_{g \star}$}
\index{functor $\Pi^{a,b}_{g \ast \bikkuri}$}
\index{functor $\psi^{(a)}_g$}
\index{functor $\Xi^{(a)}_g$}
\index{functor $\phi^{(a)}_g$}

Let us consider a meromorphic flat bundle
$\IItilde^{a,b}:=\bigoplus_{a\leq j<b}
 \nbigo_{\cnum_z}(\ast z)\,s^j$
on $(\cnum_z,0)$
with a connection
$z\del_zs^j=s^{j+1}$.
We have the non-degenerate hermitian pairing
$\IItilde^{-b+1,-a+1}$ and $\IItilde^{a,b}$
with values in 
$\nbigc^{\infty\,\moderate 0}_{\cnum_z}$
given as follows:
\[
 \Ctilde(s^i,\sbar^j)
=\frac{(\log|z|^2)^{-i-j}}{(-i-j)!}
 \chi_{i+j\leq 0}.
\]
Here,
$\chi_{i+j\leq 0}$ is $1$ if $i+j\leq 0$,
or $0$ otherwise.
For a holomorphic function $g$ on $X$,
we set
$\IItilde^{a,b}_g:=g^{\ast}\IItilde^{a,b}$.
For a $D_X$-module $M$,
we have a $D_X(\ast g)$-module
$\Pi^{a,b}_gM:=M\otimes\IItilde^{a,b}_g$.
We set
$\Pi^{a,b}_{g\star}M:=
 (\Pi^{a,b}_gM)[\star g]$
for $\star=\ast,!$.

Let $\nbigt=(M_1,M_2,C)\in\Dtriplecat(X)$.
We obtain objects 
\[
\Pi_{g!}^{a,b}\nbigt=\bigl(
 \Pi_{g\ast}^{-b+1,-a+1}M_1,
 \Pi_{g!}^{a,b}M_2,
 \Pi^{a,b}_{g!}C\bigr),
\]
\[
 \Pi_{g\ast}^{a,b}\nbigt=\bigl(
 \Pi_{g!}^{-b+1,-a+1}M_1,
 \Pi_{g\ast}^{a,b}M_2,
 \Pi^{a,b}_{g\ast}C\bigr)
\]
in $\Dtriplecat(X)$.
As in the case of $\nbigr$-triples,
we set 
\begin{multline}
\Pi^{a,b}_{g\ast!}(\nbigt):=
\varprojlim_{N\to\infty}
 \Cok\Bigl(
 \Pi^{b,N}_{g!}(\nbigt)
\lrarr
 \Pi^{a,N}_{g\ast}(\nbigt)
 \Bigr)
 \\
\simeq
 \varinjlim_{N\to\infty}
 \Ker\Bigl(
 \Pi^{-N,b}_{g!}(\nbigt)
\lrarr
 \Pi^{-N,a}_{g\ast}(\nbigt)
 \Bigr)
\end{multline}
If $\nbigt$ is non-degenerate,
$\Pi^{a,b}_{g\star}(\nbigt)$ ($\star=\ast,!$)
and $\Pi^{a,b}_{g\ast !}(\nbigt)$
are also non-degenerate.

In particular, we define
$\psi^{(a)}_g(\nbigt):=\Pi^{a,a}_{g\ast!}(\nbigt)$
and $\Xi^{(a)}_g(\nbigt):=
 \Pi^{a,a+1}_{g\ast!}(\nbigt)$.
As in the case of $\nbigr$-triples,
we define
$\phi_g^{(0)}(\nbigt)$ as the cohomology of
the naturally obtained complex:
\[
 \nbigt[!g]
\lrarr
 \Xi^{(0)}_g(\nbigt)\oplus\nbigt
\lrarr
 \nbigt[\ast g]
\]
We can recover $\nbigt$ as the cohomology of
the following complex:
\[
 \psi^{(1)}(\nbigt)
\lrarr
 \Xi^{(0)}(\nbigt)\oplus\phi^{(0)}(\nbigt)
\lrarr
 \psi^{(0)}(\nbigt)
\]
If $\nbigt$ is non-degenerate,
$\psi^{(a)}(\nbigt)$,
$\Xi^{(a)}(\nbigt)$
and $\phi^{(0)}(\nbigt)$
are also non-degenerate.
It is easy to see the following.
\begin{lem}
We have
$\DDD^{\herm}\circ\Pi^{a,b}_{g\ast!}
=\Pi^{-b+1,-a+1}_{g\ast!}\circ\DDD^{\herm}$.
In particular,
we have
$\DDD^{\herm}\circ\Xi^{(a)}_{g}
=\Xi^{(-a)}_{g}\circ\DDD^{\herm}$
and
$\DDD^{\herm}\circ\psi^{(a)}_{g}
=\psi^{(-a+1)}_{g}\circ\DDD^{\herm}$.
We also have
$\DDD^{\herm}\circ\phi^{(0)}_{g}
=\Xi^{(0)}_{g}\circ\DDD^{\herm}$.
\hfill\qed
\end{lem}

\begin{rem}
For regular holonomic $D$-modules,
the compatibilities of the functors
in {\rm\S\ref{section;13.8.16.1}}
are essentially contained in
{\rm\cite{sabbah-hermitian-duality}}.
\hfill\qed
\end{rem}

\section{De Rham functor}

\subsection{$\cnum_X$-complex-triples}

\index{$\cnum_X$-complex-triple}

Let $I_X^{\bullet}$ be a $\cnum_X$-complex.
A $(\cnum_X,I_X^{\bullet})$-complex-triple is a tuple
$(\nbigf_1^{\bullet},\nbigf_2^{\bullet},C)$,
where $\nbigf_i^{\bullet}$ are $\cnum_X$-complexes,
and a pairing
$C:\nbigf_1^{\bullet}\otimes\nbigf_2^{\bullet}
\lrarr I_X^{\bullet}$
such that 
$d C(x^i\otimes y^j)
=C(dx^i\otimes y^j)+(-1)^iC(x^i\otimes dy^j)$.

Let $I_X^{\bullet,\bullet}$ be a $\cnum_X$-double-complex.
A $(\cnum_X,I_X^{\bullet,\bullet})$-double-complex-triple is a tuple 
$(\nbigf_1^{\bullet,\bullet},\nbigf_2^{\bullet,\bullet},C)$,
where $\nbigf_i^{\bullet,\bullet}$ are
$\cnum_X$-double complexes,
and $C:\nbigf^{\bullet,\bullet}_1
\otimes\nbigf^{\bullet,\bullet}_2
\lrarr I_X^{\bullet,\bullet}$
be a morphism of double complexes,
i.e.,
\[
 d_1C(x^{\vecp},y^{\vecq})
=C(d_1x^{\vecp},y^{\vecq})
+(-1)^{p_1}C(x^{\vecp},d_1y^{\vecq})
\]
\[
 d_2C(x^{\vecp},y^{\vecq})
=C(d_2x^{\vecp},y^{\vecq})
+(-1)^{p_2}C(x^{\vecp},d_2y^{\vecq})
\]

Let $J_X^{\bullet}$ be the total complex of
$I_X^{\bullet,\bullet}$.
For any $(\cnum_X,I_X^{\bullet,\bullet})$-double-complex-triple
$(\nbigf_1^{\bullet,\bullet},\nbigf_2^{\bullet,\bullet},C)$,
we have the associated
$(\cnum_X,J_X^{\bullet})$-complex-triple,
given as follows.
The underlying $\cnum_X$-complexes
are the total complexes of
$\nbigf^{\bullet,\bullet}_i$.
We set
$\Ctilde(x^{\vecp},y^{\vecq}):=
 (-1)^{p_2q_1}C(x^{\vecp},y^{\vecq})$.
We can check that they give
a $(\cnum_X,J^{\bullet}_X)$-complex-triple
by a direct computation, 
as follows:
{\footnotesize
\begin{multline}
 d\Ctilde(x^{\vecp},y^{\vecq})
=d_1\Ctilde(x^{\vecp},y^{\vecq})
+(-1)^{p_1+q_1}d_2\Ctilde(x^{\vecp},y^{\vecq})
 \\
=(-1)^{p_2q_1}d_1C(x^{\vecp},y^{\vecq})
+(-1)^{p_1+q_1+p_2q_1}d_2C(x^{\vecp},y^{\vecq})
 \\
\!=\!(-1)^{p_2q_1}\Bigl(\!
 C(d_1x^{\vecp},y^{\vecq})
 +\!(-1)^{p_1}C(x^{\vecp},d_1y^{\vecq})
 \!\Bigr)
\!+\!(-1)^{p_1+q_1+p_2q_1}\Bigl(\!
 C(d_2x^{\vecp},y^{\vecq})
+\!(-1)^{p_2}C(x^{\vecp},d_2y^{\vecq})
 \!\Bigr)
 \\
=(-1)^{p_2q_1}C(d_1x^{\vecp},y^{\vecq})
+(-1)^{p_1+p_2+(q_1+1)p_2}
 C(x^{\vecp},d_1y^{\vecq})
 \\
+(-1)^{(p_2+1)q_1}C\bigl(
 (-1)^{p_1}d_2x^{\vecp},y^{\vecq}
 \bigr)
+(-1)^{p_1+p_2+p_2q_1}
 C(x^{\vecp},(-1)^{q_1}d_2y^{\vecq})
\\
=\Ctilde(d_1x^{\vecp},y^{\vecq})
+\Ctilde((-1)^{p_1}d_2x^{\vecp},y^{\vecq})
+(-1)^{p_1+p_2}\Bigl(
 \Ctilde(x^{\vecp},d_1y^{\vecq})
+\Ctilde(x^{\vecp},(-1)^{q_1}d_2y^{\vecq})
 \Bigr)
 \\
=\Ctilde(dx^{\vecp},y^{\vecq})
+(-1)^{p_1+p_2}
 \Ctilde(x^{\vecp},dy^{\vecq})=0
\end{multline}
}

For any $(\cnum_X,I^{\bullet}_X)$-triple
$(\nbigf_1^{\bullet},\nbigf_2^{\bullet},C)$
and any $\ell\in\seisuu$,
the shift
$(\nbigf_1^{\bullet}[-\ell],\nbigf_2^{\bullet}[\ell],C[\ell])$
is given as in \S\ref{subsection;13.8.3.10},
which is denoted by
$\nbigs_{\ell}(\nbigf_1^{\bullet},\nbigf_2^{\bullet},C)$.

\subsection{De Rham functor for $D_X$-complex-triples}
\label{subsection;13.4.5.1}
\index{de Rham functor}
Put $\omega_X^{top}:=\distribution_X^{\bullet}[2d_X]$.
We shall use the identification
$\Tot\Bigl(
 \Omega_X^{\bullet}[d_X]
 \otimes
 \Omega_{\Xbar}^{\bullet}[d_X]
 \otimes\distribution_X
 \Bigr)
\simeq
 \omega_X^{top}$
given as in \S\ref{subsection;13.8.17.1}.
For a $D_X$-triple $(M_1,M_2,C)$,
we define the pairing
$\DR(C):\DR(M_1)\otimes \overline{\DR(M_2)}\lrarr
\omega_X^{top}$
as follows:
\[
 \DR(C)(\eta^{d_X+p}\,m_1,\overline{\eta^{d_X+q}\,m_2}):=
 \eta^{d_X+p}\etabar^{d_X+q}\,
 C(m_1,\overline{m}_2)\,
 \epsilon(d_X)
 (-1)^{pd_X}
\]
A $(\cnum_X,\omega_X^{top})$-complex-triple
$\DR(M_1,M_2,C)$
is obtained.

Let $(M_1^{\bullet},M_2^{\bullet},C)$
be a $D_X$-complex-triple.
By the de Rham functor,
we have the induced
 $(\cnum_X,\omega_X^{top})$-double-complex-triple
$\bigl(
 \DR(M_1^{\bullet}),\DR(M_2^{\bullet}),
 \DR(C) \bigr)$.
Thus, we obtain a
$(\cnum_X,\omega_X^{top})$-complex-triple
\[
 \DR(M_1^{\bullet},M_2^{\bullet},C):=
 \Tot\bigl(
 \DR(M_1^{\bullet}),\DR(M_2^{\bullet}),\DR(C)
 \bigr).
\]

We have the following compatibility
with the Hermitian adjoint.
\begin{lem}
\label{lem;11.3.25.1}
Let $(M_1^{\bullet},M_2^{\bullet},C)$ be a $D_X$-complex-triple.
We have 
\[
 (-1)^{d_X}\overline{\Tot\DR(C)}\circ\exchange
=\Tot\DR(C^{\ast}),
\]
where
$\exchange$ denotes the exchange map
$A^{\bullet}\otimes B^{\bullet}\lrarr
 B^{\bullet}\otimes A^{\bullet}$
given by
$a\otimes b\longmapsto
 (-1)^{\deg(a)\deg(b)}b\otimes a$.
\end{lem}
\pf
Let $\eta^{a}$ be local sections of $\Omega^a$,
and $m^{\ell}_i$ be local sections of $M^{\ell}_i$.
We have the following equalities:
{\small
\begin{multline}
 \Tot\bigl(\DR(C^{\ast})\bigr)\bigl(
 \eta^{d_X+p}m^{\ell}_2,\,\overline{\eta^{d_X+q}m_1^{-\ell}}
 \bigr)
=(-1)^{\ell q}
\DR(C^{\ast})\bigl(
 \eta^{d_X+p}m^{\ell}_2,\,\overline{\eta^{d_X+q}m_1^{-\ell}}
 \bigr)
 \\
=\eta^{d_X+p}\etabar^{d_X+q}\,
 C^{\ast}\bigl(m^{\ell}_2,\overline{m}^{-\ell}_1\bigr)\,
 \epsilon(d_X)\,(-1)^{pd_X+\ell q}
 \\
=\etabar^{d_X+q}\eta^{d_X+p}
 \overline{C(m_1^{-\ell},\overline{m}_2^{\ell})}
\epsilon(d_X)\,(-1)^{(p+d_X)(q+d_X)+pd_X+\ell q+\ell}
 \\
=\overline{\Tot\DR(C)(\eta^{d_X+q}m_1^{-\ell},
 \overline{\eta^{d_X+p}m_2^{\ell}})}
 \, (-1)^{qd_X+pd_X+(p+d_X)(q+d_X)+\ell p+\ell q+\ell}
 \\
=
\overline{\Tot\DR(C)(\eta^{q+d_X}m_1^{-\ell},
 \overline{\eta^{p+d_X}m_2^{\ell}})}
 \, (-1)^{d_X+(\ell+p)(-\ell+q)}
\end{multline}
}
Thus, we are done.
\hfill\qed

\vspace{.1in}
We have natural identifications
$\DR(M_i[\ell])\simeq
 \DR(M_i)[\ell]$
for any $\ell\in\seisuu$.
We can check the following lemma
by using Lemma \ref{lem;13.4.5.30}.

\begin{lem}
The identifications induce
$\nbigs_{\ell}\DR\simeq\DR\circ\nbigs_{\ell}$.
\hfill\qed
\end{lem}

\subsection{The de Rham functor and the push-forward}

Let $f:X\lrarr Y$ be a morphism of complex manifolds.
We have the trace morphism
$\tr:f_!\omega_X^{top}\lrarr \omega_Y^{top}$,
which is given by the natural integral
$f_!\distribution_X^{\bullet}[2d_X]
\lrarr
 \distribution_Y^{\bullet}[2d_Y]$
multiplied with
$(2\pi\sqrt{-1})^{d_X-d_Y}$.
Let $(\nbigf_1,\nbigf_2,C)$ be a $c$-soft 
$(\cnum_X,\omega_X^{top})$-complex-triple.
The following naturally induced pairing
is denoted by $f_{\ast}C$:
\[
\begin{CD}
 f_{!}\nbigf_1\otimes f_!\nbigf_2
@>>>
 f_!\omega_X^{top}
@>{\tr}>>
 \omega_Y^{top}
\end{CD}
\]
Thus, we obtain 
a $(\cnum_Y,\omega_Y^{top})$-complex triple
$f_!(\nbigf_1,\nbigf_2,C):=(f_!\nbigf_1,f_{!}\nbigf_2,f_{\ast}C)$.

\begin{prop}
Take $\gbigt=(M_1^{\bullet},M_2^{\bullet},C)\in\Dtriplecat(X)$
such that the restriction of $f$ to the support of
$\gbigf$ is proper.
We have a natural quasi-isomorphism
\[
 \DR_Y f^{(0)}_{\dagger}(M^{\bullet}_1,M^{\bullet}_2,C)
\simeq
f_{!}\circ \DR_X(M^{\bullet}_1,M^{\bullet}_2,C).
\]
\end{prop}
\pf
We have the standard quasi-isomorphisms
of the underlying complexes.
We have only to compare the pairings.
We have only to consider the issue
for a $D_X$-triple $(M_1,M_2,C)$.
Let us argue the case that $f$ is a closed immersion.
Let $\eta_N$ be a local generator of $\omega_Y/\omega_X$.
The quasi-isomorphism
$f_{\ast}\DR_X(M)\lrarr
 \DR_Y(f_{\dagger}M)$
is given by
$\xi\otimes m \longmapsto
\eta_N\cdot \xi\otimes (m/\eta_N)$.
Let $\ell=d_Y-d_X$.
Then, we have the following:
{\small
\begin{multline}
 \DR f^{(0)}_{\dagger}C\bigl(
 \eta_N\,\xi^p\otimes(m_1/\eta_N),\,
 \overline{\eta_N\,\xi^q(m_2/\eta_N)}
 \bigr) 
 \\
=\eta_N\etabar_N \xi^p\xibar^q
 f^{(0)}_{\dagger}\bigl(
 C(m_1/\eta_N,\overline{m_2/\eta_N})
\bigr)
 \epsilon(d_Y)(-1)^{(p+\ell-d_Y)d_Y+\ell p}
 \\
=f_{\ast}\Bigl(\xi^p\xibar^q
 C(m_1,\overline{m}_2)
 \Bigr)
 \epsilon(d_Y)\epsilon(d_Y-d_X)
 (-1)^{(p+\ell-d_Y)d_Y+\ell p}
 \left(\frac{1}{2\pi\sqrt{-1}}\right)^{-\ell}
 \\
=
 f_{\ast}\DR(C)
 \bigl(\xi^p\,m_1,\overline{\xi^q\,m_2}\bigr)
\epsilon(d_X)(-1)^{(p-d_X)d_X}
 \epsilon(d_Y)(-1)^{(p+\ell-d_Y)d_Y+\ell p}
\epsilon(d_Y-d_X)
\end{multline}
}
It is easy to see
$\epsilon(d_X) \epsilon(d_Y)\epsilon(d_Y-d_X)
(-1)^{(p-d_X)d_X+(p+\ell-d_Y)d_Y+\ell p}
=1$.
Hence, we are done in this case.

\vspace{.1in}

Let us consider the case $f:X=Z\times Y\lrarr Y$.
Let $\dim Z=m$ and $\dim Y=n$.
We have the following:
{\small
\begin{multline}
 \label{eq;11.3.19.1}
 \DR_Y(f^{(0)}_{\dagger}C)\bigl(
 [m_1\eta^{m-p}]\,\omega^r,\,
 \overline{[m_2\eta^{m+p}]\,\omega^s}
 \bigr)\,
 (-1)^{p(n-r)}
 \\
=f^{(0)}_{\dagger}C
 \bigl(
 [m_1\eta^{m-p}],\,
 \overline{[m_2\eta^{m+p}]}
 \bigr)\cdot
 \omega^r\omegabar^s
 \epsilon(n)(-1)^{(r-n)n+p(n-r)}
 \\
=\int\eta^{m-p}\etabar^{m+p}
 C(m_1,\overline{m}_2)\,
 \omega^r\omegabar^s
 \epsilon(m)
 \epsilon(n)
 (-1)^{pm+(r-n)n+p(n-r)}
\left(\frac{1}{2\pi\sqrt{-1}}\right)^{m}
\\
=\int\eta^{m-p}\omega^r\etabar^{m+p}\omegabar^s
 C(m_1,\overline{m}_2)\,
 \epsilon(m)
 \epsilon(n)
(-1)^{r(m+p)-pm+n(r-n)+p(n-r)}
\left(\frac{1}{2\pi\sqrt{-1}}\right)^{m}
\end{multline}
}
We also have the following:
{\small
\begin{multline}
\label{eq;11.3.19.2}
f_{\ast}\DR(C)\bigl(
 m_1\eta^{m-p}\omega^r,\,
 \overline{m_2\eta^{m+p}\omega^s}
 \bigr) =
\\
  \int_f\bigl(
 \eta^{m-p}\omega^r
 \overline{\eta^{m+p}\omega^s}
C(m_1,\overline{m}_2)
 \bigr)
 \epsilon(m+n)
 (-1)^{(m+n)(r+m-p-m-n)}
 \left(\frac{1}{2\pi\sqrt{-1}}\right)^{m}
\end{multline}
}
It is easy to check that (\ref{eq;11.3.19.1})
and (\ref{eq;11.3.19.2}) are equal.
Hence, we have
$\DR(F_{\dagger}C)=F_{\ast}\DR(C)$.
\hfill\qed

\section{Dual of $D$-triples}
\label{subsection;13.4.12.30}

\subsection{Dual for non-degenerate $D$-triples}
\label{subsection;11.1.22.10}

\index{dual}

We set $d_X:=\dim X$.
For a $D_X$-module $M$, we have an object
\[
 \DDD_XM:=
 \nrhom_{D_X}(M,D_X\otimes\omega_X^{-1})[d_X]
\]
in the derived category $D^b(D_X)$.
\index{functor $\DDD_X$}
If $M$ is holonomic,
we have a natural identification
$\DDD_XM\simeq\nbigh^0\DDD_XM$,
and it is holonomic.
Because
\[
 \DR_{\Xbar}\nbigc_X(M)
\simeq
 \nhom_{\nbigd_X}
 \bigl(M,\distribution_X^{0,\bullet}\bigr)[d_X],
\quad
\DR_X\DDD_XM
\simeq
 \nrhom_{\nbigd_X}(M,\nbigo_X)[d_X]
\]
we have a natural isomorphism
$\nu_M:\DR_{\Xbar}\nbigc_X(M)\simeq
 \DR_X\DDD_XM$ in
$D^b_c(\cnum_X)$,
induced by the canonical isomorphism
$\distribution_X^{0,\bullet}\simeq\nbigo_X$ 
in $D^b(D_X)$.

Let $(M_1,M_2,C)\in\Dtriplecat^{nd}(X)$.
We have the induced isomorphism
$\varphi_{C,1,2}:\Mbar_2\simeq \nbigc_X(M_1)$.
Hence, we have an isomorphism
\[
 \nu_{M_1}\circ\DR\varphi_{C,1,2}:
 \DR_{\Xbar}(\Mbar_2)
\simeq\DR_X\DDD_X(M_1).
\]
We have an induced isomorphism
$\nu_{\Mbar_2}\circ\DR(\varphi_{C,2,1}):
 \DR_X(M_1)\simeq
 \DR_{\Xbar}\DDD_{\Xbar}(\Mbar_2)$
similarly.
We will prove the following proposition
in \S\ref{subsection;11.1.21.40}.

\begin{thm}
\label{thm;10.9.14.2}
There uniquely exists a non-degenerate 
hermitian pairing
 $\DDD C:\DDD_XM_1\times 
 \DDD_{\Xbar}\Mbar_2\lrarr
 \distribution_X$ such that
 the following diagram is commutative.
\begin{equation}
 \label{eq;10.9.14.1}
 \begin{CD}
 \DR_X\DDD_XM_1\otimes
 \DR_{\Xbar}\DDD_{\Xbar}\Mbar_2
 @>{\DR (\DDD C)}>>
 \omega_X^{top}
\\
 @V{\simeq}V{\kappa}V @V{=}VV \\
 \DR_{\Xbar}(\Mbar_2)\otimes
 \DR_X(M_1)
 @>{\DR(C)\circ\exchange}>> 
\omega_X^{top}
 \end{CD}
\end{equation}
Here, $\kappa$ is induced by
$\nu_{M_1}\circ\DR_X(\varphi_{C,1,2})$
and $\nu_{\Mbar_2}\circ\DR_X(\varphi_{C,2,1})$,
and $\exchange$ denotes the isomorphism
$A^{\bullet}\otimes B^{\bullet}
\simeq
 B^{\bullet}\otimes A^{\bullet}$
given by $x\otimes y\longmapsto
 (-1)^{\deg(x)\deg(y)}y\otimes x$.

The correspondence
$\nbigt=(M_1,M_2,C)\longmapsto
\DDD\nbigt:=(\DDD M_1,\DDD M_2,\DDD C)$
gives a contravariant functor 
on $\Dtriplecat^{nd}(X)$.
We have $\DDD(\DDD C)=C$
under the natural identifications
$\DDD(\DDD M_i)=M_i$.
\end{thm}
Note that we have only to show the existence
of a hermitian pairing $\DDD C$
such that the diagram (\ref{eq;10.9.14.1})
is commutative.
The other claims immediately follow.

\begin{example}
\label{example;13.4.1.11}
Let us consider the case that $M_i$ $(i=1,2)$
are smooth $D$-modules,
i.e., the case in which
$M_i$ are locally free $\nbigo_X$-modules.
Let $L_i$ $(i=1,2)$ be the corresponding local systems.
Then, $\DR_X(M_i)\simeq L_i[d_X]$.
The non-degenerate hermitian pairing 
$C:M_1\times \Mbar_2\lrarr \distribution_X$
corresponds to a non-degenerate hermitian pairing
$C_0:L_1\otimes\Lbar_2\lrarr\cnum_X$.
We have the dual 
$C_0^{\lor}:
 L_1^{\lor}\otimes\Lbar_2^{\lor}
 \lrarr\cnum_X$.

Let us look at the dual $\DDD C$ corresponds to 
$(-1)^{d_X}C_0^{\lor}$.
The pairing $\DR(C)$ on 
$L_1[d_X]\otimes \Lbar_2[d_X]$
is given by
$\epsilon(d_X)(-1)^{d_X}C_0$.
The signature comes from that of
$\Omega_X^{\bullet}[d_X]\otimes
 \Omega_{\Xbar}^{\bullet}[d_X]
\lrarr \distribution_X^{\bullet}[2d_X]$.
Hence,
$\DR_X(C)\circ\exchange$
on $\Lbar_2[d_X]\otimes L_1[d_X]$
is given by
$(\bbar,a)\longmapsto
 \epsilon(d_X)C_0(a,\bbar)$.
Let $\Psi_{C_021}:L_1^{\lor}\simeq \Lbar_2$
and $\Psi_{C_012}:L_2^{\lor}\simeq \Lbar_1$
denote the morphisms induced by 
the non-degenerate pairing $C_0$.
The isomorphism 
$\DR_{\Xbar}(\Mbar_2)\simeq
 \DR_X\DDD_XM_1$
is the composition of
\[
 \DR_{\Xbar}(\Mbar_2)
\simeq
 \nrhom_{D_X}(M_1,
 \Omega_{\Xbar}[d_X]
 \otimes_{\nbigo_{\Xbar}}\distribution_X) 
\simeq
 \nrhom_{D_X}(M_1,\nbigo_X[d_X]).
\]
The induced map
$\Lbar_2[d_X]\simeq L_1^{\lor}[d_X]$
is the shift of $\Psi_{C_021}$.
Similarly, the isomorphism
$\DR_X(M_1)\simeq
 \DR_{\Xbar}(\DDD M_2)$ is the composite of
\[
 \DR_{X}(M_1)\simeq
 \nrhom_{D_{\Xbar}}(\Mbar_2,
 \Omega_X^{\bullet}[d_X]\otimes_{\nbigo_{X}}\distribution_X)
\simeq
 \nrhom_{D_{\Xbar}}(\Mbar_2,\nbigo_{\Xbar}[d_X]).
\]
The induced map
$L_1\simeq \Lbar_2^{\lor}$
is the shift of $\Psi_{C_012}$.
Then, we obtain that
$(-1)^{d_X}C_0^{\lor}$
gives $\DDD C_0$.
\hfill\qed
\end{example}

\begin{rem}
\label{rem;13.4.13.1}
We have the induced pairing
$C\circ\exchange:
 \Mbar_2\otimes M_1\lrarr \distribution_X$.
Note that
$\DR(C\circ\exchange)
=\DR(C)\circ\exchange$
as pairings
$\DR_{\Xbar}(\Mbar_2)\times \DR_X(M_1)
\lrarr \omega^{top}$,
if we use the above identification
$\Tot\bigl(
 \Omega_{X}^{\bullet}[d_X]
 \otimes
 \Omega_{\Xbar}^{\bullet}[d_X]
\otimes\distribution_X
 \bigr)
\simeq
 \distribution_X[2d_X]$
given by
$\xi^{d_X+p}\otimes\xibar^{d_X+q}\otimes\omega
\longleftrightarrow
 \epsilon(d_X)(-1)^{pd_X}
 \xi^{d_X+p}\wedge
 \xibar^{d_X+q}\wedge\omega$.
Note that the induced identification
$\Tot\bigl(
 \Omega_{\Xbar}^{\bullet}[d_X]
 \otimes
 \Omega_{X}^{\bullet}[d_X]
\otimes\distribution_X
 \bigr)
\simeq
 \distribution_X[2d_X]$
is given by
\[
 \xibar^{d_X+q}\otimes\xi^{d_X+p}\otimes\omega
\longleftrightarrow
 \epsilon(d_X)(-1)^{qd_X+d_X}
 \xibar^{d_X+q}\wedge
 \xi^{d_X+p}\wedge\omega.
\]
\hfill\qed
\end{rem}

\subsubsection{Another formulation}

Let $(M_1,M_2,C)\in\Dtriplecat^{nd}(X)$.
We have an induced isomorphism
$\varphi_C:\Mbar_2\simeq\nbigc_X(M_1)$,
which induces
$\DDD\varphi_C:
 \DDD_{\Xbar}\nbigc_X(M_1)
\simeq
 \DDD_{\Xbar}\Mbar_2$.
According to Theorem \ref{thm;10.9.14.2},
we have a dual
$(\DDD_X M_1,\DDD_X M_2,\DDD C)$,
from which we obtain an isomorphism
$\varphi_{\DDD C}:\DDD_{\Xbar} \Mbar_2
\simeq
 \nbigc_X(\DDD_XM_1)$.
By the composition, we obtain an isomorphism
\[
\Xi_M:
 \nbigc_X(\DDD_XM_1)
\simeq
 \DDD_{\Xbar}\nbigc_X(M_1).
\]
\begin{prop}
The following diagram is commutative:
\[
 \begin{CD}
 \DR_{\Xbar}\nbigc_X\DDD_XM_1
 @>{\DR\Xi_M}>>
 \DR_{\Xbar}\DDD_{\Xbar}\nbigc_XM_1
 \\
 @V{a_1}VV @V{a_2}VV \\
 \DDD\DR_X\DDD_XM_1
@>{a_3}>>
 \DDD \DR_{\Xbar}\nbigc_XM_1
 \end{CD}
\]
Here, $a_i$ $(i=1,2,3)$ are naturally induced morphisms.
\end{prop}
\pf
Let us consider the composite 
of the following morphisms:
\begin{multline}
\label{eq;13.8.17.2}
 \DR_X(\DDD M_1)
\stackrel{\DR\varphi_{\DDD C}}{\lrarr}
 \DR_X(\nbigc_{\Xbar}\DDD_{\Xbar}\Mbar_2)
\simeq
 \DDD\DR_{\Xbar}\DDD_{\Xbar}\Mbar_2
 \\
\simeq
 \DDD \DR_X\nbigc_{\Xbar}\Mbar_2
\stackrel{\DDD \DR\varphi_C}{\lrarr}
\DDD \DR M_1
\end{multline}
By the condition of $\DDD C$,
it is equal to the composite of the following morphisms:
\begin{equation}
\label{eq;13.8.17.3}
 \DR_X(\DDD M_1)
\simeq
 \DR_{\Xbar}(\nbigc_XM_1)
\stackrel{\DR_X\varphi_C^{-1}}{\lrarr}
 \DR_{\Xbar}(\Mbar_2)
\stackrel{b}{\lrarr}
 \DDD \DR_X(M_1)
\end{equation}
Here, $b$ is the morphism
induced by $\DR(C)\circ\exchange$.
We can check that the composite of
(\ref{eq;13.8.17.3}) is equal to 
the naturally induced morphism.
Then, the claim of the proposition follows.
\hfill\qed

\subsubsection{Appendix}

Let us recall the dual functor for $D$-modules,
very briefly.
See \cite{kashiwara_text}, \cite{saito1},
\cite{saito3} and \cite{saito4} for more details.
Let $X$ be a complex manifold,
and let $H$ be a hypersurface.

\vspace{.1in}
Let $N$ be a left-$D_{X(\ast H)}$-bi-module,
i.e., it is equipped with mutually commuting 
two $D_{X(\ast H)}$-actions $\rho_i$ $(i=1,2)$.
The left $D_{X(\ast H)}$-module by $\rho_i$
is denoted by $(N,\rho_i)$.
For a $D_{X(\ast H)}$-module $L$,
let $\nhom_{D_{X(\ast H)}}(L^{\bullet},N^{\rho_1,\rho_2})$
denote  the sheaf of $D_{X(\ast H)}$-homomorphisms
$L\lrarr (N,\rho_1)$.
It is equipped with a $D_{X(\ast H)}$-action
induced by $\rho_2$.
Note that,
for a $D_{X(\ast H)}$-complex $L$,
we have the naturally defined
$D_{X(\ast H)}$-homomorphism
\begin{equation}
 \label{eq;11.1.22.2}
 L^{\bullet}\lrarr \nhom_{D_{X(\ast H)}}\Bigl(
 \nhom_{D_{X(\ast H)}}\bigl(L^{\bullet},N^{\rho_2,\rho_1}
 \bigr),
 N^{\rho_1,\rho_2}
 \Bigr)
\end{equation}
given by
$x\longmapsto
 \bigl(
 F\longmapsto
 (-1)^{|x|\,|F|}F(x)
 \bigr)$.

\vspace{.1in}

We have the natural two $D_{X(\ast H)}$-action
on $D_{X(\ast H)}\otimes\omega_X^{-1}$.
The left multiplication is denoted by $\ell$.
Let $r$ denote the action induced by the right
multiplication.
More generally, for a left $D_{X(\ast H)}$-module $N$,
we have two induced left $D_{X(\ast H)}$-module
structure on 
$N\otimes D_{X(\ast H)}\otimes\omega_X^{-1}$.
One is given by $\ell$ and
the left $D$-action on $N$,
which is denoted by $\ell$.
The other is induced by the right multiplication,
denoted by $r$.
We have the automorphism
of $\Phi_N:N\otimes D_{X(\ast H)}\otimes\omega_X^{-1}$,
which exchanges $\ell$ and $r$,
as in \cite{saito1}.

Let $\nbigg^{\bullet}$ be a left-$D_{X(\ast H)}$-injective
resolution of $\nbigo_X(\ast H)[d_X]$.
Let $M$ be a coherent $D_{X(\ast H)}$-module.
We define
\[
 \DDD M:=
 \nhom_{D_{X(\ast H)}}\Bigl(M,
 \bigl(\nbigg^{\bullet}\otimes
 D_{X(\ast H)}\otimes\omega_X^{-1}\bigr)^{\ell,r}
 \Bigr). 
\]
The morphism as in (\ref{eq;11.1.22.2})
with $\Phi_{\nbigg^{\bullet}}$ induces
an isomorphism
$M\lrarr \DDD_{X(\ast H)}\circ\DDD_{X(\ast H)}(M)$.

\vspace{.1in}

Let us consider the case $H=\emptyset$.
For any perverse sheaf $\nbigf$,
put $\DDD_{\cnum_X}(\nbigf):=
 \nrhom_{\cnum_X}(\nbigf,\omega_X^{top})$.
If $M$ is holonomic,
we have the naturally defined isomorphism
$\DR_X\DDD_X(M)\simeq \DDD_{\cnum_X}\DR_X(M)$.
See \cite{saito4} for more details on this isomorphism.
We can deduce that
the following diagram is commutative:
\[
 \begin{CD}
 \DR_X\DDD_X\circ\DDD_X M
@>{\simeq}>>
 \DDD_{\cnum_X}\circ\DDD_{\cnum_X}\DR_X M\\
 @AAA @AAA \\
 \DR_X M@>{=}>>
 \DR_X M
 \end{CD}
\]
The vertical arrows are induced by
$\DDD\circ\DDD\simeq\id$,
and the upper horizontal arrow is induced 
by the exchange of the functors $\DR$ and $\DDD$.

\subsection{Dual of complexes of non-degenerate 
$D_X$-triples}
\label{subsection;13.4.5.20}

%We consider
%$\nbigt^{\bullet}\in \nbigc(\Dtriplecat(X))$.
Assuming Theorem \ref{thm;10.9.14.2},
we define an object
$\DDD(\nbigt^{\bullet})$
in $\nbigc(\Dtriplecat^{nd}(X))$
for any $\nbigt^{\bullet}\in \nbigc(\Dtriplecat(X))$.
Its $p$-th member is given by
\[
 \DDD(\nbigt^{\bullet})^p
=\bigl(
 \DDD M_1^p,\DDD M_2^{-p}, (-1)^p\DDD C^{-p}
 \bigr).
\]
The differential
$\deltatilde^p:\DDD(\nbigt^{\bullet})^p
\lrarr
 \DDD(\nbigt^{\bullet})^{p+1}$
is given by 
\[
\begin{CD}
 \DDD(M_1^p)
 @<{(-1)^{p}\DDD(\delta_{-p-1,1})}<<
  \DDD(M_1^{p+1}),
\end{CD}
\quad
\mbox{\rm and}
\quad
\begin{CD}
 \DDD(M_2^{-p})
 @>{(-1)^{p+1}\DDD(\delta_{-p-1,2})}>>
  \DDD(M_2^{-p-1}).
\end{CD}
\]
Note that the signature is adjusted 
so that the underlying $D$-complexes
are $\DDD(M_i^{\bullet})$.

\vspace{.1in}
We have the saturated full subcategory
$\Dcomplextriplecat^{nd}_0(X)
\subset
 \Dcomplextriplecat(X)$
corresponding to
$\nbigc\bigl(
 \Dtriplecat^{nd}(X)
 \bigr)\subset
 \nbigc\bigl(
 \Dtriplecat(X)
 \bigr)$
by $\Psi_1$.
\index{category $\Dcomplextriplecat^{nd}_0(X)$}
We obtain a contravariant auto equivalence
$\DDD$ on $\Dcomplextriplecat^{nd}_0(X)$.
We have
\[
  \DDD(M_1^{\bullet},M_2^{\bullet},C)
=(\DDD M_1^{\bullet},\DDD M_2^{\bullet},
 \DDD C),
\]
where 
$(\DDD C)^p=\DDD (C^{-p})$.

\subsection{Compatibility of the push-forward and the dual}
\label{subsection;11.4.6.1}

To state the stability for push-forward,
we introduce the notion of dual for $D$-complex-triples,
which is more general than that in 
\S\ref{subsection;13.4.5.20}.
Let $\gbigt=(M_1^{\bullet},M_2^{\bullet},C)
\in\Dcomplextriplecat(X)$ be such that
each $H^p(\gbigt)\in\Dtriplecat^{nd}(X)$.
We say that $\gbigt$ has a dual
$\DDD\gbigt$,
if there exists a pairing
$\DDD C:\DDD M_1^{\bullet}\times
 \overline{\DDD M_2^{\bullet}}
\lrarr
 \distribution_X$
such that the following is commutative
in $D^b_c(\cnum_X)$:
\[
  \begin{CD}
 \Tot\DR_X\DDD_XM^{\bullet}_1\otimes
 \Tot\DR_{\Xbar}\DDD_{\Xbar}\Mbar^{\bullet}_2
 @>{\DR (\DDD C)}>>
 \omega_X^{top}
\\
 @V{\simeq}VV @V{=}VV \\
 \Tot\DR_{\Xbar}(\Mbar_2^{\bullet})\otimes
 \Tot\DR_X(M_1^{\bullet})
 @>{\DR(C)\circ\exchange}>> 
\omega_X^{top}
 \end{CD}
\]
Here, the left vertical arrow is 
induced as in the diagram (\ref{eq;10.9.14.1}).
The full subcategory of such $\gbigt$
in $\Dcomplextriplecat(X)$
is denoted by 
$\Dcomplextriplecat^{nd}_1(X)$.
\index{category $\Dcomplextriplecat^{nd}_1(X)$}

\vspace{.1in}
We have the following stability
of the category of the non-degenerate pairings
with respect to the push-forward,
and the compatibility of the dual and the push-forward,
which we will prove in
\S\ref{subsection;11.1.21.40}.
\begin{thm}
\label{thm;11.3.30.1}
Let $F:X\lrarr Y$ be a morphism
of complex manifolds.
Let $\nbigt=(M_1,M_2,C)\in\Dtriplecat^{nd}(X)$.
Assume that the restriction of $F$
to the support of $\nbigt$ is proper.
Then, we have 
$F^{(0)}_{\dagger}\nbigt
\in\Dcomplextriplecat_1^{nd}(Y)$,
and 
$F^{(0)}_{\dagger}\DDD\nbigt$
gives a dual of $F^{(0)}_{\dagger}\nbigt$
in the above sense
under the natural quasi-isomorphism of
the $D$-complexes.
In particular,
$\vecH^jF^{(0)}_{\dagger}\DDD\nbigt$
is a dual of
$\DDD \vecH^{-j}F^{(0)}_{\dagger}\nbigt$
under the natural isomorphism
of the $D$-modules.
\end{thm}
(See \S\ref{subsection;11.3.23.1}
for $F^{(0)}_{\dagger}$.)
As a consequence,
we obtain the following
for $D$-triples.

\begin{thm}
\label{thm;10.9.14.3}
Let $F$ and $\nbigt$
be as in Theorem {\rm\ref{thm;11.3.30.1}}.
Then, we have
$F^j_{\dagger}\nbigt\in \Dtriplecat^{nd}(Y)$
for any $j$.
Moreover,
we have
$F^j_{\dagger}\DDD\nbigt[-j]
\simeq
 \DDD \bigl(
 F^{-j}_{\dagger}\nbigt[j]
 \bigr)$,
i.e.,
we have
$F^{-j}_{\dagger}\DDD C=(-1)^j\DDD F^j_{\dagger}C$
under the natural identifications
$F^j_{\dagger}\DDD M_1\simeq
 \DDD F^{-j}_{\dagger}M_1$
and 
$F^{-j}_{\dagger}\DDD M_2\simeq
 \DDD F^{j}_{\dagger}M_2$.
\end{thm}

\begin{cor}
Suppose $\nbigt^{\bullet}\in\Dcomplextriplecat^{nd}_1(X)$.
Then, 
$F_{\dagger}^{(0)}(\nbigt^{\bullet})$ is an object
 in $\Dcomplextriplecat^{nd}_1(Y)$,
and we have
$\DDD F_{\dagger}^{(0)}\nbigt^{\bullet}
\simeq
 F_{\dagger}^{(0)}\DDD\nbigt^{\bullet}$
under the natural identifications 
of the underlying $D_Y$-complexes.
\hfill\qed
\end{cor}

\subsection{Compatibility with other functors}

For any $\nbigt\!=\!(M_1,M_2,C)\in\Dtriplecat(X)$,
let $\nbigt^{\ast}:=(M_2,M_1,C^{\ast})$,
where $C^{\ast}(m_2,\overline{m_1})=
 \overline{C(m_1,\overline{m_2})}$.
Once we know Theorem \ref{thm;10.9.14.2},
it is easy to obtain the compatibility with
the Hermitian adjoint and the dual.
\begin{prop}
\label{prop;11.3.23.21}
For any $\nbigt\in\Dtriplecat^{nd}(X)$,
we have $\DDD (C^{\ast})=(\DDD C)^{\ast}$,
i.e., 
$\DDD(\nbigt^{\ast})
\simeq
 (\DDD\nbigt)^{\ast}$
under the natural identification of 
the underlying $D_X$-modules.
\hfill\qed
\end{prop}

Once we know Theorem \ref{thm;10.9.14.2},
it is easy to obtain the compatibility
with the localization
in  \S\ref{subsection;10.12.27.22}.
\begin{prop}
\label{prop;10.11.17.10}
Let $H$ be a hypersurface of $X$.
For $\nbigt\in\Dtriplecat^{nd}(X)$,
we have natural isomorphisms
\[
 \DDD\bigl(\nbigt[!H]\bigr)
\simeq
 \bigl(\DDD\nbigt\bigr)[\ast H],
\quad\quad
 \DDD\bigl(\nbigt[\ast H]\bigr)
\simeq
 \bigl(\DDD\nbigt\bigr)[!H].
\]
\end{prop}
\pf
Let us show $\DDD\bigl(\nbigt[!H]\bigr)
\simeq \bigl(\DDD\nbigt\bigr)[\ast H]$.
We have only to compare the pairings
$\DDD\bigl(C[!H]\bigr)$ and
$\bigl(\DDD C\bigr)[\ast H]$
under the natural identifications
$\DDD\bigl(M_1[\ast H]\bigr)\simeq
\bigl(\DDD M_1\bigr)[!H]$
and 
$\DDD\bigl(M_2[!H]\bigr)\simeq
\bigl(\DDD M_2\bigr)[\ast H]$.
It is easy to compare their restrictions to $X\setminus H$.
We obtain the comparison on $X$
by the uniqueness of the extension.
\hfill\qed

\vspace{.1in}

Let us consider the compatibility with
the tensor product of a non-degenerate smooth $D$-triple
in \S\ref{subsection;11.3.21.10}.
Let $\nbigt_V=(V_1,V_2,C)$ be a smooth $D_X$-triple.
If $C$ is non-degenerate,
the tuple $\nbigt_V$
is locally isomorphic to a direct sum of
the trivial pairing
$(\nbigo,\nbigo,C_0)$,
where $C_0(f,\overline{g})=f\,\overline{g}$.
We put
$V_i^{\lor}:=
 \nhom_{\nbigo_X}\bigl(V_i,\nbigo_X\bigr)$.
If $C$ is non-degenerate,
we have the induced pairing
$C^{\lor}:V_1^{\lor}\otimes
\overline{V_2^{\lor}}\lrarr
 \nbigc^{\infty}_X\lrarr\distribution_X$.
In this case,
for $\nbigt=(V_1,V_2,C)$,
the dual $(V_1^{\lor},V_2^{\lor},C^{\lor})$
is denoted by $\nbigt_V^{\lor}$.
Once we know Theorem \ref{thm;10.9.14.2},
the following lemma is obvious.
\begin{lem}
\label{lem;11.2.20.20}
Let $\nbigt_i\in\Dtriplecat^{nd}(X)$ $(i=1,2)$.
Assume that
$\nbigt_2$ is smooth.
Then,
we naturally  have
\[
 \DDD(\nbigt_1\otimes\nbigt_2)
\simeq
 \DDD(\nbigt_1)\otimes\nbigt_{2}^{\lor}. 
\]
\end{lem}
\pf
We have the natural isomorphisms
of the underlying holonomic $D$-modules,
which is reviewed in Remark \ref{rem;13.4.3.10} below.
For the comparison of the pairings,
we have only to consider the case
$\nbigt_2=(\nbigo,\nbigo,C_0)$,
and then it is clear.
\hfill\qed

\vspace{.1in}

Let us consider a more general case.
Let $H$ be a hypersurface of $X$.
Let $\nbigt_V=(V_1,V_2,C)$ be 
a smooth $D_{X(\ast H)}$-triple.
We define the dual $\nbigt_V^{\lor}$
as in the case of smooth $D_X$-triple.
Once we know Theorem \ref{thm;10.9.14.2},
the following lemma is easy to see.
\begin{lem}
\label{lem;10.12.27.23}
Let $\nbigt\in\Dtriplecat^{nd}(X)$.
Let $\nbigt_V$ be a smooth non-degenerate
$D_{X(\ast H)}$-triple.
Then, we have natural isomorphisms:
\[
\DDD\bigl((\nbigt\otimes\nbigt_V)[!t]\bigr)
\simeq
 \bigl(
 \DDD(\nbigt)\otimes\nbigt_{V}^{\lor}
\bigr)[\ast t],
\quad\quad
\DDD\bigl((\nbigt\otimes\nbigt_V)[\ast t]\bigr)
\simeq
 \bigl(
 \DDD(\nbigt)\otimes\nbigt_{V}^{\lor}
\bigr)[!t].
\]
\end{lem}
\pf
The underlying $D$-modules 
of 
$\DDD\bigl((\nbigt\otimes\nbigt_V)[\ast t]\bigr)$
and
$\bigl(
 \DDD(\nbigt)\otimes\nbigt_{V}^{\lor}
\bigr)[!t]$ are naturally isomorphic.
(See Remark \ref{rem;13.4.3.10} below.)
For comparison of the pairings,
we have only to compare their restrictions
to $X\setminus H$,
which follows from Lemma \ref{lem;11.2.20.20}.
The other case can be checked similarly.
\hfill\qed

\begin{rem}
\label{rem;13.4.3.10}
Let $X$ and $H$ be as above.
Let $M$ be a holonomic $D_X(\ast H)$-module.
Let $V$ be a smooth $D_X(\ast H)$-module.
We have the natural isomorphism
$\DDD_{X(\ast H)}(M\otimes_{\nbigo_X(\ast H)} V)\simeq
 \DDD_{X(\ast H)}(M)\otimes_{\nbigo_X(\ast H)} V^{\lor}$,
given as follows.
We have two naturally induced
left $D$-module structure
$\ell$ and $r$
on $V^{\lor}\otimes_{\nbigo_X}
 (D_{X(\ast H)}\otimes\omega_X^{-1})$,
as explained in the appendix of 
{\rm\S\ref{subsection;11.1.22.10}}.
We have
\begin{multline}
 \DDD_{X(\ast H)}(M\otimes V)
\simeq
 \nrhom_{D_X(\ast H)}\bigl(M\otimes V,
 (D_{X(\ast H)}\otimes\omega_X^{-1})^{\ell,r}
 \bigr)[d_X] \\
\simeq 
 \nrhom_{D_X(\ast H)}\bigl(M,
 (V^{\lor}\otimes 
 D_{X(\ast H)}\otimes\omega_X^{-1})^{\ell,r}
 \bigr)[d_X]
\end{multline}
By $\Phi_{V^{\lor}}$, it is isomorphic to
\begin{multline}
 \nrhom_{D_X(\ast H)}\bigl(
 M, (V^{\lor}\otimes 
 D_{X(\ast H)}\otimes\omega_X^{-1})^{r,\ell}
 \bigr)[d_X]
\simeq \\
 \nrhom_{D_X(\ast H)}\bigl(
 M,(D_{X(\ast H)}
 \otimes\omega_X^{-1})^{r,\ell}
 \bigr)
\otimes V^{\lor}[d_X]. 
\end{multline}
By $\Phi_{\nbigo_X(\ast H)}$,
it is isomorphic to
$\nrhom_{D_X(\ast H)}\bigl(
 M,(D_{X(\ast H)}\otimes\omega_X^{-1})^{\ell,r}
 \bigr)
\otimes V^{\lor}[d_X]$.
Thus, we obtain the desired isomorphism
$\DDD_{X(\ast H)}(M\otimes V)\simeq
 \DDD_{X(\ast H)}(M)\otimes V^{\lor}$.
We can easily deduce
$\DDD_X(M\otimes V)
\simeq
\bigl(\DDD_X(M)\otimes V\bigr)[!H]$
and
$\DDD_X\bigl((M\otimes V)[!H]\bigr)
\simeq
\bigl(\DDD_X(M)\otimes V\bigr)$.
\hfill\qed
\end{rem}

Let us observe the compatibility of the dual
and Beilinson's functors.
We have the perfect pairing
$\IItilde^{a,b}\times\IItilde^{-b+1,-a+1}\lrarr
 \nbigo_{\cnum_z}(\ast z)$
given as follows:
\[
 \bigl\langle f(s),g(s)
 \bigr\rangle
=\underset{s=0}{\Res}
 \Bigl(f(s)g(-s)\frac{ds}{s}\Bigr)
\]
It gives an identification of the dual
of $\IItilde^{a,b}$ with $\IItilde^{-b+1,-a+1}$.
As in \S\ref{subsection;11.2.1.10},
we can  check
\[
 \bigl(
 \IItilde^{-b-1,-a-1},\IItilde^{a,b},\Ctilde_{\II}
 \bigr)^{\lor}
\simeq
 \bigl(
 \IItilde^{a,b},\IItilde^{-b+1,-a+1},\Ctilde_{\II}
 \bigr). 
\]
Then, once we know Theorem \ref{thm;10.9.14.2},
we obtain the compatibility of the dual
and Beilinson's functors
in \S\ref{subsection;11.3.21.11}.
\begin{prop}
\mbox{{}}
Let $\nbigt\in\Dtriplecat^{nd}(X)$.
\begin{itemize}
\item
We have
$\DDD\bigl(
 \Pi^{a,b}_{g\ast!}(\nbigt)\bigr)
\simeq
 \Pi^{-b+1,-a+1}_{g\ast!}(\DDD\nbigt)$.
\item
In particular,
$\psi_g^{(a)}(\nbigt)$,
$\Xi_g^{(a)}(\nbigt)$
and $\phi_g^{(0)}(\nbigt)$
are non-degenerate,
and we have natural isomorphisms:
\[
 \DDD\psi^{(a)}_g(\!\nbigt\!)
\simeq
 \psi^{(-a+1)}_g(\!\DDD\nbigt\!),
\,\,\,
 \DDD\Xi^{(a)}_g(\!\nbigt\!)
\simeq
 \Xi^{(-a)}_g(\!\DDD\nbigt\!),
\,\,\,
 \DDD\phi^{(0)}_g(\!\nbigt\!)
\simeq
 \phi^{(0)}_g(\!\DDD\nbigt\!).
\]
\end{itemize}
\end{prop}
\pf
By Lemma \ref{lem;10.12.27.23},
$\Pi^{a,N}_{g\star}(\nbigt)$ are 
non-degenerate.
Hence, we obtain that
$\Pi^{a,b}_{g\ast!}\nbigt$ are non-degenerate.
We have the following diagram:
\[
 \begin{CD}
 \Pi^{-N+1,-b+1}_{g!}(\DDD\nbigt)
 @>>>
 \Pi^{-N+1,-a+1}_{g\ast}(\DDD\nbigt)
 @>>>
 \Pi^{-b+1,-a+1}_{g\ast!}(\DDD\nbigt)
 \\
 @V{\simeq}VV @V{\simeq}VV @VVV \\
 \DDD\Pi^{b,N}_{g\ast}(\nbigt)
 @>>>
 \DDD\Pi^{a,N}_{g!}(\nbigt)
 @>>>
 \DDD\Pi^{a,b}_{g\ast!}(\nbigt)
 \end{CD}
\]
Hence, we obtain
$\Pi^{-b+1,-a+1}_{g\ast!}(\DDD\nbigt)
\simeq \DDD\Pi^{a,b}_{g\ast!}(\nbigt)$.
\hfill\qed

\subsection{Push-forward and dual of $D$-modules (Appendix)}

We recall some details on the compatibility
of the push-forward and the dual.
(See \cite{kashiwara_text} and \cite{saito4}.)

\subsubsection{Preliminary}

Let $Y$ be a complex manifold.
Let $\ell$ and $r$ denote the left and right multiplication
of $D_Y$ on itself.
We have a unique $\cnum$-linear isomorphism
$\Psi:D_Y\otimes\omega_Y^{-1}\lrarr D_Y\otimes\omega_Y^{-1}$
determined by
$\Psi\circ r(Q)=\ell(Q)\circ\Psi$
and $\Psi\circ \ell(Q)=r(Q)\circ\Psi$
for any $Q\in D_Y$.

We set
$\nbigm_{Y}:=(D_Y\otimes\omega_Y^{-1})
(\lefttop{\ell}\otimes^\ell_{\nbigo_Y})
 (D_Y\otimes\omega_Y^{-1})$,
where we use the $\nbigo$-action induced by
$\ell$ for the tensor product.
We have the natural left $D_Y$-action
$\ell\otimes\ell$:
For $v\in\Theta_Y$, we have
$(\ell\otimes\ell)(v)(P_1\otimes P_2)
=vP_1\otimes P_2+P_1\otimes vP_2$.
We also have the right $D_Y$-actions $r_i$ $(i=1,2)$:
For $Q\in D_Y$,
we have 
$r_1(Q)(P_1\otimes P_2)=\bigl(r(Q)(P_1)\bigr)\otimes P_2$
and 
$r_2(Q)(P_1\otimes P_2)=P_1\otimes \bigl(r(Q)(P_2)\bigr)$

We have a unique $\cnum$-linear isomorphism
$F:\nbigm_{Y}\lrarr \nbigm_{Y}$
determined by the conditions
$r_2(Q)\circ F=F\circ r_2(Q)$
and
$(\ell\otimes\ell)(Q)\circ F
=F\circ r_1(Q)$
for any $Q\in D_Y$.
We have the formula
$F(m_1\otimes m_2)
=(\ell\otimes\ell)\bigl(\Psi(m_1)\bigr)
 (1\otimes m_2)$.

\begin{lem}
We have
$r_1\circ F=F\circ(\ell\otimes\ell)$.
\end{lem}
\pf
We have only to show the formula
for the action of $v\in \Theta_Y$.
We have
\begin{multline}
F(vm_1\otimes m_2+m_1\otimes vm_2)
=(\ell\otimes\ell)\bigl(
 \Psi(vm_1)
 \bigr)(1\otimes m_2)
+(\ell\otimes\ell)\bigl(\Psi(m_1)\bigr)(1\otimes vm_2)
\\
=(\ell\otimes\ell)
 \bigl(\Psi(m_1)\cdot (-v)\bigr)(1\otimes m_2)
+(\ell\otimes\ell)\bigl(\Psi(m_1)\bigr)(1\otimes vm_2)
 \\
=(\ell\otimes\ell)\bigl(\Psi(m_1)\bigr)
(-v\otimes m_2-1\otimes vm_2)
+(\ell\otimes\ell)\bigl(\Psi(m_1)\bigr)
(1\otimes vm_2) \\
=(\ell\otimes\ell)\bigl(\Psi(m_1)\bigr)
 \circ
 r_1(v)(1\otimes m_2)
=r_1(v)\circ F(m_1\otimes m_2). 
\end{multline}
Thus, we are done.
\hfill\qed

\vspace{.1in}
We set
$\nbigm_{0,Y}:=
 (D_Y\otimes\omega_Y^{-1})
 \lefttop{\ell}\otimes^{\ell}
 (D_Y\otimes\omega_Y^{-1})$.
For the tensor product,
the $\nbigo_Y$-actions on the both factors
are induced by $\ell$.

We set $\nbigm_{1,Y}:=(D_Y\otimes\omega_Y^{-1})
 (\lefttop{r}\otimes^{\ell}_{\nbigo_Y})
 (D_Y\otimes\omega_Y^{-1})$.
For the tensor product,
the $\nbigo_Y$-action on the first factor
is induced by $r$,
and that on the second factor
is induced by $\ell$.
We have the natural left actions
$r\otimes\ell$,
$\ell_1$
and $r_2$:
For $v\in\Theta_Y$,
we have
$(r\otimes\ell)(v)(m_1\otimes m_2)
=\bigl(r(v)m_1\bigr)\otimes m_2
+m_1\otimes v\,m_2$,
$\ell_1(v)(m_1\otimes m_2)=vm_1\otimes m_2$
and 
$r_2(v)(m_1\otimes m_2)
=m_1\otimes \bigl(r(v)(m_2)\bigr)$.
We have a unique $\cnum$-linear map
$G:\nbigm_{0,Y}\lrarr \nbigm_{1,Y}$ 
determined by the conditions
$G\circ r_2=r_2\circ G$
and $G\circ (\ell\otimes\ell)
=\ell_1\circ G$.
\begin{cor}
\label{cor;11.3.18.2}
We have
$G\circ r_1
=(r\otimes\ell)\circ G$.
\hfill\qed
\end{cor}

\begin{rem}
Let $G$ be the morphism in 
Corollary {\rm\ref{cor;11.3.18.2}}.
We take $\omega_Y\otimes^L_{D_Y}$
by using $r_2$.
Then, $G$ induces 
the identity
$D_Y\otimes\omega_Y^{-1}
\lrarr D_Y\otimes \omega_Y^{-1}$.
\hfill\qed
\end{rem}

\subsubsection{A morphism}

Let $f:X\lrarr Y$ be a morphism of complex manifolds.
We have two natural left $f^{-1}(D_Y)$-actions
on
\[
 D_{Y\larr X}\otimes^L_{D_X}
 \bigl(
 \nbigo_X[d_X]\otimes_{f^{-1}\nbigo_Y}
 f^{-1}(D_Y\otimes\omega_Y^{-1})
 \bigr),
\]
where we consider the left multiplication of
$f^{-1}\nbigo_Y$ on $f^{-1}(D_Y\otimes\omega_Y^{-1})$
for the tensor product $\otimes_{f^{-1}}$.
One is induced by that on $D_{Y\larr X}$,
denoted by $\kappa_1$.
The other is induced by that
on $f^{-1}(D_Y\otimes\omega_Y^{-1})$,
denoted by $\kappa_2$.
We have two natural left $f^{-1}(D_Y)$-actions on
\[
 \bigl(
 D_{Y\larr X}\otimes_{D_X}^L
 \nbigo_X[d_X]
\bigr)\otimes_{f^{-1}\nbigo_Y}
 f^{-1}(D_Y\otimes\omega_Y^{-1}),
\]
where we consider the left multiplication of
$f^{-1}\nbigo_Y$ on $f^{-1}(D_Y\otimes\omega_Y^{-1})$
for the tensor product $\otimes_{f^{-1}}$.
One is induced by those on
$D_{Y\larr X}\otimes_{D_X}^L
 \nbigo_X[d_X]$
and $f^{-1}(D_Y\otimes\omega_Y^{-1})$
denoted by $\kappa_1'$,
and the other is induced by
that on $f^{-1}(D_Y\otimes\omega_Y^{-1})$,
denoted by $\kappa_2'$.

\begin{lem}
\label{lem;11.3.18.1}
We have a natural $\cnum$-linear isomorphism
\begin{multline}
 F:
 D_{Y\larr X}\otimes^L_{D_X}
 \bigl(
 \nbigo_X[d_X]\otimes_{f^{-1}\nbigo_Y}
 f^{-1}(D_Y\otimes\omega_Y^{-1})
 \bigr)
 \\
\simeq
 \bigl(
 D_{Y\larr X}\otimes_{D_X}^L
 \nbigo_X[d_X]
\bigr)\otimes_{f^{-1}\nbigo_Y}
 f^{-1}(D_Y\otimes\omega_Y^{-1})
\end{multline}
such that $F\circ \kappa_i=\kappa_i'\circ F$
$(i=1,2)$.
\end{lem}
\pf
We have
$D_{Y\larr X}=
 \omega_X\otimes
 f^{-1}(D_Y\otimes \omega_Y^{-1})
\simeq
 \bigl(
 \Omega_X^{\bullet}[d_X]\otimes D_X
 \bigr)\otimes_{f^{-1}\nbigo_Y}
 f^{-1}(D_Y\otimes\omega_Y^{-1})$.
We have the following:
\begin{multline}
 \label{eq;11.3.19.10}
 D_{Y\larr X}\otimes_{D_X}^L
 \bigl(
 \nbigo_X[d_X]\otimes_{f^{-1}\nbigo_Y}
 f^{-1}(D_Y\otimes\omega_Y^{-1})
 \bigr) \\
\simeq
 \Bigl(
 \Omega_X^{\bullet}[d_X]\otimes D_X\otimes_{f^{-1}\nbigo_Y}
 f^{-1}(D_Y\otimes\omega_Y^{-1})
 \Bigr)
\otimes_{D_X}
 \Bigl(
 \nbigo_X[d_X]\otimes_{f^{-1}\nbigo_Y}
 f^{-1}(D_Y\otimes\omega_Y^{-1})
 \Bigr)
 \\
\simeq
\Omega_X^{\bullet}[d_X]
 \otimes_{f^{-1}\nbigo_Y}
 f^{-1}\bigl(
 \nbigm_{0,Y}
 \bigr)[d_X]
\end{multline}
In the last term, 
we consider 
the $D_Y$-action $\ell\otimes\ell$-action 
on $\nbigm_{0,Y}$
for the tensor product $\otimes_{f^{-1}\nbigo_Y}$.
We also have the following:
\begin{multline}
 \bigl(
D_{Y\larr X}\otimes_{D_X}^L
 \nbigo_X[d_X]
\bigr)
\otimes_{f^{-1}\nbigo_Y}
 f^{-1}(D_Y\otimes\omega_Y^{-1})
 \\
\simeq
 \Bigl(
 \bigl(
 \Omega_X^{\bullet}[d_X]
 \otimes D_X\otimes_{f^{-1}\nbigo_Y}
 f^{-1}(D_Y\otimes\omega_Y^{-1})
\bigr)
 \otimes_{D_X}\nbigo_X[d_X]
 \Bigr)
\otimes_{f^{-1}\nbigo_Y}
 f^{-1}\bigl(D_Y\otimes\omega_Y^{-1}\bigr)
 \\
\simeq
 \Omega_X^{\bullet}[d_X]
\otimes_{f^{-1}\nbigo_Y}
 f^{-1}\Bigl(
 \nbigm_{1,Y}
 \Bigr)[d_X]
\end{multline}
In the last term,
we consider 
the $f^{-1}(\nbigo_Y)$-action
induced by the $D_Y$-action $\ell_1$-action
on $\nbigm_{1,Y}$
for the tensor product 
$\otimes_{f^{-1}\nbigo_Y}$.
Then, we obtain Lemma \ref{lem;11.3.18.1}
from Corollary \ref{cor;11.3.18.2}.
\hfill\qed

\subsubsection{Compatibility of the push-forward and the dual}

Recall that we have the trace morphism
$f_{\dagger}(\nbigo_{X}[d_X])
\lrarr
 \nbigo_Y[d_Y]$,
induced by the trace
$f_!\distribution_X^{\bullet}[2d_X]
\lrarr
 \distribution_Y^{\bullet}[2d_Y]$.
Indeed,
\begin{multline}
\label{eq;11.3.21.20}
f_{\dagger}\bigl(
 \nbigo_X[d_X]
 \bigr)
\lrarr
 f_{!}\bigl(
 \distribution^{\bullet}_X[2d_X]
  \bigr)
\otimes_{f^{-1}\nbigo_Y}
 f^{-1}(D_Y\otimes\omega_Y^{-1})
\\
\lrarr
 \distribution_Y^{\bullet}[2d_Y]
\otimes 
 D_Y\otimes\omega_Y^{-1}
\simeq\nbigo_Y[d_Y]
\end{multline}
From Lemma \ref{lem;11.3.18.1}
and (\ref{eq;11.3.21.20}),
we obtain the following trace morphism
\begin{multline}
\label{eq;11.3.21.21}
 f_{!}\bigl(D_{Y\larr X}\otimes_{D_X}^L
 (\nbigo_X[d_X]\otimes_{f^{-1}\nbigo_Y}
 f^{-1}(D_Y\otimes\omega_Y^{-1}))
 \bigr) \\
\simeq
 f_{!}\Bigl(
 (D_{Y\larr X}\otimes_{D_X}^L\nbigo_X[d_X])
 \Bigr)
 \otimes_{\nbigo_Y}(D_Y\otimes\omega_Y^{-1})
 \\
\lrarr
 \nbigo_Y[d_Y]\otimes D_Y\otimes\omega_Y^{-1}
=D_Y\otimes\omega_Y^{-1}[d_Y]
\end{multline}
Suppose that $M$ is a $D_X$-module
whose support is proper relative to $f$.
By using this trace morphism (\ref{eq;11.3.21.21}),
we obtain
{\small
\begin{multline}
 f_{\dagger}\circ \DDD M
\simeq
 f_{!}\Bigl(
 \nrhom_{D_X}\bigl(M,\,
 \nbigo_X[d_X]\otimes_{f^{-1}\nbigo_Y}
 f^{-1}(D_Y\otimes\omega_Y^{-1})\bigr)
 \Bigr) \\
\lrarr
 f_{!}\Bigl(
 \nrhom_{f^{-1}(D_Y)}\Bigl(
 D_{Y\larr X}\otimes_{D_X}^LM,\,
 D_{Y\larr X}\otimes_{D_X}^L\bigl(
 \nbigo_X[d_X]\otimes_{f^{-1}\nbigo_Y}
 f^{-1}(D_Y\otimes\omega_Y^{-1})
 \bigr)
 \Bigr)
 \Bigr) \\
\lrarr 
 \nrhom_{D_Y}\bigl(
 f_{\dagger}M,\,
 f_{!}(D_{Y\larr X}\otimes_{D_X}^L(\nbigo_X[d_X]
 \otimes_{f^{-1}\nbigo_Y}
 f^{-1}(D_Y\otimes\omega_Y^{-1})))
 \bigr)
 \\
\lrarr
 \nrhom_{D_Y}\bigl(f_{\dagger}M,\,
 D_Y\otimes\omega_Y^{-1}[d_Y]\bigr)
\simeq
 \DDD f_{\dagger}M
\end{multline}
}
The composite is an isomorphism,
if $M$ is good relatively to $f$.
See \cite{kashiwara_text} for example.

\begin{rem}
\label{rem;11.3.18.11}
We obtain the following morphism
by applying 
``$\omega_Y\otimes^L_{D_Y}$''
to {\rm(\ref{eq;11.3.21.21})}
with the $D_Y$-action $r_2$:
\begin{multline}
 \omega_Y\otimes^L_{D_Y}
 \Bigl(
 f_{!}\bigl(D_{Y\larr X}\otimes^L_{D_X}\nbigo_X[d_X]
 \bigr)
  \otimes D_Y\otimes\omega_Y^{-1}
 \Bigr)
\lrarr \\
 \omega_Y\otimes^L_{D_Y}
 \bigl(\nbigo_Y[d_Y]\otimes D_Y\otimes\omega_Y^{-1}\bigr)
=\nbigo_Y[d_Y]
\end{multline}
It is the same as the trace morphism
$f_{\dagger}\nbigo_X[d_X]
\lrarr
 \nbigo_Y[d_Y]$.
\hfill\qed
\end{rem}

\subsubsection{Compatibility of the de Rham functor, 
the push-forward and the dual}

Let $f:X\lrarr Y$ be a morphism of complex manifolds.
For simplicity,
we assume that $f$ is proper.
We have the following diagram
of the functors
from the category of holonomic $D_X$-modules
to derived category of the cohomologically constructible complexes
on $Y$:
\begin{equation}
\label{eq;11.3.18.10}
 \begin{CD}
 \DR_Y\circ f_{\dagger}\circ\DDD
 @>{\simeq}>>
 \DR_Y\circ\DDD\circ f_{\dagger}\\
 @V{\simeq}VV @V{\simeq}VV \\
 f_{\ast}\circ \DDD\circ \DR_X
 @>{\simeq}>>
 \DDD\circ f_{\ast}\circ \DR_X
 \end{CD}
\end{equation}
Here, the left vertical arrow is
given by
$\DR_Y\circ f_{\dagger}\circ \DDD
\simeq
 f_{\ast}\circ\DR_X\circ\DDD
\simeq
 f_{\ast}\circ \DDD\circ \DR_X$,
and right vertical arrow is given by
$\DR_Y\circ \DDD\circ f_{\dagger}
\simeq
 \DDD\circ\DR_X\circ f_{\dagger}
\simeq
 \DDD\circ f_{\ast}\circ \DR_X$.

\begin{prop}
\label{prop;11.3.21.22}
The diagram {\rm (\ref{eq;11.3.18.10})} is commutative.
\end{prop}
\pf
Let $N=\nbigo_X\otimes_{f^{-1}\nbigo_Y}
 f^{-1}(D_Y\otimes\omega_Y^{-1})[d_X]$.
We have
$\DR_Y f_{\dagger}\DDD M
=\omega_Y\otimes^L_{D_Y}
 f_{\ast}\nrhom_{D_X}(M,N)$.
The morphism
$\DR_Yf_{\dagger}\DDD M
\lrarr
 \DR_Y\DDD f_{\dagger}M
\lrarr
 \DDD f_{\ast}\DR_XM$
is expressed as follows
($\nrhom_{\nbiga}(Q_1,Q_1)$
are denoted by
$(Q_1,Q_2){\nbiga}$):
{\tiny
\[
 \begin{array}{ccccc}
 \omega_Y\otimes_{D_Y}^L
 f_{\ast}(M,N)_{D_X}\!\!\!
 & \lrarr & 
 \omega_Y\otimes_{D_Y}^L
 (f_{\dagger}M,f_{\dagger}N)_{D_Y}
 & \lrarr &
 \!\!\!\omega_Y\otimes^L_{D_Y}
 (\DR_Yf_{\dagger}M,\DR_Yf_{\dagger}N)_{\cnum_Y}
 \!\!\!
 \\
 & &\darr & &\darr \\
 & & 
 \!\!\!\omega_Y\otimes^L_{D_Y}
 (f_{\dagger}M,D_Y\otimes\omega_Y^{-1}[d_Y])_{D_Y}
 \!\!\!
 &\lrarr &
 (\DR_Yf_{\dagger}M,
 \DR_Y\nbigo_Y[d_Y])_{\cnum_Y}
 \\
 & & & &\darr  \\
 & & & &
 (f_{\ast}\DR_XM,
 \DR_Y\nbigo_Y[d_Y])_{\cnum_Y}
 \end{array}
\]
}
The morphism
$\DR_Yf_{\dagger}\DDD M
\lrarr
 f_{\ast}\DDD \DR_XM
\lrarr
 \DDD f_{\ast}\DR_XM$
is expressed as follows:
{\tiny
\[
\begin{array}{ccccc}
 \!\!\!\!
 \omega_Y\otimes_{D_Y}^Lf_{\ast}(M,N)_{D_X}\!\!\!\!
 & \rarr &
 \omega_Y\otimes_{D_Y}^Lf_{\ast}(\DR M,\DR N)_{\cnum_X}
 & \rarr &
 \!\!\!
 \omega_Y\otimes_{D_Y}^L
 (f_{\ast}\DR M,f_{\ast}\DR N)_{\cnum_Y}
 \!\!\!
 \\
 & &\darr & &\darr \\
 & & 
 \!\!\!
 f_{\ast}(\DR M,f^{-1}\omega_Y\otimes_{f^{-1}D_Y}^L
 \!\DR N)_{\cnum_X}
  \!\!\!
 & \rarr &
 \!\!\!(f_{\ast}\DR M,
 \omega_Y\otimes_{D_Y}^Lf_{\ast}\DR N)_{\cnum_Y}
 \!\!\!
 \\
 & & & & \darr \\
 & & & &
 \!\!(f_{\ast}\DR M,\DR_Y\nbigo_Y[d_Y])_{\cnum_Y}\!\!
\end{array}
\]
}
Then,
we can deduce the claim of the lemma.
\hfill\qed

\begin{cor}
Let $\nbigt\in\Dtriplecat^{nd}(X)$.
The following diagram of natural isomorphisms
is commutative:
\[
 \begin{CD}
 \DR_Y\circ H^j f_{\dagger}\circ \DDD(\nbigt)
 @>>>
 \DR_Y \DDD \circ H^{-j}f_{\dagger}(\nbigt)
 \\ 
 @VVV @VVV \\
 H^jf_{\ast}\circ\DDD\circ\DR_X(\nbigt)
 @>>>
\DDD H^{-j}f_{\ast}\circ \DR_X(\nbigt)
 \end{CD}
\]
Here, the cohomology for $\cnum$-complexes
are taken with respect to the middle perversity.
\end{cor}
\pf
We have only to check it
for the underlying $D$-complexes
and $\cnum$-complexes,
which have already been done
in Proposition \ref{prop;11.3.21.22}.
\hfill\qed

\section{Proof of Theorem \ref{thm;10.9.14.2} 
and Theorem \ref{thm;11.3.30.1}}
\label{subsection;11.1.21.40}

\subsection{Preliminary}

We recall the following lemma,
which will be used implicitly.
\begin{lem}
Let $\nbiga$ and $\nbigb$
be sheaves of $\cnum$-algebras on $X$.
Let $\nbigl_1$ be an $\nbiga$-module,
and let $\nbigl_2$ be 
an $(\nbiga\otimes\nbigb)$-injective module.
Then, 
$\nhom_{\nbiga}(\nbigl_1,\nbigl_2)$
is $\nbigb$-injective.
\end{lem}
\pf
Let $J_1\rarr J_2$ be a monomorphism of
$\nbigb$-modules.
Let $J_1\rarr \nhom_{\nbiga}(\nbigl_1,\nbigl_2)$
be any $\nbigb$-morphism.
Note that 
$J_1\otimes_{\cnum}\nbigl_1\lrarr 
 J_2\otimes_{\cnum}\nbigl_1$
is a monomorphism.
Hence, we have a morphism
$J_2\otimes\nbigl_1\lrarr\nbigl_2$
whose restriction to $J_1\otimes\nbigl_1$
is equal to the given morphism.
It means we obtain
$J_2\lrarr \nhom_{\nbiga}(\nbigl_1,\nbigl_2)$
whose restriction to $J_1$ is equal to
the given morphism.
\hfill\qed

\subsection{Push forward and the functor $\nbigc_X$}
\label{subsection;10.9.14.30}

Let $F:X\lrarr Y$ be a morphism of complex manifolds.
Let $M$ be a holonomic $D_X$-module
such that the restriction of $F$ to $\Supp M$
is proper.
We have the natural morphism
$F_{\dagger}\nbigc_X(M)\lrarr
 \nbigc_YF_{\dagger}(M)$ given as follows:
{\small
\begin{multline}
RF_{!}\nrhom_{D_X}\Bigl(M,\,
 D_{\Ybar\larr \Xbar}\otimes_{D_{\Xbar}}^L
 \distribution_X\Bigr)
\lrarr
 \\
RF_{!}\nrhom_{F^{-1}D_Y}\Bigl(
 D_{Y\larr X}\otimes^L_{D_X}M,\,
 D_{Y\larr X,\Ybar\larr \Xbar}
 \otimes^L_{D_{X,\Xbar}}
 \distribution_X
\Bigr)
 \\
\lrarr
 \nrhom_{D_Y}\Bigl(
 F_{\dagger}M,\,
 RF_!\bigl(
 D_{Y\larr X,\Ybar\larr \Xbar}
 \otimes_{D_{X,\Xbar}}^L
 \distribution_X
 \bigr)
 \Bigr)
\lrarr
 \nrhom_{D_Y}\bigl(
 F_{\dagger}M,\distribution_Y
 \bigr)
\end{multline}
}
Here, we have used the trace map 
(\ref{eq;11.3.23.10}) 
in the last map.
Recall the compatibility
of the push-forward and the de Rham functor
$\DR_{Y}\circ F_{\dagger}\simeq
F_{!}\circ \DR_{X}$
given as follows:
\[
 \omega_{Y}\otimes_{D_{Y}}^L
 F_!(D_{Y\larr X}\otimes_{D_{X}}^L N)
\simeq
 F_!\bigl(
  F^{-1}(\omega_{Y})\otimes_{F^{-1}(D_{Y})}^L
   D_{Y\larr X}\otimes_{D_{X}}^L N
 \bigr)
\simeq
 F_!\bigl(
 \omega_{X}\otimes_{D_{X}}^L N
 \bigr)
\]
Hence, we have the induced morphism
\begin{equation}
 \label{eq;10.9.14.10}
 F_{!}\DR_{\Xbar}\nbigc_XM
\simeq
\DR_{\Ybar}F_{\dagger}\nbigc_XM
\lrarr
 \DR_{\Ybar}\nbigc_YF_{\dagger}M.
\end{equation}
We also have the following isomorphism:
\begin{equation}
\label{eq;10.9.14.11}
 F_{!}\DR_{X}\DDD_XM
\simeq
 \DR_YF_{\dagger}\DDD_XM
\lrarr
 \DR_{Y}\DDD_YF_{\dagger}M.
\end{equation}
(See \cite{kashiwara_text}
for the compatibility
$\DDD_X\circ\DR_X
\simeq \DR_X\circ\DDD_X$.)
We shall show the following lemma.
\begin{lem}
\label{lem;10.9.10.7}
The following diagram 
in $D^b_c(\cnum_Y)$ 
is commutative:
\begin{equation}
 \label{eq;10.9.14.12}
 \begin{CD}
 F_{!}\DR_X\DDD_XM_1
 @>{\simeq}>>
 F_{!}\DR_{\Xbar}\nbigc_X M_1\\
 @VVV @VVV \\
 \DR_Y\DDD_YF_{\dagger}M_1
 @>{\simeq}>>
 \DR_{\Ybar}\nbigc_{Y}F_{\dagger}M_1
 \end{CD}
\end{equation}
The vertical arrows are given by
{\rm(\ref{eq;10.9.14.10})}
and {\rm(\ref{eq;10.9.14.11})}.
In particular, the induced morphism
$F_{\dagger}\nbigc_X(M)\lrarr 
\nbigc_YF_{\dagger}(M)$ is an isomorphism.
\end{lem}
\pf
We consider objects
$C_{1}^X:=\nbigo_X[d_X]$
and $C_{2}^X:=
 \omega_{\Xbar}\otimes^L_{D_{\Xbar}}\distribution_X$
in $D^b(D_X)$.
We have expressions
\[
 \DR_X\DDD_X(M)=
 \nrhom_{\nbigd_X}\bigl(M,\,
 C_1^X\bigr),
\quad\quad
\DR_{\Xbar}\nbigc_X(M)
=\nrhom_{\nbigd_X}\bigl(M,\,
 C_{2}^X
 \bigr).
\]
We have the natural isomorphism
$\eta_X:C_{1}^X\simeq C_{2}^X$
in $D^b(\nbigd_X)$
which induces the isomorphism
$\DR_X\DDD_X(M)\simeq
 \DR_{\Xbar}\nbigc_X(M)$.
We have the trace morphisms
$\tr:F_{\dagger}C_{i}^X\lrarr
 C_i^{Y}$ 
in $D^b(D_Y)$,
and the following diagram is commutative
in $D^b(D_Y)$:
\begin{equation}
\label{eq;10.9.14.24}
 \begin{CD}
 F_{\dagger}C_1^X
 @>{F_{\dagger}\eta_X}>>
 F_{\dagger}C_2^X \\
 @V{\tr}VV @V{\tr}VV \\
 C_1^Y @>{\eta_Y}>> C_{2}^Y
 \end{CD}
\end{equation}

\begin{lem}
\label{lem;10.9.14.23}
The natural isomorphisms
$F_{!}\DR_{\Xbar}\nbigc_X(M)\rarr
 \DR_{\Ybar}\nbigc_YF_{\dagger}(M)$
and
$F_{!}\DR_X\DDD_X(M)\rarr
 \DR_Y\DDD_YF_{\dagger}(M)$
are given as the composition
of the following natural morphisms:
\begin{multline}
\label{eq;10.9.14.20}
F_{!}\nrhom_{D_X}\bigl(M,C_i^X\bigr)
\lrarr
 F_{!}\nrhom_{F^{-1}(D_Y)}
 \bigl(
 D_{Y\larr X}\otimes_{D_X}^LM,\,
 D_{Y\larr X}\otimes_{D_X}^L
 C_i^X
 \bigr)\\
\lrarr
 \nrhom_{D_Y}\bigl(
 F_{\dagger}(M),\, 
 F_{\dagger}(C_i^X)
 \bigr)
\lrarr
 \nrhom_{D_Y}\bigl(
 F_{\dagger}(M),\, 
 C_i^Y
 \bigr)
\end{multline}
\end{lem}
\pf
Let us consider the claim for $\nbigc$.
The claim for $\DDD$ can be argued similarly.
Let $\nbigi_1^{\bullet}$ be a
$D_X\otimes F^{-1}(D_{\Ybar})$-injective
resolution of 
$D_{\Ybar\larr\Xbar}\otimes_{D_{\Xbar}}^L
 \distribution_X$
obtained by the Godement construction.
Then,
$F_!\nhom_{D_X}\bigl(
 M,\nbigi_1^{\bullet}
 \bigr)$ represents
\[
 F_!\nrhom_{D_X}(M,D_{\Ybar\larr\Xbar}\otimes^L_{D_{\Xbar}}
 \distribution_X).
\]

We remark that 
$\omega_{\Ybar}$ is represented 
by $\omega_{\Ybar}^{\sim}:=
 \Omega_{\Ybar}^{\bullet}\otimes_{\nbigo_{\Ybar}}
 D_{\Ybar}[d_Y]$,
and $\Omega_{\Ybar}^{j}\otimes_{\nbigo_{\Ybar}}
 D_{\Ybar}$ are right locally $D_{\Ybar}$-free modules.
Then, 
$F^{-1}(\omega^{\sim}_{\Ybar})
 \otimes_{F^{-1}D_{\Ybar}}\nbigi_1^{\bullet}$
is a $D_X$-injective resolution of
$F^{-1}\omega^{\sim}_{\Ybar}
 \otimes_{F^{-1}D_{\Ybar}}
 D_{\Ybar\larr\Xbar}\otimes_{D_{\Xbar}}^L
 \distribution_X
\simeq
 \omega_{\Xbar}\otimes^L_{D_{\Xbar}}
 \distribution_X$.
Hence, the isomorphism 
$\DR_{\Ybar}F_{\dagger}\nbigc_XM\simeq
F_{!}\DR_{\Xbar}\nbigc_XM$
is expressed as follows:
\begin{multline}
\label{eq;10.9.15.1}
 \omega^{\sim}_{\Ybar}\otimes_{D_{\Ybar}}
 F_{!}\nrhom_{D_X}\bigl(M,\,
 D_{\Ybar\larr \Xbar}\otimes^L_{D_{\Xbar}}
 \distribution_X
 \bigr)
\simeq \\
 F_{!}\nrhom_{D_X}\bigl(M,\,
 F^{-1}\omega^{\sim}_{\Ybar}
 \otimes_{F^{-1}D_{\Ybar}}
 D_{\Ybar\larr \Xbar}\otimes^L_{D_{\Xbar}}
 \distribution_X
 \bigr) \\
\simeq
 F_!\nrhom_{D_X}\bigl(M,\omega_{\Xbar}
 \otimes^L_{D_{\Xbar}}\distribution_X\bigr)
\end{multline}
The morphism
$\DR_{\Ybar}F_{\dagger}\nbigc_X(M)
\lrarr
 \DR_{\Ybar}\nbigc_YF_{\dagger}(M)$
is the composition of the following:
\begin{multline}
\label{eq;10.9.14.21}
\omega_{\Ybar}^{\sim}\otimes_{D_{\Ybar}}
F_!\nrhom_{D_X}\bigl(M,\,
 D_{\Ybar\larr \Xbar}
 \otimes_{D_{\Xbar}}^L
 \distribution_X
 \bigr)
\lrarr \\
\omega_{\Ybar}^{\sim}\otimes_{D_{\Ybar}}
F_!\nrhom_{F^{-1}D_Y}\bigl(
 D_{Y\larr X}\otimes^L_{D_X}M,\,
 D_{Y\larr X,\Ybar\larr \Xbar}\otimes^L_{D_{X,\Xbar}}
 \distribution_X
 \bigr) \lrarr
 \\
 \omega_{\Ybar}\otimes_{D_{\Ybar}}^L
 \nrhom_{D_{Y}}\Bigl(
 F_{\ast}(D_{Y\larr X}\otimes^L_{D_X}M),\,
 F_!\bigl(
 D_{Y\larr X,\Ybar\larr X}\otimes^L_{D_{X,\Xbar}}
 \distribution_X
 \bigr)
 \Bigr)
\lrarr \\
 \omega_{\Ybar}
 \otimes^L_{D_{\Ybar}}
 \nrhom_{D_Y}\bigl(F_{\dagger}M,\,
 \distribution_Y
 \bigr)
\end{multline}
Let $A_i$ denote the $i$-th term
in (\ref{eq;10.9.14.20}).
Let $B_i$ denote the $i$-th term
in (\ref{eq;10.9.14.21}).
For each $i=1,2,3$,
we have a natural morphism
$B_i\lrarr A_i$ given as in (\ref{eq;10.9.15.1}).
We also have a natural isomorphism
$B_4\simeq A_4$.
We have only to show the commutativity
of the following diagrams
for $i=1,2,3$:
\begin{equation}
 \label{eq;10.9.14.22}
 \begin{CD}
 B_i @>>> B_{i+1}\\
 @VVV @VVV \\
 A_i @>>> A_{i+1} \\
 \end{CD}
\end{equation}

Let us consider the diagram (\ref{eq;10.9.14.22}) with $i=1,2$.
Let $\nbigi^{\bullet}_1$ and $\omega_{\Ybar}^{\sim}$
be as above.
We take a resolution $D^f_{Y\larr X}$ of $D_{Y\larr X}$,
which is $D_X$-flat and $c$-soft with respect to $F$.
We take a $F^{-1}(D_{Y,\Ybar})$-injective resolution
$\nbigi_2^{\bullet}$
of $D_{Y\larr X,\Ybar\larr \Xbar}
 \otimes^L_{D_{X,\Xbar}}
 \distribution_X$
(see the proof of Lemma \ref{lem;10.9.10.5},
for example),
and a morphism
$D^f_{Y\larr X}
 \otimes_{D_X}\nbigi^{\bullet}_1\lrarr
 \nbigi^{\bullet}_2$,
which represents the natural isomorphism
\[
 D_{Y\larr X}\otimes_{D_X}^L
 \bigl(
  D_{\Ybar\larr \Xbar}\otimes^L_{D_{\Xbar}}
 \distribution
 \bigr)
\lrarr
 D_{Y\larr X,\Ybar\larr \Xbar}
 \otimes^L_{D_{X,\Xbar}}
 \distribution_X. 
\]
Note that
$F^{-1}(\omega^{\sim}_{\Ybar})
 \otimes_{F^{-1}D_{\Ybar}}
 \nbigi^{\bullet}_2$
is an $F^{-1}(D_Y)$-injective resolution of
$D_{Y\larr X}\otimes^L_{D_X}\bigl(
\omega_{\Xbar}\otimes^L_{D_{\Xbar}}
 \distribution_X\bigr)$.
Hence, the diagram (\ref{eq;10.9.14.22}) with $i=1$
is represented by the following diagram,
which is commutative:
{\tiny
\[
 \begin{array}{ccc}
 \omega^{\sim}_{\Ybar}
 \otimes_{D_{\Ybar}}
 F_!\nhom_{D_X}\bigl(
 M,\,\nbigi^{\bullet}_1
 \bigr)
 & \lrarr &
 \omega^{\sim}_{\Ybar}
 \otimes_{D_{\Ybar}}
 F_{!}
 \nhom_{F^{-1}(D_Y)}
 \bigl(
 D^f_{Y\larr X}\otimes M,\nbigi^{\bullet}_2
 \bigr)
 \\
 \darr & & \darr \\
 F_!\nhom_{D_X}\bigl(M,
 F^{-1}(\omega^{\sim}_{\Ybar})
 \otimes_{F^{-1}D_{\Ybar}}
 \nbigi^{\bullet}_1
 \bigr)
 & \lrarr &
 F_!\nhom_{F^{-1}D_Y}\bigl(
 D^{\sim}_{Y\larr X}\otimes_{D_X}M,
 F^{-1}(\omega^{\sim}_{\Ybar})
 \otimes_{F^{-1}D_{\Ybar}}
 \nbigi^{\bullet}_2
 \bigr)
 \end{array}
\]
}
The diagram (\ref{eq;10.9.14.22})
with $i=2$ is expressed by the following diagram,
which is commutative:
{\tiny
\[
 \begin{array}{ccc}
 \omega^{\sim}_{\Ybar}
 \otimes_{D_{\Ybar}}
 F_{!}
 \nhom_{F^{-1}(D_Y)}
 \bigl(
 D^f_{Y\larr X}\otimes M,\nbigi^{\bullet}_2
 \bigr)
 & \rarr & 
 \omega^{\sim}_{\Ybar}
 \otimes_{D_{\Ybar}}
 \nhom_{D_Y}\bigl(
 F_{\dagger}(M),
 F_!\nbigi^{\bullet}_2
 \bigr)
 \\
 \darr & & \darr \\
 F_!\nhom_{F^{-1}D_Y}\bigl(
 D^{f}_{Y\larr X}\otimes_{D_X}M,
 F^{-1}(\omega^{\sim}_{\Ybar})
 \otimes_{F^{-1}D_{\Ybar}}
 \nbigi^{\bullet}_2
 \bigr)
 & \rarr &
 \nhom_{D_Y}\Bigl(\!
 F_{\dagger}M,
 F_!\bigl(
 F^{-1}(\omega^{\sim}_{\Ybar})
 \otimes_{F^{-1}D_{\Ybar}}
 \nbigi^{\bullet}_2\bigr)\!
 \Bigr)
 \end{array}
\]
}
The diagram (\ref{eq;10.9.14.22})
with $i=3$ is commutative
because of the construction of
the trace morphism.
Thus, the proof of Lemma \ref{lem;10.9.14.23}
is finished.
\hfill\qed

\vspace{.1in}
Now, Lemma \ref{lem;10.9.10.7}
follows from Lemma \ref{lem;10.9.14.23}
and the commutativity of (\ref{eq;10.9.14.24}).
\hfill\qed

\subsection{Pairing on the push-forward}

Let $(M_1,M_2,C)\in \Dtriplecat^{nd}(X)$.
Let $F:X\lrarr Y$ be a morphism
such that the restriction of $F$
to $\Supp M_i$ are proper.
We have 
$F_{\dagger}\varphi_{C,1,2}:
 F_{\dagger}\Mbar_2
\lrarr
 F_{\dagger}\nbigc_{X}M_1$,
and the morphism
$b:F_{\dagger}\nbigc_{X}M_1
\simeq
\nbigc_{Y}F_{\dagger}M_1$
in \S\ref{subsection;10.9.14.30}.
We also have $\varphi_{F^{(0)}_{\dagger}C,1,2}:
 F_{\dagger}\Mbar_2\lrarr
 \nbigc_YF_{\dagger}M_1$.

\begin{lem}
\label{lem;10.9.10.5}
We have
$b\circ F_{\dagger}\varphi_{C,1,2}=
\varphi_{F^{(0)}_{\dagger}C,1,2}$.
\end{lem}
\pf
We set $D^{f}_{\Ybar\larr \Xbar}:=
 F^{-1}(D_{\Ybar}\otimes\omega_{\Ybar}^{-1})
 \otimes_{F^{-1}\nbigo_{\Ybar}}
 \nbige_{\Xbar}^{\bullet}
 \otimes_{\nbigo_{\Xbar}}
 D_{\Xbar}[d_X]$.
We use the notation
$D^f_{Y\larr X}$ in a similar meaning.
It naturally gives a resolution 
of $D_{\Ybar\larr \Xbar}$,
which is $D_{\Xbar}$-flat and $c$-soft with respect to $F$.
Then, $\varphi_{F^{(0)}_{\dagger}C,1,2}$
is represented as the composition
of the following morphisms:
{\small
\begin{multline}
F_!\bigl(
 D^f_{\Ybar\larr \Xbar}\otimes_{D_{\Xbar}}
 \Mbar_2
 \bigr) \\
\lrarr 
F_!\nhom_{F^{-1}D_Y} 
 \bigl(
 D^f_{Y\larr X}\otimes_{D_X} M_1,
 F^{-1}(D_{Y,\Ybar}\otimes\omega_{Y,\Ybar}^{-1})
 \otimes_{F^{-1}(\nbigo_{Y,\Ybar})}
 \distribution^{\bullet}_{X}
 \bigr)[2d_X]  \\
\lrarr
 \nhom_{D_Y} 
 \Bigl(
 F_{\ast}\bigl(
 D^f_{Y\larr X}\otimes_{D_X} M_1\bigr),
 D_{Y,\Ybar}\otimes\omega_{Y,\Ybar}^{-1}
 \otimes_{\nbigo_{Y,\Ybar}} 
 F_!
 \distribution^{\bullet}_{X}
 \Bigr)[2d_X]  \\
\lrarr
  \nhom_{D_Y} 
 \Bigl(
 F_{\ast}\bigl(
 D^f_{Y\larr X}\otimes_{D_X} M_1\bigr),
 \distribution_Y
 \Bigr)
\end{multline}
}

We take an $F^{-1}(D_{Y,\Ybar})$-injective
resolution $\nbigj_1^{\bullet}$ of 
$F^{-1}(D_{Y,\Ybar}\otimes\omega_{Y,\Ybar}^{-1})
 \otimes_{F^{-1}\nbigo_{Y,\Ybar}}
 \distribution_X^{\bullet}[2d_X]$
which is isomorphic to
$D_{Y\larr X,\Ybar\larr \Xbar}
 \otimes^L_{D_{X,\Xbar}}
 \distribution_X$
in $D^b(F^{-1}(D_{Y,\Ybar})\!)$.
We take a $D_{Y,\Ybar}$-injective resolution
$\nbigj_2^{\bullet}$ of $\distribution_Y$,
and a morphism
$F_!\nbigj_1^{\bullet}\lrarr 
 \nbigj_2^{\bullet}$ such that
the following diagram in $D^b(D_{Y,\Ybar})$
is commutative:
\[
 \begin{CD}
 F_{!}\Bigl(
 F^{-1}(D_{Y,\Ybar}\otimes\omega^{-1}_{Y,\Ybar})
 \otimes_{F^{-1}\nbigo_{Y,\Ybar}}
 \distribution_X^{\bullet}[2d_X]
 \Bigr)
 @>>>
 \distribution_Y \\
 @VVV @VVV \\
 F_!\nbigj^{\bullet}_1
 @>>>\nbigj_2^{\bullet}
 \end{CD}
\]
Then, 
$\varphi_{F^{(0)}_{\dagger}C,1,2}$
is represented by the composition of the following
morphisms:
{\small
\begin{multline}
 F_{!}(D_{\Ybar\larr\Xbar}^{f}
 \otimes_{D_{\Xbar}}\Mbar_2)
\lrarr 
  F_{!}\nhom_{F^{-1}(D_Y)}
 \Bigl(
 D_{Y\larr X}^{f}\otimes_{D_X}M_1,
 \nbigj_1^{\bullet}
 \Bigr) \\
\lrarr
 \nhom_{D_Y}
 \Bigl(
 F_{!}\bigl(
 D_{Y\larr X}^{f}\otimes_{D_X}M_1\bigr),\,
 F_{!}\nbigj_1^{\bullet}
 \Bigr) 
\lrarr
 \nhom_{D_Y}
 \Bigl(
 F_{!}\bigl(
 D_{Y\larr X}^{f}\otimes_{D_X}M_1\bigr),\,
 \nbigj_2^{\bullet}
 \Bigr) 
\end{multline}
}
Hence, we obtain the following factorization
of $\varphi_{F^{(0)}_{\dagger}C,1,2}$:
{\small
\begin{multline}
\label{eq;10.9.14.32}
 F_{\dagger}\Mbar_2
\stackrel{A}{\lrarr}
 F_{!}\nrhom_{F^{-1}(D_Y)}
 \bigl(
 D_{Y\larr X}\otimes_{D_X}^LM_1,\,
 D_{Y\larr X,\Ybar\larr\Xbar}
 \otimes_{D_{X,\Xbar}}^L
 \distribution_X\,
 \bigr) \\
\stackrel{A_1}{\lrarr}
 \nrhom_{D_Y}\Bigl(
 F_{\dagger}M_1,\,
 F_{!}\bigl(
 D_{Y\larr X,\Ybar,\larr\Xbar}
 \otimes^L_{D_{X,\Xbar}}
 \distribution_X
 \bigr)
 \Bigr)
\stackrel{A_2}{\lrarr}
 \nrhom_{D_Y}
 \bigl(F_{\dagger}M_1,\distribution_Y\bigr)
\end{multline}
}
Let us look at the morphism $A$.
It is obtained as the push-forward 
of the composition $A'$ of the following morphisms:
{\scriptsize
\begin{multline}
 D^f_{\Ybar\larr \Xbar}\otimes_{D_{\Ybar}}\Mbar_2
\lrarr
 \nhom_{F^{-1}(D_Y)}
 \Bigl(
 D^f_{Y\larr X}\otimes_{D_X}M_1,\,
 F^{-1}\bigl(
 D_{Y,\Ybar}\otimes\omega^{-1}_{Y,\Ybar}
 \bigr)\otimes_{F^{-1}\nbigo_{Y,\Ybar}}
 \distribution_X^{\bullet}
 \Bigr)[2d_X]\\
\lrarr
 \nhom_{F^{-1}(D_Y)}\bigl(
 D^f_{Y\larr X}\otimes_{D_X}M_1,\nbigj_1^{\bullet}
 \bigr)
\end{multline}
}
We set
$D^g_{\Ybar\larr \Xbar}:=
 F^{-1}(D_{\Ybar}\otimes\omega_{\Ybar}^{-1})
 \otimes_{F^{-1}\nbigo_{\Ybar}}
 \Omega_{\Xbar}^{\bullet}\otimes D_{\Xbar}[d_X]$
and
$D^g_{Y\larr X}:=
 F^{-1}(D_{Y}\otimes\omega_{Y}^{-1})
 \otimes_{F^{-1}\nbigo_Y}
 \Omega_X^{\bullet}\otimes D_{X}[d_X]$.
We have
$\bigl(
 D^g_{Y\larr X}\otimes_{\cnum}
 D^g_{\Ybar\larr\Xbar}\bigr)
 \otimes_{D_{X,\Xbar}}\distribution_X
\simeq
 F^{-1}(D_{Y,\Ybar}\otimes\omega_{Y,\Ybar}^{-1})
 \otimes_{F^{-1}\nbigo_{Y,\Ybar}}
 \distribution_X^{\bullet}[2d_X]$,
where 
we use the identification of
$\Omega_X^{\bullet}[d_X]
\otimes
 \Omega_{\Xbar}^{\bullet}[d_X]
 \otimes\distribution_X
\simeq
 \distribution^{\bullet}[2d_X]$.
In $D^b(F^{-1}D_{\Ybar})$,
the morphism $A'$ is represented by 
the composition $A''$ of the following morphisms:
{\scriptsize
\begin{multline}
 D^g_{\Ybar\larr \Xbar}\otimes_{D_{\Ybar}}\Mbar_2
\stackrel{B_0}\lrarr
 \nhom_{F^{-1}(D_Y)}
 \Bigl(
 D^g_{Y\larr X}\otimes_{D_X}M_1,\,
 F^{-1}\bigl(
 D_{Y,\Ybar}\otimes\omega_{Y,\Ybar}
 \bigr)\otimes_{F^{-1}\nbigo_{Y,\Ybar}}
 \distribution_X^{\bullet}
 \Bigr)[2d_X]\\
\lrarr
 \nhom_{F^{-1}(D_Y)}\bigl(
 D^g_{Y\larr X}\otimes_{D_X}M_1,\nbigj_1^{\bullet}
 \bigr)
\end{multline}
}

We take a $D_X\otimes F^{-1}(D_{\Ybar})$-injective
resolution $\nbigj_3^{\bullet}$ of
$D_{\Ybar\larr\Xbar}^g\otimes_{D_{\Xbar}}
 \distribution_X$.
We have a morphism
$D^g_{Y\larr X}\otimes_{D_X}\nbigj_3^{\bullet}
\lrarr
 \nbigj_1^{\bullet}$,
which represents the natural isomorphism
\[
D^g_{Y\larr X}\otimes_{D_X}
 \bigl(
 D^g_{\Ybar\larr\Xbar}\otimes_{D_{\Xbar}}
 \distribution
 \bigr)
\simeq
 F^{-1}(D_{Y,\Ybar}\otimes\omega_{Y,\Ybar}^{-1})
 \otimes_{F^{-1}\nbigo_{Y,\Ybar}}
 \distribution_{X}^{\bullet}[2d_X]. 
\]
The composition of 
the following morphisms is denoted by 
$B_1$:
\[
 D^g_{\Ybar\larr\Xbar}\otimes_{D_{\Xbar}}\Mbar_2
\lrarr
 \nhom_{D_X}\bigl(M_1,D^g_{\Ybar\larr\Xbar}\otimes
 \distribution_X\bigr)
\lrarr
 \nhom_{D_X}(M_1,\nbigj_3^{\bullet})
\]
Here, the first one is
induced by
$\varphi_{C,1,2}:\Mbar_2\lrarr
\nhom_{D_X}(M_1,\distribution_X)$.
The composite of the following morphisms
is denoted by $B_2$:
\begin{multline}
 \nhom_{D_X}(M_1,\nbigj_3^{\bullet})
\lrarr
 \nhom_{F^{-1}D_Y}\bigl(
 D^g_{Y\larr X}\otimes_{D_{X}}M_1,\,
 D^g_{Y\larr X}\otimes_{D_X}\nbigj_3^{\bullet}
 \bigr) \\
\lrarr
 \nhom_{F^{-1}D_Y}\bigl(
 D^g_{Y\larr X}\otimes_{D_{X}}M_1,\,
 \nbigj_1^{\bullet}
 \bigr)
\end{multline}
The composite $B_2\circ B_1$
factors through $B_0$.
Hence, $B_2\circ B_1=A''$
in $D^b(F^{-1}D_Y)$.
It means that the morphism $A$
is factorized into
{\scriptsize
\[
\begin{CD}
 F_{\dagger}\Mbar_2
 @>{F_{\dagger}\varphi_{C,1,2}}>>
 F_{\dagger}\nbigc_XM_1
 @>{B_3}>>
 F_!\nrhom_{F^{-1}(D_Y)}\Bigl(
 D_{Y\larr X}\otimes^L_{D_X}M_1,\,
 D_{Y\larr X,\Ybar\larr \Xbar}
 \otimes^L_{D_{X,\Xbar}}
 \distribution_X
 \Bigr)
\end{CD}
\]
}
We have 
$A_2\circ A_1\circ B_3=b$.
Thus, the proof of Lemma \ref{lem;10.9.10.5}
is finished.
\hfill\qed

\vspace{.1in}
We have the following commutative diagram
as the compatibility of the de Rham functor
and the push-forward:
\begin{equation}
\label{eq;10.9.10.6}
\begin{CD}
 F_{!}\DR_{\Xbar}\Mbar_2
@>{F_!\DR\varphi_{C,1,2}}>{\simeq}>
 F_{!}\DR_{\Xbar}\nbigc_XM_1\\
 @V{\simeq}VV @V{\simeq}VV \\
 \DR_{\Ybar}F_{\dagger}\Mbar_2
@>{\DR F_{\dagger}\varphi_{C,1,2}}>{\simeq}>
 \DR_{\Ybar}F_{\dagger}\nbigc_XM_1
\end{CD} 
\end{equation}

\begin{cor}
\label{cor;10.9.15.23}
We have $F^j_{\dagger}(M_1,M_2,C)\in
 \Dtriplecat^{nd}(Y)$.
\end{cor}
\pf
By Lemma \ref{lem;10.9.10.7},
Lemma \ref{lem;10.9.10.5}
and the commutative diagram (\ref{eq;10.9.10.6}),
we obtain that the induced morphism
$\DR \varphi_{F^{(0)}_{\dagger}C,1,2}:
 \DR_{\Ybar}F_{\dagger}\Mbar_2
\lrarr
 \DR_{\Ybar}\nbigc_YF_{\dagger}M_1$
is an isomorphism.
Because $F_{\dagger}^{(0)}C$ and
$F_{\dagger}C$ are equal up to the signature,
we obtain the corollary.
\hfill\qed

\subsection{Duality and push-forward}
\label{subsection;13.4.4.1}

Let $\nbigt=(M_1,M_2,C)\in\Dtriplecat^{nd}(X)$.
Let $F:X\lrarr Y$ be a morphism
such that its restriction to $\Supp M_i$
are proper.
Assume that 
there exists a non-degenerate hermitian pairing
$\DDD_X C:\DDD_X M_1\times \DDD_{\Xbar} \Mbar_2
 \lrarr \distribution_X$,
which is a dual of $C$.
We have the induced pairings:
\[
 F^{(0)}_{\dagger}C:
 F_{\dagger}M_1\times
 F_{\dagger}\Mbar_2
\lrarr
\distribution_Y
\]
\[
 F^{(0)}_{\dagger}\DDD_X C:
 F_{\dagger}\DDD_X M_1\times
 F_{\dagger}\DDD_{\Xbar} \Mbar_2
\lrarr
\distribution_Y
\]
Recall we have a natural equivalence
$F_{\dagger}\DDD M_i\simeq 
 \DDD F_{\dagger}M_i$,
under which we have the induced pairing
\[
  F^{(0)}_{\dagger}\DDD_X C:
 \DDD_Y F_{\dagger}M_1
 \times
 \DDD_{\Ybar}F_{\dagger}\Mbar_2
\lrarr
\distribution_Y.
\]

\begin{lem}
\label{lem;10.9.14.41}
The dual $\DDD_Y F^{(0)}_{\dagger}C$ exists,
and we have
$\DDD_Y F^{(0)}_{\dagger}C
=F^{(0)}_{\dagger}\DDD_X C$.
\end{lem}
\pf
We shall compare the following morphisms
\begin{equation}
\label{eq;10.9.10.3}
  \DR F^{(0)}_{\dagger}\DDD C:
 \DR_Y \DDD_YF_{\dagger}M_1
 \otimes 
 \DR_{\Ybar}
 \DDD_{\Ybar}F_{\dagger}\Mbar_2
\lrarr
 \omega_Y^{top}
\end{equation}
\begin{equation}
\label{eq;10.9.10.4}
 \DR 
 F^{(0)}_{\dagger}C\circ\exchange:
 \DR_{\Ybar}\bigl(
 F_{\dagger}\Mbar_2\bigr)
\otimes
 \DR_{Y}\bigl(
 F_{\dagger}M_1\bigr)
\lrarr 
\omega_Y^{top}
\end{equation}
We have an identification $\Upsilon_1$ of
$\DR_Y \DDD_YF_{\dagger}M_1
 \otimes 
 \DR_{\Ybar}
 \DDD_{\Ybar}F_{\dagger}\Mbar_2$
and 
$\DR_{\Ybar}\bigl(
 F_{\dagger}\Mbar_2\bigr)
\otimes
 \DR_{Y}\bigl(
 F_{\dagger}M_1\bigr)$
given by the following isomorphisms:
\[
 \nu_{F_{\dagger}M_1}\circ
 \DR(\varphi_{F^{(0)}_{\dagger}C,1,2}):
\DR_{\Ybar}F_{\dagger}\Mbar_2\simeq
 \DR_Y\DDD_YF_{\dagger}M_1
\]
\[ 
 \nu_{F_{\dagger}\Mbar_2}
 \circ\DR(\varphi_{F^{(0)}_{\dagger}C,2,1}):
\DR_{Y}F_{\dagger}M_1\simeq
 \DR_{\Ybar}\DDD_{\Ybar}F_{\dagger}\Mbar_2
\]
We have only to show that
(\ref{eq;10.9.10.3})
and (\ref{eq;10.9.10.4}) are equal
under the identification $\Upsilon_1$.

We have another identification 
\[
\Upsilon_2:
\DR_Y \DDD_YF_{\dagger}M_1
 \otimes 
 \DR_{\Ybar}
 \DDD_{\Ybar}F_{\dagger}\Mbar_2
\simeq
\DR_{\Ybar}\bigl(
 F_{\dagger}\Mbar_2\bigr)
\otimes
 \DR_{Y}\bigl(
 F_{\dagger}M_1\bigr).
\]
Indeed, we have the isomorphism
$\kappa_{1,2}:
 \DR_{\Ybar}F_{\dagger}\Mbar_2
 \simeq
 \DR_Y\DDD_YF_{\dagger}M_1$
induced by
$\nu_{M_1}\circ \DR(\varphi_{C,1,2})$,
given as follows:
\begin{equation}
\label{eq;10.9.10.2}
 \DR_{\Ybar}F_{\dagger}\Mbar_2
\simeq
 F_{!}\DR_{\Xbar}\Mbar_2
\simeq
 F_{!}\DR_{\Xbar}\nbigc_XM_1
\simeq
 F_{!}\DR_X\DDD_XM_1
\simeq
 \DR_Y\DDD_YF_{\dagger}M_1
\end{equation}
Similarly, we have the isomorphism
$\kappa_{2,1}:
\DR_{Y}F_{\dagger}M_1\simeq
 \DR_{\Ybar}\DDD_{\Ybar}F_{\dagger}\Mbar_2$.
The identification $\Upsilon_2$ is induced by
$\kappa_{1,2}$ and $\kappa_{2,1}$.
We can check that
(\ref{eq;10.9.10.3}) and
(\ref{eq;10.9.10.4}) are equal
under the identification $\Upsilon_2$.
Hence, we have only to show
$\nu_{F_{\dagger}M_1}\circ
 \DR\varphi_{F^{(0)}_{\dagger}C,1,2}
=\kappa_{1,2}$
and 
$\nu_{F_{\dagger}\Mbar_2}
 \circ
 \DR\varphi_{F^{(0)}_{\dagger}C,2,1}
=\kappa_{2,1}$.
As for the first,
we have only to show the commutativity
of the following diagrams:
\begin{equation}
 \label{eq;10.9.14.40}
 \begin{CD}
 F_{!}\DR_{\Xbar}\Mbar_2 
 @>{F_{!}\DR\varphi_C}>>
 F_{!}\DR_{\Xbar}\nbigc_XM_1 
 @>{F_{!}\nu_{M_1}}>>
 F_{!}\DR_{X}\DDD_XM_1\\ 
 @V{\simeq}VV @V{\simeq}VV @V{\simeq}VV \\
 \DR_{\Ybar}F_{\dagger}\Mbar_2 
 @>{\DR \varphi_{F^{(0)}_{\dagger}C}}>>
 \DR_{\Ybar}\nbigc_YF_{\dagger}M_1 
 @>{\nu_{F_{\dagger}M_1}}>>
 \DR_{\Xbar}\DDD_{\Xbar}F_{\dagger}M_1
 \end{CD}
\end{equation}
The commutativity of the left square 
follows from Lemma \ref{lem;10.9.10.5}
and the commutativity of (\ref{eq;10.9.10.6}).
The commutativity of the right square
is given in Lemma \ref{lem;10.9.10.7}.

We can prove the second one
by considering the non-degenerate $D$-triple
$C':\Mbar_2\times M_1\lrarr\distribution_X$
on the conjugate $\Xbar$ 
given by
$C'(\mbar_2,m_1)=C(m_1,\mbar_2)$.
We remark
$F_{\dagger}^{(0)}C'=(F_{\dagger}^{(0)}C)\circ\exchange$.
Note that the signature $(-1)^{d_X-d_Y}$ appears twice.
One is caused by the effect of the change of the orientation
on the trace map.
The other is the difference of the identifications
$\Tot\bigl(
\Omega_X^{\bullet}[d_X]\otimes\Omega_{\Xbar}[d_X]
\otimes\distribution_X\bigr)
\simeq
 \distribution_X^{\bullet}[2d_X]$
and 
$\Tot\bigl(
\Omega_Y^{\bullet}[d_Y]\otimes\Omega_{\Ybar}[d_Y]
\otimes\distribution_Y\bigr)
\simeq
 \distribution_Y^{\bullet}[2d_Y]$.
Then, we can apply the argument for the first.
Thus, we obtain Lemma \ref{lem;10.9.14.41}.
\hfill\qed

\subsection{Canonical prolongation of good meromorphic flat bundles}
\label{subsection;10.12.14.2}

Let $V_1$ be a good meromorphic flat bundle
on $(X,D)$.
We obtain a good meromorphic flat bundle
$\nbigc_X(V_1)(\ast D)$ on $(\Xbar,\Dbar)$,
and a good meromorphic flat bundle
$V_2:=\overline{\nbigc_X(V_1)(\ast D)}$
on $(X,D)$.
\begin{lem}
\label{lem;10.11.17.1}
The natural isomorphism
$\overline{\nbigc_X(V_1)_{|X\setminus D}}
\simeq
 V_{2|X\setminus D}$
is uniquely extended to isomorphisms
\[
 \overline{\nbigc_X(V_{1})}
 \simeq V_{2}[!D],
\quad\quad
\overline{\nbigc_X(V_1[!D])}\simeq
 V_{2}. 
\]
\end{lem}
\pf
We set 
$\nbigm:=\overline{\nbigc_X(V_{1})}$.
We have $\nbigm(\ast D)=V_2$.
Let $\nbigm'$ be any holonomic $D_X$-module
such that $\nbigm'(\ast D)=V_2$.
We have a unique morphism
$\varphi_1:
 \overline{\nbigc_X(\nbigm')}
\lrarr V_1$
such that $\varphi_{1|X\setminus D}=\id$.
Hence, we have 
a morphism $\varphi_2:\nbigm\lrarr\nbigm'$
such that $\varphi_{2|X\setminus D}=\id$.
It means $\nbigm=\nbigm[!D]=V_2[!D]$
by the universal property of $V_2[!D]$.
We obtain the other isomorphism similarly.
\hfill\qed

\vspace{.1in}
Hence, for integers $m$,
we obtain non-degenerate
holonomic hermitian pairings
\[
 \nbigt_!(V_1,V_2,m):=
 \bigl(
 V_1[\ast D],V_2[!D],(-1)^mC_1
 \bigr),
\]
\[
\nbigt_{\ast}(V_1,V_2,m):=
 \bigl(
 V_1[! D],V_2[\ast D],(-1)^mC_2
 \bigr).
\]
Here $C_i$ are non-degenerate hermitian pairings
corresponding to the isomorphisms in Lemma 
\ref{lem;10.11.17.1}.

Let us consider their dual.
By using the description of the Stokes structure,
we have
$V_2^{\lor}=
 \overline{\nbigc_X(V_1^{\lor})(\ast D)}$.
Hence, we have
$\nbigt_{\star}(V_1^{\lor},V_2^{\lor},m)$
for $\star=\ast,!$.
We have the canonical isomorphisms
$\DDD_X(V_i)\simeq
 V_i^{\lor}[!D]$
and $\DDD_X(V_i[!D])\simeq
 V_i^{\lor}$.

\begin{lem}
$\nbigt_{\star}(V_1,V_2,m)$
have duals,
and we have
\[
 \DDD\nbigt_{\ast}(V_1,V_2,m)
=\nbigt_{!}(V_1^{\lor},V_2^{\lor},m+d_X),
\quad
 \DDD\nbigt_{!}(V_1,V_2,m)
=\nbigt_{\ast}(V_1^{\lor},V_2^{\lor},m+d_X).
\]
\end{lem}
\pf
The case $D=\emptyset$
has already been proved
in Example \ref{example;13.4.1.11}.
Let $\pi:\Xtilde(D)\lrarr X$ be the real blow up.
Let $\nbigl_i$ be the local system on
$\Xtilde(D)$ with the Stokes structure
associated to $V_i$.
Let $\nbigl_i^{\lor}$
be the local system on
$\Xtilde(D)$ with the Stokes structure
associated to $V_i^{\lor}$.

We have the associated constructible sheaves
$\nbigl_i^{<D}$ and $\nbigl_i^{\leq D}$,
from $\nbigl_i$.
(See \cite{mochi9} for the notation.)
Similarly, we have
$\nbigl_i^{\lor\,<D}$ and $\nbigl_i^{\lor\,\leq D}$,
from $\nbigl_i^{\lor}$.
Let $k:X\setminus D\lrarr\Xtilde(D)$
be the inclusion.
We set $\omega_{\Xtilde}(D):=
 k_{!}\cnum_{X\setminus D}[2d_X]$,
which is a dualizing complex.
We set
$\DDD_{\Xtilde(D)}(\nbigf)
:=\nrhom_{\cnum_{\Xtilde(D)}}(\nbigf,\omega_{\Xtilde(D)})$.
We have natural identifications
$\DDD_{\Xtilde(D)}\nbigl_i^{<D}[d_X]
\simeq
 \nbigl_i^{\lor\,\leq D}[d_X]$
and
$\DDD_{\Xtilde(D)}\nbigl_i^{\leq D}[d_X]
\simeq
 \nbigl_i^{\lor\,< D}[d_X]$.
We have naturally defined non-degenerate pairings
$\nbigl_{1}^{<D}[d_X]\times
 \nbiglbar_2^{\leq D}[d_X]
 \lrarr\omega_{\Xtilde(D)}$
and
$\nbigl_{1}^{\leq D}[d_X]\times
 \nbiglbar_2^{< D}[d_X]\lrarr\omega_{\Xtilde(D)}$.
We also have
$\nbigl_{1}^{\lor <D}[d_X]\times
 \nbiglbar_2^{\lor \leq D}[d_X]\lrarr\omega_{\Xtilde(D)}$
and
$\nbigl_{1}^{\lor \leq D}[d_X]\times
 \nbiglbar_2^{\lor < D}[d_X]
 \lrarr\omega_{\Xtilde(D)}$.

\vspace{.1in}

Look at 
$\DR_XV_1[!D]\lrarr
 \DDD \DR_{\Xbar}\Vbar_2$
and 
$\DR_{\Xbar}\Vbar_2\lrarr
 \DDD \DR_{X}V_1[!D]$
induced by the pairing
$V_1[!D]\times \Vbar_2
\lrarr\distribution_X$.
The morphism
$\DR_XV_1[!D]\lrarr
 \DDD \DR_{\Xbar}\Vbar_2$
is factorized as follows:
\[
\DR_X V_1[!D]
\lrarr
\nhom_{\cnum_X}\bigl(
 \DR_{\Xbar}\Vbar_2,
 \distribution_X^{\bullet}
 \bigr)[2d_X]
\lrarr
 \DDD \DR_{\Xbar}\Vbar_2
\]
We have the following commutative diagram:
\[
 \begin{CD}
 \DR_X^{<D}V_1
 @>>>
 \nhom_{\cnum_X}\bigl(
 \DR^{\leq D}_{\Xbar}\Vbar_2,
 \Omega^{\bullet\,<D}_X
 \bigr)[2d_X]\\
 @VVV @VVV \\
 \DR_XV_1[!D]
 @>>>
 \nhom_{\cnum_X}\bigl(
 \DR^{\leq D}_{\Xbar}\Vbar_2,
 \distribution^{\bullet}_X
 \bigr)[2d_X]
 \end{CD}
\]
Hence, the isomorphism
$\DR_XV_1[!D]\lrarr
 \DDD \DR_{\Xbar}\Vbar_2$
is obtained as 
the push-forward of the unique isomorphism
$\nbigl_1^{<D}[d_X]\simeq
 \DDD_{\Xtilde(D)}
 \nbiglbar_2^{\leq D}[d_X]$
whose restriction to $X\setminus D$ is the identity.

The morphism
$\DR_{\Xbar}\Vbar_2\lrarr
 \DDD \DR_{\Xbar}V_1[!D]$
is represented as the composition
of the following morphisms
\begin{multline}
 \DR_{\Xbar}\Vbar_2
\lrarr
 \nhom_{\cnum_X}\bigl(
 \DR^{<D}_XV_1,\nbige_X^{\bullet}
 \bigr)[2d_X]
 \\
\lrarr
 \nhom_{\cnum_X}\bigl(
 \DR^{<D}_XV_1,\distribution_X^{\bullet}
 \bigr)[2d_X]
\lrarr
 \DDD \DR^{<D}_{X}V_1
\end{multline}
Hence, the isomorphism
$\DR_{\Xbar}\Vbar_2\lrarr
 \DDD \DR_XV_1[!D]$
is obtained as the push-forward of
the unique isomorphism
$\nbiglbar_2^{\leq D}[d_X]
\simeq
 \DDD_{\Xtilde(D)}\nbigl_1^{<D}[d_X]
=\nbigl_1^{\lor\leq D}[d_X]$
whose restriction to $X\setminus D$ is the identity.

Let us consider the diagram:
\begin{equation}
 \label{eq;10.9.15.10}
 \begin{CD}
 \DR_X(V_1^{\lor})
 \times
 \DR_{\Xbar}(\Vbar_2[!D])
 @>{A_2}>>
 \omega_X^{top}\\
 @V{A_1}VV @V{=}VV \\
 \DR_{\Xbar}(\Vbar_2)
\times
 \DR_X(V_1[!D])
 @>{(-1)^{d_X}A_3}>>
 \omega_X^{top}
 \end{CD}
\end{equation}
The identification $A_1$
and the pairings $A_i$ $(i=2,3)$
are obtained as the push forward of
the identification and the pairings
on the real blow up.
Hence, the diagram (\ref{eq;10.9.15.10})
is commutative.
\hfill\qed

\vspace{.1in}

Let $g$ be a holomorphic function  on $X$
such that $g^{-1}(0)=D$.
Combined with Beilinson construction,
we obtain the following 
non-degenerate hermitian pairings,
which are mutually dual:
\[
 \Pi^{a,b}_{g!}C_1:
 \Pi^{-b+1,-a+1}_{g\ast}V_1
 \times
 \overline{
 \Pi^{a,b}_{g!}V_2}
\lrarr
\distribution_X
\]
\[
(-1)^{d_X}
\Pi^{-b+1,-a+1}_{g\ast}C_2:
 \Pi^{a,b}_{g!}V_1^{\lor}\times
 \overline{
 \Pi^{-b+1,-a+1}_{g\ast}V_2^{\lor}}
\lrarr
\distribution_X
\]
We obtain the following non-degenerate
hermitian pairings, which are mutually dual:
\[
 \Xi_{g}^{(a)}C_1:
 \Xi^{(-a)}_{g}V_1
 \times 
 \overline{
 \Xi^{(a)}_{g}V_2}
\lrarr
\distribution_X
\]
\[
 (-1)^{d_X}
\Xi_{g}^{(-a)}C_2:
 \Xi^{(a)}_{g}V_1^{\lor}\times 
 \overline{
 \Xi^{(-a)}_{g}V_2^{\lor}}
\lrarr
\distribution_X
\]
We also obtain the following hermitian pairings
which are mutually dual:
\[
 \psi^{(a)}_gC_1:
 \psi^{(-a+1)}_{g}V_1
 \times 
 \overline{\psi^{(a)}_{g}V_2}
\lrarr
\distribution_X
\]
\[
 (-1)^{d_X}
 \psi^{(-a+1)}_gC_2:
 \psi^{(a)}_{g}V_1^{\lor}\times 
 \overline{
 \psi^{(-a+1)}_{g}V_2^{\lor}}
\lrarr
\distribution_X
\]
Namely,
there exist 
$\DDD\psi_g^{(a)}\nbigt_{\star}(V_1,V_2)$
and 
$\DDD\Xi_g^{(a)}\nbigt_{\star}(V_1,V_2)$.

\subsection{Special case}
\label{subsection;10.9.15.21}

Let $\nbigt=(M_1,M_2,C)\in
\Dtriplecat^{nd}(X)$.
Assume that there exists a good cell 
$\nbigc=(Z,U,\varphi,V_1)$ for $M_1$
at $P\in\Supp\nbigt$.
Let us show that,
on a small neighbourhood $X_P$ of $P$,
there exists a dual pairing
$\DDD C:\DDD M_{1|X_P}
 \times
 \DDD M_{2|X_P}
 \lrarr\distribution_{X_P}$.

Let $V_2:=\overline{\nbigc_Z(V_1)}(\ast D_Z)$.
Because $M_2\lrarr \overline{\nbigc_X(M_1)}$
is an isomorphism,
the tuple
$(Z,U,\varphi,V_2)$ is a good cell of $M_2$.
We have the objects
$\nbigt_{V!}:=
 (V_{1\ast},V_{2!},C_{1!})
 \in\Dtriplecat^{nd}(Z)$
and 
$\nbigt_{V\ast}:=
 (V_{1!},V_{2\ast},C_{1\ast})
 \in\Dtriplecat^{nd}(Z)$.
We have 
$\varphi_{\dagger}V_{1!}
 \lrarr
 \nbigm_1\lrarr
 \varphi_{\dagger}V_{1\ast}$.
By applying the functor 
$\overline{\nbigc_X}$,
we obtain 
$\varphi_{\dagger}V_{2!}
\lrarr\nbigm_2\lrarr
 \varphi_{\dagger}V_{2\ast}$.
By construction of the pairings,
we obtain the morphisms
$\varphi_{\dagger}\nbigt_{V!}
\lrarr\nbigt\lrarr
 \varphi_{\dagger}\nbigt_{V\ast}$.
Let $g$ be a cell function for $(Z,U,\varphi,V_1)$,
and we put $g_Z:=g\circ\varphi$.
Recall that we obtain 
$\phi^{(0)}_g\nbigt\in \Dtriplecat^{nd}(X)$
as the cohomology of the complex
$\varphi_{\dagger}\nbigt_{V!}
 \lrarr
 \varphi_{\dagger}\Xi^{(0)}_{g_Z}\nbigt_{V}
 \oplus\nbigt
 \lrarr
 \varphi_{\dagger}\nbigt_{V\ast}$.

Suppose that the claim of 
Theorem {\rm\ref{thm;10.9.14.2}}
holds for non-degenerate pairings
such that the dimensions of the support
are strictly less than
$\dim\Supp\nbigt$.
Assume also that
$\dim(g^{-1}(0)\cap\Supp(\nbigt))<
 \dim\Supp(\nbigt)$.
Then, we have the dual
$\DDD \phi^{(0)}_g\nbigt$
and 
$\DDD \psi^{(a)}_g\nbigt$
in $\Dtriplecat^{nd}(X)$
of 
$\phi^{(0)}_g\nbigt$
and 
$\psi^{(a)}_g\nbigt$,
respectively.
We have the naturally induced morphisms
$\DDD\psi^{(0)}_g\nbigt
\lrarr
 \DDD\phi^{(0)}_g\nbigt
\lrarr
 \DDD\psi^{(1)}_g\nbigt$.
We also have the duals
$\DDD\varphi_{\dagger}\psi^{(a)}_{g_Z}\nbigt$
and $\DDD\varphi_{\dagger}\Xi^{(a)}_{g_Z}\nbigt$
(\S\ref{subsection;10.12.14.2}).
We have natural isomorphisms
$\DDD\psi^{(a)}_g\nbigt\simeq
 \DDD\varphi_{\dagger}\psi^{(a)}_{g_Z}\nbigt$
(\S\ref{subsection;13.4.4.1}).
We define a non-degenerate hermitian pairing
$\DDD\nbigt$
as the cohomology of the complex:
\[
 \DDD\varphi_{\dagger}\psi^{(0)}_{g_Z}\nbigt_V
 \lrarr
 \DDD\varphi_{\dagger}\Xi^{(0)}_{g_Z}\nbigt_{V}
 \oplus\DDD \phi^{(0)}_g\nbigt
 \lrarr
 \DDD\varphi_{\dagger}\psi^{(1)}_{g_Z}\nbigt_V
\]
It is easy to check
$\DDD\nbigt$ gives a dual of $\nbigt$.

\begin{cor}
\label{cor;10.9.15.20}
Suppose that the claim of 
Theorem {\rm\ref{thm;10.9.14.2}}
holds for any non-degenerate pairings
such that the dimensions of the support
are less than $n$.

Let $\nbigt\in\Dtriplecat^{nd}(X)$
with $\dim\Supp\nbigt=n$.
If there exists a good cell at any point 
of $P\in\nbigt$,
there exists a dual of $\nbigt$.
\hfill\qed
\end{cor}

\subsection{End of the proof}
\label{subsection;10.12.27.31}

Let $V_1$ be a meromorphic flat bundle
on $(X,D)$, which is not necessarily good.
We have a non-degenerate pairing
$\nbigt=\bigl(V_1,\overline{\nbigc_X(V_1)},C\bigr)$,
where $C$ is the naturally defined pairing.
By applying Corollary \ref{cor;10.9.15.20}
and resolution of turning points 
(\cite{kedlaya}, \cite{kedlaya2},
 see also \cite{mochi6} and \cite{mochi7}),
we obtain the existence of a dual of $\nbigt$.
Similarly, we have a dual of
$(V_{1!},\overline{\nbigc_X(V_{1!})},C')$.
Then, by using the argument in 
\S\ref{subsection;10.9.15.21}
for cells which are not necessarily good,
we can finish the proof of
Theorem \ref{thm;10.9.14.2}.
Once we obtain Theorem \ref{thm;10.9.14.2},
we obtain Theorem \ref{thm;11.3.30.1}
from Corollary \ref{cor;10.9.15.23}
and Lemma \ref{lem;10.9.14.41}.
\hfill\qed

\section{Real structure}

\subsection{Real structure of non-degenerate $D$-triple}
\label{subsection;11.3.29.1}

We define an auto equivalence
$\gammatilde^{\ast}:=\DDD^{\herm}\circ\DDD
=\DDD\circ\DDD^{\herm}$
on $\nbigc\bigl(\Dtriplecat^{nd}(X)\bigr)$.
\index{functor $\gammatilde^{\ast}$}
Similarly,
we define 
$\gammatilde^{(0)\ast}:=
 \DDD^{(0)\herm}\circ\DDD^{(0)}
=\DDD^{(0)}\circ\DDD^{(0)\herm}$
on $\Dcomplextriplecat_1^{nd}(X)$.
\index{functor $\gammatilde^{(0)\ast}$}
We have 
$\Psi_1\circ\gammatilde^{(0)\ast}
=\gammatilde^{\ast}\circ\Psi_1$.
In the following,
we will often omit to denote ``$(0)$''.

We obtain the following proposition
from Lemma \ref{lem;11.3.23.20}
and Theorem \ref{thm;10.9.14.2}.
\begin{prop}
Let $f:X\lrarr Y$ be a morphism of complex manifolds.
Let $\nbigt^{\bullet}\in
 \nbigc\bigl(\Dtriplecat^{nd}(X)\bigr)$
such that the restriction of $f$ to
the support of $\nbigt^{\bullet}$ is proper.
Then,
we have the following natural isomorphism
in $\nbigc(\Dtriplecat^{nd}(Y))$:
\[
H^j\bigl(
 F_{\dagger}\bigl(
 \gammatilde^{\ast}\nbigt^{\bullet}
 \bigr)
 \bigr)[-j]
\simeq
 \gammatilde^{\ast}\Bigl(
 H^{j}F_{\dagger} \bigl(\nbigt^{\bullet}\bigr)[-j]
 \Bigr)
\]
Correspondingly,
for $\gbigt\in\Dcomplextriplecat^{nd}(X)$,
we have the following natural isomorphism
in $\Dcomplextriplecat^{nd}(Y)$:
\[
 \vecH^j\bigl(
 F^{(0)}_{\dagger}(\gammatilde^{\ast}\gbigt)
 \bigr)
\simeq
 \gammatilde^{\ast}
 \vecH^{j}F^{(0)}_{\dagger}(\gbigt)
\]
\hfill\qed
\end{prop}

A real structure of $\nbigt^{\bullet}
\in\nbigc\big(\Dtriplecat^{nd}(X))$
is defined to be an isomorphism
$\kappa:\gammatilde^{\ast}\nbigt^{\bullet}\simeq
 \nbigt^{\bullet}$
such that
$\kappa\circ\gammatilde^{\ast}\kappa=\id$.
Similarly,
a real structure of $\gbigt\in\Dcomplextriplecat^{nd}(X)$
is defined to be an isomorphism
$\kappa:\gammatilde^{\ast}\gbigt
\simeq\gbigt$
such that
$\kappa\circ\gammatilde^{\ast}\kappa=\id$.
A real structure of
$\nbigt^{\bullet}\in\nbigc\bigl(\Dtriplecat^{nd}(X)\bigr)$
naturally corresponds to
a real structure of
$\Psi_1\bigl(\nbigt^{\bullet}\bigr)
 \in\Dcomplextriplecat^{nd}(X)$.

\vspace{.1in}
For our later purpose,
we give the compatibility with the shift functor.
\begin{lem}
\label{lem;11.3.24.3}
We have
$\kappa_{1,\ell}:
\nbigs_{-\ell}\circ\DDD^{\herm}
\simeq
 \DDD^{\herm}\circ\nbigs_{\ell}$
and
$\kappa_{2,\ell}:
 \nbigs_{-\ell}\circ\DDD
\simeq\DDD\circ\nbigs_{\ell}$
on $\nbigc\bigl(\Dtriplecat^{nd}(X)\bigr)$
and $\Dcomplextriplecat_1^{nd}(X)$.
As a consequence,
we have 
$\kappa_{3,\ell}:
 \nbigs_{-\ell}\circ\gammatilde^{\ast}
\simeq
 \gammatilde^{\ast}\circ\nbigs_{\ell}$
on the categories.
They satisfy $\kappa_{a,\ell+m}=\kappa_{a,\ell}\circ\kappa_{a,m}$.
\end{lem}
\pf
We have already observed
the compatibility of the shift and $\DDD^{\herm}$
in \S\ref{subsection;13.8.5.1}.
We have 
\[
 \DDD\nbigs_{\ell}(\gbigt)
=\bigl(
 \DDD(M_1^{\bullet}[-\ell]),
 \DDD(M_2^{\bullet}[\ell]),
 \DDD(C[\ell])
 \bigr),
\]
\[
 \nbigs_{-\ell}\DDD(\gbigt)
=\bigl(
 \DDD(M_1^{\bullet})[\ell],
 \DDD(M_2^{\bullet})[-\ell],
 \DDD(C)[-\ell]
 \bigr).
\]
We have natural identifications
\[
 \DDD(M_1^{\bullet}[-\ell])^{-p}
=\DDD(M_1^{\bullet}[-\ell]^{p})
=\DDD(M_1^{p-\ell}),
\]
\[
 \DDD(M_1^{\bullet})[\ell]^{-p}
=\DDD(M_1^{\bullet})^{-p+\ell}
=\DDD(M_1^{p-\ell}).
\]
We also have
$\DDD(M_2^{\bullet}[\ell])^{p}
=\DDD(M_2^{-p+\ell})$
and
$\DDD(M_2^{\bullet})[-\ell]^{p}
=\DDD(M_2^{-p+\ell})$.
Under the natural identifications,
we have
$\DDD(C[\ell])^p
=(-1)^{\ell}
   \bigl(\DDD(C)[-\ell]\bigr)^p$.
Indeed, we have the following equalities:
\begin{equation}
 \DDD(C[\ell])^p
=\DDD\bigl( (C[\ell])^{-p} \bigr)
=(-1)^{\ell(\ell-1)/2+\ell p}\DDD(C^{-p+\ell})
\end{equation}
\begin{equation}
  \bigl(\DDD(C)[-\ell]\bigr)^p
=(-1)^{\ell(\ell+1)/2+\ell p}(\DDD C)^{p-\ell}
=(-1)^{\ell(\ell+1)/2+\ell p}\DDD(C^{-p+\ell})
\end{equation}
The isomorphism $\kappa_{2,\ell}$
is given by the following morphisms:
\[
\begin{CD}
 \DDD(M_1^{\bullet})[\ell]^{-p}
 @<{(-1)^{\ell p}\epsilon(-\ell)}<<
 \DDD(M_1^{\bullet}[-\ell])^{-p},
\end{CD}
\quad\quad
\begin{CD}
 \DDD(M_2^{\bullet})[-\ell]^p
 @>{(-1)^{\ell p}\epsilon(\ell)}>>
 \DDD(M_2^{\bullet}[\ell])^p
\end{CD}
\]
(See (1.1.5.2) of \cite{deligne-SGA4}.)
It is easy to check that they
give the desired isomorphism.
\hfill\qed

\begin{cor}
A real structure of $\gbigt\in\Dcomplextriplecat(X)$
naturally induces
a real structure of
$\nbigs_{\ell}(\gbigt)$.
A similar claim holds
for real structures and shifts
of objects in $\nbigc\bigl(\Dtriplecat^{nd}(X)\bigr)$.
\hfill\qed
\end{cor}

\subsubsection{Functoriality of the real structure}

Let $\Dtriplecat(X,\real)$
be the category of non-degenerate $D$-triples
with real structures.
\index{category $\Dtriplecat(X,\real)$}
We have the following functoriality.
\begin{prop}
\label{prop;11.3.24.10}
\mbox{{}}
\begin{itemize}
\item
$\Dtriplecat(X,\real)$
is equipped with 
the functors $\DDD$ and $\DDD^{\herm}$.
\item
For any hypersurface $H\subset X$,
we have the exact functors
$[\ast H]$ and $[!H]$ 
on $\Dtriplecat(X,\real)$,
which are compatible with
the corresponding functors
on $\Dtriplecat(X)$.
\item
For a proper morphism $f:X\lrarr Y$,
we have the cohomological functor
$H^jf_{\dagger}:
 \Dtriplecat(X,\real)
\lrarr
 \Dtriplecat(Y,\real)$
compatible with
$H^jf_{\dagger}:
 \Dtriplecat(X)\lrarr
 \Dtriplecat(Y)$.
\hfill\qed
\end{itemize}
\end{prop}

Let $H$ be any hypersurface of $X$.
Let $\nbigt_1=(V_1,V_2,C)$ be a non-degenerate
smooth $D_{X(\ast H)}$-triple.
We have its dual
$\nbigt_1^{\lor}=(V_1^{\lor},V_2^{\lor},C^{\lor})$
as a smooth $D_{X(\ast H)}$-triple.
We set
$\gammatilde_{sm}^{\ast}(\nbigt_1):=
 (\nbigt_1^{\lor})^{\ast}$.
A real structure of $\nbigt_1$
as a smooth $D_{X(\ast H)}$-triple
is an isomorphism
$\kappa:\gammatilde_{sm}^{\ast}(\nbigt_1)
\simeq\nbigt_1$
such that
$\kappa\circ\gammatilde_{sm}^{\ast}(\kappa)=\id$.

\begin{lem}
\label{lem;13.8.30.10}
Let $\nbigt\in\Dtriplecat^{nd}(X,\real)$.
Let $\nbigt_1$ be a smooth $D_{X}$-triple
with real structure.
Then,
$\nbigt\otimes\nbigt_1$
has a naturally induced real structure.
\end{lem}
\pf
We naturally have
$\gammatilde^{\ast}(\nbigt\otimes\nbigt_1)
\simeq
 \gammatilde^{\ast}(\nbigt)
\otimes
 \gammatilde_{sm}^{\ast}(\nbigt_1)$.
Hence, the claim follows.
\hfill\qed

\begin{prop}
\label{prop;13.8.30.11}
Let $\nbigt\in\Dtriplecat^{nd}(X,\real)$.
Let $\nbigt_1$ be a smooth $D_{X(\ast H)}$-triple
with real structure.
Then,
$(\nbigt\otimes\nbigt_1)[\star H]
\in\Dtriplecat^{nd}(X)$
has a unique $\real$-structure
such that 
its restriction to $X\setminus H$
is the natural one as in Lemma 
{\rm\ref{lem;13.8.30.10}}.
\end{prop}
\pf
As $D_{X(\ast H)}$-triple,
we have a natural isomorphism
$\gammatilde^{\ast}(\nbigt\otimes\nbigt_1)
\simeq
 \gammatilde^{\ast}(\nbigt)
\otimes\gammatilde_{sm}^{\ast}(\nbigt_1)$.
Hence, the claim of the proposition is easy to see.
\hfill\qed

%\vspace{.1in}
%By using Proposition \ref{prop;13.8.30.11},
\begin{cor}
For any holomorphic function $g$ on $X$,
we have Beilinson's functors
$\Pi^{a,b}_{g,\ast!}$
on $\Dtriplecat(X,\real)$,
compatible with
those on 
$\Dtriplecat(X)$.
In particular,
we have
the nearby cycle functor $\psi^{(a)}_g$ 
and the maximal functor $\Xi^{(a)}_g$.
Moreover,
we have the vanishing cycle functor
$\phi^{(a)}_g$.
\hfill\qed
\end{cor}

\subsection{Descriptions of real perverse sheaves}
\label{subsection;13.8.29.20}

\subsubsection{Preliminary}

Let $X$ be a $d_X$-dimensional complex manifold.
Let $D_c^b(\cnum_X)$ denote the derived
category of cohomologically bounded constructible
$\cnum_X$-complexes.
We set $a_X^!\cnum:=\cnum_X[2d_X]$ in
$D_c^b(\cnum_X)$
with a naturally defined real structure.
\index{complex $a_X^{\bikkuri}\cnum$}
For any $\nbigc\in D^b_c(\cnum_X)$,
we set $\DDD(\nbigf):=
\nrhom_{\cnum_X}(\nbigf,a_X^!\cnum)$.
The real structure of $a_X^!\cnum$ induces
an identification
$\DDD(\overline{\nbigf})\simeq
 \overline{\DDD\nbigf}$.

\vspace{.1in}

We have a quasi-isomorphism
$a_X^{!}\distribution_X^{\bullet}[2d_X]$
given by 
$s\longmapsto (2\pi\sqrt{-1})^{d_X}\iota(s)$
at the degree $-d_X$,
where $\iota:\cnum_X\lrarr \distribution_X$
is the inclusion.
We shall identify them in $D^b_c(\cnum_X)$.
We have the real compatible real structure
on $\distribution_X^{\bullet}[2d_X]$,
which is twisted by $(2\pi\sqrt{-1})^{d_X}$
from the natural one.
\begin{rem}
It is also induced by the identification
$\Omega_X^{\bullet}[d_X]
 \otimes
 \Omega_{\Xbar}^{\bullet}[d_X]
 \otimes\distribution_X^{\bullet}
\simeq
 \distribution_X^{\bullet}[2d_X]$.
Namely, the composite of the following
identifications
\begin{multline}
 \overline{\distribution_X^{\bullet}[2d_X]}
\simeq
 \overline{\Omega_X^{\bullet}[d_X]
 \otimes \Omega_{\Xbar}^{\bullet}[d_X]
 \otimes\distribution_X}
\simeq
 \Omega^{\bullet}_{\Xbar}[d_X]\otimes
 \Omega^{\bullet}_X[d_X]\otimes\distribution_X
 \\
\simeq
 \Omega^{\bullet}_X[d_X]\otimes
 \Omega^{\bullet}_{\Xbar}[d_X]\otimes\distribution_X
\simeq
 \distribution_X^{\bullet}[2d_X]
\end{multline}
$\omega$ is mapped to
$(-1)^{d_X}\overline{\omega}$.
Here, 
\hfill\qed
\end{rem}

Let $F:X\lrarr Y$ be a proper morphism.
We use the trace morphism
$F_{\ast}(a_X^{!}\cnum)\lrarr
 a_Y^!\cnum$,
induced by the identification
$a_X^!\cnum\simeq
 \distribution_X^{\bullet}[2d_X]$.
The trace morphism is compatible
with the real structures of
$a_X^!\cnum$ and $a_Y^!\cnum$.
We have the Verdier duality
$F_{\ast}\DDD(\nbigf)\simeq
 \DDD F_{\ast}\nbigf$
for any $\nbigf$ in $D_c^b(\cnum_Y)$,
and the following natural diagram is commutative:
\[
 \begin{CD}
 F_{\ast}\DDD\nbigfbar
 @>{\simeq}>>
 \DDD F_{\ast}\nbigfbar\\
 @V{\simeq}VV @V{\simeq}VV \\
 \overline{F_{\ast}\DDD\nbigf}
 @>{\simeq}>>
 \overline{\DDD F_{\ast}\nbigf}
 \end{CD}
\]

\vspace{.1in}
Let $\nbigf_i$  $(i=1,2)$ 
be objects in $D^b_c(\cnum_X)$.
A morphism
$C:\nbigf_1\otimes\nbigf_2\lrarr a_X^!\cnum$
is equivalent to a morphism
$\Psi_C:\nbigf_1\lrarr\DDD\nbigf_2$.
We have the induced pairings:
\[
 \overline{C}:
 \overline{\nbigf}_1\otimes
 \overline{\nbigf}_2\lrarr
 a_X^!\cnum,
\quad
 C\circ\exchange:
 \nbigf_2\otimes\nbigf_1\lrarr a_X^!\cnum
\]
We have
$\Psi_{\overline{C}}=\overline{\Psi_C}$
under the natural isomorphism
$\DDD\nbigfbar_2\simeq \overline{\DDD\nbigf_2}$.
We have 
$\Psi_{C\circ\exchange}
=\DDD \Psi_C$
under the natural isomorphism
$\nbigf_2\simeq \DDD\bigl(\DDD\nbigf_2\bigr)$.
For any integer $\ell$,
the shift $C[\ell]:\nbigf_1[-\ell]\otimes\nbigf_2[\ell]\lrarr
a_X^!\cnum$ is induced.
(See \S\ref{subsection;11.4.6.20}.)
We have $\Psi_{C[\ell]}=\Psi$
under the natural identification
$\DDD(\nbigf_2[\ell])\simeq
 \DDD(\nbigf_2)[-\ell]$,
as remarked in Lemma \ref{lem;13.4.5.30}.

We say that $C$ is non-degenerate,
if $\Psi_C$ is a quasi-isomorphism.
In that case,
we have the induced non-degenerate pairing
$\DDD' C:\DDD\nbigf_2\otimes\DDD\nbigf_1\lrarr
 a_X^!\cnum$ in $D_c^b(\cnum_X)$
given by 
$\DDD\nbigf_2\otimes\DDD\nbigf_1\lrarr
 \nbigf_1\otimes\nbigf_2\lrarr a_X^!\cnum$.
By construction,
we have $\Psi_{\DDD' C}=\Psi_C^{-1}$.

\subsubsection{Category of non-degenerate pairings}

Let $\cnumcomplextriplecat(X)$
denote the category of tuples of
$\nbigf_i\in D_b^c(\cnum_X)$ $(i=1,2)$
with a pairing
$C:\nbigf_1\otimes\nbigfbar_2
\rarr
 a_X^!\cnum$
in $D^b_c(\cnum_X)$.
A morphism
$(\nbigf_1,\nbigf_2,C)
\rarr
 (\nbigf_1',\nbigf_2',C')$
is a pair of morphisms
$\varphi_1:\nbigf_1'\rarr\nbigf_1$
and 
$\varphi_2:\nbigf_2\rarr\nbigf_2'$
such that
$C'\circ(\id\otimes\varphi_2)
=C\circ(\varphi_1\otimes\id)$.
Let $\cnumcomplextriplecat^{nd}(X)
\subset\cnumcomplextriplecat(X)$ be the full subcategory of
the objects
$(\nbigf_1,\nbigf_2,C)$
such that $C$ is non-degenerate.

We define contravariant auto equivalences
$\DDD^{\herm}$ and $\DDD$
on $\cnumcomplextriplecat^{nd}(X)$
as follows:
\[
 \DDD^{\herm}(\nbigf_1,\nbigf_2,C):=
 (\nbigf_2,\nbigf_1,\overline{C\circ\exchange}),
\quad
 \DDD(\nbigf_1,\nbigf_2,C):=
 (\DDD\nbigf_1,\DDD\nbigf_2,
 \DDD' C\circ\exchange)
\]
We define an auto equivalence $\gammatilde^{\ast}$
on $\cnumcomplextriplecat^{nd}(X)$
by 
$\gammatilde^{\ast}:=\DDD\circ\DDD^{\herm}
=\DDD^{\herm}\circ\DDD$.
We have
\[
 \gammatilde^{\ast}\bigl(
 \nbigf_1,\nbigf_2,C
 \bigr)
=(\DDD\nbigf_2,\DDD\nbigf_1,\overline{\DDD'C}).
\]
For any integer $\ell$,
the shift is defined by 
$\nbigs_{\ell}(\gbigf)=(\nbigf_1[-\ell],\nbigf_2[\ell],C[\ell])$.

\begin{lem}
\label{lem;11.3.24.2}
We have
$\nu_{1,\ell}:
 \nbigs_{-\ell}\circ\DDD^{\herm}
\simeq \DDD^{\herm}\circ\nbigs_{\ell}$
and
$\nu_{2,\ell}:
 \nbigs_{-\ell}\circ\DDD
\simeq \DDD\circ\nbigs_{\ell}$.
As a consequence,
we have
$\nu_{3,\ell}:\nbigs_{-\ell}\circ\gammatilde^{\ast}
\simeq\gammatilde^{\ast}\circ\nbigs_{\ell}$.
They satisfy
$\nu_{a,\ell+m}=\nu_{a,\ell}\circ\nu_{a,m}$.
\end{lem}
\pf
We have
\[
\DDD^{\herm}\nbigs_{\ell}(\gbigf)
=(\nbigf_2[\ell],\nbigf_1[-\ell],
 \overline{C[\ell]\circ\exchange}), 
\]
\[
 \nbigs_{-\ell}\DDD^{\herm}(\gbigf)
=(\nbigf_2[\ell],\nbigf_1[-\ell],
 \overline{C\circ\exchange}[-\ell]).
\]
The isomorphism $\nu_{1,\ell}$
is given by the natural identifications
of the underlying $\cnum_X$-complexes.
We have 
$\Psi_{\overline{C[\ell]\circ\exchange}}
=\overline{\DDD\Psi_{C[\ell]}}
=\overline{\DDD\Psi_C}
=\Psi_{\overline{C\circ\exchange}}$.
Hence, the pairings are equal.

Let us consider the case of $\DDD$.
We have the following:
\[
 \nbigs_{-\ell}\DDD(\gbigf)
=\bigl((\DDD\nbigf_1)[\ell],(\DDD\nbigf_2)[-\ell],
 (\DDD' C\circ\exchange)[-\ell] \bigr)
\]
\[
 \DDD\nbigs_{\ell}(\gbigf)
=\bigl(
 \DDD(\nbigf_1[-\ell]),
 \DDD(\nbigf_2[\ell]),
 \DDD'\bigl(C[\ell]\bigr)\circ\exchange
 \bigr)
\]
We have natural isomorphisms
$\DDD(\nbigf_1)[\ell]
\simeq
 \DDD(\nbigf_1[-\ell])$
and 
$\DDD(\nbigf_2)[-\ell]
\simeq
 \DDD(\nbigf_2[\ell])$,
as in \S\ref{subsection;11.4.6.20}.
Under the identification,
we obtain the equality of the pairings
from
$\Psi_{(\DDD'C\circ\exchange)[-\ell]}
=\Psi_{\DDD'(C[\ell])\circ\exchange}$,
which can be deduced easily
as in the case of $\DDD^{\herm}$.
\hfill\qed

\vspace{.1in}

A real structure of an object
$\gbigf\in\cnumcomplextriplecat^{nd}(X)$ is
an isomorphism
$\kappa:\gammatilde^{\ast}\gbigf
\simeq
 \gbigf$
such that
$\kappa\circ\gammatilde^{\ast}\kappa=\id$.
We have the following corollary.
\begin{cor}
A real structure of
$\gbigf\in\cnumcomplextriplecat(X)$
naturally induces
a real structure of
$\nbigs_{\ell}(\gbigf)$.
\hfill\qed
\end{cor}

Let $F:X\lrarr Y$ be a proper map.
For any $\gbigf\in\cnumcomplextriplecat^{nd}(X)$
we have a naturally induced object
$F_{\ast}\gbigf:=
 (F_{\ast}\nbigf_1,F_{\ast}\nbigf_2,F_{\ast}C)$.
It gives a functor
$\cnumcomplextriplecat^{nd}(X)
\lrarr
 \cnumcomplextriplecat^{nd}(Y)$.
We have natural transformations
$F_{\ast}\DDD^{\herm}
\simeq
 \DDD^{\herm}F_{\ast}$,
$F_{\ast}\DDD
\simeq
 \DDD F_{\ast}$
and 
$F_{\ast}\DDD
\simeq
 \DDD F_{\ast}$.
If an object $\gbigf$ in $\cnumcomplextriplecat^{nd}(X)$
is equipped with a real structure $\kappa$,
then $F_{\ast}\gbigf$ is equipped with 
a naturally induced real structure $F_{\ast}\kappa$.

\subsubsection{Pairs of perverse sheaves with a pairing}

Let $\cnumtriplecat^{nd}(X)$
denote the full subcategory of 
objects $(\nbigf_1,\nbigf_2,C)$
such that $\nbigf_i$ are perverse,
in $\cnumcomplextriplecat^{nd}(X)$ 
\index{category $\cnumtriplecat^{nd}(X)$}
We have the naturally defined contravariant
auto equivalences
$\DDD$ and $\DDD^{\herm}$
on $\cnumtriplecat^{nd}(X)$.
We also have $\gammatilde^{\ast}$
on $\cnumtriplecat^{nd}(X)$.
A real structure of an object in $\cnumtriplecat^{nd}(X)$
is defined to be a real structure 
as an object in $\cnumcomplextriplecat^{nd}(X)$.

Let $\gbigf=(\nbigf_1,\nbigf_2,C)\in
 \cnumcomplextriplecat^{nd}(X)$.
Let $H^{j}(\nbigf_i)$ $(i=1,2)$ denote 
the $j$-th perverse cohomology.
We have the induced pairing
$H^j(C):
 H^{-j}(\nbigf_1)
 \otimes 
 H^j(\nbigf_2)\lrarr a_X^!\cnum$.
We set
$H^j\bigl(\gbigf\bigr):=
 \bigl(H^{-j}(\nbigf_1), H^j(\nbigf_2), H^j(C) \bigr)$.
Thus, we obtain a functor
$H^j:
 \cnumcomplextriplecat^{nd}(X)
\lrarr
 \cnumtriplecat^{nd}(X)$.
We have natural isomorphisms
$H^j\circ\DDD
\simeq
 \DDD\circ H^{-j}$,
$H^j\circ\DDD^{\herm}
\simeq
 \DDD^{\herm}\circ H^{-j}$
and
$H^j\circ\gammatilde^{\ast}
\simeq
 \gammatilde^{\ast}\circ H^{-j}$.
Therefore,
if $\gbigf$ is equipped with a real structure,
$H^j(\gbigf)$ are equipped with induced real structures.

Let $\cnumcomplextriplecat^{nd}(X,\real)$
denote the category of
objects $\gbigf\in\cnumcomplextriplecat^{nd}(X)$
with a real structure $\kappa$.
Morphisms are defined in a natural way.
Let $\cnumtriplecat^{nd}(X,\real)
\subset\cnumcomplextriplecat^{nd}(X,\real)$
denote the full subobject
$(\gbigf,\kappa)$
such that $\gbigf\in\cnumtriplecat^{nd}(X)$.
\index{category $\cnumcomplextriplecat^{nd}(X,\real)$}

\subsubsection{Perverse sheaves with real structure}

For a field $K$,
let $\Per(K_X)$ be the category of 
$K_X$-perverse sheaves.
\index{category $\Per(K_X)$}
Let $\Per(\cnum_X,\real)$ be the category of
$\cnum_X$-perverse sheaves $\nbigf$
with real structure $\rho$,
i.e., $\rho:\nbigfbar\simeq\nbigf$
such that $\rho\circ\overline{\rho}=\id$.
A morphism
$(\nbigf_1,\rho_1)\lrarr(\nbigf_2,\rho_2)$
is a morphism
$\varphi:\nbigf_1\lrarr\nbigf_2$
in $\Per(\cnum_X)$
such that $\varphi\circ\rho_1=\rho_2\circ\varphi$.
We have a natural equivalence
$\Per(\real_X)\simeq\Per(\cnum_X,\real)$
given by 
$\nbigf_{\real}\longmapsto
 \nbigf_{\real}\otimes\cnum$
with $\rho(x\otimes \alpha)=x\otimes\alphabar$.

We have a natural functor
$\Lambda:\Per(\cnum_X,\real)\lrarr
 \cnumtriplecat^{nd}(X,\real)$
For any $(\nbigf,\rho)$
in $\Per(\cnum_X,\real)$,
We have the induced non-degenerate pairing
$C_{\rho}:
 \DDD\nbigf\otimes\overline{\nbigf}
\lrarr
 a_X^!\cnum$:
\[
\begin{CD}
 \DDD\nbigf\otimes\overline{\nbigf}
@>{\DDD\rho\otimes\id}>>
 \DDD\nbigfbar\otimes\overline{\nbigf}
@>{\kappa}>>
 a_X^!\cnum
\end{CD}
\]
The morphism $\kappa$ is naturally given one.
It is easy to check that
$(\DDD\nbigf,\nbigf,C_{\rho})$
with $\kappa=(\id,\id)$
is an object in $\cnumtriplecat^{nd}(X,\real)$,
which is denoted by
$\Lambda(\nbigf,\rho)$.
The following lemma is easy to check.

\begin{lem}
The functor $\Lambda$ is an equivalence.
It is compatible with the proper push-forward,
the shift operator,
$\DDD$ and $\DDD^{\herm}$.
\hfill\qed
\end{lem}

\subsubsection{Example}
Let us look at an example.
Let $\nbigf_i=\cnum[d_X]$  $(i=1,2)$.
Let $e_i$ correspond to $1[d_X]$.
We consider a pairing
$C:\nbigf_1\otimes\overline{\nbigf_2}
\lrarr
 a_X^!\cnum$
given by
$C(e_1,\overline{e}_2)=u$,
where $u$ correspond to $1[2d_X]$
in $a_X^!\cnum=\cnum[2d_X]$.
We have
$\DDD\nbigf_i=\cnum\,e_i^{\lor}$,
where $e_i^{\lor}\otimes e_i\longmapsto u$.
Then, 
$\overline{\DDD' C}:
 \DDD\nbigf_2\otimes\overline{\DDD\nbigf_1}
\lrarr
 a_X^!\cnum$
is given by
$\overline{\DDD' C}(e_2^{\lor},e_1^{\lor})
=(-1)^{d_X}u$.
Note that
$\nbigf_1\lrarr\DDD\overline{\nbigf}_2$
is given by
$e_1\longmapsto \overline{e}_2^{\lor}$,
and 
$\overline{\nbigf}_2\lrarr\DDD\nbigf_1$
is given by
$\overline{e}_2\longmapsto
 (-1)^{d_X}e_1^{\lor}$.
Let 
$\kappa=(\kappa_1,\kappa_2):
(\DDD\nbigf_2,\DDD\nbigf_1,\overline{\DDD' C})
\lrarr
 (\nbigf_1,\nbigf_2,C)$ be given by
$\kappa_1(e_1)=e_2^{\lor}$
and 
$\kappa_2(e_1^{\lor})=(-1)^{d_X}e_2$.
The induced map
\[
\gamma^{\ast} \kappa:
 \Bigl(
 \DDD(\DDD\nbigf_1),\DDD(\DDD\nbigf_2),
 \overline{(\DDD'(\overline{\DDD' C}))}
 \Bigr)
\lrarr
 (\DDD\nbigf_1,\DDD\nbigf_2,\overline{\DDD' C})
\]
given by
$(\gamma^{\ast}\kappa)_1(e_2^{\lor})
=(-1)^{d_X}(e_1^{\lor})^{\lor}$
and
$(\gamma^{\ast}\kappa)_2(e_1^{\lor})
=(e_2^{\lor})^{\lor}$.
Note the natural identification
$(e_i^{\lor})^{\lor}=(-1)^{d_X}e_i$.
Hence, we have
$\kappa\circ\gamma^{\ast}\kappa=\id$.

\subsection{Compatibility with the de Rham functor}
\label{subsection;13.4.13.10}

For any object
$\nbigt=(\nbigm_1^{\bullet},\nbigm_2^{\bullet},C)
 \in\Dcomplextriplecat(X)$,
we define
$\DR_X(\nbigt)\in \cnumcomplextriplecat(X)$
as in \S\ref{subsection;13.4.5.1}.
It induces an equivalence of the categories
$\Dtriplecat^{nd}(X)$
and 
$\cnumtriplecat^{nd}(X)$.

\begin{prop}
\label{prop;11.3.24.1}
For $\nbigt\in\Dtriplecat^{nd}(X)$,
we have natural isomorphisms
\[
 \DDD\circ \DR_X(\nbigt)
\simeq
 \DR_X\circ\DDD(\nbigt),
\quad\quad
 \DDD^{\herm}\circ\DR_X(\nbigt)
\simeq
 \DR_X\circ\DDD^{\herm}(\nbigt),
\]
which are given by
the natural identifications of 
underlying complexes.
\end{prop}
\pf
Let us consider the case of $\DDD^{\herm}$.
We remark that the conjugate 
$\overline{\DR_X(C)}$ in Lemma \ref{lem;11.3.25.1}
is taken for the natural real structure of
$\distribution^{\bullet}_X[2d_X]$.
Because the quasi-isomorphism
$a_X^{!}\cnum\simeq \distribution_X[2d_X]$
is given by the multiplication of $(2\pi\sqrt{-1})^{d_X}$,
the conjugate with respect to the real structure of
$a_X^{!}\cnum$ is $(-1)^{d_X}\overline{\DR_X(C)}$.
Then, by using Lemma \ref{lem;11.3.25.1},
we can compare the pairings.

Let us consider the case of $\DDD$.
We have the natural identifications:
\[
 \DDD (\overline{\DR_XM})
=\DDD (\DR_{\Xbar}\Mbar)
\simeq
 \DR_{\Xbar}(\DDD_{\Xbar}\Mbar)
\simeq
 \overline{\DR_X\DDD M}
\simeq
 \overline{\DDD \DR_X M}
\]
Here, we use the identifications
$\Omega^{\bullet}_X[d_X]\otimes
 \Omega^{\bullet}_{\Xbar}[d_X]\otimes
 \distribution_X
\simeq
 \distribution_X^{\bullet}$
for the exchanges of
$\DR_X\circ\DDD_X\simeq
 \DDD\circ\DR_X$
and $\DR_{\Xbar}\circ\DDD_{\Xbar}\simeq
 \DDD\circ\DR_{\Xbar}$,
the composite
$\DDD \overline{\DR_XM}
\simeq
 \overline{\DR_X M}$
is induced by the real structure of 
$\distribution_X^{\bullet}[2d_X]\simeq
a_X^!\cnum$.
Then, it follows from the the definition of 
$\DDD C$
that the pairings $\DR\DDD (C)$
and $\DDD \DR(C)$
are equal under the natural identifications.
\hfill\qed

\begin{cor}
We have 
$\gammatilde^{\ast}\DR_X(\nbigt)
\simeq
 \DR_X\gammatilde^{\ast}\nbigt$
given by the natural isomorphisms
of the underlying complexes.
In particular,
the de Rham functor gives a functor
$\Dtriplecat^{nd}(X,\real)
\lrarr
 \cnumtriplecat^{nd}(X,\real)$.
\hfill\qed
\end{cor}

Let $\Dcomplextriplecat^{nd}(X)\subset\Dcomplextriplecat(X)$
denote the full subcategory of
objects 
$\gbigt=(\nbigm^{\bullet}_1,\nbigm^{\bullet}_2,C)$
such that
(i) $H^i(\gbigt)\in\Dtriplecat^{nd}(X)$,
(ii) the dual $\DDD C$ exist
in the sense of \S\ref{subsection;11.4.6.1}.
The functor $\DR$ gives a functor
$\Dcomplextriplecat^{nd}(X)
\lrarr
 \cnumcomplextriplecat^{nd}(X)$.

\begin{lem}
\label{lem;13.4.5.2}
By the rule of the signatures,
we have a natural transform
 $\nbigs_{\ell}\circ\DR\simeq
 \DR\circ\nbigs_{\ell}$
of the functors from
$\Dcomplextriplecat(X)$
to $\cnumcomplextriplecat(X)$.
\end{lem}
\pf
We have the natural isomorphism
$\Tot\DR_X(M_i^{\bullet}[\pm\ell])
\simeq
 \Tot\DR_X(M_i^{\bullet})[\mp\ell]$,
which are given by 
the multiplication of $(-1)^{p\ell}$
on $\Omega^{d_X+p}\otimes M_i^{r\pm \ell}$.
We can compare the pairings
by a direct computation.
\hfill\qed

\begin{prop}
We have natural isomorphisms
$\DR\DDD^{\herm}\simeq
 \DDD^{\herm}\DR$
and
$\DR\DDD\simeq
 \DDD\DR$,
which are given by the natural identifications
of the underlying complexes.
As a consequence,
we have 
$\DR\circ\gammatilde^{\ast}
 \simeq
 \gammatilde^{\ast}\circ\DR$.
\end{prop}
\pf
Let us consider the case of $\DDD^{\herm}$.
We have the natural isomorphism
of the underlying complexes
in the case of $\DDD^{\herm}$.
The comparison of the pairings 
can be reduced to Proposition \ref{prop;11.3.24.1}
by Lemma \ref{lem;11.3.24.3},
Lemma \ref{lem;11.3.24.2}
and 
Lemma \ref{lem;13.4.5.2}.

Let us consider the case of $\DDD$.
Let $M^{\bullet}$ be a complex of 
holonomic $D$-modules.
Let us observe that 
the following diagram of the natural 
quasi-isomorphisms is commutative:
\begin{equation}
 \label{eq;13.4.5.4}
 \begin{CD}
\DR\nbigs_{-\ell}\DDD M^{\bullet}
@>>>
 \DR\DDD\nbigs_{\ell}(M^{\bullet}) \\
 @VVV @VVV \\
 \nbigs_{-\ell}\DDD \DR(M^{\bullet})
 @>>>
 \DDD\nbigs_{\ell}\DR(M^{\bullet})
 \end{CD}
\end{equation}
Indeed, we take an $D_X$-injective resolution $I_X^{\bullet}$
of $\nbigo_X[d_X]$.
We take an $\cnum_X$-injective resolution $J^{\bullet}$ of
$\Tot\Omega_X^{\bullet}[d_X]\otimes I_X^{\bullet}$.
We have the following commutative diagram of the natural morphisms:
{\scriptsize
\[
 \begin{array}{ccc}
 \Tot
 \Omega_X^{\bullet}[d_X]
 \otimes
 \bigl(
 \nhom_{D_X}(M^{\bullet},I_X^{\bullet})[-\ell]
 \bigr)
 & \lrarr &
 \Tot
 \Omega_X^{\bullet}[d_X]
 \otimes
 \nhom_{D_X}(M^{\bullet}[\ell],I_X^{\bullet})
 \\
 \darr & & \darr \\
\Tot
 \nhom_{\cnum_X}(
 \Omega_X^{\bullet}[d_X]\otimes
 M^{\bullet},
 \Omega_X^{\bullet}[d_X]\otimes I_X^{\bullet})[-\ell]
 & \lrarr &
 \Tot
 \nhom_{\cnum_X}(
 \Omega_X^{\bullet}[d_X]\otimes
 M^{\bullet}[\ell],
 \Omega_X^{\bullet}[d_X]\otimes I_X^{\bullet})
 \\
 \darr & & \darr \\
\Tot
 \nhom_{\cnum_X}(
 \Omega_X^{\bullet}[d_X]\otimes
 M^{\bullet},
 J^{\bullet}_X)[-\ell]
 & \lrarr &
 \Tot
 \nhom_{\cnum_X}(
 \Omega_X^{\bullet}[d_X]\otimes
 M^{\bullet}[\ell],
 J^{\bullet}_X)
 \end{array}
\]
}
Then, we obtain the commutativity of
(\ref{eq;13.4.5.4}).

The comparison of the parings 
can be reduced to Proposition \ref{prop;11.3.24.1}
by 
Lemma \ref{lem;11.3.24.3},
Lemma \ref{lem;11.3.24.2}
and 
Lemma \ref{lem;13.4.5.2},
as in the case of $\DDD^{\herm}$.
\hfill\qed

\vspace{.1in}
We remark the commutativity of the following
diagrams of the natural isomorphisms:
\[
 \begin{CD}
 \nbigs_{-\ell}\DDD \DR(\nbigt)
 @>>>
 \DDD\nbigs_{\ell} \DR(\nbigt)
 \\
 @VVV @VVV \\
 \DR\nbigs_{-\ell}\DDD(\nbigt)
 @>>>
 \DR\DDD\nbigs_{\ell}(\nbigt)
 \end{CD}
\]
We have similar diagrams
for $\DDD^{\herm}$ and $\gammatilde^{\ast}$.
In particular,
the shift functor of real structures are compatible
with the de Rham functor.

We have the commutativity of the following diagram:
\[
 \begin{CD}
 \DDD\gammatilde^{\ast}\DR\nbigt
 @>>>
 \gammatilde^{\ast}\DDD \DR\nbigt \\
 @VVV @VVV \\
 \DR\DDD\gammatilde^{\ast}\nbigt
 @>>>
 \DR\gammatilde^{\ast}\DDD\nbigt
 \end{CD}
\]
Hence, the dual of real structures
is compatible with the de Rham functor.
Similarly,
the functors $\DDD^{\herm}$ and $\gammatilde^{\ast}$
of real structures are compatible
with the de Rham functor.

Let $\nbigt\in\Dcomplextriplecat^{nd}(X)$.
Let $F:X\lrarr Y$ be a morphism of complex manifolds
such that $Supp(\nbigt)$ is proper with respect to $F$.
By Proposition \ref{prop;11.3.21.22},
the following is commutative:
\[
 \begin{CD}
 \DR_Y\circ F_{\dagger}\circ\gammatilde^{\ast}(\nbigt)
 @>>>
 \DR_Y\circ \gammatilde^{\ast}\circ F_{\dagger}(\nbigt)\\
 @VVV @VVV \\
 F_{\ast}\circ\gammatilde^{\ast}\circ\DR_X(\nbigt)
 @>>>
 \gammatilde^{\ast}\circ F_{\ast}\circ \DR_X(\nbigt)
 \end{CD}
\]
In particular,
the push-forward functor of real structures
is compatible with the de Rham functor.

\subsection{Regular case}

Let $\Dtriplecat^{nd}_{reg}(X,\real)\subset 
 \Dtriplecat^{nd}(X,\real)$ 
denote the full subcategory of
the objects whose underlying $D$-modules
are regular holonomic.
By the Riemann-Hilbert correspondence,
the de Rham functor induces an equivalence
\[
\DR:\Dtriplecat^{nd}_{reg}(X,\real)
\lrarr
 \cnumtriplecat^{nd}(X,\real)
\simeq
 \Per(\real_X).
\]
According to the results in \S\ref{subsection;13.4.13.10},
it is compatible with the push-forward
by proper morphisms.
It is also compatible with the dual.

\subsubsection{Description}
\label{subsection;13.8.6.1}

Let $M$ be any regular holonomic $D_X$-module
such that $\DR_X(M)$ has a real structure,
i.e., $\iota:\overline{\DR_X(M)}\simeq \DR_X(M)$
such that $\iota\circ\iotabar=\id$.
It induces a pairing
$C_{\iota}:
 \DDD_XM\otimes
 \Mbar\lrarr\distribution_X$
corresponding to
\[
 \DR_X(\DDD_XM)\otimes
 \DR_{\Xbar}(\Mbar)
\lrarr
 \DDD_X\DR_X(M)\otimes
 \DR_X(M)
\lrarr
 \omega^{top}.
\]
We put $\nbigt:=(\DDD M,M,C_{\iota})$.
\begin{prop}
$\nbigt$ has a real structure given by
$\kappa=(\id,\id)$,
i.e.,
$(\nbigt,\kappa)\in
 \Dtriplecat^{nd}_{reg}(X,\real)$.
Conversely, any object in
$\Dtriplecat^{nd}_{reg}(X,\real)$
has such a description
up to isomorphisms.
\end{prop}
\pf
Let us prove the first claim.
We have only to check that $(\id,\id)$
gives a morphism
$\gammatilde^{\ast}\nbigt\lrarr\nbigt$.
Namely,
we have only to check that
$(\DDD C_{\iota})^{\ast}=C_{\iota}$
under the natural isomorphism
$\DDD_X\DDD_X(M)\simeq M$.
The morphism
$\Mbar\lrarr \nbigc_X(\DDD M)$
given by $C_{\iota}$
induces 
$\iota: 
 \DR_{\Xbar}(\Mbar)
\simeq
 \DR_X(M)$,
if we use the identification
$\Omega_X^{\bullet}[d_X]
\otimes
 \Omega_{\Xbar}^{\bullet}[d_X]
\otimes\distribution_X
\simeq
 \distribution_X^{\bullet}[2d_X]$.
See Remark \ref{rem;13.4.13.1}
for the signature.
The morphism
$\DR_X(\DDD_X M)\simeq
 \DR_{\Xbar}( \DDD_{\Xbar}\Mbar)$
induced by $C_{\iota}$
is identified with the morphism
induced by
$\iota:\DR_{\Xbar}\Mbar\simeq\DR_XM$.
Then, we obtain the commutativity of the following diagram:
\[
 \begin{CD}
 \DR M
 \otimes \DR \DDD \Mbar
 @>{\DR(\DDD C_{\iota})}>>
 a_X^!\cnum\\
 @VV{\iota\otimes\iota^{\ast}}V @VV{\id}V \\
 \DR_{\Xbar}(\Mbar)
 \otimes
 \DDD\DR_XM
 @>{\DR(C_{\iota})\circ\exchange}>>
 a_X^!\cnum
 \end{CD}
\]
Then, we obtain the commutativity of the following diagram:
\[
 \begin{CD}
 \DR \DDD M
 \otimes\DR \Mbar 
 @>{\overline{\DR(\DDD C_{\iota})\circ\exchange}}>>
  a_X^!\cnum \\
 @VV{\id\otimes\id}V @VV{\id}V \\
 \DR_{X}(\DDD M)
 \otimes
 \DR_{\Xbar}M
 @>{\DR(C_{\iota})}>>
 a_X^!\cnum
 \end{CD}
\]
It means $\DDD(C_{\iota})^{\ast}=C_{\iota}$.
As for the second claim,
we have only to care the signatures
which can be checked as in the case of first claim.
\hfill\qed

\subsubsection{Localizations}

Let $H$ be a hypersurface of $X$.
Let $j:X\setminus H\lrarr X$ be the inclusion.
We have the functor $[\ast H]$
on $\Dtriplecat^{nd}_{reg}(X,\real)$.
It is equivalent to $j_{\ast}j^{\ast}$
on $\Per(\real_{X})$.
The essential image of the functor $[\ast H]$
is naturally equivalent to
$\Dtriplecat^{nd}_{\reg}(X\setminus H,\real)
\simeq
 \Per(\real_X)$.

Similarly, the functor $[!H]$ on 
$\Dtriplecat^{nd}_{\reg}(X,\real)$
is equivalent to
$j_!j^{\ast}$ on 
$\Per(\real_X)$,
and the essential image 
is naturally equivalent to
$\Dtriplecat^{nd}_{\reg}(X\setminus H,\real)
\simeq
 \Per(\real_{X\setminus H})$.

\subsubsection{Beilinson functors}

For any integers $a\leq b$,
let $\nbigi^{a,b}$ be the $\real$-local system
on $\cnum^{\ast}$
underlying the smooth 
$\nbigr_{\cnum(\ast 0)}$-triple $\II^{a,b}$
with the real structure.
By using $\nbigi^{a,b}$,
we obtain the nearby cycle functors $\psi^{(a)}$,
the maximal functors $\Xi^{(a)}$
and the vanishing cycle functors $\phi^{(a)}$
on $\Per(\real_X)$.
By the compatibility of the localizations
with the de Rham functor
$\DR_X:\Dtriplecat^{nd}_{reg}(X,\real)
\lrarr \Per(\real_X)$,
we obtain that $\DR_X$ is compatible with
the functors $\psi^{(a)}$,
$\Xi^{(a)}$ and $\phi^{(a)}$.

\subsection{$\real$-Betti structure}
\label{subsection;13.8.29.40}

Let $\nbigt=(M_1,M_2,C)\in\Dtriplecat^{nd}(X,\real)$.
By the construction above,
we have a real structure of $\DR_X(M_2)$,
i.e.,
a pre-$\real$-Betti structure of $M_2$
in the sense of \cite{mochi9}.
\begin{prop}
\label{prop;13.8.30.1}
The pre-$\real$-Betti structure is
an $\real$-Betti structure of $M_2$
in the sense of {\rm\S7.2}
in {\rm\cite{mochi9}}.
\end{prop}
\pf
We have only to check the claim
locally around any points $P$ of the support of $M_2$.
We shall shrink $X$ without mention.
We use a Noetherian induction on the support.
If the support of $M_2$ is a point,
then the claim is clear.
We take a cell $\nbigc=(Z,U,\varphi,V)$
of $M_2$ at $P$
with a cell function $g$.
We have the $\real$-structure
of the meromorphic flat bundle $V$
in the sense of \S6.4 of \cite{mochi9}.

\begin{lem}
\label{lem;13.8.29.31}
The $\real$-structure is good 
in the sense of 
{\rm \S6.4} of {\rm\cite{mochi9}}.
\end{lem}
\pf
We consider the special case $X=Z=\Delta=\{|z|<1\}$
with $\dim X=1$.
We may assume $M_i(\ast P)=V_i$.
We have the pairing
$M_1(\ast P)\times 
\overline{M_2(\ast P)}
\lrarr
 \distribution_{X}(\ast P)$.
By Lemma 12.6.10 of \cite{mochi7},
we can check that it comes from
a pairing on the real blow up.
Hence, the real structure is compatible
with the Stokes structure.
Moreover,
we may replace 
$\distribution_X(\ast P)$
with $\nbigc^{\infty\moderate P}_X$.

Let us consider the case $\dim Z=1$.
We may assume that 
$X=\Delta_t\times\Delta^m$ with $g=t$.
Let $h:X\lrarr \Delta_t$
denote the projection.
We obtain
$h_{\dagger}(M_1,M_2,C)$
with an induced real structure.
Then, $V$ is a direct summand of
$(h\circ\varphi)^{\ast}h_{\dagger}M_2(\ast g)$.
Hence, by the result in the previous paragraph,
we obtain that the $\real$-structure of $V$ is good.

For any irreducible smooth curve 
$C\subset Z$ with
$C\not\subset Z\setminus U$,
we can check that
the $\real$-structure of $V_{|C}$ is good,
by using the result in the previous paragraph.

Let $Q$ be any point of $Z$.
Let $(Z_Q,\psi_Q)$ be a local resolution of $Z$
around $Q$.
(See \S6.4 of \cite{mochi9}.)
We can easily deduce that
the $\real$-structure of $\psi_Q^{\ast}V$ is good,
by using the result in the previous paragraph.
Hence, the $\real$-structure of $V$ is good.
\hfill\qed

\vspace{.1in}
We have the induced $\real$-structures of
$\nbigt[\star g]$,
and the morphisms
$\nbigt[!g]\lrarr\nbigt\lrarr\nbigt[\star g]$
are compatible with the $\real$-structures.
We have two pre-$\real$-Betti structures on
$M_2[\star g]=\varphi_{\star}(V)$.
One is induced by the $\real$-structure of
$\nbigt[\star g]$.
The other is induced by the good $\real$-structure
of $V$.

\begin{lem}
\label{lem;13.8.29.30}
The two pre-$\real$-Betti structures are the same.
\end{lem}
\pf
Let $V_2:=V$.
Let $(Z,U,\varphi,V_1)$ be a cell of $M_1$.
We have the naturally defined sesqui-linear pairing
$C^U_V:
 V_{1|U}\times\Vbar_{2|U}\lrarr \nbigc^{\infty}_{U}$.
For any irreducible smooth curve $C\subset Z$
with $C\not\subset D:=Z\setminus U$,
the pairing of $V_{i|U\cap C}$
is extended to
$V_{1|C}\times
 \Vbar_{2|C}
\lrarr
 \nbigc^{\infty\moderate D}_{Z}$,
which can be shown by an argument
in the proof of Lemma \ref{lem;13.8.29.31}.
For any $Q\in Z$,
we take a local resolution 
$(Z_Q,\psi_Q)$ of $V_i$.
By considering the Stokes filtrations of $V_i$,
and by using Lemma 12.6.10 of \cite{mochi7},
we obtain a unique extension of the pairing
$\psi_Q^{\ast}V_1
\times
 \overline{\psi_Q^{\ast}V_2}
\lrarr
 \nbigc^{\infty\moderate \check{D}_Q}_{\check{Z}_Q}$.
Then, we obtain a unique extension 
of the pairing
$C_V^{\moderate}:V_1\times
 \overline{V_2}
\lrarr
 \nbigc^{\infty\moderate D}_{Z}$.
Hence, we have
$C_{V\ast}:
V_{1!}\times
\overline{V_{2\ast}}
\lrarr
 \distribution_Z$
and 
$C_{V!}:
V_{1\ast}\times
\overline{V_{2!}}
\lrarr
 \distribution_Z$.
We have the $\real$-structure of
$\nbigt_V^U:=
 (V_{1|U},V_{2|U},C^U_V)$.
Let 
$\nbigt_{V\ast}
:=(V_{1!},V_{2\ast},C_{V\ast})$
and 
$\nbigt_{V!}
:=(V_{1\ast},V_{2!},C_{V!})$.
For any irreducible curve
$C\subset Z$,
the $\real$-structure of
$\nbigt_{V|C\cap U}^U$
is extended to those of
$\nbigt_{V\star|C}$,
which can be shown by an argument
in the proof of Lemma \ref{lem;13.8.29.31}.
Hence, the $\real$-structure
of $\nbigt_V^U$
is extended to those of
$\nbigt_{V\star}$.
The induced pre-$\real$-Betti structures of $V_{2\star}$
are the $\real$-Betti structures.
Because the $\real$-structures of
$\nbigt[\star g]$
is obtained as the push-forward of the $\real$-structures
of $\nbigt_{V\star}$,
we obtain the claim of Lemma \ref{lem;13.8.29.30}.
\hfill\qed

\vspace{.1in}
By Lemma \ref{lem;13.8.29.31}
and Lemma \ref{lem;13.8.29.30},
the cell $\nbigc$ is compatible with 
the pre-$\real$-Betti structure of $M_2$
in the sense of \S7.1.3 in \cite{mochi9}.
We have the induced $\real$-structure
on $\phi^{(0)}(\nbigt)$.
The induced pre-$\real$-Betti structure
of $\phi^{(0)}(M_2)$ is a $\real$-Betti structure,
by the hypothesis of the induction.
Thus, the pre-$\real$-Betti structure of
$M_2$ is a $\real$-Betti-structure.
\hfill\qed

\vspace{.1in}
Let $\Hol(X,\real)$ denote the category of
holonomic $D_X$-modules with $\real$-Betti structure.
By Proposition \ref{prop;13.8.30.1},
we obtain a functor
$\Upsilon:
 \Dtriplecat^{nd}(X,\real)\lrarr\Hol(X,\real)$.
\begin{prop}
The functor $\Upsilon$
is an equivalence.
\end{prop}
\pf
It is clearly faithful.
Let us observe that the functor is full.
Let $\nbigt^{(i)}=(M_1^{(i)},M_2^{(i)},C^{(i)})$
with $\kappa^{(i)}$
be objects in $\Dtriplecat^{nd}(X,\real)$.
Let $f:M_2^{(1)}\lrarr M_2^{(2)}$
be a morphism in $\Hol(X,\real)$.
The following induced morphism is commutative:
\begin{equation}
 \label{eq;13.8.30.2}
 \begin{CD}
 \DR(\DDD M_2^{(2)})
 @>{\DR(\DDD f)}>>
 \DR(\DDD M_2^{(1)})
 \\
 @VVV @VVV \\
 \overline{\DR(\DDD M_2^{(2)})}
 @>{\overline{\DR(\DDD f)}}>>
 \overline{\DR(\DDD M_2^{(1)})}
 \end{CD}
\end{equation}
Here, the vertical arrows are induced by
the real structure of
$\DR\DDD(M_2^{(i)})$.

We set $f_2:=f$,
and let $f_1:M_1^{(2)}\lrarr M_1^{(1)}$
be the morphism induced by
$\DDD(f):
 \DDD M_2^{(2)}\lrarr \DDD M_2^{(1)}$
with isomorphisms
$M_1^{(i)}\simeq \DDD M_2^{(i)}$
underlying $\kappa^{(i)}$.
Let us consider the following diagram:
\begin{equation}
 \label{eq;13.8.30.3}
 \begin{array}{ccc}
 \DR\DDD M_2^{(2)}
 &
 \simeq \DR M_1^{(2)} \simeq
 &
 \DR\nbigc_{\Xbar}\Mbar_2^{(2)}
 \\
 \darr & & \darr \\
\DR\DDD M_2^{(1)}
 &
\simeq \DR M_1^{(1)}\simeq
 &
 \DR\nbigc_{\Xbar}\Mbar_2^{(1)}
 \end{array}
\end{equation}
Here, the vertical arrows are given as
$\DR\DDD f$
and $\DR\nbigc_{\Xbar}\fbar$,
respectively.
Under the natural isomorphism
$\DR_X\nbigc_{\Xbar}\Mbar_2^{(i)}
\simeq
 \DR_{\Xbar}\DDD_{\Xbar}\Mbar_2^{(i)}$,
we have
$\DR\nbigc_{\Xbar}(\fbar)
=\overline{\DR \DDD f}$.
The horizontal arrows are induced by the real
structures of
$\DR \DDD M_2^{(i)}$.
Hence, (\ref{eq;13.8.30.3}) is commutative.
Then, we obtain that
$(f_1,f_2)$ gives a morphism
$\nbigt^{(1)}\lrarr\nbigt^{(2)}$
such that 
$\Upsilon(f_1,f_2)=f$.
Hence, $\Upsilon$ is full.

Let us observe that $\Upsilon$ is essentially surjective.
By the fully faithfulness,
we have only to check it locally around any point $P$ of $X$.
We may shrink $X$ without mention.
We use a Noetherian induction on 
the support of $M$.
Let $M\in\Hol(X,\real)$.
We take a cell $\nbigc=(Z,U,\varphi,V)$ of $M$ at $P$
with a cell function $g$.
We have a description of $M$
as the cohomology of
\begin{equation}
 \label{eq;13.8.30.4}
 \psi_g(\varphi_{\dagger}V)
\lrarr
 \Xi_g(\varphi_{\dagger}V)
\oplus
 \phi_g(M)
\lrarr
 \psi_g(\varphi_{\dagger}V)
\end{equation}
The $\real$-structure of $V$ is good.
As in the proof of Lemma \ref{lem;13.8.29.30},
we obtain $\nbigt_{V\star}\in\Dtriplecat^{nd}(Z,\real)$ 
$(\star=\ast,!)$.
By using the Beilinson construction,
we obtain
$\Xi(\varphi_{\dagger}V)
 \in
 \Dtriplecat^{nd}(X,\real)$.
By the hypothesis of the induction,
we have $\nbigt_i\in\Dtriplecat^{nd}(X,\real)$
$(i=1,2)$
such that
$\Upsilon(\nbigt_1)\simeq 
 \psi_g(\varphi_{\dagger}V)$
and 
$\Upsilon(\nbigt_2)\simeq 
 \phi_g(\varphi_{\dagger}V)$.
By the fully faithfulness of $\Upsilon$,
(\ref{eq;13.8.30.4})
comes from a complex 
in $\Dtriplecat^{nd}(X,\real)$.
Hence, $M$ comes from an object in
$\Dtriplecat^{nd}(X,\real)$,
i.e.,
we obtain the essential surjectivity of $\Upsilon$.
\hfill\qed

\vspace{.1in}

\begin{prop}
\label{prop;13.8.30.12}
The functor $\Upsilon$
is compatible with the following functors:
\begin{itemize}
\item
The push-forward by any proper morphisms.
\item
The dual $\DDD$.
\end{itemize}
\end{prop}
\pf
It follows from the compatibilities
in \S\ref{subsection;13.4.13.10}.
\hfill\qed

\vspace{.1in}
For any hypersurface $H$ on $X$,
we have the localization functors
$[\star H]$ $(\star=\ast,!)$
on $\Hol(X,\real)$.
\begin{prop}
\label{prop;13.8.30.13}
The localizations are compatible with $\Upsilon$,
i.e.,
for any $M\in \Dtriplecat^{nd}(X,\real)$,
we naturally have
$\Upsilon(M[\star H])\simeq
 \Upsilon(M)[\star H]$.
\end{prop}
\pf
It follows from the characterization
of the localization
as in Theorem 8.1.4 of \cite{mochi9}.
\hfill\qed

\vspace{.1in}
Let $H$ be any hypersurface of $X$.
Let $\nbigt_1=(V_1,V_2,C)$ be a smooth $D_{X(\ast H)}$-triple
with a real structure.
Then, $V_2$ is a meromorphic flat bundle
on $(X,H)$ with a good real structure.
\begin{prop}
\label{prop;13.8.30.14}
Let $\nbigt\in\Dtriplecat^{nd}(X,\real)$.
We have
$\Upsilon(\nbigt\otimes\nbigt_1[\star H])
\simeq
 \bigl(
 \Upsilon(\nbigt)
 \otimes V_2
 \bigr)[\star H]$
in $\Hol(X,\real)$.
\end{prop}
\pf
If $H=\emptyset$,
the claim is clear.
The general case easily follows.
\hfill\qed

\begin{cor}
We naturally have
$\Upsilon\circ\Xi
\simeq \Xi\circ\Upsilon$,
$\Upsilon\circ\psi
\simeq
 \psi\circ\Upsilon$
and 
$\Upsilon\circ\phi
\simeq
 \phi\circ\Upsilon$.
\hfill\qed
\end{cor}

\subsection{Basic examples}
\label{subsection;11.3.30.2}

\subsubsection{Real structures
on the simplest $D$-triple}

Let $Z$ be a complex manifold.
We consider the $D$-triple 
$\nbigu_Z=(\nbigo_Z,\nbigo_Z,C_0)$,
where $C_0(f,\gbar)=f\cdot \gbar$.
We have
$\gammatilde^{\ast}\nbigu_Z=\bigl(
 \DDD\nbigo_Z,\DDD\nbigo_Z,\DDD C_0\bigr)$.
Let
$\nu:\DDD\nbigo_Z\simeq \nbigo_Z$
be the isomorphism
determined by the following condition.
We have the natural identification of
the cohomology sheaf $H^{-d_Z}\bigl(\DR\nbigo_Z\bigr)$
and the sheaf of flat sections of $\nbigo_Z$.
The isomorphisms
$\DR\DDD \nbigo_Z
\simeq
 \nhom_{D_Z}(\nbigo_Z,\!\nbigo_Z[d_Z])
\simeq
 \nhom_{\cnum_Z}(\DR\nbigo_Z,\!\DR\nbigo_Z[d_Z])$
induce the following isomorphism:
\begin{multline}
 H^{-d_Z}\bigl(
 \DR \DDD\nbigo_Z
 \bigr)
\simeq
  \nhom_{\cnum_Z}\bigl(
 H^{-d_Z}(\DR\nbigo_Z),
 H^{-2d_Z}(\DR\nbigo_Z[d_Z])
 \bigr)
 \\
\simeq
 \nhom_{\cnum_Z}\bigl(
 \cnum_Z,\cnum_Z
 \bigr)
\end{multline}
Then, $\nu$ induces
$H^{-d_Z}(\DR\DDD\nbigo_Z)
\simeq
 H^{-d_Z}(\DR\nbigo_Z)$
given by
$\id\longmapsto 1$.
We remark that the composition of
the natural isomorphisms
\[
\begin{CD}
 \nbigo_Z@>>>
 \DDD(\DDD\nbigo_Z)
 @<{\DDD\nu}<<
 \DDD\nbigo_Z
 @>{\nu}>>
\nbigo_Z
\end{CD}
\]
is $(-1)^{d_Z}$.
In other words,
$(-1)^{d_Z}\nu$ can be identified with $\DDD\nu$.
See also the appendix below.
The following proposition is clear.
\begin{prop}
Let $a$ be a complex number such that $|a|=1$.
Then, 
\[
 (a\,\nu^{-1},(-1)^{d_Z}\overline{a}\,\nu):
 \gamma^{\ast}\nbigu_Z\simeq\nbigu_Z
\]
gives a real structure.
\hfill\qed
\end{prop}
Let us consider more specific real structures
of the shifts of $\nbigu_Z$.

\subsubsection{The $D$-triple
corresponding to the simplest
$\real$-perverse sheaf}

We have the $D$-triple
$(\DDD\nbigo_Z,\nbigo_Z,C_{\iota})$
corresponding to 
the simplest $\real$-perverse sheaf $\real[d_Z]$,
as in \S\ref{subsection;13.8.6.1}.
Here, $C_{\iota}$ is the pairing
induced by the natural real structure
of $\DR_Z(\nbigo_Z)$.
We can check 
$(v^{-1},1):
(\DDD\nbigo_Z,\nbigo_Z,C_{\iota})
\simeq
 \nbigu_Z$.
The $(\id,\id)$ gives a real structure of
$(\DDD\nbigo_Z,\nbigo_Z,C_{\iota})$,
and the corresponding real structure of
$\nbigu_Z$ is given by $(1,(-1)^{d_Z})$.

\subsubsection{Natural real structures of $\nbigu_Z[d_Z]$
and $\nbigu_Z[-d_Z]$}

Let us observe that
$\nbigu_Z[d_Z]$
and $\nbigu_Z[-d_Z]$
in $\nbigc\bigl(\Dtriplecat(Z)\bigr)$,
or equivalently 
in $\Dcomplextriplecat(Z)$,
have canonical real structures.
We assume $Z$ is connected for simplicity.

We have the natural identification
$H^0\bigl(\DR\nbigo_Z[-d_Z]\bigr)$
with the sheaf of flat sections of $\nbigo_Z$,
i.e.,
$H^0\bigl(\DR\nbigo_Z[-d_Z]\bigr)=\cnum$.
We have the natural identifications
\begin{multline}
 H^0\bigl(
 \DR
 \DDD(\nbigo_Z[d_Z])
 \bigr)
\simeq 
 \nhom_{D_Z}\bigl(
 \nbigo_Z[d_Z],\nbigo_Z[d_Z]
 \bigr)
 \\
\simeq
 \nhom_{\cnum_Z}\bigl(
 \DR\nbigo_Z[d_Z],
 \DR\nbigo_Z[d_Z]
 \bigr).
\end{multline}
We have the unique isomorphism
$\mu:\DDD(\nbigo_Z[d_Z])\simeq\nbigo_Z[-d_Z]$
such that it induces
$H^0\DR\DDD(\nbigo_Z[d_Z])
\simeq
 H^0(\DR\nbigo_Z[-d_Z])=\cnum$
given by
$\id\longmapsto 1$.
We have the induced morphism
$\DDD\mu:\DDD\bigl(\nbigo_Z[-d_Z]\bigr)
\simeq
 \nbigo_Z[d_Z]$.
Note the following commutative diagram:
\begin{equation}
 \label{eq;11.3.27.1}
 \begin{CD}
 \DDD(\nbigo_Z[d_Z]) 
 @>{\mu}>>
 \nbigo_Z[-d_Z]\\
 @VVV @V{=}VV \\
 \DDD(\nbigo_Z)[-d_Z]
 @>{\epsilon(-d_Z)\nu}>>
 \nbigo_Z[-d_Z]
 \end{CD}
\end{equation}
Here, the left vertical arrow is
given by the exchange of the dual and the shift.
(It is given by the multiplication of
$(-1)^{pd_Z}\epsilon(d_Z)$
on the degree $p$-part.)
We also have the following commutative diagram:
\begin{equation}
\label{eq;13.8.6.2}
 \begin{CD}
 \DDD(\nbigo_Z[-d_Z])
 @>{\DDD\mu}>>
 \nbigo_Z[d_Z]\\
 @VVV @V{=}VV \\
 \DDD(\nbigo_Z)[d_Z]
 @>{\epsilon(d_Z)\,\nu}>>
 \nbigo_Z[d_Z]
 \end{CD}
\end{equation}
Here, the left vertical arrow is given by 
the exchange of the dual and the shift.
The commutativity follows from
the commutativity of (\ref{eq;11.3.27.1}).

\begin{prop}
The pair
$\bigl(\mu^{-1},\DDD\mu\bigr)$
gives a real structure of 
the object
$\nbigu_Z[d_Z]$
in $\Dcomplextriplecat(Z)$
obtained as the shift of the real structure
$\bigl(\epsilon(-d_Z)\,\nu^{-1},\epsilon(d_Z)\,\nu\bigr)$
of $\nbigu_Z$.
The pair $\bigl((\DDD\mu)^{-1},\mu\bigr)$
is a real structure of $\nbigu_Z[-d_Z]$
in $\Dcomplextriplecat(Z)$
obtained from the real structure
$\bigl(\epsilon(d_Z)\nu^{-1},\epsilon(-d_Z)\,\nu\bigr)$
of $\nbigu_Z$
by the shift.
\end{prop}
\pf
The commutativity of the diagrams
(\ref{eq;11.3.27.1}) and (\ref{eq;13.8.6.2})
implies the claims.
\hfill\qed

\begin{cor}
We have the corresponding real structures
on the objects $\nbigu_Z[d_Z]$ and $\nbigu_Z[-d_Z]$
in $\nbigc(\Dtriplecat(Z))$.
\hfill\qed
\end{cor}

\subsubsection{Trace morphism}

For simplicity, we assume that $Z$ is proper.
Let us look at the push-forward
$a_{Z\dagger}\nbigu_Z
=(a_{Z\dagger}\nbigo_Z,a_{Z\dagger}\nbigo_Z,
 a_{Z\dagger}C_0)$
as the complex of $D$-triples.
Let us observe that
we have natural morphisms
$a_{Z\dagger}^0\bigl(\nbigu_Z[d_Z]\bigr)
\lrarr
 \nbigu_{\pt}$
and 
$\nbigu_{\pt}
\lrarr
 a_{Z\dagger}^0\bigl(
 \nbigu_Z[-d_Z]
 \bigr)$.
We have
\[
 a_{Z\dagger}C_0(\eta^{d_Z-p}m_1,\,
 \overline{\eta^{d_Z+p}m_2})
:=\left(
 \frac{1}{2\pi\sqrt{-1}}
 \right)^{d_Z}
 \int\eta^{d_Z-p}\etabar^{d_Z+p}
 C_0(m_1,\overline{m}_2)
 \,
 \epsilon(p+d_Z)
\]
In particular,
$a_{Z\dagger}^{d_Z}\nbigu_Z=
 \bigl(
 a_{Z\dagger}^{-d_Z}\nbigo_Z,\,
 a_{Z\dagger}^{d_Z}\nbigo_Z,\,C_1
 \bigr)$,
where 
\begin{equation}
 \label{eq;13.8.6.3}
 C_1(1,\etabar^{2d_Z})
\!\!=\!\!
 \left(\frac{1}{2\pi\sqrt{-1}}\right)^{d_Z}\!\!\!\!
 \int \etabar^{2d_Z}(-1)^{d_Z}
\!\!=\!\!\overline{
 \left(\frac{1}{2\pi\sqrt{-1}}\right)^{d_Z}
 \!\!\!\!\!
 \int \!\!\eta^{2d_Z}}\!
=C_0(1,\tr(\etabar^{2d_Z}))
\end{equation}
We also have
$a_{Z\dagger}^{-d_Z}\nbigu_Z=\bigl(
 a_{Z\dagger}^{d_Z}\nbigo_Z,\,
 a_{Z\dagger}^{-d_Z}\nbigo_Z,\,
 C_2
 \bigr)$,
where
\begin{equation}
 \label{eq;13.8.6.4}
 C_2(\eta^{2d_Z},1)
=\left(\frac{1}{2\pi\sqrt{-1}}\right)^{d_Z}
 \int \eta^{2d_Z}
=C_0\bigl(\tr(\eta^{2d_Z}),1\bigr)
\end{equation}
The equalities (\ref{eq;13.8.6.3})
and (\ref{eq;13.8.6.4}) mean that
we have the following natural morphisms:
\begin{equation}
 \label{eq;11.3.27.2}
 (a_Z^{-1},\tr):
 a_{Z\dagger}^0\bigl(
 \nbigu_Z[d_Z]
 \bigr)
\lrarr
 \nbigu_{\pt},
\quad\quad
 (\tr,a_Z^{-1}):
 \nbigu_{\pt}
\lrarr
 a_{Z\dagger}^{0}\bigl(
 \nbigu_Z[-d_Z]
 \bigr)
\end{equation}

The natural real structure of
$\nbigu_{\pt}$ is given by $(\id,\id)$.
\begin{prop}
\label{prop;11.3.27.11}
The morphisms in {\rm(\ref{eq;11.3.27.2})}
are compatible with the real structures.
\end{prop}
\pf
Let us consider the dual of
$\tr:a_{Z\dagger}^{d_Z}\nbigo_Z\lrarr\cnum$.
We consider $\Phi_1$ induced as follows:
\[
 \begin{CD}
 \bigl(
 a_{Z\dagger}^{0}(\nbigo_Z[d_Z])
 \bigr)^{\lor}
 @<{\tr^{\lor}}<<
 \cnum^{\lor} \\
 @V{B}V{\simeq}V @V{A}V{\simeq}V\\
 a_{Z\dagger}^{0}\bigl(\DDD(\nbigo_Z[d_Z])\bigr)
 @<{\Phi_1}<<
 \cnum
 \end{CD}
\]
Here, $A$ is induced by the natural multiplication
$\cnum\times\cnum\lrarr\cnum$,
and $B$ is given by the compatibility
of the dual and the push-forward.
We have the identification:
\begin{equation}
 \label{eq;11.3.27.3}
 a_{Z\dagger}^0\bigl(\DDD\nbigo_Z[d_Z]\bigr)
\simeq
 a_{Z\ast}^0\nhom_{D_Z}\bigl(\nbigo_Z[d_Z],
 \nbigo_Z[d_Z]\bigr)
\end{equation}
We have $\Phi_1(1)=\id$
under the identification (\ref{eq;11.3.27.3}).
Indeed, the identity is mapped to
the trace by the following morphisms:
{\small
\begin{multline}
 a_{Z\ast}^{0}\nrhom_{D_Z}(\nbigo_Z[d_Z],\nbigo_Z[d_Z])
\lrarr
 a_{Z\ast}^0\nrhom_{\cnum_Z}\bigl(
 \DR\nbigo_Z[d_Z],\DR\nbigo_Z[d_Z]
 \bigr) \\
\lrarr
 \bigl[
 \nrhom_{\cnum}\bigl(
 a_{Z\ast}\DR(\nbigo_Z)[d_Z],\,
 a_{Z\ast}\DR(\nbigo_Z)[d_Z]
 \bigr)
 \bigr]^{0}
\lrarr 
 \nrhom_{\cnum}\bigl(
 a_{Z\ast}^{d_Z}\DR(\nbigo_Z),\,\cnum
 \bigr)
\end{multline}
}
Hence, the following diagram is commutative,
by the construction of $\mu$:
\begin{equation}
 \begin{CD}
 a_{Z\dagger}^0(\nbigo_Z[d_Z])^{\lor}
 @>{B}>>
 a_{Z\dagger}^0\DDD\bigl(\nbigo_Z[d_Z]\bigr)
 @>{a_{Z\dagger}\mu}>>
 a_{Z\dagger}^0\nbigo_Z[-d_Z]
 \\
 @A{\tr^{\lor}}AA  @A{\Phi_1}AA @A{a_Z^{-1}}AA \\
 \cnum^{\lor} @>{A}>> \cnum @>{\id}>> \cnum
 \end{CD}
\end{equation}
Let us consider the dual of 
$a_Z^{-1}:\cnum\lrarr a_{Z\dagger}^0\bigl(\nbigo_Z[-d_Z]\bigr)$.
Let $\Phi_2$ be determined by 
the commutativity of the left square in the 
following diagram:
\begin{equation}
\label{eq;11.3.27.10}
 \begin{CD}
  a_{Z\dagger}^0\bigl(\nbigo_Z[-d_Z]\bigr)^{\lor}
 @>{B'}>>
  a_{Z\dagger}^0\bigl(\DDD(\nbigo_Z[-d_Z])\bigr)
 @>{a_{Z\dagger}\DDD\mu}>>
 a_{Z\dagger}^0\bigl(\nbigo_Z[d_Z]\bigr)
 \\
 @V{(a_Z^{-1})^{\lor}}VV @V{\Phi_2}VV @V{\tr}VV \\
 \cnum^{\lor} @>{A'}>>\cnum @>{=}>>\cnum
 \end{CD}
\end{equation}
Here, $A'$ is induced by the natural pairing
$\cnum\times\cnum\lrarr\cnum$,
and $B'$ is induced by the exchange of
the push-forward and the dual.
Let us observe that 
the right square in (\ref{eq;11.3.27.10})
is also commutative.
For that purpose,
we have only to show 
$(a_{Z\dagger}\mu\circ B)^{\lor}
=a_{Z\dagger}\DDD\mu\circ B'$
under the natural identification
$(a_{Z\dagger}^0\nbigo_Z[d_Z])^{\lor\lor}
\simeq
 a_{Z\dagger}^0\nbigo_Z[d_Z]$.
By a simple diagram chasing,
it can be reduced to the following:
\begin{itemize}
\item
$\mu\circ(\DDD\mu)^{-1}$
is the natural isomorphism
$\DDD\DDD\nbigo_Z[d_Z]\simeq \nbigo_Z[d_Z]$.
\item
The natural isomorphisms
$a_{Z\dagger}\simeq a_{Z\dagger}\DDD\circ\DDD$
and 
$a_{Z\dagger}\simeq \DDD\circ\DDD a_{Z\dagger}$
are compatible with
the identification
$a_{Z\dagger}\DDD\circ\DDD
\simeq
 \DDD\circ a_{Z\dagger}\circ\DDD
\simeq
 \DDD\circ\DDD a_{Z\dagger}$.
\end{itemize}
Then,
we obtain the claim of Proposition \ref{prop;11.3.27.11}.
\hfill\qed

\subsubsection{The real structure of the underlying perverse sheaves}

Let us look at the induced real structure
of $\cnum$-complex-triples.
We have
\[
 \Psi_1\nbigu_Z[d_Z]
=\bigl(
 \nbigo_Z[-d_Z],\,
 \nbigo_Z[d_Z],\,
 (-1)^{d_Z(d_Z+1)/2}C_0
 \bigr)
\]
The double complex
$\DR \Psi_1\nbigu_Z[d_Z]$
is given as follows:
we have
$\DR(\nbigo_Z[-d_Z])\simeq
 \cnum[0]$,
$\DR(\nbigo_Z[d_Z])\simeq
 \cnum[2d_Z]$
and
{\small
\begin{multline}
 \DR\bigl(
 (-1)^{d_Z(d_Z+1)/2}C_0
 \bigr)
 (\alpha\otimes\overline{\beta\,e})
=\alpha\betabar
 (-1)^{d_Z(d_Z+1)/2}
 \left(\frac{1}{2\pi\sqrt{-1}}\right)^{d_Z}
 (-1)^{d_Z(d_Z+1)/2}
 \\
=\alpha\overline{\left(
 \betabar\,\Bigl(\frac{1}{2\pi\sqrt{-1}}\Bigr)^{d_Z}
 \right)}\,(-1)^{d_Z}
\end{multline}
}
Then, we have
\[
 \Tot\DR\bigl(\Psi_1\nbigu_Z[d_Z]\bigr)
\simeq
 \bigl(\cnum[0],\cnum[2d_Z],C_1\bigr)
\]
where
$C_1\bigl(
 \alpha,\,\beta(2\pi\sqrt{-1})^{d_Z}e
 \bigr)=\alpha\betabar$.
We remark that $(-1)^{d_Z}$ appears
when we take the total complex.
The real structure
$\gamma^{\ast}\Tot\DR(\Psi_1\nbigu_Z[d_Z]) 
\simeq
 \Tot\DR(\Psi_1\nbigu_Z[d_Z]) $
is given by $(\id,\id)$.
It induces the real structure
$\real(2\pi\sqrt{-1})^{d_Z}\,e\subset
 \cnum[2d_Z]$.

\vspace{.1in}

We have
$\Psi_1\nbigu_Z[-d_Z]
=\bigl(
 \nbigo_Z[d_Z],\,\nbigo_Z[-d_Z],\,
 (-1)^{d_Z(d_Z-1)/2}C_0
 \bigr)$.
The double complex
$\DR\Psi_1\nbigu_Z[d_Z]$ is given as follows:
$\DR(\nbigo_Z[d_Z])
\simeq
 \cnum[2d_Z]$,
$\DR\bigl(\nbigo_Z[-d_Z]\bigr)
\simeq
 \cnum[0]$ and
{\small
\begin{multline}
 \DR\bigl(
 (-1)^{d_Z(d_Z-1)/2}C_0
 \bigr)(\alpha e\otimes\betabar)
=\alpha\betabar\cdot
 (-1)^{d_Z(d_Z-1)/2}
 \left(\frac{1}{2\pi\sqrt{-1}}\right)^{d_Z}
 (-1)^{d_Z(d_Z+1)/2}
 \\
=\alpha\betabar
 \left(\frac{1}{2\pi\sqrt{-1}}\right)^{d_Z}
 (-1)^{d_Z}
\end{multline}
}
Then, we obtain 
\[
 \Tot\DR(\Psi_1\nbigu_Z[-d_Z])
\simeq 
 \Bigl(
 \cnum[2d_Z],\,\cnum[0],\,C_2
 \Bigr),
\quad\quad
 C_2\bigl(
 \alpha(2\pi\sqrt{-1})^{d_Z}e,\,\beta
 \bigr)
=\alpha\betabar.
\]
The real structure
$\gamma^{\ast}\Tot\DR(\Psi_1\nbigu_Z[-d_Z])
\simeq
 \Tot\DR(\Psi_1\nbigu_Z[-d_Z])$
is given by $(\id,\id)$.
It gives the real structure
$\real\cdot 1\subset \cnum[0]$.

\vspace{.1in}

Let us directly compute the real structures of
the underlying perverse sheaves of
$\nbigu_Z$ associated to
\[
\kappa_1=\bigl(\epsilon(-d_Z)\nu^{-1},\epsilon(d_Z)\nu\bigr),
\quad\quad
\kappa_2=\bigl(\epsilon(d_Z)\nu^{-1},\epsilon(-d_Z)\nu\bigr),
\]
although we can describe them as the shift.
Note that under
$\cnum[2d_Z]\simeq \distribution_X^{\bullet}[2d_Z]$
with the twisted trace morphisms,
the real part of
$\cnum[2d_Z]$ is given by
$(2\pi\sqrt{-1})^{d_Z}\real[2d_Z]\subset
 \cnum[2d_Z]$.
It gives the identification
\[
 H^{-d_Z}\Hom\bigl(\overline{\cnum[d_Z]},\cnum[2d_Z]\bigr)
\simeq
 \overline{H^{-d_Z}\Hom\bigl(\cnum[d_Z],\cnum[2d_Z]\bigr)}.
\]
Let $f\in H^{-d_Z}\Hom\bigl(\cnum[d_Z],\cnum[2d_Z]\bigr)$
be the identity.
It sends $(2\pi\sqrt{-1})^{d_Z}[d_Z]$
to $(2\pi\sqrt{-1})^{d_Z}[2d_Z]$.

The pairing and the first component of $\kappa_1$
maps $f$ to
\[
 \fbar\in
 \overline{H^{-d_Z}\Hom\bigl(\cnum[d_Z],\cnum[2d_Z]\bigr)}.
\]
Hence, the real structure of
$H^{-d_Z}\Hom\bigl(\cnum[d_Z],\cnum[2d_Z]\bigr)$
is $\real\,f$,
and the real structure of
$\cnum[d_Z]$ is given by
$(2\pi\sqrt{-1})^{d_Z}\real[d_Z]
\subset\cnum[d_Z]$.
The pairing of the first component of $\kappa_2$
maps to $f$ to $(-1)^{d_Z}f$.
Hence, the real structure of 
$H^{-d_Z}\Hom\bigl(\cnum[d_Z],\cnum[2d_Z]\bigr)$
is $\real(2\pi\sqrt{-1})^{d_Z}f$,
and the real structure of 
$\cnum[d_Z]$ is given by
$\real[d_Z]\subset\cnum[d_Z]$.

\subsubsection{Appendix}

We use the identification
$\DDD\nbigo_Z\simeq \nbigo_Z$
given as follows:
\[
 D_X\otimes\Theta^{\bullet}_X
\lrarr
 \nhom_{D_X}\Bigl(
 D_X\otimes\Theta_X^{\bullet},\,
 \bigl(D_X\otimes\omega_X^{-1}[d_X]\bigr)^{\ell,r}
 \Bigr)
\]
\[
 P\otimes\tau
\longmapsto
 \Bigl(
 Q\otimes\omega
\longmapsto
 (-1)^{|\tau|\,|\omega|}
 C(Q\otimes\omega,P\otimes\tau)
 \Bigr)
\]
Here,
$C(Q\otimes\omega,P\otimes\tau)
=\ell(Q)r(P)(\omega\wedge\tau)\,(-1)^{|\tau|}$.
Note that the induced morphism
$\nbigo_X\lrarr
 \nhom_{D_X}\bigl(
 D_X\otimes\omega_X^{-1},
 D_X\otimes\omega_X^{-1}
 \bigr)$
sends $1$ to the identity.
By taking tensor product 
$\omega_X\otimes_{D_X}$,
we obtain
\[
 \omega_X\otimes \Theta_X^{\bullet}
\lrarr
 \nhom_{D_X}\bigl(
 D_X\otimes\Theta_X^{\bullet},
 \nbigo_X[d_X]
 \bigr)
\lrarr
 \nhom_{D_X}\bigl(\nbigo_X,\nbigo_X[d_X]\bigr)
\]
Let $\eta$ be a local generator of $\omega_X$.
Then, $\eta\otimes\eta^{-1}$
is sent to $(-1)^{d_X}$.
Hence,
the induced isomorphism
$\cnum[d_X]\lrarr
\nhom\bigl(
 \cnum[d_X],\cnum[2d_X]
 \bigr)$
maps $1[d_X]$
to $(-1)^{d_X}f$.

\chapter{Good systems of ramified irregular values}

In \S\ref{section;13.5.6.100},
we prepare some terminology 
which are convenient to control the irregularity.
In \S\ref{section;13.5.13.10},
we study a resolution of turning points
for Lagrangian covers.
It simplifies and clarifies
the definition of wild harmonic bundles.

\section{Good system of ramified irregular values}
\label{section;13.5.6.100}

\subsection{Good set of irregular values}
\label{subsection;13.5.12.40}
\index{good set of irregular values}

Let $X$ be a complex manifold
with a normal crossing hypersurface $D$.
Let $P$ be any point of $D$.
We take a holomorphic coordinate neighbourhood
$(X_P;z_1,\ldots,z_n)$ around $P$
such that 
$X_P\cap D=\bigcup_{i=1}^{\ell}\{z_i=0\}$.
For any $\vecm=(m_i)\in\seisuu^{n}$,
we set $\vecz^{\vecm}=\prod_{i=1}^nz_i^{m_i}$.
We naturally regard $\seisuu^{\ell}\subset\seisuu^{n}$
by $\vecm\longmapsto(\vecm,0,\ldots,0)$.
A finite subset 
$\nbigi\subset\nbigo_X(\ast D)_P/\nbigo_{X,P}$
is called a good set of irregular values at $P$,
if the following holds
(\cite{mochi7}, \cite{sabbah4}):
\begin{itemize}
\item
 For any $\gminia\in\nbigi\setminus\{0\}$,
 we take a lift $\gminiatilde\in\nbigo_{X}(\ast D)_P$.
 Then, there exists 
 $\ord(\gminia)\in \seisuu^{\ell}_{\leq 0}$
 such that
 $\gminiatilde\vecz^{-\ord(\gminia)}$ is 
 an invertible element of $\nbigo_{X,P}$.
\item
 For any pair $\gminia_i\in\nbigi$ $(i=1,2)$
 with $\gminia_1\neq\gminia_2$,
 we takes lift $\gminiatilde_i\in\nbigo_X(\ast D)_P$.
 Then, there exists
 $\ord(\gminia_1-\gminia_2)\in\seisuu^{\ell}_{\leq 0}$
 such that
 $(\gminiatilde_1-\gminiatilde_2)
 \vecz^{-\ord(\gminia_1-\gminia_2)}$
 is an invertible element of $\nbigo_{X,P}$.
\item
 For any such pairs 
 $(\gminia_1,\gminia_2)$ and
 $(\gminib_1,\gminib_2)$
 in $\nbigi$, 
 either
 $\vecz^{\ord(\gminia_1-\gminia_2)-\ord(\gminib_1-\gminib_2)}
 \in \nbigo_{X,P}$
 or 
 $\vecz^{-\ord(\gminia_1-\gminia_2)+\ord(\gminib_1-\gminib_2)}
 \in \nbigo_{X,P}$
 holds.
\end{itemize} 
 The condition is independent of the choices of lifts.

\subsection{Good system of ramified irregular values}

Let $X$ be a complex manifold
with a normal crossing hypersurface $D$.
We generalize the notion of good system of irregular values
in \S2.4.1 of \cite{mochi7}.

Let $P$ be any point of $D$.
We introduce a category $C_P(X,D)$ as follows.
Objects in $C_P(X,D)$ are holomorphic maps
$\varphi:(Z,Q)\lrarr (X,P)$
of smooth complex manifolds
which are coverings with ramification along $D$.
We set $D_Z:=\varphi^{-1}(D)$.
Morphisms 
$F:\bigl(Z,Q,\varphi\bigr)
\lrarr
 \bigl(Z',Q',\varphi'\bigr)$
are holomorphic maps $F:(Z,Q)\lrarr (Z',Q')$
such that $\varphi'\circ F=\varphi$.
Such morphisms induce the morphisms
$\nbigo_{Z'}(\ast D_{Z'})_{Q'}
\lrarr
\nbigo_{Z}(\ast D_Z)_{Q}$
over $\nbigo_X(\ast D)_Q$.
Let $\nbigotilde_{X}(\ast D)_P$
denote an inductive limit of
$\nbigo_{Z}(\ast D_Z)_Q$.
Similarly,
let $\nbigotilde_{X,P}$
denote the inductive limit of
$\nbigo_{Z,Q}$.
\index{ring $\nbigotilde_{X}(\ast D)_P$}
\index{ring $\nbigotilde_{X,P}$}

We have another description of the rings.
Let $\cnum\{z_1,\ldots,z_n\}$
denote the ring of convergent power series.
Let $\cnum\{z_1,\ldots,z_n\}_{z_1\ldots z_{\ell}}$
denote its localization with respect to
$z_1\cdots z_{\ell}$.
For a coordinate $(z_1,\ldots,z_n)$ such that
$D=\bigcup_{i=1}^{\ell}\{z_i=0\}$,
we have natural isomorphisms
\[
\nbigotilde_{X,P}\simeq
\varinjlim_e
\cnum\bigl\{z_1^{1/e},\ldots,z_{\ell}^{1/e},z_{\ell+1},\ldots,z_n\bigr\},
\]
\[
 \nbigotilde_{X}(\ast D)_P
 \simeq
\varinjlim_e
\cnum\bigl\{z_1^{1/e},\ldots,z_{\ell}^{1/e},z_{\ell+1},\ldots,z_n\bigr\}
 _{z_1^{1/e}\cdots z_{\ell}^{1/e}}.
\]

\begin{df}
\label{df;13.5.7.201}
A finite subset 
$\nbigi\subset
 \nbigotilde_X(\ast D)_P\big/\nbigotilde_{X,P}$
can be regarded as
$\nbigi\subset\nbigo_{Z}(\ast D_Z)_Q
 \big/\nbigo_{Z,Q}$
for some $\bigl((Z,Q),\varphi\bigr)\in C_P(X,D)$.
It is called a good set of ramified irregular values at $P$,
if (i) it is a good set of irregular values on small $(Z,D_Z)$,
(ii) it is stable under the action of 
the Galois group of $\varphi$.
\index{good set of ramified irregular values}
\hfill\qed
\end{df}
Note that,
 if $P_1$ is close to $P$,
 we choose $Q_1\in \varphi^{-1}(P_1)$,
 and we obtain a natural map
 $\nbigi_P\lrarr
 \nbigo_{Z}(\ast D_Z)_{Q_1}
 \big/\nbigo_{Z,Q_1}
\lrarr 
 \nbigotilde_{X}(\ast D)_{P_1}
 \big/\nbigotilde_{X,P_1}$.
The image is well defined,
and it gives a good set of ramified irregular
values in
 $\nbigotilde_{X}(\ast D)_{P_1}
 \big/\nbigotilde_{X,P_1}$.

\begin{df}
\label{df;13.5.7.100}
A good system of ramified irregular values
on $(X,D)$ is a family of good sets of 
ramified irregular values
$\vecnbigi=\bigl\{\nbigi_P\,\big|\,P\in D
 \bigr\}$ satisfying the following condition.
\begin{description}
\item[\bf (A)]
 If $P_1$ is sufficiently close to $P$,
 the image of
 $\nbigi_P$ in 
$\nbigotilde_{X}(\ast D)_{P_1}
 \big/\nbigotilde_{X,P_1}$
 is equal to
 $\nbigi_{P_1}$.
\hfill\qed
\end{description}
\end{df}
\index{good system of ramified irregular values}

The easiest example of good system of
ramified irregular values is given by
setting $\nbigi_P=0$ for any $P\in D$.
It is called the trivial good system of ramified irregular values,
and denoted by $\veczero$.

\begin{rem}
Let $\vecnbigi=(\nbigi_P\,|\,P\in D)$
be a good system of ramified irregular values
on $(X,D)$.
\begin{itemize}
\item
We set
$-\vecnbigi:=(-\nbigi_P\,|\,P\in D)$,
where
$-\nbigi_P:=\bigl\{
 -\gminia\,\big|\,
 \gminia\in\nbigi_P
 \bigr\}$,
which is a good system of irregular values on $(X,D)$.
\item
Let $F:X'\lrarr X$ be a morphism of complex manifolds
such that
$D':=F^{-1}(D)$ is normal crossing.
For any $P'\in D'$, we have a naturally defined morphism
$\nbigotilde_{X}(\ast D)_{F(P')} \big/
\nbigotilde_{X,F(P')}
 \lrarr 
\nbigotilde_{X'}(\ast D')_{P'}\big/
\nbigotilde_{X',P'}$.
The image of $\nbigi_{F(P')}$ is denoted by
$F^{-1}(\nbigi)_{P'}$,
which is a good set of irregular values at $P'$.
Thus, we obtain a good system of irregular values
$F^{-1}(\vecnbigi)$ on $(X',D')$.
\item
In particular,
for any open subset $X'\subset X$,
we have the naturally defined good system 
$\vecnbigi_{|X'}:=(\nbigi_P\,|\,P\in X'\cap D)$
on $(X',D\cap X')$.
\hfill\qed
\end{itemize}
\end{rem}

\subsection{Specialization of good set of ramified irregular values}
\label{subsection;13.5.6.30}
\index{specialization}

Let $X$ and $D$ be as above.
Let $D=\bigcup_{i\in \Lambda}D_i$
be the irreducible decomposition.,
Let $I\subset \Lambda$ be any subset.
We set $I^c:=\Lambda\setminus I$.
We put $D_I:=\bigcap_{i\in I}D_i$,
$D(I):=\bigcup_{i\in I}D_i$,
$\del D_I:=D_I\cap D(I^c)$
and $D_I^{\circ}:=D_I\setminus \del D_I$.
For any $P\in D_I$,
we have the naturally defined maps:
\[
 \varphi_{I1}:
 \nbigotilde_{X}(\ast D)_P\big/
 \nbigotilde_{X,P}
\lrarr
 \nbigotilde_{X}(\ast D)_P\big/
 \nbigotilde_{X}(\ast D(I^c))_P
\]
\[
 \varphi_{I2}:
 \nbigotilde_{X}(\ast D(I^c))_P\big/\nbigotilde_{X,P}
\lrarr
 \nbigotilde_{D_I}(\ast \del D_I)_P\big/
 \nbigotilde_{D_I,P}
\]
For any subset 
$\nbigj_P\subset
 \nbigotilde_{X}(\ast D)_P\big/
 \nbigotilde_{X,P}$,
we set
$\nbigj_P(I)_{|D_I}:=
 \varphi_{I2}
 \bigl(
 \varphi_{I1}^{-1}(0)
 \bigr)$.

Let $\vecnbigi$ be a good system of ramified irregular values
on $(X,D)$.
By applying the above procedure to any $P\in D_I$,
we obtain a good system of ramified irregular values
on $(D_I,\del D_I)$,
denoted by
$\vecnbigi(I)_{|D_I}$.

\vspace{.1in}
In the local case,
we may also use the following procedure.
Let $\nbigj$ be any finite subset of
$\nbigotilde_X(\ast D)_P/\nbigotilde_{X,P}$.
For any $\gminia\in\nbigo_{X}(\ast D)_P/\nbigo_{X,P}$,
we define
$\nbigj(-\gminia):=\bigl\{
 \gminib-\gminia\,\big|\,
 \gminib\in\nbigj
 \bigr\}$.
We obtain the set
$\nbigj(-\gminia,I)_{|D_I}:=
 \nbigj(-\gminia)(I)_{|D_I}
\subset
 \nbigo_{D_I}(\ast \del D_I)_P\big/
 \nbigo_{D_I,P}$
in a similar way.
If $\nbigj\subset\nbigo_X(\ast D)_P/\nbigo_X$
is a good set of irregular values,
$\nbigj(-\gminia,I)_{|D_I}$
is a good set of irregular values
for any $\gminia\in\nbigj$.

\subsection{Resolution}
\index{resolution for good set of irregular values}

Let $X$ be a complex manifold with 
a simply normal crossing hypersurface $D$.
For simplicity,
we assume that the number of the irreducible components of $D$
is finite,
and that the numbers of the connected components of
$D_I^{\circ}$ are finite.
Let $\vecnbigi=(\nbigi_P\,|\,P\in D)$
be a family of finite subsets
$\nbigi_P\subset\nbigotilde_{X}(\ast D)_P\big/\nbigotilde_{X,P}$
such that 
(i) each $\nbigi_P$ is Galois invariant,
(ii) it satisfies the condition {\bf(A)}
in Definition \ref{df;13.5.7.100}.
The following proposition is essentially proved in \cite{mochi7}.
We will give a proof in 
\S\ref{subsection;13.5.10.50}--\S\ref{subsection;13.5.10.51}.

\begin{prop}
\label{prop;13.5.7.110}
There exists a projective morphism of complex manifolds
$F:X'\lrarr X$ such that
(i) $D':=F^{-1}(D)$ is simply normal crossing,
(ii) $X'\setminus D'\simeq X\setminus D$,
(iii) $F^{-1}\vecnbigi$ is a good set of irregular values
 on $(X',D')$.
\end{prop}

\begin{rem}
Proposition {\rm\ref{prop;13.5.7.110}}
shall be used typically as follows.
Let $\vecgbigi_i$ $(i=1,2)$
be good system of ramified irregular values.
We set
$(\gbigi_1\oplus\gbigi_2)_P
:=\gbigi_{1P}\cup\gbigi_{2P}$
and 
$(\gbigi_1\otimes\gbigi_2)_P
:=\bigl\{
 \gminia_1+\gminia_2\,\big|\,
 \gminia_i\in\gbigi_{iP}
 \bigr\}$.
Then, the systems
$\vecgbigi_1\oplus\vecgbigi_2$
and 
$\vecgbigi_1\otimes\vecgbigi_2$
are not good systems of ramified irregular values,
in general.
We use Proposition {\rm\ref{prop;13.5.7.110}}
to resolve the points at which 
they are not good.
\hfill\qed
\end{rem}

\subsubsection{System of tuples of ramified ideals}
\label{subsection;13.5.10.50}

Let $X$ and $D$ be as above.
Let $D=\bigcup_{i\in\Lambda}D_i$
be the irreducible decomposition.
We consider a family of finite tuples of
 finitely generated ideals 
 $\gbigi_P=(\gbigi_{P,1},\ldots,\gbigi_{P,s(P)})$
 of $\nbigotilde_{X,P}$
for $P\in D$.
We assume that $\{s(P)\}$ is bounded.
We impose the following condition:
\begin{itemize}
\item
 Each $\gbigi_P$ is Galois invariant.
\item
 If $P_1$ is sufficiently close to $P$,
 we have
\[
  \{\gbigi_{P,1}^{(P_1)},\ldots,\gbigi_{P,s(P)}^{(P_1)},
 \nbigotilde_{X,P_1}
 \}
 =\{\gbigi_{P_1,1},\ldots,\gbigi_{P_1,s(P_1)},
 \nbigotilde_{X,P_1}
 \}
\]
 as sets.
 Here, $\gbigi_{Pi}^{(P_1)}$ denotes
 the ideal in $\nbigotilde_{X,P_1}$
 induced by $\gbigi_{Pi}$.
\end{itemize}
Such a family
$\vecgbigi=\bigl\{
 \gbigi_P\,|\,P\in D
 \bigr\}$
is called a system of tuples of ramified ideals
on $(X,D)$.
It is called principal,
if the following holds:
\begin{itemize}
\item
 For any $P$,
 $\gbigi_{Pi}$ are principal ideals,
 and 
 $\{\gbigi_{Pi}\}$ is totally ordered
 with respect to the inclusion.
\end{itemize}
The following lemma is clear.
\begin{lem}
There exists
$(e_i\,|\,i\in\Lambda)\in\seisuu_{>0}^{\Lambda}$
such that the following holds:
\begin{itemize}
\item
Let $P$ be any point of $D_I$.
We take a holomorphic coordinate
$(z_1,\ldots,z_n)$ around $P$
 such that $D_j$ $(j\in I)$ are expressed as
 $\{z_{\alpha(j)}=0\}$.
 Then, 
 $\gbigi_{Pi}$ are contained in the extension
 $\nbigo_{X,P}[z_{\alpha(j)}^{1/e_j}\,|\,j\in I]$.
\hfill\qed
\end{itemize}
\end{lem}

\begin{lem}
\label{lem;13.5.7.200}
Let $\vecgbigi$ be any family of 
tuples of ramified ideals on $(X,D)$.
There exists a projective morphism
$F:X'\lrarr X$ such that
(i) $D':=F^{-1}(D)$ is simply normal crossing,
(ii) $X'\setminus D'\simeq X\setminus D$,
(iii) $F^{-1}\vecgbigi=\bigl(
 F^{-1}(\gbigi_{F(P')})\nbigotilde_{\Xtilde'}(\ast D')_{P'}
 \,\big|\,
 P'\in D'
 \bigr)$ is principal.
\end{lem}
\pf
For each $P\in D$,
we take a small neighbourhood $X_P$ of $P$,
and a ramified covering
$\varphi_P:Z_P\lrarr X_P$
whose ramification indexes along $D_i$
are $e_i$.
We have the tuple of the ideals
$\vecgbigi_{Z_P}:=
 \bigl(
 \gbigi_{Z_P,i}\,\big|\,
 i=1,\ldots,s(P) \bigr)$
which induces
$\vecgbigi_P$,
and it is invariant under the action
of the Galois group $G(\varphi_P)$
of $\varphi_P$.
By applying the construction
in \S15.1 of \cite{mochi7},
we canonically obtain a normal complex variety $X_P'$
equipped with a $G(\varphi)$-action,
and a $G(\varphi)$-equivariant
projective birational morphism $G_P:Z_P'\lrarr Z_P$
such that
(i) $Z_P'\setminus(\varphi_P\circ G_P)^{-1}(D)
\simeq
 Z_P\setminus \varphi_P^{-1}(D)$,
(ii) each ideal $G_P^{-1}(\gbigi_{Z_P,i})\nbigo_{Z'_P}$ 
is principal,
(iii) the ideals are totally ordered with respect to
 the inclusion.
In particular, we obtain $X^{(1)}_P=Z_P'/G(\varphi)$
with an induced birational morphism 
$F^{(1)}_P:X^{(1)}_P\lrarr X_P$
which satisfies 
$X^{(1)}_P\setminus (F^{(1)}_P)^{-1}(D)
\simeq X_P\setminus D$.
We can glue $(X^{(1)}_P,F^{1}_P)$ $(P\in D)$,
and we obtain a complex variety
$X^{(1)}$ with a morphism
$F^{(1)}:X^{(1)}\lrarr X$.
Then, by taking an appropriate projective birational morphism
$X'\lrarr X^{(1)}$,
we obtain the desired one.
\hfill\qed

\subsubsection{Proof of Proposition \ref{prop;13.5.7.110}}
\label{subsection;13.5.10.51}

The following lemma is clear.
\begin{lem}
There exists
$(e_i\,|\,i\in\Lambda)\in\seisuu_{>0}^{\Lambda}$
such that the following holds:
\begin{itemize}
\item
Let $P$ be any point of $D_I$.
We take a holomorphic coordinate
$(z_1,\ldots,z_n)$ around $P$
 such that $D_j$ $(j\in I)$ are expressed as
 $\{z_{\alpha(j)}=0\}$.
 Then, 
 $\nbigi_P$ is contained in 
 $\nbigo_{X}(\ast D)_P[z_{\alpha(j)}^{1/e_j}\,|\,j\in I]
 \big/\nbigo_{X,P}[z_{\alpha(j)}^{1/e_j}\,|\,j\in I]$.
\hfill\qed
\end{itemize}
\end{lem}

For each $P\in D$,
we take a small neighbourhood $X_P$
and a ramified covering 
$\varphi_P:(Z_P,Q)\lrarr (X_P,P)$
whose ramification indexes along $D_i$ are $e_i$.
We set $D_{Z_P}:=\varphi_P^{-1}(D)$.
We have
$\nbigi_P\subset
 \nbigo_{Z_P}(\ast D_{Z_P})_Q\big/
 \nbigo_{Z_P,Q}$.
For each $\gminia\in\nbigi_P$,
we take lift $\gminiatilde
\in\nbigo_{Z_P}(\ast D_{Z_P})_Q$.
We take a coordinate 
$(\xi_1,\ldots,\xi_n)$ of $Z_P$
such that
$D_{Z_P}=\bigcup_{i=1}^{\ell}\{\xi_i=0\}$.
For $\gminia\in\nbigi_P$,
let $m_i(\gminia)\geq 0$ be the pole order
of $\gminia$ along $\{\xi_i=0\}$.
We set
$\vecxi^{\vecm(\gminia)}:=
 \prod_{i=1}^{\ell}\xi_i^{m_i(\gminia)}$.
We consider the ideal
$\gbigi_{P,\gminia}$
generated by
$\vecxi^{\vecm(\gminia)}$
and 
$\gminiatilde\vecxi^{\vecm(\gminia)}$.
It is independent of the choice of $\gminiatilde$.
Note that
the ideal $\gbigi_{P,\gminia}$ is principal
means either one of the following holds;
(i) $\gminia=0$ in
 $\nbigotilde_{X}(\ast D)_P\big/\nbigotilde_{X,P}$,
(ii) there exists 
 $\ord(\gminiatilde)$
 in $\seisuu_{\leq 0}^{\ell}$.

They induce a tuple of finitely generated ideals
$\vecgbigi_P=(\gbigi_{P,\gminia})$
in $\nbigotilde_{X,P}$,
which is Galois invariant.
Note that $\vecgbigi_P$ is principal
implies that,
for each pair $(\gminia,\gminib)$ in $\nbigi_P$,
we have either
$\vecxi^{\vecm(\gminia)-\vecm(\gminib)}
 \in\nbigo_{Z_P,Q}$
or 
$\vecxi^{-\vecm(\gminia)+\vecm(\gminib)}
 \in\nbigo_{Z_P,Q}$.
(We may use it below to deduce the condition (iii)
in Definition \ref{df;13.5.7.201}.)

The family
$\vecgbigi=(\vecgbigi_P\,|\,P\in D)$
is a system of tuples of ramified ideals.
By applying Lemma \ref{lem;13.5.7.200},
we obtain 
a projective birational morphism
$F_1:X_1\lrarr X$
such that
(i) $D_1:=F_1^{-1}(D)$ is simply normal crossing,
(ii) $X_1\setminus D_1\simeq X\setminus D$,
(iii) $F^{-1}\nbigi$  satisfies the condition 
(i) in Definition \ref{df;13.5.7.201}.

By applying similar arguments to tuples
$\nbigj_P:=
 \{\gminia-\gminib\,|\,\gminia,\gminib\in\nbigi_P\}$,
we obtain a projective birational morphism
$F_2:X_2\lrarr X$
such that
(i) $D_2:=F_2^{-1}(D)$ is simply normal crossing,
(ii) $X_2\setminus D_2\simeq X\setminus D$,
(iii) $F^{-1}\nbigi$  satisfies the condition 
(i) and (ii) in Definition \ref{df;13.5.7.201}.
The condition (iii) in Definition \ref{df;13.5.7.201}
is also satisfied.
Note that the above remark.
\hfill\qed

\section{Resolution of turning points for 
Lagrangian covers}
\label{section;13.5.13.10}

\subsection{Statement}
\label{subsection;13.5.12.20}
Let $X$ be an $n$-dimensional complex manifold 
with a normal crossing hypersurface $H$.
We consider a reduced complex analytic closed subset
$\Sigma\subset T^{\ast}(X\setminus H)$
such that 
(i) it is finite over $X\setminus H$,
(ii) the smooth part of $\Sigma$ is Lagrangian
with respect to the natural symplectic structure of
$T^{\ast}(X\setminus H)$.
We call such $\Sigma$ by a Lagrangian cover of $X\setminus H$.
If we are given any section $\omega\in\Omega^1_{X\setminus H}$,
then $\omega+\Sigma\subset T^{\ast}(X\setminus H)$ denotes
the translation of $\Sigma$ by $\omega$.
\index{Lagrangian cover}

Let $H=\bigcup_{j\in \Lambda}H_j$ be 
the irreducible decomposition.
For any $\vecN\in\seisuu^{\Lambda}$,
let $\vecN H$ denote the divisor $\sum N_iH_i$.
For simplicity,
we assume that
$\Lambda$ is finite,
and that the connected components of
$H_I^{\circ}$ are finite
for any $I\subset\Lambda$.
\begin{df}
A Lagrangian cover $\Sigma$ is called meromorphic, if
the closure of $\Sigma$ in 
$T^{\ast}X(\log D)\otimes \nbigo_X(\vecN H)$
is a complex analytic subset which is finite over $X$
for some $\vecN\in\seisuu_{>0}^{\Lambda}$.
\hfill\qed
\end{df}
\index{meromorphic Lagrangian cover}

We introduce some conditions
for the behaviour of meromorphic 
Lagrangian cover around $H$.

\begin{df}
A meromorphic Lagrangian cover 
is called logarithmic,
if the closure in $T^{\ast}X(\log H)$
is finite over $X$.
\hfill\qed
\end{df}
\index{logarithmic Lagrangian cover}

For any $P\in X$,
let $X_P\subset X$ denote a small neighbourhood of $P$.
We set $H_P:=H\cap X_P$ in that situation.

\begin{df}
A meromorphic Lagrangian cover $\Sigma$
is called unramifiedly good at $P\in H$,
if there exist a good set of irregular values
$\nbigi_P\subset M(X_P,H_P)/H(X_P)$
and logarithmic 
$\Sigma_{P,\gminia}\subset
 T^{\ast}(X_P\setminus H_P)$ $(\gminia\in \nbigi_P)$
such that 
\begin{equation}
 \label{eq;13.5.13.1}
 \Sigma_{|X_P\setminus H_P}
=
 \bigsqcup_{\gminia\in\nbigi_P}
 \bigl(
 d\gminiatilde+\Sigma_{P,\gminia}
 \bigr).
\end{equation}
Here, $\gminiatilde\in M(X_P,H_P)$ are 
lifts of $\gminia$.
The meromorphic Lagrangian cover
is called unramifiedly good on $(X,H)$,
if it is unramifiedly good at any $P\in H$.
\hfill\qed
\end{df}
\index{unramifiedly good Lagrangian cover}

\begin{df}
A meromorphic Lagrangian cover $\Sigma$
is called good at $P\in H$,
if there exists a ramified covering
$\psi_P:(X_P',H_P')\lrarr (X_P,H_P)$
such that
$\psi_P^{-1}(\Sigma)$ is unramifiedly good 
on $(X_P',H_P')$.
It is called good,
if it is good at any $P\in H$.
\hfill\qed
\end{df}
\index{good Lagrangian cover}

We will prove the following theorem.

\begin{thm}
\label{thm;13.5.11.1}
For any meromorphic Lagrangian cover $\Sigma$,
there exists a projective morphism
$F:X'\lrarr X$ such that
(i) $H':=F^{-1}(H)$ is normal crossing,
(ii) $X'\setminus H'\simeq X\setminus H$,
(iii) $F^{-1}(\Sigma)$ is good on $(X',H')$.
\end{thm}
\index{resolution for Lagrangian cover}

\subsubsection{Wild harmonic bundle}

Before going to the proof,
we give consequences on wildness of harmonic bundles.
Let $X$ be any complex manifold
with a normal crossing hypersurface 
$H=\bigcup_{i\in\Lambda}H_i$.
For simplicity,
we assume that
$\Lambda$ is finite,
and that the connected components of $H_I^{\circ}$ 
are finite for any $I\subset\Lambda$.

Let $(E,\delbar_E,\theta,h)$ be a harmonic bundle
on $X\setminus H$.
Let $\Sigma(\theta)\subset T^{\ast}(X\setminus H)$
denote the spectral variety of $\theta$.
It is easy to observe that $\Sigma(\theta)$ is 
a Lagrangian cover of $X\setminus H$, 
by using Gabber's theorem.
(See \S A:III.3 of \cite{bjork}, for example.)
The following corollary is
an easy characterization of wildness of
harmonic bundles.

\begin{cor}
\label{cor;13.5.11.2}
There exists a projective morphism
$F:X'\lrarr X$ such that
(i) $H':=F^{-1}(H)$ is simply normal crossing,
(ii) $X'\setminus H'\simeq X\setminus H$,
(iii) $F^{-1}\harmonicbundle$
be a good wild harmonic bundle on $(X',H')$,
if and only if
$\Sigma(\theta)$ is meromorphic.
\end{cor}
\pf
The only if part is clear.
The if part follows from 
Theorem \ref{thm;13.5.11.1}.
\hfill\qed

\vspace{.1in}
We call $F$ satisfying the condition 
in Corollary {\rm\ref{cor;13.5.11.2}}
by a resolution for 
the wild harmonic bundle
$\harmonicbundle$ on $(X,H)$.
The following corollary means that
we do not have to care with the choice
of an open subset.

\begin{cor}
\label{cor;13.5.13.31}
There exists a resolution 
for $\harmonicbundle$ on $(X,H)$,
if and only if
there exists a resolution 
for $\harmonicbundle$ on $(X,H_1)$
for some $H_1\supset H$.
\hfill\qed
\end{cor}

\begin{cor}
\label{cor;13.5.13.30}
There exists a resolution for 
$\harmonicbundle$ on $(X,H)$,
if and only if
$(E,\delbar_E,\theta)$
is extended to a meromorphic Higgs sheaf
on $(X,H)$.
\end{cor}
\pf
If $(E,\delbar_E,\theta)$ is extended to
a meromorphic Higgs sheaf,
$\Sigma(\theta)$ is meromorphic.
Hence, there exists a resolution.
Conversely, if we have a resolution
$F:(X',H')\lrarr (X,H)$,
According to \cite{mochi7},
$F^{-1}(E,\delbar_E,\theta)$
is extended to a good filtered Higgs bundle
on $(X',H')$.
By taking the push-forward,
$(E,\delbar_E,\theta)$ is extended
to a meromorphic Higgs sheaf on $(X,H)$.
\hfill\qed

\subsection{Separation of ramification and polar part}

\subsubsection{Separation of zeroes and poles}

Let $Z$ be any irreducible normal complex analytic space
with a hypersurface $H=\bigcup_{i\in\Lambda} H_i$.
Suppose that $(Z,H)$ is equipped with 
an action of a finite group $G$.
Let $\nbige$ be any $G$-equivariant
locally free $\nbigo_Z$-module on $Z$.
Let $f_1,\ldots,f_m$ be sections of 
$\nbige(\vecN H)$ 
for some $\vecN\in\seisuu_{>0}^{\Lambda}$,
such that the $G$-action induces a permutation on
$\{f_1,\ldots,f_m\}$.
Let $\gbigi_{f_i}$ denote the ideal sheaf
of the $0$-set of $f_i$ regarded as a section of
$\nbige(\vecN H)$.
We set $\gbigi^{(0)}_{f_i}:=
 \gbigi_{f_i}+\nbigo(-\vecN H)$
in $\nbigo_Z$.
By applying the construction in \S15.1 of \cite{mochi7},
we obtain the following lemma.

\begin{lem}
\label{lem;13.5.12.10}
There exist an irreducible normal complex
analytic space $Z^{(1)}$ with a $G$-action,
and a $G$-equivariant projective
morphism $F:Z^{(1)}\lrarr Z$
such that the following holds:
\begin{itemize}
\item
$Z^{(1)}\setminus F^{-1}(H)\simeq Z\setminus H$.
(It implies that there exists a closed analytic subset
$A\subset H$ such that
$Z^{(1)}\setminus F^{-1}(A)\simeq
 Z\setminus A$.)
\item
 The ideals
 $\gbigi^{(1)}_{f_i}:=
 \nbigo_{Z^{(1)}}\,
 F^{-1}(\gbigi^{(0)}_{f_i})$
 are principal.
\hfill\qed
\end{itemize}
\end{lem}

Let $Z(F^{\ast}(f_i))\subset Z^{(1)}\setminus F^{-1}(H)$
be the set of the points $P\in Z^{(1)}\setminus F^{-1}(H)$
such that $F^{\ast}(f_i)(P)=0$.
We can reword the second condition
in Lemma \ref{lem;13.5.12.10}
as follows.

\begin{lem}
\label{lem;13.5.11.10}
For any $Q\in F^{-1}(H)$ and for any 
$i=1,\ldots,m$,
one of the following holds:
\begin{description}
\item[\bf (A1)]
$Q$ is not contained
in the closure of  $Z\bigl(F^{-1}(f_i)\bigr)$.
\item[\bf (A2)]
$F^{-1}(f_i)$ is a section of
$F^{\ast}\nbige$ around $Q$.
\end{description}
\end{lem}
\pf
Fix $Q\in F^{-1}(H)$ and $f_i$.
We set $P:=F(Q)$.
We take a frame $e_1,\ldots,e_r$ of $\nbige$.
We have the expression
$f_i=\sum_{p=1}^r A_p e_p$.
Let $h_j$ be defining functions of $D_j$ around $P$.
If $P\not\in D_j$, we put $h_j:=1$.
We set $\alpha:=\prod h_j^{N_j}$.
Then, the ideal $\gbigi_{f_i}$
is generated by
$\alpha A_{p}$ $(p=1,\ldots,r)$ and $\alpha$.
By the construction,
$\gbigi^{(1)}_{f_i}$ is principal.
If $F^{-1}\alpha$ is a generator,
we have that
$F^{-1}(A_p)\in \nbigo_{Z^{(1)},Q}$ for any $p$,
which means that {\bf (A2)} holds.
If one of $\alpha A_p$, say $\alpha A_{p_0}$, is a generator,
$A_{p_0}$ is invertible outside $F^{-1}(H)$,
which means that {\bf (A1)} holds.
Thus, Lemma \ref{lem;13.5.11.10} is proved.
\hfill\qed

\subsubsection{A general construction}
\label{subsection;13.5.12.21}

Let $W$ be any irreducible normal complex analytic space.
Let $D$ be an effective divisor on $W$.
The support of $D$ is denoted by $|D|$.
Let $\pi:E(D)\lrarr W$ be a holomorphic vector bundle on $W$.
Let $\Upsilon\subset E(D)$ 
be a reduced irreducible closed analytic 
subset of $E(D)$
such that
(i) the induced map $\pi:\Upsilon\lrarr W$ is finite and flat,
(ii) $\dim \Upsilon=\dim W$.
Let $d$ be the degree of $\Upsilon\lrarr W$.
We set
\[
 A(\Upsilon):=\Bigl\{
 (x_1,\ldots,x_d)\in
 \overbrace{
 \Upsilon\times_W\cdots\times_W\Upsilon}^d\,\Big|\,
 \pi(x_i)\not\in D,\,\,
 x_i\neq x_j\,\,(i\neq j)
 \Bigr\}.
\]
Let $\Abar(\Upsilon)$ denote the closure of $A(\Upsilon)$
in $E(D)\times_W\cdots\times_WE(D)$
with the reduced structure.
Let $Y(\Upsilon)$ denote the normalization of
$\Abar(\Upsilon)$.
Let $\pi_{Y(\Upsilon)}$ denote the naturally induced morphism
$Y(\Upsilon)\lrarr W$.

Let $\gbigs_d$ denote the $d$-th symmetric group.
We have a naturally induced action on
$\gbigs_d$ on $Y(\Upsilon)$.

\begin{lem}
The natural morphism 
$[\pi_{Y(\Upsilon)}]:
 Y(\Upsilon)/\gbigs_d\lrarr W$
is an isomorphism
of complex analytic spaces.
\end{lem}
\pf
Let $W_0(\Upsilon)$ denote the set of the points
$P\in W\setminus D$
such that the number of the points
in $\Upsilon_W\times\{P\}$ is $d$.
Then, 
$A(\Upsilon)$ is a $\gbigs_d$-principal bundle
over $W_0(\Upsilon)$.
We naturally have
$A(\Upsilon)/\gbigs_d\simeq
 W_0(\Upsilon)$.

Because 
$\Abar(\Upsilon)$ is contained in
$\Upsilon\times_W\cdots\times_W\Upsilon$,
it is finite over $W$.
Hence,
$\pi_{Y(\Upsilon)}$ is finite.
By the construction,
the dimension of each irreducible component
$Y(\Upsilon)$ is $\dim W$.
Hence, we obtain that $\pi_{Y(\Upsilon)}$ is an open map
(see \cite{Grauert-Remmert}).
Then, it is easy to check that
$[\pi_{Y(\Upsilon)}]$ gives a homeomorphism
of topological spaces.
Both of 
$Y(\Upsilon)/\gbigs_d$
and $W$ are normal.
They are isomorphic on the open
subset which is the complement
of a nowhere dense closed analytic subset.
Hence, we obtain that they are isomorphic
as complex analytic spaces.
\hfill\qed

\vspace{.1in}

Let $D_{Y(\Upsilon)}$ be the pull back of $D$
by $Y(\Upsilon)\lrarr W$.
The projection of
$E(D)\times_W\times\cdots\times_W E(D)$
onto the $i$-th component
induces a section
of $\pi^{\ast}_{Y(\Upsilon)}E(D)$ on $Y(\Upsilon)$,
which is denoted by $s_i$.
By applying the construction in Lemma \ref{lem;13.5.12.10}
to $\pi_{Y(\Upsilon)}^{\ast}E$ with 
the tuple of the sections 
$\bigl\{s_i\,\big|\,i=1,\ldots,m\bigr\}\cup
\bigl\{s_i-s_j\,\big|\,
 i,j=1,\ldots,m
 \bigr\}$,
we obtain a reduced irreducible normal complex analytic
space $Y_1$ with a $\gbigs_d$-action,
and an $\gbigs_d$-equivariant projective morphism 
$F_Y:Y_1\lrarr Y(\Upsilon)$.

We set $W_1:=Y_1/\gbigs_d$
which is an irreducible normal complex analytic space.
We have a naturally induced morphism
$F_W:W_1\lrarr W$.
It satisfies 
$W_1\setminus |D_1|
\simeq
 W\setminus|D|$.
We set $D_1:=F^{\ast}_W(D)$,
$E_1:=F_W^{\ast}E$,
and $\Upsilon_1:=F_W^{\ast}\Upsilon$ over $W_1$.
\begin{lem}
$Y_1$ is naturally isomorphic to
$Y(\Upsilon_1)$.
\end{lem}
\pf
By the construction,
we naturally have 
$Y_1\lrarr
 E_1(D_1)\times_{W_1}\cdots\times_{W_1}E_1(D_1)$.
The image is contained in
$\Abar(\Upsilon_1)$.
Hence, 
we have the naturally defined morphism
$\psi:
 Y_1\lrarr Y(\Upsilon_1)$.
The restriction of $\psi$ over
$W_1\setminus |D_1|$
is an isomorphism.
Because $Y_1\lrarr W_1$ is finite,
$\psi$ is also finite.
Both of $Y_1$ and $Y(\Upsilon)$ is normal.
Then, we obtain that $\psi$ is an isomorphism.
\hfill\qed

\vspace{.1in}

Let $Q$ be any point of $|D_{Y_1}|$.
For each $i$, one of the following holds:
 \begin{description}
 \item[$\pmb{(A_i1)}$]
 $Q$ is not contained in the closure of $Z(F_Y^{\ast}(s_i))$.
 \item[$\pmb{(A_i2)}$]
  $F_Y^{\ast}(s_i)$ are sections of 
 $\pi_{Y(\Upsilon_1)}^{\ast}E_1$.
 \end{description}
For each pair $(i,k)$, one of the following holds:
 \begin{description}
 \item[$\pmb{(A_{i,k}1)}$]
 $Q$ is not contained in the closure of $Z(F_Y^{\ast}(s_i-s_k))$.
 \item[$\pmb{(A_{i,k}2)}$]
  $F_Y^{\ast}(s_i-s_k)$ are sections of 
 $\pi_{Y(\Upsilon_1)}^{\ast}E_1$.
 \end{description}

\subsection{Proof of Theorem \ref{thm;13.5.11.1}}

\subsubsection{A criterion for goodness}

Let $X=\Delta^n$ and $H=\bigcup_{i=1}^{\ell}\{z_i=0\}$.
Let $\Sigma$ be a meromorphic Lagrangian cover
on $(X,H)$.
Let $\Sigmabar\subset 
 T^{\ast}X(\log H)\otimes\Omega(\vecN H)$
be the closure which is finite over $X$.
We obtain a normal complex space
$Y(\Sigmabar)$
with a morphism
$\pi_{Y(\Sigmabar)}:Y(\Sigmabar)
\lrarr X$
as in \S\ref{subsection;13.5.12.21}.
We also have the sections $s_i$ $(i=1,\ldots,d)$
of $\pi_{Y(\Sigmabar)}^{\ast}
 T^{\ast}X(\log D)\otimes\nbigo(\vecN H)$.
Let $O=(0,\ldots,0)$.

\begin{prop}
\label{prop;13.5.13.4}
Fix $Q\in\pi_{Y(\Sigmabar)}^{-1}(O)$.
Suppose that for any $i$,
at least one of the following holds:
\begin{description}
\item[$\pmb{(A_i1)}$]
 $Q$ is not contained in $F(s_i)$.
\item[$\pmb{(A_i2)}$]
 $s_i$ is a section of
 $\pi_{Y(\Sigmabar)}^{\ast} T^{\ast}X(\log D)$
\end{description}
Moreover, suppose that for any $(i,k)$,
\begin{description}
\item[$\pmb{(A_{ik}1)}$]
 $Q$ is not contained in $F(s_i)$.
\item[$\pmb{(A_{ik}2)}$]
 $s_i-s_k$ is a section of
 $\pi_{Y(\Sigmabar)}^{\ast} T^{\ast}X(\log D)$
\end{description}
Then, there exists a finite subset
$\nbigi_O\subset 
 \nbigotilde_{X}(\ast H)_P\big/
 \nbigotilde_{X,P}$
such that
(i) it satisfies the first two conditions in 
{\rm\S\ref{subsection;13.5.12.40}},
(ii) we have the decomposition 
{\rm(\ref{eq;13.5.13.1})} around $O$.
\end{prop}
\pf
By taking an appropriate ramified covering,
we have only to consider the case that
$Y(\Sigmabar)\lrarr X$ is not ramified
along $\{z_i=0\}$ $(i=1,\ldots,\ell)$.
We will shrink $X$ around $O$
without mention in the following argument.

We define an equivalence relation on 
$\{s_1,\ldots,s_d\}$ as follows:
$s_i\sim s_j$
if $s_i-s_j$ is a section of
$\pi_{Y(\Sigmabar)}^{\ast}T^{\ast}X(\log D)$
around $Q$.
We obtain a decomposition
by the equivalence relation:
\begin{equation}
 \label{eq;13.5.13.2}
 \{s_1,\ldots,s_d\}
=\bigsqcup_{\kappa\in T}
 B_{\kappa}
\end{equation}
Let $G_Q\subset\gbigs_d$ be the stabilizer of $Q$.
Then, the action of $G_Q$
on $\{s_1,\ldots,s_d\}$ preserves 
the decomposition (\ref{eq;13.5.13.2})
Indeed, if $g^{\ast}s_i=s_j$ for some $g\in G_Q$,
$Q$ is contained in the closure of $Z(s_i-s_j)$,
and hence $s_i,s_j\in B_{\kappa}$ for some $\kappa$.

For each $\kappa$,
we have a section $\omega_{\kappa}$ of
$\pi^{\ast}_{Y(\Sigmabar)}\Bigl(
 T^{\ast}X(\log D)\otimes
 \bigl(\nbigo_X(\vecN H)\big/\nbigo_X\bigr)
 \Bigr)$
induced by $s_i\in B_{\kappa}$,
which is independent of the choice of $s_i$.
By the consideration in the previous paragraph,
$\omega_{\kappa}$ is invariant
under the action of $G_Q$.
Hence, it is the pull back of a section
$\omega_{\kappa,0}$
of 
$T^{\ast}X(\log H)\otimes
 \nbigo_X(\vecN H)\big/\nbigo_X$.
Note that we have the induced exterior derivative
$d$
on $\Omega^{\bullet}_X(\log H)
 \otimes\nbigo_X(\vecN H)\big/\nbigo_X$.
Because 
we have $d\omega_{\kappa}=0$
around general points of $H$,
we have 
$d\omega_{\kappa}=0$.
Hence, we can take a meromorphic function
$\gminia_i\in M(X,H)$
such that $d\gminia_i$ induces $\omega_{\kappa,0}$.
We set $t_i:=s_i-d\gminia_{\kappa}$,
which are sections of
$\pi_{Y(\Sigmabar)}^{\ast}T^{\ast}X(\log D)$.
We obtain decompositions
$s_i=\pi_{Y(\Sigmabar)}^{\ast}d\gminia_{\kappa}+t_i$.

By the conditions 
$\pmb{(A_i1)}$ and $\pmb{(A_{ik1})}$,
we obtain that
$\ord\gminia_{\kappa}$ 
and $\ord(\gminia_{\kappa}-\gminia_{\kappa'})$
exist.
We have the Lagrangian covers
given by $\{t_i\,|\,s_i\in B_{\kappa}\}$,
denoted by $\Sigma_{\kappa}$,
and we have
$\Sigma=\bigsqcup \bigl(
 d\gminia_{\kappa}+\Sigma_{\kappa}
\bigr)$.
\hfill\qed

\subsubsection{Proof of Theorem \ref{thm;13.5.11.1}}

Let us return to the situation in
\S\ref{subsection;13.5.12.20}.
By applying  the construction in
\S\ref{subsection;13.5.12.21}
to $T^{\ast}X(\log H)\otimes\nbigo(\vecN H)$
with the closure $\Sigmabar$ of $\Sigma$,
and by applying the resolution of singularity,
we obtain a projective morphism $F:X_1\lrarr X$
of complex manifolds,
satisfying the following conditions:
\begin{itemize}
\item
 $H_1:=F^{-1}(H)$ is normal crossing.
\item
 We have $F^{\ast}\Sigmabar$
 in $F^{\ast}\bigl(
 T^{\ast}_X(\log H)
 \otimes\nbigo(\vecN H)
 \bigr)$. 
 For any $Q\in Y(F^{\ast}\Sigmabar)$
 and for any $i$ (resp. $(i,k)$)
 we have $\pmb{(A_i1)}$ or $\pmb{(A_i2)}$
 (resp.
 $\pmb{(A_{ik}1)}$ or $\pmb{(A_{ik}2)}$).
\end{itemize}
We obtain an effective divisor
$\vecN_1 H_1=F^{\ast}(\vec N H)$.
We have a naturally defined morphism
$F^{\ast}T^{\ast}X(\log H)\otimes\nbigo_X(\vecN H)
\lrarr
 T^{\ast}X_1(\log H_1)\otimes\nbigo_{X_1}(\vecN_1H_1)$.
Let $\Sigmabar_1$ denote the image of
$F^{\ast}\Sigmabar$.
Because
$X_1\setminus H_1\simeq
 X\setminus H$,
for any 
$Q\in Y(\Sigmabar_1)$,
 and for any $i$ (resp. $(i,k)$)
 we have $\pmb{(A_i1)}$ or $\pmb{(A_i2)}$
 (resp.  $\pmb{(A_{ik}1)}$ or $\pmb{(A_{ik}2)}$).
Then, by applying Proposition \ref{prop;13.5.13.4},
for each $P\in H_1$,
we obtain a finite subset
$\nbigi_P\subset
 \nbigotilde_{X_1}(\ast H_1)_P\big/
 \nbigotilde_{X_1,P}$
such that (i) it satisfies the first two conditions in 
{\rm\S\ref{subsection;13.5.12.40}},
(ii) we have the decomposition 
{\rm(\ref{eq;13.5.13.1})} around $P$.
It is easy to observe that the system
$\vecnbigi=\bigl\{
 \nbigi_P\,|\,P\in H_1
 \bigr\}$
satisfies the assumption in 
Proposition \ref{prop;13.5.7.110}.
Hence, we can take 
a projective morphism of complex manifolds
$F_1:X'\lrarr X_1$ such that
(i) $H':=F_1^{-1}(H_1)$ is simply normal crossing,
(ii) $X'\setminus H'\simeq X_1\setminus H_1$,
(iii) $F^{-1}\vecnbigi$ is a good set of irregular values
 on $(X',H')$.
Then, the induced projective morphism
$X'\lrarr X$ has the desired property.
Thus, the proof of Theorem \ref{thm;13.5.11.1}
is finished.
\hfill\qed

\backmatter%%%%%%%%%%%%%%%%%%%%%%%%%%%%%%%%%%%%%%%%%%%%%%%%%%%%%%%

\printindex

%%%%%%%%%%%%%%%%%%%%%%%%%%%%%%%%%%%%%%%%%%%%%%%%%%%%%%%%%%%%%%%%%%%%%%

\noindent
{\em Address\\
Research Institute for Mathematical Sciences,
Kyoto University,
Kyoto 606-8502, Japan,\\
takuro@kurims.kyoto-u.ac.jp

\end{document}